%% file: proba_02X_ang_arxiv.tex
\documentclass{amsbook}

\usepackage{amssymb}
\usepackage{amsfonts}
\usepackage{amsmath}
\usepackage{lastpage}
\usepackage[legalpaper,bookmarks=true,colorlinks=true,linkcolor=blue,citecolor=blue]{hyperref}
\usepackage{graphicx}%
\usepackage{fancyhdr}
\usepackage{color}
\usepackage[mathlines]{lineno}
\usepackage{lscape}
\usepackage{epsfig}
\usepackage{natbib}
\usepackage[]{listings}

\newtheorem{theorem}{Theorem}
\theoremstyle{plain}

\newtheorem{corollary}{Corollary}

\newtheorem{definition}{Definition}

\newtheorem{exercise}{Exercise}
\newtheorem{lemma}{Lemma}

\newtheorem{proposition}{Proposition}
\newtheorem{remark}{Remark}

\numberwithin{equation}{section}

\begin{document}
\Large
\pagenumbering{roman}
\begin{center}

\huge \textbf{Gane Samb LO}\\
\vskip 6cm
\Huge \textbf{Mathematical Foundations of Probability Theory} \\

\vskip 7cm

\huge \textit{\textbf{Statistics and Probability African Society (SPAS) Books Series}.\\
 \textbf{Saint-Louis, Calgary, Alberta. 2018}}.\\

\bigskip \Large  \textbf{DOI} : http://dx.doi.org/10.16929/sbs/2016.0008\\
\bigskip \textbf{ISBN} 978-2-9559183-8-8
\end{center}

\newpage

\begin{center}
\huge \textbf{SPAS TEXTBOOKS SERIES}
\end{center}

\bigskip \bigskip

\Large

 \begin{center}
 \textbf{GENERAL EDITOR of SPAS EDITIONS}
 \end{center}

\bigskip
\noindent \textbf{Prof Gane Samb LO}\\
gane-samb.lo@ugb.edu.sn, gslo@ugb.edu.ng\\
Gaston Berger University (UGB), Saint-Louis, SENEGAL.\\
African University of Sciences and Technology, AUST, Abuja, Nigeria.\\

\bigskip

\begin{center}
\Large \textbf{ASSOCIATED EDITORS}
\end{center}

\bigskip
\noindent \textbf{KEhinde Dahud SHANGODOYIN}\\
shangodoyink@mopipi.ub.bw\\
University of Botswana (Botswana)\\

\noindent \textbf{Blaise SOME}\\
some@univ-ouaga.bf\\
Chairman of LANIBIO, UFR/SEA\\
Joseph Ki-Zerbo University (Ouaga I), Burkina-Faso.\\

\bigskip
\begin{center}
\Large \textbf{ADVISORS}
\end{center}

\bigskip

\noindent \textbf{Ahmadou Bamba SOW}\\
ahmadou-bamba.sow@ugb.edu.sn\\
Gaston Berger University, Senegal.\\

\noindent \textbf{Tchilabalo Abozou KPANZOU}\\
kpanzout@yahoo.fr\\
Kara University, Togo.\\

\newpage

\Large \textbf{List of published books}\\

\Large \textbf{List of published or scheduled books in English}\\

\noindent $\square$ Weak Convergence (IA) - Sequences of Random Vectors. Gane Samb LO, Modou NGOM and Tchilabalo A. KPANZOU. 2016.\\
Doi : 10.16929/sbs/2016.0001. ISBN 978-2-9559183-1-9\\

\noindent $\square$ A Course on Elementary Probability Theory. Gane Samb LO. 2017.\\
Doi : 10.16929/sbs/2016.0003. ISBN 978-2-9559183-3-3\\

\noindent $\square$ Measure Theory and Integration By and For the Learner. Gane Samb LO.\\
Doi :  http://dx.doi.org/10.16929/sbs/2016.0005. ISBN  978-2-9559183-5-7\\

\newpage

\noindent \textbf{Library of Congress Cataloging-in-Publication Data}\\

\noindent Gane Samb LO, 1958-\\

\noindent Mathematical Foundations of Probability Theory.\\

\noindent SPAS Books Series, 2018.\\

\noindent \textit{DOI} : 10.16929/sbs/2016.0008\\

\noindent \textit{ISBN} 978-2-9559183-8-8\\

\newpage

\noindent \textbf{Author : Gane Samb LO}\\
\bigskip

\bigskip
\noindent \textbf{Emails}:\\
\noindent gane-samb.lo@ugb.edu.sn, ganesamblo@ganesamblo.net.\\

\bigskip
\noindent \textbf{Url's}:\\
\noindent www.ganesamblo@ganesamblo.net\\
\noindent www.statpas.net/cva.php?email.ganesamblo@yahoo.com.\\

\bigskip \noindent \textbf{Affiliations}.\\
Main affiliation : University Gaston Berger, UGB, SENEGAL.\\
African University of Sciences and Technology, AUST, ABuja, Nigeria.\\
Affiliated as a researcher to : LSTA, Pierre et Marie Curie University, Paris VI, France.\\

\noindent \textbf{Teaches or has taught} at the graduate level in the following universities:\\
Saint-Louis, Senegal (UGB)\\
Banjul, Gambia (TUG)\\
Bamako, Mali (USTTB)\\
Ouagadougou - Burkina Faso (UJK)\\
African Institute of Mathematical Sciences, Mbour, SENEGAL, AIMS.\\
Franceville, Gabon\\

\bigskip \noindent \textbf{Dedicatory}.\\

\noindent \textbf{To my first and tender assistants, my daughters Fatim Zahr\`a Lo and Maryam Majiguèn Azr\`a Lo}

\bigskip \noindent \textbf{Acknowledgment of Funding}.\\

\noindent The author acknowledges continuous support of the World Bank Excellence Center in Mathematics, Computer Sciences and Intelligence Technology, CEA-MITIC. His research projects in 2014, 2015 and 2016 are funded by the University of Gaston Berger in different forms and by CEA-MITIC.

\title{Mathematical Foundations of Probability Theory}

\begin{abstract} \large (\textbf{English}) In the footsteps of the book \textit{Measure Theory and Integration By and For the Learner} of our series in Probability Theory and Statistics, we intended to devote a special volume of the very probabilistic aspects of the first cited theory. The book might have assigned the title : From Measure Theory and Integration to Probability Theory. The fundamental aspects of Probability Theory, as described by the keywords and phrases below, are presented, not from experiences as in the book \textit{A Course on Elementary Probability Theory}, but from a pure mathematical view based on Measure Theory. Such an approach places Probability Theory in its natural frame of Functional Analysis and constitutes a firm preparation to the study of Random Analysis and Stochastic processes. At the same time, it offers a solid basis towards Mathematical Statistics Theory. The book will be continuously updated and improved on a yearly basis.

\noindent (\textbf{Fran\c{c}ais}) \noindent \\

\noindent \textbf{Keywords}. Measure Theory and Integration; Probabilistic Terminology of Measure Theory and Applications;
Probability Theory Axiomatic; Fundamental Properties of Probability Measures; Probability Laws of Random Vectors; Usual Probability Laws review; Gaussian Vectors; Probability Inequalities; Almost sure and in Probability Convergences; Weak convergences; Convergence in Lp; Kolmogorov Theory on sequences of independent real-valued random variables; Central Limit Theorem, Laws of Large Numbers, Berry-Essen Approximation, Law of the iterated logarithm for real valued independent random variables; Existence Theorem of Kolmogorov and Skorohod for Stochastic processes; Conditional Expectations; First examples of stochastic process : Brownian and Poisson Processes.\\  

\noindent \textbf{AMS 2010 Classification Subjects :} 60-01; 60-02; 60-03;G0Axx; 62GXX.
\end{abstract}
\maketitle

\frontmatter
\tableofcontents
\mainmatter
\Large
\include{preface_gen_en}

\include{proba_02_00_a_preface_ang}
\include{proba_02_00_intro_ang}

\include{proba_02_01_ang}

\include{proba_02_02_ang}

\include{proba_02_03_ang}

\include{proba_02_03_exo_ang}
\include{proba_02_04_ang}

\include{proba_02_05_ang}

\include{proba_02_06_ang}

\include{proba_02_07_ang}

\include{proba_02_07_kol_ang}

\include{proba_02_07_Gauss_series}

\include{proba_02_07_clt_ang}
\include{proba_02_07_be_ang}

\include{proba_02_07_LIL_ang}

\include{proba_02_08_ang}

\include{proba_02_09_ang}

\include{proba_02_09_proof_ang}

\include{proba_02_appendix}

\include{proba_02_biblio_ang}
\end{document}

%% file: preface_gen_en.tex
\chapter*{General Preface}

\noindent \textbf{This textbook} is one of the elements of a series whose ambition is to cover a broad part of Probability Theory and Statistics.  These textbooks are intended to help learners and readers, of all levels, to train themselves.\\

\noindent As well, they may constitute helpful documents for professors and teachers for both courses and exercises.  For more ambitious  people, they are only starting points towards more advanced and personalized books. So, these textbooks are kindly put at the disposal of professors and learners.

\bigskip \noindent \textbf{Our textbooks are classified into categories}.\\

\noindent \textbf{A series of introductory  books for beginners}. Books of this series are usually destined to students of first year in universities and to any individual wishing to have an initiation on the subject. They do not require advanced mathematics.  Books on elementary probability theory (See \cite{ips-probelem-ang}, for instance) and descriptive statistics are to be put in that category. Books of that kind are usually introductions to more advanced and mathematical versions of the same theory. Books of the first kind also prepare the applications of those of the second.\\

\noindent \textbf{A series of books oriented to applications}. Students or researchers in very related disciplines  such as Health studies, Hydrology, Finance, Economics, etc.  may be in need of Probability Theory or Statistics. They are not interested in these disciplines  by themselves.  Rather, they need to apply their findings as tools to solve their specific problems. So, adapted books on Probability Theory and Statistics may be composed to focus on the applications of such fields. A perfect example concerns the need of mathematical statistics for economists who do not necessarily have a good background in Measure Theory.\\

\noindent \textbf{A series of specialized books on Probability theory and Statistics of high level}. This series begins with a book on Measure Theory, a book on its probability theory version, and an introductory book on topology. On that basis, we will have, as much as possible,  a coherent presentation of branches of Probability theory and Statistics. We will try  to have a self-contained approach, as much as possible, so that anything we need will be in the series.\\

\noindent Finally, a series of \textbf{research monographs} closes this architecture. This architecture should be so diversified and deep that the readers of monograph booklets will find all needed theories and inputs in it.\\

\bigskip \noindent We conclude by saying that, with  only an undergraduate level, the reader will  open the door of anything in Probability theory and statistics with \textbf{Measure Theory and integration}. Once this course validated, eventually combined with two solid courses on topology and functional analysis, he will have all the means to get specialized in any branch in these disciplines.\\

\bigskip \noindent Our collaborators and former students are invited to make live this trend and to develop it so  that  the center of Saint-Louis becomes or continues to be a re-known mathematical school, especially in Probability Theory and Statistics.

%% file: proba_02_00_intro_ang.tex
\chapter*{Introduction} \label{proba_02_intro}

\noindent \textbf{Mathematical Foundation of Probability Theory}.\\

\noindent In the introduction to the book \textit{Measure Theory and Integration By and For The Learner}, we said : \\

\noindent \textit{Undoubtedly, Measure Theory and Integration is one of the most important part of Modern Analysis, with Topology and Functional Analysis for example. Indeed, Modern mathematics is based on functional analysis, which is a combination of the Theory of Measure and Integration, and Topology}.\\

\noindent \textit{The application of mathematics is very pronounced in many fields, such as finance (through stochastic calculus), mathematical economics (through stochastic calculus), econometrics [which is a contextualization of statistical regression to economic problems], physic statistics. Probability Theory and Statistics has become an important tool for the analysis of biological phenomena and genetics} modeling.\\

\noindent This quotation already stressed the important role played by Probability Theory in the application of Measure Theory. So, Probability Theory seems to be one of the most celebrated extensions of Measure Theory and Integration when it comes to apply it to real life problems.\\

\noindent Probability Theory itself may be presented as the result of modeling of stochastic phenomena based on random experiences. This way is illustrated in the element of this series : \textit{A Course on Elementary Probability Theory}.\\

\noindent But for theoretical purposes, it may be presented as a mathematical theory, mainly based on Measure Theory and Integration, Topology and Functional Analysis. This leads to impressive tools that reveal themselves very powerful in dealing real-life problems.\\

\noindent In this book, we tried to give the most common elements of the Theory as direct rephrasing and adaptation of results Measure Theory according to the following scenario.\\

\noindent Chapter \ref{chap1} is devoted to a complete rephrasing of the Measure Theory and Integration Terminology to that of Probability Theorem, moving from a general measures to normed measures called Probability Measures.\\

\noindent Chapters \ref{proba_02_rv}, \ref{proba_02_upl} and \label{proba_02_gauss} deal with a simple fact in Measure Theory and Integration, namely the image-measure, which becomes the most important notion in Probability Theory and called under the name of Probability Laws.  Chapter \ref{proba_02_rv} includes a wide range of characterizations for Probability Laws of Random vectors we might need in research problems in Probability Theory and Mathematical Statistics. In particular, the concept of independence is visited from various angles, which leads to a significant number of important characterizations of it. In Chapter \ref{proba_02_upl}, usual and important probability Laws are given and reviewed in this chapter in connection with the their generations described made in \cite{ips-probelem-ang}. Finally Chapter \ref{proba_02_upl} presents the so important Gaussian random vectors.\\

\noindent Chapter \ref{proba_02_conv} is concerned with the theory of convergence of sequences of (real-valued, mainly) random variables. The three types of Convergence : Almost-sure, in Probability and in $L^p$.\\

\noindent It is important to notice the the book \textit{Weak Convergence (IA). Sequences of random vectors} (See \cite{ips-wcia-ang}) has its place exactly here, within the global frame of the series. Due to its importance and its size, we preferred to devote a booklet of medium size (about two hundred pages) to an introduction to weak convergence.\\
 
\noindent Because of the importance of Inequalities in Probability Theory, we devote Chapter \ref{proba_02_ineg} to them. This chapter will
continuously updated and augmented on a yearly basis.\\

\noindent In Chapter \ref{probab_02_indep}, we presented the main results of the study of sequence of independent random variables which occupied the researchers in a great part of the 19th century. The laws that were studied are until now the most important ones of the theory, although they are exented to the non-independent cases nowadays. But there is no way to join the current studies if the classical main tools and proofs are not mastered. We introduce to the Kolmogorov Strong Law of Large numbers, the Central Limit Theorem, the Berry-Essen Approximation and the Law of the Iterated Logarithm.\\

\noindent Chapter \ref{proba_02_ce} uses the Radon-Nikodym Theorem to found the important notion of Mathematical Expectation which is the main tool form moving to independent to dependent data.\\

\noindent Finally, Chapter \ref{thfondamentalKolm} presents the Fundamental Theorem of Kolmogorov which is considered as the foundation of Modern Probability Theory. Versions of the Theorem are given, among them, the Skorohod Theorem. This chapter is the bridge with the course on Stochastic processes.\\

\noindent \textbf{The place of the book within the series}.\\

\noindent While the book \textit{A Course on Elementary Probability Theory} may  read at any level, the current one should no be read before the full exposition of Measure Theory and Integration (\cite{ips-mestuto-ang} or a similar book). Indeed, the latter book is cited in any couple of pages. The demonstrations in that book are quoted in the current one. Without assuming those demonstration, this textbook would have a very much greater number of pages.\\

\noindent Reading the textbook \cite{ips-wcia-ang} is recommended after Chapter \ref{proba_02_conv} of the current book.\\

\noindent Now, this book combined with \cite{ips-wcrv-ang} open the doors of many other projects of textbooks, among whom we cite :\\

\noindent (a) Asymptotics of Sequences of Random Vectors\\

\noindent (b) Stochastic Processes\\

\noindent (c) Mathematical Statistics\\

\noindent (d) Random Measures\\
 
\noindent (e) Times Series\\

\noindent (f) etc.\\
 
\noindent Consequently, the series will expand to those areas.\\


%% file: proba_02_01_ang.tex
\chapter[From Measure Theory to Probability Theory]{An update of the Terminology from Measure Theory to Probability Theory} \label{chap1}\label{proba_02_intotp}

\section{Introduction} \label{proba_02_intotp_01}

\noindent This course of Probability Theory is the natural continuation of the one on Measure Theory and Integration. Its constitutes the very minimal basis for a fundamental course which enables to prepare for more advanced courses on  Probability Theory and Statistics, like Stochastic Processes, Stochastic Calculus or to prepare specialized Mathematical statistics, etc.\\

\noindent The book \textit{A Course on Elementary Probability Theory} (\cite{ips-probelem-ang}) of this series concentrated on discrete probability measures and focused on random experiences, urn models, generation of random variables and associated computations. The reader will not find such results here. We recommend him to go back to this book or to similar ones which directly deal with  Probability Theory related to real experiences. This textbook treats the mathematical aspects of Probability Theory, as a branch of Measure Theory and Integration as exposed in \cite{{ips-mestuto-ang}}, where the Measure Theory terminology can be found.\\

\noindent This course begins with new expressions and names of concepts introduced in Measure Theory and Integration. Next, a specific orientation will be taken to present the base of modern Probability Theory.\\

\noindent
\section{Probabilistic Terminology} \label{proba_02_intotp_02}

\subsection{Probability space} $ $\\

\noindent A probability space is a measure space $(\Omega,\mathcal{A},m)$ where the measure
assigns the unity value to the whole space $\Omega$, that is, 
\begin{equation*}
m(\Omega )=1.
\end{equation*}

\bigskip \noindent Such a measure is called a probability measure. Probability measures are generally denoted in blackboard font : $\mathbb{P}$, $\mathbb{Q}$, etc.\\

\bigskip \noindent We begin with this definition : \\

\noindent

\begin{definition}
Let  $(\Omega ,\mathcal{A})$ be a measurable space. The mapping
\begin{equation*}
\begin{array}{cccc}
\mathbb{P}: & \mathcal{A} & \rightarrow & \mathbb{R} \\ 
& A & \hookrightarrow & \mathbb{P}(A)%
\end{array}%
\end{equation*}

\noindent is a probability measure if and only if $\mathbb{P}$ is a measure and 
$\mathbb{P}(\Omega)=1$, that is :\\

\noindent (a)  $0\leq \mathbb{P}\leq \mathbb{P}(\Omega )=1$.\\

\noindent (b) For any countable collection of measurable sets $\{A_{n},n\geq 0\}\subset \mathcal{A}$, pairwise disjoints, we have
 
\begin{equation*}
\mathbb{P}(\sum_{n\geq 0}A_{n})=\sum_{n\geq 0}\mathbb{P}(A_{n}).
\end{equation*}
\end{definition}

\bigskip

\noindent We adopt a special terminology in Probability Theory.\\

\bigskip  \noindent (1) The whole space $\Omega$ is called \textit{universe}.\\

\noindent (2) Measurable sets are called \textit{events}. Singletons are elementary events whenever they are measurable.\\

\noindent \textbf{Example}. Let us consider a random experience in which we toss two dies and get the outcomes as the ordered pairs $(i,j)$, where $i$ and $j$ are respectively the number of the first and next the second face of the two dies which come out.  Here, the universe is $\Omega=\{1,...,6 \}^2$. An ordered pair $\{(i,j)\}$ is an elementary event. As an other example, the event : \textit{the sum of the faces is less or equal to 3} is exactly  
$$
A=\{(1,1),(1,2),(2,1)\}.
$$

\bigskip \noindent (3) \textit{Contrary event}. Since $\mathbb{P}(\Omega )=1$, the probability of the complement
of an event $A$, also called the contrary event to $A$ and denoted $\overline{A}$, is computed as

\begin{equation*}
\mathbb{P}(\overline{A})=1-\mathbb{P}(A).
\end{equation*}

\bigskip \noindent The previous facts form simple transitions from Measure Theory and Integration terminology to that of Probability Theory. We are going to continue to do the same in more elaborated transitions in the rest of that chapter.\\

\subsection{Properties of a Probability measure} $ $\\

\noindent Probability measures inherit all the properties of a measure.\\

\noindent (P1) A probability measure is sub-additive, that is, for any countable collection of events $\{A_{n},n\geq 0\}\subset \mathcal{A}$, we have
 
\begin{equation*}
\mathbb{P}(\bigcup_{n\geq 0}A_{n}) \leq \sum_{n\geq 0}\mathbb{P}(A_{n}).
\end{equation*}  

\noindent (P2) A probability measure $\mathbb{P}$ is non-decreasing, that is, for any ordered pair of events $(A,B) \in \mathcal{A}^2$ such that $A \subset B$, we have

$$
\mathbb{P}(A) \leq \mathbb{P}(B)
$$

\bigskip \noindent and more generally for any ordered pair of events $(A,B) \in \mathcal{A}^2$, we have

$$
\mathbb{P}(B \setminus A) = \mathbb{P}(B)- \mathbb{P}(A \cap B).
$$

\bigskip \noindent (P3) A probability measure $\mathbb{P}$ is continuous below, that is, for any non-decreasing sequence of events $(A_{n})_{n\geq 0} \subset \mathcal{A}$, we have
$$
\mathbb{P}(\bigcup_{n\geq 0} A_n)=\lim_{n\rightarrow +\infty} \mathbb{P}(A_n),
$$

\bigskip \noindent and is continuous above, that is, for any non-increasing sequence of events $(A_{n})_{n\geq 0} \subset \mathcal{A}$, we have
$$
\mathbb{P}(\bigcap_{n\geq 0} A_n)=\lim_{n\rightarrow +\infty} \mathbb{P}(A_n)
$$

\noindent The \textsl{continuity above} in Measure Theory requires that the values of the measures of the $A_n$'s be finite for at least one integer $n\geq 0$. Here, we do not have to worry about this, since all $\mathbb{P}(A_n)$'s are bounded by one.\\

\subsection{Random variables} $ $\\

\noindent Measurable mappings are called random variables. Hence, a mapping 
\begin{equation}
\begin{array}{cccc}
X: & (\Omega ,\mathcal{A)} & \rightarrow & (E,\mathcal{B})%
\end{array}
\label{0110}
\end{equation}

\noindent is a random variable, with respect to the $\sigma$-algebras $\mathcal{A}$ and $\mathcal{B}$ if and only if it is measurable with respect to the same $\sigma$-algebras.\\

\bigskip \noindent \textbf{Probability law}.\\

\noindent There is not a more important phrase in Probability Theory that \textit{Probability law}. I dare say that the essence of probability Theory is finding probability laws of random phenomena by intellectual means and the essence of Statistical theory is the same but by means of inference from observations or data.\\

\bigskip \noindent Suppose that we have a probability measure $\mathbb{P}$ on the measurable space $(\Omega,\mathcal{A})$ in Formula \ref{0110}. We have the following definition.

\begin{definition} \label{02_def_01} The Probability law of the random variable $X$ in Formula \ref{0110} is the image-measure of $\mathbb{P}$ by $X$, denoted as $\mathbb{P}_X$, which is a probability measure on $E$ given by

$$
\mathcal{B} \ni \mathcal{B} \mapsto \mathbb{P}_X(B)=\mathbb{P}(X \in B). \  \Diamond 
$$
\end{definition}

\bigskip \noindent Such a simple object holds everything in Probability Theory.\\

\bigskip \noindent \textbf{Classification of random variables}.\\

\noindent Although the space $E$ in Formula (\ref{0110}) is arbitrary, the following cases are usually and commonly studied :\\

\noindent (a) If $E$ is $\overline{\mathbb{R}}$, endowed with the usual Borel $\sigma$-algebra, the random variable is called a \textit{real random variables} (\textit{rrv}).\\

\noindent (b) If $E$ is $\overline{\mathbb{R}}^{d}$ $(d\in \mathbb{N}^{\ast })$, endowed with the usual Borel $\sigma$-algebra, $X$ is called a $d$-random vector or a random vector of dimension $d$, denoted $X=(X_{1},X_{2},...,X_{d})^t$, where $X^t$ stands for the transpose of $X$.\\

\noindent (c) More generally if $E$ is of the form $\mathbb{R}^{T}$, where $T$ is a non-empty set, finite or countable infinite or non-countable infinite, $X$ is simply called a stochastic process. The $\sigma$-algebra on $\mathbb{R}^{T}$, which is considered as the collections of mapping from $T$ to $\mathbb{R}$ is constructed by using the fundamental theorem of Kolmogorov, which generally is stated in the first chapter of a course on Stochastic processes, and which is extensively stated in Chapter \ref{thfondamentalKolm}.\\

\noindent In the special case where $T=\mathbb{N}$, $X$ is a sequence of real random variables $X=\{X_1, X_2, ....\}$.\\

\noindent (d) If $E$ is some metric space $(S,d)$ endowed with the Borel $\sigma$-algebra denoted as $\mathcal{B}(S)$, the term \textit{random variable} is simply used although some authors prefer using \textit{random element}.\\

\subsection{Mathematical Expectation} $ $ \\

\noindent It is very important to notice that, at the basic level, the mathematical expectation, and later the conditional mathematical expectation, is defined for a real random variable.\\

\noindent \textbf{(a) Mathematical expectation of \textbf{rrvs}'s}.\\

\noindent Let $X \ : \ (\Omega ,\mathcal{A}, \mathbb{P}) \rightarrow (\overline{\mathbb{R}}, \mathcal{B}_{\infty}(\overline{\mathbb{R}})$ be a \textit{real random variable}. Its mathematical expectation with respect to the probability measure $\mathbb{P}$ or its $\mathbb{P}$-mathematical expectation, denoted by $\mathbb{E}_{\mathbb{P}}(X)$ is simply its integral with respect to $$\mathbb{P}$$ whenever it exists and we denote :

$$
\mathbb{E}_{\mathbb{P}}(X)=\int_{\Omega} X \ d\mathbb{P}.
$$

\bigskip \noindent The full notation $\mathbb{E}_{\mathbb{P}}$ of the mathematical expectation reminds us to which probability measure the mathematical expectation is relative to. In many examples, it may be clear that all the mathematical expectations are relative to only one probability measure so that we may drop the subscript and only write

$$
\mathbb{E}(X)=\int_{\Omega} X \ d\mathbb{P}.
$$

\noindent Also, the parentheses may also be removed and we write $\mathbb{E}X$.\\

\noindent \textbf{(b) Mathematical expectation of a function of an arbitrary random variable}.\\

\noindent For an arbitrary random variable as defined in Formula (\ref{0110}) and for any real-valued  measurable mapping

\begin{equation}
\begin{array}{cccc}
h: & (E,\mathcal{B)} & \rightarrow & (\overline{\mathbb{R}}, \mathcal{B}_{\infty}(\overline{\mathbb{R}}),
\end{array}
\end{equation}

\bigskip \noindent the composite mapping
$$
h(X)=h\circ X:  (\Omega ,\mathcal{A,}\mathbb{P}\mathcal{)}  \rightarrow \overline{\mathbb{R}}
$$

\bigskip \noindent is a real random variable. We may define the mathematical expectation of $h(X)$ with respect to $\mathbb{P}$ by 

$$
\mathbb{E} h(X) =\mathbb\int_{\Omega} h(X) \ d\mathbb{P},
$$

\bigskip \noindent  whenever the integral exists.\\

\noindent \textbf{(c) Use of the probability law for computing the mathematical expectation}.\\

\noindent We already know from the properties of image-measures (See page \cite{ips-mestuto-ang}, Doc 04-01 Point (V)), that we may compute the mathematical expectation of $h(X)$, if it exists, by

\begin{equation}
\mathbb{E}(h(X))=\int_{\Omega}h \ d\mathbb{P}_X=\int_{\Omega} h(x) \ d\mathbb{P}_X(x). \ \label{FIG}
\end{equation}

\bigskip \noindent If $X$ is itself a real random variable, its expectation, if it exists, is 

\begin{equation}
\mathbb{E}(X)=\int_{\mathbb{R}} x \ d\mathbb{P}_X(x). \ \label{FIGS}
\end{equation}

\bigskip \noindent (d) \textbf{Mathematical expectation of a vector}.\\

\noindent The notion of mathematical expectation may be extended to random vectors by considering the vector of the mathematical expectations of the coordinates. Let us consider the random vector $X$ such that $X^t=(X_{1},X_{2},...,X_{d})$. The Mathematical vector expectation $\mathbb{E}(X)$ is defined by 

$$
(\mathbb{E}(X))^t=(\mathbb{E}X_{1}, \mathbb{E}X_{2},..., \mathbb{E}X_{d}).
$$

\bigskip \noindent A similar extension can be operated for random matrices.\\

\noindent \textbf{(e) Properties of the Mathematical expectation}.\\

\noindent As an integral of real-valued measurable application, the mathematical expectation inherits all the properties of integrals we already had in Measure Theory. Here, we have to add that :
\textit{constant real random variables and bounded random variables have finite expectations}. The most important legacy to highlight is the following.\\

\begin{theorem} \label{proba_02_expect} On the class of all random variables with finite mathematical expectation denoted $\mathcal{L}^1(\Omega, \mathcal{A}, \mathbb{P})$, the mathematical expectation operator :\\

\noindent (a) is linear, that is for all $(\alpha,\beta)\in \mathbb{R}^2$, for all $(X,Y)\in \mathcal{L}^1(\Omega, \mathcal{A}, \mathbb{P})$,

$$
\mathbb{E}(\alpha X + \beta Y)=\alpha \mathbb{E}(X) + \beta \mathbb{E}(Y),
$$

\bigskip \noindent  (b) is non-negative, that for all non-negative $X$ random variable, we have $\mathbb{E}(X)\geq 0$\\

\bigskip \noindent  (c) and satisfies for all non-negative $X$ random variable : $\mathbb{E}(X)=0$ if and only if $X=0$, $\mathbb{P}$-\textit{a.e}.

\bigskip \noindent  (c) Besides, we have for all real-valued random variables $X$ and $Y$ defined on $(\Omega, \mathcal{A})$,

$$
\biggr( |X| \leq Y, \ Y \in \mathcal{L}^1(\Omega, \mathcal{A}, \mathbb{P})\biggr) \Rightarrow X \in \mathcal{L}^1(\Omega, \mathcal{A}, \mathbb{P}). 
$$

\bigskip \noindent and

$$
\biggr| \int X \ d\mathbb{P} \biggr| \leq \int |X| \ d\mathbb{P} \leq \int Y \ d\mathbb{P}. 
$$
\end{theorem}

\bigskip \noindent The first formula in the following Lemma is often used to computing the mathematical expectation of non-negative real-valued randoms variables. We generalize with respect to the counting measure. For example, this will render much comprehensible the proof the Kolmogorov Theorem \ref{sllnEPKol} (Chapter \ref{probab_02_indep}, page \pageref{sllnEPKol}) on strong laws of large numbers.\\

\noindent Let us define, for a real-valued random variable $X$, its lower endpoint $lep(X)$ and the upper endpoint $uep(X)$ respectively by

$$
lep(X)=\inf \{t \in \mathbb{R}, \ \mathbb{P}(X\leq t) >0 \}, \  uep(X)=\sup \{t \in \mathbb{R}, \ \mathbb{P}(X\leq t) < 1 \}.
$$ 

\bigskip \noindent This means that $\mathbb{P}(X\leq t)=0$ for all $t>uep(X)$ and similarly, we have $\mathbb{P}(X\leq t)=0$  for all $t\leq lep(X)$. Actually, we have $uep(X)=\|X\|_{\infty}$ in the $L^{\infty}$ space. The values space of $X$ becomes $\mathcal{V}_X=[lep(X), \ uep(X)]$.\\

\noindent We have : 

\begin{proposition} \label{proba_02_sec_01_esp} Let $X$ be any real-valued and non negative random variable, we have \label{formulaCF}
$$
\ \ \mathbb{E}(X)=\int_{0}^{uep(X)} \mathbb{P}(X>t) dt, \ \ (CF)
$$

\bigskip \noindent and

$$
1 + \mathbb{E}\left([X]_{+}\right)=\sum_{n \in [0, \ [uep(X)]_{+}]} \mathbb{P}(X\geq n), \ \ (DF1)
$$

\noindent where $[x]_{+}$ (resp. $[x]^{+}$) stands for the greatest (resp. smallest) integer less or equal (resp. greater or equal) to $x \in \overline{\mathbb{R}}$. Also, for any \textit{a.s} finite real-valued random variable, we have

$$
-1 + \sum_{n \in [0, \ [uep(X)]_{+}]} \mathbb{P}(|X|\geq n) \leq \mathbb{E}|X|\leq \sum_{n \in [0, \ [uep(X)]_{+}]} \mathbb{P}(|X|\geq n). \ \ (DF2)
$$
\end{proposition}

\bigskip \noindent Since $\mathbb{P}(X>t)=0$ for all $t>uep(X)$, extending the integration domain to $+\infty$ does not effect the value of the integral.\\

\bigskip \noindent \textbf{Proof}. \\

\noindent \textbf{Proof of (CF)}. The function $t \mapsto \mathbb{P}(X>t)$ is bounded and has at most a countable number of discontinuity. So its improper Riemann integral is a Lebesgue's one and we may apply the Tonelli's Theorem (See Chapter 8, Doc 07-01 in \cite{ips-mestuto-ang}) at Line (L13) below as follows :

\begin{eqnarray*}
\int_{0}^{uep(X)} \mathbb{P}(X>t) dt &=&\int_{0}^{uep(X)} \mathbb{P}(X>t) d\lambda(t)\\
&=&  \int_{(t\in ]0, uep(X)])} \biggr( \int_{\mathcal{V}_X} 1_{(x>t)} \ d\mathbb{P}_X(x) \biggr) d\lambda(t)\\
&=&  \int_{(t\in ]0, uep(X)])}  1_{(x>t)} \ d\mathbb{P}_X(x) d\lambda(t)\\
&=& \int_{\mathcal{V}_X} \biggr( \int_{(t\in ]0, uep(X)])}  1_{(x>t)} \biggr) d\mathbb{P}_X(x) \ \ (L13)\\
&=& \int_{\mathcal{V}_X} \biggr( \int_{0}^{\min(x,uep(X)} d\lambda(t) \biggr) d\mathbb{P}_X(x) \\
&=& \int_{\mathcal{V}_X} \min(x,uep(X) d\mathbb{P}_X(x) \\
&=& \mathbb{E}(min(X,uep(X)))=\mathbb{E}(X),
\end{eqnarray*}

\noindent  since $X \leq uep(X)$ \textit{a.s.}\\

\noindent \textbf{Proof of (DF1)}. We use the counting measure $\nu$ on $\mathbb{N}$ and say

\begin{eqnarray*}
\sum_{n\geq 0, \ n\leq [uep(X)]^{+}} \mathbb{P}(X\geq n)&=&  \int_{[0,\ [uep(X)]^{+}]} \mathbb{P}(X\geq n) \ d\nu(n)\\
&=&\int_{[0,\ [uep(X)]^{+}]} \biggr(\int_{\mathcal{V}_X} 1_{(x\geq n)} \ d\mathbb{P}_X(x)\biggr) d\nu(n) \ \ (L22)\\
&=&\i\int_{\mathcal{V}_X} \biggr( \int_{[0,\ [uep(X)]^{+}]} 1_{(x\geq n)} d\nu(n)\biggr) \ d\mathbb{P}_X(x) \ \ (L23)\\
&=&\int_{\mathcal{V}_X} \nu([0, \max(x,[uep(X)]^{+})]) \ d\mathbb{P}_X(x) \ (L24)\\
&=&\int_{0}^{+\infty} \biggr(\biggr[\max(x,[uep(X)]^{+})\biggr]_{+}+1\biggr) \ d\mathbb{P}_X(x)\\
&=&\mathbb{E}\biggr(\biggr[\max(X,[uep(X)]^{+})\biggr]_{+}\biggr)+1.\\
\end{eqnarray*}

\noindent We conclude that 
\begin{eqnarray*}
\sum_{n\geq 0, \ n\leq [uep(X)]^{+}} \mathbb{P}(X\geq n)&=&\mathbb{E}[X]_{+}+1.\\
\end{eqnarray*}


\bigskip \noindent  \textbf{Proof of (DF2)}. The left-hand inequality is derived from (DF1) when applied to non-negative random variabe $|X|$. To establish the right-hand inequality, suppose that $X$ is non-negative. Let us denote $A_n=(X\geq n)$, $n\geq 0$ with $A_0=\Omega$ clearly. We have for any $n\geq 1$, $A_{n-1}\setminus A_n)=(n-1\leq X <n)$. If $uep(X)=+\infty$, the sets $]n-1,n]$, $n\geq 1$, form a partition of $\mathbb{R}_+$. If  $uep(X)$ if finite, the sets $]n-1,n]$, $1\leq n \leq N=[uep(X)]_{+}+1$ for a partition of $[0, N]$ which covers $X$ \textit{a.s.}. So we have  

$$
\sum_{1\leq n < N+1} \biggr(A_{n-1}\setminus A_n\biggr) =\Omega.
$$

\bigskip \noindent  By the Monotone Convergence Theorem when $N=+\infty$, but by finite additivity for $N$ finite, we have

\begin{eqnarray*}
\mathbb{E}(X)&=& \mathbb{E}\biggr(X \sum_{1\leq n < N+1} 1_{A_{n-1}\setminus A_n} \biggr) \ \ (L31)\\
&=& \sum_{1\leq n < N+1} \mathbb{E}\biggr(X  1_{A_{n-1}\setminus A_n} \biggr)\\
&=& \sum_{1\leq n < N+1} \mathbb{E}\biggr(X  1_{(n-1\geq |X| <n)} \biggr)\\
&\leq & \sum_{1\leq n < N+1} n \mathbb{E}\biggr(1_{(n-1\geq X <n)} \biggr)\\
&\leq & \sum_{1\leq n < N+1} n \biggr(\mathbb{P}(A_{n-1})-\mathbb{P}(A_{n})\biggr) \ \  (L35)\\
\end{eqnarray*}

\bigskip \noindent  Suppose that $N$ is infinite. By developing the last line, we have for 

\begin{eqnarray*}
\sum_{1\leq n \leq k+1} n (\mathbb{P}(A_{n-1})-\mathbb{P}(A_{n}))&=& \sum_{0\leq n \leq k} \mathbb{P}(A_n) - (k+1) \mathbb{P}(A_{k+1})\\
&\leq & \sum_{0\leq n \leq k} \mathbb{P}(A_n). \ \ (L42)
\end{eqnarray*}

\bigskip \noindent  By letting $k \rightarrow +\infty$ in Line (42) and by combining the results with Lines (L31) and (L36), we get the inequality.\\

\noindent If $N$ is finite, the last line is exactly

$$
\sum_{0\leq n \leq N} \mathbb{P}(A_n) - N \mathbb{P}(A_{N} \leq \sum_{0\leq n \leq N} \mathbb{P}(A_n),
$$

\bigskip \noindent  and hence, is less or equal to $\sum_{0\leq n \leq [uep(X)]_{+}} \mathbb{P}(A_n)$ and we conclude that

$$
\mathbb{E}(X) \leq \sum_{0\leq n \leq [uep(X)]_{+}} \mathbb{P}(A_n). 
$$

\bigskip \noindent  To get the right-hand inequality in (DF2), we just apply the last formula to $|X|$. $\blacksquare$\\

\bigskip \noindent  \textbf{An easy example}. Suppose that $X$ is Bernoulli random variable with $\mathbb{P}(X=1)=1-\mathbb{P}(X=0)=p$, $0<p<1$. We have
$\mathbb{E}(X)=p$, $uep(X)=1$, $[uep(X)]_{+}=1$, $\mathbb{P}(A_0)=1$, $\mathbb{P}(A_1)=p$ and we exactly have

$$
\sum_{0\leq n \leq [uep(X)]_{+}} \mathbb{P}(A_n)=1+p=\mathbb{E}(X)+1.
$$

\bigskip
\subsection{Almost-sure events} $ $\\

\noindent In Measure Theory, we have studied null-sets and the notion of almost-everywhere  \textit{(a.e)} properties. In the context of probability theory, for any null-set $N$, we have
\begin{equation*}
\mathbb{P}(\overline{N})=1-\mathbb{P}(N)=1.
\end{equation*}

\bigskip \noindent So, the complement of any null-set is an almost-sure (\textit{a.s.}) event. Then, a random property $\mathbb{P}$ holds \textit{a.s.} if and only if

\begin{equation*}
\mathbb{P}(\{\omega \in \Omega ,\text{ }\mathbb{P}(\omega )\text{ true } \})=1.
\end{equation*}

\bigskip \noindent An almost-everywhere (\textsl{a.e.}) property is simply called an almost-sure (\textsl{a.s.}) property. Let us recall some properties of \textit{a.e.} properties.\\

\noindent (P1) If A is an \textit{a.e.} event and if the event $B$ is contained $A$, then $B$ is an \textit{a.s.} event.\\

\noindent (P2) A countable union of \textit{a.s.} events is an \textit{a.s.} event.\\

\noindent (P3) If each assertion of a countable family of assertions holds \textit{a.s.}, then all the assertions of the family hold simultaneously \textit{a.s.}\\

\subsection{Convergences of real-valued random variables} $ $\\

\noindent We may also rephrase the convergence results in Measure Theory as follows.\\

\noindent \textbf{(a) Almost-sure convergence}.\\

\noindent Let  $X$ and $X_{n}$, $n\geq 0$, be random variables defined on the probability space $(\Omega ,\mathcal{A},\mathbb{P})$. The sequence $(X_{n})_{n\geq 0}$ converges to $X$ almost-surely, denoted

$$
X_{n}\rightarrow X, \ a.s.,
$$

\bigskip
\noindent as $n \rightarrow +\infty$ if and only if $(X_{n})_{n\geq 0}$ converges to $X$ \textit{a.e.}, that is

$$
\mathbb{P}(\{\omega \in \Omega, \ X_n(\omega) \rightarrow X(\omega)\})=1.
$$

\bigskip \noindent \textbf{(b) Convergence in Probability}.\\

\noindent The convergence in measure becomes the convergence in Probability. Let $X$ be an \textit{a.s.}-finite real random variable and $(X_{n})_{n\geq 0}$ be a sequence of \textit{a.e.}-finite real random variables. We say that $(X_{n})_{n\geq 0}$ converges to $X$ in probability, denoted

$$
X_n \overset{\mathbb{P}}{\longrightarrow} X,
$$

\bigskip \noindent if and only if, for any $\varepsilon>0$,

$$
\mathbb{R}(|X_n-X|>\varepsilon) \rightarrow 0 \ as \ n\rightarrow +\infty.
$$

\bigskip \noindent We remind that the \textit{a.s.} limit and the limit in probability are \textit{a.s.} unique.\\

\bigskip \noindent We re-conduct the comparison result from Measure Theory.\\

\noindent \textbf{(c) Comparison between these two Convergence types}.\\

\noindent Let $X$ be an \textit{a.s.}-finite real random variable and $(X_{n})_{n\geq 0}$ be a sequence of \textit{a.e.} finite real random variables. Then we have the following implications, where all the unspecified limits are done as $n\rightarrow +\infty$ .\\
 
\noindent (1) If $X_{n} \rightarrow X$, \textit{a.s.}, then $X_n \overset{\mathbb{P}}{\longrightarrow} X$.\\

\noindent (2) If $X_n \overset{\mathbb{P}}{\longrightarrow} X$, then there exists a sub-sequence $(X_{n_k})_{k\geq 0}$ of $(X_{n})_{n\geq 0}$ such that $X_{n_{k}}\rightarrow X$, $a.s.$, as $k\rightarrow +\infty$.\\

\noindent Later, we will complete the comparison theorem by adding the convergence in the space

$$
L^{p}=\{X \in L^0, \  \mathbb{E}\left\vert X\right\vert^{p}<\infty\}, \ p\geq 1.
$$

\bigskip \noindent The weak convergence also will be quoted here while its study is done in \cite{ips-wcrv-ang}.

\section{Independence} \label{proba_02_intotp_03} 

\noindent The notion of independence is extremely important in Probability Theory and its applications. The main reason that the theory, in its earlier stages, has been hugely developed in the frame of independent random variables. Besides, a considerable number methods of handling dependent random variables are still generalizations of techniques used in the independence frame. In some dependence studies, it is possible to express the dependence from known functions of independent objects. In others, approximations based on how the dependence is near the independence are used.\\

\bigskip \noindent So, mastering the independence notion and related techniques is very important. In the elementary book (\cite{ips-probelem-ang}), we introduced the independence of events in the following way :\\

\noindent \textbf{Definition}. Let $A_{1}$, $A_{2}$, . . . , $A_{n}$ be events in a probability space $(\Omega , \mathbb{P}(\Omega), \mathbb{P})$. We have the following definitions.\\

\noindent (A)  The events $A_{1}$, $A_{2}$, . . . , $A_{n-1}$ and $A_{n}$ are pairwise independent if and only if 
\begin{equation*}
\mathbb{P}(A_{i}\cap A_{j})=\mathbb{P}(A_{i})\text{ }\mathbb{P}(A_{j}),\text{
for all }1\leq i\neq j\leq n.
\end{equation*}

\noindent (B)  The events $A_{1}$, $A_{2}$, . . ., $A_{n-1}$ and $A_{n}$ are mutually independent if and only if for any subset $\left\{
i_{1,}i_{2,}...,i_{k}\right\} $ of  $\left\{ 1,2,...,n\right\}$, with $2\leq k\leq n$, we have  
\begin{equation*}
\mathbb{P}(A_{i_{1}}\cap A_{i_{2}}\cap ...\cap A_{i_{k}})=\mathbb{P}(A_{i_{1}})\text{ }%
\mathbb{P}(A_{i_{2}}) ... \mathbb{P}(A_{i_{k}})\text{.}
\end{equation*}

\bigskip \noindent (C)  Finally, the events $A_{1}$, $A_{2}$, . . ., $A_{n-1}$ and $A_{n}$ fulfills the global factorization formula if and only if  
\begin{equation*}
\mathbb{P}(A_{1}\cap A_{2}\cap ...\cap A_{n})=\mathbb{P}(A_{1})\text{ }%
\mathbb{P}(A_{2}) ... \mathbb{P}(A_{n})\text{.}
\end{equation*}

\noindent We showed with examples that none two definitions from the three definitions (A), (B) and (C) are equivalence. It is important to know that, without any further specification, \textit{independence}  refers to Definition (B).\\

\noindent Measure Theory and Integration (MTI) \textit{gives the nicest and most perfect way} to deal with the notion of independence and, by the way, with the notion of dependence with copulas.\\

\subsection{Independence of random variables}$ $\\

\noindent Let $X_{1},...,X_{n}$ be $n$ random variables defined on the same probability space
\begin{equation*}
\begin{array}{cccc}
X_{i} & (\Omega ,\mathcal{A},\mathbb{P}) & \mapsto & (E_{i},\mathcal{B}_{i}).%
\end{array}%
\end{equation*}

\bigskip \noindent Let $(X_{1},...,X_{n})$ be a the $n$-tuple defined by 
\begin{equation*}
\begin{array}{cccc}
(X_{1},...,X_{n})^t : & (\Omega ,\mathcal{A}) & \mapsto & (E,\mathcal{B})%
\end{array}%
\end{equation*}

\noindent where  $E=\Pi _{1\leq i\leq n}E_{i}$ is the product space of the $E_{i}$'s endowed with the product $\sigma$-algebra, $\mathcal{B}=\otimes _{1\leq i\leq n}\mathcal{B}_{i}$. On 
each $(E_{i},\mathcal{B}_{i})$, we have the probability law $\mathbb{P}_{X_{i}}$ of $X_{i}$.\\

\noindent Each of the $\mathbb{P}_{X_{i}}$'s is called a marginal probability law of $(X_1,...,X_n)^t$.\\

\bigskip \noindent  On $(E,\mathcal{B})$, we have the following product probability measure 
\begin{equation*}
\mathbb{P}_{X_{1}}\otimes ...\otimes \mathbb{P}_{X_{n}},
\end{equation*}

\noindent characterized on the semi-algebra 
$$
S=\{\Pi _{1\leq i\leq n} \text{ } A_{i},A_{i}\in \mathcal{B}_{i}\}
$$

\bigskip \noindent of measurable rectangles by 

\begin{equation}
\mathbb{P}_{X_{1}}\otimes ...\otimes \mathbb{P}_{X_{n}}\left(\prod_{1\leq i\leq
n}A_{i})=\prod_{1\leq i\leq n}\mathbb{P}_{X_{i}}(A_{i}\right).  \label{0115}
\end{equation}

\noindent Now, we have two probability measures

$$
\mathbb{P}_{X_{1}}\otimes ...\otimes \mathbb{P}_{X_{n}}
$$ 

\noindent that is the product probability measure of the marginal probability measures and the probability law  
\begin{equation*}
\mathbb{P}_{(X_{1},...,X_{n})}(B)=\mathbb{P}((X_{1},...,X_{n})\in B).
\end{equation*}

\bigskip \noindent of the $n$-tuple $(X_{1},...,X_{n})$ on $(E,\mathcal{B})$, with is the image-measure of $\mathbb{P}$ by $(X_{1},...,X_{n})$. The latter probability measure is called the joint probability measure.\\

\bigskip \noindent By the $\lambda$-$\pi$ Lemma (See \cite{ips-mestuto-ang}, Exercise 11 of Doc 04-02, Part VI, page 228), these two probability measures are equal whenever they agree on the 
semi-algebra $\mathcal{S}$.\\

\noindent Now, we may give the most general definition of the independence of random variables :

\begin{definition}\label{proba_02_01_defIndependence}
The random variables $X_{1},...,$ and $X_{n}$ are independent if and only if the joint probability law  $\mathbb{P}_{(X_{1},...,X_{n})}$ of the
vector $(X_{1},...,X_{n})$ is the product measure of its marginal probability laws $\mathbb{P}_{X_{i}}$, that is :\\

\noindent For any  $B_{i}\in \mathcal{B}_{i},$ $1\leq i\leq n$, 
\begin{equation}
\mathbb{P}(X_{1}\in B_{1},X_{2}\in B_{2},...,X_{n}\in B_{n})=\prod_{1\leq
i\leq n}\mathbb{P}_{X_{i}}(B_{i}).  \label{0117}
\end{equation}
\end{definition}

\bigskip

\noindent For an ordered pair of random variables, the two random variables 
\begin{equation*}
{\Large 
\begin{array}{cccc}
X: & (\Omega ,\mathcal{A}) & \mapsto & (E,\mathcal{B})%
\end{array}
}
\end{equation*}

\noindent and

\begin{equation*}
{\Large 
\begin{array}{cccc}
Y & (\Omega ,\mathcal{A}) & \mapsto & (F,\mathcal{G})%
\end{array}%
}
\end{equation*}

\bigskip \noindent are independent if and only if $A\in \mathcal{B}$ et B$\in \mathcal{G}$,
\begin{equation*}
\mathbb{P}(X\in B,Y\in G)=\mathbb{P}(X\in A)\times \mathbb{P}(Y\in B).
\end{equation*}

\bigskip \noindent \textbf{Important Remark}. The independence is defined for random variables defined on the same probability space.
The space in which they take values may differ.\\

\bigskip \noindent Formula (\ref{0117}) may be rephrased by means of measurable functions in place of measurable subsets. We have

\begin{theorem} \label{proba_02_intotp_th02}
The random variables $X_{1},...,$ and $X_{n}$ are independent if and only if, for all non-negative and measurable real-valued functions $h_{i}:(E_{i},\mathcal{B}_{i})\mapsto \mathbb{R}$, we have

\begin{equation}
\mathbb{E}\biggr( \prod_{1\leq i\leq n}h_{i}(X_{i}) \biggr)=\prod_{1\leq i\leq n}\mathbb{E}%
(h_{i}(X_{i})).  \label{0118}
\end{equation}
\end{theorem}

\bigskip \noindent \textbf{Proof}.\\

\noindent We have to show the equivalence between Formulas (\ref{0117}) and (\ref{0118}). Let us begin to suppose that Formula (\ref{0118}) holds. Let us prove Formula  (\ref{0117}). Let $A_{i}\in \mathcal{B}$ and set $h_{i}=1_{A_{i}}$. Each $h_{i}$ is non-negative and measurable. Further 
\begin{equation*}
h_{i}(X_{i})=1_{A_{i}}(X)=1_{(X_{i}\in A_{i})}.
\end{equation*}

\bigskip \noindent and then
\begin{equation}
\mathbb{E}(h_{i}(X_{i}))=\mathbb{E}(1_{(X_{i}\in A_{i})})=\mathbb{P}%
(X_{i}\in A_{i}).  \label{0120}
\end{equation}

\bigskip \noindent As well, we have
\begin{equation*}
\prod_{1\leq i\leq n}h_{i}(X_{i})=\prod_{1\leq i\leq n}1_{(X_{i}\in
A_{i})}=1_{(X_{1}\in A_{1},...,X_{n}\in A_{n})}
\end{equation*}

\bigskip \noindent and then
\begin{eqnarray}
\mathbb{E}(\prod_{1\leq i\leq n}h_{i}(X_{i}))&=&\mathbb{E}(1_{(X_{1}\in A_{1},...,X_{n}\in A_{n})}) \label{0121}\\
&=&\mathbb{P}(X_{1}\in A_{1},...,X_{n}\in A_{n}). \notag
\end{eqnarray}

\bigskip \noindent By putting together (\ref{0120}) and (\ref{0121}), we get (\ref{0117}).\\

\bigskip \noindent Now, assume that (\ref{0117}) holds. Let  $h_{i}:(E_{i},\mathcal{B}_{i})\mapsto 
\mathbb{R}$ be measurable functions. Set 
\begin{equation*}
\mathbb{E}\biggr(\prod_{1\leq i\leq n}h_{i}(X_{i})\biggr)=\mathbb{E}(h(X_{1},...,X_{n})),
\end{equation*}

\noindent where $h(x_{1},...,x_{n})=h_{1}(x_{1})h_{2}(x_{2})...h_{n}(x_{n})$. The equality between the joint probability law and the product margin probability measures leads to
\begin{equation*}
\mathbb{E}(h(X_{1},...,X_{n}))=\int h(x_{1},...,x_{n})\text{ }d\mathbb{P}%
_{(X_{1},...,X_{n})}(h(x_{1},...,x_{n})
\end{equation*}

\begin{equation*}
=\int h(x_{1},...,x_{n})\text{ }d\left\{ \mathbb{P}_{X_{1}}\otimes
...\otimes \mathbb{P}_{X_{n}}\right\} (x_{1},...,x_{n}).
\end{equation*}

\bigskip \noindent From there, we apply Fubini's theorem, 
\begin{eqnarray*}
&&\mathbb{E}(h(X_{1},...,X_{n}))\\
&=&\int_{\Omega _{1}}d\mathbb{P}_{X_{1}}(x_{1})\int_{\Omega _{2}}d\mathbb{P}_{X_{2}}(x_{2})\int ...\\
&...& d\mathbb{P}_{X_{n-1}}(x_{n-1})\int h(x_{1},...,x_{n})\text{ }d\mathbb{P}_{X_{n}}(x_{1})\\
&=&\int_{\Omega _{1}}d\mathbb{P}_{X_{1}}(x_{1})\int_{\Omega _{2}}d\mathbb{P}_{X_{2}}(x_{2})\int ...d\mathbb{P}_{X_{n-1}}(x_{n-1})\int
h_{1}(x_{1})h_{2}(x_{2})...h_{n}(x_{n})\text{ }d\mathbb{P}_{X_{n}}(x_{1})\\
&=&\int_{\Omega _{1}}h_{1}(x_{1})\text{ }d\mathbb{P}_{X_{1}}(x_{1})\int_{\Omega _{2}}h_{2}(x_{2})\text{ }d\mathbb{P}_{X_{2}}(x_{2})...\int_{\Omega_{n}}.h_{n}(x_{n})\text{ }d\mathbb{P}_{X_{n}}(x_{1})\\
&=&\prod_{1\leq i\leq n}\int_{\Omega _{i}}h_{i}(x_{i})\text{ }d\mathbb{P}_{X_{i}}(x_{i})=\prod_{1\leq i\leq n}\int_{\Omega _{i}}h_{i}(X_{i})\text{ }d\mathbb{P}\\
&=&\prod_{1\leq i\leq n}\mathbb{E}(h_{i}(X_{i})). \ \blacksquare
\end{eqnarray*}

\bigskip \noindent This demonstration says that we have independence if and only if Formula (\ref{0118}) holds for all measurable functions $h_{i}:(E_{i},%
\mathcal{B}_{i})\mapsto \mathbb{R}$, $\mathbb{P}_{X_{i}}$-integrable or simply for all measurable and bounded functions $h_{i}:(E_{i},\mathcal{B}_{i})\mapsto \mathbb{R}$ or for all
non-negative measurable functions $h_{i}$, $i\in \{1,...,n\}$.\\

\noindent Let us come back to independence of events.\\

\subsection{Independence of events} \label{proba_02_ssec_indep_ev}$ $\\

\noindent Independence of events is obtained from independence of random variables.\\

\noindent \textbf{(a) Simple case of two events}.\\

\noindent We say that two events  $A\in \mathcal{A}$ and $B\in \mathcal{A}$ are independent if and only if the random variables  $1_{A}$ and $1_{B}$ are independent, that is, for all 
$h_{i}:\mathbb{R}\rightarrow \mathbb{R}$ (i=1,2) non-negative and measurable  
\begin{equation}
\mathbb{E}h_{1}(1_{A})h_{2}(1_{B})=\mathbb{E}h_{1}(1_{A})\mathbb{E}%
h_{2}(1_{B})  \label{0131}
\end{equation}

\noindent \noindent As a direct consequence, we have for $h_{i}(x)=x$, that Formula \ref{0131} implies that  
\begin{equation*}
\mathbb{E(}1_{A}1_{B})=\mathbb{E(}1_{AB})=\mathbb{E}(1_{A})\mathbb{E}(1_{B}),
\end{equation*}

\noindent that is 
\begin{equation}
\mathbb{P}(AB)=\mathbb{P(}A)\times \mathbb{P}(B).  \label{0132}
\end{equation}

\bigskip \noindent Now, we are going to prove that (\ref{0132}), in its turn, implies (\ref{0131}). First, let us show
that (\ref{0132}) implies 
\begin{equation}
\mathbb{P(}A^{c}B)=\mathbb{P(}A^{c})\times \mathbb{P}(B),  \label{0132a}
\end{equation}

\begin{equation}
\mathbb{P(}AB^{c})=\mathbb{P(}A)\times \mathbb{P}(B^{c})  \label{0132b}
\end{equation}

\bigskip \noindent and\\

\begin{equation}
\mathbb{P(}A^{c}B^{c})=\mathbb{P(}A^{c})\times \mathbb{P}(B^{c}).
\label{0132c}
\end{equation}

\bigskip \noindent Assume that (\ref{0132}). Since, 
\begin{equation*}
B=AB+A^{c}B,
\end{equation*}

\noindent we have
\begin{equation*}
\mathbb{P}(B)=\mathbb{P}(AB)+\mathbb{P}(A^{c}B)=\mathbb{P}(A)\mathbb{P}(B)+%
\mathbb{P}(A^{c}B).
\end{equation*}

\noindent Then 
\begin{equation*}
\mathbb{P}(A^{c}B)=\mathbb{P}(B)-\mathbb{P}(A)\mathbb{P}(B)=\mathbb{P}(B)(1-%
\mathbb{P}(A))=\mathbb{P}(A^{c})\mathbb{P}(B).
\end{equation*}

\noindent Hence (\ref{0132a}) holds. And (\ref{0132b}) is derived in the same manner by exchanging the role of $A$ and $B$. Now, to prove (\ref%
{0132c}), remark that
\begin{equation*}
A^{c}B^{c}=(A\cup B)^{c}=(AB^{c}+A^{c}B+AB)^{c}.
\end{equation*}

\noindent Then, we get

\begin{eqnarray*}
\mathbb{P}(A^{c}B^{c})&=&1-\mathbb{P}(AB^{c})-\mathbb{P}(A^{c}B)-\mathbb{P}(AB)\\
&=&1-\mathbb{P}(A)\mathbb{P}(B^{c})-\mathbb{P}(A^{c})\mathbb{P}(B)-\mathbb{P}(A)\mathbb{P}(B)\\
&=&1-\mathbb{P}(A)-\mathbb{P}(A^{c})\mathbb{P}(B)\\
&=&1-\mathbb{P}(A)-\mathbb{P}(A^{c})(1-\mathbb{P}(B^{c}))\\
&=&1-\mathbb{P}(A)-\mathbb{P}(A^{c})+\mathbb{P}(A^{c})\mathbb{P}(B^{c})\\
&=&\mathbb{P}(A^{c})\mathbb{P}(B^{c}).
\end{eqnarray*}

\bigskip \noindent Hence (\ref{0132c}) holds.\\

\noindent Finally, let us show that Formula (\ref{0132}) ensures Formula  (\ref{0131}). Consider two non-negative and measurable mappings $h_{i}:\mathbb{R}\rightarrow \mathbb{R}$,  (i=1,2). We have
 
\begin{equation*}
h_{1}(1_{A})=h_{1}(1)1_{A}+h_{1}(0)1_{A^{c}}
\end{equation*}

\noindent and

\begin{equation*}
h_{2}(1_{B})=h_{2}(1)1_{B}+h_{2}(0)1_{B^{c}}.
\end{equation*}

\noindent As well, we have

\begin{eqnarray*}
h_{1}(1_{A})h_{2}(1_{B})&=&h_{1}(1)h_{2}(1)1_{AB}+h_{1}(1)h_{2}(0)1_{AB^{c}}\\
&+&h_{1}(0)h_{2}(1)1_{A^{c}B}+h_{1}(0)h_{2}(0)1_{A^{c}B^{c}}.
\end{eqnarray*}

\noindent Then, we have

\begin{equation*}
\mathbb{E(}h_{1}(1_{A}))=h_{1}(1)\mathbb{P}(A)+h_{1}(0)\mathbb{P}(A^{c})
\end{equation*}

\bigskip \noindent and\\

\begin{equation*}
\mathbb{E(}h_{2}(1_{B}))=h_{2}(1)\mathbb{P}(B)+h_{2}(0)\mathbb{P}(B^{c}).
\end{equation*}

\bigskip \noindent We also have

\begin{equation*}
\mathbb{E}h_{1}(1_{A})h_{2}(1_{B})=h_{1}(1)h_{2}(1)\mathbb{P}(A)\mathbb{P}%
(B)+h_{1}(1)h_{2}(0)\mathbb{P}(A)\mathbb{P}(B^{c})
\end{equation*}%
\ 
\begin{equation*}
+h_{1}(0)h_{2}(1)\mathbb{P}(A^{c})\mathbb{P}(B)+h_{1}(0)h_{2}(0)\mathbb{P}%
(A^{c})\mathbb{P}(B^{c}).
\end{equation*}

\bigskip \noindent By comparing the three last formulas, we indeed obtain that
\begin{equation*}
\mathbb{E}h_{1}(1_{A})h_{2}(1_{B})=\mathbb{E}h_{1}(1_{A})\mathbb{E}%
h_{2}(1_{B}).
\end{equation*}

\bigskip \noindent The previous developments lead to the definition (and theorem).

\begin{definition} (Definition-Theorem). The events $A$ and $B$ are independent if and only if $1_{A}$
and $1_{B}$ are independent if and only if 
\begin{equation}
\mathbb{P=}(AB)=\mathbb{P}(A)\times \mathbb{P}(B).  \label{0133}
\end{equation}
\end{definition}

\bigskip \noindent \textbf{(b) Case of an arbitrary finite number $k\geq 2$ of events}.\\

\noindent Let us extend this definition to an arbitrary number $k$ of events and compare it with the definition (B) in the preliminary remarks of this section.\\

\noindent Let $A_{i}$, $1\leq i\leq k,$ be $k$ events and $h_{i}:\mathbb{R} \rightarrow \mathbb{R}$, $(i=1,...,k)$, be $k$ non-negative and measurable mappings. The events $A_{i}$ are independent if and only if the mappings
$1_{A_{i}}$, $1\leq i \leq k$, are independent if and only if for all measurable finite mappings $h_i$, $1\leq i \leq k$, we have 
\begin{equation}
\mathbb{E(}\prod_{1\leq i\leq k}h_{i}(1_{A_{i}}))=\prod_{1\leq i\leq k}\mathbb{E(}h_{i}(1_{A_{i}})).  \label{0134}
\end{equation}

\bigskip \noindent Let us put for each $s$-tuple of non-negative integers $1\leq i_{1}\leq i_{2}\leq ...\leq i_{s}\leq k$, $1 \leq s \leq  k$,
\begin{equation*}
h_{i_{j}}(x)=x,\text{ }j=1,...,s
\end{equation*}

\bigskip \noindent and

\begin{equation*}
h_{i}(x)=1 \ for \ i\notin \{i_{1},i_{2},...,i_{s}\}.
\end{equation*}

\bigskip \noindent Hence, by Formula (\ref{0134}), we get for any subset $\left\{i_{1,}i_{2,}...,i_{k}\right\} $ of  $\left\{ 1,2,...,n\right\}$, with $2\leq k\leq n$ 
\begin{equation}
\mathbb{P}(\bigcap_{1\leq j\leq s}A_{i_{j}})=\prod_{1\leq j\leq s}\mathbb{P}
(A_{i_{j}}).  \label{0135}
\end{equation}

\bigskip \noindent This is Definition (B) in the preliminary remarks of this section. By the way, it is also a generalization of Formula (\ref{0132}) for two events ensembles. We may, here again, use  straightforward computations similar to those done for the case $k\geq 2$, to show that Formula \ref{0135} also implies Formula \ref{0134}. This leads to the definition below.

\begin{definition} (Definition-Theorem) The events $A_{i},$ $1\leq i\leq k,$ are independent if and only if the mappings $1_{A_{i}}$ are independent if and only if for each $s$-tuple $1\leq i_{1}\leq i_{2}\leq ...\leq i_{s}\leq k$, of non-negative integers, 
\begin{equation}
\mathbb{P}\left(\bigcap_{1\leq j\leq s}A_{i_{j}}\right)=\prod_{1\leq j\leq s}\mathbb{P}%
(A_{i_{j}}).
\end{equation}
\end{definition}

\bigskip \noindent \textbf{(c) An interesting remark}.\\

\noindent A useful by-product of Formula (\ref{0134}) is that if $\{A_i, \ 1\leq i \leq n\}$, is a collection of independent events, then any elements of any collection of events $\{B_i, \ 1\leq i \leq n\}$, with $B_i=A_i$ or $B_i=A_i^c$, are also independent.\\

\noindent To see this, it is enough to establish Formula (B). But for any $\left\{i_{1,}i_{2,}...,i_{k}\right\} $ of  $\left\{ 1,2,...,n\right\}$, with $2\leq k\leq n$, we make take
$h_{i_j}(x)=x$ if $B_{i_j}=A_{i_j}$ or $h_{i_j}(x)=1-x$ if $B_{i_j}=A_{i_j}^c$ for $j=1,...,k$ and $h_i(x)=1$ for $i \notin \{i_1,...,i_k\}$ in Formula \ref{0134} and use the independence of the $A_i$'s.\\

\noindent We get, for $\left\{i_{1,}i_{2,}...,i_{k}\right\} \subset \left\{ 1,2,...,n\right\}$, with $2\leq k\leq n$, 
that
$$ 
\mathbb{P}(\bigcap_{1\leq j\leq s}B_{i_{j}})=\prod_{1\leq j\leq s}\mathbb{P}(B_{i_{j}}). \ \square
$$

\bigskip 

\subsection{Transformation of independent random variables} $ $\\

\noindent 
Consider the independent random variables 
$$
\begin{array}{cccc}
X_{i}: & (\Omega ,\mathcal{A},\mathbb{P}) & \mapsto & (E_{i},\mathcal{B}_{i}),
\end{array}
$$

\bigskip \noindent $i=1,...,n$ and  $g_{i}:(E_{i},\mathcal{B}_{i})\mapsto (F_{i},F_{i})$, $n$ measurable mappings.\\

\bigskip \noindent Then, the random variables $g_{i}(X_{i})$ are also independent.\\

\noindent Indeed, if $h_{i}:F_{i}\rightarrow \mathbb{R}$, $1\leq i \leq n$, are measurable and bounded real-valued mappings, then the $h_{i}(g_{i})$ are also real-valued bounded and measurable mappings. Hence, the $h_{i}(g_{i}(X_{i}))$'s are $\mathbb{P}$-integrable. By independence of the $X_{i}$, we get
 
\begin{equation*}
\mathbb{E}(\prod_{1\leq i\leq n}h_{i}\circ g_{i}(X_{i}))=\prod_{1\leq i\leq
n}\mathbb{E}(h_{i}\circ g_{i}(X_{i})),
\end{equation*}

\noindent \noindent and this proves the independence of the $h_{i}\circ g_{i}(X_{i})$'s. We have the proposition :\\

\begin{proposition}
Measurable transformations of independent random variables are independent
\end{proposition}

\subsection{Family of independent random variables}. $ $\\

\noindent Consider a family of random variables

$$
\begin{array}{cccc}
X_{t} & (\Omega ,\mathcal{A},\mathbb{P}) & \mapsto & (E_{t},\mathcal{B}_{t}), \ (t\in T).
\end{array}%
$$

\bigskip \noindent This family $\{X_{t},t\in T\}$ may be finite, infinite and countable or infinite and non countable. It is said that the random variables of this family are independent if and only the random variables in any finite sub-family of the family are independent, that is, for any subfamily $\{t_{1},t_{2},...,t_{p}\}\subset T$, $2\leq p <+\infty$, the mappings $X_{t_{1}},X_{t_{2}},...,X_{t_{p}}$ are independent.\\

\noindent The coherence of this definition will be a consequence of the Kolmgorov Theorem.\\

\section{Pointcarr\'e and Bonferroni Formulas} \label{proba_02_intotp_04}

\bigskip
\noindent \textbf{Poincarr\'e or Inclusion-exclusion Formula}.\\

\noindent In \cite{ips-mestuto-ang}, we already proved these following formulas for subsets $A_1$, ..., $A_n$ of $\Omega$, $n\geq 2$ :\\

\begin{eqnarray}
&&1_{ \bigcup_{(1 \leq j \leq j}A_{j})}  \notag \\
&=&\sum_{1\leq j\leq n} 1_{A_{j}} +\sum_{r=2}^{n}(-1)^{r+1} \sum_{1\leq i_{1}<...<t_{r}\leq n}1_{A_{i_{1}}...A_{i_r}} \ \ \label{hp01}
\end{eqnarray}

\bigskip \noindent and

\begin{equation}
Card(\bigcup_{1\leq j\leq n}A_j)=\sum_{1\leq j\leq n}Card(A_{j})+\sum_{r=2}^{n}(-1)^{r+1} \sum_{1\leq i_{1}<...<t_{r}\leq n} Card(A_{i_{1}}...A_{ir}). \ \label{hp02}
\end{equation}

\bigskip \noindent In the cited book, Formula (\ref{hp01}) is proved and very similar techniques may be repeated to have Formula (\ref{hp02}). The same techniques also lead the formula

\begin{equation}
\mathbb{P}(\bigcup_{1\leq j\leq n}A_j)=\sum_{1\leq j\leq n}\mathbb{P}(A_{j})+\sum_{r=2}^{n}(-1)^{r+1} \sum_{1\leq i_{1}<...<t_{r}\leq n} \mathbb{P}(A_{i_{1}}...A_{ir}),  \label{hp03}
\end{equation}

\noindent if $A_1$, ..., $A_n$ are events.\\

\noindent These three formula are different versions of the Pointcarr\'e's Formula, also called Inclusion-Exclusion Formula.\\ 

\noindent \textbf{Bonferroni's Inequality}.\\

\noindent Let $A_1$, ..., $A_n$ be measurable subsets of $\Omega$, $n\geq 2$. Define

\begin{eqnarray*}
\alpha_0&=&\sum_{1\leq j\leq n}\mathbb{P}(A_{j})\\
\alpha_1&=&\alpha_0 - \sum_{1\leq i_{1}<t_{2}\leq n} \mathbb{P}(A_{i_{1}}A_{i_2})\\
\alpha_2&=&\alpha_1 +  \sum_{1\leq i_{1}<...<t_{3}\leq n} \mathbb{P}(A_{i_{1}}...A_{i_3})\\
\cdots &=& \cdots\\
\alpha_r&=&\alpha_{r-1} + (-1)^{r+1} \sum_{1\leq i_{1}<...<t_{r}\leq n} \mathbb{P}(A_{i_{1}}...A_{i_r})\\
\cdots &=& \cdots\\
\alpha_r&=&\alpha_{r-1} + (-1)^{n+1}  \mathbb{P}(A_1 A_2 A_3 ... A_{n})
\end{eqnarray*}

\bigskip \noindent Let $p = n \text{ mod } 2$, that is $n=2p+1+h$, $h\in \{0,1\}$. We have the Bonferroni's inequalities : if $n$ is odd,

$$
\alpha_{2k+1}\leq \mathbb{P}\biggr( \bigcup_{1\leq j\leq n}A_n\biggr) \leq \alpha_{2k}, \ k=0,...,p  \ (BF1)
$$

\bigskip \noindent and if $n$ is even,

$$
\alpha_{2k+1}\leq \mathbb{P}\biggr( \bigcup_{1\leq j\leq n}A_j \biggr) \leq \alpha_{2k}, \ k=0,...,p-1. \ \ (BF2)
$$

\bigskip \noindent We may easily extend this formula to cardinalities in the following way. Suppose the $A_i$'s are finite subsets of $\Omega$ and one of them at least is non-empty. Denote by $M$ the
cardinality of $\Omega_0=\bigcup_{1\leq j\leq n}A_n$. Hence

$$
\mathbb{P}(\Omega_0) \ni A \mapsto \mathbb{P}(A)=\frac{1}{M} Card (A),
$$ 

\noindent is a probability measure and the Bonferroni inequalities hold. By multiplying the formulas by $M$, we get

$$
\beta_{2k+1}\leq Card \biggr( \bigcup_{1\leq j\leq n}A_n \biggr) \leq \alpha_{2k}, \ p=0,1,...
$$

\noindent where the sequence $(\beta_s)_{0\leq s \leq n}$ is defined sequentially by

$$
\beta_0=\sum_{1\leq j\leq n} Card(A_{j})
$$

\bigskip \noindent and for $r>0$, 

$$
\beta_r=\alpha_{r-1} + (-1)^{r+1} \sum_{1\leq i_{1}<...<t_{r}\leq n} Card(A_{i_{1}}...A_{i_r}).
$$

\bigskip \noindent The extension has been made in the case where one of $A_i$'s is non-empty. To finish, we remark that all inequalities hold as equalities of null terms if all the sets $A_i$'s are empty.\\

\bigskip \noindent Remark also that for $0<s\leq n$, we have

$$
\alpha_s=\sum_{1\leq j\leq n}\mathbb{P}(A_{j})+\sum_{r=2}^{s}(-1)^{r+1} \sum_{1\leq i_{1}<...<t_{r}\leq n} \mathbb{P}(A_{i_{1}}...A_{ir})
$$

\bigskip \noindent and

$$
\beta_s=\sum_{1\leq j\leq n} Card(A_{j})+\sum_{r=2}^{s}(-1)^{r+1} \sum_{1\leq i_{1}<...<t_{r}\leq n} Card(A_{i_{1}}...A_{ir}).
$$

%% file: proba_02_02_ang.tex
\chapter{Random Variables in $\mathbb{R}^d$, $d\geq 1$} \label{proba_02_rv} 

\noindent This chapter will focus on the basic important results of Probability Theory concerning random vectors. Most of the properties exposed
here and relative to discrete \textit{real} random variables are already given and proved in the textbook \cite{ips-probelem-ang} of this series. The new features
are the extensions of those results to vectors and the treatment of the whole thing as applications of the contents of Measure Theory and integration.\\

\noindent Three important results of Measure Theory and Integration, namely $L^p$ spaces, Lebesgue-stieljes measures and Radon-Nokodym's Theorem are extensively used.\\

\noindent First, we will begin with specific results for \textit{real} random variables.\\

\section{A review of Important Results for Real Random variables} \label{proba_02_rv_sec_01}
\noindent First, let us recall inequalities already established in Measure Theory. Next, we will introduce the new and important Jensen's one and give some of its applications.\\

\noindent \textbf{Remarkable inequalities}.\\

\noindent The first three inequalities are results of Measure Theory and Integration (See Chapter 10 in \cite{ips-mestuto-ang}).\\

\noindent \textbf{(1) H\"{o}lder Inequality}. Let $p>1$ and  $q>1$ be two conjugated positive rel numbers, that is, $1/p+1/q=1$ and let  

\begin{equation*}
\begin{array}{cccc}
X,Y: & (\Omega ,\mathcal{A},\mathbb{P}) & \mapsto & \mathbb{R}
\end{array}
\end{equation*}

\bigskip \noindent  be two random variables $X \in L^{p}$ and $Y \in L^{q}$. Then $XY$ is integrable and we have
 
\begin{equation*}
\left\vert \mathbb{E}(XY)\right\vert \leq \left\Vert X\right\Vert _{p}\times
\left\Vert Y\right\Vert _{q},
\end{equation*}

\bigskip \noindent  where for each $p\geq 1$, $\left\Vert X\right\Vert _{p}=(\inf |X|^p)^{1/p}$.\\

\bigskip \noindent  \textbf{(2) Cauchy-Schwartz's Inequality}. For $p=q=2$, the H\"{o}lder inequality becomes the Cauchy-Schwartz one :
\begin{equation*}
\left\vert \mathbb{E}(XY)\right\vert \leq \left\Vert X\right\Vert _{2}\times
\left\Vert Y\right\Vert _{2}.
\end{equation*}

\bigskip \noindent  \textbf{(3) Minskowski's Inequality}. Let $p\geq 1$ (including $p=+\infty$). If $X$ and $Y$ are in $L^{p}$, then we have 
\begin{equation*}
\left\Vert X+Y\right\Vert _{p}\leq \left\Vert X\right\Vert _{p}+\left\Vert
Y\right\Vert _{p}.
\end{equation*}

\bigskip \noindent  \textbf{(4) $C_p$  Inequality}. Let $p \in [1,+\infty[$. If $X$ and $Y$ are in $L^{p}$, then for $C_p=2^{p-1}$, we have 
\begin{equation*}
\|X+Y\|_p^p\leq C_p (\|X\|_p^p+\|Y\|_p^p).
\end{equation*}

\bigskip \noindent  \textbf{(5) Jensen's Inequality}.\\

\noindent \textbf{(a) Statement and proof of the inequality}.\\

\begin{proposition} \label{proba-02.jensen} (Jensen's inequality). Let $\phi $ be a convex function defined from a closed interval $I$ of $\mathbb{R}$ to $\mathbb{R}$. Let $X$ be a \textit{rrv} with values in $I$ such that $\mathbb{E}(X)$ is finite. Then $\mathbb{E}(X)\in I$ and 
\begin{equation*}
\phi (\mathbb{E}(X))\leq \mathbb{E}(\phi (X)).
\end{equation*}
\end{proposition}

\bigskip \noindent  \textbf{Proof}. Here, our proof mainly follows the lines of the one in \cite{partha}.\\

\noindent  Suppose that the hypotheses hold with $0 \in I$ and $\phi(0)=0$. That $\mathbb{E}(X)\in I$ is obvious. First, let us assume that $I$ is a compact interval, that is, $I=[a,b]$, with $a$ and $b$ finite and $a<b$. A convex function has left-hand and right-hand derivatives and then, is continuous (See Exercise 6 on Doc 03-09 of Chapter 4, page 191). Thus, $\phi$ is uniformly continuous on $I$. For $\varepsilon >0$, there exists $\delta >0$ such that
 
\begin{equation}
\left\vert x-y\right\vert \leq \delta \Rightarrow \left\vert \phi (x)-\phi
(y)\right\vert \leq \varepsilon .  \label{jensen01}
\end{equation}

\bigskip \noindent  We may cover $I$ with a finite number of disjoint intervals $E_{j}$ $(1\leq j\leq k)$, of diameters not greater than $\delta$. By using the Choice's Axiom, let us pick one $x_{j}$ in each $E_{j}$. Let $\mu$ be a une probability measure on $I$. We have 

\begin{eqnarray*}
\left\vert \int_{I}\phi (x) \ d\mu -\sum_{1\leq j\leq k}\phi (x_{j})\right\vert &=&\left\vert \sum_{1\leq j\leq k}\int_{E_{j}}\phi (x)\text{ }d\mu-\sum_{1\leq j\leq k}\phi (x_{j})\text{ }\mu (E_{j}) \right\vert \ (J02)\\
&=&\left\vert\sum_{1\leq j\leq k}\int_{E_{j}}\phi (x) \ d\mu -\sum_{1\leq j\leq k}\int_{E_{j}}\phi (x_{j}) \ d\mu \right\vert \\
&\leq& \sum_{1\leq j\leq k}\int_{E_{j}}\left\vert \phi (x)-\phi (x_{j})\right\vert \ d\mu \leq \sum_{1\leq j\leq k}\varepsilon \mu (E_{j})\leq \varepsilon . 
\end{eqnarray*}

\bigskip \noindent  We also have 

\begin{eqnarray*}
\left\vert \int_{I}x\text{ }d\mu -\sum_{1\leq j\leq k}x_{j}\text{ }\mu (E_{j})\right\vert
&=&\left\vert \sum_{1\leq j\leq k}\int_{E_{j}}x\text{ }d\mu -\sum_{1\leq j\leq k}x_{j}\text{ }\mu (E_{j})\right\vert\\
&=&\left\vert \sum_{1\leq j\leq k}\int_{E_{j}}x\text{ }d\mu -\sum_{1\leq j\leq k}\int_{E_{j}}x_{j} \right\vert\\
&\leq& \sum_{1\leq j\leq k}\int_{E_{j}}\left\vert x-x_{j}\right\vert d\mu\\
&\leq &\sum_{1\leq j\leq k}\delta \text{ }\mu (E_{j})\leq \delta .
\end{eqnarray*}

\bigskip \noindent  Then, by uniform continuity, we get 
\begin{equation}
\left\vert \phi \left(\int_{I}xd\mu \right)-\phi \left(\sum_{1\leq j\leq k}x_{j}\mu
(E_{j})\right)\right\vert \leq \varepsilon .  \label{jensen03}
\end{equation}

\bigskip \noindent  By applying the convexity of  $\phi$, we have 
\begin{equation*}
\phi \left(\int_{I}x\text{ }d\mu \right)\leq \varepsilon +\phi \left(\sum_{1\leq j\leq
k}x_{j}\text{ }\mu (E_{j})\right)\leq \varepsilon +\sum_{1\leq j\leq k}\phi (x_{j})%
\text{ }\mu (E_{j}).
\end{equation*}

\bigskip \noindent  By applying Formula (J02) to last term of the right-hand side, we have

\begin{equation*}
\phi `\left(\int_{I}x\text{ }d\mu \right)\leq 2\varepsilon +\int_{I}\phi (x)\text{ }d\mu,
\end{equation*}

\noindent for any $\varepsilon >0$. This implies 
\begin{equation}
\phi \left(\int_{I}x\text{ }d\mu \right)\leq \int_{I}\phi (x)\text{ }d\mu .
\label{jensen04}
\end{equation}

\bigskip \noindent  Now let $I$ be arbitrary and $\mu$ be a probability measure on $\mathbb{R}$%
. Put, for each $n\geq 1$, $I_{n}=[a_{n},b_{n}]$ with $(a_{n},b_{n})\rightarrow (-\infty
,+\infty )$ as $n\rightarrow \infty$ and $\mu (I_{n})>0$ for large values of $n$. Let us consider the probability measures
 $\mu _{n}$ on  $I_{n}$ defined by  
\begin{equation*}
\mu _{n}(A)=\mu (A)/\mu (I_{n}),\text{ }A\subset I_{n}.
\end{equation*}

\bigskip \noindent  Let us apply the inequality (\ref{jensen04}) to have 
\begin{equation*}
\phi \left(\int_{I_{n}}x\text{ }d\mu _{n}\right)\leq \int_{I_{n}}\phi (x)\text{ }d\mu
_{n}.
\end{equation*}

\bigskip \noindent  But, by the Monotone Convergence Theorem, we get  
\begin{equation*}
\int x\text{ }d\mu =\lim_{n\uparrow \infty }\int_{I_{n}}x\text{ }d\mu
=\lim_{n\uparrow \infty }\mu (I_{n})\int_{I_{n}}x\text{ }d\mu _{n}
\end{equation*}

\noindent and
\begin{equation*}
\lim_{n\uparrow \infty }\mu (I_{n})\int_{I_{n}}\phi (x)\text{ }d\mu
_{n}=\int \phi (x)\text{ }d\mu .
\end{equation*}

\bigskip \noindent  By using the  continuity of $\phi$, and the the Monotone Convergence Theorem, and the fact that $\int xd\mu$ exists, we conclude by  
\begin{eqnarray*}
\phi \left(\int_{I}xd\mu \right)&=&\lim_{n\rightarrow \infty }\phi (\int_{I_{n}}x \ d\mu)\\
&=&\lim_{n\rightarrow \infty }\phi (\mu (I_{n})\int_{I_{n}}x\text{ }d\mu _{n})\\
&=&\lim_{n\rightarrow \infty }\phi (\mu(I_{n})\int_{I_{n}}x \ d\mu _{n} + (1-\mu(I_{n})) \times 0)\\
&\leq& \lim_{n\rightarrow \infty }\mu (I_{n})\text{ }\phi \left(\int_{I_{n}}x \ d\mu _{n}\right)+(1-\mu (I_{n}))\text{ }\phi (0) \ (By \ convexity)\\
&\leq& \lim_{n\rightarrow \infty }\mu (I_{n})\text{ }\phi \left(\int_{I_{n}}xd\mu _{n}\right)\leq \lim_{n\rightarrow \infty }\mu (I_{n})\text{ }\phi\left(\int_{I_{n}}xd\mu _{n}\right)\\
&=& \lim_{n\rightarrow \infty }\mu (I_{n})\int_{I_{n}}\phi (x)\text{ }d\mu _{n} \ (Since \ \phi(0)=0)\\
&=&\int \phi (x)\text{ }d\mu .  \ (J03)
\end{eqnarray*}

\bigskip \noindent  The proof above is valid for any probability measure on $\mathbb{R}$. Since $X$ is integrable, $X$ is \textit{a.e.} finite and hence the support of $\mathbb{P}_X$ is a subset of $\mathbb{R}$. Hence, by applying (J3) to $\mathbb{P}_X$, we have the Jensen's inequality with the restrictions $0 \in I$, $\phi(0)=0$. We remove them as follows :\\

\noindent If $0\notin I$, we may enlarge $I$ to contains $0$ without any change of the inequality. If $\phi(0)\neq 0$, we may still apply the inequality to the convex function $\psi(x)=\phi(x)-\phi(0)$ which satisfies $\psi(0)=0$ and get the result.\\

\bigskip \noindent  \textbf{(b) Some applications of the Jensen's Inequality}.\\

\noindent The following stunning results on $L^p$ hold when the measure is a probability measure. They do not hold in general.\\

\noindent \textbf{(b1) Ordering the spaces $L^p$}.\\

\noindent Let $1<p<q$, $p$ finite but $q \in [1, +\infty]$. Let $X\in L^{q}$. Then $X\in L^{p}$ and 

$$
\left\Vert X\right\Vert _{p}\leq \left\Vert X\right\Vert _{q}.
$$

\bigskip \noindent  For $q=+\infty$, the inequality holds for any finite measure.\\

\bigskip \noindent  \textbf{Proof}. We consider two cases.\\

\noindent \textbf{Case $q$ finite}. Set $g_{1}(x)=x^{p}$, $g_{2}(x)=x^{q}$. Then the function $g_{2}\circ g_{1}^{-1}(x)=x^{q/p\text{ }}$ is convex on $(0,+\infty )$ since its second derivative is non-negative on $(0,+\infty)$. Let us set $X=g_{1}^{-1}(Y)$. In order to stay on $(0,+\infty)$, put $Z=|X|$ and take $Z=g_1^{-1}(Y)$, $Y\in (0,+\infty)$. The application of Jensen's Inequality leads to 
\begin{equation*}
g_{2}\circ g_{1}^{-1}(\mathbb{E}(Y))\leq \mathbb{E}(g_{2}\circ
g_{1}^{-1}(Y)).
\end{equation*}

\bigskip \noindent  Then we have
\begin{equation*}
g_{1}^{-1}(\mathbb{E}(Y))\leq g_{2}^{-1}(\mathbb{E}(g_{2}\circ
g_{1}^{-1}(Y)),
\end{equation*}

\bigskip \noindent  that is
\begin{equation*}
g_{1}^{-1}(\mathbb{E}(g_{1}(Z))\leq g_{2}^{-1}(\mathbb{E}(g_{2}(Z)).
\end{equation*}

\bigskip \noindent  This is exactly :
\begin{equation*}
\left\Vert X\right\Vert _{p}\leq \left\Vert X\right\Vert _{q}.
\end{equation*}

\bigskip \noindent  \textbf{Case $q=+\infty$}. By definition, $X \in L^{\infty}$ means that the set 
$$
\{M \in [0,+\infty[, \  |X| \leq M, \ \mathbb{P}-a.e\}
$$
*\noindent  is not empty and the infimum of that set is $\|X\|_{\infty}$. But for any $0\leq M<+\infty$ such that $|X| \leq M$, $\mathbb{P}$-\textit{a.e}. By taking the power and integrating, we get that

$$
\biggr(\int |X|^p \ d\mathbb{P}\biggr)^{1/p}\leq M.
$$ 

\bigskip \noindent  By taking the minimum of those values $M$, we get $\|X\|_{p}\leq \|X\|_{\infty}$.\\

\noindent \textbf{Conclusion}. If we have two real and finite numbers $p$ and $q$ such that $1\leq p \leq q$, we have the following ordering for $L^p$ spaces associated to a probability measure :

$$
L^{\infty} \subset L^q \subset L^p \subset L^1. 
$$

\bigskip \noindent  \textbf{(b2) Limit of the sequence of $L^p$-norm}.\\

\noindent We have

$$
\|X\|_{p} \nearrow \|X\|_{\infty} \text{ as } p\nearrow +\infty. \ (LN)
$$

\bigskip \noindent  \textbf{Proof}. If $\|X\|_{p_0}=+\infty$ for some $p_0\geq 1$, the results of Point (b2) above imply that $\|X\|_{\infty}=+\infty$ and $\|X\|_{p}=+\infty$ for all $p\geq p_0$ and the Formula (LN) holds.\\

\noindent Now suppose that $\|X\|_{p}<+\infty$ for all $p\geq 1$. By definition, $\|X\|_{\infty}=+\infty$ if the set
$$
\{M \in [0,+\infty[, \  |X| \leq M, \ \mathbb{P}-a.e\}
$$

\noindent is empty and is its infimum in the either case. In both cases, we have $\mathbb{P}(|X|>c)>0$ for all $c < \|X\|_{\infty}$ (as a consequence of the infimum). We get the following inequalities, which first exploit the relation : $|X| \leq \|X\|_{\infty}$, \textit{a.e.}. Taking the powers in that inequality and integrating yield, for $c<\|X\|_{\infty}$,

$$
\|X\|_{\infty} \geq \biggr( \int |X|^p \ d\mathbb{P} \biggr)^{1/p}\geq \biggr( \int_{(|X|\geq c)} |X|^p \ d\mathbb{P} \biggr)^{1/p} \geq c \biggr( \mathbb{P}(|X|\geq c)\biggr)^{1/p}.
$$

\bigskip \noindent  By letting first $p\rightarrow +\infty$, we get

$$
c \leq \liminf_{p\rightarrow +\infty} \|X\|_{p} \leq \limsup_{p\rightarrow +\infty} \|X\|_{p} \leq \|X\|_{\infty}.
$$

\bigskip \noindent  By finally letting $c\nearrow \|X\|_{\infty}$, we get the desired result.\\

\bigskip

\section{Moments of Real Random Variables} \label{proba_02_rv_sec_02}

\noindent \textbf{(a) Definition of the moments}.\\

\noindent The moments play a significant role in Probability Theory and in Statistical estimation. In the sequel, $X$ and $Y$ are two \textit{rrv}'s, $X_{1}$, $X_{2}$, ... and  $Y_{1}$, $Y_{2}$ are finite sequences of \textit{rrv}'s, $\alpha_{1}$, $\alpha_{2}$, ... and  $\beta_{1}$, $\beta_{2}$ are finite sequences of real numbers.\\

\noindent Let us define the following parameters, whenever the concerned expressions make sense.\\

\noindent \textbf{ (a1) Non centered moments of order $k\geq 1$} :

$$
m_{k}(X)=E\left\vert X\right\vert ^{k},
$$

\bigskip \noindent  which always exists as the integral of a non-negative random variable.\\

\bigskip \noindent  \textbf{ (a2) Centered Moment of order $k\geq 1$}. 
$$
\mu _{k}(X)=E\left\vert X-m_{1}\right\vert ^{k},
$$

\bigskip \noindent  which is defined if If $m_1(X)=\mathbb{E}X$ exists and is finite.\\

\bigskip \noindent  \textbf{(b) Focus on the centered moment of order 2}.\\

\noindent \textbf{(b1) Definition}.\\

\noindent If $\mathbb{E}X$ exists and is finite, the centered moment of second order 
$$
\mu_{2}(X)=\mathbb{E}\biggr(X-\mathbb{E}(X)\biggr)^{2},
$$

\bigskip \noindent  is called the variance of $X$. Throughout the textbook, we will use the notations

\begin{equation*}
\mu_{2}(X)=:\mathbb{V}ar(X)=:\sigma _{X}^{2}.
\end{equation*}

\bigskip \noindent  The number  $\sigma_{X}=\mathbb{V}ar(X)^{1/2}$ is called the standard deviation of $X$.\\

\bigskip \noindent  \textbf{(b2) Covariance between $X$ and $Y$}.\\

\noindent If $\mathbb{E}X$ and $\mathbb{E}Y$ exist and are finite, we may define the covariance between $X$ and $Y$ by

$$
\mathbb{C}ov(X,Y)=\mathbb{E}\biggr( (X-\mathbb{E}(X))\mathbb{((Y}-\mathbb{E}(Y))\biggr).
$$

\bigskip \noindent  \textbf{Warning}. It is important to know that the expectation operator is used in the Measure Theory and Integration frame, that is, $\mathbb{E}h(X)$ exists and is finite if and only if
$\mathbb{E}|h(X)|$ is finite. Later, when using Radon-Nikodym derivatives and replacing Lebesgue integrals by Riemann integrals, one should always remember this fact.\\

\bigskip \noindent  \textbf{Warning}. From now on, we implicitly assume the existence and the finiteness of the first moments of the concerned real random variables when using the variance or the covariance.

\bigskip \noindent  \textbf{(b3) Expansions of the variance and covariance}.\\

\noindent By expanding the formulas of the variance and the covariance and by using the linearity of the integral, we get, whenever the expressions make sense, that

\begin{equation*}
\mathbb{V}ar(X)=\mathbb{E}(X^{2})-\mathbb{E}(X)^{2},
\end{equation*}

\bigskip \noindent  (In other words, the variance is the difference between the non centered moment of order 2 and the square of the expectation), and

\begin{equation*}
\mathbb{C}ov(X,Y)=\mathbb{E}(XY)-\mathbb{E}(X)\mathbb{E}(Y).
\end{equation*}

\bigskip \noindent  \textbf{(b4) Two basic inequalities based on the expectation and the variance}.\\

\noindent The two first moments, when they exist, are basic tools in Statistical estimation. In turn, two famous inequalities are based on them. The first is the :\\

\noindent \textbf{Markov's inequality} : For any random variable $X$, we have for any $\lambda>0$
\begin{equation*}
\mathbb{P}(\left\vert X\right\vert >\lambda )\leq \frac{\mathbb{E}\left\vert
X\right\vert }{\lambda }.
\end{equation*}

\bigskip \noindent  (See Exercise 6 in Doc 05-02 in Chapter 6 in \cite{ips-mestuto-ang}). Next we have the : \\

\noindent \textbf{Tchebychev's inequality} : If  $X-\mathbb{E}(X)$ is defined \textit{a.e.}, then for any $\lambda >0$,
\begin{equation*}
\mathbb{P}(\left\vert X-\mathbb{E}(X)\right\vert >\lambda )\leq \frac{\mathbb{V}ar(X)}{\lambda ^{2}}.
\end{equation*}

\bigskip \noindent  This inequality is derived by applying the Markov's inequality to $|X-\mathbb{E}(X)|$ and by remarking that $(|X-\mathbb{E}(X)|>\lambda)=((X-\mathbb{E}(X))^2>\lambda^2)$, for any $\lambda>0$.\\

\bigskip \noindent  \textbf{(c) Remarkable properties on variances and covariances}.\\

\noindent Whenever the expressions make sense, we have the following properties.\\

\noindent (P1) $\mathbb{V}ar(X)=0$ if and only if $X=\mathbb(X)$ $a.s.$\\

\noindent (P2) For all $\lambda>0$, $\mathbb{V}ar(\lambda X)=\lambda ^{2}\mathbb{V}ar(X)$\\

\noindent (P3) We have

$$
\mathbb{V}ar(\sum_{1\leq i\leq k} \alpha_i X_{i})=\sum_{1\leq i\leq k} \mathbb{V}ar(X_{i}) \alpha_i^2 + 2\sum_{i<j} \mathbb{C}ov(X_{i},Y_{j}) \alpha_i  \alpha_j .
$$

\bigskip \noindent  (P4) We also have

$$
\mathbb{C}ov\biggr(\sum_{1\leq i\leq k}\alpha_i  X_{i}), \sum_{1\leq i\leq \ell} \beta_i Y_{i})) =\sum_{1\leq i\leq k} \sum_{1\leq j\leq \ell} \mathbb{C}ov(X_{i},Y_{j}) \alpha_i  \beta_j .
$$

\bigskip \noindent  (P5)] If $X$ and $Y$ are independent, then  $\mathbb{C}ov(X,Y)=0$.\\

\noindent (P6) Si $X_{1},...,X_{k}$ are pairwise independent, then  

$$
\mathbb{V}ar(\sum_{1\leq i\leq k} \alpha_i X_{i})=\sum_{1\leq i\leq k} \mathbb{V}ar(X_{i}) \alpha_i^2.
$$

\bigskip \noindent  (P7)  If none of $\sigma_{X}$ and  $\sigma_{Y}$ is null, then the coefficient

\begin{equation*}
\rho _{XY}=\frac{\mathbb{C}ov(X,Y)}{\sigma _{X}\sigma _{Y}},
\end{equation*}

\bigskip \noindent  is called the linear correlation coefficient between $X$ and $Y$ and satisfies
 
\begin{equation*}
\left\vert \rho _{XY}\right\vert \leq 1.
\end{equation*}

\bigskip \noindent  \textbf{Proofs or comments}. Most of these formulas are proved in the textbook \cite{ips-probelem-ang} of this series. Nevertheless we are going to make comments of the proofs at the light of Measure Theory and Integration and prove some of them.\\

\noindent (P1) We suppose that $\mathbb{E}(X)$ exists and is finite. We have $Y=(X-\mathbb{E}(X))^2\geq 0$ and $\mathbb{V}ar(X)=\mathbb{E}Y$. Hence, $\mathbb{V}ar(X)=0$ if and only if $Y=0$ \textit{a.e.} $\square$\\

\noindent (P2) This is a direct application of the linearity of the integral as recalled in Theorem \ref{proba_02_expect}. $\square$\\ 

\noindent (P3) This formula uses (P2) and the following the identity :

$$
\biggr(\sum_{1\leq i\leq k} a_{i}\biggr)^2=\sum_{1\leq i\leq k} a_i^2 +2\sum_{i<j} a_i a_j,
$$

\bigskip \noindent  where $a_i$, $1\leq i \leq k$, are real and finite numbers. Developing the variance and applying this alongside the linearity of the mathematical expectation together lead to the result. $\square$\\

\noindent (P4) This formula uses the following identity

$$
\biggr(\sum_{1\leq i\leq k} a_{i}\biggr)\biggr(\sum_{1\leq i\leq \ell} b{i}\biggr)=\sum_{1\leq i\leq k}\sum_{1\leq j\leq \ell} a_i b_j,
$$

\bigskip \noindent  where the $a_i$, $1\leq i \leq k$, and the $b_i$, $1\leq i \leq \ell$, are real and finite numbers. By developing the covariance and applying this alongside the linearity of the mathematical expectation lead to the result. $\square$\\

\noindent (P5) Suppose that $X$ and $Y$ are independent. Since $X$ and $Y$ are \textit{real} random variables, Theorem \ref{0118} implies that :  $\mathbb{E}(XY)=\mathbb{E}(X)\mathbb{E}(Y)$ by. Hence, by Point (b3) above, we get

$$
\mathbb{C}ov(X,Y)=\mathbb{E}(XY)-\mathbb{E}(X)\mathbb{E}(Y)=\mathbb{E}(X)\mathbb{E}(Y)-\mathbb{E}(X)\mathbb{E}(Y)=0. \square
$$

\bigskip \noindent  (P6) If the $X_i$'s are pairwise independent, the covariances in the formula in (P3) vanish and we have the desired result. $\square$\\

\noindent (P7) By applying the Cauchy-Schwartz inequality to $X-\mathbb{E}(X)$ and to $Y-\mathbb{E}(Y)$, that is the H\H{o}lder inequality for $p=q=2$, we get

$$
|\mathbb{C}ov(X,Y)| \leq \sigma_X \sigma_Y.
$$

\bigskip \noindent  If none of $\sigma_X$ and $\sigma_Y$ is zero, we get $|\rho_{XY}|\leq 1$. $\square$\\

\section{Cumulative distribution functions} \label{proba_02_rv_sec_03}

\bigskip \noindent  An important question in Probability Theory is to have parameters or functions which characterize probability laws. In Mathematical Statistics, these characteristics may be used in statistical tests. For example, if $X$ is a real value random variable having finite moments of all orders, that is : for all $k\geq 1$, $\mathbb{E}|X|^k<+\infty$. Does the sequence 
$(\mathbb{E}|X|^k)_{k\geq 1}$ characterize the probability law $\mathbb{P}_X$? This problem, named after the moment problem, will be addressed in a coming book. \\

\noindent The first determining function comes from the Lebesgue-Stieljes measure studied in Chapter 11 in \cite{ips-mestuto-ang}. We will use the results of that chapter without any further recall. \\

\noindent \textbf{(a) The cumulative distribution function of a real-random variable}.\\

\noindent Let $X : (\Omega ,\mathcal{A},\mathbb{P}) \mapsto  \mathbb{R}$ be a random real-valued random variable. Its probability law $\mathbb{P}_X$ satisfies :

$$
\forall x \in \mathbb{R}, \ \mathbb{P}_X(]-\infty,x])<+\infty.
$$

\bigskip \noindent  Hence, the function

$$
\mathbb{R} \ni x \mapsto F_X(x)=\mathbb{P}_X(]-\infty,x]),
$$

\bigskip \noindent  is a distribution function and  $\mathbb{P}_X$ is the unique probability-measure such that

$$
\forall (a,b) \in \mathbb{R}^2 \text{ such that } a\leq b, \ \mathbb{P}_X(]a,b])=F_X(b)-F_X(a).
$$

\bigskip \noindent  Before we go further, let us give a more convenient form of $F_X$ by writing for any $x \in \mathbb{R}$,

\begin{eqnarray*}
F_X(x)&=&\mathbb{P}_X(]-\infty,x])=\mathbb{P}(X^{-1}(]-\infty,x]))\\
&=&\mathbb{P}(\{\omega \in \Omega, \ X(\omega)\leq x\})=\mathbb{P}(X\leq x). 
\end{eqnarray*}

\bigskip \noindent  Now, we may summarize the results of the Lebesgue-Stieljes measure in the context of probability Theory.\\

\noindent \textbf{Definition}. For any real-valued random variable $X : (\Omega ,\mathcal{A},\mathbb{P}) \mapsto  \mathbb{R}$, the function defined by

$$
\mathbb{R} \ni \mapsto F_X(x)=\mathbb{P}(X \leq x),
$$

\bigskip \noindent  is called the cumulative distribution (cdf) function of $X$.\\

\noindent It has the two sets of important properties.\\

\noindent \textbf{(b) Properties of $F_X$}.\\

\noindent (1) It assigns non-negative lengths to intervals, that is 

$$
\forall (a,b) \in \mathbb{R}^2 \text{ such that } a\leq b, \ \Delta_{a,b}F=F_X(b)-F_X(a)\geq 0.
$$

\bigskip \noindent  (2) It is right-continuous at any point $t\in \mathbb{R}$.\\

\noindent (3) $F(-\infty )=\lim_{x \rightarrow -\infty} F(x)=0$ \ and $F(+\infty)=\lim_{x \rightarrow +\infty} F(x)=1$.\\

\noindent \textbf{Warning}. Point (1) means, in the case of one-dimension, that $F_X$ is non-decreasing. \textit{So, it happens that the two notions of non-negativity of lengths by $F$ and non-drecreasingness of  $F$ coincide in dimension one}. However, we will see that this is not the case in higher dimensions, and that  non-decreasingness is not enough to have a \textit{cdf}.\\

\noindent \textbf{(c) Characterization}.\\

\noindent The \textit{cdf} is a characteristic function of the probability law of a random variable with values in $\mathbb{R}$ from the following fact, as seen Chapter 11 in \cite{ips-mestuto-ang} of this series :\\

\noindent \textit{There exists a one-to-one correspondence between the class of Probability Lebesgue-Stieljes measures $\mathbb{P}_F$ on $\mathbb{R}$ and the class of \textbf{cfd}'s 
$F_{\mathbb{P}}$ on $\mathbb{R}$ according the relations}

$$
\biggr(\forall x \in \mathbb{R}, \ \ F_{\mathbb{P}}(x)=\mathbb{P}(]-\infty,x]) \biggr), \ \ \biggr(\forall (a,b) \in \left(\mathbb{R}\right), \ a\leq b, \ \ 
\mathbb{P}_F(]a,b])=\Delta_{a,b}F \biggr)
$$

\bigskip \noindent The \textit{cdf} is a characteristic function of the probability law of random variables. The means that two random real variable $X$ and $Y$ with the same distribution function have the same probability law.\\

\noindent \textbf{(d) How Can we Define a Random Variable Associated to a \textit{Cdf}}.\\

\noindent Let us transform the properties in Point (b) into a definition.\\

\noindent \textbf{(d1) Definition}. A function $F : \mathbb{R} \rightarrow [0,1]$ is \textit{cdf} if and only if conditions (1), (2) and (3) of Point (b) above are fulfilled.\\

\noindent Once we know that $F$ is (\textit{cdf}), can you produce a random variable $X : (\Omega ,\mathcal{A},\mathbb{P}) \mapsto  \mathbb{R}$ such that $F_X=X$? meaning : can we construct a probability space $(\Omega ,\mathcal{A},\mathbb{P})$ holding a random variable such that for all $x \in \mathbb{R}$, $F(x)=\mathbb{P}(X \leq x)$?\\

\noindent This is the simplest form the Kolmogorov construction. A solution is the following.\\

\bigskip \noindent  \textbf{(d2) A Simple form of Kolmogorov construction}.\label{kolmconst_01}\\

\noindent  Since $F$ is a \textit{cdf}, we may define the Lebesgue-Stieljes measure $\mathbb{P}$ on $(\mathbb{R} ,\mathcal{B}(\mathbb{R})$ defined by

$$
\mathbb{P}(]y,x])=\Delta_{y,x}F=F(x)-F(y), \ \ -\infty <y<x<+\infty. \ (LS11)
$$

\bigskip \noindent  By Conditions (3) in the definition of a \textit{cdf} in Point (b) above, $\mathbb{P}$ is normed and hence, is a probability measure. By letting $y \downarrow -\infty$ in (LS1), we get

$$
\forall x \in \mathbb{R}), \ F(x)=\mathbb{P}(]-\infty,x]). \ \ (LS12)
$$

\bigskip \noindent  Now take $\Omega=\mathbb{R}$, $\mathcal{A}=\mathcal{B}(\mathbb{R})$ and let $X : (\Omega ,\mathcal{A},\mathbb{P}) \mapsto  \mathbb{R}$ be the identity function
$$
\forall \omega \in \Omega, \ X(\omega)=\omega.
$$
  
\bigskip \noindent  It is clear that $X$ is a random variable and we have for $x\in \mathbb{R}$, we have

\begin{eqnarray*}
F_X(x)&=&\mathbb{P}(\{\omega \in \mathbb{R}, \ X(\omega)\leq x\})\\
&=&\mathbb{P}(\{\omega \in \mathbb{R}, \ \omega\leq x\})\\
&=&\mathbb{P}(]-\infty, x)=F(x),
\end{eqnarray*}

\bigskip \noindent  where we used (LS12). We conclude the $X$ admits $F$ as a \textsl{cdf}.\\

\noindent \textbf{Warning}. This construction may be very abstract at a first reading. If you feel confused with it, we may skip it and wait a further reading to catch it.

\bigskip \noindent  \textbf{(e) Decomposition of \textit{cdf} in discrete and continuous parts}.\\

\noindent Let $F$ be a \textit{cdf} on $\mathbb{R}$ and let us denote by $\mathbb{P}_{F}$ the associated Lebesgue-measure. We already know from Measure Theory that : $x \in \mathbb{R}$ is a continuity point of $F$ if and only if 

$$
\mathbb{P}_{F}(\{x\})=F(x)-F(x-0)=0, \ \ (CC)
$$

\bigskip \noindent  where for each $x \in \mathbb{R}$
$$
F(x+)\equiv F(x+0)=\lim_{h \searrow 0}F(x+h)
$$

\noindent and\\

$$
F(x-)\equiv F(x-0)=\lim_{h \searrow 0}F(x+h)
$$

\bigskip \noindent are the right-limit hand and the the left-limit hand of $F(\circ)$ at $x$, whenever they exist. In the present case, they do because of the monotonicity of $F$.\\

\noindent So, a \textsl{cdf} is continuous if and only if Formula (CC) holds for each $x \in \mathbb{R}$. In the general case, we are able to decompose the \textsl{cdf} into two non-negative distributions functions $F_c$ and $F_d$, where $F_c$ is continuous and $F_d$ is discrete in a sense we will define. As a reminder, a distribution function (\textit{df}) on $\mathbb{R}$ is a function satisfying only Conditions (1) and (2) in Point (b) above.\\

\noindent Let us define a discrete \textit{df} $F_d$ on $\mathbb{R}$ as a function such that there exists a countable number of distinct real numbers $\mathcal{D}=\{x_j,  \ j\in J\}$, 
$J \subset \mathbb{N}$ and a family of finite and positive real numbers $(p_j)_{j\in J}$ such that \label{ddf1}

$$
\forall x \in \mathbb{R}, F_d(x)=\sum_{x_j \leq x, j\in J} p_j <+\infty. \ (DDF1)
$$

\bigskip \noindent  Let $\nu$ be the discrete measure $\mathcal{D}$ defined by \label{ddf2}

$$
\forall (y,x)\in \mathbb{R}^2 \text{such that }y\leq x, \nu(]y,x])=\sum_{x_j \in ]y,x], j\in J} p_j <+\infty. \ (DDF2)
$$

\bigskip \noindent  By combining (DDF1) and (DDF2), we see that$F_d$ a discrete \textit{df} is a \textsl{df} of a counting measure which is finite on  bounded above intervals. It follows that for each $j \in J$,

$$
\nu(\{x_j\})=F_d(x_j)-F_d(x_j-0)=p_j>0. 
$$

\bigskip \noindent  This implies that a discrete \textsl{df} is never continuous at all points. We still may call $\mathcal{D}$ the support of $F_d$ by extension of the support of $\nu$.\\

\noindent We know that $F$, as a non-decreasing function, has at most a countable number of discontinuity points. Let us denote the set of those discontinuity points by $\mathcal{D}=\{x_j,  \ j\in J\}$, 
$J \subset \mathbb{N}$ and put $p_j=F(x_j)-F(x_j-0)>0$. Going Back to Measure Theory (see Solution of Exercise 1, Doc 03-06, Chapter 4 in \cite{ips-mestuto-ang} of this series), we have that

$$
\forall (y,x)\in \mathbb{R}^2 \text{such that }y \leq x,\ \sum_{x_j \in ]y,x], j\in J} F(x_j)-F(x_j-0)\leq F(x)-F(y).
$$

\bigskip \noindent  By letting $y \downarrow-\infty$, we have 
 
$$
\forall x \in \mathbb{R}, F_d(x)=\sum_{x_j \leq x, j\in J} p_j \leq F(x)<+\infty.
$$

\bigskip \noindent  Besides, the set discontinuity points of $F_d$ is $\mathcal{D}$ since discontinuity points $x$ of $F_d$ must satisfy $\nu(\{x\})>0$.\\
   
\noindent Next, let us define $F_c=F-F_d$. It is clear that $F_c$ is right-continuous and non-negative. Let us prove that $F_c$ is continuous. By the developments above, $F_c$ is continuous outside $\mathcal{D}$ and for each $j \in J$, we have

\begin{eqnarray*}
F_c(x_j)-F_c(x_j-0)&=&\biggr(F(x_j)-F(x_j-0)\biggr)-\biggr(F_d(x_j)-F_d(x_j-0)\biggr)\\
&=&p_j-p_j=0.
\end{eqnarray*}

\bigskip \noindent It remains to show $F_c$ is a \textsl{df} by establishing that : it assigns to intervals non-negative lengths. For each $x\in \mathbb{R}$, $h>0$, we have

$$
\Delta_{x,x+h}F_c=\biggr(F(x+h)-F(x)\biggr) -\biggr(F_d(x+h)-F_d(x)\biggr).
$$

\bigskip \noindent But, by definition $F_d(x+h)-F_d(x)$ is the sum of the jumps of $F$ at discontinuity points in $]x , \ x+h]$. We already know (otherwise, get help from a simple drawing) that this sum of jumps is less than  $F(x+h+0)-F(x)$ which is $F(x+h)-F(x)$ by right-continuity of $F$. Hence for all $x\in \mathbb{R}$, for all $h>0$, $\Delta_{x,x+h}F_c\geq 0$. In total, $F$ is a \textit{df}.\\

\bigskip \noindent  We get the desired decomposition : $F=F_c+F_d$. Suppose that we have another alike decomposition $F=F_c^{\ast}+F_d^{\ast}$. Since the functions are bounded, we get
$F_c-F_c^{\ast}=F_d^{\ast}-F_d$. Let us denote by $\mathcal{D}$ and $\mathcal{D}^{\ast}$ and by $p_x$ and $p_x^{\ast}$ the supports and the discontinuity jumps (at $x$) of $F_d$ and $F_d^{\ast}$ respectively.\\

\noindent If the supports are not equal, thus for $x \in \mathcal{D} \Delta \mathcal{D}^{\ast}$, $F_d^{\ast}-F_d$ is discontinuous at $x$.\\

\noindent If $\mathcal{D}=\mathcal{D}^{\ast}$ and $p_x \neq p_x^{\ast}$, the discontinuity jump of $F_d^{\ast}-F_d$ at $x$ is $p_x^{\ast}-p_x>0$.\\

\bigskip \noindent  Since none of the two last conclusions is acceptable, we get that the equation $F_c-F_c^{\ast}=F_d^{\ast}-F_d$ implies that $F_d^{\ast}$ and $F_d$ have the same support and the same discontinuity jumps, and hence are equal and then so are $F_c^{\ast}$ and $F_c$.\\

\noindent We get the following important result.\\

\noindent \textbf{Proposition}. A \textit{cdf} $F$ is decomposable into the addition of two non-negative distribution functions (\textit{df}) $F_c$ and $F_d$, where $F_c$ is continuous and $F_d$ is discrete and $F_c(-\infty)=F_d(-\infty)=0$. The decomposition is unique.\\

\noindent \textbf{Warning}. Such a result is still true for a \textit{df} but the condition $F_c(-\infty)=F_d(-\infty)=0$ is not necessarily true.\\

\noindent \textbf{NB}. We did not yet treat the probability density existence and its use for real random variables. This will be done in the next section which is concerned with random vectors.

\section{Random variables on $\overline{\mathbb{R}}^d$ or Random Vectors}  \label{proba_02_rv_sec_04} 

\noindent  \textbf{(a) Introduction}.\\

\noindent Random vectors are generalizations of real random variables. A random vector of dimension $d\geq 1$ is a random variable 
$$
X : (\Omega ,\mathcal{A},\mathbb{P}) \mapsto  (\overline{\mathbb{R}}^d, \mathcal{B}_{\infty}(\overline{\mathbb{R}}^d).
$$

\noindent  with values in $\overline{\mathbb{R}}^d$.\\

\noindent \textbf{Important Remarks}. In general, it is possible to have $\left(\overline{\mathbb{R}}^d, \mathcal{B}_{\infty}\left(\overline{\mathbb{R}}^d\right)\right)$ as the set of values of random vectors, especially when we are concerned with general probability laws. But, the most common tools which are used for the study of random vectors such as the cumulative random vectors, the characteristic functions, the absolute probability density function are used for finite component random vectors with values in $\left(\mathbb{R}^d, \mathcal{B}(\mathbb{R}^d\right)$.\\

\noindent Throughout this section, we use a random vector with $d$ components as follows. Let 
$$
X=\left[ 
\begin{array}{c}
X_{1} \\
X_{2} \\ 
\cdot \cdot \cdot \\
X_{d-1} \\ 
X_{d}
\end{array}
\right]
$$

\bigskip \noindent  From Measure Theory, we know that $X$ is a random variable if and only if each $X_i$, $1\leq i \leq d$, is a real random variable.\\

\noindent If $d=1$, the random vector becomes a real random variable, abbreviated (\textit{rrv}).\\

\noindent \textbf{Notation}. To save space, we will rather use the transpose operator and write $X^t=(X_1, X_2, \cdots, X_d)$ or $X=(X_1, X_2, \cdots, X_d)^t$. Let $(Y_1, Y_2, \cdots, X_d)^t$ another be $d$-random vector and two other random vectors $(Z_1, Z_2, \cdots, Z_r)^t$ and $(T_1, T_2, \cdots, T_s)^t$ of dimensions $r\geq 1$ and $s\geq 1$,  all of them being defined on $(\Omega ,\mathcal{A},\mathbb{P})$.\\

\noindent \textbf{Matrix Notation}. \noindent To prepare computations on the matrices, let us denote any real matrix $A$ of $r\geq 1$ lines and $s\geq 1$ columns in the form $A=(a_{ij})_{1\leq i \leq r, \ 1\leq j \leq s}$, where the lowercase letter $a$ is used to denote the elements of the matrix whose name is the uppercase letter $A$. As well, we will use the notation $(A)_{ij}=a_{ij}$, $1\leq i \leq r, \ 1\leq j \leq s$.\\

\noindent A matrix of $r$ lines and $s$ columns is called a $(r \times s)$-matrix, a square matrix with $r$ lines and $r$ columns is a $r$-matrix and a vector of $r$ components is a $d$-vector.\\ 

\noindent The $s$ columns of a matrix $A$ are elements of $\mathbb{R}^r$ and are denoted by $A^1$, $A^2$, ..., $A^s$. The $r$ lines of the matrix $A$ are $(1 \times s)$-matrices denoted $A_1$, ...,$A_r$, that is $A_1^t$, ...,$A_r^t$ belong to $\mathbb{R}^s$.\\

\bigskip \noindent  So, for $1\leq j \leq s$, $A^j=(a_{1j}, a_{2j}, \cdots, a_{rj})^t$ and for  $1\leq i\leq r$, $A_i=(a_{i1}, a_{i2}, \cdots, a_{is})$.\\

\noindent We also have $A=[A^1, A^2, ..., A^d]$ and $A^t=[A_1^t, A_2^t, ..., A_p^t]$.\\

\noindent Introduce the scalar product in $\mathbb{R}^s$ in the following way. Let $x$ and $y$ be two elements $\mathbb{R}^s$ with $x^t=(x_1,\cdots, x_s)$ and $y^t=(y_1,\cdots, y_s)$.\\

\noindent We define the scalar product $<x,y>$ of $x$ and $y$ as the matrix product of the $(1 \times s)$-matrix $x^t$ by the $(s \times 1)$-matrix $y$ which results in the real number

$$
<x,y> \ = \ x^ty=\sum_{i=1}^{s} x_iy_i.
$$ 

\bigskip \noindent  With the above notation, the matrix operations may be written in the following way.\\

\noindent If (1) \textit{Sum of matrices of same dimensions}. If $A=(a_{ij})_{1\leq i \leq r, \ 1\leq j \leq s}$ and $B=(b_{ij})_{1\leq i \leq r, \ 1\leq j \leq s}$ are two $(r \times s)$-matrix, then $A+B$ is the $(r \times s)$-matrix : $A+B=(a_{ij}+b_{ij})_{1\leq i \leq r, \ 1\leq j \leq s}$, that is : $(A+B)_{ij}=(A)_{ij}+(B)_{ij}$ for $1\leq i \leq r, \ 1\leq j \leq s$.\\

\noindent If (2) \textit{Multiplication by a scalar}. If $\lambda$ is a real number and if $A=(a_{ij})_{1\leq i \leq r, \ 1\leq j \leq s}$, then $\lambda A$ is the $(r \times s)$-matrix : 
$\lambda A=(\lambda a_{ij})_{1\leq i \leq r, \ 1\leq j \leq s}$, that is :$(\lambda A)_{ij}=\lambda (A)_{ij}$ for $1\leq i \leq r, \ 1\leq j \leq s$.\\

\noindent If (3) \textit{Product of Matrices}. If $A=(a_{ij})_{1\leq i \leq r, \ 1\leq j \leq s}$ and $B=(b_{ij})_{1\leq i \leq s, \ 1\leq j \leq q}$ such that the number of columns of $A$ (the first matrix) is equal to the number of lines of $B$ (the second of the second matrix), the product matrix $AB$ is a $(r,q)$-matrix defined by
$$
AB=(A_iB^j)_{1\leq i \leq r, \ 1\leq j \leq q},\ (PM1)
$$

\bigskip \noindent  that is, for $1\leq i \leq r, \ 1\leq j \leq q$,

$$
(AB)_{ij}=(A_iB^j)=\sum_{k=1}^d a_{ik} b_{kj}. \ (MP2)
$$

\bigskip \noindent  \textbf{(b) Variance-covariance and Covariance Matrices}.\\

\bigskip \noindent  {(b1) Definition of Variance-covariance and Covariance Matrices}.\\

\noindent We suppose that the components of our random vectors have finite second moments. We may define\\

\noindent (i) \textbf{the  mathematical expectation vector $\mathbb{E}(X)$ of $X\in \mathbb{R}^d$} by the vector

$$
\mathbb{E}(X)^t=(\mathbb{E}(X_1), \mathbb{E}(X_2), \cdots, \mathbb{E}(X_d)),
$$

\bigskip \noindent  (ii) \textbf{the covariance matrix $\mathbb{C}ov(X,Y)$ between $X \in \mathbb{R}^d$ and $Z \in \mathbb{R}^r$ } by the $(d \times r)$-matrix

$$
\mathbb{C}ov(X,Y)=\Sigma_{XY}=\mathbb{E}\biggr((X-\mathbb{E}(X)) (Z-\mathbb{E}(Z))^t \biggr),
$$

\bigskip \noindent  in an other notation

$$
\mathbb{C}ov(X,Y)=\Sigma_{XY}=\biggr(\mathbb{E}(X_i-\mathbb{E}(X_i)) (Z_j-\mathbb{E}(Z_j)) \biggr)_{1\leq i \leq d, \ 1\leq j \leq r},
$$

\bigskip \noindent  (iii) \textbf{the variance-covariance matrix $\mathbb{V}ar(Y)$ of $X \in \mathbb{R}^d$} by the $(d \times d)$-matrix

\begin{eqnarray*}
\mathbb{V}ar(X)&=&\Sigma_{X}=\mathbb{E}\biggr((X-\mathbb{E}(X)) (X-\mathbb{E}(X))^t \biggr)\\
&=&\biggr( \mathbb{E}(X_i-\mathbb{E}(X_i)) \mathbb{E}(X_j-\mathbb{E}(X_j)) \biggr)_{1\leq j \leq d, \ 1\leq j \leq d}. \ \square
\end{eqnarray*}

\bigskip \noindent  Let us explain more the second definition. The matrix 
$$
(X-\mathbb{E}(X)) (Z-\mathbb{E}(Z))^t
$$ 

\bigskip \noindent   is the product of the $(d \times 1)$-matrix $X-\mathbb{E}(X)$, with

$$
(X-\mathbb{E}(X))^t=(X_1-\mathbb{E}(X_1), X_2-\mathbb{E}(X_2), \cdots, X_d-\mathbb{E}(X_d)),
$$ 

\bigskip \noindent  by the $(1 \times r)$-matrix with 
$$
(Z-\mathbb{E}(Z))^t=(Z_1-\mathbb{E}(Z_1), Z_2-\mathbb{E}(Z_2), \cdots, X_r-\mathbb{E}(Z_r)).
$$

\bigskip \noindent    The $(ij)$-element of the product matrix, for ${1\leq i \leq d, \ 1\leq j \leq r}$, is

$$
(X_i-\mathbb{E}(X_i)) (Z_j-\mathbb{E}(Z_j)). 
$$

\bigskip \noindent  By taking the mathematical expectations of those elements, we get the matrix of covariances

$$
\mathbb{C}ov(X,Y)=\biggr(\mathbb{C}ov(X_i,Y_j) \biggr)_{1\leq i \leq d, \ 1\leq j \leq r}
$$

\bigskip \noindent  For $X=Z$ (and then $d=r$), we have

$$
\mathbb{V}ar(X) \equiv \mathbb{C}ov(X,X) =\biggr(\mathbb{C}ov(X_i,X_j) \biggr)_{1\leq i \leq d, \ 1\leq j \leq d}.
$$

\bigskip \noindent  We have the following properties.\\

\bigskip \noindent  \textbf{(b2) Properties}.\\

\noindent Before we state the properties, it is useful to recall that linear mappings from $\mathbb{R}^d$ to $\mathbb{R}^r$ are of the form 
$$
\mathbb{R}^d \ni X \mapsto AX \in \mathbb{R}^r, 
$$

\bigskip \noindent  where $A$ is a  $(r \times d)$-matrix of real scalars. Such mappings are continuous (uniformly continuous, actually) and then measurable with respect to the usual $\sigma$-algebras on  
$\mathbb{R}^d$ and $\mathbb{R}^p$.\\

\noindent  Here are the main properties of the defined parameters.\\

\noindent (P1) For any $\lambda \in \mathbb{R}$, 

$$\mathbb{E}(\lambda X)=\lambda \mathbb{E}(X).$$

\bigskip \noindent   (P2) For two random vectors $X$ and $Y$ of the same dimension $d$, $AX \in \mathbb{R}^p$ and

$$\mathbb{E}(X+Y)=\mathbb{E}(X)+\mathbb{E}(Y).$$

\bigskip \noindent   (P3) For any $(p \times d)$-matrix $A$ and any $d$-random vector $X$, 

$$
\mathbb{E}(AX)=A \mathbb{E}(X) \in \mathbb{R}^p.
$$

\bigskip \noindent   (P4) For any $d$-random vector $X$  and any $s$-random vector $Z$, 
$$
\mathbb{C}ov(X,Z)=\mathbb{C}ov(Z,X)^t.
$$

\bigskip \noindent   (P5) For any $(p \times d)$-matrix $A$, any $(q \times s)$-matrix $B$, any $d$-random vector $X$ and any $s$-random vector $Z$,

$$
\mathbb{C}ov(AX, BZ)= A \mathbb{C}ov(X,Z) B^t, 
$$

\bigskip \noindent  which is a $(p,q)$-matrix.\\

\bigskip \noindent  \textbf{Proofs}.\\

\noindent We are just going to give the proof of (P3) and (P4) to show how work the computations here.\\

\noindent \textit{Proof of (P3)}. The $i$-th element of the column vector $AX$ of $\mathbb{R}^p$, for $1\leq i \leq p$, is 

$$
(AX)_i=A_iX=\sum_{1\leq j \leq d} a_{ij}X_j
$$ 

\bigskip \noindent  and its real mathematical expectation, is

$$
\mathbb{E}(AX)_i=\sum_{1\leq j \leq d} a_{ij}\mathbb{E}(X_j).
$$ 

\bigskip \noindent  But the right-hand member is, for $1\leq i \leq p$, the $i$-th element of the column vector $A\mathbb{E}(X)$. Since $\mathbb{E}(AX)$ and $A\mathbb{E}(X)$ have the same components, we get

$$
\mathbb{E}(AX)=A\mathbb{E}(X). \ \square
$$

\bigskip \noindent  \textit{Proof of (P5)}. We have

$$
\mathbb{C}ov(AX, BZ)=\mathbb{P}\left( (AX-\mathbb{E}(AX))(BZ-\mathbb{E}(BZ))^t  \right). \ (COV1)
$$

\bigskip \noindent  By (P5), we have

\begin{eqnarray*}
(AX-\mathbb{E}(AX))(BX-\mathbb{E}(BZ))&=&A(X-\mathbb{E}(X)) (B(Z-\mathbb{E}(Z)))^t\\
&=&A \biggr( (X-\mathbb{E}(X)) (Z-\mathbb{E}(Z))^t\biggr) B^t.
\end{eqnarray*}

\bigskip \noindent  Let us denote $C=(X-\mathbb{E}(X)) (Z-\mathbb{E}(Z))^t$. We already know that 

$$
c_{ij}=\biggr((X_i-\mathbb{E}(X_i))(Z_j-\mathbb{E}(Z_j))\biggr), \ (i,j)\in \{1,...,d\}^2.
$$

\bigskip \noindent   Let us fix $(i,j) \in \{1,...,p\} \times \{1,...,q\}$. The $ij$-element of the $(p,q)$-matrix $ACB^t$ is

$$
(ACB^t)_{ij}=(AC)_i (B^t)^j.
$$ 

\bigskip \noindent  But  the elements of $i$-th line of $AC$ are $\{A_iC^1,A_2C^2, ..., A_iC^p\}$ and the column $(B^t)^j$ contains the elements of the $j$-th line of $B$, that is $b_{j1}$, $b_{j2}$, ...,$b_{js}$. We get

\begin{eqnarray*}
(ACB^t)_{ij}&=&\sum_{1\leq k\ leq s} \biggr((AC)_i\biggr)_k \biggr((B^t)^j\biggr)_k\\
&=& \sum_{1\leq k\ \leq s} \biggr(A_iC^k \biggr) b_{jk}\\
&=& \sum_{1\leq k\ \leq s} \sum_{1\leq p\ \leq p} a_{ih} (X_h-\mathbb{E}(X_h))(Z_k-\mathbb{E}(Z_k)) b_{jk}. \ (COV2) \\
\end{eqnarray*}

\bigskip \noindent  Hence, by applying Formula (COV1), the $ij$-element of $\mathbb{C}ov(AX, BZ)$ is

$$
\mathbb{E}\biggr((ACB^t)_{ij}\biggr) = \sum_{1\leq k\ \leq s} \sum_{1\leq p\ \leq p} a_{ih} \mathbb{C}ov(X_h,Z_k) b_{jk}. \ (COV3) \\
$$

\bigskip \noindent  Actually we have proved that for any $(p\times d)$-matrix $A$, for any $(d \times s)$-matrix and for any $(q \times)$-matrix, the $ij$-element of $ACB^t$ is given by

$$
\sum_{1\leq k\ \leq s} \sum_{1\leq p\ \leq p} a_{ih} c_{hk} b_{jk}. \ (ACBT)
$$

\bigskip \noindent  When applying this to Formula (COV2), we surely have that

$$
\mathbb{C}ov(AX, BZ)=A \mathbb{C}ov(X, Z) B^t.
$$

\bigskip \noindent  \textbf{(b3) Focus on the Variance-covariance matrix}.\\

\noindent (P6) Let $A$ be a $(p\times d)$-matrix and $X$ be a $d$-random vector. We have

$$
\Sigma_{AX}= A \Sigma_{X} A^t.
$$

\bigskip \noindent  (P7) The Variance-covariance matrix $\Sigma_X$ of $X \in \mathbb{R}^d$ is a positive matrix as a quadratic form, that is :

$$
\forall u \in \mathbb{R}^d, \ u \Sigma_X u^t \geq 0.
$$

\bigskip \noindent  If $\Sigma_X$ is invertible, then it is definite-positive that is

$$
\forall u \in \mathbb{R}^d\setminus \{0\}^d, \ u \Sigma_X u^t > 0.
$$

\bigskip \noindent  (P8) $\Sigma_X$ is symmetrical and there exists an  orthogonal $d$-matrix $T$ such that $T\Sigma_X T^t$ is a diagonal matrix

$$
T\sigma_X T^t=diag(\delta_1, \delta_2, ..., \delta_d),
$$

\bigskip \noindent  with non-negative eigen-values $\delta_j\geq 0$, $1\leq j\leq d$. Besides we have the following facts :\\

\noindent (P8a) The columns of $T$ are eigen-vectors of $\Sigma_X$ respectively associated the eigen-values $\delta_j\geq 0$, $1\leq j\leq d$ respectively.\\

\noindent (P8b) The columns, as well as the lines, of $T$ form an orthonormal basis of $\mathbb{R}^d$.\\

\noindent (P8c) $T^{-1}=T$.\\

\noindent (P8d) The number of positive eigen-values is the rank of $\Sigma_X$ and $\Sigma_X$ is invertible if and only all the eigen-values are positive.\\

\noindent (P8e) The determinant of $\Sigma_X$ is given by

$$
\left|\Sigma_X\right| \equiv det(\Sigma_X) =\prod_{j=1}^{d} \delta_j.
$$

\bigskip \noindent  \textbf{Proofs}. Property (P6) is a consequence of (P5) for $A=B$. Formula (P8) and its elements are simple reminders of Linear algebra and diagonalization of symmetrical matrices. The needed reminders are gathered in Subsection \ref{proba_02_appendix_ortho} in Section \ref{proba_02_appendix}.  The only point to show is (P7). And by (P5), we have for any $u \in \mathbb{R}^d$

\begin{eqnarray*}
u \Sigma_X u^t&=& u \mathbb{E}\biggr( (X-\mathbb{E}(X))(X-\mathbb{E}(X))^t \biggr) u^t\\
&=& \mathbb{E}\biggr( u(X-\mathbb{E}(X))(X-\mathbb{E}(X))^t u^t\biggr). 
\end{eqnarray*}

\bigskip \noindent  But $u(X-\mathbb{E}(X))(X-\mathbb{E}(X))^t u^t=\biggr(u(X-\mathbb{E}(X)\biggr) \biggr(u(X-\mathbb{E}(X))\biggr)^t$. Since $u(X-\mathbb{E}(X))$ is $d$-vector, we have

$$
\biggr(u(X-\mathbb{E}(X)\biggr) \biggr(u(X-\mathbb{E}(X))\biggr)^t=\|u(X-\mathbb{E}(X)\|^2\geq 0.
$$

\bigskip \noindent  Hence $u \sigma_X u^t\geq 0$. $\blacksquare$.\\

\bigskip \noindent  \textbf{(c) Cumulative Distribution Functions}.\\

\noindent The presentation of \textit{cdf}'s on $\mathbb{R}^d$ follows that lines we already used  for \textit{cdf}'s on $\mathbb{R}$. But the notations are heavier.\\

\noindent Let us recall the notion of volume we already introduced in Chapter 11 in \cite{ips-mestuto-ang}. \\

\noindent \textbf{(c1) Notion of Volume of cuboids by $F$}.\\

\noindent \textbf{Simple case}. Let us begin by the case $d=2$. Consider a rectangle 
$$
]a,b]=]a_1,b_1] \times ]a_2,b_2]=\prod_{i=1}^{2} ]a_i,b_i],
$$

\noindent  for $a=(a_{1},a_{2})\leq b=(b_{1},b_{2})$ meaning $a_i \leq b_i$, $1\leq i \leq 2$. The volume of $]a,b]$ by $F$ is denoted

\begin{equation*}
\Delta F(a,b)=F(b_{1},b_{2})-F(b_{1},a_{2})-F(a_{1},b_{2})+F(a_{1},a_{2}).
\end{equation*}

\bigskip \noindent   \textbf{Remark}. In the sequel we will use both notations $\Delta_{a,b}F$ and $\Delta F(a,b)$ equivalently. The function $\Delta F(a,b)$ is obtained according to the following rule :\\

\noindent  \textbf{Rule of forming $\Delta F(a,b)$}. First consider $F(b_{1},b_{2})$ the value of the distribution function at the right endpoint $b=(b_{1},b_{2})$ of the interval $]a,b]$. Next proceed to the replacements of each $b_{i}$ by $a_{i}$ by replacing exactly one of them, next two of them etc., and add each value of $F$ at the formed points, with a sign \textit{plus} $(+)$ if the number
of replacements is even and with a sign \textit{minus} $(-)$ if the number of replacements is odd.\\

\noindent  We also may use a compact formula. Let $\varepsilon =(\varepsilon _{1},\varepsilon_{2})\in \{0,1\}^{2}.$ We have four elements in $\{0,1\}^{2}:$ $(0,0),$ $%
(1,0),$ $(0,1),$ $(1,1).$ Consider a particular $\varepsilon _{i}=0$ or $1$, we have

\begin{equation*}
b_{i}+\varepsilon _{i}(a_{i}-b_{i})=\left\{ 
\begin{tabular}{lll}
$b_{i}$ & $if$  & $\varepsilon _{i}=0$ \\ 
$a_{i}$ & $if$ & $\varepsilon _{i}=1$%
\end{tabular}%
\right. .
\end{equation*}%

\bigskip \bigskip \noindent  So, in
\begin{equation*}
F(b_{1}+\varepsilon _{1}(a_{1}-b_{1}),b_{1}+\varepsilon _{2}(a_{2}-b_{2})),
\end{equation*}

\bigskip \noindent  the number of replacements of the $b_{i}$'s by the corresponding $a_{i}$ is the number of the coordinates of $\varepsilon =(\varepsilon _{1},\varepsilon
_{2})$ which are equal to the unity $1$. Clearly, the number of replacements is 
\begin{equation*}
s(\varepsilon )=\varepsilon _{1}+\varepsilon _{2}=\sum_{i=1}^{2}\varepsilon_{i}
\end{equation*}

\bigskip \noindent  We way rephrase the Rule of forming  $\Delta F(a,b)$ into this formula
\begin{equation*}
\Delta F(a,b)=\sum_{\varepsilon =(\varepsilon _{1},\varepsilon _{2})\in
\{0,1\}}(-1)^{s(\varepsilon )}F(b_{1}+\varepsilon
_{1}(a_{1}-b_{1}),b_{1}+\varepsilon _{2}(a_{2}-b_{2})).
\end{equation*}

\bigskip \noindent  We may be more compact by defining the product of vectors as the vector of the products of coordinates as

\begin{equation*}
(x,y)\ast (X,Y)=(x_1X_1, ..., y_d Y_d), \ d=2.
\end{equation*}

\bigskip \noindent  The formula becomes 
\begin{equation*}
\Delta F(a,b)=\sum_{\varepsilon \in \{0,1\}}(-1)^{s(\varepsilon
)}F(b+\varepsilon \ast (a-b)).
\end{equation*}

\bigskip \noindent  Once the procedure is understood for $d=2$, we may proceed to the general case.\\

\bigskip
\noindent \textbf{General case, $d\geq 1$}.\\

\noindent Let $a=(a_{1},...,a_{d})\leq b=(b_{1},...,b_{d})$ two points of $\mathbb{R}^{d}$. The volume of the cuboid

$$
]a,b]=\prod_{i=1}^{d} ]a_i,b_i],
$$

\noindent  by $F$, is defined by

\begin{equation*}
\Delta F(a,b)=\sum_{\varepsilon =(\varepsilon _{1},...,\varepsilon _{d})\in
\{0,1\}^{k}}(-1)^{s(\varepsilon )}F(b_{1}+\varepsilon
_{1}(a_{1}-b_{1}),...,b_{d}+\varepsilon _{k}(a_{d}-b_{d}))
\end{equation*}

\bigskip \noindent  or
\begin{equation*}
\Delta F(a,b)=\sum_{\varepsilon \in \{0,1\}^{d}}(-1)^{s(\varepsilon)}F(b+\varepsilon \ast (a-b)).
\end{equation*}

\bigskip  \bigskip \noindent  Similarly to the case $d=2$, we have the \\ 

\noindent \textbf{General rule of forming $\Delta F(a,b)$}. $\Delta F(a,b)$ in formed as follows. First consider $F(b_{1},b_{2},...,b_{d})$ the value of $F$ at right endpoint $b=(b_{1},b_{2},...,b_{d})$ of the interval $]a,b]$. Next proceed to the replacement of each $b_{i}$ by $a_{i}$ by replacing exactly one of them, next two of them etc., and add the each value of $F$ at these points with a sign plus $(+)$ if the number of replacements is even and with a sign minus $(-)$ if the number of replacements is odd.\\

\bigskip \noindent  \textbf{(c2) Cumulative Distribution Function}.\\

\noindent In this part, we study finite components vectors.\\

\noindent \textbf{Definition}. For any real-valued random variable $X : (\Omega ,\mathcal{A},\mathbb{P}) \mapsto  \mathbb{R}^d$, the function defined by

$$
\mathbb{R}^d \ni x \mapsto F_X(x)=\mathbb{P}(X \leq x),
$$

\bigskip \noindent  where $x^t=(x_1,...,x_d)$ and

$$
F_X(x_1,...,x_d)=\mathbb{P}(X_1 \leq x_1, X_2 \leq x_2, ..., X_d \leq x_d)=\mathbb{P}_X\biggr( \prod_{i=1}^{d}]-\infty, x_i]\biggr).
$$

\bigskip \noindent  is called the cumulative distribution (cdf) function of $X$.\\

\noindent It has the two sets of important properties.\\

\noindent \textbf{Properties of $F_X$}.\\

\noindent (1) It assigns non-negative volumes to cuboids, that is 

$$
\forall (a,b) \in (\mathbb{R}^d)^2 \text{ such that } a\leq b, \ \Delta_{a,b}F\geq 0.
$$

\bigskip \noindent  (2) It is right-continuous at any point $t\in \mathbb{R}^d$, that is,
 
\begin{equation*}
F_{X}(t^{(n)})\downarrow F_{m}(t)
\end{equation*}

\bigskip \noindent  as 
\begin{equation*}
(t^{(n)}\downarrow t) \Leftrightarrow (\forall 1\leq i\leq d, \ t_{i}^{(n)}\downarrow t_{i}).
\end{equation*}

\bigskip \noindent  (3) $F_X$ satisfies the limit conditions :\\

\noindent Condition (3-i)%
\begin{equation*}
\lim_{\exists i,1\leq i\leq k,t_{i}\rightarrow -\infty
}F_{X}(t_{1},...,t_{k})=0
\end{equation*}%

\bigskip \noindent  and Condition(3-ii)

\begin{equation*}
\lim_{\forall i,1\leq i\leq k,t_{i}\rightarrow +\infty
}F_{X}(t_{1},...,t_{k})=1.
\end{equation*}

\bigskip \noindent  As we did in one dimension, we have :\\

\noindent \textbf{Definition}. A function $F : \mathbb{R}^d \rightarrow [0,1]$ is \textbf{cdf} if and only if Conditions (1), (2) and (3) above hold.\\

\noindent \textbf{(c3) Characterization}.\\

\noindent The \textit{cdf} is a characteristic function of the probability law of random variables of $\mathbb{R}^d$ from the following fact, as seen in Chapter 11 in \cite{ips-mestuto-ang} of this series :\\

\noindent \textit{There exists a one-to-one correspondence between the class of Probability Lebesgue-Stieljes measures $\mathbb{P}_F$ on $\mathbb{R}^d$ and the class of \textbf{cfd}'s 
$F_{\mathbb{P}}$ on $\mathbb{R}^d$ according the relations}

$$
\forall x \in \mathbb{R}^d, \ \ F_{\mathbb{P}}(x)=\mathbb{P}(]-\infty,x])
$$

\noindent and

$$
\forall (a,b) \in \left(\mathbb{R}^d\right), \ a\leq b, \ \ \mathbb{P}_F(]a,b])=\Delta_{a,b}F.
$$

\bigskip \noindent  This implies that two $d$-random vectors $X$ and $Y$ having the same distribution function have the same probability law.\\

\noindent \textbf{(c4) Joint \textit{cdf}'s and marginal \textit{cdf}'s}.\\

\noindent Let us begin by the sample case where $d=2$. Let $X \ : \ (\Omega ,\mathcal{A},\mathbb{P}) \rightarrow \mathbb{R}^2$ be a random coupe, with $X^t=(X_1,X_2)$. We have,

\begin{eqnarray*}
(X_1\leq x)&=&X_1^{-1}(]-\infty, x])= X_1^{-1}(]-\infty, x]) \cap X_1^{-1}(]-\infty, +\infty])\\
&=&\lim_{y\uparrow +\infty} X_1^{-1}(]-\infty, x]) \cap X_2^{-1}(]-\infty, y])
=\lim_{y\uparrow +\infty} (X_1\leq x, X_2\leq y)
\end{eqnarray*}

\noindent  and by applying the Monotone Convergence Theorem, we have
 
\begin{eqnarray*}
\forall x\in \mathbb{R}, \ F_{X_1}(x)&=&\mathbb{P}(X_1\leq x)\\
&=&\lim_{y\uparrow +\infty} \mathbb{P}(X_1\leq x, X_2\leq y)=\lim_{y\uparrow +\infty} F_{(X_1,X_2)}(x,y).
\end{eqnarray*}

\noindent  We write, for each $x\in \mathbb{R}$,

$$
F_{X_1}(x)=F_{(X_1,X_2)}(x,+\infty).
$$

\bigskip \noindent  The same thing could be done for the $X_2$. We may now introduce the following terminology.\\

\bigskip \noindent \textbf{Definition}. $F_{(X_1,X_2)}$ is called the joint \textsl{cdf} of the ordered pair $(X_1,X_2)$. $F_{X_1}$ and $F_{X_2}$ are called the marginal \textsl{cfd}'s of the couple. The marginal \textsl{cdf}'s may be computed directly but they also may be derived from the joint \textsl{cdf} by

$$
F_{X_1}(x_1)=F_{(X_1,X_2)}(x_1,+\infty)  \text{ and } F_{X_2}(x_2)=F_{(X_1,X_2)}(+\infty,x_2), \ \ (x_1,x_2) \in \mathbb{R}^2.
$$

\bigskip \noindent The extension to higher dimensions is straightforward. Let $X \ : \ (\Omega ,\mathcal{A},\mathbb{P}) \rightarrow \mathbb{R}^d$ be a random vector with  $X^t=(X_1,...,X_d)$.\\

\noindent (i) Each marginal \textsl{cdf} $F_{X_i}$, $1 \leq i \leq d$, is obtained from the joint \textsl{cdf} $F_X=:F_{(X_1,...,X_d)}$ by

$$
F_{X_i}(x_i)=F_{(X_1,...,X_d)}\left( +\infty, ..., +\infty,\underset{i-th \ argument}{\underbrace{x_i}}, +\infty, ..., +\infty\right), \ x_i \in \mathbb{R},
$$

\bigskip \noindent or

$$
F_{X_i}(x_i)=\lim_{(\forall j\in\{1,...,d\}\setminus \{i\}, \ x_j\uparrow +\infty)} F_{(X_1,...,X_d)}(x_1,...,x_d), \ \ x_i \in \mathbb{R}.
$$

\bigskip \noindent (ii) Let $(X_{i_1}, ...,X_{i_r})^t$ be a  sub-vector of $X$ with $1\leq r <d$, $1\leq i_1 <i_2 <...< i_r)$. Denote $I=\{i_1,...,i_r\}$, the marginal \textsl{cfd} of $(X_{i_1}, ...,X_{i_r})$ is given by

$$
F_{(X_{i_1}, ...,X_{i_r})}\left(x_{i_1}, ...,x_{i_r}\right)=\lim_{\forall j\in\{1,...,d\}\setminus I, \ x_j\uparrow +\infty} F_{(X_1,...,X_d)}(x_1,...,x_d), \ (x_{i_1}, ...,x_{i_r}) \in \mathbb{R}^r.
$$

\bigskip \noindent (iii) Let $X^{(1)}=(X_1,...,X_r)^t$ and $X^{(2)}=(X_{r+1},...,X_b)^t$ be two sub-vectors which partition $X$ into consecutive blocs. The marginal \textsl{cdf}'s of $X^{(1)}$ and $X^{(2)}$ are respectively given by

$$
F_{X^{(1)}}(x)=F_{(X_1,...,X_d)}\left(x, \underset{(d-r) \text{ times}}{\underbrace{+\infty, ...,+\infty}}\right), \ x \in \mathbb{R}^{r}
$$

\noindent and

$$
F_{X^{(2)}}(y)=F_{(X_1,...,X_d)}\left(\underset{r \text{ times}}{\underbrace{+\infty, ...,+\infty}}, y\right), \ y \in \mathbb{R}^{d-r}.
$$

\bigskip \noindent After this series of notation, we have this important theorem concerning a new characterization of the independence.

\begin{theorem} \label{proba_02_rv_th03}  
Let $X \ : \ (\Omega ,\mathcal{A},\mathbb{P}) \rightarrow \mathbb{R}^d$ be a random vector. Let us adopt the notation above. The following equivalences hold.\\

\noindent (i) The margins $X_i$, $1\leq i \leq d$ are independent if and only if the joint \textsl{cdf} of $X$ is factorized in the following way :

$$
\forall (x_1,...,x_d) \in \mathbb{R}^d, \ F_{(X_1,...,X_d)}(x_1,...,x_d)=\prod_{j=1}^{d} F_{X_i}(x_i). \ (FLM01)
$$

\bigskip \noindent (i) The two marginal vectors $X^{(1)}$ and $X^{(2)}$ are independent if and only if the joint \textsl{cdf} of $X$ is factorized in the following way : for 
$(x^{(1)},x^{2)}) \in \mathbb{R}^d$, we have 

$$
\ F_{(X_1,...,X_d)}(x^{(1)},x^{2)})=F_{X^{(1)}}(x^{(1)}) F_{X^{(2)}}(x^{(2)}). \ \ (FLM02)
$$
\end{theorem}

\bigskip \noindent \textbf{Proof}. This important characterization follows as a simple result of Measure Theory and Integration. The proof of the two points are very similar. So, we only give the proof
of the first one.\\

\noindent Suppose that the components of $X$ are independent. By Theorem \ref{proba_02_intotp_th02} in Section \ref{proba_02_intotp_03} in Chapter \ref{proba_02_intotp}, we have for any 
$(x_1,...,x_d) \in \mathbb{R}^d$,

\begin{eqnarray*}
F_{(X_1,...,X_d)}(x_1,...,x_d)&=&\mathbb{P}(X_1\leq x_1,...,X_d\leq x_d)=\mathbb{E}\biggr(\prod_{j=1}^{d} 1_{]-\infty,x_j]}(X_j)\biggr)\\
&=& \prod_{j=1}^{d} \mathbb{E} \biggr(1_{]-\infty,x_j]}(X_j)\biggr)=\prod_{j=1}^{d} F_{X_j}(x_j).
\end{eqnarray*}
 
\noindent Conversely, if Formula (FLM01) holds, the Factorization Formula (FACT02) in Part (10.03) in Doc 10-01 in Chapter 11 in \cite{ips-mestuto-ang} of this series, we have :  for any $a=(a_{1},...,a_{k})\leq b=(b_{1},...,b_{k})$,

$$
\Delta_{a,b}F_X= \prod_{1 \leq i \leq k} (F_{X_i}(b_i)-F_{X_i}(a_i)).
$$

\bigskip \noindent By using the Lebesgue-Stieljes measures and exploiting the product measure properties, we have for any $(a,b)\in \mathbb{R}^2$, $a\leq b$,

$$
\mathbb{P}_X(]a,b])= \prod_{1 \leq i \leq k} \mathbb{P}_{X_i}(]a_i,b_i])=\biggr(\otimes_{j=1}^{d} \mathbb{P}_{X_i}\biggr)(]a,b]).
$$

\bigskip \noindent So the probability measures $\mathbb{P}_X$ and $\otimes_{j=1}^{d} \mathbb{P}_{X_i}$ coincide on the $\pi$-system of rectangles of the form $]a,b]$ which generates 
$\mathcal{B}(\mathbb{R}^d)$. Hence they simply coincide. Thus the components of $X$ are independent.\\

\noindent One handles the second point similarly by using Formula (FACT05) in the referred book at the same part, in the same document and the same section.\\

\bigskip \noindent  \textbf{(c5) How Can we Define a Random Variable Associated to a \textit{Cdf}}. \label{kolmconst_02}\\

\noindent As on $\mathbb{R}$, the Kolmogorov construction on $\mathbb{R}^d$, $d\geq 2$, is easy to perform.\\

\noindent For any \textsl{cdf} $F$ on $\mathbb{R}^d$, we may define the Lebesgue-Stieljes measure $\mathbb{P}$ on $(\overline{\mathbb{R}}^d ,\mathcal{B}_{\infty}(\overline{\mathbb{R}}^d)$ defined by

$$
\mathbb{P}(]y,x])=\Delta_{y,x}F, \ \ (y,x) \in (\mathbb{R}^d)^2, \ y \leq x. \ (LS21)
$$

\bigskip \noindent  Now take $\Omega=\mathbb{R}^d$, $\mathcal{A}=\mathcal{B}(\mathbb{R}^d)$ and let $X : (\Omega ,\mathcal{A},\mathbb{P}) \mapsto  \mathbb{R}^d$ be the identity function
$$
\forall \omega \in \Omega, \ X(\omega)=\omega.
$$

\bigskip \noindent  Thus we have :

$$
\forall x \in \mathbb{R}^d, \ F(x)=\mathbb{P}(]-\infty,x]). \ (LS22)
$$

\bigskip \noindent \textbf{Particular case}. In may situations, the above construction is stated as following : let $n\geq 1$ and $F_1$, $F_2$, ...., $F_n$ be $n$ \textsl{cdf}'s respectively defined on
$\mathbb{R}^{d_i}$, $d_i\geq 1$. Can we construct a probability space $(\Omega ,\mathcal{A},\mathbb{P})$ holding $n$ independent random vectors $X_1$, ..., $X_n$ such that for any 
$1\leq i \leq n$, $F_{X_i}=F_i$.\\

\noindent The answer is yes. It suffices to apply the current result to the $cdf$ $F$ defined on $\mathbb{R}^d$, with $d=d_1+...+d_n$ and defined as follows :

$$
\forall (x_1,...,x_n)^t \prod_{j=1}^{d} \mathbb{R}^{d_j}, \ F(x_1,...,x_n)=\prod_{j=}^{d} F_j(x_j).
$$

\bigskip \noindent Using Formula (FACT05) in Part (10.03) in Doc 10-01 in Chapter 11 in \cite{ips-mestuto-ang} of this series, we see that $F$ is a \textsl{cdf}. We may consider the identity function on $\mathbb{R}^d$ as above, form $X_1$ by taking the first $d_1$ components, $X_2$ by the next $d_2$ components, ..., $X_d$ by the last $d_n$ components. These subvectors are independent and respectively have the \textsl{cdf}'s $F_1$, ..., $F_n$.\\

\section[Laws of Random vectors]{Probability Laws and Probability Density Functions of Random vectors} \label{proba_02_rv_sec_05}

\bigskip \noindent  Throughout this section we deal with random vectors, like the $d$-random vector ($d\geq 1$) 
$$
X : (\Omega ,\mathcal{A},\mathbb{P}) \mapsto  (\overline{\mathbb{R}}^d, \mathcal{B}_{\infty}(\overline{\mathbb{R}}^d),
$$

\bigskip \noindent  with $X^t=(X_1, X_2, \cdots, X_d)$.\\

\bigskip \noindent  \textbf{ A- Classification of Random vectors}.\\

\noindent \textit{(a) Discrete Probability Laws}.\\

\noindent \textbf{Definition}. The random variable $X$ is said to be discrete if it takes at most a countable number of values in $\overline{\mathbb{R}}$ denoted $\mathcal{V}_X=\{x^{(j)}, j\in J\}$,
$\emptyset \neq J \subset \mathbb{N}$.\\

\noindent \textbf{NB}. On $\mathbb{R}$, we denote the values taken by such a random variable by sub-scripted sequences $x_j$, $j\in J$. In $\mathbb{R}^d$, $d\geq 2$, we use super-scripted sequences 
in the form $x^{(j)}$, $j\in J$, to avoid confusions with notation of components or powers.\\

\noindent Next, we give a set of facts from which we will make a conclusion on how to work with such random variables.\\

\noindent We already know from Measure Theory that $X$ is measurable (Se Chapter 4, Doc 08-03, Criterion 4) if and only if
$$
\forall j \in J, (X=x^{(j)}) \in \mathcal{A}.
$$

\noindent  Besides, we have for any $B  \in \mathcal{B}_{\infty}(\overline{\mathbb{R}}^d)$,

$$
(X \in B)= \sum_{j \in J, x^{(j)} \in B} (X=x^{(j)}). \ (DD01)
$$

\bigskip \noindent  Now, we clearly have 

$$
\sum_{j \in J} \mathbb{P}(X=x^{(j)})=1. \ (DD02)
$$ 

\bigskip \noindent  From (DD01), the probability law $\mathbb{P}_X$ of $X$ is given by

$$
\mathbb{P}_X(B)=\sum_{j \in J, x^{(j)} \in B} \mathbb{P}(X=x^{(j)}), \ (DD03)
$$

\bigskip \noindent  for any $B \in \mathcal{B}_{\infty}(\overline{\mathbb{R}}^d)$. Let us denote the function defined on $\mathcal{V}_X$ by

$$
\mathcal{V}_X \in x \mapsto f_X(x)=\mathbb{P}_X(\{x\})=\mathbb{P}(X=x^{(j)}).
$$

\bigskip \noindent Next, let us consider the counting measure $\nu$ on $\mathbb{R}^d$ with support $\mathcal{V}_X$. Formulas (DD02) and (DD03) imply that 

$$
\int f_X \ d\nu=1. \ (RD01)
$$

\bigskip \noindent  and for any $B \in \mathcal{B}_{\infty}(\overline{\mathbb{R}}^d)$, we have

$$
\int_{B} d\mathbb{P}_X = \int_B f_X \ d\nu. \ (RD02) 
$$

\bigskip \noindent  We conclude that $f_X$ is the Radon-Nikodym derivative of $\mathbb{P}_X$ with respect to the $\sigma$-finite measure $\nu$. Formula (RD02) may be written in the form

$$
\int h d\mathbb{P}_X = \int h f_X \ d\nu. \ (RD03) 
$$

\bigskip \noindent  where $h=1_B$. By using the four steps method of the integral construction, Formula (RD03) becomes valid whenever $\mathbb{E}h(X)=\int h \ d\mathbb{P}_X$ make senses.\\

\noindent We may conclude as follows.\\

\noindent \textbf{Discrete Probability Laws}.\\

\noindent If $X$ is discrete, that is, it takes a countable number of values in $\overline{\mathbb{R}}^d$ denoted $\mathcal{V}_X=\{x^{(j)}, j\in J\}$, its probability law $\mathbb{P}_X$ is also said to be discrete. It has a probability density function \textit{pdf} with respect to the counting measure on $\mathbb{R}^d$ supported by $\mathcal{V}_X$ and defined by

$$
f_X(x)=\mathbb{P}(X=x), \ x \in \overline{\mathbb{R}}^d,
$$

\bigskip \noindent  which satisfies

$$
f_X(x^{(j)})=\mathbb{P}(X=x^{(j)}) \text{ for } j\in J \text{ and } f_X(x)=0 \text{ for } x \notin \mathcal{V}_X.
$$

\bigskip \noindent  As a general rule, integrating any measurable function $h : \mathcal{B}_{\infty}(\overline{\mathbb{R}}^d) \rightarrow \mathcal{B}_{\infty}(\overline{\mathbb{R}})$ with respect to the probability law $\mathbb{P}_X$ is performed through the \textit{pdf} $f_X$ in the Discrete Integral Formula \label{dif1}

$$
\mathbb{E}h(X)=\int h f_X d\nu= \sum_{j \in J} h(x_j) f_X(x^{(j)}). \ \ (DIF1)
$$

\bigskip \noindent  which becomes for $h=1_B$, $B \in \mathcal{B}_{\infty}(\overline{\mathbb{R}}^d)$,

$$
\mathbb{P}_X(B)=\mathbb{P}(X\in B) = \sum_{j \in J, x^{(j)}}  f_X(x^{(j)}). (DIF2)
$$

\bigskip \bigskip \bigskip \noindent  Some authors name \textit{pdf}'s with respect to counting measures as \textit{mass pdf's}. For theoretical purposes, they are Radon-Nikodym derivatives.

\bigskip \noindent  \textbf{ (b) Absolutely Continuous Probability Laws}.\\

\noindent \textbf{(b1) Lebesgue Measure on $\mathbb{R}^d$}.\\

\noindent We already have on $\overline{\mathbb{R}}^d$ the $\sigma$-finite Lebesgue measures $\lambda_d$, which is the unique measure defined by the values

$$
\lambda_d\biggr(\prod_{i=1}{2} ]a_i,b_i]\biggr)= \prod_{i=1}^{d} (b_i-a_i), \ (LM01)
$$

\bigskip \noindent  for any points $a=(a_{1},...,a_{d})^t \leq b=(b_{1},...,b_{d})^t$ of $\mathbb{R}^{d}$. This formula also implies

$$
\lambda_d\biggr(\prod_{i=1}^{d} ]a_i,b_i]\biggr)= \prod_{i=1}^{d} \lambda_1(]a_i, b_i]), \ (LM02)
$$

\bigskip \noindent  Let us make some Measure Theory reminders. Formula (LM02) ensures that $\lambda_d$ is the product measure of the Lebesgue measure $\lambda_1=\lambda$, that is

$$
\lambda_d=\lambda^{\otimes d}.
$$

\bigskip \noindent  Hence, we may use Fubini's Theorm for integrating a measurable function function $h : \mathcal{B}_{\infty}(\overline{\mathbb{R}}^d) \rightarrow \mathcal{B}_{\infty}(\overline{\mathbb{R}})$ through the formula

$$
\int h d\lambda_d =\int d\lambda(x_1) \int ....\int d\lambda(x_{d-1}) \int h(x_1,...,x_d) f_X(x_1,...,x_d) d\lambda(x_{d}),
$$

\bigskip \noindent  when applicable (for example, when $h$ is non-negative or $h$ is integrable).\\

\noindent \textbf{(b2) Definition}.\\

\noindent The probability Law $\mathbb{P}_X$ is said to be absolutely continuous if it is continuous with respect to $\lambda_d$. By extension, the random variable itself is said to be 
absolutely continuous.\\

\noindent \textbf{NB}. It is important to notice that the phrase \textbf{absolutely continuous} is specifically related to the continuity with respect to Lebesgue measure.\\

\bigskip \noindent  In the rest of this Point (b), we suppose that $X$ is absolutely continuous.\\

\noindent \textbf{(b3) Absolutely Continuous \textsl{pdf}'s}.\\

\noindent By Radon-Nikodym's Theorem, there exists a Radon-Nikodym derivative denoted $f_X$ such that for any $B \in \mathcal{B}_{\infty}(\overline{\mathbb{R}}^d)$,

$$
\int_B d\mathbb{P}_X = \int_B f_X \ d\lambda_d.
$$

\bigskip \noindent  The function $f_X$ satisfies

$$
f_X \geq 0 \text{ and } \int_{\mathbb{R}} f_X \ d\lambda_d=1.
$$

\bigskip \noindent  Such a function is called a \textsl{pdf} with respect to the Lebesgue measure. Finally, we may conclude as follows.\\

\noindent As a general rule, integrating any measurable function $h : \mathbb{R}^d \rightarrow \overline{\mathbb{R}}$ with respect to the probability law $\mathbb{P}_X$, which is absolutely continuous, is performed through the \textit{pdf} $f_X$ with the Absolute Continuity Integral Formula \label{acif}

$$
\mathbb{E}h(X)=\int h f_X \ d\lambda_d. \ (ACIF)
$$

\bigskip \noindent  Since $\lambda_d$ is the product of the Lebesgue measure on $\mathbb{R}$ $d$ times, we may use Fubini's Theorem when applicable to have

$$
\mathbb{E}h(X)= \int d\lambda_1(x_1) \int ....\int d\lambda_1(x_{d-1}) \int h(x_1,...,x_d) f_X(x_1,...,x_d) d\lambda_1(x_{d}).
$$

\bigskip
\noindent In particular, the \textit{cdf} of $X$ becomes

\begin{eqnarray*}
F_X(x)&=&\int_{-\infty}{x_d} d\lambda(x_1) \int_{-\infty}{x_2} \\
&....& \int_{-\infty}{x_{d-1}} d\lambda(x_{d-1}) \int_{-\infty}{x_d} h(x_1,...,x_d) f_X(x_1,...,x_d) d\lambda(x_{d})
\end{eqnarray*}

\bigskip \noindent  for any $x=(x_1,...,x_d)^t \in \mathbb{R}^d$.\\

\bigskip \noindent  \textbf{(b4) Criterion for Absolute Continuity from the Cdf}.\\

\noindent In practical computations, a great deal of Lebesgue integrals on $\mathbb{R}$ are Riemann integrals. Even integrals with respect to the multidimensional Lebesgue Measure can be multiple Riemann ones. But we have to be careful for each specific case (See Points (b5) and (b6) below).\\
  
\noindent Let be given the \textsl{dcf} $F_X$ of a random vector, the absolute continuity of $X$ would give for any $x \in \mathbb{R}^d$ \label{ac01}

$$
F_{X}(x)=\int_{-\infty }^{x_{1}} d\lambda(x_{1})\int_{-\infty }^{x_{2}}d\lambda(x_{2})...d\lambda(x_{d-1})\int_{-\infty }^{x_{k}}f_{X}(x_{1},...,x_{d}) \ d\lambda(x_{d}). \ (AC01)
$$

\bigskip \noindent  If $f_X$ is locally bounded and locally Riemann integrable (LLBRI), we have

$$
f_{X}(x_{1},x_{2},...,x_{k})=\frac{\partial ^{k}F_{X}(x_{1},x_{2},...,x_{k})}{\partial x_{1}\partial x_{2}...\partial x_{x}}, \ (AC02)
$$

\bigskip \noindent  $\lambda_d$-\textit{a.e.}. (See Points (b5) and (b6) below for a more detailed explanation of \textit{LLBRI} functions and for a proof).\\

\noindent From a computational point of view, the above Formula quickly helps to find the \textit{pdf}, if it exists.\\

\noindent \noindent \textbf{(b5) Cautions to be taken when replacing Lebesgue integral by Riemann ones}. \label{LLBRI}\\

\noindent Let us consider that we are on $\mathbb{R}$ and let $X$ be a real random variable with an absolutely  \textit{pdf} $f$. For any measurable function $h$ from $\mathbb{R}$ to $\mathbb{R}$, the expectation

$$
\mathbb{E}(h(X))=\int_{\mathbb{R}} h(x) f(x) \ d\lambda(x), \ (EL)
$$

\bigskip \noindent   is defined with respect to the Lebesgue measure. It happens that for computation such an integral, we lean to use the improper Riemann integral

$$
\mathbb{E}(h(X))=\int_{-\infty}^{+\infty} h(x) f(x) \ dx. \ (ER)
$$

\bigskip \noindent   Although this works for a lot of cases, we cannot use the just mentioned formula without a minimum of care, since in Riemann
integration we may have that $\int_{\mathbb{R}} h(x) f(x) \ dx$ is finite and $\int_{\mathbb{R}} |f(x)| \ dx$ infinite, a situation that cannot occur with Lebesgue integration.\\

\noindent We may use the results of Doc 06-07 in Chapter 7 in \cite{ips-mestuto-ang} of this series to recommend the following general rule that we will follow in this book.\\

\noindent Let us suppose that the function $h f$ is \textsl{LLBRI} (implying that $hf$ is $\lambda$-\textit{a.e.} continuous on $\mathbb{R}$). We have :\\

\noindent (a) If $\mathbb{E}(h(X))$ exists and is finite, then Formula (ER) holds as an improper Riemann integral (as an application of the Dominated Convergence Theorem), that is

$$
\mathbb{E}(h(X))=\lim_{n\rightarrow +\infty} \int_{a_n}^{b_n} h(x) f(x) \ dx, \ (ER02)
$$

\bigskip \noindent   for any sequence $(a_n,b_n)_{n\geq 0}$ such that $(a_n,b_n)\rightarrow (-\infty, \ +\infty)$ as $n\rightarrow +\infty$. In such a case, we may chose a particular alike sequence to compute 
$\mathbb{E}(h(X))$.\\

\noindent To check whether $\mathbb{E}(h(X))$ is finite, we may directly use Riemann integrals (which are based on the Monotone Convergence Theorem)

$$
\mathbb{E}(h^-(X))=\int_{-\infty}^{+\infty}  (hf)^+(x) f(x) \ dx \ \ (ENPa) 
$$

\noindent and

$$
\mathbb{E}(h^-(X))=\int_{-\infty}^{+\infty}  (hf)^-(x) f(x) \ dx, (ENPb) 
$$

\bigskip \noindent   and apply the classical Riemann integrability criteria.\\

\noindent (b) If the Riemann improper integral of $|hf|$ exists and is finite, then the Lebesgue integral of $hf$ exists (by using the Monotone Convergence Theorem on the positive and negative parts) and Formula (ER) holds.\\

\noindent (c) Even if $\mathbb{E}(h(X))$ exists and is infinite, Formula (ER) still holds, by using the Monotone Convergence Theorem on the positive and negative parts and exploiting Formula (ENP).\\

\noindent Finally, such results are easily extended in dimension $d\geq 2$, because of the Fubini's integration formula.\\

\bigskip \noindent  \textbf{(b6) Back to Formula (AC01)}.\\

\noindent \textbf{Dimension one}. If $f$ is \textsl{LLBRI}, we surely have that $f$ is $\lambda$-\textit{a.e.} continuous and we may treat the integrals $\int_{-\infty}^{x} f_X(t) \ d\lambda(t)$ as a Riemann ones. By the known results for indefinite Riemann integrals, we have

$$
\biggr( \forall x \in \mathbb{R}, F_X(x)=\int_{-\infty}^{x} f_X(t) \ d\lambda(t) \biggr) \Leftrightarrow \frac{dF_X}{dx}=f_X \ \lambda-\textit{a.e.}.
$$

\bigskip \noindent   Remark that the constant resulting in the solution of the differential equation in the right-hand assertion is zero because of $F_X(-\infty)=0$.\\

\bigskip 
\noindent \textbf{Dimension $d\geq 2$}. Let $d=2$ for example. Let $f_X$ be \textit{LLBRI}. By Fubini's theorem,

$$
\forall (x,y) \in \mathbb{R}^2, F_X(x,y)=\int_{-\infty}^{x} d\lambda(s) \biggr( \int_{-\infty}^{y} f_X(s,t) \ d\lambda(t) \biggr).
$$

\bigskip \noindent   The function,  

$$
t \rightarrow \int_{-\infty}^{y} f_X(s,t) \ d\lambda(t)
$$

\bigskip \noindent is bounded (by the unity) and continuous. By, returning back to Riemann integrals, we have

$$
\forall (x,y) \in \mathbb{R}^2, \ \ \frac{\partial F_X(x,y)}{\partial x}=\int_{-\infty}^{y} f_X(x,t) \ d\lambda(t).
$$ 

\bigskip \noindent By applying the results for dimension one to the partial function $f_X(x,t)$, for $x$ fixed, which is (LBLI), we get

$$
\forall (x,y) \in \mathbb{R}^2, \frac{\partial^2 F_X(x,y)}{\partial y \partial x}=\int_{-\infty}^{y} f_X(x,t)\  \lambda-a.e.
$$ 

\bigskip \noindent   The order of derivation may be inverted as in the Fubini's Theorem.\\

\noindent The general case $d\geq 2$ is handled by induction. $\blacksquare$\\

\bigskip \noindent  \textit{(c) General case}.\\

\noindent Let us be cautious! Later, we will be concerned by practical computations and applications of this theory. We will mostly deal with discrete or absolutely continuous random variables. But, we should be aware that these kind of probability laws form only a small part of all the possibilities, as we are going to see it.\\

\noindent By the Lebesgue Decomposition Theorem (Doc 08-01, Part III, Point (05-06), Chapter 9), there exists a unique decomposition of $\mathbb{P}_X$ into an absolutely continuous measure $\phi_{ac}$, associated to a non-negative Radon-Nikodym $f_{r,X}$ and $\lambda_d$-singular measure $\phi_{s}$, that is, for any $B \in \mathcal{B}_{\infty}(\overline{\mathbb{R}}^d)$,

$$
\mathbb{P}_X(B)=\int_B f_{r,X} \ d\lambda_d+ \phi_s(B).
$$ 

\bigskip \noindent  The $\lambda$-singularity of $\phi_{s}$ means that there exists a $\lambda_d$-null set $N$ such that for all $B \in \mathcal{B}_{\infty}(\overline{\mathbb{R}}^d)$,

$$
\phi_s(B)=\phi_s(B\cap N)
$$  

\bigskip \noindent  Suppose that none of $\phi_s$ and $\phi_s$  is the null measure. If $N$ is countable that is $N$ may be written as $N=\{x^{(j)}, j\in J\}$, $\emptyset \neq J \subset \mathbb{N}$, the measure is discrete and for any $B \in \mathcal{B}_{\infty}(\overline{\mathbb{R}}^d)$, we have

$$
\phi_s(B)=\sum_{j\in J} f_{d,X}(x^{(j)}),
$$  

\noindent  where 

$$
f_{d,X}(x^{(j)})=\phi_s(x^{(j)}), \ \ j\in J.
$$

\bigskip \noindent  We have $0<a=\phi_{ac}(\overline{\mathbb{R}}^d), \ b=\phi_{s}(\overline{\mathbb{R}}^d)\leq 1$, and $a+b=1$. Let us denoting by $\nu$ the counting measure with support $N$. Then $f_X^{(1)}=f_{r,X}/a$ is an absolutely continuous \textit{pdf} and $f_X^{(2)}=f_{d,X}/b$ is a discrete \textsl{pdf} and we have for all $B \in \mathcal{B}_{\infty}(\overline{\mathbb{R}}^d)$,

$$
\mathbb{P}_X(B)= a \int_B  f_X^{(1)} \ d\lambda_d + (a-1) \int_B f_X^{(2)} \ d\nu.
$$

\bigskip \noindent  Hence, $\mathbb{P}_X$ is mixture of two probability laws, the first being absolutely continuous and the second being discrete.\\

\noindent We may be more precise in dimension one.\\

\noindent \textbf{More detailed decomposition on $\mathbb{R}$}. We already saw that a real \textsl{cdf} $F $may be decomposed into two \textsl{df}'s :

$$
F=F_c +F_d,
$$

\bigskip \noindent  where $F_c$ is continuous and $F_d$ is discrete. Surely, the Lebesgue-Stieljes measure associated with $F_d$, denoted by $\phi_d$, is discrete. 
The Lebesgue-Stieljes measure associated with $F_c$, denoted by $\phi_c$, may decomposed as above into

$$
\phi_c=\phi_{ac} + \phi_s
$$ 

\bigskip \noindent  where $\phi_{ac}$ is absolutely continuous and  $\phi_s$ is singular. Since $F_c(-\infty)=F_d(-\infty)=0$, we may go back to the \textsl{df}'s to have :

$$
F=F_{ac}+F_c+F_d,
$$
 
\bigskip \noindent  where $F_{ac}$ is \textsl{df} of measure absolutely continuous, $F_d$ is a discrete \textsl{df} and $F_c$ is a continuous and, unless it is equal to the null measure, is neither discrete nor absolutely continuous.\\

\noindent This fact is obvious since $F_c$ is continuous and cannot be discrete. Also, it is singular and cannot be absolutely continuous.\\

\noindent We have the following conclusion.\\

\noindent \textbf{Position of any probability law with respect to the Lebesgue measure}. Any probability law is a mixture of an absolutely continuous probability measure $\mathbb{P}_{ac,X}$, associated to a $pdf$ $f_{ac,X}$, a discrete distribution probability measure $\mathbb{P}_{d,X}$, which is a $\lambda$-singular measure $\mathbb{P}_{d,X}$ which has a countable strict support $\mathcal{V}_{d,X}$ and of a $\lambda$-singular probability measure $\mathbb{P}_{c,X}$ which has a non-countable $\lambda$-null set support, respectively associated to $p_1\geq 0$, $p_2\geq 0$ and $p_2\geq 0$, with $p_1+p_2+p_3=1$, such that

$$
\mathbb{P}_{X}=p_1 \mathbb{P}_{ac,X}+ p_2 \mathbb{P}_{d,X}+ p_3 \mathbb{P}_{c,X}.
$$  
 
\bigskip \noindent  The probability measures are respectively associated to the \textsl{df}'s $F_{ac}$, $F_d$, $F_c$ so that we have

$$
F_X=F_{ac} +F_d + F_c,
$$

$$
\frac{dF_{ac}(x)}{dx}=f_{ac,X}, \ \lambda-a.e,
$$

\bigskip \noindent  $\mathcal{V}_{d,X}$ is the set of discontinuity points of $F$, and $F_c$ is continuous but not $\lambda$-\textit{a.e.} differentiable.\\

\noindent By \textit{strict countable support of $\mathbb{P}_{d,X}$}, we mean a support such that for any point $x$ in, we have $\mathbb{P}_{d,X}(\{x\})>0$.\\

\noindent \textbf{Warning}. If the decomposition has more that one term, the corresponding functions among $F_{ac}$, $F_d$ and $F_c$ are not \textit{cdf}'s but only \textsl{df}'s.\\  

\bigskip \noindent \textbf{(b7) Marginal Probability Density functions}.\\

\noindent Let us begin, as usual, by the simple case where $d=2$. Let $X \ : \ (\Omega ,\mathcal{A},\mathbb{P}) \rightarrow \mathbb{R}^2$ be a random couple, with $X^t=(X_1,X_2)$. Let us suppose that 
$X$ has a \textit{pdf} $f_{(X_1,X_1)}$ with respect to a $\sigma$-finite product measure $m=m_1 \otimes m_2$ on $\mathbb{R}^2$. Let us show that each $X_i$, $i \in \{1,2\}$, has a \textsl{pdf} with respect to $m$. We have, for any Borel set $B$,

\begin{eqnarray*}
\mathbb{P}(X_1 \in B)&=&\mathbb{P}((X_1,X_2) \in B\times \mathbb{R})\\
&=& \int 1_{B\times \mathbb{R}} f_{(X_1,X_2)}(x,y) \ dm(x,y)\\
&=&\int_B \biggr(\int_{\mathbb{R}} f_{(X_1,X_2)}(x,y) \ dm_2(y)\biggr) \ dm_1(x).
\end{eqnarray*}

\bigskip \noindent By definition, the function

$$
f_{X_1}(x) =\int_{\mathbb{R}} f_{(X_1,X_1)}(x,y) \ dm_2(y), \ m-a.e \ in \ x \in \mathbb{R},
$$ 

\bigskip \noindent is the \textsl{pdf} of $X$ with respect of $m_1$, named as the \textit{marginal} \textsl{pdf} of $X_1$. We could do the same for $X_2$. We may conclude as follows.

\bigskip \noindent \textbf{Definition}. Suppose that the random order pair  $X^t=(X_1,X_2)$ has a  \textit{pdf} $f_{(X_1,X_1)}$ with respect to a $\sigma$-finite product measure $m=m_1 \otimes m_2$ on $\mathbb{R}^2$. Then each $X_i$, $i \in \{1,2\}$, has the marginal \textsl{pdf}'s $f_{X_i}$ with respect to $m_i$, and

$$
f_{X_1}(x) =\int_{\mathbb{R}} f_{(X_1,X_2)}(x,y) \ dm_2(y), m_1-a.e. \in \ x \in \mathbb{R}
$$ 

\bigskip \noindent  and

$$
f_{X_2}(x) =\int_{\mathbb{R}} f_{(X_1,X_2)}(x,y) \ dm_1(x), \ m_2-a.e. \in x \in \mathbb{R}.
$$

\bigskip \noindent The extension to higher dimensions is straightforward. Let $X \ : \ (\Omega ,\mathcal{A},\mathbb{P}) \rightarrow \mathbb{R}^d$ be a random vector with  $X^t=(X_1,...,X_d)$. Suppose that
$X$ has a \textsl{pdf} $f_{(X_1,...,X_d)}$ with respect to a $\sigma$-finite product measure $m=\otimes_{j=1}^{d} m_j$.\\ 

\bigskip  \noindent (i) Then each $X_j$, $j \in \{1,d\}$, has the marginal \textsl{pdf}'s $f_{X_j}$ with respect to $m_j$ given $m_i$-\textit{a.e.}, for $x \in \mathbb{R}$, by

$$
f_{X_j}(x) =\int_{\mathbb{R}^{d-1}} f_{(X_1,...,X_d)}(x_1,...,x_d) \ d\biggr(\otimes_{i\leq i \leq d,\ i\neq j} m_i\biggr)(x_1,...,x_{j-1},x_{j+1},...,x_d).
$$ 

\bigskip \noindent (ii) Let $(X_{i_1}, ...,X_{i_r})^t$ be a  sub-vector of $X$ with $1\leq r <d$, $1\leq i_1 <i_2 <...< i_r)$. Denote $I=\{i_1,...,i_r\}$, the marginal \textsl{pdf} of $(X_{i_1}, ...,X_{i_r})$ with respect to $m=\otimes_{i=1}^{r} m_{i_j}$ is given for $(x_1,...,x_r) \in \mathbb{R}^r$ by

\begin{eqnarray*}
&&f_{(X_{i_1}, ...,X_{i_r})}(x_1,...,x_r) =\\
&&\int_{\mathbb{R}^{d-r}} f_{(X_1,...,X_d)}(x_1,...,x_d) \ d\biggr(\otimes_{1\ i \leq d,\ i\notin I} m_i\biggr)(x_j, \j\in \{1,...,n\}\setminus I).
\end{eqnarray*}

\bigskip \noindent Let $X^{(1)}=(X_1,...,X_r)^t$ and $X^{(2)}=(X_{r+1},...,X_b)^t$ be two sub-vectors which partition $X$ into two consecutive blocs. Then  $X^{(1)}$ and $X^{(2)}$ have the \textit{pdf} 
$f_{X^{(1)}}$ and $f_{X^{(2)}}$ with respect to $\otimes_{j=1}^{r} m_j$ and $m=\otimes_{j=r+1}^{d} m_j$ respectively, and given for $x\in \mathbb{R}^{r}$ by

$$
f_{X^{(1)}}(x) =\int_{\mathbb{R}^{d-r}} f_{(X_1,...,X_d)}(x_1,...,x_d) \ d\biggr(\otimes_{r+1 \leq i \leq d} m_i\biggr)(x_{r+1},...,x_d), 
$$ 

\bigskip \noindent  and for $ x \in \mathbb{R}^{d-r}$ by

$$
f_{X^{(2)}}(x) =\int_{\mathbb{R}^{r}} f_{(X_1,...,X_d)}(x_1,...,x_d) \ d\biggr(\otimes_{1 \leq i \leq r} m_i\biggr)(x_1,...,x_r). 
$$ 

\bigskip \noindent After this series of notations, we have this important theorem for characterizing the independence.

\begin{theorem} \label{proba_02_rv_th04}
Let $X \ : \ (\Omega ,\mathcal{A},\mathbb{P}) \rightarrow \mathbb{R}^d$ be a random vector. Let us adopt the notation above. Suppose that we are given a $\sigma$-finite product measure $m=\otimes_{j=1}^{d} m_j$ on $\mathbb{R}^d$, and $X$ has a \textsl{pdf} $f_X$ with respect to $m$. We have the following facts.\\

\noindent (i) The margins $X_i$, $1\leq i \leq d$ are independent if and only if the joint \textsl{pdf} of $X$ is factorized in the following way :

$$
\forall (x_1,...,x_d) \in \mathbb{R}^d, \ f_{(X_1,...,X_d)}(x_1,...,x_d)=\prod_{j=1}^{d} f_{X_i}(x_i), \ m.a.e. \ (DLM01)
$$

\bigskip \noindent (i) The two marginal vectors $X^{(1)}$ and $X^{(2)}$ are independent if and only if the joint \textsl{pdf} of $X$ is factorized in the following way : for all $(x^{(1)},x^{2)}) \in \mathbb{R}^d$, 

$$
f_{(X_1,...,X_d)}(x^{(1)},x^{2)})=f_{X^{(1)}}(x^{(1)}) f_{X^{(2)}}(x^{(2)}), \ m.a.e. \  (FLM02)
$$
\end{theorem}

\bigskip  \noindent \textbf{Proof}. It will be enough to prove the first point, the proof of the second being very similar. Suppose that the $X_i$ are independent. It follows that for any Borel rectangle
$B=B_1 \times ... \times B_d$, we have

\begin{eqnarray*}
\mathbb{P}(X\in B)&=& \int_{B} f_X(x) \ dm(x)\\
&=&\mathbb{P}(X_1 \in B_1,...,X_d \in B_d)\\
&=&\prod_{j=1}^{d} \mathbb{P}(X_j \in B_j)\\
&=&\prod_{j=1}^{d} \int_{B_j} f_{X_j}(x_i) \ dm_j(x_j)\\
&=&\int_{B_1 \times ... \times B_d} \biggr(\prod_{j=1}^{d} f_{X_j}(x_i) \biggr) dm(x)
\end{eqnarray*}

\noindent Thus the two finite measures

$$
B \mapsto \int_{B} f_X(x) \ dm(x) \text{ and } B \mapsto \int_{B} \biggr(\prod_{j=1}^{d} f_{X_j}(x_i) \biggr) dm(x)
$$

\noindent coincide on a $\pi$-system generating the whole $\sigma$-algebra. Thus, they coincide. Finally, we get two finite indefinite integrals with respect to the same $\sigma$-finite measure $m$. By 
the Radon-Nikodym Theorem, the two Radon-Nikodym derivatives are equal $m$-\textit{a.e}.\\

\noindent Suppose now that Formula (DLM01) holds. Thanks to Fubini's Theorem, we readily get the factorization of the joint \textsl{cdf} and get the independence through Theorem \ref{proba_02_rv_th03}.\\

\section{Characteristic functions} \label{proba_02_rv_sec_06}

\bigskip \noindent  After the \textsl{cdf}'s, are going to see a second kind of characterization function for probability laws.\\

\noindent \textbf{I - Definition and first properties}.\\

\noindent It is important to say that, in this section, we only deal with finite components random vectors with values in spaces $\mathbb{R}^{d}$, $d\geq 1$, endowed with the Borel $\sigma$-algebra 
$\mathcal{B}(\mathbb{R}^{d})=\mathcal{B}(\mathbb{R})^{\otimes d}$ which is the product $\sigma$-algebra of $\mathcal{B}(\mathbb{R})$ $d$ times.\\

\noindent \textbf{(a) Characteristic function}.\\

\begin{definition} For any random variable $X:(\Omega ,\mathcal{A},\mathbb{P})\mapsto 
\mathbb{R}^{d}$, the function 
\begin{equation*}
u\mapsto \phi _{X}(u)=\mathbb{E}(e^{i<X,u>}),
\end{equation*}

\bigskip \noindent  is called the characteristic function of $X$. Here, $i$ is the complex number with positive imaginary part such that $i^2=-1$.
\end{definition}

\noindent This function always exists since we interpret the integral in the following way

\begin{equation*}
\mathbb{E}(e^{i<X,u>})=\mathbb{E}(\cos <X,u>)+i\text{ }\mathbb{E}(\sin<X,u>),
\end{equation*}

\bigskip \noindent  which is defined since the integrated real and imaginary parts are bounded.\\

\noindent The role played by the characteristic function in Probability Theory may also be played by  a few number of functions called moment generating functions. These functions do not always exist, and if they do, they may be defined only on a part of $\mathbb{R}^d$. The most used of them is defined as follows.\\

\noindent \textbf{(a) Moment Generated Function (\textit{mgf})}.\\

\noindent The following function 
\begin{equation*}
u\mapsto \varphi_{X}(u)=\mathbb{E}(e^{<X,u>}), \ u \in \mathbb{R}^d,
\end{equation*}

\bigskip \noindent  when defined on a domain $D$ of $\mathbb{R}^d$ containing the null vector as an interior point, is called the moment generating function (\textit{mfg}) of $X$.\\

\noindent \noindent If $\varphi_{X}$ exists on some domain $D$ to which zero is interior, we will prefer it to $\Phi_{X}(u)$, to avoid to use the \textit{complex} number $i$ involved in $\Phi_{X}$. Non-mathematician users of Probability Theory would like this.\\

\noindent Besides, we may find the characteristic function by using the moment generating functions as follows :
\begin{equation*}
\Phi_X (u)=\varphi_X(iu), u \in  \mathbb{R}^d.
\end{equation*}

\noindent  The characteristic function has these two immediate properties.

\begin{proposition}
For all $u\in \mathbb{R}^{d}$,

$$
\left\Vert \phi _{X}(u)\right\Vert \leq 1=\left\Vert \phi _{X}(0)\right\Vert.
$$

\bigskip \noindent Besides $\phi _{X}(u)$ is uniformly continuous at any point $u \in \mathbb{R}^d$.
\end{proposition}

\noindent This proposition is easy to prove. In particular, the second point is an immediate application to the Dominated Convergence Theorem.

\bigskip \noindent  Here are the :\\

\noindent \textbf{II - Main properties of the characteristic function}.

\begin{theorem} \label{proba02_rv_prop_cf} We have the following facts.\\

\noindent (a) Let $X$ be a random variable with value in $\mathbb{R}^{d}$, $A$ a $(k\times d)$-matrix of real scalars, $B$ a vector of $\mathbb{R}^{k}$. Then the characteristic function of  $Y=AX+B \in \mathbb{R}^{k}$ is given,  
\begin{equation*}
\mathbb{R}^{k}\ni u\mapsto \phi _{Y}(u)=e^{<B,u>}\phi _{X}(A^t u), \ u \in \mathbb{R}^{k}.
\end{equation*}

\bigskip \noindent   (b) Let $X$ and $Y$ be two independent random variables with values in $\mathbb{R}^{d}$, defined on the same probability space. The for any
$u\in \mathbb{R}^{d}$, we have

\begin{equation*}
\phi _{X+Y}(u)=\phi _{X}(u)\text{ }\times \phi _{Y}(u).
\end{equation*}

\bigskip \noindent (c) Let $X$ and $Y$ be two random variables respectively with values in $\mathbb{R}^{d}$ and in $\mathbb{R}^{k}$ and defined on the same probability measure. If the random variables $X$ and $Y$ are independent, then for any $u\in \mathbb{R}^{d}$ and for $v \in \mathbb{R}^{k}$, we have
\begin{equation}
\phi _{(X,Y)}(u,v)=\phi _{X}(u)\times \phi _{Y}(v).  \label{prodfc}
\end{equation}
\end{theorem}

\bigskip \noindent  Let us make some remarks before we give the proof of the theorem. In Part A, Section 3, Chapter 6, the characterization (c) was stated and quoted as (CI4), and admitted without proof. Here, the proof will be based on a characterization of product measure.\\

\noindent Point (c) provides a characterization of the independence between $X$ and $Y$. But the decomposition in Point (b) is not enough to ensure the independence. You may consult counter-examples book of \cite{stoyanov} or the monograph \cite{ips-probelem-ang} of this series, Part A, Section 3, Chapter 6, where is reported a counter-example from \cite{stoyanov}.\\

\bigskip \noindent  \textbf{Proof of Theorem \ref{proba02_rv_prop_cf}}.\\

\noindent Point (a). By definition, we have $<AX+B,u>=$ $^{t}(AX+B)u=$ $^{t}X(A^{T}u)+B^{T} u$. Hence, 
\begin{eqnarray*}
\phi _{AX+B}(u)&=&\mathbb{E}(e^{^{t}X(A^{t} u)+B ^{t}u})=e^{<B,u>}\times \mathbb{%
E}(e^{<X,A^{t}u)})\\
&=&e^{<B,u>}\phi _{X}(A^{t} u).
\end{eqnarray*}

\bigskip \noindent  Point (b). Let $X$ and $Y$ be independent. We may form $X+Y$ since they both have their values in $\mathbb{R}^{d}$, and they are defined on the same probability space. We have for any 
$u \in \mathbb{R}^{d}$,

\begin{equation*}
\phi _{X+Y}(u)=\mathbb{E}\left(e^{<X+Y,u>}\right)=\mathbb{E}\left(e^{<X,u>}e^{<Y,u>}\right)=%
\mathbb{E}\left(e^{<X,u)}\right)\times \mathbb{E}\left(e^{<Y,u)}\right).
\end{equation*}

\bigskip \noindent  Point (c). Let $X$ and $Y$ be two independent random variables with values in $\mathbb{R}^{d}$ and $\mathbb{R}^{k}$. Let $u$ and $v$ be two respectively elements of $\mathbb{R}^{d}$ and $\mathbb{R}^{k}$. We have 
\begin{equation*}
<\left( 
\begin{array}{c}
X \\ 
Y%
\end{array}%
\right) ,\left( 
\begin{array}{c}
u \\ 
v%
\end{array}%
\right) >=<X,u>+<Y,v>.
\end{equation*}

\bigskip \noindent  Then
\begin{equation*}
\phi _{(X,Y)}(u,v)=E\left(\exp <\left( 
\begin{array}{c}
X \\ 
Y%
\end{array}%
\right) ,\left( 
\begin{array}{c}
u \\ 
v%
\end{array}%
\right) >\right)
\end{equation*}%
\begin{equation*}
=\mathbb{E}\left(e^{<X,u>+<Y,v>})=E(e^{<X,u>})\mathbb{E}(e^{<Y,v>}\right)
\end{equation*}
\begin{equation*}
=\phi _{X}(u)\times \phi _{Y}(v).
\end{equation*}

\bigskip \noindent  The proof is over. $\blacksquare$.\\

\bigskip \noindent  Now, we want to move to next very important other characterization. When $d=1$, we have an explicit inversion formula which expresses the \textit{cdf} of a probability law on $\mathbb{R}^d$ by means of its characteristic function. The characterization of a probability law on $\mathbb{R}$ by its characteristic function follows from this inversion formula.\\

\noindent But when $d>1$, things are more complicated and we may need a non-standard version of the Theorem of Stone-Weierstrass Theorem. In that case a more general characterization of probability measures in metric spaces may be useful. So we begin with general characterizations.\\

\bigskip
\noindent \textbf{III - Characterization of a probability law on a metric space}.\\

\noindent Let us suppose that we are working on a metric space $(E,\rho)$ endowed with the metric $\rho$. We are going to use the class $C_b(E)$ of real-valued continuous and bounded functions defined on $E$. Let us begin by reminding that, by the $\lambda$-$\pi$ Lemma (See \cite{ips-mestuto-ang}, Exercise 11 of Doc 04-02, Part VI, page 228), the class of open sets $\mathcal{O}$ is a determining class of probability measures since it is a $\pi$-system, containing $E$ and generating $\mathcal{B}(E)$, that is, for two probability measures $\mathbb{P}_j$ ($j\in \{1,2\}$) on $(E,\mathcal{B}(E))$, we have

\begin{equation} 
(\mathbb{P}_1=\mathbb{P}_2) \Leftrightarrow (\forall G \in \mathcal{O}, \ \mathbb{P}_1(G)=\mathbb{P}_2(G)). \label{caractOpens}
\end{equation}

\bigskip \noindent Actually, this characterization can be extended to integrals of $f \in C_b(E)$. For this, we need the following tool.

\begin{lemma} \label{probab_02_02_approxOpensBoundedFunctions} Let  $G$ be a non-empty open in $E$. There exists a non-decreasing sequence of functions $(f_m)_{m\geq 1}$ such that :\\

\noindent (1) for each $m\geq 1$, $f_m$ is a Lipschitz function of coefficient $m$ and
$$
0 \leq f_m \leq 1_G.
$$

\noindent and $f_m=0$ on $\partial G$ and\\

\noindent (2) we have
$$
f_m \uparrow 1_G, \ as \ m \uparrow +\infty.
$$
\end{lemma}

\noindent The proof is given in the Appendix Chapter \ref{proba_02_appendix} in Lemma \ref{proba02_rv_charac_opensets} (page \pageref{proba02_rv_charac_opensets}).\\

\noindent This lemma may be used to get the following characterization : for two probability measures $\mathbb{P}_j$ ($j\in \{1,2\}$) on $(E,\mathcal{B}(E))$, we have

\begin{equation} 
(\mathbb{P}_1=\mathbb{P}_2) \Leftrightarrow \left(\forall f \in C_b(E), \ \int f \ d\mathbb{P}_1=\int f \ d\mathbb{P}_2\right). \label{caractBoundedFunctions}
\end{equation}

\noindent To establish this, we only need to show the indirect implication. Suppose that right-hand assertion holds. For any $G \in \mathcal{O}$, we consider the the sequence $(f_m)_{n\geq 1}$ in Lemma
\ref{probab_02_02_approxOpensBoundedFunctions} and we have

\begin{equation*} 
\forall m\geq 1, \ \int f_m \ d\mathbb{P}_1=\int f_m \ d\mathbb{P}_2.
\end{equation*}

\noindent By letting $m\uparrow +\infty$ and by applying the Monotone Convergence Theorem, we get $\mathbb{P}_1(G)=\mathbb{P}_2(G)$. Since this holds for any $G \in \mathcal{O}$, we get 
$\mathbb{P}_1=\mathbb{P}_2$ by Formula \label{caractOpens}.\\

\noindent \textbf{IV - Characterization of a probability law on $\mathbb{R}^d$ by its characteristic function}.\\

\noindent We are going to prove that characteristic functions also determine probability laws on $\mathbb{R}^d$. 

\begin{theorem} \label{proba02_rv_cf02}
Let $X$ and $Y$ be two random variables with values in $\mathbb{R}^{d}$. Their characteristic functions coincide on $\mathbb{R}^{d}$ if and only if do their probability laws on 
$\mathcal{B}(\mathbb{R}^{d})$, that is
\begin{equation*}
\Phi _{X}=\Phi _{Y} \Leftrightarrow \mathbb{P}_{X}=\mathbb{P}_{Y}.
\end{equation*}
\end{theorem}

\bigskip \noindent  \textbf{Proof}. We are going to use an approximation based on a version of the theorem of Stone-Weierstrass. Let us begin by reminding that the class of intervals of $\mathbb{R}^d$ 
$$
\mathcal{I}_d=\{]a,b[=\prod_{j=1}^{d} ]a_j,b_j[, \ a\leq b, (a,b) \in \left(\mathbb{R}^d\right)^2\}
$$

\noindent is a $\pi$-system, contains $E=\mathbb{R}^d$ and generates $\mathcal{B}\left(\mathbb{R}^d\right)$. By the the $\lambda$-$\pi$ Lemma, it constitutes a determining class for probability measures.\\

\noindent Fix $G=]a,b[$ with $a_j<b_j$, for all $1\leq j \leq d$. For any $j \in \{1,\cdots,d\}$ and consider the sequence $(f_{j,m})_{m\geq 0}\subset C_b(\mathbb{R}^d)$ constructed for $G_j=]a_j,b_j[$ in 
Lemma \ref{probab_02_02_approxOpensBoundedFunctions}. The numbers $f_{j,m}(a_j)$ and $f_{j,m}(a_j)$ are zero. So the functions

$$
f_m(x)=\prod_{j=1}^{d} f_{j,m}(x_j), \ x=(x_1,...,x_d)^t \in \mathbb{R}^d, \ m\geq 1,
$$

\noindent vanish on the border $\partial G$ of $G$ since 

$$
\partial G=\{x \in G, \ \exists j\in\{1,\cdots,d\}, \ x_j=a_j \ or \ x_j=b_j\}
$$
 
\noindent It becomes clear that for any probability measure $\mathbb{L}$ on $(\mathbb{R}^d,\mathcal{B}(\mathbb{R}^d))$, we have

$$
\forall ]a,b[ \in \mathcal{I}_d, \ f_m \uparrow 1_{]a,b[} \ and \ \int f_m \ d\mathbb{L} \uparrow \mathbb{L}(]a,b[), \ as \ m\uparrow +\infty.
$$

\noindent We may seize the opportunity to state a new characterization of probability measures of $\mathbb{R}^d$. Let $C_{b,0}(\mathbb{R}^d)$ be the class of functions $f$ for which there exists $]a,b[ \in \mathcal{I}_d$ such that $0\leq f\leq 1$ and $f=0$ outside $]a,b[$. We get that : \\

\noindent For two probability measures $\mathbb{P}_j$ ($j\in \{1,2\}$) on $(\mathbb{R}^d,\mathcal{B}(\mathbb{R}^d))$ :

\begin{equation} 
(\mathbb{P}_1=\mathbb{P}_2) \Leftrightarrow  \left(\forall f \in C_{b,0}(\mathbb{R}^d), \ \int f \ d\mathbb{P}_1=\int f \ d\mathbb{P}_2\right). \label{caractBoundedFunctionSpec}
\end{equation}

\noindent Now fix $f \in C_{b,0}(\mathbb{R}^d)$ associated with  $[a,b]$. Let $\varepsilon \in ]0,1[$ . Fix $r>0$ and $K_r=[-r,r]^d$. We choose $r$ such that 

\begin{equation}
-r\leq \min(a_1,...,a_d) \ and \  r \geq \max(b_1,...,b_d) \label{approx18}
\end{equation}

\noindent and

\begin{equation}
\mathbb{P}_X(K_r^c)+\mathbb{P}_Y(K_r^c)\leq \frac{\varepsilon}{2(2+\varepsilon)}. \label{approx18b}
\end{equation}

\noindent Now consider the class $\mathcal{H}$ of finite linear combinations of functions of the form

\begin{equation}
\prod_{j=1}^{d} \exp\biggr(i n_j\pi x_j/r\biggr), \label{EX}
\end{equation}

\noindent where $n_j\in \mathbb{Z}$ is a constant and $i$ is the normed complex of angle $\pi/2$ and let $\mathcal{H}_{r}$ be the class of the restrictions $h_{r}$ of elements $h \in \mathcal{H}$ on 
$K_r=[-r,r]^d$.\\

\noindent It is clear that $\mathcal{H}_{r}$ is a sub-algebra of $C_b(K_r)$ with the following properties.\\

\noindent (a) for each $h \in \mathcal{H}$, the uniform norm of $h$ on $\mathbb{R}^d$ is equal to the uniform norm of $h$ on $K_r$, that is

$$
\|f\|_{\infty}=\sup_{x \in \mathbb{R}^d} |h(x)|= \sup_{x \in K_r} |h(x)|=\|f\|_{K_r}. 
$$ 

\noindent This comes from that remark that $h$ is a finite linear combination  of functions of the form in Formula \ref{EX} above and each factor $\exp\left(i n_j\pi x_j/r\right)$ is a 
$2r$-periodic function.\\

\noindent (b) $\mathcal{H}_{r}$ separates the points of $K_r\setminus \partial K_r$ and separates points of $K_r\setminus \partial K_r$ from points of $\partial K_r$. Indeed, if $x$ and $y$ are two points in $K_r$, at the exception where both of them are edge points of $K_r$ of the form

$$
(x,y) \in \{(s_1,...,s_d) \in K_r, \ \forall j\in\{1,...,d\}, \ s_j=r \ or \ s_j=r\}^2,
$$

\noindent \noindent there exists $j_0 \in \{1,\cdots,d\}$ such that $0<|x_{j_0}-y_{j_0}|<2r$ that is $|(x_{j_0}-y_{j_0})/r|<2$ and the function 
$$
h_{r}(x)=\exp(i\pi x_{j_0}/r)
$$ 

\noindent separates $x$ and $y$ since $h_{r}(x)=h_{r}(y)$ would imply $\exp( i\pi (x_{j_0}-x_{j_0})/r)=1$, which in term would imply $x_{j_0}-x_{j_0}=2\ell r$, $\ell \in \mathbb{Z}$. The only possible value of $\ell$ would be zero and this is impossible since $x_{j_0}-y_{j_0} \neq 0$.\\

\noindent (c) For all the points in $t \in \partial K_r$, the function $g(t)\equiv 0 \in \mathcal{H}_{r}$ converges to $f(t)=0$.\\

\noindent (d) $\mathcal{H}_{r}$ contains all the constant functions.\\

\noindent We may then apply Corollary 2 in \cite{loSW2018} (Corollary \ref{sec_EF_cor_05} in the appendix, page \pageref{sec_EF_cor_05}) to get that : there exists $h_{r} \in \mathcal{H}_{r}$ such that

\begin{equation}
\|f-h_{r}\|_{K_r} \leq \varepsilon/4. \label{approx19}
\end{equation}

\bigskip \noindent and by Point (a) above (using also that the norm of $f\in C_{b,0}$ less or equal to $1$), we have

\begin{equation}
\|h\|_{\infty}=\|h_{r}\|_{K_r]} \leq \|f\|_{\infty} + \varepsilon/4\leq 1+\varepsilon. \label{approx20}
\end{equation}

\bigskip \noindent Now, y the assumption of equality of the characteristic functions, we have 
\begin{equation*}
\mathbb{E}(h(X))=\mathbb{E}(h(Y)).
\end{equation*}

\bigskip \noindent  We have have 
\begin{eqnarray*}
\mathbb{E}(f(X))-\mathbb{E}(f(Y)) &=&\left( \int f\text{ }d\mathbb{P}_{X}-\int h\text{ }d\mathbb{P}_{X}\right) +\left( \int h\text{ }d\mathbb{P}_{X}-\int h\text{ }d\mathbb{P}_{Y}\right)\\ 
&+&\left( \int h\text{ }d\mathbb{P}_{Y}-\int f\text{ }d\mathbb{P}_{Y}\right)  \\
&=&\left( \int f\text{ }d\mathbb{P}_{X}-\int h\text{ }d\mathbb{P}_{X}\right)
+\left( \int h\text{ }d\mathbb{P}_{Y}-\int f\text{ }d\mathbb{P}_{Y}\right). 
\end{eqnarray*}

\bigskip \noindent  The first term satisfies

\begin{eqnarray}
\mathbb{E}\left\vert \int f\text{ }d\mathbb{P}_{X}-\int h_{r}\text{ }d\mathbb{P}_{X}\right\vert  &\leq &\int_{K_r]}\left\vert f-h_{r}\right\vert \text{ }d\mathbb{P}_{X}  \label{cv34d} \\
&&+\int_{K_r^c}\left\vert f-h\right\vert \text{ }d\mathbb{P}_{X} \notag \\
&\leq &\varepsilon/4+(\left\Vert f\right\Vert +\left\Vert h\right\Vert)\mathbb{P}_{X}(K_r^c), \notag \\
&\leq &\varepsilon/4+(2+\varepsilon) \mathbb{P}_{X}(K_r^c), \notag 
\end{eqnarray}

\noindent where we used Formulas \ref{approx19} and \ref{approx20}.\\

\bigskip \noindent  By treating the second term in the same manner, we also get 
\begin{equation}
 \mathbb{E}\left\vert \int f\text{ }d\mathbb{P}_{Y}-\int h\text{}d\mathbb{P}_{Y}\right\vert \leq \varepsilon/4+(2+\varepsilon) \mathbb{P}_{Y}(K_r^c). \label{cv34e}
\end{equation}

\bigskip \noindent  By putting together Formulas (\ref{cv34d}) and (\ref{cv34e}) and by remembering Formulas \eqref{approx18} and \eqref{approx18b}, we get 

\begin{equation*}
\left\vert \mathbb{E}(f(X))-\mathbb{E}(f(Y))\right\vert \leq \varepsilon/2+(2+\varepsilon)(\mathbb{P}_{X}(K_r^c)+\mathbb{P}_{Y}(K_r^c))\leq \varepsilon.
\end{equation*}

\bigskip \noindent for any $\varepsilon \in ]0,1[$. So, for all $f \in C_{b,0}(\mathbb{R}^d)$,
\begin{equation*}
\int f d\mathbb{P}_X=\int f d\mathbb{P}_Y.
\end{equation*}

\noindent We close the proof by applying Formula \eqref{caractBoundedFunctionSpec} above.\\


\bigskip \noindent  \textbf{V - Inversion Formula on $\mathbb{R}$ and applications}.\\

\noindent Here, we consider the characteristic function of a Lebesgue-Stieljes measures on $\mathbb{R}$, not necessarily a probability measure. After the proof of the following proposition, we will get another characterization of probability laws by characteristic functions by means of \textit{cdf}'s. Let us begin to state the 
 
\begin{proposition} \label{proba02_rv.inv01} 
Let $F$ be an arbitrary distribution function. Let

$$
\Phi(x)=\int exp(itx) \ d\lambda_F(x), \ x \in \mathbb{R},
$$

\bigskip \noindent  where $\lambda_{F}$ denotes the Lebesgue-Stieljes measure associated with $F$. Set for two
reals numbers $a$ and $b$ such that $a<b$,

\begin{equation}
J_{U}=:J_{U}(a,b)=\frac{1}{2\pi }\int_{-U}^{U}\frac{e^{-iau}-e^{-ibu}}{iu}\Phi _{X}(u)\text{ }du. \label{proba02_rv.inv02}
\end{equation}

\bigskip \noindent (a) Then, we have, as U $\rightarrow +\infty$, $J_{U}$ converges to

\begin{eqnarray}
F(b-)-F(a)+\frac{1}{2}\biggr(F_{X}(a)-F(a-)+F(b)-F(b-)\biggr)). \ \label{proba02_rv.inv03} 
\end{eqnarray}
 
\bigskip \noindent (b) If $a$ and $b$ are continuity points of $F$, then

\begin{eqnarray}
F(b)-F(a)= \lim_{U \rightarrow +\infty} J_{U}. \label{proba02_rv.inv04}
\end{eqnarray}

\bigskip  \noindent  If $F$ is absolutely continuous, that is there exists a measurable $\lambda$-\textsl{a.e.} finite function $f$ such that for $x\in \mathbb{R}$,
\begin{equation}
F(x)= \int_{-\infty}^{x} f(t) d(x), \label{proba02_rv.inv05}
\end{equation}

\bigskip \noindent then, we have $\lambda$-a.e.,

$$
f(x)=\frac{1}{2\pi} \int_{-\infty}^{+\infty} e^{-ixu} \Phi(u)\text{ }du.
$$
\end{proposition}

\bigskip \noindent  \textbf{Proof}. Recall Dirichlet's Formula
\begin{eqnarray*}
\int_{-\infty }^{0}\frac{\sin x}{x}dx&=&\int_{0}^{+\infty }\frac{\sin x}{x}dx\\
&=&\lim_{b\rightarrow +\infty }\int_{0}^{b}\frac{\sin x}{x}dx=\pi/2,
\end{eqnarray*}

\noindent which can be proved, for example, using complex integration based on residues. We deduce from it that
the numbers
$$
\int_{a}^{b}\frac{\sin x}{x}dx=\int_{a}^{0}\frac{\sin x}{x}dx+\int_{0}^{b}\frac{\sin x}{x}dx, \ a\leq 0 \leq b
$$ 

\noindent are uniformly bounded in $a$ and $b$, say by $M$. By using Fubini's theorem, we have
\begin{eqnarray*}
J_{U}&=&\frac{1}{2\pi }\int_{-U}^{U}\frac{e^{-iau}-e^{-ibu}}{iu}\left( \int e^{iux}d\mathbb{P}_{X}(x)\right) \text{ }du\\
&=&\int d\mathbb{P}_{X}(dx)\times \frac{1}{2\pi }\int_{-U}^{U}\frac{e^{-i(a-x)u}-e^{-i(b-x)u}}{iu}du\\
&=&\int J(U,x) \ d\mathbb{P}_{X}(x),
\end{eqnarray*}

\noindent where

\begin{eqnarray*}
J(U,x)&=&\frac{1}{2\pi }\int_{-U}^{U}\frac{e^{-i(a-x)u}-e^{-i(b-x)u}}{iu}du\\
&=&\frac{1}{2\pi i}\int_{-U}^{U}\frac{\cos (u(a-x))-\cos (u(b-x))}{u}du\\
&+&\frac{1}{2\pi }\int_{-U}^{U}\frac{\sin (u(b-x))-\sin (u(a-x))}{u}du.
\end{eqnarray*}

\noindent But, we also have
 
\begin{equation*}\normalfont
\int_{-U}^{U}\frac{\cos (u(a-x))-\cos (u(b-x))}{u}du=0.
\end{equation*}

\bigskip \bigskip \noindent  Since the integrated functions are odd and the integration is operated on a symmetrical compact interval with respect to zero. We get

\begin{eqnarray*}
J(U,x)&=&\frac{1}{2\pi }\int_{-U}^{U}\frac{\sin (u(b-x))-\sin (u(a-x))}{u}du\\
&=&\frac{1}{2\pi }\int_{-U(b-x)}^{U(b-x)}\frac{\sin v}{v}dv-\frac{1}{2\pi } \int_{-U(a-x)}^{U(a-x)}\frac{\sin v}{v}dv.
\end{eqnarray*}

\bigskip \noindent  Thus,  $J(U,x)$ uniformly bounded bounded by $M/\pi$. Next by considering the position of $x$ with respect of the interval $(a,b)$ and by handling accordingly the signs of $(b-x)$ and $(a-x)$, we easily arrive at the following set of implications :
\begin{equation*}
(x<a\text{ }or \text{ }x>b)\Rightarrow J(U,x)\rightarrow 0\text{ as }%
U\rightarrow +\infty,
\end{equation*}%
\begin{equation*}
(x=a\text{ }or\text{ }x=b)\Rightarrow J(U,x)\rightarrow 1/2\text{ as }%
U\rightarrow +\infty
\end{equation*}%
\begin{equation*}
(a<x<b)\Rightarrow J(U,x)\rightarrow 1\text{ as }U\rightarrow +\infty .
\end{equation*}

\bigskip \noindent Then
 
\begin{equation*}
J(U,x)\rightarrow 1_{]a,b[}+\frac{1}{2}1_{\{a\}}+\frac{1}{2}1_{\{b\}}.
\end{equation*}

\bigskip \noindent  From there, we apply the Fatou-Lebesgue Theorem to get

\begin{eqnarray*}
&&J_{U}\rightarrow \int \left( 1_{]a,b[}+\frac{1}{2}1_{\{a\}}+\frac{1}{2}%
1_{\{b\}}\right) d\mathbb{P}_{X}(x)\\
&=&F(b-)-F(a)+\frac{1}{2}\biggr(F(a)-F(a-)+F(b)-F(b-)\biggr).
\end{eqnarray*}

\bigskip \noindent  This proves Point (a). $\square$\\
 
\noindent Point (b) If $a$ and $b$ are continuity points of $F$, the limit in (\ref{proba02_rv.inv03}) reduces to $F(b)-F(a)$. $\square$\\

\noindent Point (c)  Now, from (\ref{proba02_rv.inv04}), we deduce that $F$ is continuous and next, the derivative of $F$ at $x$ is $f(x)$ when $f$ is continuous. But a measurable function that is integrable is $\lambda$-\textit{a.e.} continuous. So, 

$$
\frac{dF(x)}{dx}=f(x), \ \lambda-a.e.
$$

\bigskip \noindent   Also, by (\ref{proba02_rv.inv04}), we have for all $h>0$,
 
\begin{equation}
\frac{F(a+h)-F(a)}{h}=\lim_{U\rightarrow +\infty} \frac{1}{2\pi }\int_{-U}^{U}\frac{e^{-iau}-e^{-i(a+h)u}}{ihu}
\Phi _{X}(u)\text{ }du. \label{proba02_rv.inv07}
\end{equation}

\noindent Then for any $a\in \mathbb{R}$,

\begin{eqnarray*} 
f(a)&=& \lim_{h\rightarrow 0} \frac{F(a+h)-F(b)}{h}\\
&=& \frac{1}{2\pi} \lim_{h\rightarrow 0} \lim_{U\rightarrow +\infty} \int_{-U}^{U}\frac{e^{-iau}-e^{-i(a+h)u}}{ihu}\Phi(u)\text{ }du\\ 
&=& \frac{1}{2\pi} \lim_{U\rightarrow +\infty} \lim_{h\rightarrow 0} \int_{-U}^{U}\frac{e^{-iau}-e^{-i(a+h)u}}{ihu}\Phi(u)\text{ }du\\
&=& \frac{1}{2\pi} \lim_{U\rightarrow +\infty} \int_{-U}^{U} \lim_{h\rightarrow 0}\frac{e^{-iau}-e^{-i(a+h)u}}{ihu}\Phi(u)\text{ }du;
\end{eqnarray*}

\bigskip \noindent  where the exchange between integration and differentiation in the last line is allowed by the use the Fatou-Lebesgue theorem based on the fact that the integrated function is bounded by the unity which is integrable on $(-U,U)$, $U$ fixed.\\

\bigskip \noindent  So, we arrive at

\begin{eqnarray*} 
f(a)&=& \frac{1}{2\pi} \lim_{U\rightarrow +\infty} \int_{-U}^{U} \lim_{h\rightarrow 0}\frac{e^{-iau}-e^{-i(a+h)u}}{ihu}\Phi(u)\text{ }du\\
&=&  \frac{1}{2\pi} \lim_{U\rightarrow +\infty} \int_{-U}^{U}  e^{-iau} \Phi(u)\text{ }du\\
&=&  \frac{1}{2\pi} \int_{-\infty}^{+\infty} e^{-iau} \Phi(u)\text{ }du.
\end{eqnarray*}

\bigskip \noindent  $\lambda$-a.e. $\blacksquare$\\

\bigskip \noindent  \textbf{Application}.\\

\noindent Now, let us use this to prove Theorem \ref{proba02_rv_cf02} for $k=1$.\\

\noindent Let $X$ and $Y$ be two \textit{rrv}'s with equal characteristic functions. By (\ref{proba02_rv.inv04}), their distribution functions $F_X$ and $F_Y$ are equal on the set
$D_{X,Y}$ of continuity points of both $F_X$ and $F_Y$. The complement of that set is at most countable. So, for $x\in D_{X,Y}$ fixed, we may find a sequence of numbers $(x_n)_{n\geq 0}$ such that 

$$
(x_n)_{n\geq 0} \subset D, \text{ such that } x_n \rightarrow x \text{ as }n\uparrow +\infty.
$$

\bigskip \noindent   So, we will have for any $n\geq 0$

$$
F_X(x)-F_X(a_n)=F_Y(x)-F_Y(a_n).
$$

\bigskip \noindent  By letting $n\rightarrow +\infty$, we get for all $x\in D_{X,Y}$

$$
F_X(x)=F_Y(x).
$$

\bigskip  \noindent  For any $x\in \mathbb{R}$, we also can can find monotone sequence $(x_n)_{n\geq 0}$ such that 

$$
(x_n)_{n\geq 0} \subset D, \text{ such that } x_n \downarrow x \text{ as }n\uparrow +\infty.
$$

\bigskip \noindent  By right-continuity at $x$ of $F_X$ and $F_Y$, we have

$$
F_X(x)=\lim_{n\uparrow +\infty} F_X(x_n)=\lim_{n\uparrow +\infty} F_Y(x_n)=F_Y(x_n).
$$

\bigskip \noindent  Conclusion $F_X=F_X$. Thus by the first characterization, $X$ and $Y$ have the same probability law.$\blacksquare$\\

\bigskip \noindent  \textbf{IV - A characterization of independence}.\\

\noindent We are going to see that Point (c) of Theorem \ref{proba02_rv_prop_cf} is a rule for independence because of Theorem \ref{proba02_rv_cf02}. We have

\begin{theorem} \label{proba02_rv_prop_cf031}  Let $X$ and $Y$ be two random variables respectively with values in $\mathbb{R}^{d}$ and in $\mathbb{R}^{k}$ and defined on the same probability measure. The random variables $X$ and $Y$ are independent if and only if for any $u\in \mathbb{R}^{d}$ and for $v \in \mathbb{R}^{k}$, we have
\begin{equation}
\phi _{(X,Y)}(u,v)=\phi _{X}(u)\times \phi _{Y}(v)  \label{prodfc03}
\end{equation}
\end{theorem}

\bigskip \noindent  \textbf{Proof}. We need only to prove that (\ref{prodfc03}) implies independence of $X$ and $Y$. Suppose that (\ref{prodfc03}) holds. It is clear that the left-hand member of (\ref{prodfc03}) is the characteristic function of the product measure $\mathbb{P}_{X}\otimes \mathbb{P}_{Y}$. Since the characteristic functions of the probability laws $\mathbb{P}_{(X,Y)}$ and $\mathbb{P}_{X}\otimes \mathbb{P}_{Y}$ coincide, we get
$$
\mathbb{P}_{(X,Y)}=\mathbb{P}_{X}\otimes \mathbb{P}_{Y},
$$

\bigskip \noindent  which is the definition of the independence between $X$ and $Y$. $\blacksquare$\\

\bigskip \noindent  \textbf{V - Characteristic functions and moments for \textit{rrv}}.\\

\noindent We are going to see how to find the moments from the characteristic function in the following. Let us write

\begin{equation*}
\Phi_{X}(u)=\int e^{iux}d\mathbb{P}_{X}(x), \ u \in \mathbb{R}.
\end{equation*}%

\bigskip \noindent  The function 
\begin{equation*}
g(u,x)=\cos (ux)+i\sin (ux)=e^{iux}
\end{equation*}

\bigskip \noindent  is differentiable with respect to $u$ and its derivative is 
\begin{equation*}
g^{\prime}(u,x)=ix(\cos (ux)+i\sin ux))=ixe^{iux}.
\end{equation*}

\bigskip \noindent  It is bounded by $Y(x)=\left\vert x\right\vert$. The integral of this function $Y(x)$ is the mathematical
expectation of $X$, that is,
 
\begin{equation*}
\int Y(x)d\mathbb{P}_{X}(x)=\int \left\vert x\right\vert d\mathbb{P}%
_{X}(X)=E\left\vert X\right\vert .
\end{equation*}

\bigskip \noindent  Suppose that the mathematical expectation is finite. Then, by the Dominated Convergence Theorem (See  Point 06.14 in Doc 06.14, Chapter 7, in \cite{ips-mestuto-ang} of this series), we may exchange integration and differentiation. The method may be repeated by a second differentiation and so forth. We conclude this quick discussion in

\begin{proposition} \label{proba_02_rv_sec_06_propMoments}
If $\mathbb{E}(X)$ exists and is finite, then the function $u\mapsto \phi _{X}(u)$ is differentiable and we have
\begin{equation*}
\phi _{X}^{\prime }(u)=\int ixe^{iux}d\mathbb{P}_{X}(x).
\end{equation*}

\bigskip \noindent  And we have  
\begin{equation*}
i\times \mathbb{E}(X)=\phi _{X}^{\prime }(0).
\end{equation*}

\bigskip \noindent  More generally, if for $k\geq 1$, $\mathbb{E}\left\vert X\right\vert
^{k}$ exists and is finite, then the function $u\mapsto \phi _{X}(u)$ is differentiable $k$ times with 
\begin{equation*}
\phi _{X}^{(k)}(u)=i^{k}\int x^{k}e^{iux}d\mathbb{P}_{X}(x)
\end{equation*}

\bigskip \noindent  and
\begin{equation*}
\mathbb{E} X^{k}=-i^{k}\quad \phi _{X}^{(k)}(0).
\end{equation*}
\end{proposition}

\bigskip
\section{Convolution, Change of variables and other properties} \label{proba_02_rv_sec_07}

\bigskip \noindent  \textbf{I - Convolution product of probability density functions on $\mathbb{R}$}.\\

\noindent Let $X$ and $Y$ be two real-valued random variables which are defined on the same probability space $(\Omega ,\mathcal{A},\mathbb{P})$ and mutually independent. 
Set $Z=X+Y$. By definition, the probability law of $Z$ is called the convolution product of the probability laws of $\mathbb{P}_X$ and of $\mathbb{P}_Y$, denoted as

\begin{equation}
\mathbb{P}_{Z}= \mathbb{P}_{X} \ast \mathbb{P}_{Y}. \label{proba02_rv_pconv01}
\end{equation}

\bigskip \noindent  Now, suppose that $X$ and $Y$ have probability density functions $f_{X}$ and $f_{Y}$ with respect to the Lebesgue measure $\lambda$.
Then $Z$ has an absolutely probability density function $f_Z$ denoted as

\begin{equation*}
f_{X+Y}=f_{X}\ast f_{Y}.
\end{equation*}

\bigskip \noindent  We have the following

\begin{proposition} Let $X$ and $Y$ be to real-valued and independent random variables, defined on the same probability measure $(\Omega ,\mathcal{A},\mathbb{P})$, and admitting the probability density functions $f_{X}$ and $f_{Y}$ with respect to a $\sigma$-finite product  measure $\nu=\nu_1 \otimes \nu_2$. Then $Z$ has a \textsl{pdf} $f_Z$ which has the two following to expressions :

\begin{equation*}
f_{X}\ast f_{Y}(z)=\int_{\mathbb{R}}f_{X}(z-x)\text{ }f_{Y}(x)\text{ }d\lambda(x).
\end{equation*}
\end{proposition}

\bigskip \noindent  \textbf{Proof}. Assume the hypotheses of the proposition hold. Let us use the joint probability law of $(X,Y)$ to have
\begin{equation*}
F_{Z}(z)=\mathbb{P}(X+Y\leq x)=\int_{(x+y\leq z)}d\mathbb{P}_{(X,Y)}(x,y).
\end{equation*}

\bigskip \noindent  Since $X$ and $Y$ are independent, we have  
\begin{equation*}
\mathbb{P}_{(X,Y)}=\mathbb{P}_{X} \otimes \mathbb{P}_{Y}.
\end{equation*}

\bigskip \noindent  We may apply Fubini's Theorem to get  
\begin{eqnarray*}
F_{Z}(z)&=&\int (x+y)d\mathbb{P}_{(X,Y)}(x,y)=\int d\mathbb{P}_{X}(x)\int_{y\leq z-x}d\mathbb{P}_{Y}(y)\\
&=&\int f_{X}(x) d\nu_1(y) \left(\int_{y\leq z-x}f_{Y}(y) \ d\nu_2(y)\right).
\end{eqnarray*}

\bigskip \noindent We recall the the Lebesgue measure is invariant by translation. Let us make the change variable $u=y+x$, to have 

\begin{eqnarray*}
F_{Z}(z)&=&\int f_{X}(x)\text{ }dx\left(\int_{u\leq z}f_{Y}(u-x)\text{ }du\right)\\
&=&\int f_{X}(x)\text{ }dx\text{ }\left(\int_{-\infty }^{z}f_{Y}(u-x)\text{ }du\right).
\end{eqnarray*}

\bigskip \noindent Let us use again the Fubini's Theorem to get 

\begin{equation*}
F_{Z}(z)=\int \int_{-\infty }^{z}f_{Y}(u-x)\text{ }f_{X}(x)\text{ }dx\text{ }%
dy=\int_{-\infty}^{z}\left(\int f_{Y}(u-x)\text{ }f_{X}(x)\text{ }dx\right)\text{ }du.
\end{equation*}

\bigskip \noindent Taking the differentiation with respect to $z$, we get 
 
\begin{equation*}
f_{Z}(z)=\int f_{Y}(z-x)\text{ }f_{X}(x)\text{ }dx.
\end{equation*}

\bigskip \noindent For such a formula for discrete random variables, the reader is referred \cite{ips-probelem-ang} of this series, Formula (3.24), Part D, Section 3, Chapter 6.\\

\bigskip \noindent  \textbf{II - Change of Variable by Diffeomorphisms and Introduction to the Gauss Random variables}.\\ \label{changeVariable}

\noindent \textbf{(a) Recall of the Change of Variable Formula for Riemann Integrals on $\mathbb{R}^{d}$} (See \cite{valiron}, page 275, for double integration).\\

\bigskip \noindent  Suppose we have the following Riemann integral on $\mathbb{R}^{d}$,

\begin{equation*}
I=\int_{D}f(x_1,x_2,\cdots,x_d)\text{ }dx_1 dx_2 \cdots dx_d,
\end{equation*}

\bigskip \noindent  where $D$ is a domain of $\mathbb{R}^{d}$. We will write for short with  $x^t=(x_1,x_2,\cdots,x_d)$, 

\begin{equation*}
I=\int_{D}f(x)\text{ }dx.
\end{equation*}

\bigskip \bigskip \noindent  Suppose that we have a diffeomorphism $h$ from an other domain $\Delta$ of $\mathbb{R}^{k}$ to $D$. This means that the function 
\begin{equation*}
h:\Delta \mapsto D
\end{equation*}

\bigskip \noindent  (a) is a bijection (one-to-one mapping).\\

\noindent (b) $h$ and its inverse function $g=h^{-1}$ have continuous partial derivatives (meaning that they are both of class $C^1$).\\
 
\noindent  Let us write $h$ as :
\begin{equation*}
D\ni x=h(y)\longleftrightarrow y\in \Delta .
\end{equation*}

\bigskip \noindent  The components of $h$ are denoted by $h_{i}$ : 
\begin{equation*}
x_{i}=h_{i}(y)=h_{i}(y_{1},...,y_{d}).
\end{equation*}

\bigskip \noindent  The $d$-square matrix of elements  
\begin{equation*}
\frac{\partial x_{i}}{\partial y_{j}}=\frac{h_{i}(y_{1},...,y_{d})}{\partial
y_{j}}
\end{equation*}

\bigskip \noindent  written also as
\begin{equation*}
M(h)=\left[ \left( \frac{\partial x_{i}}{\partial y_{j}}\right) _{ij}\right].
\end{equation*}

\bigskip \noindent  is called the Jacobian matrix of the transformation. The absolute value of its determinant is called the \textbf{Jacobian coefficient} of the change of variable. We way write it as

\begin{equation*}
J(h,y)=\det \left( \left[ \left( \frac{\partial x_{i}}{\partial y_{j}}%
\right) _{ij}\right] \right) .
\end{equation*}

\bigskip \noindent  The change of variable formula is the following

\begin{equation*}
I=\int_{\Delta }f(h(y))\text{ }\left\vert J(h,y)\right\vert \text{ }dy.
\end{equation*}

\bigskip \noindent  We replace $x$ by $h(y)$, the domain $D$ by $\Delta$, but we multiply the integrated function by the Jacobian coefficient (depending on $y$).\\

\noindent \textbf{(b) An example leading to the Gaussian probability Law}.\\

\noindent Let us give a classical example. Suppose we want to compute 
\begin{equation*}
I=\int_{[0,+\infty \lbrack \times \lbrack 0,+\infty \lbrack
}e^{-(x^{2}+y^{2})}dx\text{ }dy.
\end{equation*}

\bigskip \noindent   Let us the polar coordinates of $(x,y)$ in $(\mathbb{R}_{+})^2$ :

\begin{equation*}
\left\{ 
\begin{array}{c}
x=r\cos \theta \\ 
y=r\sin \theta%
\end{array}
\right .
\end{equation*}%

\bigskip \noindent  with
\begin{equation*}
(x,y)\in D=[0,+\infty \lbrack \times \lbrack 0,+\infty \lbrack
\longleftrightarrow (r,\theta )\in \lbrack 0,+\infty \lbrack \times \lbrack
0,\pi /2].
\end{equation*}

\bigskip \noindent  The Jacobian coefficient of the transformation is  
\begin{equation*}
J(r,\theta )=\left\vert 
\begin{array}{cc}
\frac{\partial x}{\partial r} & \frac{\partial x}{\partial \theta } \\ 
\frac{\partial y}{\partial r} & \frac{\partial y}{\partial \theta }%
\end{array}%
\right\vert =\left\vert 
\begin{array}{cc}
\cos \theta & -r\sin \theta \\ 
\sin \theta & r\cos \theta%
\end{array}%
\right\vert =r\cos ^{2}\theta +r\sin ^{2}\theta =r.
\end{equation*}

\bigskip \noindent  We apply the change of variable formula to have 
\begin{equation*}
I=\int_{[0,+\infty \lbrack \times \lbrack 0,\pi /2]}re^{-r^{2}}dr\text{ }%
d\theta =\int_{[0,\pi /2]}d\theta \int_{\lbrack 0,+\infty \lbrack
}re^{-r^{2}}dr=\frac{\pi }{4}.
\end{equation*}

\bigskip \noindent  By the Fubini's Formula, we have 
\begin{equation*}
I=\int_{[0,+\infty \lbrack }e^{-x^{2}}dx\text{ }\int_{[0,+\infty \lbrack
}e^{-y^{2}}dy=(\int_{0}^{+\infty }e^{-u^{2}}du)^{2}.
\end{equation*}

\bigskip \noindent  Then, we have

\begin{equation*}
\int_{0}^{+\infty }e^{-u^{2}}du=\frac{\sqrt{\pi }}{2}.
\end{equation*}

\bigskip \noindent  Finally, by a new change of variable, where we take the evenness of the function $u \mapsto exp(-u^{2}/2)$, leads to  
\begin{equation*}
\frac{1}{\sqrt{2\pi }}\int_{-\infty }^{+\infty }e^{-u^{2}/2}du=1.
\end{equation*}

\bigskip \noindent  This is a probability density function. Compare this with the lengthy proof in Section 5, Chapter 7, in \cite{ips-probelem-ang} of this series.\\

\noindent Let us apply this Formula to finding new density functions.\\

\bigskip \noindent  \textbf{(c) Finding a probability density function by change of variables}.\\

\noindent Let $X$ be a random variable in $\mathbb{R}^d$ of probability density function $f_{X}$ with respect to the Lebesgue measure on  $\mathbb{R}^{k}$, still denoted by 
$\lambda_k(x)=dx$. Suppose that $D$ is the support of $X$. Let 
\begin{equation*}
h:\Delta \mapsto D
\end{equation*}

\bigskip \noindent  be a diffeomorphism and  
\begin{equation*}
Y=h^{-1}(X)
\end{equation*}

\bigskip \noindent  be another random vector. Then, the probability density function of $Y$ exists and is given by

\begin{equation*}
f_{Y}(y)=f_{X}(h(y))\text{ }\left\vert J(h)\right\vert \text{ }1_{\Delta}(y). \ (CVF)
\end{equation*}%

\bigskip \noindent  This follows from an immediate application of the variable change formula. Let $B$ be a borel set of $\mathbb{R}^{d}$, we have

\begin{eqnarray*}
\int_{x\in h(B)}f_{X}(x)\text{ }dx=\int_{h^{-1}(x)\in B}f_{X}(x)\text{ }dx.
\end{eqnarray*}

\bigskip \noindent  Let us apply the variable change formula as follows :

\begin{eqnarray*}
\mathbb{P}(Y\in B) &=& \int_{y\in B}f_{X}(h(y))\text{ }1_{\Delta}(y)\left\vert J(h,y)\right\vert dy\\
&=&\int_{B}\left\{ f_{X}(h(y))\text{ }1_{\Delta }(y)\left\vert
J(h,y)\right\vert \right\} dy.
\end{eqnarray*}

\bigskip \noindent  We deduce from this that 
\begin{equation*}
f_{Y}(y)=f_{X}(h(y))\text{ }1_{\Delta }(y)\left\vert J(h,y)\right\vert
\end{equation*}

\bigskip \noindent is the probability density function of $Y$.\\

\noindent In Mathematical Statistics, this tool is extensively used, especially for Gaussian random variables.\\

\bigskip \noindent  \textbf{(d) Important example}.\\

\noindent This example is important for two reasons. First, we will have to apply many of the techniques used in the this chapter and secondly, the object of the example is the starting point of the study of stable laws.\\

\noindent Let us consider two independent $\mathcal{E}(\lambda)$-random variables $X_1$ and $X_2$, $\lambda>0$ on a same probability space (such a construction is achieved through the Kolmogorov construction method) and let us set $X_s=X_1-X_2$. The \textit{pdf} of $X_s$ is the convolution product of $f_{X_1}$ and $f_{-X_2}$. The \textit{pdf} $f_{-X_0}$ is

$$
f_{-X_2}(y)=\lambda \exp(\lambda y), \ y\leq 0.
$$

\bigskip \noindent So, we have for all $x \in \mathbb{R}$,

\begin{eqnarray*}
f_{X_s}(x)&=&\left(f_{X_1} * f_{-X_2}\right)(x)\\
&=&\int f_{X_1}(x-y) f_{-X_2}(y) \ dy\\
&=&\lambda^2 \int \biggr(\exp(-\lambda (x-y)) 1_{(x-y\geq 0)}\biggr) \biggr(\exp(\lambda y) 1_{(y\leq  0)}\biggr) \ dy\\
&=&\lambda^2 \int \biggr(\exp(-\lambda (x-y)) 1_{(y \leq  x)}\biggr) \biggr(\exp(\lambda y) 1_{(y\leq  0)}\biggr) \ dy.
\end{eqnarray*}

\bigskip \noindent If $x\leq 0$, we have

\begin{eqnarray*}
f_{X_s}(x)&=&\lambda^2 \int_{-\infty}{x} \exp(-\lambda (x-y)) \exp(\lambda y) \ dy\\
&=&\lambda^2 \exp(-\lambda x) \int_{-\infty}{x} \exp(2\lambda y) \ dy\\
&=&\lambda^2 \exp(-\lambda x) \left[\frac{e^{2\lambda y}}{2\lambda}\right]_{-\infty}^{x}\\
&=&\frac{\lambda}{2} \exp(\lambda x).
\end{eqnarray*}

\bigskip \noindent If $x\geq 0$, we have

\begin{eqnarray*}
f_{X_s}(x)&=&\lambda^2 \int_{-\infty}{0} \exp(-\lambda (x-y)) \exp(\lambda y) \ dy\\
&=&\lambda^2 \exp(-\lambda x) \int_{-\infty}^{0} \exp(2\lambda y) \ dy\\
&=&\frac{\lambda}{2} \exp(-\lambda x).
\end{eqnarray*}

\bigskip \noindent In total, we have

\begin{equation}
f_{X_s}(x)=\frac{\lambda}{2} \exp(-\lambda |x|), \ x \in \mathbb{R}. \label{symm}
\end{equation}

\bigskip \noindent Next, let us see an interesting application of the inversion formula. The characteristic function of $X_s$ is 

\begin{eqnarray*}
\Phi_{X_s}(u)&=&\Phi_{X_1-X_2}(u)=\Phi_{X_1-X_2}(u)\\
&=&\Phi_{X_1}(u) \Phi_{X_2}(-u)\\
&=&\frac{1}{1-it/\lambda}\frac{1}{1+it/\lambda}
\end{eqnarray*}

\bigskip \noindent which leads to

\begin{equation}
\Phi_{X_s}(u)=\frac{\lambda^2}{\lambda^2+u^2}, \ u \in \mathbb{R}.
\end{equation}

\bigskip \noindent Now let us apply the inversion formula to this characteristic function. We have  $\lambda$-a.e. for all $x \in \mathbb{R}$

\begin{eqnarray*}
\frac{\lambda}{2} \exp(-\lambda |x|)&=&\frac{1}{2\pi} \int e^{-iux} \Phi_{X_s}(u) \ du\\
&=&\frac{1}{2\pi} \int e^{-iux} \frac{\lambda^2}{\lambda^2+u^2} \ du,\\
\end{eqnarray*}

\bigskip \noindent and by dividing both members by by ($\lambda/2)$ we get

\begin{equation}
\exp(-\lambda |x|)= \int e^{-iux} \frac{\lambda}{\pi(\lambda^2+u^2)} \ du, \label{fc_cauchy00}
\end{equation}

\bigskip \noindent and by replacing $x$ by $-x$, we conclude that we have $\lambda$-a.e. for all $x \in \mathbb{R}$, 

\begin{equation*}
\int e^{iux} \frac{\lambda}{\pi(\lambda^2+u^2)} \ du = \exp(-\lambda |x|). \label{fc_cauchy}
\end{equation*}

\bigskip \noindent It happens that

$$
f_{C(0,\lambda)}=\frac{\lambda}{\pi(\lambda^2+u^2)}, \ x \in \mathbb{R},
$$

\bigskip \noindent is the \textit{pdf} of a Cauchy random variable of parameters $0$ and $\lambda>0$ (see Chapter \ref{proba_02_upl}, Section \ref{proba_02_upl_sec_02}, \pageref{cauchyDis}). We just found the characteristic of a Cauchy random variable, which is not easy to find by direct methods.\\

\section{Copulas} \label{proba_02_rv_07a}

\noindent The lines below should form a part of Section \ref{proba_02_rv_07a} which was devoted to \textit{cdf}'s. But, nowadays, the notion of copula is central in Statistics theory, although copulas are simply particular \textit{cdf}'s in Probability. So we think that introducing to copulas in a section might serve for references.\\

\noindent A very recurrent source on copulas is \cite{nelsen}. However, the lines below will use the note of \cite{losklar2018}.\\

\noindent \textbf{Definition} A copula on $\mathbb{R}^d$ is a \textit{cdf} $C$ whose marginal \textit{cdf}'s defined by, for $1\leq i \leq d$,
$$
\mathbb{R} \ni s \mapsto C_{i}(s)=C\left( +\infty, ..., +\infty,\underset{i-th \ argument}{\underbrace{s}}, +\infty, ..., +\infty\right),
$$
 
\bigskip  \noindent are all equal to the $(0,1)$-uniform \textit{cdf} which in turn is defined by
$$
x \mapsto x 1_{[0,1[} + 1_{[1,+\infty[},
$$

\noindent and we may also write, for all $s \in [0,1]$, 

\begin{equation}
C_{i}(s)=C\left(1, ..., 1,\underset{i-th \ argument}{\underbrace{s}}, 1, ...,1 \right)=s. \label{probcop}
\end{equation}

\bigskip \noindent The copula became very popular with following the important theorem of \cite{sklar1959}

\begin{theorem} \label{theo1} For any \textit{cdf} $F$ on $\mathbb{R}^d$, $d\geq 1$, there exists a copula $C$ on 
$\mathbb{R}^d$ such that

\begin{equation}
\forall x \in \mathbb{R}^d, \ F(x)=C(F_1(x),...,F_d(x)). \label{sklar} 
\end{equation}
\end{theorem} 

\bigskip \noindent This theorem is now among the most important tools in Statistics since it allows to study the dependence between the components of a random vector through the copula, meaning that the intrinsic dependence does not depend on the margins.\\

\noindent We are going to provide a recent proof due to \cite{losklar2018}. Fortunately, the tools we need are available in the current series, in particular in \cite{ips-wcrv-ang}.\\

\noindent \textbf{Proof of \cite{sklar1959}'s Theorem}.\\

\noindent \textbf{(A) - Complements}. We first need some complements to the properties of the generalized inverse function given in \cite{ips-wcrv-ang}. Let us begin by defining generalized functions. Let $[a,b]$ and $[c,d]$ be non-empty intervals of $\mathbb{R}$ and let $G : [a,b] \mapsto [c,d]$ be a non-decreasing mapping such that
\begin{eqnarray*}
c&=& \inf_{x \in [a,b]} G(x),\ \ \ (L11)\\
d&=& \sup_{x \in [a,b]} G(x). \ \ \ (L12)\\
\end{eqnarray*}

\noindent Since $G$ is a mapping, this ensures that
\begin{eqnarray*}a&=&\inf \{x\in \mathbb{R}, \ G\left(x\right) > c \}, \ \ \ (L13)\\
b&=&\sup \{x\in \mathbb{R}, G\left( x\right) < d\}. \ \ \ (L14)\\
\end{eqnarray*}

\noindent If $x=a$ or $x=b$ is infinite, the value of $G$ at that point is meant as a limit. If $[a,b]$ is bounded above or below in $\mathbb{R}$, $G$ is extensible on $\mathbb{R}$ by taking $G\left( x\right)=G\left( a+\right)$ for $x\leq a$ and $G\left( x\right) =G\left( b-0\right) $ for $x\geq b$. As a general rule, we may consider $G$ simply as defined on $\mathbb{R}$. In that case, $a=lep(G)$ and $b=uep(G)$ are called \textit{lower end-point} and \textit{upper end-point} of $G$.\\

\noindent The generalized inverse function of $G$ is given by

$$
\forall u \in [lep(G),uep(G)],  \ G^{-1}\left(u\right) =\inf \left\{ x\in \mathbb{R}, \ G\left( x\right) \geq u\right\}.
$$

\bigskip \noindent The properties of $G^{-1}$ have been thoroughly studied, in particular in \cite{billingsley}, \cite{resnick}. The results we need in this paper are gathered and proved in wcrv or in \cite{ips-wcia-fr} (Chapter 4, Section 1) and reminded as below.

\begin{lemma} \label{lemA} Let $G$ be a non-decreasing right-continuous function with the notation above. Then  $G^{-1}$ is left-continuous and we have

\begin{equation*}
\forall u \in [c,d], \ G(G^{-1}(u))\geq u \ (A) \  and \  \forall x \in [a,b], \ G^{-1}(G(x))\leq x \ (B) 
\end{equation*}

\bigskip \noindent and 

\begin{equation}
\forall x \in [lep(G), uep(G)], \ G^{-1}(G(x)+0)=x. \label{FF}
\end{equation}
\end{lemma}

\noindent \textit{Proof}. The proof of Formulas (A) and (B) are well-known and can be found in the cited books above. Let us prove Formula \eqref{FF} for any $x \in [a,b]$.\\

\noindent On one side, we start by the remark that $G^{-1}(G(x)+0)$ is the limit of $G^{-1}(G(x)+h)$ as $h \searrow 0$. But for any $h>0$, $G^{-1}(G(x)+h)$ is the infimum of the set of $y\in [a,b]$ such that $G(y)\geq G(x)+h$. Any these $y$ satisfies $y\geq x$. Hence $G^{-1}(G(x)+0)\geq x$.\\

\noindent On the other side $G(x+h) \searrow G(x)$ by right-continuity of $G$, and by the existence of the right-hand limit of the non-decreasing function $G^{-1}(\circ)$, $G^{-1}(G(x+h)) \searrow G^{-1}(G(x)+0)$. Since $G^{-1}(G(x+h))\leq x+h$ by Formula (B), we get that $G^{-1}(G(x)+0)\leq x$ as $h \searrow 0$. The proof is complete. $\square$\\

\noindent \textbf{(B) - Proof of Sklar's Theorem}.  Define for $s=(s_1,s_2,\cdots,s_d) \in [0,1]^d$,

\begin{equation}
C(s)=F( F_1^{-1}(s_1+0), F_2^{-1}(s_2+0), \cdots, F_d^{-1}(s_d+0)). \label{cop}
\end{equation}
 
\bigskip \noindent It is immediate that $C$ assigns non-negative volumes to cuboids of $[0,1]^d$, since according to Condition (DF2), Formula \eqref{VP} for $C$ derives from the same for $F$ where the arguments are the form $F_i^{-1}(\circ+0)$, $1\leq i\leq d$.\\

\noindent Also $C$ is right-continuous since $F$ is right-continuous as well as each $F_i^{-1}(\circ+0)$, $1\leq i\leq d$. By passing, this explains why we took the right-limits because the $F_i^{-1}(\circ)$'s are left-continuous.\\

\noindent Finally, by combining Formulas \eqref{FF} and \eqref{cop}, we get the conclusion of Sklar in Formula \eqref{sklar}. The proof is finished. $\square$\\

\section{Conclusion} \label{proba_02_rv_sec_08}

\bigskip \noindent \textbf{(A) Back to independence of Random vectors}.\\
 
\noindent Because of the importance of the notion of independence and since several characterizations of the independence are scattered this chapter and in Chapter \ref{proba_02_intotp}, we think that a summary on this point may be useful to to reader.\\

\noindent (1) The most general definition of a finite family of random variables is given in Definition \ref{proba_02_01_defIndependence} (page \pageref{proba_02_01_defIndependence}). This definition covers all type of random variables and uses the finite product measure. Random variables of an infinite family are independent if and only if the elements each finite sub(family are independence.\\

\noindent In this general case, Theorem \label{proba_02_intotp_th02} (page \pageref{proba_02_intotp_th02}) gives a general characterization.\\

\noindent (2) When we have a random real-valued vector in $\mathbb{R}^d$, $d\leq 1$, the independence of the coordinates and the independence of sub-vectors are characterized :\\

\noindent (2a) in Theorem \ref{proba_02_rv_th03} (page \pageref{proba_02_rv_th03}), using the cumulative distribution functions,\\

\noindent (2b) in Theorem \ref{proba02_rv_prop_cf031} (page \pageref{proba02_rv_prop_cf031}), using the characteristic functions,\\

\noindent (2c) in Theorem \ref{proba_02_rv_th04} (page \pageref{proba_02_rv_th04}), using the probability density functions with respect to the measure.\\

\bigskip \noindent \textbf{(B) General advices to determine probability laws}.\\

\noindent Now, we have the means to characterize the usual probability laws by their distribution functions or their characteristic functions. Its is also important
to know the parameters of the usual laws. In the next two chapters, we will be dealing with them. Estimating these from data is one of the most important
purposes of Statistics.\\

\bigskip \noindent  In trying to find the probability laws, the following ideas may be useful.\\

\noindent (A) Using the convolution product to find the probability law of the sum of two independent real-value random variables.\\

\noindent (B) Using the product of characteristic function to find the probability law of the sum of two independent random variables of equal dimension.\\

\noindent (C) Finding the distribution function of the studied random variable and differentiate it if possible, and try to identify a known probability law.\\

\noindent (D) Directly finding the characteristic function of the studied random variable and trying to identify a known probability law.\\

\noindent (E) Using the Change of Variable Formula to derive \textsl{pdf}'s if applicable.\\

\noindent (F) In particular, the following easy \textit{stuff} may be useful : \\

\noindent \textbf{A useful stuff}. Suppose that two random elements $X$ and $Y$, are defined on the same probability space and take their values in the same measure space $(E,\mathcal{B}, \nu)$, which is endowed with a measure $\nu$. Suppose that $X$ and $Y$ have \textsl{pdf}'s $f_X$ and $f_Y$ with respect to $\nu$ and that these two \textit{pdf}'s a common support $\mathcal{V}$ and have a common variable part, meaning that there exist a non-negative function $h : E \rightarrow \mathbb{R}$ and constants $C_1>0$ and 
$C_2>0$ such that

$$
\forall x\in E, \ f_X(x)=C_1 h(x) \text{ and } f_Y(x)=C_2 h(x).
$$

\bigskip \noindent Then $f_X=f_Y$, $\nu$-\textit{a.s.} and $C_1=C_2$. $\Diamond$\\

\noindent The proof is obvious since

$$
1=\int_{\mathcal{V}} f_X \ d\nu =C_1 \int_{\mathcal{V}} h \ d\nu = \int_{\mathcal{V}} f_X \ d\nu =C_2 \int_{\mathcal{V}} h \ d\nu.
$$

\bigskip \noindent which leads to

$$
C_1=C_2=1/\left(\int_{\mathcal{V}} h \ d\nu\right).
$$

\bigskip \noindent Despite its simplicity, this \textit{stuff} is often used and allows to get remarkable Analysis formulas, some of them being extremely difficult, even impossible, to establish by other methods.\\

%% file: proba_02_03_ang.tex
\chapter{Usual Probability Laws} \label{proba_02_upl}

\noindent We begin to focus on real random variables. Later, we will focus on Random vectors in Chapter \ref{proba_02_gauss}.\\

\noindent Actually, the researchers have discovered a huge number of probability laws. A number of dictionaries of probability laws exist (See for example, \cite{kotz1}, which is composed of 13 volumes at least). Meanwhile, people are still continuing to propose new probability laws and their properties (see \cite{okorie} for a recent example).\\

\noindent This chapter is just a quick introduction to this wide area. A short list among the most common laws is given. Some others concern new important probability laws (Skewed normal, hyperbolic, etc.).\\

\noindent \textbf{I - Review of usual probability law on $\mathbb{R}$}. \label{proba_02_upl_od}\\

\noindent We begin with discrete random variables. For such random variables, the discrete integration formula is used to find the parameters and the characteristic functions. This has already been done in the monograph of \cite{ips-probelem-ang}. We will not repeat the computations here.

\section{Discrete probability laws} \label{proba_02_upl_sec_01}

\bigskip \noindent For each random variable $X$, the values set or support $\mathcal{V}_X$, the probability density function with respect to the appropriate counting measure, the characteristic function and/or the moment generating function and the moments are given.\\

\noindent \textbf{(1) Constant random variable $X=a$, a.s, $a\in \mathbb{R}$}.\\

\noindent $X$ takes only one value, the value $a$.\\

\noindent Discrete probability density function on $\mathcal{V}_X=\{a\}$ :\\

\begin{equation*}
\mathcal{V}_X=\{a\} \text{ and } \mathbb{P}(X=a)=1.
\end{equation*}

\bigskip \noindent Distribution function : 

\begin{equation*}
F_{X}(x)=1_{[a,+\infty \lbrack }, \ x \in \mathbb{R}.
\end{equation*}

\bigskip \noindent Characteristic function :
 
\begin{equation*}
\Phi_{X}(u)=e^{iau}, \ t\in \mathbb{R}.
\end{equation*}

\bigskip \noindent Moment generating function :
 
\begin{equation*}
\varphi_{X}(u)=e^{au}, \ t\in \mathbb{R}.
\end{equation*}

\bigskip \noindent Moments of order $k\geq 1$
\begin{equation*}
\mathbb{E}X^{k}=a^{k}, \mathbb{E}(X-a)^{k}=0.
\end{equation*}

\bigskip \noindent \textbf{A useful remark}. A constant random variable is independent from any other random variable defined on the same probability space. Indeed let $X=a$ and $Y$ be another any other random variable defined on the same probability space. The joint characteristic function of $(X,Y)$ is given by

\begin{eqnarray*}
\Phi_{(X,Y)}(u,v)&=&\mathbb{E}\exp(iXu+iYv)=\mathbb{E}\biggr(\exp(iau) \exp(iYv)\biggr)\\
&=&\exp(iau) \mathbb{E} \exp(iYv)=\Phi_{X}(u) \Phi_{Y}(v),
\end{eqnarray*}

\bigskip  \noindent for any $(u,v)\in \mathbb{R}^2$. By Theorem \ref{proba02_rv_prop_cf031} in Chapter \ref{proba_02_rv}, $X$ and $Y$ are independent.\\

\bigskip \noindent \textbf{(2) Uniform Random variable on $\{1,2,...,n\}$, $n\geq 1$}.\\

\noindent $X\sim \mathcal{U}(1,2,...,n)$ takes each value in  $\{1,2,...,n\}$ with the same probability.\\

\noindent Discrete probability density function on $\mathcal{V}_X=\{1,2,...,n\}$ :

\begin{equation*}
\mathbb{P}(X=k)=1/n, \ k\in \{1,...,n\}
\end{equation*}

\bigskip \noindent Distribution function :

\begin{equation*}
F(x)=\left\{ 
\begin{tabular}{lll}
0 & if  & $x<1,$\\
$\frac{i-1}{n}$ & if  & $\frac{i-1}{n}\leq x< \frac{i}{n}$, $1\leq i \leq n,$\\
1 & if & $x\geq n.$ \\
\end{tabular}
\right. 
\end{equation*}

\bigskip \noindent Characteristic function :

\begin{equation*}
\Phi _{X}(u)=\frac{1}{n}\sum_{j=1}^{n}e^{iju}, \ u \in \mathbb{R}.
\end{equation*}

\bigskip \noindent Moments of order $k\geq 1$ :

\begin{equation*}
\mathbb{E}X^{k}=\frac{1}{n}\sum_{j=1}^{n}j^{k}.
\end{equation*}

\bigskip \noindent Mathematical expectation and variance :

$$
\mathbb{E}(X)=\frac{n+1}{2}, \ \mathbb{V}ar(X)=\frac{(n-1)(n+1)(4n+3)}{12}.
$$

\bigskip \noindent \textbf{(3) Bernoulli Random Variable with parameter $0<p<1$}.\\

\noindent $X\sim \mathcal{B}(p)$ takes two values : $1$ (Success) and $0$ (failure).\\

\noindent Discrete probability density function on $\mathcal{V}_X=\{0,1\}$ :

\begin{equation*}
\mathbb{P}(X=1)=p=1-\mathbb{P}(X=0).
\end{equation*}

\bigskip \noindent Distribution function : 

\begin{equation*}
F(x)=0\times 1_{]-\infty ,0[}+p\times 1_{[0,1[}+1_{[1,+\infty \lbrack }, \ x \in \mathbb{R}.
\end{equation*}

\bigskip \noindent Characteristic function :
 
\begin{equation*}
\Phi _{X}(u)=q+pe^{iu}, \ u \in \mathbb{R}.
\end{equation*}

\bigskip \noindent Moments of order $k\geq 1$ :

\begin{equation*}
\mathbb{E}X^{k}=p.
\end{equation*}

\bigskip \noindent Mathematical expectation and variance :

$$
\mathbb{E}(X)=p, \ \mathbb{V}ar(X)=pq.
$$

\bigskip \noindent \textbf{(4) Binomial random variable with parameters $0<p<1$ and $n\geq 1$}.\\

\noindent  $X\sim \mathcal{B}(n,p)$  takes its values in $\{0,1,...,n\}$.\\

\noindent Discrete probability density function on $\mathcal{V}_X=\{0,1,...,n\}$ :

\begin{equation*}
\mathbb{P}(X=k)=C_{n}^{k}\text{ }p^{k}(1-p)^{n-k}, \ k=0,...,n.
\end{equation*}

\bigskip \noindent Characteristic function. Since $X$ is the sum of $n$ independent Bernoulli $\mathcal{B}(p)$ random variables, Point (b) and Theorem \ref{proba02_rv_prop_cf} and the value of the characteristic function of a Bernoulli random variable, yield

\begin{equation*}
\Phi _{X}(u)=(q+pe^{iu})^{n}, \ u \in \mathbb{R}.
\end{equation*}

\bigskip \noindent Mathematical expectation and variance :

\begin{equation*}
\mathbb{E}(X)=np,\text{ and  }Var(X)=np(1-p).
\end{equation*}

\bigskip \noindent The above parameters are computed by still using the decomposition of Binomial random variable by into a sum of independent Bernoulli random variables.\\
 
\bigskip \noindent \textbf{(5) Geometric Random Variable with parameter $0<p<1$}.\\

\noindent  $X\sim \mathcal{G}(p)$ takes its values in $\mathbb{N}$.\\

\noindent Discrete probability density function on $\mathcal{V}_X=\mathbb{N}$ :

\begin{equation*}
\mathbb{P}(X=k)=p(1-p)^{k}, \ k\in \mathbb{N}.
\end{equation*}

\bigskip \noindent Characteristic function :

\begin{equation*}
\Phi _{X}(u)=p/(1-qe^{iu}), \ u \in \mathbb{R}. 
\end{equation*}

\bigskip \noindent Mathematical expectation and variance :

\begin{equation*}
\mathbb{E}(X)=q/p,\text{ }Var(X)=q/p^{2}.
\end{equation*}

\bigskip \noindent \textbf{(6) Negative Binomial Random Variable with parameters $r\geq 1$ and $0<p<1$}.\\

\noindent  $X_r\sim \mathcal{B}N(r,p)$  takes the values in $\{r,r+1,... \}$.\\

\bigskip \noindent Discrete probability density function on $\mathcal{V}_X=\{r,r+1,... \}$ :

\begin{equation*}
\mathbb{P}(X=k)=C_{k-1}^{r-1}p^{k}(1-p)^{r-k},\ \ k\geq r.
\end{equation*}

\bigskip \noindent Characteristic function.  Since $X_r$ is the sum of $r$ independent Geometric $\mathcal{G}(p)$ random variables, Theorem  and the value of the characteristic function of a Bernoulli random variable, yield

\begin{equation*}
\Phi _{X}(u)=\left\{ pe^{iu}/(1-qe^{iu}\right\} ^{r}, \ u<-\log(1-p).
\end{equation*}

\bigskip \noindent Mathematical expectation and variance :

\begin{equation*}
\mathbb{E}(X)=rq/p,\text{ }Var(X)=rq/p^{2}.
\end{equation*}

\bigskip \noindent \textbf{(7) Poisson Random variable of parameter $\lambda>0$}.\\

\noindent $X\sim \mathcal{P}(\lambda)$ takes its values in $\mathbb{N}$.\\

\noindent Discrete probability density function on $\mathcal{V}_X=\mathbb{N}$ :

\begin{equation*}
\mathbb{P}(X=k)=\frac{\lambda ^{k}}{k!}e^{-\lambda },\text{ }k\geq 0.
\end{equation*}

\bigskip \noindent Characteristic function
\begin{equation*}
\Phi _{X}(u)=\exp (\lambda (e^{iu}-1)), \ u \in \mathbb{R}.
\end{equation*}

\bigskip \noindent Mathematical expectation and variance :

\begin{equation*}
\mathbb{E}(X)=Var(X)=\lambda .
\end{equation*}

\bigskip \noindent \textbf{(8) Hyper-geometric Random Variable}.\\

\noindent  $X\sim \mathcal{H}(N,\theta ,n)$ or  $H(N,M,n),$ $1\leq n\leq N,$ $0<\theta <1,$ $\theta =M/N,$ takes its values in $\{0,1,...,min(n,M)\}$.\\

\bigskip \noindent Discrete probability density function on $\mathcal{V}_X=\{0,1,...,min(n,M)\}$ :

\begin{equation*}
(X=k)=\frac{C_{M}^{k}\times C_{N-M}^{n-k}}{C_{N}^{n}},\text{ }k=0,...,\min(n,M).
\end{equation*}

\bigskip \noindent Characteristic function  of no use.\\

\noindent Mathematical expectation and variance :

\begin{equation*}
\mathbb{E}(X)=rM/n, \text{ and } V(X)=rM(n-M)(n-r)/\{n^{2}(n-1)\}.
\end{equation*}

\bigskip \noindent \textbf{(9) Logarithmic Random Variable}.\\

\noindent  $X\sim Log(p)$ takes its values in $\{1,2,...\}$.\\

\bigskip \noindent Discrete probability density function on $\mathcal{V}_X=\{1,2,...\}$ :

\begin{equation*}
\mathbb{P}(X=k)=-qk/(k \log p),\text{ }k\geq 1.
\end{equation*}

\bigskip \noindent Characteristic function :
\begin{equation*}
\Phi _{X}(u)=\log (1-qe^{iu})/\log (p),  \ u \in \mathbb{R}.
\end{equation*}

\bigskip \noindent Moment Generating function :
\begin{equation*}
\Phi _{X}(u)=\log (1-qe^{u})/\log (p),  \ u<-\log(1-p).
\end{equation*}

\bigskip \noindent Mathematical expectation and variance :

\begin{equation*}
\mathbb{E}(X)=-q/(p\log (p)),V(X)=-q(q+\log (p))/(p\log (p)).
\end{equation*}

\newpage
\section{Absolutely Continuous Probability Laws} \label{proba_02_upl_sec_02}

\bigskip \noindent For each random variable $X$, the support $\mathcal{V}_X$, the probability density function with respect to the Lebesgue measure, the characteristic function and/or the moment generating function, the moments are given. By definition, the support $\mathcal{V}_X$ of $X$ is given by

$$
\mathcal{V}_X=\overline{\{x \in \mathbb{R}, f_X(x) \neq 0\}}
$$  

\bigskip \noindent We also have
\begin{equation*}
\mathbb{P}(X \in \mathcal{V}_X)=1.
\end{equation*}

\noindent \noindent For any real-valued random variable, we may define

\begin{equation*}
lep(F)=\inf \{x,F(x)>0\}
\end{equation*}%

\bigskip \noindent and

\begin{equation*}
uep(F)=\sup \{x,F(x)<1\}.
\end{equation*}

\bigskip \noindent where $lep(F)$ and $uep(F)$ respectively stand for \textit{lower end-point of $F$} and \textit{upper end-point of $F$}. As a result we have

$$
X \in [lep(F), uep(F)], \ \textit{a.e.}
$$

\bigskip \noindent The first examples given without computations are done in \cite{ips-mestuto-ang}.\\

\bigskip \noindent \textbf{(1) Continuous uniform Random variable on a bounded compact set}.\\

\noindent  Let $a$ and $b$ be  two real numbers such that $a<b$. $X\sim \mathcal{U}(a,b)$.\\

\noindent Domain : $\mathcal{V}_X=[a,b]$.\\

\noindent Absolutely continuous probability density function on $\mathcal{V}_X=[a,b]$ :

\begin{equation*}
f_{X}(x)=\frac{1}{b-a}1_{[a,b]}(x), \ x \in \mathbb{R}.
\end{equation*}

\noindent Distribution function : 

\begin{equation*}
F_{X}(x)=\left\{ 
\begin{array}{c}
1 \ \ if \ \ x\geq b, \\ 
(x-a)/(b-a) \ \  if \ \ a\leq x\leq b, \\ 
0 \ \ if \ \ x\leq a.%
\end{array}%
\right.
\end{equation*}

\bigskip \noindent Characteristic function :

\begin{equation*}
\Phi_{X}(u)=\frac{e^{ibu}-e^{iau}}{iu(b-a)}, u \in \mathbb{R}.
\end{equation*}

\bigskip \noindent Moments of order $k\geq 1$ :

$$
\mathbb{E}X^k=\frac{b^{k+1}-b^{k+1}}{(k+1)(b-a)}.
$$

\bigskip \noindent Mathematical expectation and variance :

\begin{equation*}
\mathbb{E}(X)=(a+b)/2\text{, }et\text{ }Var(X)=(b-a)^{2}/12.
\end{equation*}

\bigskip \noindent \textbf{(2) Exponential Random Variable of parameter $b>0$}.\\

\noindent  $X\sim \mathcal{E}(b)$ is supported on $\mathbb{R}_{+}.$\\

\noindent Absolutely continuous probability density function on $\mathcal{V}_X=\mathbb{R}_{+}$ :

\begin{equation*}
f_{X}(x)=b e^{-b x}1_{(x\geq 0)}.
\end{equation*}

\bigskip \noindent Distribution function : 

\begin{equation*}
F_{X}(x)=(1-e^{-b x})1_{(x\geq 0)}.
\end{equation*}

\bigskip \noindent Characteristic function :

\begin{equation*}
\Phi_{X}(u)=(1-iu/b )^{-1}, 
\end{equation*}

\bigskip \noindent Moment Generating Function :

\begin{equation*}
\phi_{X}(u)=(1-u/b)^{-1}, \ u <b.  
\end{equation*}

\bigskip \noindent Moments of order $k\geq 1$

$$
\mathbb{E}(X^k)=\frac{k!}{\lambda^k}.
$$

\bigskip \noindent Mathematical expectation and variance :

\begin{equation*}
\mathbb{E}(X)=1/\lambda ,Var(X)=1/\lambda ^{2}.
\end{equation*}

\bigskip \noindent \textbf{(3) Gamma Random variable with Parameter $a>0$ and $b>0$}.\\

\noindent  $X\sim \gamma (a,b)$  is defined ob $\mathbb{R}_{+}$.\\

\noindent Absolutely continuous probability density function on $\mathcal{V}_X=\mathbb{R}_{+}$ :
\begin{equation*}
f_{X}(x)=\frac{b^{a}}{\Gamma (a)}x^{a-1}e^{-bx}1_{(x\geq 0)}
\end{equation*}

\bigskip \noindent with

\begin{equation*}
\Gamma (a)=\int_{0}^{\infty }x^{a-1}\text{ }e^{-x}\text{ }dx.
\end{equation*}

\bigskip \noindent Characteristic function :

\begin{equation*}
\Phi _{X}(u)=(1-iu/b)^{-a}.
\end{equation*}

\bigskip \noindent Moments of order $k\geq 1$ :

$$
\mathbb{E}(X^k)= \frac{1}{b^k} \prod_{j=0}^{k} a+j.
$$

\bigskip \noindent Mathematical expectation and variance :
\begin{equation*}
\mathbb{E}(X)=a/b,Var(X)=a/b^{2}.
\end{equation*}

\bigskip \noindent \textbf{Be careful}. Some authors, many of them in North America, take $\gamma(a,1/b)$ as the \textit{gamma law}. If you read somewhere that
$\mathbb{E}(X)=ab$ for $X \sim \gamma(a,b)$, be aware that in our definition we have $X \sim \gamma(a,1/b)$.\\

\bigskip \noindent \textbf{(4) Symmetrized Exponential random  variable with $\lambda>0$}.\\

\noindent $X\sim \mathcal{E}_s(\lambda)$ in defined on $\mathbb{R}$.\\

\noindent From the non-negative random variable $X$ , it is always possible to define a symmetrized random variable $X_s$ by considering two independent $\mathcal{E}(\lambda)$-random variables $X_1$ and $X_2$ on a same probability space (such a construction is achieved through the Kolmogorov construction method) and by setting $X_s=X_1-X_2$. Another way to define it is to have an $\mathcal{E}(\lambda)$-random variable $X$ and a $(0,1)$-uniform random variable $U$ independent of $X$ and to set $X_s = -X 1_{(U\leq 0.5)} + X 1_{(U > 0.5)}$. We are going to use the first method. It is clear that $X_s$ is a symmetric random variable. Further if $X$ admits an absolutely continuous \textit{pdf}, $X_s$ has the $pdf$

$$
f_{X_x}(x)=\frac{1}{2} f_{X}(|x|), \ x\in \mathbb{R}.
$$

\noindent By applying this to the exponential random variable, a Symmetrized Exponential random  following $X\sim \mathcal{E}_s(\lambda)$ has the following \textit{pdf}.\\
 
\noindent Absolutely continuous probability density function on $\mathcal{V}_X=\mathbb{R}$ :

\begin{equation*}
f_{X}(x)=\frac{\lambda}{2} \exp(-\lambda |x|), \ x \in \mathbb{R}.\\
\end{equation*}

\bigskip \noindent Distribution function : \\

$$
F_X(x)=\frac{1}{2} e^{\lambda x} 1_{(x<0)} + \left(1- \frac{1}{2}e^{-\lambda x}\right) 1_{(x\geq 0)}, \ x\in \mathbb{R}.
$$

\bigskip \noindent Characteristic function : (See Formula \ref{fc_cauchy}, Chapter \ref{proba_02_rv}, page \pageref{fc_cauchy})

$$
\Phi_X(u)=\exp(-\lambda |u|), \ u \in \mathbb{R}.
$$

\bigskip \noindent Mathematical expectation and variance :\\

$$
\mathbb{E}X_s=0 \ \ and \ \ \mathbb{V}ar(X_s)=\frac{2}{\lambda^2}. 
$$

\bigskip \noindent To justify the variance, we may remark that $X_s=(X_1-1/\lambda)-(X_2-1/\lambda)$, that is, 
$X_s=(X_1-\mathbb{E}(X_1))-(X_2-\mathbb{E}(X_1=2))$ and exploit that $X_s$ is the difference between two independent and centered random variables.\\

\bigskip \noindent \textbf{Remark}.  For $\gamma=1$, this law holds the name of Laplace random variable.\\

\bigskip \noindent \textbf{(5) Beta Random variables of parameter $a>0$ and $b>0$}.\\

\noindent $X\sim B(a,b)$ is defined on $(0,1)$.\\

\noindent Absolutely continuous probability density function on $\mathcal{V}_X$ :

\begin{equation*}
f_{X}(x)=\frac{1}{B(a,b)}x^{a-1}(1-x)^{b-1}1_{(0,1)}(x),
\end{equation*}

\bigskip \noindent where

$$
B(a,b)=\int_{1}^{1} x^{a-1}(1-x)^{b-1} \ dx.
$$

\bigskip \noindent Mathematical expectation and variance :

\begin{equation*}
\mathbb{E}(X)=a/(a+b) \text{ and } Var(X)=ab/[(a+b)^{2}(a+b+1)].
\end{equation*}

\bigskip \noindent \textbf{(6) Pareto Random Variable of parameter $a>0$}.\\

\noindent $X\sim Par(a,\alpha )$, with parameters $\alpha >0$ and $a\geq 0$, is supported by $]a,+\infty]$.\\

\noindent Absolutely continuous probability density function on $\mathcal{V}_X=]a, +\infty[$ :

\begin{equation*}
f_{X}(x)=\alpha a^{\alpha} x^{-\alpha-1} 1_{(x>a)}.
\end{equation*}

\bigskip \noindent \textbf{(7) Cauchy random variable with $\lambda>0$ and $a\in \mathbb{R}$}. \label{cauchyDis}\\

\noindent $X\sim C(a,\lambda)$ in defined on $\mathbb{R}$.\\

\noindent Absolutely continuous probability density function on $\mathcal{V}_X=\mathbb{R}$ :

\begin{equation*}
f_{X}(x)=\frac{\lambda}{ \pi (\lambda^2+(x-a)^{2})}, \ x \in \mathbb{R}. \ \label{cauchyDisPdf}
\end{equation*}

\bigskip \noindent Distribution function : \\

$$
F_X(x)=\frac{1}{\pi} \left(\arctan \left(\frac{x-a}{\lambda}\right)-\frac{\pi}{2}\right), \ \in \mathbb{R}.
$$

\bigskip \noindent Characteristic function : \\

$$
\Phi_X(u)=\exp(iua - \lambda |u|), \ u \in \mathbb{R}.
$$

\bigskip \noindent The mathematical expectation does not exist.\\

\noindent \textbf{Proof}.  We have to prove that for $a=0$ and $\lambda=1$, Formula \ref{cauchyDisPdf} gives a \textit{pdf}. Indeed, using that the primitive of $(1+x^2)$ is $\arctan x$, the inverse of the tangent function $\tan x$, we have

$$
\int_{-\infty}^{+\infty} \frac{dx}{\pi(1+x^2)}=\frac{1}{\pi} \left[\arctan x\right]_{-\infty}^{+\infty}=1.
$$

\bigskip \noindent Next, setting $X=\lambda Z + a$, where $Z$ follows a $C(0,1)$ law leads to the general case in \ref{cauchyDisPdf} by differentiating
$F_Z(x)=F_X((x-a)/\lambda)$, $x \in \mathbb{R}$.\\

\noindent The expression of the characteristic function of a \textit{rrv} $Z$ following a standard Cauchy law is given by Formula \ref{fc_cauchy00} (Chapter \ref{proba_02_rv}, page \pageref{fc_cauchy00}). By the transform $X=\lambda Z + a$, we have the general characteristic function of a Cauchy distribution.\\

\noindent Finally, we have for $a=0$ and $\lambda=1$,

$$
\mathbb{E}(X^+)=\int_{0}^{+\infty} \frac{x}{\pi(1+x^2)} dx=+\infty,
$$

\bigskip \noindent and 

$$
\mathbb{E}(X^-)=\int_{-\infty}^{0} \frac{x}{\pi(1+x^2)} dx=-\infty,
$$

\bigskip \noindent and then $\mathbb{E}(X)$ is not defined. Concerning that point, we recommend to go back to the remark concerning the caution to take while using the improper Riemann integration at the place of the Lebesgue integral (See Point (b5) in Section \ref{proba_02_rv_sec_05} in Chapter \ref{proba_02_rv}, page \pageref{LLBRI}).\\

\bigskip \noindent \textbf{(8) Logistic Random Variable with parameters $a\in \mathbb{R}$ and $b>0$}.\\

\noindent  $X\sim \ell(a,b)$ is supported by the whole real line.\\

\noindent Absolutely continuous probability density function on $\mathcal{V}_X=\mathbb{R}$ :

\begin{equation*}
f_{X}(x)=b^{-1}e^{-(x-a)/b}/(1+e^{-(x-a)/b}),\ x\in \mathbb{R}.
\end{equation*}

\bigskip \noindent Characteristic function :

\begin{equation*}
\Phi _{X}(u)=e^{iau}\pi b\text{ }cosec(i\pi bu).
\end{equation*}

\bigskip \noindent Mathematical expectation and variance :
 
\begin{equation*}
E(X)=a;V(X)=b^{2}\pi^{2}/3.
\end{equation*}

\bigskip \noindent \textbf{(9) Weibull Random Variable  with parameters $a>0$ and $b>0$}.\\

\noindent  $X\sim W(a,b)$ is supported by $\mathbb{R}_{+}$.\\

\noindent Absolutely continuous probability density function on $\mathcal{V}_X=\mathbb{R}_{+}$ :\\

\begin{equation*}
f_{X}(x)=ab\text{ }x^{b-1}\exp (-ax^{-b})1_{(x>0)}.
\end{equation*}

\bigskip \noindent Characteristic function :
\begin{equation*}
\Phi _{X}(u)=a^{-iu/b}\Gamma (1+iu/b), \ u \in \mathbb{R}.
\end{equation*}

\bigskip \noindent Mathematical expectation and variance :

\begin{equation*}
\mathbb{E}(X)=(1/a)^{1/b}\Gamma (1+1/b);V(X)=a^{-2/b}(\Gamma (1+2/b)-\Gamma (1+1/b)).
\end{equation*}

\bigskip \noindent \textbf{(10) Gumbel Random Variable $a\in \mathbb{R}$ and $b>0$}.\\

\noindent  $X\sim Gu(a,b)$ is supported by the whole line $\mathbb{R}$.\\

\noindent Absolutely continuous probability density function on $\mathcal{V}_X=\mathbb{R}$ :

\begin{equation*}
f_{X}(x)=(u/b)e^{-u},\text{ }with\text{ }u=e^{-(x-a)/b}.
\end{equation*}

\bigskip \noindent Characteristic function :

\begin{equation*}
\Phi _{X}(u)=e^{iua}\Gamma (1-ibu), \ u \in \mathbb{R}.
\end{equation*}

\bigskip \noindent Mathematical expectation and variance :

\begin{equation*}
\mathbb{E}(X)=a+\gamma b.
\end{equation*}

\bigskip \noindent where $\gamma=  $ is the Euler's number and 

\begin{equation*}
Var(X)=\pi ^{2}b^{2}/2.
\end{equation*}%

\bigskip \noindent \textbf{(11) Double-exponential Random Variable with parameter $b>0$}.\\

\noindent See Point (4) above.\\

\noindent  $X\sim \mathcal{E}_{d}(b)$ is defined on the whole real line.\\

\noindent Absolutely continuous probability density function on $\mathcal{V}_X=\mathbb{R}$ :

\begin{equation*}
f_{X}(x)=\frac{b}{2}\exp (-b \left\vert x\right\vert ),x\in \mathbb{R}.
\end{equation*}

\bigskip \noindent Characteristic function :

\begin{equation*}
\Phi _{X}(u)=(1+(u/b)^{2})^{-1}, \ \ u \in \mathbb{R}.
\end{equation*}

\bigskip \noindent Moments of order $k\geq 1$ :

\bigskip \noindent Mathematical expectation and variance :

\begin{equation*}
\mathbb{E}(X)=0, \text{ and }Var(X)=2b^{-2}.
\end{equation*}

\bigskip \noindent \textbf{(12) Gaussian Random Variable with parameters $m\in \mathbb{R}$ and $\sigma>0$}. $X\sim \mathcal{N}(m,\sigma ^{2})$ is supported by the whole real line $\mathcal{V}_X=\mathbb{R}$.\\

\noindent Absolutely continuous probability density function on $\mathcal{V}_X=\mathbb{R}$ :

\begin{equation*}
f_{X}(x)=\frac{1}{\sqrt{2\pi }}\exp (-(x-m)^{2}/\sigma ^{2}),x\in \mathbb{R}.
\end{equation*}

\bigskip \noindent Characteristic function :

\begin{equation*}
\Phi _{X}(u)=e^{-um}\exp (-\sigma ^{2}u^{2}/2).
\end{equation*}

\bigskip \noindent Moments of order $k\geq 1$.

$$
\mathbb{X} \left( \frac{X-m}{\sigma^2} \right)=\frac{2^{k} \ k!}{(2k)}.
$$

\bigskip \noindent Mathematical expectation and variance :
\begin{equation*}
\mathbb{E}(X)=m\text{, }Var(X)=\sigma ^{2}.
\end{equation*}

\bigskip \noindent Because of its importance in the history of Probability Theory, as explained by its name of \textit{normal} probability law, we will devote a special study to it in Chapter
\ref{proba_02_gauss}.

\bigskip \noindent \textbf{(13) Chi-square Probability law of parameters $d\geq 1$}.\\

\noindent  $X\sim \chi_d^2$ is supported by $\mathcal{V}_X=\mathbb{R}_+$.\\

\noindent \textbf{Definition} A Chi-square Probability law of $d\geq 1$ degrees of freedom is simply a Gamma law of parameters $a=d/2$ and $b=1/2$, that is

$$
\chi_d^2=\gamma(d/2,1/2).
$$

\bigskip \noindent By reporting the results of $\gamma$-laws, we have the following facts.\\

\noindent Absolutely continuous probability density function on $\mathcal{V}_X=\mathbb{R}_{+}$ :\\

\noindent Absolutely continuous probability density function on $\mathcal{V}_X=\mathbb{R}_{+}$ :
\begin{equation*}
f_{X}(x)=\frac{1}{2^{d/2}\Gamma (d/2)}x^{\frac{d}{2}-1}e^{-\frac{x}{2}}1_{(x\geq 0)}.
\end{equation*}

\bigskip \noindent Characteristic function :

\begin{equation*}
\Phi _{X}(u)=(1-i2u)^{-d/2}.
\end{equation*}

\bigskip \noindent Moments of order $k\geq 1$.

$$
\mathbb{E}(X^k)= \frac{1}{2^{-k}} \prod_{j=0}^{k} \left(\frac{d}{2}+j\right).
$$

\bigskip \noindent Mathematical expectation and variance :
\begin{equation*}
\mathbb{E}(X)=d,\ Var(X)=2d.
\end{equation*}

\bigskip 
\noindent \textbf{Important properties}. Chi-square distributions are generated from Gaussian random variables as follows.\\

\noindent \textbf{Fact 1}. If $Z$ follows a standard Gaussian probability law, $Z^2$ follows a Chi-square law of one degree of freedom :

$$
Z \sim \mathcal{N}(0,1) \Rightarrow Z^2 \sim \chi_1^2.
$$

\bigskip \noindent \textbf{Proof}. Suppose that $Z \sim \mathcal{N}(0,1)$ and put $X=Z^2$. It is clear that the domain of $Y$ is $\mathcal{V}_X=\mathbb{R}_+$. For any $y \geq 0$,

\begin{eqnarray*}
F_X(x)&=&\mathbb{P}(Z^2\leq x)\\
&=&\mathbb{P}(|Z| \leq \sqrt{x})\\
&=&\mathbb{P}(Z \in ]-\infty, \sqrt{x}] \setminus ]-\infty, \sqrt{x}[)\\
&=&\mathbb{P}(Z \leq  \sqrt{x}) - \mathbb{P}(Z \leq  -\sqrt{x}).
\end{eqnarray*}

\bigskip \noindent Remind that $Z$ has an even absolutely continuous \textsl{pdf} $f_Z$. This implies that  $\mathbb{P}(Z=t)=0$ for any $t \in \mathbb{R}$ and we get for any $y \geq 0$,

$$
F_Y(x)=F_Z(\sqrt{x}) - F_Z(-\sqrt{x}).
$$

\bigskip \noindent By differentiating by $x$, we get the absolutely continuous \textsl{pdf} of $X$ for any any $x \in \mathcal{V}_X$

$$
f_X(x)=\frac{1}{2\sqrt{x}}\biggr(f_Z(\sqrt{x})+f_Z(-\sqrt{x})\biggr)=\frac{1}{\sqrt{x}} f_Z(\sqrt{x}),
$$

\bigskip \noindent which leads to

$$
f_X(x)=\frac{\left(\frac{1}{2}\right)^{1/2}}{\sqrt{\pi}} x^{1-1/2} \exp(-x/2), \ \ x \in \mathbb{R}_+.
$$

\bigskip  \noindent By comparing with the absolutely continuous \textit{pdf} $f_{\chi^1_2}$ of a Chi-square probability law, we see that $f_{\chi^1_2}$ and $f_X$ are two absolutely continuous \textbf{pdf}'s with the same support $\mathcal{V}$ and a common variable part 

$$
h(x)=x^{1-1/2} \exp(-x/2),  \ \ x\in \mathcal{V}.
$$

\bigskip \noindent By the Easy Stuff remark in Section \ref{proba_02_rv_sec_08} in Chapter Section \ref{proba_02_rv_sec_08}, it follows that they are equal and by the way, we get the stunning equality

$$
\Gamma(1/2)=\sqrt{\pi}.
$$

\bigskip \noindent \textbf{Fact 2}. Let $d\geq 2$. The convolution product of $d$ Chi-square law of one degree of freedom is a Chi-square law probability of $d$ degrees of freedom. In particular, if $X_1$, ..., $X_d$ are $d$ independent real-valued random variables, defined on the same probability space, identically following a Chi-square law of one degree of freedom, we have

$$
X_1^2+ \cdots +  X_d^2 =\sum_{1\leq i \leq d} X_i^2 \sim \chi_d^2. 
$$

\bigskip \noindent Indeed, if $X_1$, ..., $X_d$ are independent and identically follow a Chi-square law of one degree of freedom, the probability law $X_1^2+ \cdots +  X_d^2$ is characterized by its characteristic function

$$
\varphi_{X_1^2+ \cdots +  X_d^2 }(t)= \prod_{j=1}^{d}=\varphi_{X_j}(t)=(1-2it)^{d/2},
$$ 

\bigskip \noindent which establishes that $X_1^2+ \cdots +  X_d^2$ follows a $\chi_d^2$ law.\\
 
\bigskip \noindent \textbf{(14) Around the Normal Variance Mixture class of random variables}\\

\noindent We are introducing some facts on this class of random variables which are important tools in financial data statistical studies. We only provide some of their simple features, not dwelling in their deep relations. It is expected to treat these random variables in completion of Chapter \label{proba_02_gauss} later.\\

\noindent \textbf{(a)}. A normal variance mixture is defined as follows :

$$
X = \mu + \sigma \sqrt{W} Z, \ (NMV)
$$

\bigskip \noindent where $Z$ is a standard random variable, $W$ is a positive random variable defined on the same space as $Z$ and independent of
$Z$, $\mu$ is a real number (the mean of $X$) and $\sigma$ is a positive random variable. Hence we have

$$
\mathbb{E}(X)=\mu + \sigma \mathbb{E}(\sqrt{W}) \mathbb{E}(Z)=\mu 
$$

\bigskip \noindent and

$$
\mathbb{V}(X)= \sigma^2 \mathbb{E}(W^2) \mathbb{E}(Z^2)=\sigma^2 \mathbb{E}(W^2). 
$$

\bigskip \noindent The following distributions of $W$ are generally used.\\

\noindent \textbf{(b)}. The inverse gamma law $W=1/Y \sim Ig(\alpha, \beta)$, where $Y \sim \gamma(\alpha,\beta)$, $\alpha>0$, $\beta>0$. \label{inversegamma}\\

\noindent Absolutely continuous probability density function on $\mathcal{V}_X=\mathbb{R}_{+}$ :\\

$$
f_{W}(x)=\frac{\beta^{\alpha}}{\Gamma(\alpha)} x^{-\alpha-1} \exp(-\beta/x), \ x>0.
$$

\bigskip \noindent Mathematical expectation and variance :

$$
\mathbb{E}(X)=\frac{\beta}{\alpha-1} \ for \ \alpha>1;  \mathbb{V}ar(X)=\frac{\beta^2}{(\alpha-1)^2(\alpha-2)} \ for \ \alpha>2.
$$

\noindent \textbf{(c)}. The Generalized Inverse Gaussian (GIG) law : $W \sim Gig(a,b,c)$, $(a,b,c)\in  \mathbb{R}_{+}^3$. \label{gig}\\

\noindent Parameters domains :

\begin{eqnarray*}
b>0 \ \ and \ \ c\geq 0 &if& a<0\\
b>0 \ \ and \ \ c>0 &if& a=0\\
b\geq 0 \ \ and \ \ c>0 &if& a>0.\\
\end{eqnarray*}

\noindent Absolutely continuous probability density function on $\mathcal{V}_X=\mathbb{R}_{+}$ :\\

$$
f_{W}(x)=\frac{b^{-a}(bc)^{a}}{2K_{a}\left(\left(bc\right)^{1/2}\right)} x^{a-1} e^{-(cx+b/x)/2}, \ x\in \mathbb{R}.
$$

\noindent 

$$
K_{a}\left(\left(bc\right)^{1/2}\right)=\frac{b^{-a}(bc)^{a}}{2} \int_{0}^{+\infty} x^{a-1} e^{-(cx+b/x)/2} dx.
$$

\bigskip \noindent This function, called a \textit{modified Bessel function}, is not directly defined in a simple argument, but on a composite argument $\left(bc\right)^{1/2}$ and one should pay a particular attention to the simultaneous domain of the parameters $(a,b,c)$.\\

\bigskip \noindent Mathematical expectation and variance :

$$
\mathbb{E}(X)=\frac{\beta}{\alpha-1} \ for \ \alpha>1;  \mathbb{V}ar(X)=\frac{\beta^2}{(\alpha-1)^2(\alpha-2)} \ for \ \alpha>2.
$$

\bigskip \noindent \textbf{(d)}. Student distribution of $\nu\geq 1$ degrees of freedom : $X \sim t(\nu)$. \label{student0}\\

\noindent If in Formula (NMV), we take $W$ as inverse Gamma random variable $Ig(\nu/2,\nu/2)$, where $\nu\geq 1$ is an integer, the \textit{pdf} of $X$ becomes :\\

\noindent for $x \in \mathcal{V}_X=\mathbb{R}$ :

$$
f(x)=\frac{\Gamma((\nu+1)/2)}{\sigma \Gamma(\nu/2)(\nu\pi)^{1/2}}\biggr(1+\frac{(x-\mu)^2/2}{\nu}\biggr)^{-(\nu+1)/2}.
$$

\bigskip \noindent \textbf{(e)}. Symmetric Generalized Hyperbolic distribution : $X \sim SGH(\mu, a,b,c)$. \label{SGH}\\

\noindent Parameters : $\mu \in \mathbb{R}$, $a$, $b$ and $c$ given in the \textit{Gig} law presentation.\\

\noindent If, in Formula (NMV), we take $W$ as the generalized inverse Gaussian random variable $Gig(a,b,c)$, the \textit{pdf} of $X$ becomes :\\

\noindent for $x \in \mathcal{V}_X=\mathbb{R}$ :

$$
f(x)=\frac{(ab)^{-a/2}c^{1/2}}{\sigma (2\pi)^{1/2}K_{a}\left((bc)^{1/2}\right)}
\frac{K_{a-1/2}\left( \left(b+c(x-mu)^2/\sigma \right)^{1/2} \right)}{\left(b+c(x-mu)^2/\sigma\right)^{1/4-a/2}}.
$$ 

\bigskip 
\noindent \textbf{(f)}. Generalized Hyperbolic distribution : $X \sim GH(\mu, a,b,c)$.\\

\noindent The latter probability law is a particular case of the following model : \label{GNMV}

$$
X = \mu + \gamma W + \sigma \sqrt{W} Z, \ (GNMV)
$$

\bigskip \noindent for $\gamma=0$. If $\gamma$ is an arbitrary real number and we take $W$ as a generalized inverse Gaussian random variable $Gig(a,b,c)$, the \textit{pdf} of $X$ is  :\\

\noindent for $x \in \mathcal{V}_X=\mathbb{R}$ :

$$
f(x)=c \frac{\exp(\gamma(x-\mu)/\sigma)K_{a-1/2}\left( \left(b+\sigma^{-1}(x-mu)^2(c+\gamma^2/\sigma) \right)^{1/2}\right)}{\left(b+\sigma^{-1}(x-mu)^2(c+\gamma^2/\sigma)\right)^{1/4-a/2}},
$$ 

\noindent where

$$
c=\frac{(ab)^{-a/2}c^{c}(c+\gamma^2/\sigma)^{1/2-a}}{\sigma (2\pi)^{1/2}K_{a}\left((bc)^{1/2}\right)}.
$$ 

\bigskip
\noindent \textit{Comments}. This part (14) was only an introduction to an interesting modern and broad topic in Statistical studies in Finance. The multivariate version has also been developed.\\

\bigskip \noindent \textbf{(15) Probabiliy Laws of the Gaussian sample}.\\

\noindent In Mathematical Statistics, the study Gaussian samples holds a special place, at least at the beginning of the exposure of the theory. The following probability laws play the major roles.\\

\noindent \textit{(a) The \textit{Chi-square} probability law of $n\geq 1$ degrees of freedom}.\\

\noindent $X \sim \chi_{2}^{n}$.\\

\noindent This law has been introduced in Point (13) above.\\

\noindent \textit{(b) The \textit{Student} probability law of $n\geq 1$ degrees of freedom}. \label{student}\\

\noindent  $X \sim t(n)$ is defined on the whole real line.\\

\noindent Absolutely continuous probability density function on $\mathcal{V}_X=\mathbb{R}$ :

\begin{equation*}
f_{X}(x)=\frac{\Gamma((n+1)/2)}{(n\pi)^{1/2}\Gamma(n/2)} \left(1 + \frac{x^2}{n}\right)^{-(n+1)/2}
\end{equation*}

\bigskip \noindent Characteristic function. No explicit form.\\

\bigskip \noindent Moments of order $k\geq 1$ :

\bigskip \noindent Mathematical expectation and variance :

\begin{equation*}
\mathbb{E}(X)=0, \text{ and }Var(X)=\frac{n}{n-2}), \ n\geq 3.
\end{equation*}

\bigskip \noindent \textit{(c) The Fisher probability law of degrees of of freedom $n\geq 1$ and $m\geq 1$}. \label{fisher}\\
\noindent  $X\sim F(n,m))$ is defined on the positive real line.\\

\noindent Absolutely continuous probability density function on $\mathcal{V}_X=\mathbb{R}_{+}$ :

\begin{equation*}
f_{X}(x)=\frac{n^{n/2}m^{m/2}\Gamma((n+m)/2)}{\Gamma(n/2)\Gamma(m/2)} \frac{x^{n/2-1}}{(m+nx)^{(n+m)/2}}.
\end{equation*}

\bigskip \noindent Characteristic function. No explicit form.\\

\noindent Mathematical expectation and variance :

\begin{equation*}
\mathbb{E}(X)=\frac{m}{m-2}, \ m\geq 3 \text{ and }Var(X)=\frac{2m^2 (n+m-2)}{n(m-2)^2(m-4)}, \ m\geq 5.
\end{equation*}

\bigskip 
\noindent We take this opportunity to propose an exercise which illustrate the change of variable formula given in page \pageref{changeVariable} and which allows to find the just given laws.

%% file: proba_02_03_exo_ang.tex
\begin{exercise} \label{exo_proba_02_01} Let $\left( X,Y\right)$ be a $2$-random vector with \textit{pdf} $f_{\left( X,Y\right)}$, on its support $D$ with respect to the Lebesgue measure on $\mathbb{R}^2$. Consider the following transform
$$
(x,y) \mapsto h\left( x,y\right) =\left(x,x+y\right) \in \Delta,
$$

\noindent that is : \\

\begin{center}
$\left\{ 
\begin{array}{c}
U=X \\ 
V=X+Y%
\end{array}
\right.$ 
\end{center}

\bigskip 
\noindent (a) Find the law of $(U,V)$ and their marginals law.\\

\noindent (b) Precise the \textit{pdf} of $V$ if $X$ and $Y$ are independent.\\
 
\noindent (c) Application : Let $X\sim \gamma \left( \alpha ,b\right)$ and $Y\sim \gamma
\left( \beta ,b\right)$. Show that $V=X+Y\sim \gamma \left(\alpha +\beta ,b\right)$.
\end{exercise}

\begin{exercise} \label{exo_proba_02_02} Let Let $\left( X,Y\right)$ be a $2$-random vector with \textit{pdf} $f_{\left( X,Y\right)}$ with respect to the Lebesgue measure on $\mathbb{R}^2$. Consider the following transfrom
$$
(x,y) \mapsto h\left( x,y\right) =\left(x/y,y\right),
$$

\noindent that is

$$
\left( X,Y\right) \mapsto (U,V)=\left(X/Y,Y\right).
$$

\bigskip \noindent (a) Apply the general change of variable formula to write the \textit{pdf} of $(X/Y,Y)$ and deduce the marginal \textit{pdf} of $U=X/Y$.\\

\noindent (b) Precise it for $X$ and $Y$ independent.\\

\noindent (c) In what follows, $X$ and $Y$ are independent. Precise the \textit{pdf} of $U$ when $X$ and $Y$
are both standard Gaussian random variables. Identify the found probability law.\\ 

\noindent (d) Let $X$ be a standard Gaussian random variable and $Y=Z^{1/2}$ the square-root of a $\chi_{n}^{2}$ random variable with $n\geq 1$. Precise the \textit{pdf} of $U=X/\sqrt{Z}$. Begin to give the \textit{pdf} of $Y$ by using its \textit{cdf}.\\

\noindent Deduce from this the probability law of

$$
t(n)=\frac{\sqrt{n}X}{Y}\equiv\frac{\mathcal{N}(0,1)}{\sqrt{\chi_{2}^{n}/n}}
$$

\bigskip \noindent where the term after the sign $\equiv$ is a rephrasing of a ratio of two independent random variables : a $\mathcal{N}(0,1)$ random variable by the square root of a chi-square random variable divided by its number of freedom degrees.\\

\noindent Conclude that a \textit{t(n)}-random variable has the same law as the ratio of two independent random variables : a $\mathcal{N}(0,1)$ random variable by by square root of a chi-square random variable divided by its number of freedom degrees.\\

\noindent (e) Let $X$ and $Y$ be two independent random variable following Chi-square laws of respective number of freedom degrees $n\geq 1$ and $m\geq 1$. Precise the \textit{pdf} of $U=X/Y$.\\

\noindent Deduce from this the probability law of

$$
F_{n,m}=\left(\frac{m}{n}\right)\frac{X}{Y}\equiv\frac{\chi_{2}^{n}/n}{\chi_{2}^{m}/m},
$$

\noindent where the term after $\equiv$ is a rephrasing of a ratio of two independent random variables : a Chi-square random variable of number of freedom degrees $n\geq 1$ by a Chi-square random variable of number of freedom degrees $m\geq 1$.\\

\noindent Conclude that a Fisher random variable with numbers of freedom degrees $n\geq 1$ and $m\geq 1$ has 
the same probability law as a ratio of two independent random variables : a Chi-square random variable of number of freedom degrees $n\geq 1$ by a Chi-square random variable of number of freedom degrees $m\geq 1$.
\end{exercise}

\noindent \textbf{Solutions of Exercise \ref{exo_proba_02_01}}. We have the transformation :\\

\begin{center}
$\left\{ 
\begin{array}{c}
u=x \\ 
v=x+y%
\end{array}%
\right. \Leftrightarrow \left\{ 
\begin{array}{c}
x=u \\ 
y=-u+v%
\end{array}%
\right.$ 
\end{center}

\bigskip \noindent The Jacobian matrix is :\\

\begin{center}
$J_{\left( u,v\right) }=\left( 
\begin{array}{cc}
1 & 0 \\ 
-1 & 1%
\end{array}%
\right)$
\end{center}

\bigskip \noindent with determinant$\det \left(J_{\left( u,v\right) }\right) =1$. By the Change of variable formula, we have
$$
f_{\left( U,V\right) }\left( u,v\right)=
f_{\left( X,Y\right) }\left(u,-u+v\right) 1_{\Delta}\left( u,v\right).
$$

\bigskip \noindent The marginal laws are

$$
f_{U}\left( u\right) =\int_{D_{v}}f_{\left( U,V\right) }\left( u,v\right) dv? \ U \in \mathcal{V}_U
$$

\noindent and

\begin{eqnarray*}
f_{V}\left( v\right) &=&\int_{D_{u}}f_{\left( U,V\right) }\left( u,v\right)
du \\
&=&\int_{D_{U}}f_{\left( X,Y\right) }\left( u,-u+v\right) \mathbb{I}%
_{D\left( U,V\right) }\left( u,v\right) du. 
\end{eqnarray*}

\noindent \textbf{Question (b)} We have
$$
f_{V}(v)=f_{X+Y}\left( v\right) =\int f_{X}\left( u\right) f_{Y}\left( v-u\right)
du=\int f_{Y}\left( u\right) f_{X}\left( v-u\right) du,
$$

\noindent which is the convolution product between $X$ and $Y$.\\
 
\noindent \textbf{Question (c)}. We recall that

$$
f_{X}\left( x\right) =\frac{b^{\alpha }}{\Gamma \left( \alpha \right) }%
x^{\alpha -1}e^{-bx}\mathbb{I}_{\mathbb{R}_{+}}\left( x\right)
$$

\noindent and 

$$
f_{Y}\left( x\right) =\frac{b^{\beta }}{\Gamma \left( \beta \right)}x^{\beta -1}e^{-bx}\mathbb{I}_{\mathbb{R}_{+}}\left( x\right).
$$

\bigskip \noindent We have 

\begin{eqnarray*}
f_{X+Y}\left( v\right)  &=&\frac{b^{\alpha }}{\Gamma \left( \alpha \right) }%
\frac{b^{\beta }}{\Gamma \left( \beta \right) }\int_{0}^{v}u^{\alpha
-1}e^{-bu}\left( v-u\right) ^{\beta -1}e^{-b\left( v-u\right) }du \\
&=&\frac{b^{\alpha +\beta }}{\Gamma \left( \alpha \right) \Gamma \left(
\beta \right) }\times v^{\alpha +\beta -2}e^{-bv}\int_{0}^{v}\left( \frac{u}{%
v}\right) ^{\alpha -1}\left( 1-\frac{u}{v}\right) ^{\beta -1}du.
\end{eqnarray*}

\bigskip \noindent By taking the further change of variables $x=u/v$, we get

$$
f_{X+Y}\left( v\right) =\left[ \frac{b^{\alpha +\beta }}{\Gamma \left(
\alpha \right) \Gamma \left( \beta \right) }\int_{0}^{1}\left( x\right)
^{\alpha -1}\left( 1-x\right) ^{\beta -1}dx\right] \times v^{\alpha +\beta
-2}e^{-bv}\mathbb{I}_{\mathbb{R}_{+}}\left( v\right).
$$

\bigskip
\noindent Since $X+Y$ has the same domain and the same variable part of a $\gamma(\alpha+\beta,b)$, they have the same 
constant and then we have $X+Y \sim \gamma(\alpha+\beta,b)$ and 
$$
\frac{b^{\alpha+\beta }}{\Gamma \left( \alpha +\beta \right) }=\frac{b^{\alpha +\beta }}{%
\Gamma \left( \alpha \right) \Gamma \left( \beta \right) }\times \beta
\left( \alpha ,\beta \right)
$$

\noindent with 

$$
B\left( \alpha ,\beta \right) =\int_{0}^{1}\left( x\right) ^{\alpha
-1}\left( 1-x\right) ^{\beta -1}dx=\frac{\Gamma \left( \alpha \right) \Gamma
\left( \beta \right) }{\Gamma \left( \alpha +\beta \right)}.
$$

\noindent $\square$\\

\noindent \textbf{Solutions of Exercise \ref{exo_proba_02_02}}.\\

\noindent Question (a). We have the transformation : \\

\begin{center}
$\left\{ 
\begin{array}{c}
u=\frac{x}{y} \\ 
v=y%
\end{array}%
\right. \Rightarrow \left\{ 
\begin{array}{c}
x=uv \\ 
y=v%
\end{array}
\right)$
\end{center}

\bigskip \noindent The the Jacobian matrix is 
$$
J_{\left( u,v\right) }=\left( 
\begin{array}{cc}
v & u \\ 
1 & 0%
\end{array}%
\right)$$

\noindent with determinant $\det \left( J_{\left( u,v\right) }\right) =v$. The \textit{pdf} of $(U,V)$ becomes

$$
f_{\left( \frac{X}{Y},Y\right) }\left( u,v\right) =f_{\left(
X,Y\right) }\left( uv,v\right) \left\vert v\right\vert \mathbb{I}_{D\left( 
\frac{X}{Y},Y\right) }\left( u,v\right),
$$

\noindent and the marginal law of $V=X/Y$ is

\begin{eqnarray*}
f_{\frac{X}{Y}}\left( u\right)  &=&\int_{D_{v}}f_{\left( U,V\right) }\left(
u,v\right) dv \\
&=&\int_{D_{v}}f_{\left( X,Y\right) }\left( uv,v\right) \left\vert
v\right\vert \mathbb{I}_{D\left( \frac{X}{Y},Y\right) }\left( u,v\right) dv%
\text{ \ \ \ \ }\left( 1\right). 
\end{eqnarray*}

\noindent Question (b) If $X$ and $Y$ are independent, we have

\begin{eqnarray*}
f_{\frac{X}{Y}}\left( u\right)  &=&\int_{D_{v}}f_{\left( U,V\right) }\left(
u,v\right) dv \\
&=&\int_{\left( uv\in D_{X},v\in D_{Y}\right) }f_{X}\left( uv\right)
f_{Y}\left( v\right) \left\vert v\right\vert dv.
\end{eqnarray*}

\noindent Question (c) Now if $X$ and $Y$ are standard Gaussian random variables, we get
\begin{eqnarray*}
f_{\frac{X}{Y}}\left( u\right)  &=&\frac{1}{2\pi }\int_{-\infty }^{\infty
}e^{-\frac{1}{2}u^{2}v^{2}}e^{-\frac{1}{2}v^{2}}\left\vert v\right\vert dv \\
&=&\frac{1}{2\pi }\int_{-\infty }^{\infty }e^{-\frac{1}{2}\left(
u^{2}+1\right) v^{2}}\left\vert v\right\vert dv \\
&=&\frac{1}{\pi }\int_{0}^{\infty }ve^{-\frac{1}{2}\left( u^{2}+1\right)
v^{2}}dv \\
&=&\frac{1}{\pi }\left[ -\frac{e^{-\frac{1}{2}\left( u^{2}+1\right) v^{2}}}{%
u^{2}+1}\right] _{0}^{\infty } \\
&=&\frac{1}{\pi \left( u^{2}+1\right)},
\end{eqnarray*}

\noindent \noindent which is the standard Cauchy probability law.\\

\noindent Question (d). The \textit{pdf} of $Y$ is easily derived from the relation

$$
\forall y\geq 0, \ F_Y(y)=F_{Z}(y^2)
$$

\noindent \noindent which, after differentiation, gives

\begin{equation*}
f_{Y}\left( y\right) =2yf_{Y_{1}}\left( y^{2}\right) =\frac{2\left( \frac{1}{%
2}\right) ^{\frac{n}{2}}}{\Gamma \left( \frac{n}{2}\right) }y^{n-1}e^{-\frac{%
1}{2}y^{2}}1_{\mathbb{R}_{+}}\left( y\right). 
\end{equation*}

\noindent From there, we have

\begin{eqnarray*}
f_{\frac{X}{Y}}\left( u\right)  &=&\int f_{X}\left( uv\right) f_{Y}\left(
v\right) \left\vert v\right\vert dv\text{ \ \ } \\
&=&\frac{1}{\sqrt{2\pi }}\frac{1}{2^{\frac{n}{2}-1}\Gamma \left( \frac{n}{2}%
\right) }\int_{0}^{\infty }v^{n}e^{-\frac{1}{2}u^{2}v^{2}}e^{-\frac{v^{2}}{2}%
}dv.
\end{eqnarray*}

\bigskip 
\noindent Let us set 

\begin{eqnarray*}
A &=&\int_{0}^{\infty }v^{n}e^{-\frac{1}{2}u^{2}v^{2}}e^{-\frac{v^{2}}{2}}dv%
\text{\ } \\
&=&\int_{0}^{\infty }v^{n}e^{-\frac{v^{2}}{2}\left( u^{2}+1\right) }dv\text{%
\ },
\end{eqnarray*}

\noindent and make the change of variable $t=\frac{v^{2}}{2}\left( u^{2}+1\right)$. Then we have

$$
v=\sqrt{\frac{2t}{u^{2}+1}}=\sqrt{\frac{2}{u^{2}+1}}t^{\frac{1}{2}}
$$

\noindent and

\begin{eqnarray}
dv&=&\frac{1}{2}\sqrt{\frac{2}{u^{2}+1}}t^{-\frac{1}{2}}dt\\
&=&\frac{1}{2}\sqrt{\frac{2}{u^{2}+1}}\times \frac{1}{v}\times \sqrt{\frac{2}{u^{2}+1}}dt\\
&=&\frac{dt}{u^{2}+1}\times \frac{\left( u^{2}+1\right) ^{\frac{1}{2}}}{\sqrt{2}}t^{-\frac{1}{2}}.
\end{eqnarray}

\noindent Next, we have

\begin{eqnarray*}
A &=&\left( \frac{2}{u^{2}+1}\right) ^{\frac{n}{2}}\frac{1}{\left[ 2\left(
u^{2}+1\right) \right] ^{\frac{1}{2}}}\int_{0}^{\infty }t^{\frac{n+1}{2}%
-1}e^{-t}dt \\
&=&\frac{2^{\frac{n}{2}}}{\sqrt{2}\left( u^{2}+1\right) ^{\frac{n}{2}+\frac{1%
}{2}}}\Gamma \left( \frac{n+1}{2}\right), 
\end{eqnarray*}

\noindent which leads to 

\begin{eqnarray*}
f_{\frac{X}{Y}}\left( u\right)  &=&\frac{1}{\sqrt{2\pi }}\frac{1}{2^{\frac{n%
}{2}-1}\Gamma \left( \frac{n}{2}\right) }\times \frac{2^{\frac{n}{2}}}{\sqrt{%
2}\left( u^{2}+1\right) ^{\frac{n}{2}+\frac{1}{2}}}\Gamma \left( \frac{n+1}{2%
}\right)  \\
&=&\frac{\Gamma \left( \frac{n+1}{2}\right) }{\sqrt{\pi }\Gamma \left( \frac{%
n}{2}\right) }\left( u^{2}+1\right) ^{-\frac{n+1}{2}}
\end{eqnarray*}

\noindent Finally, by taking $W=\sqrt{n}U=\sqrt{n}X/Y$, the pdf of $W$ is

\begin{equation*}
f_{\frac{\sqrt{n}X}{Y}}\left( u\right) =\frac{\Gamma \left( \frac{n+1}{2}\right) }{%
\sqrt{n\pi }\Gamma \left( \frac{n}{2}\right) }\left(1 + \frac{u^2}{n}\right) ^{(n+1)/2}.
\end{equation*}

\bigskip 
\noindent Question (e). Using the right expressions for the Chi-square random variable leads to

\begin{eqnarray*}
f_{\frac{X}{Y}}\left( u\right)  &=&\int f_{X}\left( uv\right) f_{Y}\left(
v\right) \left\vert v\right\vert dv\text{ } \\
&=&\frac{\left( \frac{1}{2}\right) ^{\frac{n}{2}}\left( \frac{1}{2}\right) ^{%
\frac{m}{2}}}{\Gamma \left( \frac{n}{2}\right) \Gamma \left( \frac{m}{2}%
\right) }\int_{0}^{\infty }\left( uv\right) ^{\frac{n}{2}-1}e^{-\frac{1}{2}%
uv}v^{\frac{m}{2}-1}e^{-\frac{1}{2}v}dv \\
&=&\frac{u^{\frac{n}{2}-1}}{2^{\frac{n+m}{2}}\Gamma \left( \frac{n}{2}%
\right) \Gamma \left( \frac{m}{2}\right) }\int_{0}^{\infty }v^{\frac{n+m}{2}%
-1}e^{-\frac{1}{2}v\left( u+1\right) }dv \\
&=&\frac{u^{\frac{n}{2}-1}}{2^{\frac{n+m}{2}}\Gamma \left( \frac{n}{2}%
\right) \Gamma \left( \frac{m}{2}\right) }\times \frac{\Gamma \left( \frac{%
n+m}{2}\right) }{\left( u+1\right) ^{\frac{n+m}{2}}} \\
&=&\frac{\Gamma \left( \frac{n+m}{2}\right) }{2^{\frac{n+m}{2}}\Gamma \left( 
\frac{n}{2}\right) \Gamma \left( \frac{m}{2}\right) }\times u^{\frac{n}{2}%
-1}\times \left( u+1\right) ^{-\left( \frac{n+m}{2}\right)}.
\end{eqnarray*}

\bigskip 
\noindent Now using the general rule $f_{aZ}(t)=|a|^{-1}f_X(t/a)$ gives the \textit{pdf}

\begin{equation*}
f_{X}(u)=\frac{n^{n/2}m^{m/2}\Gamma((n+m)/2)}{\Gamma(n/2)\Gamma(m/2)} \frac{u^{n/2-1}}{(m+nux)^{(n+m)/2}}1_{(u\geq 0)},
\end{equation*}

\bigskip 
\noindent which is the \textit{pdf} of a Fisher $F_{n,m}$ random variable.
 

%% file: proba_02_04_ang.tex
\chapter{An Introduction to Gauss Random Measures} \label{proba_02_gauss}

 This chapter focuses of Gaussian probability measures on $\mathbb{R}$ first and next on $\mathbb{R}^d$, $d\geq 2$ exclusively. This is explained by the role of such probability laws in the history of Probability Theory and its presence in a great variety of sub-fields of Mathematics and in a considerable number of Science domains. Knowing that law and its fundamental properties is mandatory.

\section{Gauss Probability Laws on $\mathbb{R}$}

\noindent \textbf{(A) Standard Gauss Gauss Probability Law}.\\

\noindent We already encounter the function

\begin{equation} \label{proba_02_gauss_prop01}
f_{0,1}(x)=\frac{1}{\sqrt{2\pi}} \exp(-x^2/2), \ x \in \mathbb{R}.
\end{equation}

\bigskip \noindent and we proved that $\int_{\mathbb{R}} f_{0,1}(x) \ dx=1$. We remark that $f$ is locally bounded and locally Riemann integrable (LLBRI). So, we may equivalently consider the Riemann integral of $f_{0,1}$ or its Lebesgue integral. Without express notification, we will use Riemann integrals as long as we stay in the case where these Riemann integrals are Lebesgue's one.\\

\begin{definition} A random variable $X\ : \ (\Omega, \mathcal{A}, \mathbb{P}) \rightarrow \mathbb{R}$ is said to follow a \textsl{standard normal} or \textit{standard Gaussian} probability law, or in other words : $X$ is a \textsl{standard normal} or \textit{standard Gaussian} random variable if and only if $f_{0,1}$ is the \textsl{pdf} of $X$, that is the Radon-Nikodym of $\mathbb{P}_X$ with respect to the Lebesgue measure. Its strict support is the whole real line $\mathbb{R}$.\\
\end{definition}

\bigskip \noindent The historical derivation of such \textit{pdf} in the earlier Wworks of \textit{de} Moivre, Laplace and Gauss (1732 - 1801) is stated in \cite{loeve} and in \cite{ips-mestuto-ang} of this series.\\

\bigskip \noindent The main properties of a standard normal random variable are the following.

\begin{theorem}  \label{proba_02_gauss_th01} If $X$ is a standard normal random variable, then :\\

\noindent (1) $\mathbb{E}(X)=0$ and $\mathbb{V}ar(X)=1$.\\

\noindent (2) $X$ has finite moments of all orders and :\\
$$
\mathbb{E}X^{2k+1}=0 \text{  and  } \mathbb{E}X^{2k}=\frac{(2k)!}{2^k k!}, \ \ k\geq 1, 
$$

\bigskip \noindent and, in particular, its kurtosis parameter $K_X$ satisfies

$$
K_X=\frac{\mathbb{E}(X^4)}{\mathbb{E}(X^{2})^2}=3.
$$

\bigskip \noindent (3) Its \textit{mgf} is

$$
\varphi(u)=\exp(u^2/2), \ u \in \mathbb{R}
$$

\bigskip \noindent and its characteristic function is

$$
\Phi(u)=\exp(-u^2/2), \ u \in \mathbb{R}.
$$

\bigskip \noindent (5) Its \textit{cdf} 

$$
G(x)=\frac{1}{\sqrt{2\pi}} \int_{-\infty}^{x} \exp(-u^2/2) \ du, \ u \in \mathbb{R}
$$

\bigskip \noindent admits the approximation, for $x>1$,

\begin{equation*}
C\left\{ \frac{1}{x}-\frac{1}{x^{2}}\right\} e^{-x^{2}/2}\leq 1-G(x)\leq 
\frac{Ce^{-x^{2}/2}}{x}, 
\end{equation*}

\bigskip \noindent where $C=1/\sqrt{2\pi}$.\\

\noindent (6) The quantile function $G^{-1}(1-s)$ is expanded as $s\downarrow 0$, according to 
\begin{eqnarray*}
&&\phi ^{-1}(1-s)\\
&=&\left\{ (2\log (1/s))^{1/2}-\frac{\log 4\pi +\log \log (1/s)%
}{2(2\log (1/s))^{1/2}}+O((\log \log (1/s)^{2}(\log 1/s)^{-1/2}))\right\} .
\end{eqnarray*}

\bigskip \noindent and the derivative of $G^{-1}(1-s)$ is, as $s\downarrow 0,$%

\begin{eqnarray*}
\biggr(G^{-1}(1-s)\biggr)^{\prime} &=&(2\log (1/s))^{1/2}-\frac{\log 4\pi +\log \log (1/s)}{%
2(2\log (1/s))^{1/2}} \\
&&+O((\log \log (1/s)^{2}(\log 1/s)^{-1/2})).
\end{eqnarray*}

\bigskip \noindent (7) The following property holds. For each $x \in \mathbb{R}$,

$$
\lim_{n \rightarrow +\infty} G\biggr((2\log n)^{1/2} x + (2\log n)^{1/2} -\frac{\log 4\pi +\log \log n}{2(2\log n)^{1/2}} \biggr)^n=\exp(-e^{-x}).
$$
\end{theorem}

\bigskip \noindent For right now, we are only concerned with the three first Points. The other points are related to the tail $1-G$ of a normal law. We will deal with this in the monograph devoted to extreme value theory.\\ 

\noindent \textbf{Proof of Theorem}.\\

\noindent Points (1) and (2). Let $k\geq 0$. By using Formula (ACIF) (See Section \ref{proba_02_rv_sec_05} in Chapter \ref{proba_02_rv}, page \pageref{acif}), we have

$$
\mathbb{E}X^{2k+1}=\int_{\mathbb{R}} x^{2k+1} f_{0,1}(x) \ d\lambda(x).
$$

\bigskip \noindent The function $|x|^{2k+1} f_{0,1}(x)$ is locally bounded and locally Riemann integrable. So, we may use the recommendations in Point (b) in Section \ref{proba_02_rv_sec_05} in Chapter \ref{proba_02_rv} to get 
$$
\int_{\mathbb{R}} x^{2k+1} f_{0,1}(x) \ d\lambda(x)=\lim_{n\rightarrow+\infty} \int_{[-n,n]} x^{2k+1} f_{0,1}(x) \ dx.
$$

\bigskip \noindent Now, Riemann integration techniques ensure that, for each $n\geq 1$, $\int_{-n}^{n} x^{2k+1} f_{0,1}(x) \ dx=0$ since the continuous function $x^{2k+1} f_{0,1}(x)$ is odd on the symmetrical interval $[-n,n]$ with respect to zero. By putting together all the previous facts, we have

$$
\mathbb{E}X^{2k+1}=0.
$$

\bigskip \noindent For even order moments, we denote $I_{k}=\mathbb{E}X^{2k}$, $k\geq 0$. We have $I_0=1$. For $k\geq 1$, let us use Riemann integrals and integrations by parties. We have

\begin{eqnarray*}
I_{k}&=&\frac{1}{\sqrt{2\pi}} \int_{-\infty}^{+\infty} x^{2k} \exp(-x^2/2) \ dx= \frac{1}{\sqrt{2\pi}} \int_{-\infty}^{+\infty} x^{2k-1} \biggr(x\exp(-x^2/2)\biggr) \ dx\\
&=&\frac{1}{\sqrt{2\pi}} \int_{-\infty}^{+\infty} x^{2k-1} d\biggr(-\exp(-x^2/2)\biggr)\\
&=&\frac{1}{\sqrt{2\pi}}   \biggr[-\exp(-x^2/2)\biggr]_{-\infty}^{+\infty} + (2k-1) \frac{1}{\sqrt{2\pi}} \int_{-\infty}^{+\infty} x^{2k-2} \exp(-x^2/2) \ dx\\
&=& (2k-1) I_{k-1}.
\end{eqnarray*}

\bigskip \noindent We get by induction that

$$
I_k=(2k-1)I_{k-1}=(2k-1)(2k-3)I_{k-3}= \cdots =(2k-1)(2k-3)(2k-5) \cdots 3 I_{0}.
$$

\bigskip \noindent Hence

$$
I_k=(2k-1)(2k-3)(2k-5) \cdots 3.
$$

\bigskip \noindent By multiplying $I_k$ by the even numbers $(2k)(2k-2)\cdots 2 =2^k k!$ and dividing it as well, we get the results.\\

\bigskip \noindent (3) By still using Riemann integrals and using Formula (AC01) (See page \pageref{ac01}), we have for all $u\in \mathbb{R}$

\begin{eqnarray*}
\varphi(u)&=&\mathbb{E}(e^{tX})=\frac{1}{\sqrt{2\pi}} \int_{-\infty}^{+\infty} e^{ux} \exp(-x^2/2) \ dx\\
&=&\frac{1}{\sqrt{2\pi}} \int_{-\infty}^{+\infty} \exp\biggr( \frac{1}{2} \biggr(x^2-2tu\biggr) \biggr) \ dx
\end{eqnarray*}

\bigskip \noindent By using $x^2-2tu=(u-t)^2-u^2$, we get

$$
\varphi(u)=\exp(u^2/2) \biggr( \int_{-\infty}^{+\infty} \exp(-(x-u)^2/2) \ dx\biggr)
$$

\bigskip \noindent By using the change of variable $y=u-1$, we get that integral between the parentheses is one, and the proof is finished. $\square$\\

\bigskip \noindent \textbf{(B) Real Gauss Probability Laws}.\\

\noindent Now given a standard random variable $Z$, $m$ a real number and $\sigma>0$, the random variable

$$
X=\sigma Z + m,
$$

\bigskip \noindent has the \textit{cdf}, for $x \in \mathbb{R}$,

$$
F_X(x)=\mathbb{P}(\sigma Z + m \leq x)=\mathbb{P}\left(Z \leq \frac{x-m}{\sigma}\right)=F_Z\left(\frac{x-m}{\sigma}\right)
$$

\bigskip \noindent which, by differentiating the extreme members, leads to

$$
f_{m,\sigma}(x)=\frac{dF_X(x)}{dx}=\frac{1}{\sigma} \frac{dF_Z((x-m)/\sigma)}{dx},
$$

\bigskip \noindent for $x \in \mathbb{R}$. Since the functions $f_{m,\sigma}$ and $f_{0,1}$ are bounded and continuous, we may apply the recommendations of Point (b) 
[Section \ref{proba_02_rv_sec_05}, Chapter \ref{proba_02_rv}] to conclude that

$$
f_{m,\sigma}(x)=\frac{1}{\sigma \sqrt{2\pi}} \exp \biggr( - \frac{(x-m)^2}{2}\biggr), \ x \in \mathbb{R}. \ \ (RG)
$$

\bigskip \noindent is the absolute \textsl{pdf} of $X$. By using the properties of expectations and variances and properties characteristic functions, we have :\\

$$
\mathbb{E}X=m \text{  and  } \mathbb{V}ar(X)=\sigma^2, 
$$

$$
\varphi_X(u)=\exp(mu + \sigma^2 u^2/2), \ u \in \mathbb{R},
$$

\bigskip \noindent and

$$
\Phi_X(u)=\exp(imi -\sigma^2 u^2/2), \ u \in \mathbb{R}.
$$

\bigskip \noindent Before we conclude, we see that if $\sigma=0$, $X=m$ and its $mgf$ is $\exp(mu)$, which is of the form $\exp(mu + \sigma^2 u^2/2)$ for $\sigma=0$. We may conclude as follows.\\

\noindent \textbf{Definition - Proposition} (DEF01).\\

\noindent A real random variable is said to follow a Gaussian or normal probability law, denoted $X \sim \mathcal{N}(m,\sigma^2)$, if and only if its \textit{mgf} is given by

$$
\varphi_X(u)=\exp(mu + \sigma^2 u^2/2), \ u \in \mathbb{R}, \ (RGM)
$$

\bigskip \noindent or, if and only of its, characteristic function given by

$$
\Phi_X(u)=\exp(imu -\sigma^2 u^2/2), \ u \in \mathbb{R}. \ (RGC)
$$

\bigskip \noindent If $X$ is not degenerate, that is $\sigma>0$, its absolutely continuous \textsl{pdf} is

$$
f_X=\frac{1}{\sigma \sqrt{2\pi}} \exp \biggr( - \frac{(x-m)^2}{2} \biggr), \ x \in \mathbb{R}. \ (RGD)
$$

\bigskip \noindent Its first parameters are

$$
m=\mathbb{E}X \text{  and  } \sigma^2=\mathbb{E}X^{2}.  
$$

\bigskip \noindent $\Diamond$.

\bigskip \noindent \textbf{(D) Some immediate properties}.\\

\noindent \textbf{(D1) Finite linear combination of independent real Gaussian randoms}.\\

\noindent Any linear combination of a finite number $d\geq 2$ of independent random variables $X_1, \cdots, X_d$ with coefficient $\delta_1, \cdots, \delta_d$ follows a normal law. Precisely,
if the $X_i$'s are independent and $X_i \sim \mathcal{N}(m_i, \sigma_{i}^2)$, $1\leq i \leq d$, if we denote $m^t=(m_1,...,m_d)$, $\delta^t=(\delta_1, \cdots, \delta_d)$ and $\Sigma=diag(\sigma_{1}, \cdots, \sigma_{d})$, we have

$$
\sum_{1\leq j \leq d} \delta_j X_j \sim \mathcal{N}(m^t\delta, \delta^t \Sigma \delta).
$$

\bigskip \noindent To see this, put

$$
Y=\sum_{1\leq j \leq d} \delta_j X_j
$$

\noindent By the factorization property formula, we have for any $u\in \mathbb{R}^d$,

\begin{eqnarray*}
\Phi_Y(u)&=&\mathbb{E} \exp \biggr( \sum_{1\leq j \leq d} u \delta_j X_j\biggr)\\
&=&\mathbb{E} \exp \biggr( \prod_{1\leq j \leq d} u \delta_j X_j\biggr)\\
&=& \prod_{1\leq j \leq d} \mathbb{E}  \exp(u \delta_j X_j)\\
&=& \prod_{1\leq j \leq d} \exp(im_j \delta_j u -\delta_j^2 \sigma_j^2 u^2/2)\\
&=& \exp \biggr( i \langle m, \delta \rangle u - \biggr(\delta^t  \Sigma \delta\biggr) \frac{u^2}{2}\biggr).
\end{eqnarray*}

\bigskip \noindent From there, we may conclude that $Y$ follows a real random vector with the given parameters.\\

\noindent \textbf{(D2) Towards Gaussian Random Vectors}.\\

\noindent Let us remain in the frame of the previous point (D1). Let $X$ be the vector defined by $X^t=(X_1,\cdots, X_d)$ with independent real Gaussian Random variables with the given parameters. We have for any $u\in \mathbb{R}^d$ with $u^t=(u_1,...,u_d)$,

\begin{eqnarray*}
\Phi_X(u)&=&\mathbb{E} \exp i \langle u, Z\rangle\\
&=&\mathbb{E} \exp \biggr( \sum_{1\leq j \leq d} u_j X_j\biggr)\\
&=& \prod_{1\leq j \leq d} \mathbb{E}  \exp(u u_j X_j)\\
&=& \prod_{1\leq j \leq d} \exp(im_j u_j - \sigma_j^2 u^2/2)\\
&=& \exp \bigg( im^tu - \frac{u^t  \Sigma u}{2} \biggr).
\end{eqnarray*}

\noindent By also using the same techniques for the \textit{mgf}, we get that for any $u\in \mathbb{R}^d$ with $u^t=(u_1,...,u_d)$,

$$
\varphi_X(u)=\exp \bigg(  \langle m, u \rangle  + \frac{u^t  \Sigma u}{2}\biggr).
$$

\bigskip \noindent A random vector $X$ whose components are independent and satisfy  $X_i \sim \mathcal{N}(m_j, \sigma_{j}^2)$, $1\leq j \leq d$ has the \textit{mgf}

$$
\varphi_X(u)=\exp \bigg(  \langle m, u \rangle  + \frac{u^t  \Sigma u}{2}\biggr). \ (RV01)
$$

\noindent for any for any $u\in \mathbb{R}^d$, where $m^t=(m_1,...,m_d)$ and  $\Sigma=diag(\sigma_{1}, \cdots, \sigma_{d})$. Besides, we have \\

$$
\mathbb{E}X=m  \text{ and } \mathbb{V}ar(X)=\Sigma. 
$$

\bigskip \noindent This offers us a good transition to the introduction of Gaussian random vectors.\\

\bigskip
\section{Gauss Probability Law on $\mathbb{R}^d$, Random Vectors}

\bigskip \noindent \textbf{(A) Introduction and immediate properties}.\\

\noindent In general, the study of random vectors relies so much on quadratic forms and orthogonal matrices topic. We advice the reader to read at least the definitions, theorems and propositions on the aforementioned topic in Section \ref{proba_02_appendix_ortho} in the Appendix Chapter \ref{proba_02_appendix}. Each time a property on orthogonal matrices is quoted, it is supposed to be found in the appendix in the aforementioned section.\\

\noindent The above formula (RV01) gave us a lead to the notion of Gaussian random vectors. we have :

\begin{definition} A random variable $X\ : \ (\Omega, \mathcal{A}, \mathbb{P}) \rightarrow \mathbb{R}^d$, $d\geq 1$, is said to follow a $d$-multivariate Gaussian  probability law, or in other words : $X$ is a $d$-Gaussian random vector if and only its \textsl{mgf} is defined by

$$
\varphi_X(u)=\exp \bigg(  \langle m, u \rangle  + \frac{u^t  \Sigma u}{2}\biggr), \ u\in \mathbb{R}^2,  \ (RV02)
$$

\bigskip \noindent  where $m$ is a $d$-vectors or real numbers and $\Sigma$ is a symmetrical and semi-positive $d$-matrix or real numbers, and we write $X \sim \mathcal{N}_d(m,\Sigma)$.
\end{definition}

\bigskip \noindent By comparing with Formula (RV01), we immediately have :

\begin{proposition} \label{proba_02_gauss_prop02} A random vectors with real-valued independent Gaussian components is a Gaussian vector.
\end{proposition}

\bigskip \noindent We also have the following properties.

\begin{proposition} \label{proba_02_gauss_prop03} $X$ admits the \textsl{mgf} in Formula (REV02), then we have
$$
\mathbb{E}(X)=m  \text{ and } \mathbb{V}ar(X)=\Sigma. 
$$
\end{proposition}

\bigskip \noindent \textbf{Proof}. We are going to construct a random vector $Y$ which has the \textsl{mgf}

$$
\exp \bigg(  \langle m, u \rangle  + \frac{u^t  \Sigma u}{2}\biggr), \ u\in \mathbb{R}^d.  \ (RV03)
$$

\noindent and next use the characterization of the probability law by the \textsl{mgf}. Since $\Sigma$ is symmetrical and semi-positive, we may find an orthogonal $d$-matrix $T$ such that

$$
T \Sigma T^t =diag(\delta_1, \cdots, \delta_d),
$$

\bigskip \noindent where $\delta_1, \cdots, \delta_d$ are non-negative real numbers. By the Kolmogorov Theorem as applied in Point (c5) in Section \ref{proba_02_rv_sec_04} in Chapter \ref{proba_02_rv}, we may find a probability space $(\Omega, \mathcal{A}, \mathbb{P})$ holding a $d$-random vector $Z$ whose components are centered  independent real-valued Gaussian random vectors with respective variances $\delta_j$, $1 \leq j \leq d$. Its follows that $Z$ is Gaussian
and hence, by Formula (REV01), we have for $D=diag(\delta_1, \cdots, \delta_d)$, for $u\in \mathbb{R}^d$,

$$
\varphi_Z(u)=\exp \bigg( \frac{u^t  D u}{2}\biggr).
$$ 

\bigskip \noindent Now let us set $Y=m+T^t Z$. We have

\begin{eqnarray*}
\varphi_Y(u)&=&\mathbb{E} \exp \bigg( (m+T^t Z)^t u \biggr)= \exp (m^t u)  \mathbb{E} \exp \bigg( Z^t (T u) \biggr)\\
&=& \exp (m^t u)  \mathbb{E} \exp \bigg(\langle Z,  T u \rangle \biggr)\\
&=& \exp (\langle m, u \rangle ) \exp \bigg( \frac{u^t T^t DT u}{2}\biggr)
\end{eqnarray*}

\bigskip \noindent But, by the properties of orthogonal matrices, we have $T\Sigma T^t=D$, which implies that $T\Sigma T^tDT=\Sigma$. Thus we have

\begin{eqnarray*}
\varphi_Y(u)=\exp \bigg(  \langle m, u \rangle  + \frac{u^t  \Sigma u}{2}\biggr), \ u\in \mathbb{R}^d.
\end{eqnarray*}

\bigskip \noindent This a direct proof, based on the Kolmogorov construction, that the function in Formula (REV03) is characteristic. An other method would rely on the Bochner Theorem we do not mention here. At the end, we have that $Y \sim \mathcal{N}_d(m,\Sigma)$. By using the properties and expectation vectors and variance-covariance properties seen in Chapter \ref{proba_02_rv}, we have

$$
\mathbb{E}Y=\mathbb{E}(m+T^t Z)=\mathbb{E}(m)+T^t \mathbb{E}(Z)=m 
$$ 

\bigskip \noindent and, since the constant vector $m$ is independent from $T^tZ$,

$$
\mathbb{V}ar(Y)=\mathbb{V}ar(m+T^Z)=\mathbb{V}ar(T^Z)=T^t DT=\Sigma. 
$$ 

\bigskip \noindent We conclude as follows : for any random variable characterized by its \textsl{mdf} given in Formula (REV02), its expectation vector and its variance-covariance matrix are given as above. $\blacksquare$\\ 

\bigskip \noindent \textbf{Important Remark}. In the notation $X \sim \mathcal{N}_d(m,\Sigma)$, $m$ and $\Sigma$ are the respective expectation vector and the variance-covariance matrix of $X$.\\

\noindent Let us now study other  important properties of Gaussian vectors.\\

\bigskip \noindent \textbf{(B) - Linear transforms of Gaussian Vectors}.\\

\begin{proposition} The following assertions hold.\\

\noindent (a) Any finite-dimension linear transform of a Gaussian random vector is a Gaussian random vector.\\

\noindent (b) Any linear combination of the components of a Gaussian random vector is a real Gaussian random variable.\\

\noindent (c) If a random vector $X\ : \ (\Omega, \mathcal{A}, \mathbb{P}) \rightarrow \mathbb{R}^d$, $d\geq 1$, follows a $\mathcal{N}_d(m,\Sigma)$ probability law and if $A$ is a $(k \times d)$-matrix and $B$ a $k$-vector, $k\geq 1$, then $AX+B$ follows a $\mathcal{N}_k(Am+B,A\Sigma A^t)$ probability law.\\
\end{proposition}

\bigskip \noindent \textbf{Proof}. It is enough to prove Point (c). Suppose that the assumption of that point hold. Thus $Y=AX$ is $k$-random vector. By Point (a) of Theorem \ref{proba02_rv_prop_cf} in 
Section \ref{proba_02_rv_sec_06} in Chapter \ref{proba_02_rv}, we have

$$
\Phi_{AX+B}(v)=\exp(B^t v) \Phi_{X}(A^t v)
$$

\bigskip \noindent and combining this with Formula (RV02) gives, for any $v\in \mathbb{R}^k$,

\begin{eqnarray*}
\varphi_{AX+B}(v)&=&\exp( B^t v) \Phi_{X}(A^t v)\\
&=&\exp( B^t v) \exp \bigg(   m^t (A^t v)  + \frac{v^t  A\Sigma A^t v}{2}\biggr)\\
&=& \exp \bigg( (B+Am)^t v + \frac{ (A^t v)^t  A\Sigma (A^t v)}{2}\biggr)\\
&=& \exp \bigg( \langle B+Am, v \rangle + \frac{ v^t \left(A\Sigma A \right) v)}{2}\biggr). \ \square.
\end{eqnarray*}

\bigskip \noindent This proves (c) which is a more precise form of (a). Point (c) is only an application of Point (c) to a $(d\times 1)$-matrix $A$.\\

\bigskip \noindent Point (c) provides a new definition of Gaussian vectors \textbf{given we already have the definition of a real-valued Gaussian random variable}. We have : \\

\noindent \textbf{Definition - Proposition}. \\

\noindent (a) (DEF01) Any $d$-random vector, $d\geq 2$, is Gaussian if an only if any linear combination of its components is a real-valued Gaussian random variable.\\

\noindent (b) (DEF02) \textit{Given we already have the definition of a real-valued Gaussian random variable}, a $d$-random vector, $d\geq 2$,  is Gaussian if any linear combination of its components is a real-valued Gaussian random variable. $\Diamond$\\

\bigskip \noindent \textit{Proof of Point (a)}. Let $X\ : \ (\Omega, \mathcal{A}, \mathbb{P}) \rightarrow \mathbb{R}^d$, $d\geq 2$, be a random vector such that any linear combination of its components is a real-valued Gaussian random variable.\\

\noindent First, each component is Gaussian and hence is square integrable and next, by Cauchy-Schwartz inequality, any product of two components is integrable. Hence the expectation vector $m$ and the variance-covariance matrix $\Sigma$ of $X$ have finite elements. Next, the characteristic function of $X$ is satisfies, for any $u\in \mathbb{R}^d$.\\

$$
\Phi_X(u)=\mathbb{E} \exp \bigg( i u^t X \biggr)=\Phi_{u^t X}(1), \ (RV04)
$$

\bigskip \noindent where $\Phi_{u^t X}$ is the characteristic of $u^t X=u_1 X_1+\cdots +  u_d X_d$ which is supposed to be a real-valued normal random variable with parameters

$$
\mathbb{E}(u^t X)=\sum_{1\leq j \leq d} u_j X_j = \langle u, m \rangle
$$

\noindent and

\begin{eqnarray*}
\mathbb{V}ar(u^t X)&=&\mathbb{V}ar(\sum_{1\leq j \leq d} u_j X_j)\\
&=&\sum_{1\leq j \leq d} \sum_{1\leq j \leq d} \mathbb{C}ov(X_i, X_j) u_i u_j\\
&=&u^t \Sigma u.
\end{eqnarray*}

\bigskip \noindent Now using the characteristic function of a $\mathcal{N}(\langle u, m \rangle, u^t \Sigma u)$ allows to conclude. $\square$.

\bigskip\bigskip \noindent \textbf{Some consequences}.\\

\noindent (a) A sub-vector of a Gaussian Vector is a Gaussian vector since it is a projection, then a finite-dimensional linear transform, of the vector.\\

\noindent (b) As particular cases of Point (a), components of a Gaussian vector are Gaussian.\\

\noindent (c) A vector whose components are \textbf{independent} and Gaussian is Gaussian.\\

\noindent (c) \textbf{But, in general}, a vector whose components are Gaussian is not necessarily Gaussian. Here is a general, using \cite{sklar1959}'s Theorem, to construct counter-examples. As stated in Section \ref{proba_02_rv_07a} of Chapter \ref{proba_02_rv}, for any random vector $X$ of dimension $d\geq 1$, the \textit{cdf} $F_X$ of $X$ satisfies

\begin{equation*}
\forall x \in \mathbb{R}^d, \ F_X(x)=C(F_{X,1}(x),...,F_{X,d}(x)), 
\end{equation*}
 
\bigskip \noindent where $C$ is a copula and $F_{X,j}$ stand for the individual marginal \textit{cdf}'s and the copula is unique if the marginal \textit{cdf}'s are continuous. By choosing the $F_{X,j}$ as \textit{cdf}'s of Gaussian random variables $X_{j}$, the vector $X=(X_{1},...,X_{d})^t$ has Gaussian components. But not any copula $C$ makes $F_X$ a \textit{cdf} of Gaussian vector.\\

\noindent For example, for $d=2$, by taking the least copula $C(u,v)=max(u + v -1,0)$, $(u,v) \in [0,1]^2$, $\Phi$ the \textit{cdf} of a $\mathcal{N}(0,1)$ random variable, a random vector 
$(X,Y)^t$ associated with the \textit{cfd}

$$
\forall x=(x,y)^t \in \mathbb{R}^2, \ F(x)=\max(\Phi(x)+\Phi(y) -1,0),
$$

\noindent is not Gaussian but has Gaussian components.\\

\bigskip \noindent \textbf{(C) - Uncorrelated and Gaussian Component}.\\

\noindent Let us begin to resume the result of this part by saying this : For a Gaussian vector, uncorrelation and independence of its sub-vectors are the same. Precisely we have :\\

\begin{proposition} Let  $Y\ : \ (\Omega, \mathcal{A}, \mathbb{P}) \rightarrow \mathbb{R}^r$ and $Z\ : \ (\Omega, \mathcal{A}, \mathbb{P}) \rightarrow \mathbb{R}^s$, $r\geq 1$, $s\geq b$, be a two random vectors such that $X^t=(Y^t, Z^t)$ is a $d$-Gaussian vector, $d=r+s$. Suppose that $Y$ and $Z$ are uncorrelated, that is, their covariance matrices are null matrices

$$
\mathbb{C}ov(Y,Z)\mathbb{C}ov(Z,Y)^t=\biggr(\mathbb{C}ov(Y_i, Z_j\biggr)_{1\leq i \leq r, \ 1\leq j \leq s}=0, \ 
$$

\bigskip \noindent that is also

$$
\forall (u,j)\in \{1,...,r\} \times \{1,...,s\}, \ \ \mathbb{C}ov(Y_i, Z_j)=0.
$$

\bigskip \noindent Then $Z$ and $Y$ and independent.
\end{proposition}

\bigskip \noindent \textbf{Proof}.\\

\noindent Since $X$ is Gaussian, its sub-vectors $Z$ and $Z$ are Gaussian and have \textit{mgf} functions 

$$
\mathbb{R}^r \ni v \mapsto \varphi_Y(v)=\exp\left(m_Y^t v + v^t \Sigma_Y v\right) \ \ (RV05a)
$$

\bigskip \noindent and

$$
\mathbb{R}^s \ni w \mapsto \varphi_Z(w)=\exp\left(m_Z^t w + w^t \Sigma_Y w\right), \ (RV05b)
$$

\bigskip \noindent where $m_Y$ and $\Sigma_Y$ (resp. $m_Z$ and $\Sigma_Z$) are the expectation vector and the variance-covariance matrix of $Y$ (resp. $Z$). The components of $X$ are
$X_i=Y_i$ for $1\leq i \leq r$ and $X_i=Z_i$ for $r+1\leq i \ d$. Suppose that $Y$ and $X$ are uncorrelated. Denote also by $m_X$ and $\Sigma_X$ the expectation vector and the variance-covariance matrix of $X$.\\

\noindent Thus for any  $v \in \mathbb{R}^r$, $w \in \mathbb{R}^r$, we have by denoting
$u^t=(v^t,w^t)$, $u \in \mathbb{R}^d$,

\begin{eqnarray*}
u^t \Sigma_X u&=&\sum_{1\leq i \leq d, \ 1\leq j \leq d} \mathbb{C}ov(X_i,X_j)\\
&=&\sum_{1\leq i \leq r, \ 1\leq j \leq r} \mathbb{C}ov(X_i,X_j)+\sum_{r+1\leq i \leq r, \ r+1\leq j \leq d} (L2)\\
&+&\sum_{1\leq i \leq r, \ r+1\leq j \leq d} \mathbb{C}ov(X_i,X_j)+\sum_{r+1\leq i \leq d, \ 1 \leq j \leq s} \ (L3)\\
\end{eqnarray*}

\bigskip \noindent The covariances of Line (L3) are covariance between a component of $Y$ and another of $Z$ and by hypothesis, the summation in that line is zero. In the first term of Line (L2), the covariances are those between components of $Y$ and the second term contains those of components of $Z$. We get   

\begin{eqnarray*}
u^t \Sigma_X u=v^t \Sigma_Y v + w^t \Sigma_Y w
\end{eqnarray*}

\noindent with the same notation, we have $u^t m_X=v^t m_Y + v^t m_Z$ and, by taking Formula (RV05) into account, we arrive at

\begin{eqnarray*}
\varphi_X(u)&=&\varphi_{(Y,Z)}(u,w)=\exp \biggr(u^t m_X + \frac{u^t \Sigma_X u}{2}\biggr)\\
&=&\exp \biggr(v^t m_Y + v^t m_Z + \frac{v^t \Sigma_Y v + w^t \Sigma_Y w}{2}\biggr)\\
&=&\varphi_Y(v) \varphi_Z(w).
\end{eqnarray*}

\bigskip \noindent We finally have for any $v \in \mathbb{R}^r$, $w \in \mathbb{R}^r$,

$$
\varphi_{(Y,Z)}(u,w)=\varphi_Y(v) \varphi_Z(w).
$$

\bigskip \noindent By Theorem \ref{proba02_rv_prop_cf031} in Section \ref{proba_02_rv_sec_06} in Chapter \ref{proba_02_rv}, we conclude that $Y$ and $Z$ are independent.\\

\bigskip \noindent \textbf{WARNING} Gaussian Random vectors do not have the exclusivity of such a property. To make it simple, this property holds for a random ordered pair $(X,Y)$ if for example, for any $(u,v)\in \mathbb{R}^2$,

$$
\biggr|\Phi_{(X,Y)}(u,v)-\Phi_{X}(u)\Phi_{Y}(v)\biggr|\leq h_{X,Y}(n,u),
$$

\bigskip \noindent where $h_{X,Y}$ is a function satisfying $h_{X,Y}(0,0)=0$.\\

\noindent \textbf{Example : Associated random variables}. A finite family of $d$ real random variables $X_j\ : \ (\Omega, \mathcal{A}, \mathbb{P}) \rightarrow \mathbb{R}$, $1\leq j \leq d$ is said to be associated if and only for any pair$(f,g)$ of bounded real-valued and measurable functions functions defined on $\mathbb{R}^d$ both coordinate-wisely non-decreasing, we have

$$
\mathbb{C}ov\biggr(f(X_1,...,X_d)g(X_1,...,X_d)\biggr)\geq 0.
$$

\noindent Let us denote $X^t=(X_1,...,X_d)$ and let $\Sigma_X$ be the variance-covariance matrix of $X$. If that sequence is associated, \cite{newman} Theorem states that for any $u\in \mathbb{R}^d$,

$$
\biggr|\Phi_{(X_1,...,X_d)}(u)-\prod_{1\leq j\leq d}\Phi_{X_j}(u_j)\biggr| \leq \frac{1}{2} \sum_{1\leq i \leq d, \ 1\leq j \leq d} |u_iu_j|\mathbb{C}ov(X_i,X_j).
$$

\bigskip \noindent It is useful to know that the covariances $\mathbb{C}ov(X_i,X_j)$ are non-negative for associated variables. Thus, associated and uncorrelated variables are independent. 

\bigskip \noindent \textbf{(D) - Density probability function of Gaussian Vectors with a positive variance-covariance matrix}.\\

\noindent Probability laws of Gaussian vectors with non-singular variance-covariance matrix may be characterized by their absolute \textit{pdf}. We have the following :\\

\noindent \textbf{Proposition - Definition} (DEF03).\\

\noindent (a) Let $X\ : \ (\Omega, \mathcal{A}, \mathbb{P}) \rightarrow \mathbb{R}^d$, $d\geq 1$, be Gaussian random vector of expectation vector $m$ and variance-covariance matrix $\Sigma$. If $\Sigma$ is invertible, then $X$ has the \textit{pdf}

$$
\frac{det(\Sigma)^{-1/2}}{(2\pi)^{d/2}} \exp \biggr( -\frac{(x-m)^t \Sigma^{-1}(x-m)}{2} \biggr), \ x\in \mathbb{R}^d. \ x \in \mathbb{R}^d. \ (RVD)
$$

\bigskip \noindent (b) (DEF03) A random vector $X\ : \ (\Omega, \mathcal{A}, \mathbb{P}) \rightarrow \mathbb{R}^d$, $d\geq 1$, whose variance-covariance is invertible is a Gaussian vector if and only if it admits the absolute density probability \textsl{pdf} (RVD) above.\\

\bigskip \noindent \textbf{Proof}. We use the same techniques as in the proof of Proposition \ref{proba_02_gauss_prop03} and based on the Kolmogorov construction of a probability space $(\Omega, \mathcal{A}, \mathbb{P})$ holding a $d$-random vector $Z$ whose components are centered  independent real-valued Gaussian random vectors having as variances the eigen-values $\delta_j$, $1 \leq j \leq d$ of 
$\Sigma$. Set $D=diag(\delta_1,...,\delta_d)$. All those eigen-value $\delta_j$, $1\leq j \leq d$, are positive and

$$
det(\Sigma)=\prod_{1\leq j \leq d} \delta_j.
$$

\bigskip \noindent We may use the \textsl{pdf}'s of each $Z_j$ and make profit of their independence to get the \textsl{pdf} of $Z$, which is for any $x\in \mathbb{R}^d$,

$$
f_Z(z)=\prod_{1\leq j \leq d} f_{Z_i}(z_i)=\prod_{1\leq j \leq d} \frac{1}{(2\pi \delta_j)} \exp\biggr(- \sum_{1\leq j \leq d} \frac{z_i^2}{2\delta_j}\biggr),
$$

\bigskip \noindent which yields

$$
f_Z(z)= \frac{det(\Sigma)^{-1/2}}{(2\pi)^{d/2}} \exp\biggr(-  \frac{z^t D^{-1} z}{2}\biggr), \ z\in \mathbb{R}^d.
$$

\bigskip \noindent Now, let $T$ be an orthogonal matrix such that $T \Sigma T^t=D$. Set $Y=T(Z+m)$, that is : $Z=T^t Y -m$,  is a diffeomorphism which preserves the whole domain $\mathbb{R}^d$ of $Z$ and the Jacobian coefficient $J(y)$ is the determinant of $T^t$ which is $\pm 1$. The change of variable formula (CVF) in Section \ref{proba_02_rv_sec_07} in Chapter \ref{proba_02_rv} leads to

$$
f_Y(y)=f_Y(T(z-m)=\frac{det(\Sigma)^{-1/2}}{(2\pi)^{d/2}} \exp\biggr(- \frac{(y-m) T^t D^{-1} T (y-m)}{2}\biggr),
$$

\bigskip \noindent where $\ y\in \mathbb{R}^d$. Since $T^t D^{-1} T=\Sigma^{-1}$, we conclude that

$$
f_Y(y)=\frac{det(\Sigma)^{-1/2}}{(2\pi)^{d/2}} \exp\biggr(-  \frac{(y-m) \Sigma^{-1} (y-m)}{2}\biggr), \ y\in \mathbb{R}^d.
$$

\bigskip \noindent By combining this with Proposition \ref{proba_02_gauss_prop03}, we conclude that the non-negative function given in formula (RVD) is an absolute \textsl{pdf} and is the \textsl{pdf} of any random vector with the \textsl{mgf} given in Formula (RV02) for a non-singular matrix $\Sigma$. $\blacksquare$\\

\noindent \textbf{Different definitions}. We provided three definitions (DEF01), (DEF02) and (DEF03) for Gaussian vectors. The first which is based of the characteristic function or the \textsl{mgf} is the most general. The second suppose we already have the definition a real Gaussian random \textbf{}variable. The last assumes that the variance-covariance is invertible.\\

\noindent \textbf{Remark}. Another way to proceed for the last  proof is to directly show that the function given Formula (RVD) is a \textsl{pdf} and to compute its \textsl{mgf} by using the orthogonal transform of $\Sigma$. By trying to do so, Formula (UID) in Section \ref{proba_02_appendix_ortho} in Chapter \ref{proba_02_appendix_ortho} may be useful.\\  

\noindent \textbf{(E) Quadratic forms of Gaussian Vectors}.\\

\noindent Let $X \sim \mathcal{N}_d(m,\Sigma)$, $d\geq 1$, be a $d$-dimensional Random Vector. We have the following sample result.\\

\begin{proposition} 
If $\Sigma$ is invertible, then the quadratic form $(X-m)^t \Sigma^{-1} (X-m)$ follows a Chi-square probability law of $d$ degrees of freedom, that is

$$
(X-m)^t \Sigma^{-1} (X-m) \sim \chi_{d}^{2}.
$$
\end{proposition}

\bigskip \noindent \textbf{Proof}. Suppose that $X \sim \mathcal{N}_d(m,\Sigma)$ and $\Sigma$ is invertible. Let $T$ be an orthogonal matrix such that 
$$
T \Sigma T^t=D=diag(\delta_1,...,\delta_d)
$$

\bigskip \noindent which entails

$$
\Sigma^{-1}= T^t D^{-1} T.
$$

\bigskip \noindent Set $Y=T (X -m)$. Thus $Y$ is a Gaussian vector. Its variance-covariance matrix is $\Sigma_Y=T \Sigma T^t=D$. Hence the components $Y_1$, ..., $Y_d$ are Gaussian and not correlated. Hence they are independent. By Fact 2 in Point (11) on the Chi-square probability law in Section \ref{proba_02_upl_sec_02} in Chapter \ref{proba_02_upl_sec_02}, we have

$$
Q=\sum_{1\leq j \leq d}\frac{Y_j^2}{\delta_j} \sim \chi_d^2.
$$
 
\bigskip \noindent Since $D^{-1}=diag(1/\delta_1,...,1/\delta_d)$, we have

$$
Q=Y^tD^{1-}Y=(X-m)^t T^t D^{1-}T (X-m)=(X-m)^t \Sigma^{-1}(X-m) \sim \chi_d^2. \ \square
$$

%% file: proba_02_05_ang.tex
\chapter{Introduction to Convergences of Random Variables} \label{proba_02_conv}

\section{Introduction} \label{proba_02_conv_sec_01}

\noindent The convergence of random variables, extended by the convergence of their probability laws, is a wide field with quite a few number of sub-fields. In Statistical terms, any kind of convergence theory of sequences of random variables is classified in the asymptotic methods area.\\

\noindent We are going to introduce some specific types of convergence.\\

\noindent Let $(X_n)_{n\geq 0}$ be a sequence of random elements with values in a Borel space $(E,\mathcal{B})$, where $\mathcal{B}$ is the $\sigma$-algebra generated by the class of open set $\mathcal{O}$, such that each $X_n$, $n\geq 0$, is defined on some probability space $(\Omega_n, \mathcal{A}_n,\mathbb{P}^{(n)})$.\\

\noindent Let also $X_{\infty} : (\Omega_{\infty}, \mathcal{A}_{\infty},\mathbb{P}^{({\infty})}) \rightarrow (E,\mathcal{O})$ be some random element.\\

\noindent \textbf{Notation}. We will simply write $X=X_{\infty}$ if no confusion is possible.\\

\noindent \textbf{Regularity Condition}. At least we suppose that the topological space $(E,\mathcal{O})$ is separated ensuring that limits are unique and for sequences of any random elements 
$X, X_n : (\Omega, \mathcal{A},\mathbb{P}) \rightarrow (E,\mathcal{O})$, $n\geq 0$, we have

$$
(X_n \rightarrow X) \in \mathcal{A}.
$$

\noindent Now let us present some the following definitions for convergence of random variables after the \\

\noindent \textbf{Warning : In this textbook, only the convergences (A), (B) , (F) and (G) will be addressed, and they will studied on $E=\mathbb{R}^d$}.\\

\bigskip \noindent \textbf{(A) Almost-sure Convergence}. \textit{Suppose that $X_{\infty}$ and all the elements of the sequence $(X_n)_{n\geq 0}$ are defined on the same probability space $(\Omega, \mathcal{A},\mathbb{P})$}.\\

\noindent The sequence $(X_n)_{n\geq 0}$ converges almost-surely to $X_{\infty}$ and we denote

$$
X_n \rightarrow X_{\infty}, \ a.s. \ as \ n\rightarrow +\infty,
$$

\bigskip 
\noindent if and only 

$$
\mathbb{P}(X_n \nrightarrow X_{\infty})=0. \ (ASC)
$$

\bigskip \bigskip 
\noindent \textbf{(B) Convergence in Probability}. \textit{Suppose that $X_{\infty}$ and all the elements of the sequence $(X_n)_{n\geq 0}$ are defined on the same probability space 
$(\Omega, \mathcal{A},\mathbb{P})$ and $E$ is a normed real linear space and its norm is denoted by $\|.\|$}.\\

\noindent The sequence $(X_n)_{n\geq 0}$ converges in probability to $X_{\infty}$ and we denote

$$
X_n \overset{\mathbb{P}}{\longrightarrow} X_{\infty},  \ as \ n\rightarrow +\infty,
$$

\bigskip 
\noindent if and only for any $\varepsilon>0$

$$
\lim_{n\rightarrow +\infty} \mathbb{P}(\|X_n-X_{\infty}\|>\varepsilon)=0. \ (CP)
$$

\bigskip \bigskip 
\noindent \textbf{(C) General Convergence in Probability}. \textit{Suppose that $X_{\infty}$ and all the elements of the sequence  $(X_n)_{n\geq 0}$ are defined on the same probability space $(\Omega, \mathcal{A},\mathbb{P})$}.\\

\noindent The sequence $(X_n)_{n\geq 0}$ generally converges in probability to $X_{\infty}$ and we denote

$$
 X_{n} \overset{\mathbb{P}^{(g)}}{\longrightarrow} X_{\infty},  \ as \ n\rightarrow +\infty,
$$

\bigskip 
\noindent if and only for any open set $G \in E$, 

$$
\lim_{n\rightarrow +\infty} \mathbb{P}(X \in G, \ X_n \notin G)=0. \ (GCP)
$$

\bigskip \bigskip 
\noindent \textbf{(D) Complete Convergence}. \textit{Suppose that $X_{\infty}$ and all the elements of the sequence  $(X_n)_{n\geq 0}$ are defined on the same probability space $(\Omega, \mathcal{A},\mathbb{P})$ and $E$ is a normed real linear space and its norm is denoted by $\|.\|$}.\\

\noindent The sequence $(X_n)_{n\geq 0}$ completely converges to $X_{\infty}$ and we denote

$$
X_{n} \overset{c.c}{\longrightarrow} X_{\infty},  \ as \ n\rightarrow +\infty,
$$

\bigskip 
\noindent if and only for any $\varepsilon>0$,

$$
\sum_{n}^{+\infty} \mathbb{P}(\|X_n-X_{\infty}\|>\varepsilon)<+\infty. \ (CC)
$$

\bigskip \bigskip 
\noindent \textbf{(E) Convergence in $p$-th moment, $p>0$}. \textit{Suppose that $X_{\infty}$ and all the elements of the sequence  $(X_n)_{n\geq 0}$  are defined on the same probability space $(\Omega, \mathcal{A},\mathbb{P})$ and $E$ is a normed real linear space and its norm is denoted by $\|.\|$}. Let $r>0$.\\

\noindent The sequence $(X_n)_{n\geq 0}$ converges to $X_{\infty}$ in the $r$-th moment and we denote

$$
X_{\infty} \overset{m^r}{\longrightarrow} X_{\infty},  \ as \ n\rightarrow +\infty,
$$

\bigskip 
\noindent if and only for 

$$
\lim_{n\rightarrow +\infty} \mathbb{E}\|X_n-X_{\infty}\|^r=0. \ \ (MR)
$$

\bigskip \bigskip 
\noindent \textbf{(F) Convergence in moment $L^{p}$, $p\geq 1$}. \textit{Suppose that $X_{\infty}$ and all the elements of the sequence  $(X_n)_{n\geq 0}$  are \textbf{real-valued} mappings defined on the same probability space $(\Omega, \mathcal{A},\mathbb{P})$ and belong all to $L^{p}(\Omega, \mathcal{A},\mathbb{P})$}.\\

\noindent The sequence $(X_n)_{n\geq 0}$ convergences to $X_{\infty}$ in $L^{p}$ and we denote

$$
X_{\infty} \overset{L^{p}}{\longrightarrow} X_{\infty},  \ as \ n\rightarrow +\infty,
$$

\bigskip \noindent if and only if 

$$
\lim_{n\rightarrow +\infty} \mathbb{E}\|X_n-X_{\infty}\|_{p}=0. \ \ (CLP)
$$

\bigskip \noindent \textbf{Important Remark}. It is of the greatest importance to notice that all the previous limits, the random variables $X_{\infty}$ and $X_{n}$, $n\geq 0$, are defined on the same
\textbf{probability space}. This will not be the case in the next definition. Each random element may be defined on its own probability space. We will come back to this remark after the definition.\\

\noindent \textbf{(G) Weak Convergence in a metric space}. \textit{Suppose $E$ is a metric space $(E,d)$ endowed with the metric Borel $\sigma$-algebra}. Denote by $\mathcal{C}_b(E)$ the class of all real-valued, bounded and continuous functions defined on $E$. Define the probability laws :\\

$$
\mathbb{P}_{\infty}= \mathbb{P}^{(\infty)} X_{\infty}^{-1}, \ \ \mathbb{P}_{n}= \mathbb{P}^{(n)} X_n^{-1}, \ n\geq 0.  
$$

\bigskip \noindent The sequence $(X_n)_{n\geq 0}$ weakly convergences to $X_{\infty}$ and we denote

$$
X_{\infty} \rightsquigarrow X_{\infty},  \ as \ n\rightarrow +\infty,
$$

\bigskip \noindent if and only for $f\in \mathcal{C}_b(E)$

$$
\lim_{n\rightarrow +\infty} \int_E f \ d\mathbb{P}_{n}=\int_E f \ d\mathbb{P}_{\infty}.\  \ (WC)
$$

\bigskip \noindent \textbf{Remark}. We effectively see that only the probability laws of $X_n$, $n\geq 0$ and $X_{\infty}$ are concerned in Formula (WC), at the exclusion of the paths
$\{X_{\infty}(\omega), \ \omega\in \Omega_{\infty} \}$ and $\{X_n(\omega), \ \omega\in \Omega_n \}$, $n\geq 0$. In general, a type of convergence which ignores the domain of elements of the sequence whose limit is considered, is called \textit{weak} or \textit{vague}.\\

\noindent As announced earlier, we are going to study convergences type (A), (B) , (F) and (G) for sequences of random vectors in $\mathbb{R}^d$, $d\geq 1$. At this step, the three remarks are should be
made.\\

\noindent \textbf{(a) Convergence (A) and (B) are already treated in the Measure Theory and Integration book}. We will give easy extensions only.\\

\noindent \textbf{(b) Convergence (G) is treated in a separate monograph}. At this step of this course of probability theory, the weak convergence theory for random vectors may be entirely treated. This is what we did in \cite{ips-wcrv-ang}, as an element of the current series. The reason we expose that theory in an independent textbook us that we want it to be a first part of the exposition of Weak convergence embracing the most general spaces, including, stochastic processes.\\

\noindent The reader is free to read it as soon as he has completed the chapters 1 to 4 of this textbook. But, for coherence's sake, we will give the needed reminders to have a comprehensive comparison between the different kinds of convergence.\\

\noindent \textbf{(c) Space $L^p$}. Convergence in $L^p$ is simply a convergence in a normed space $L^p$. We already know for the Measure Theory and Integration book  that this space is a Banach one.\\

\noindent \textbf{(d) Convergences of real sequences}. When dealing with random vectors, a minimum prerequisite is to master the convergence theory for non-random sequences of real number. This is why we always include a related appendix in our monographs dealing with it. In this book, the reminder is exposed in Section \ref{proba_02_appendix_limits} in the Appendix chapter 
\ref{proba_02_appendix}.\\

\bigskip \noindent After the previous remarks, we see that this chapter is rather a review one with some additional points. In particular, the equi-continuity notion will be introduced for the comparison between the convergence in measure and the $L^p$-convergence.\\

\newpage
\noindent \textbf{Part A : Convergences of real-valued random variables}. \label{proba_02_convergence_A}.\\

\section{Almost-sure Convergence, Convergence in probability} \label{proba_02_conv_sec_02}

\noindent As recalled previously, such convergences have been studied in probability and Integration [Chapter 7 in \cite{ips-mestuto-ang}]. We are just going to report the results.

\bigskip \noindent \textbf{(a) Almost-everywhere convergence}.\\

\bigskip \noindent A sequence of random variables $X_{n})_{n\geq 1}$\ defined from ($\Omega ,\mathcal{A}, \mathbb{P})$\ to $\overline{\mathbb{R}}$ converges almost-surely to a random variable $X:(\Omega ,\mathcal{A},m)$\ $\longmapsto \overline{\mathbb{R}}$ and we
denote
\begin{equation*}
X_{n}\longrightarrow f,\text{ }a.s.,
\end{equation*}

\bigskip \noindent if and only if 
\begin{equation*}
\mathbb{P}(X_n \ \nrightarrow \ X)=0.
\end{equation*}

\bigskip \noindent If the elements of the sequences $X_n$ are finie \textsl{a.s.}, we have :\\

\noindent \textbf{Characterization}. A sequence of \textit{a.s.} finite random variables $(X_{n})_{n\geq 1}$\ defined from ($\Omega ,\mathcal{A}, \mathbb{P})$\ to $\overline{\mathbb{R}}$ converges almost-surely to a random variables $X:(\Omega ,\mathcal{A},m)$\ $\longmapsto \overline{\mathbb{R}}$ if and only if

$$
\mathbb{P}\left( \bigcap_{k\geq 1}\bigcap_{N\geq 1}\bigcup_{n\geq N}(\left\vert X_{n}-f\right\vert <1/k)  \right) = 0,
$$

\bigskip \noindent if and only if, for any $k\geq 1$

$$
\mathbb{P}\left( \bigcup_{N\geq 1}\bigcap_{n\geq N}(\left\vert X_{n}-f\right\vert <1/k)  \right) = 0
$$

\bigskip  \noindent if and only if, for any $\varepsilon >0$

$$
\mathbb{P}\left( \bigcup_{N\geq 1}\bigcap_{n\geq N}(\left\vert X_{n}-f\right\vert \geq \varepsilon)  \right) = 0.
$$

\bigskip \noindent \textbf{(b) Convergence in Probability}.\\

\noindent A sequence of \textbf{a.s. finite} random variables $(X_{n})_{n\geq 1}$\ defined from $(\Omega ,\mathcal{A},\mathbb{P})$\ to $\mathbb{R}$
converges in Probability with respect to the probability measure $\mathbb{P}$ to an \textit{a.s} finite random variable $X:(\Omega ,\mathcal{A},m)$\ $\longmapsto \mathbb{R}$, denoted
\begin{equation*}
X_{n}\rightarrow _{\mathbb{P}}X
\end{equation*}

\bigskip \noindent if and only for any $\varepsilon >0,$
\begin{equation*}
\mathbb{P}(\left\vert X_{n}-X\right\vert >\varepsilon)\longrightarrow 0\text{ as }
n\longrightarrow +\infty .
\end{equation*}

\bigskip \noindent \textbf{Remark}. The convergence in probability is only possible if the limit $f$ and the $X_n$'s are \textbf{a.s.} since we need to get the differences $X_n-X$. The \textit{a.s.} finiteness justifies this.\\

\noindent \textbf{NB}. It is important to notice that the inequality in $(| X_{n}-X|>\varepsilon)$ may be strict or not.\\

\bigskip \noindent \textbf{(c) - Properties of the \textit{a.s.} convergence}.\\

\noindent (c1) The \textit{a.s.} limit is \textit{a.s.} unique.\\

\noindent (c2) We have the following operations on \textit{a.s} limits :\\

\noindent Let $(X_{n})_{n\geq 1}$ and $(Y_{n})_{n\geq 1}$ be sequences of \textit{a.s.} finite functions. Let $a$ and $b$ be finite real numbers. Suppose that  $X_{n}\rightarrow X$ $a.s.$ and $Y_{n}\rightarrow Y$ $a.s.$. Let $H(x,y)$ a continuous function of $(x,y) \in D$, where $D$ is an open set of $\mathbb{R}^2$. We have : \\

\noindent \textbf{(1)} $aX_{n}+bY_{n}\rightarrow aX+Yg$ $a.s$.\\

\noindent \textbf{(2)} $X_{n}Y_{n}\rightarrow XY$ $a.s$\\

\noindent \textbf{(3)} If $\mathbb{P}(Y=0)=0$ (that is $Y$ is $a.s$ nonzero), then

$$
X_{n}/Y_{n}\rightarrow X/Y, \ a.s.
$$

\bigskip \noindent \textbf{(4)} If $(X_{n}, Y_n)_{n\geq 1} \subset D$ \textit{a.s.} and $(X,Y) \in D$ \textit{a.s.}, then

$$
H(X_{n},Y_{n})\rightarrow H(X,Y), \ a.s.\\
$$

\bigskip \noindent \textbf{(d)- \textit{a.s.} Cauchy sequences}. If we deal with \textit{a.s.} finite functions, it is possible to consider Cauchy Theory on sequences of them. And we have the following definition and characterizations.\\

\noindent \textbf{Definition}. A sequence $(X_{n})_{n\geq 1}$ of \textit{a.s.} finite functions is an $\mathbb{P}$-\textit{a.s.} Cauchy sequence if and only if

$$
\mathbb{P}\biggr(X_p-Y_q \nrightarrow 0, \ as \ (p,q) \rightarrow (+\infty,+\infty) \biggr)=0,
$$

\bigskip \noindent that is, the $\omega$ for which the real sequence $(X_n(\omega))_{n\geq 0}$ is a Cauchy sequence on $\mathbb{R}$ form an \textit{a.s.} event.\\

\bigskip \noindent \textbf{Other expressions}. A sequence $(X_{n})_{n\geq 1}$ of \textit{a.s.} finite functions is an $\mathbb{P}$-\textit{a.s.} Cauchy sequence :\\

\noindent if and only if for any $k\geq 1$,

$$
\mathbb{P}\left( \bigcap_{n\geq 1} \bigcup_{p\geq n} \bigcup_{q \geq n} (|f_p-f_q| > 1/k) \right)=0
$$

\bigskip \noindent if and only if for any $k\geq 1$,

$$
\mathbb{P}\left( \bigcap_{n\geq 1} \bigcup_{p\geq 0}  (|f_{p+n}-X_n| > 1/k) \right)=0
$$

\bigskip \noindent if and only if for any $\varepsilon >0$,

$$
\mathbb{P}\left( \bigcap_{n\geq 1} \bigcup_{p\geq n} \bigcup_{q \geq n} (|f_p-f_q| > \varepsilon) \right)=0
$$

\bigskip \noindent if and only if for any $\varepsilon>0$,

$$
\mathbb{P}\left( \bigcap_{n\geq 1} \bigcup_{p\geq 0}  (|X_{p+n}-X_n| > \varepsilon) \right)=0
$$

\bigskip \noindent \textbf{Property}. Let  $(X_{n})_{n\geq 1}$ be a sequence of \textit{a.s.} finite functions.\\ 

\noindent $(X_{n})_{n\geq 1}$ is an $\mathbb{P}$-\textit{a.s.} Cauchy sequence if and only if $(X_{n})_{n\geq 1}$ converges \textit{a.s.} to an \textit{a.s.} finite function.

\bigskip  \bigskip \noindent \textbf{(e) - Properties of the convergence in probability}.\\

\noindent \textbf{(e1)} The limit in probability is \textit{a.s.} unique.\\

\noindent \textbf{(e2)} Operation on limits in Probability.\\

\noindent The operations of limits in probability are not simple as those for \textit{a.s.} limits. The secret is that such operations are related to weak convergence. The concepts of
tightness or boundedness are needed to handle this. But we still have some general laws and complete results on operations on constant and non-random limits.\\

\noindent \noindent Let $X_{n}\rightarrow _{\mathbb{P}} X$ and $Y_{n}\rightarrow _{\textbf{}}Y$, $a\in \mathbb{R}$. We have :\\

\noindent (1) In the general case where $X$ and $Y$ are random and \textit{a.s.} finite, we have :\\

\noindent \textbf{(1a)} $X_{n}+Y_{n}\rightarrow_{\mathbb{P}}X+Y$.\\

\noindent \textbf{(2b)} $aX_{n} \rightarrow_{\mathbb{P}} aX$\\

\noindent \textbf{(2)} - Finite and constant limits in probability.\\

\noindent Let $X=A$ and $Y=B$ be constant and non-random. we have \\

\noindent \textbf{(2a)} $aX_{n}+bY_{n}\rightarrow_{\mathbb{P}} aA+bB$.\\

\noindent \textbf{(2b)} $X_{n}Y_{n}\rightarrow_{\mathbb{P}} AB$.\\

\noindent \textbf{(3c)} If $B \neq 0$,  then

$$
X_{n}/Y_{n}\rightarrow_{\mathbb{P}} A/B.
$$

\bigskip \noindent \textbf{(3d)} If $(X_{n}, Y_n)_{n\geq 1} \subset D$ \textit{a.s.} and $(A,B) \in D$, then

$$
H(X_{n},Y_{n})\rightarrow_{\mathbb{P}} H(A,B).
$$

\bigskip \bigskip \noindent \textbf{(f) - Cauchy sequence in probability or mutually convergence in probability}.\\

\noindent Here again, we deal with \textit{a.s.} finite random variables and consider a Cauchy Theory on sequences of them. And we have the following definition and characterizations.\\

\noindent \textbf{Definition}. A sequence $(X_{n})_{n\geq 1}$ of \textit{a.s.} random variables is a Cauchy sequence in probability if and only if, for any $\varepsilon$,

$$
\mathbb{P}( |X_p-X_q| > \varepsilon) \rightarrow 0 \ as \ (p,q) \rightarrow (+\infty,+\infty).
$$

\bigskip \noindent \textbf{Properties}. Let  $(X_{n})_{n\geq 1}$ be a sequence of \textit{a.s.} random variables. We have :\\ 

\noindent \textbf{P1} \noindent $(X_{n})_{n\geq 1}$ is a Cauchy sequence in probability if and only if $(X_{n})_{n\geq 1}$ converges in probability to an \textit{a.s.} random variable.\\

\noindent \textbf{P2} \noindent If $(X_{n})_{n\geq 1}$ is a Cauchy sequence in probability, then $(X_{n})_{n\geq 1}$ possesses a subsequence $(X_{n_k})_{k\geq 1}$ and an \textit{a.s.} random variable such that $f$ such that 
$$
X_{n_k} \rightarrow X \ \textit{a.s.} \ as \ k\rightarrow +\infty, 
$$

\bigskip \noindent and 

$$
X_n \rightarrow_{\mathbb{P}} X \ as \ n\rightarrow +\infty.
$$

\bigskip \noindent \textbf{(g) - Comparison between \textit{a.e.} convergence and convergence in probability}.\\

\noindent \textbf{(1)}. If $X_{n}\rightarrow X$ $a.s.$, then $X_{n}\rightarrow _{\mathbb{P}}f$.\\

\noindent The reverse implication is not true. It is only true for a sub-sequence as follows.\\

\noindent \textbf{(2)}. Let $X_{n}\rightarrow_{\mathbb{P}}X$. Then, there exists a sub-sequence $(X_{n_{k}})_{k\geq 1}$ of $(X_{n})_{n\geq 1}$ converging \textit{a.s} to $X$.\\

\noindent \textbf{Terminology}. Probability Theory results concerning a \textit{a.s.} limit is qualified as \textsl{strong}. Since such results imply versions with limits in probability which are called \textit{weak}.\\

\section{Convergence in $L^p$} \label{proba_02_conv_sec_03}

\noindent We already know that $L^p(\Omega, \mathcal{A}, \mathbb{P})$ is a Banach space, with for $X \in L^p$,

$$
\|X\|_p=\left(\mathbb{E}(|X|^p\right)^{1/p}, \ p\in [1,\infty[
$$

\bigskip \noindent and

$$
\|X\|_{\infty}=\inf \{M>0, \ |X| \leq M, \ \textsl{a.s.}\}, \ p=+\infty.
$$

\bigskip \noindent In this section, we are going to compare $L^p$ convergence and the \textsl{a.s.} convergence or the convergence in probability.\\

\noindent We restrict ourselves to the case where $p$ is finite.\\

\noindent \textbf{(a) Immediate implications}. \\

\noindent We have the following facts.

\begin{proposition}  \label{proba_02_con_sec_lp_01} $(X_n)_n\subset L^p$ and $X \in L^p$ and let $X_n \overset{L^{p}}{\longrightarrow} X$. Then :\\

\noindent (i) $X_n \longrightarrow_{\mathbb{P}} X$\\

\noindent and\\

\noindent (ii) $\|X_n\|_p \rightarrow \|X\|_p$,\\

\noindent meaning that : the convergence in $L^p$ implies the convergence in probability and the convergence of $p$-th absolute moments.
\end{proposition}

\bigskip \noindent \textbf{Proof}. $(X_n)_n\subset L^p$ and $X \in L^p$ and let $X_n \overset{L^{p}}{\longrightarrow} X$.\\

\noindent Proof of Point (i). For any $\varepsilon>0$ and by the Markov inequality, we have

$$
\mathbb{P}(|X_n-X|> \varepsilon)=\mathbb{P}(|X_n-X|^p> \varepsilon^p)\leq \frac{\|X_n-X\|_p^p}{\varepsilon} \rightarrow 0.
$$

\bigskip \noindent Thus the convergence in $L^p$ implies the convergence in probability.\\

\noindent Proof of Point (ii). This is immediate from the second triangle inequality

$$
\biggr| \|X_n\|_p - \|X\|_p \biggr| \leq \|X_n-X\|_p \rightarrow 0. \ \square
$$

\bigskip \noindent On can the question : does one of Points (i) and (ii) implies the convergence in $L^p$? We need the concepts of continuity of a sequence of real random variables. Most of the materials used below comes from \cite{loeve}.\\

\noindent \textbf{(b) Continuity of a sequence of random variables}. \\

\noindent We have already seen the notion of continuity for a real-valued $\sigma$-additive application defined on the $\sigma$-algebra $\mathcal{A}$ with respect to the probability measure 
$\mathbb{P}$ pertaining to the probability space $(\Omega, \mathcal{A}, \mathbb{P})$, which holds whenever as follows :\\

$$
\forall A \in \mathcal{A}, \ \mathbb{P}(A)=0 \Rightarrow \phi(A).
$$

\bigskip 
\noindent Such a definition may be extended to the situation where we replace $\mathbb{P}(A)=0$ by a limit of the form :

$$
\phi(A) \rightarrow 0 \ as  \ \mathbb{P}(A) \rightarrow 0, \ \ (AC01)
$$    

\bigskip 
\noindent which may be discretized in the form :

$$
\biggr((A_p)_{p\geq 0}\subset\mathcal{A} \ and \ \mathbb{P}(A_p)\rightarrow 0 \biggr) \Rightarrow
\biggr(\phi(A_p)\rightarrow 0\biggr),\ \ (AC02)
$$

\bigskip \noindent where the limits are meant as $p \rightarrow +\infty$.\\

\noindent Let $\phi=i=\phi_{X}$ be an indefinite integral associated to the absolute value of random variable $X$, that is

$$
\phi_X(A)=\int_A |X| \ d\mathbb{P}, \ A \in \mathcal{B}(\mathbb{R}).
$$
 
\bigskip \noindent We denote $B(X,c)=(|X|>c)$ for any $c>0$ and introduce the condition

$$
\lim_{c \uparrow +\infty} \phi_X(B(X,c))=0,
$$

\bigskip \noindent that is

 $$
\lim_{c \uparrow +\infty} \int_{(|X|>c)} |X| \ d\mathbb{P}=0. \ \ (CI) 
$$
 
\bigskip \noindent Let us introduce the following :\\

\noindent \textbf{Definitions}. \\

\noindent (a) A random variable $X \in \overline{\mathbb{R}}$ is $\mathbb{P}$-absolutely continuous if and only if Formula (AC01) holds.\\

\noindent (b) A random variable $X \in \overline{\mathbb{R}}$ is $\mathbb{P}$-continuously integrable  if and only if Formula (CI01) holds. $\Diamond$\\

\bigskip \noindent We have the following first result.\\

\begin{proposition} \label{proba_02_con_sec_lp_02} If $X$ is integrable, then it is $\mathbb{P}$-absolutely continuous and $\mathbb{P}$-continuously integrable.
\end{proposition}

\bigskip \noindent \textbf{Proof}. Let $X$ be integrable. Now, since $(|X|>c) \downarrow (|X|=+\infty)$ as $c\uparrow +\infty$, we get by the monotone convergence theorem (Do not forget that any limit is achieved through a discretized form)

$$
\lim_{c \uparrow +\infty} \int_{(|X|>c)} X \ d\mathbb{P}=\int_{(|X|=+\infty)} X \ d\mathbb{P}.
$$

\bigskip \noindent Since $X$ in integrable, it is \textit{a.s.} finite, that is $\mathbb{P}(|X|=+\infty)=0$, which leads to $\int_{(|X|=+\infty)} X \ d\mathcal{P}$ since the indefinite integral of the integrable random variable $X$ is continuous with respect to $\mathcal{P}$. Hence $X$ is $\mathbb{P}$-continuous integrable.\\

\noindent Now, suppose that $(A_p)_{p\geq 0}\subset\mathcal{A}$   and $\mathbb{P}(A_p)\rightarrow 0$ as $p\rightarrow +\infty$. We have for any $c>0$, $p\geq 0$,

\begin{eqnarray*}
\int_{A_p} |X| \ d\mathbb{P}&=&\int_{A_p \cap B(X,c)} |X| \ d\mathbb{P} + \int_{A_p \cap B(X,c)^c} |X| \ d\mathbb{P}\\
&\leq& \int_{B(X,c)} |X| \ d\mathbb{P} + c \mathbb{P}(A_p).
\end{eqnarray*}

\bigskip \noindent By letting $p\rightarrow +\infty$ first and next $c\uparrow +\infty$, we get Formula (AC02). Hence $X$ is $\mathbb{P}$-absolutely continuous. $\square$\\

\bigskip \noindent Now we may extend the definitions above to a sequence of integrable random variables by requiring that Formulas (AC02), page \pageref{ac02} or (CI01), page \pageref{ac02}, to hold uniformly. This gives :\\

\noindent \textbf{Definitions}. \\

\noindent (a) A sequence of integrable random variables $(X_n)_{n\geq 0} \subset L^1$, is  $\mathbb{P}$- uniformly and absolutely continuous (\textit{uac}) if and only if

$$
\lim_{\mathbb{P}(A)\rightarrow 0} \sup_{n\geq 0} \int_{A} |X_n| \ \mathbb{P}=0,  \ \ (UAC1)
$$

\bigskip \noindent which is equivalent to

$$
\forall \varepsilon>0, \exists \eta>0, \ \forall A\in \mathcal{A}, \mathbb{P}(A)<\eta \Rightarrow  \forall n\geq 0, \  \int_{A} |X_n| \ d\mathbb{P} <\varepsilon. \ \ (UAC2)
$$

\bigskip \noindent (b) A sequence of integrable random variables $(X_n)_{n\geq 0} \subset L^1$, is  $\mathbb{P}$-uniformly continuously integrable (\textit{uci}) if and only if

$$
\lim_{c\uparrow +\infty} \sup_{n\geq 0}, \ \int_{(|X_n|>c)} |X_n| \ \mathbb{P}. \ \ (UCI)
$$

\bigskip \noindent. $\Diamond$\\

\noindent \textbf{Example}. As in \cite{billingsley}, let us consider a sequence of random variables $(X_n)_{n\geq 0} \subset L^{1+r}$, $r>0$ such that

$$
\sup_{n\geq 0} \mathbb{E} |X_n|^{1+r}=C<+\infty.
$$

\bigskip \noindent Such a sequence is $\mathbb{P}$-\textit{uci} since for all $c>0$,

$$
\int_{(|X_n|>c)} |X_n| \ \mathbb{P}=\int_{(|X_n|>c)} \frac{|X_n|^{1+r}}{|X_n|^r} \ \mathbb{P}\leq c^{-r}C,
$$

\bigskip \noindent and next

$$
\sup_{n\geq 0} \int_{(|X_n|>c)} |X_n| \ \mathbb{P}\leq c^{-r}C \rightarrow 0 \ as \ c\uparrow +\infty.
$$

\bigskip \noindent Unlike the situation where we had only one integrable random variable, the two notions of $\mathbb{P}$-\textit{uac} and $\mathbb{P}$-\textit{uci} do not coincide for sequences. We have :

\begin{proposition} \label{proba_02_con_sec_lp_03}  A sequence of integrable random variables $(X_n)_{n\geq 0} \subset L^1$\\

\noindent (i) is $\mathbb{P}$-\textit{uci} \\

\noindent if and only if \\

\noindent (ii) it is $\mathbb{P}$-\textit{uci} and the sequence of integrals $(\mathbb{E}|X_n|)_{n\geq 0}$ is bounded. $\Diamond$\\
\end{proposition}

\bigskip \noindent \textbf{Proof}. Let us consider a sequence of integrable random variables $(X_n)_{n\geq 0} \subset L^1$.\\

\noindent Let us suppose that is $\mathbb{P}$-\textit{uci}. Hence by definition, by the classical results of limits in $\mathbb{R}$, where

$$
\sup_{c>0} \sup_{n\geq 0} \int_{(|X_n|>c)} |X_n| \ \mathbb{P}=C< +\infty,
$$

\bigskip \noindent next, for any $n\geq 0$, for any $c_0>0$,

$$
\mathbb{E}|X_n|=\int_{(|X_n|>c)} |X_n| \ \mathbb{P} +\int_{(|X_n|\leq c)} |X_n| \ \mathbb{P}\leq C +c_0 (|X_n|\leq c)\leq C+c_0,
$$

\bigskip \noindent and thus

$$
\sup_{n\geq 0} \mathbb{E}|X_n| \leq C+c_0<+\infty.
$$

\bigskip \noindent Besides, if we are given $(A_p)_{p\geq 0}\subset\mathcal{A}$ and $\mathbb{P}(A_p)\rightarrow 0$ as $p\rightarrow +\infty$, we have for any $c>0$, $p\geq 0$,

\begin{eqnarray*}
\int_{A_p} |X_n| \ d\mathbb{P}&=&\int_{A_p \cap B(X_n,c)} |X_n| \ d\mathbb{P} + \int_{A_p \cap B(X_n,c)^c} |X_n| \ d\mathbb{P}\\
&\leq& \int_{B(X_n,c)} |X_n| \ d\mathbb{P} + c \mathbb{P}(A_p)\\
&\leq& \int_{B(X_n,c)} |X_n| \ d\mathbb{P} + c \mathbb{P}(A_p)\\
&\leq& \sup_{n\geq 0} \int_{B(X_n,c)} |X_n| \ d\mathbb{P} + c \mathbb{P}(A_p).
\end{eqnarray*}

\bigskip \noindent By letting $p\rightarrow +\infty$ first and next $c\uparrow +\infty$, we get Formula (UCA). Hence the sequence $\mathbb{P}$-\textit{uac}.\\

\noindent Suppose now that the sequence is $\mathbb{P}$-\textit{uac} and the sequence of integrals $(\mathbb{E}|X_n|)_{n\geq 0}$ is bounded. Put

$$
\sup_{n\geq 0} \mathbb{E}|X_n|=C<+\infty.
$$

\bigskip \noindent By the Markov inequality, we have

$$
\sup_{n\geq 0} \mathbb{P}(|X_n|>c)\leq  \sup_{n\geq 0} \frac{\mathbb{E}(|X_n|}{c}\leq Cc^{-1}. \ \ (MK)
$$

\bigskip \noindent Let us apply Formula (UAC2). Let $\varepsilon>0$ and let $\eta>0$ such that 

$$
\mathbb{P}(A)<\eta \Rightarrow  \forall n\geq 0,  \  \int_{A} |X_n| \ d\mathbb{P} <\varepsilon. \ \ (MK1)
$$

\bigskip \noindent Let $c_0>0$ such that $Cc_0^{-1}<\eta/2$. By Formula (MK) above we have for all $c\leq c_0$, for all $n\geq 0$, $\mathbb{P}(|X_n|>c)\leq \eta/2<\eta$, and by Formula (MK1),

$$
\forall n\geq 0,  \  \int_{|X_n|>c} |X_n| \ d\mathbb{P} <\varepsilon
$$ 

\bigskip \noindent that is

$$
\forall c\leq c_0, \ \ \sup_{n\geq 0} \int_{|X_n|>c} |X_n| \ d\mathbb{P} \leq \varepsilon.
$$ 

\bigskip \noindent This means that the sequence is $\mathbb{P}$-\textsl{uci}. $\blacksquare$\\

\bigskip \noindent Now, we are able to give the converse of Proposition \ref{proba_02_con_sec_lp_01}.\\

\begin{theorem} \label{proba_02_con_sec_th_01} Let $(X_n)_{n\geq 0} \subset L^p$ be e sequence of elements of $L^p$, and $X$ some random variable $X \in \overline{\mathbb{R}}$. We have :\\

\noindent (a) If $X_n \overset{L^{p}}{\longrightarrow} X$, then $X \in L^p$.\\

\noindent (b) If If $X_n \overset{L^{p}}{\longrightarrow} X$, then If $X_n \longrightarrow_{\mathbb{P}} X$.\\

\noindent (c) Suppose that $X_n \longrightarrow_{\mathbb{P}} X$ and one of the three conditions holds.\\

\noindent (c1) The sequence $(|X_n|^p)_{n\geq 0}$ is $\mathbb{P}$-uniformly and absolutely integrable.\\

\noindent (c2) The sequence $(|X_n-X|^p)_{n\geq 0}$ is $\mathbb{P}$-uniformly and absolutely integrable.\\

\noindent (c3) The sequence $(|X_n|^p)_{n\geq 0}$ is $\mathbb{P}$-uniformly and continuously integrable.\\

\noindent (c4) $\|X_n\|_p \rightarrow \|X\|_p<+\infty$.\\

\noindent Then $X_n \overset{L^{p}}{\longrightarrow} X$.\\

\noindent (All the limits are meant when $n\rightarrow +\infty$).
\end{theorem}

\bigskip \noindent \textbf{Proof of Theorem \ref{proba_02_con_sec_th_01}}.\\

\noindent \textit{Proof of (a)}. First remark that the random variables  $|X_n-X|$ are \textit{a.s.} defined since the $X_n$'s are \textit{a.s.} finite. Next, by Minkowski's inequality, for any $n\geq 0$.

$$
\|X\|_p \leq \|X_n\|_p + \|X_n-X\|_p
$$

\bigskip \noindent By $X_n \overset{L^{p}}{\longrightarrow} X$, there exists $n_0\geq 0$ such that $\|X_n-X\|_p \leq 1$ and thus, $\|X\|_p \leq 1 +\|X_{n_0}\|_p <+\infty$.\\

\noindent \textit{Proof of (b)}. It is done in Proposition \ref{proba_02_con_sec_lp_01}.\\

\noindent \textit{Proof of (c)}. Suppose that $X_n \longrightarrow_{\mathbb{P}} X$.\\

\noindent \textit{Let (c1) hold}. Let $\varepsilon>0$. By Point (f) in Section \ref{proba_02_conv_sec_02}, the sequence $(X_n)_{n\geq 0}$ is of Cauchy in probability. Hence for  
$B_{r,s}(\varepsilon)=(|X_r-X_s|>(\varepsilon/2)^{1/p})$, we have

$$
\mathbb{P}(B_{r,s}(\varepsilon)) \rightarrow 0 \ as \ (r, s)\rightarrow(+\infty,+\infty). \ \ (LP1)
$$

\bigskip \noindent Since (c1) holds, we use Formula (UAC2) to find a value $\eta>0$ such that, for $C_p=2p^{p-1}$,

$$
\mathbb{P}(A)<\eta \Rightarrow  \forall n\geq 0,  \  \int_{A} |X_n|^p \ d\mathbb{P} <(\varepsilon/2C_p). \ \ (LP2)
$$

\bigskip \noindent From Formula (LP1), we can find an integer $r_0$ such that for any $r\geq r_0$, for any $s\geq 0$, 

$$
\mathbb{P}(B_{r,r+s}(\varepsilon))<\eta.
$$

\bigskip \noindent Now, based on the previous facts and the $C_p$ inequality, we have for all $r\geq r_0$ and $s\geq 0$,

\begin{eqnarray*}
\int |X_r-X_{r+s}|^p \ d\mathbb{P} &\leq& \int_{B_{r,r+s}(\varepsilon)} |X_r-X_{r+s}|^p \ d\mathbb{P}\\
&+& \int_{B_{r,s}(\varepsilon)^c} |X_r-X_{r+s}|^p \ d\mathbb{P} \ \ (L1)\\
 &\leq& C_p \biggr(\int_{B_{r,r+s}(\delta)} |X_r|^p + \int_{B_{r,r+s}(\delta)} |X_{r+s}|^p\biggr) + \varepsilon/2\\
&\leq&  \varepsilon/2+ \varepsilon/2=\varepsilon.\\
\end{eqnarray*}

\bigskip \noindent This implies that the sequence $(X_n)_{n\geq 0}$ us a Cauchy sequence in $L^p$ and since $L^p$ is a Banach space, it converges in $L^p$ to $Y$. By Point (b), we also have that 
$X_n$ converges in Probability to $Y$. Thus $X=Y$ \textit{a.s.}. Finally $\|X-n-Y\|_p$=$\|X-n-X\|_p \rightarrow 0$ and $X_n$ converges to $X$ in $L^p$.\\

\noindent \textit{Let (c2) hold}. The same method may used again. The form of $B_{r,s}(\varepsilon)$ does not change since $X$ is dropped in the difference. When concluding in Line (L1) in the last group of formulas, we use

$$
\int_{B_{r,r+s}(\varepsilon)} |X_r-X_{r+s}|^p \leq \int_{B_{r,r+s}(\varepsilon)} |X_r-X|^p \  d\mathbb{P}+ \int_{B_{r,r+s}(\varepsilon)} |X_{r+s}-X|^p \ d\mathbb{P},
$$

\bigskip \noindent and the conclusion is made similarly.\\

\noindent \textit{Let (c3) hold}. By Proposition \ref{proba_02_con_sec_lp_03}, (c1) holds and we have the results.\\

\noindent \textit{Let (c4) hold}. We are going to use the Young version of the Dominated Convergence Theorem [YCDT] (See \cite{ips-mestuto-ang}, Chapter 7, Doc 06-02, Point (06.07c)). We have

$$
x|X_n-X|^p \leq C_p (|X_n|^p + |X|^p)=Y_n.
$$

\bigskip \noindent Hence $|X_n-X|^p$ converges to zero in probability and is bounded, term by term, by a sequence $(Y_n)_{n\geq 0}$ of non-negative and integrable random variables such that :\\

\noindent (i) $Y_n$ converges to $Y=2C_p |X|^p$\\

\noindent and \\

\noindent (ii) $\int Y_n \ \mathbb{P}$ convergences to $\int Y \ \mathbb{P}$. 

\bigskip \noindent By the YDCT, we get the conclusion, that is $\|X_n-X\|_p^p=\int |X_n-X|^p \ d\mathbb{P}\rightarrow 0$ as $n\rightarrow 0$.\\

\bigskip \noindent We still have to expose a simple review of weak convergence on $\mathbb{R}$. But we prefer stating it, for once, on $\mathbb{R}^d$ in the next part.\\

\newpage
\noindent \textbf{Part B : Convergence of random vectors}. \label{proba_02_convergence_B}.\\

\section{A simple review on weak convergence} \label{proba_02_conv_sec_04}

\noindent A general introduction of the theory of weak convergence is to be found in \cite{ips-wcrv-ang}. The main fruits of that theory on $\mathbb{R}^d$ are summarized below. \\

\bigskip \noindent First of all, it is interesting that characteristic elements of probability laws on $\mathbb{R}^d$ (\textit{cdf}'s, \textit{pdf}'s, \textit{mgf}'s, characteristic functions, etc.) still play the major roles in weak convergence. \\

\noindent The main criteria for weak convergence are stated here :

\begin{theorem} \label{proba_02_conv_sec_04_thportmanteau} (A particular version of Portmanteau Theorem) Let $d$ be a positive integer. The sequence of random vectors $X_{n}:(\Omega _{n },\mathcal{A}_{n },
\mathbb{P}^{(n)})\mapsto (\mathbb{R}^d,\mathcal{B}(\mathbb{R}^d))$, $\geq 1$, weakly converges to the random vector $X:(\Omega _{\infty},\mathcal{A}_{\infty},\mathbb{P}_{\infty})\mapsto (\mathbb{R}^d,\mathcal{B}(\mathbb{R}^d))$ if and only if one of these assertions holds.\\

\bigskip \noindent (i) For any real-valued continuous and bounded function $f$ defined on $\mathbb{R}^d$,
\begin{equation*}
\lim_{n\rightarrow +\infty} \mathbb{E}f(X_n) =\mathbb{E}f(X).
\end{equation*}

\bigskip \noindent (ii) For any open set $G$ in $\mathbb{R}^d$, 
\begin{equation*}
\liminf_{n\rightarrow +\infty} \mathbb{P}_n(X_{n}\in G)\geq \mathbb{P}_{\infty}(X \in G).
\end{equation*}

\bigskip \noindent (iii) For any closed set $F$ of  $\mathbb{R}^d$, we have
\begin{equation*}
\limsup_{n\rightarrow +\infty} \mathbb{P}_{n}(X_{n}\in F)\leq \mathbb{P_{\infty}}(X \in F).
\end{equation*}

\bigskip \noindent (iv) For any Borel set $B$ of $\mathbb{R}^d$ that is $\mathbb{P}_{X}$-continuous, that is $\mathbb{P}_{\infty}(X \in \partial B)=0$, we have
\begin{equation*}
\lim_{n\rightarrow +\infty} \mathbb{P}_n(X_{n }\in B)=\mathbb{P}_{X}(B)=\mathbb{P}_{\infty}(X \in B).
\end{equation*}

\bigskip \noindent (v) For any continuity point $t=(t_{1},t_{2},...,t_{d})$ of $F_{X}$, we have,

\begin{equation*}
F_{X_n}(t) \rightarrow F_{X}(t) \text{ as } n\rightarrow +\infty.
\end{equation*}

\noindent where for each $n\geq 1$, $F_{X_n}$ is the distribution function of $X_n$ and $F_{X}$ that of $X$.\\

\noindent (vi) For any point $u=(u_{1},u_{2},...,u_{d})\in \mathbb{R}^{k}$, 
\begin{equation*}
\Phi _{X_{n}}(u)\mapsto \Phi _{X}(u) \text{ as }n\rightarrow +\infty,
\end{equation*}

\noindent where for each $n\geq 1$, $\Phi_{X_n}$ is the characteristic function of $X_n$ and $\Phi_{X}$ is that of $X$.\\

\bigskip \noindent \textbf{(c)} If the moment functions $\varphi_{X_{n}}$ exist on  $B_{n}$, $n\geq 1$ and $\varphi _{X}$ exists on $B$, where the $B_{n}$ and $B$ are neighborhoods of $0$ and
$B \cap_{n\geq 1}$, and if for any $x\in B$, 
\begin{equation*}
\Psi _{X_{n}}(x)\rightarrow \Psi _{X}(x)\text{ as }n\rightarrow +\infty, 
\end{equation*}

\noindent then $X_{n}$ weakly converges to $X$.\newline

\end{theorem}

\bigskip \noindent The characteristic function as a tool of weak convergence is also used through the following criteria.\\

\noindent \textbf{Wold Criterion}. The sequence $\{X_n, \ \ n\geq 1\} \subset \mathbb{R}^d$ weakly converges to $X \in \mathbb{R}^d$, as $n \rightarrow +\infty$ if and only if for any $a \in \mathbb{R}^d$, the sequence $\{<a,X_n>, \ \ n\geq 1\} \subset \mathbb{R}$ weakly converges to $X \in \mathbb{R}$ as $n \rightarrow +\infty$.\\

\noindent We also have :\\

\noindent \textbf{The Continuous mapping Theorem}. Assume that the sequence $\{X_n, \ \ n\geq 1\} \subset $ weakly converges to $X \in \mathbb{R}^d$, as $n \rightarrow +\infty$. Let $k\geq 1$ and let $f : \mathbb{R}^d \rightarrow \mathbb{R}^d$ be a continuous function. Then $\{f(X_n), \ \ n\geq 1\} \subset \mathbb{R}^k$ weakly converges to $f(X) \in \mathbb{R}^k$.\\

\noindent The \textit{pdf}'s may be used in the following way.\\

\begin{proposition} \label{cv.scheffeExt} These two assertions hold.\\

\noindent (A) Let $X_{n}: (\Omega _{n },\mathcal{A}_{n},\mathbb{P}^{(n)}) \mapsto (\mathbb{R}^{d},\mathcal{B}(\mathbb{R}^{d}))$ be random vectors
and $X : (\Omega_{\infty},\mathcal{A}_{\infty},\mathbb{P}_{\infty})\mapsto (\mathbb{R}^{d},\mathcal{B}(\mathbb{R}^{d}))$ another random vector, all of them absolutely continuous with respect to the Lebesgue measure denoted as $\lambda_d$. Denote $f_{X_{n}}$ the probability density function of $X_n$,  
$n\geq 1$ and by $f_{X}$ the probability density function of $X$. Suppose that we have 

\begin{equation*}
f_{X_{n}}  \rightarrow f_{X}, \text{ }\lambda_k-a.e., \text{ as } n\rightarrow +\infty.
\end{equation*}

\bigskip
\noindent Then $X_{n}$ weakly converges to $X$ as $n\rightarrow +\infty$.\\

\noindent (B) Let $X_{n }: (\Omega _{n },\mathcal{A}_{n},\mathbb{P}^{(n)})\mapsto (\mathbb{R}^{d},\mathcal{B}(\mathbb{R}^{d}))$ be discrete random vectors and $X$ : $(\Omega_{\infty},\mathcal{A}_{\infty},\mathbb{P}_{\infty})\mapsto (\mathbb{R}^{d},\mathcal{B}(\mathbb{R}^{d}))$ another discrete random vector. For each $n$, define $D_n$ the countable support of $X_n$, that 
$$
\mathbb{P}^{(n)}(X_n \in D_n)=1 \text{ and for each } x \in D_n, \text{ } \mathbb{P}^{(n)}(X_n=x)\neq 0,
$$

\bigskip
\noindent and $D_{\infty}$ the countable support of $X$. Set  $D=D_{\infty} \cup (\cup_{n \geq 1} D_n)$ and denote by $\nu$ as the counting measure on $D$. Then the probability densities of the $X_n$ and of $X$ with respect to $\nu$ are defined on $D$ by
$$
f_{X_n}(x)=\mathbb{P}^{(n)}(X_n=x), \text{ } n\geq 1, \text{ } f_{X}(x)=\mathbb{P}^{(\infty)}(X=x), \text{ }x\in D.
$$

\bigskip
\noindent If 
\begin{equation*}
(\forall x\in D), f_{X_{n }}(x)  \rightarrow f_{X}(x),
\end{equation*}

\bigskip
\noindent then $X_{n}$ weakly converges to $X$.
\end{proposition}

\bigskip \noindent \textbf{In summary}, the weak convergence in $\mathbb{R}^{d}$ holds when the distribution functions, the characteristic functions, the moment functions (if they exist) or the probability density functions (if they exist) with respect to the same measure $\nu$, point-wisely converge to the distribution function, or to the characteristic function or to moment function (if it exists), or to the probability density function (if it exists) with respect to $\nu$ of a probability measure in $\mathbb{R}^{d}$. In the case of point-wise convergence of the distribution functions, only matters the convergence for continuity points of the limiting distribution functions.\\

\bigskip \noindent  In Chapter 1 in \cite{ips-wcrv-ang}, a number of direct applications are given and a review of some classical weak convergence results are stated.

\section{Convergence in Probability and \textit{a.s.} convergence on $\mathbb{R}^d$} \label{proba_02_conv_sec_05}

\bigskip Let us denote by $\|.\|$ one of the three equivalent usual norms on  $\mathbb{R}^d$. Because of the continuity of the norm, $\|X\|$ becomes a real-valued random variable for any random vector. From this simple remark, we may extend the \textit{a.s.} convergence and the convergence in probability on $\mathbb{R}^d$ in the following way.\\

\noindent \textbf{Definitions}.\\

\noindent Let $X$ and $(X_n)_{n\geq 0}$ be, respectively, a random vector and a sequence of random vectors defined on the same on the same probability space $(\Omega, \mathcal{A},\mathbb{P})$ with values $\overline{\mathbb{R}}^d$. Let us denote by $X_n^{(j)}$ the $j$-th component of $X_n$ for each $1\leq j \leq d$, $n\geq 1$.\\

\noindent (a) The $(X_n)_{n\geq 0}$ converges \textsl{a.s.} to $X$ as $n\rightarrow +\infty$ if and only if, each sequence of components $(X_n^{(j)})_{n\geq 0}$ converges to $X_j$ as $n\rightarrow +\infty$.\\

\noindent (b) Let  $X$ and the elements of sequences $(X_n)_{n\geq 0}$ have \textsl{a.s.}-finite components. Then $(X_n)_{n\geq 0}$ converges \textsl{a.s.} to $X$ if and only if
$$
\|X_n-X\| \rightarrow 0, \ \textsl{a.s.} \ \ as  \ n\rightarrow +\infty.
$$

\noindent (c) Let  $X$ and the elements of sequences $(X_n)_{n\geq 0}$ have \textsl{a.s.}-finite components. Then $(X_n)_{n\geq 0}$ converges  to $X$ in probability if and only if
$$
\|X_n-X\| \rightarrow_{\mathbb{P}} 0, \  \ as  \ n\rightarrow +\infty.
$$

\noindent $\Diamond$\\

\bigskip \noindent For the coherence of the definition, we have to prove the equivalence between Points (a) and (b) above in the case where the random vectors have \textit{a.s.} finite components. This is let as an easy exercise.\\

\noindent We have the following properties.\\

\begin{proposition} \label{proba_02_conv_sec_05_prop_01} Let $X$, $Y$, $(X_n)_{n\geq 0}$ and $(Y_n)_{n\geq 0}$ be, respectively, two random vectors and two sequences of random vectors defined on the same on the same probability space $(\Omega, \mathcal{A},\mathbb{P})$ with values $\overline{\mathbb{R}}^d$. Let $f : \mathbb{R}^d \rightarrow \mathbb{R}^k$ be a continuous function. Finally let $a$ and $b$ tow real numbers. The limits in the proposition are meant as $n\rightarrow +\infty$.\\

\noindent (1) If $X_n \rightarrow X$, \textit{a.s}, then $X_n \rightarrow_{\mathbb{P}} X$.\\

\noindent (2) Let $X_n \rightarrow X$, \textit{a.s} and $Y_n \rightarrow Y$, \textit{a.s}. Then, we have\\

\noindent (2a) $aX_n+bY_n \rightarrow aX+bY$, \textit{a.s.}\\

\noindent and\\

\noindent (2b) $f(X_n) \rightarrow f(X)$, \textit{a.s.}\\

\noindent (3) Let $X_n \rightarrow_{\mathbb{P}} X$, and $Y_n \rightarrow_{\mathbb{P}} Y$. Then, we have\\

\noindent (3a) $aX_n+bY_n \rightarrow_{\mathbb{P}} aX+bY$.\\

\noindent and, if $X=A$ is a non-random constant vector, we have\\

\noindent (3b) $f(X_n) \rightarrow_{\mathbb{P}}  f(A)$.\\

\noindent But in general, if $f$ is a Lipschitz function, we have\\

\noindent (3c) $f(X) \rightarrow_{\mathbb{P}}  f(X)$.
\end{proposition}

\bigskip \noindent \textbf{Proofs}. By going back to the original versions on $\mathbb{R}$ for \textit{a.s.} and convergence in probability, all these results become easy to prove except Points
(3b) and (3c). But a proof of Point (2v) is given in the proof of Lemma 8 in \cite{ips-wcrv-ang} of this series. Point (3c) is proved as follows.\\

\noindent Let $X_n \rightarrow_{\mathbb{P}} X$ and let $f$ be a Lipschitz function associated to a coefficient $\rho>0$, that is

$$
\forall (x,y) \in \mathbb{R}^d, \ \|f(x)-f(y)\| \leq \rho \|x-y\|.
$$

\bigskip
\noindent Hence for any $\varepsilon>0$,

$$
(\|f(X_n)-f(X)\| >\varepsilon) \subset (\|f(X_n)-f(X)\| >\varepsilon/\rho)
$$  

\bigskip
\noindent and hence

$$
\mathbb{P}(\|f(X_n)-f(X)\| >\varepsilon) \subset \mathbb{P}(\|f(X_n)-f(X)\| >\varepsilon/\rho) \rightarrow 0.
$$  

\bigskip
\noindent Thus $\|f(X_n)-f(X)\| \rightarrow_{\mathbb{P}} 0$. $\square$\\

\noindent \textbf{Immediate implications}. Since projections are Lipschitz functions, we get the if $X_n \rightarrow_{\mathbb{P}} X$, then we also get the convergence in probability component-wise, that is : each sequence of components $(X_n^{(j)})_{n\geq 0}$ converges to $X_j$ in probability. Conversely, the convergence in probability implies the convergence if probability of the vectors. Indeed, take for example

$$
\|x\|=\max_{1\leq i \ d} |x_i|, \ x=(x_1,...,x_d)^t \in \mathbb{R}^d.
$$

\bigskip
\noindent We have, for each $n\geq 0$,

$$
(\|X_n-X\| > \varepsilon) \subset \bigcup_{1\leq i \ d} (\|X_n^{(j)}-X_j\| > \varepsilon),
$$

\bigskip \noindent which leads, for each $n\geq 0$, to

$$
\mathbb{P}(\|X_n-X\| > \varepsilon) \leq \sum_{1\leq i \ d} \mathbb{P}(\|X_n^{(j)}-X_j\| > \varepsilon).
$$

\bigskip
\noindent Since $d$ is fixed, the conclusion is obvious.\\

\section[Comparison between convergence types]{Comparison between convergence in probability and weak convergence} \label{proba_02_conv_sec_06}

\bigskip \noindent This section is reduced to the statements of results concerning the comparison between the weak convergence and the convergence in probability.\\

\noindent We remember that in the definition of weak convergence, the elements of the sequence $(X_n)_{n\geq 0}$ may have their own probability spaces. So, in general, the comparison with convergence in probability does not make sense unless we are in the particular case where all the elements of the sequence $(X_n)_{n\geq 0}$ and the limit random variable $X$ are defined on the same probability space.\\

\noindent Before, we state the results, let us give this definition.\\

\noindent \textbf{Definition}. Let $(X_n)_{n\geq 0}$ and $(Y_n)_{n\geq 0}$ be two random vectors and two sequences of random vectors defined on the same on the same probability space $(\Omega, \mathcal{A},\mathbb{P})$ with values $\overline{\mathbb{R}}^d$.\\

\noindent They are equivalent in probability if and only if :

$$
\|X_n-Y_n\| \rightarrow_{\mathbb{P}} 0, \ as \ n\rightarrow+\infty. 
$$

\bigskip \noindent They are \textit{a.s.} equivalent with respect to their \textit{a.s.} convergence or divergence if and only if

$$
\|X_n-Y_n\| \rightarrow 0 \ a.s., \ as \ n\rightarrow+\infty. \ \Diamond
$$

\bigskip \noindent We have :\\

\begin{proposition} \label{proba_02_conv_sec_06_prop_01} Let $X$, $Y$, $(X_n)_{n\geq 0}$ and $(Y_n)_{n\geq 0}$ be, respectively, two random vectors and two sequences of random vectors defined on the same on the same probability space $(\Omega, \mathcal{A},\mathbb{P})$ with values $\overline{\mathbb{R}}^d$. We have :\\

\noindent (a) The convergence in probability implies the weak convergence, that is : \\

\noindent If $X_{n}\rightarrow_{\mathbb{P}} X$, that d $Y_{n}\rightsquigarrow c$, then $(X_{n},Y_{n})\rightsquigarrow (X,c)$.\\

\noindent (b) The weak convergence and convergence in probability to a constant are equivalent, that is : \\

\noindent $X_{n}\longrightarrow _{\mathbb{P}}c$ as $n \rightarrow +\infty$ if and only if $X_{n}\rightsquigarrow c$ as $n \rightarrow +\infty$.\\

\noindent (c) Two equivalent sequences  in probability weakly converge to the same limit if one of them does.\\

\noindent (d) (Slutsky's Theorem) If \ $X_{n}\rightsquigarrow X$ \ and $Y_{n}\rightsquigarrow c$, then $(X_{n},Y_{n})\rightsquigarrow (X,c)$.\\

\noindent (e) (Coordinate-wise convergence in probability) $X_{n}\longrightarrow _{\mathbb{P}}X$ \ and $Y_{n}\longrightarrow _{\mathbb{P}}Y$ if and only if  $(X_{n},Y_{n})\longrightarrow _{\mathbb{P}}(X,Y)$.
\end{proposition}

\bigskip \noindent The proofs of all these facts are given in \cite{ips-wcrv-ang} of this series.\\

\bigskip \noindent \textbf{A comment}. We know that the result (e) does not holds in general for the weak convergence. This means that the convergence in probability implies the weak convergence but not the contrary.\\

\noindent The \textit{a.s.} equivalence takes a special shape for partial sums. Let us consider a sequence $(X_n)_{n\geq 0}$ of real random variables defined on the same on the same probability space $(\Omega, \mathcal{A},\mathbb{P})$. For a sequence of positive numbers $(c_n)_{n\geq 1}$, let us consider the truncated random variables $X_n^{(t)}=X_n 1_{|X_n|\leq c_n}$, $n\geq 1$, that is each $X_n^{(t)}$, $n\geq 1$, remains unchanged for $|X_n|\leq c_n$ but vanishes otherwise. Let us form the the partial sums $S_n=\sum_{1\leq k \leq n} X_n$ and $S_n^{(t)}=\sum_{1\leq k \leq n} X_n^{(t)}$, $n\geq 1$.\\

\noindent Let $(b_n)_{n\geq 1}$ be a sequence of real numbers converging to $+\infty$. We are going to see that the \textit{a.s.} equivalence between $(S_n^{(t)}/b_n)_{n\geq 1}$ and $(S_{n}/b_n)_{n\geq 1}$ is controlled by the series 
$$
\sum_{n\geq 0} \mathbb{P}(|X_n|>c_n). 
$$

\bigskip \noindent Indeed, since the event $(X_k \neq X_k^{(t)})$ occurs only if $(|X_k|\geq c_k)$, we have

$$
\mathbb{P}( X_n \neq X_n^{(t)}, \ i.o)=\lim_{n\uparrow +\infty} \mathbb{P}\biggr( \bigcup_{k\geq n} (X_k \neq X_k^{(c)}) \biggr)\leq
\lim_{n\uparrow +\infty} \sum_{k\geq n}\mathbb{P} (|X_k|\geq c).
$$

\noindent So if the series $\sum_{n\geq 0} \mathbb{P}(|X_n|>c_n)$ is convergent, we have $\mathbb{P}( X_n \neq X_n^{(t)}, \ i.o)$ which implies that there exists a null-set $\Omega_0^c$ such that for any $\omega \in \Omega_0$, we can find $N(\omega)$ such that for any $n\geq N(x)$, $X_n=X_n^{(t)}$ so that for $n>N(\omega)$,

$$
\biggr|\frac{S_n-S_n^{(t)}}{b_n}\biggr|=\frac{1}{b_n} |S_N(\omega)^{(t)}-S_N(\omega)|\rightarrow 0, \ as \ n\rightarrow +\infty.
$$

\bigskip \noindent This proves the claim. $\square$

%% file: proba_02_06_ang.tex
\chapter{Inequalities in Probability Theory} \label{proba_02_ineg}

\noindent Here, we are going to gather a number of some inequalities we may encounter and use in Probability Theory. Some of them are already known from the first chapters.\\

\noindent The reader may skip this chapter and comes back to it only when using, later, an inequality which is is stated here and especially when he/she wants to see the proof.\\

\noindent Unless an express specification is given, the random variables $X$, $Y$, $X_i$, $Y_i$, $i\geq 1$, which used below, are defined on the same probability space $(\Omega, \mathcal{A}, \mathbb{P})$.\\

\noindent Readers who want to read this chapter in the first place will need an earlier introduction to the notion of conditional expectation right now, instead of waiting Chapter \ref{proba_02_ce} where this notion is studied.\\

\section{Conditional Mathematical Expectation} \label{proba_02_ineg_sec_01}
 
\noindent We are going to use the Radon-Nikodym Theorem as stated in Doc 08-01 in Chapter 9 in \cite{ips-mestuto-ang}.\\

\noindent Let be given a sub-$\sigma$-algebra $\mathcal{B}$  of $\mathcal{A}$, $Y$ a measurable mapping from $(\Omega, \mathcal{A})$ to a measurable space $(E, \mathcal{F})$ and finally a measurable 
mapping $h$ from $(E, \mathcal{F})$ to $\overline{\mathbb{R}}$, endowed with the usual $\sigma$-algebra $\mathcal{B}_{\infty}(\overline{\mathbb{R}})$. We always suppose that $h(Y)$ is defined, and quasi-integrable, that is :  $\mathbb{E}(h(Y)^+)$ or $\mathcal{E}(h(Y)^-)$ is finite.\\

\noindent Now the mapping

$$
\mathcal{B} \ni B \mapsto \phi_{\mathcal{B}}(B)=\int_{B} h(Y) \ d\mathbb{P}
$$
 
\bigskip \noindent is $\sigma$-additive and is continuous with respect to $\mathbb{P}$. By the Radon-Nikodym Theorem as recalled earlier, $\phi_{\mathcal{B}}$ possesses a Radon-Nikodym derivative with respect to $\mathbb{P}$, we denoted as

$$
\frac{d\phi_{\mathcal{B}}}{d\mathbb{P}}=:\mathbb{E}(h(Y)/\mathcal{B}).
$$

\bigskip \noindent By the properties of Radon-Nikodym derivatives (please, visit again the aforementioned source if needed), we may define. \\

\noindent \textbf{Definition} \label{proba_02_ineg_sec01_defEC} If the mathematical of expectation $h(Y)$ exists, the conditional mathematical expectation of $h(Y)$ denoted as 
$$
\mathbb{E}(h(Y)/\mathcal{B}),
$$

\bigskip \noindent is the $\mathbb{P}$.\textbf{a.s} unique real-valued and $\mathcal{B}$-measurable random variable such that

$$
\forall B \in \mathcal{B}, \ \ \int_{B} h(Y) \ d\mathbb{P} = \int_{B} \mathbb{E}(h(Y)/\mathcal{B}) \ d\mathbb{P}.  \ (CE01)  
$$ 

\bigskip \noindent Moreover, $\mathbb{E}(h(Y)/\mathcal{B})$ is \textit{a.s.} finite if $h(Y)$ is integrable. $\Diamond$\\

\noindent \textbf{Extension of the Definition}. By putting $Z=1_B$ in (CE01), and by using the classical three steps method of Measure Theory and Integration, we easily get that when 
$\mathbb{E}|h(Y)|<+\infty$, Formula (CE01) is equivalent to any one of the two following others :\\

\noindent (a) For any non-negative and $\mathcal{B}$-measurable random variable $Z$, 

$$
\int Z h(Y) \ d\mathbb{P} = \int Z \mathbb{E}(h(Y)/\mathcal{B}) \ d\mathbb{P}.  \ (CE02) \ \Diamond 
$$ 

\bigskip \noindent (b) For any $\mathcal{B}$-measurable and integrable random variable $Z$,

$$
\int Z h(Y) \ d\mathbb{P} = \int Z \mathbb{E}(h(Y)/\mathcal{B}) \ d\mathbb{P}.  \ (CE02) \ \Diamond 
$$ 

\bigskip \noindent In the extent of this chapter, we will directly utilize (CE02) as a definition each time we need it.\\

\noindent The following exercises will be proved in Chapter \ref{proba_02_ce} as properties of the mathematical expectation.\\

\noindent \textbf{Exercise}. Show that the following properties.\\

\noindent (1) If $X$ is a real-valued and quasi-integrable random variable, then

$$
\mathbb{E}\biggr(\mathbb{E}(X)/\mathcal{B})\biggr)=\mathbb{E}(X).
$$

\bigskip \noindent (2) If $X$ is a real-valued and quasi-integrable random variable $\mathcal{B}$-measurable, then 

$$
\mathbb{E}(X)/\mathcal{B})=X. \ \textit{a.s.}
$$

\bigskip \noindent (3) Let $X$ be a real-valued and quasi-integrable random variable independent of $\mathcal{B}$ in the following sense : for all real-valued and quasi-integrable random variable $Z$ 
$\mathcal{B}$-measurable,

$$
\mathbb{E}(ZX) = \mathbb{E}(Z) \mathbb{E}(X).
$$

\bigskip \noindent Then, we have
$$
\mathbb{E}(X)/\mathcal{B})=\mathbb{E}(X)). \ \textit{a.s.} 
$$

\bigskip \noindent (4) If $X$ is a real-valued and quasi-integrable random variable independent $\mathcal{B}$ and if $Z$ is a real-valued, quasi-integrable  and $\mathcal{B}$-measurable random variable, we have

$$
\mathbb{E}(ZX)/\mathcal{B})= Z \mathbb{E}(X)/\mathcal{B}).
$$

\bigskip \noindent (5) If $X$ and $Y$ are real-valued random variables both non-negative or both integrable, then

$$
\mathbb{E}\biggr((X+Y)/\mathcal{B}\biggr)=\mathbb{E}(X/\mathcal{B}) + \mathbb{E}(Y/\mathcal{B}). \ \ \Diamond
$$

\section{Recall of already known inequalities} \label{proba_02_ineg_sec_0}.

\noindent \textbf{(1) Inequality of Markov}. If $X\geq 0$, then for all $x>0$,

$$
\mathbb{P}(X\geq x) \leq \frac{1}{x}. \ \Diamond
$$

\bigskip \noindent \textbf{(2) Inequality of Chebychev}. If $\mathbb{E}(X)$ exists and $X-\mathbb{E}(X)$ is \textit{a.s.} defined,  then for all $x>0$,

$$
\mathbb{P}(|X-\mathbb{E}(X)|>x) \leq \frac{\mathbb{E}(X-\mathbb{E}(X))^2}{x^2}. \ \Diamond
$$

\bigskip \noindent These two inequalities are particular forms of the following one.\\

\bigskip \noindent \textbf{(3) Basic Inequality}. (As in \cite{loeve})  Let $X$ be any real-valued random variable and $g$ be a non-null, non-decreasing and non-negative mapping from $\mathbb{R}$ to 
$\mathbb{R}$. Then for any $a \in \mathbb{R}$, we have

$$
\frac{\mathbb{E}g(X)-g(a)}{\|g(X=\|_{\infty}}\leq \mathbb{P}(X \geq a) \leq \frac{\mathbb{E}g(X)}{g(a)}. \ \ (BI01)
$$

\noindent If, in addition, $g$ is even or if $g$ satisfies

$$
\forall a\geq 0, \ g(\max(-a, a)) \geq g(a), \ \ (AI01)
$$

\noindent 

$$
\frac{\mathbb{E}g(X)-g(a)}{\|g(X=\|_{\infty}}\leq \mathbb{P}(|X| \geq a) \leq \frac{\mathbb{E}g(X)}{g(a)}. \ \ (BI02) \ \Diamond
$$

\bigskip \noindent \textbf{Proof of Formula (BI01)}. The mathematical expectation $\mathbb{E}(g(X))$ exists since $g$ is of constant sign. By using the same method of establishing the Markov inequality, we have

$$
\mathbb{E}g(X) \geq \int_{X \geq a} \ d\mathbb{P} \geq g(a) \mathbb{P}(X\geq a), \ \ (BI03)
$$

\bigskip \noindent where we used the non-decreasingness of $g$. So, we get the right-hand inequality of Formula (BI01) even if $g(a) = 0$. We also gave

$$
\int_{X \geq a} \ d\mathbb{P} \leq \|g(X=\|_{\infty} \mathbb{P}(|X| \geq a)
$$

\noindent and

$$
\int_{X < a} \ d\mathbb{P} \leq g(a) \mathbb{P}(|X| < a)\leq g(a)
$$

\bigskip \noindent and by these formulas,

\begin{eqnarray*}
\mathbb{E}g(X)&=&\int_{X\geq a} g(X) \ d\mathbb{P} +\int_{X<a} g(X) \ d\mathbb{E}\\
&=&\|g(X=\|_{\infty} \mathbb{P}(|X| \geq a)+g(a),
\end{eqnarray*}

\noindent that is
$$
\mathbb{E}g(X) \leq \|g(X=\|_{\infty} \mathbb{P}(|X| \geq a)+g(a),
$$

\noindent which gives the left-and inequality in Formula (BI01) even if $\|g(X=\|_{\infty}=+\infty$ (it cannot be zero by assumption). $\square$\\

\noindent \textbf{Proof of Formula (BI02)}. Since $\mathbb{P}(X\leq a)\leq \mathbb{P}(|X-\leq a)$, we pnly have to justify the right-hand inequality of (BI02). But we may use the simple remark that 
$X\geq \max(-a,a)$ on $(|X|\geq a)$ to modify (BI03) as follows 
$$
\mathbb{E}g(X) \geq \int_{|X| \geq a} \ d\mathbb{P} \geq g(\max(-a,a)) \mathbb{P}(|X|\geq a), \ \ (BI04)
$$

\noindent So, using Assumption (AI01) - which holds if $g$ is even - allows to conclude. $\square$\\

\bigskip \noindent \textbf{(4) H\"{o}lder Inequality}. Let $p>1$ and  $q>1$ be two conjugated positive rel numbers, that is, $1/p+1/q=1$ and let  
\begin{equation*}
{\Large 
\begin{array}{cccc}
X,Y: & (\Omega ,\mathcal{A},\mathbb{P}) & \mapsto & \mathbb{R}%
\end{array}%
,}
\end{equation*}

\bigskip \noindent  be two random variables $X \in L^{p}$ and $Y \in L^{q}$. Then $XY$ is integrable and we have
 
\begin{equation*}
\left\vert \mathbb{E}(XY)\right\vert \leq \left\Vert X\right\Vert _{p}\times
\left\Vert Y\right\Vert _{q},
\end{equation*}

\bigskip \noindent  where for each $p\geq 1$, $\left\Vert X\right\Vert _{p}=(\mathbb{E}(|X|^p)^{1/p}$.\\

\bigskip \noindent  \textbf{(5) Cauchy-Schwartz's Inequality}. For $p=q=2$, the H\"{o}lder inequality becomes the Cauchy-Schwartz one :
\begin{equation*}
\left\vert \mathbb{E}(XY)\right\vert \leq \left\Vert X\right\Vert _{2}\times
\left\Vert Y\right\Vert _{2}.
\end{equation*}

\bigskip \noindent  \textbf{(6) Minskowski's Inequality}. Let $p\geq 1$ (including $p=+\infty$). If $X$ and $Y$ are in $L^{p}$, then we have 
\begin{equation*}
\left\Vert X+Y\right\Vert _{p}\leq \left\Vert X\right\Vert _{p}+\left\Vert
Y\right\Vert _{p}.
\end{equation*}

\bigskip \noindent  \textbf{(7) $C_p$  Inequality}. Let $p \in [1,+\infty[$. If $X$ and $Y$ are in $L^{p}$, then for $C_p=2^{p-1}$, we have 
\begin{equation*}
\|X+Y\|_p^p\leq C_p (\|X\|_p^p+\|Y\|_p^p).
\end{equation*}

\bigskip \noindent \textbf{(8) Ordering the spaces $L^p, \ p\geq 1$}.\\

\noindent Let $1<p<q$, $p$ finite but $q \in [1, +\infty]$. Let $X\in L^{q}$. Then $X\in L^{p}$ and 

$$
\left\Vert X\right\Vert _{p}\leq \left\Vert X\right\Vert _{q} \leq \left\Vert X\right\Vert _{+\infty}.
$$

\bigskip \noindent  \textbf{(9) Jensen's Inequality}.\\

\noindent Let $\phi $ be a convex function defined from a closed interval $I$ of $\mathbb{R}$ to $\mathbb{R}$. Let $X$ be a \textit{rrv} with values in $I$ such that $\mathbb{E}(X)$ is finite. Then $\mathbb{E}(X)\in I$ and 
\begin{equation*}
\phi (\mathbb{E}(X))\leq \mathbb{E}(\phi (X)).
\end{equation*}

\bigskip \noindent  \textbf{(10) Inequality for two convex functions a random variable}. Let $g_i$, $i\in \{1,2\}$ be two finite real-valued convex and increasing functions (then invertible function as increasing and continuous functions) such that $g_2$ is convex in $g_1$ meaning that $g_2g_1^{-1}$ is convex. For any real-valued random variable $X$ such that
$X$ and $g_1(X)$ are integrable, we have

\begin{equation*}
g_{1}^{-1}(\mathbb{E}(g_{1}(Z))\leq g_{2}^{-1}(\mathbb{E}(g_{2}(Z)).
\end{equation*}

\bigskip
\noindent \textbf{(11) Bonferroni's Inequality}.\\

\noindent Let $A_1$, ..., $A_n$ be measurable subsets of $\Omega$, $n\geq 2$. Define

\begin{eqnarray*}
\alpha_0&=&\sum_{1\leq j\leq n}\mathbb{P}(A_{j})\\
\alpha_1&=&\alpha_0 - \sum_{1\leq i_{1}<t_{2}\leq n} \mathbb{P}(A_{i_{1}}A_{i_2})\\
\alpha_2&=&\alpha_1 +  \sum_{1\leq i_{1}<...<t_{3}\leq n} \mathbb{P}(A_{i_{1}}...A_{i_3})\\
\cdots &=& \cdots\\
\alpha_r&=&\alpha_{r-1} + (-1)^{r+1} \sum_{1\leq i_{1}<...<t_{r}\leq n} \mathbb{P}(A_{i_{1}}...A_{i_r})\\
\cdots &=& \cdots\\
\alpha_r&=&\alpha_{r-1} + (-1)^{n+1}  \mathbb{P}(A_1 A_2 A_3 ... A_{n}).
\end{eqnarray*}

\bigskip \noindent Let $p = n \text{ mod } 2$, that is $n=2p+1+h$, $h\in \{0,1\}$. We have the Bonferroni's inequalities : if $n$ is odd,

$$
\alpha_{2k+1}\leq \mathbb{P}\biggr( \bigcup_{1\leq j\leq n}A_n\biggr) \leq \alpha_{2k}, \ k=0,...,p  \ (BF1)
$$

\bigskip \noindent and if $n$ is even,

$$
\alpha_{2k+1}\leq \mathbb{P}\biggr( \bigcup_{1\leq j\leq n}A_j \biggr) \leq \alpha_{2k}, \ k=0,...,p-1. \ \ (BF2)
$$

\bigskip
\section{Series of Inequalities} \label{proba_02_ineg_sec_03}

\bigskip \noindent  \textbf{(12) Order relations for conditional expectations}. Let $X$ and $Y$ be two real-valued random variables such that $X \leq Y$. Let $\mathcal{B}$ be a $\sigma$-sub-algebra 
of $\mathcal{A}$. Then, whenever the expressions in the two sides make sense and are finite, we have

$$
\mathbb{E}(X/\mathcal{B}) \leq \mathbb{E}(Y/\mathcal{B}) \ \textit{a.s.} \ (CE03)
$$

\bigskip \noindent Besides, the conditional expectation is a contracting operator in the following sense : for any real-valued and quasi-integrable random variable $X$, we have

$$
|\mathbb{E}(X/\mathcal{B})| \leq \mathbb{E}(|X|/\mathcal{B}). \ \ (CE04) \ \Diamond
$$

\bigskip \noindent \textbf{Proof}. Suppose that all the assumptions hold. We have for all $B \in\mathcal{B}$,

$$
\int_B \mathbb{E}(X/\mathcal{B}) \ d\mathbb{P}=\int_B X \ d\mathbb{P} \leq \int_B Y=\int_B \mathbb{E}(Y/\mathcal{B}) \ d\mathbb{P}
$$

\bigskip \noindent Take an arbitrary $\varepsilon>0$ and set $B(\varepsilon)=(\mathbb{E}(X/\mathcal{B} > \mathbb{E}(Y/\mathcal{B})+\varepsilon)$. It it sure that $B_0 \in\mathcal{B}$ and we have

\begin{eqnarray*}
\int_{B(\varepsilon)} \mathbb{E}(X/\mathcal{B}) \ d\mathbb{P}&\geq&\biggr(\int_{B(\varepsilon)} \mathbb{E}(Y/\mathcal{B})+\varepsilon\biggr) \ d\mathbb{P}\\
 &\geq&\biggr(\int_{B(\varepsilon)} \mathbb{E}(Y/\mathcal{B}) \ d\mathbb{P}\biggr)+\varepsilon \mathbb{P}(B(\varepsilon)).
\end{eqnarray*}

\bigskip \noindent The two last formulas cannot hold together unless $\mathbb{P}(B(\varepsilon))=0$ for all $\varepsilon>0$. By the Monotone convergence Theorem, we get that 
$\mathbb{P}(\mathbb{E}(X/\mathcal{B} > \mathbb{E}(Y/\mathcal{B})=0$, which proves Inequality (CE03). To prove Inequality (CE04), we apply (CE03) and Point (4) in the exercise in Section \ref{proba_02_ineg_sec_01} to $X \leq |X|=X^{+} + X^{-}$ and to $-X \leq |X|=X^{+} + X^{-}$, we get

\begin{eqnarray*}
|\mathbb{E}(X/\mathcal{B})|&=&\max(-\mathbb{E}(X/\mathcal{B}),\mathbb{E}(X/\mathcal{B})\\
&=&\max(\mathbb{E}(-X/\mathcal{B}),\mathbb{E}(X/\mathcal{B})\leq \mathbb{E}(|X|/\mathcal{B}).
\end{eqnarray*}

\bigskip \noindent  \textbf{(13) Jensen's Inequality for Conditional Mathematical Expectations}. Let $\mathcal{B}$ be a $\sigma$-sub-algebra of $\mathcal{A}$. Let $\phi $ be a convex function defined from a closed interval $I$ of $\mathbb{R}$ to $\mathbb{R}$. Let $X$ be a \textit{rrv} with values in $I$ such that $\mathbb{E}(X)$ is finite. Then $\mathbb{E}(X)\in I$ and 
\begin{equation*}
\phi (\mathbb{E}(X/\mathcal{B}))\leq \mathbb{E}(\phi(X)/\mathcal{B}).
\end{equation*}

\noindent \textbf{Proof}. It will be given on Chapter \ref{proba_02_ce}, Theorem \ref{proba_02_jensenME} (See page \pageref{proba_02_jensenME})
$\square$\\

\bigskip \noindent  \textbf{(14) Kolmogorov's Theorem for sums independent random variables}.\\

\noindent Let $X_1$, ..., $X_n$ be independent centered and square integrable random variables. We denote $\mathbb{V}ar(X_i)=\sigma_i^2$, $1\leq i \leq n$.  Let $c$ be a non-random number
(possibly infinite) satisfying
$$
\sup_{1\leq k \leq n} |X_k| \leq c \ a.s.
$$

\bigskip \noindent Denote the partial sums by

$$
S_0=0, \ S_k=\sum_{i=1}^k X_i, \ k\geq 1  \text{ and } s_0=0, \ s_k^2=\sum_{i=1}^k \sigma_i^2.
$$

\bigskip \noindent We have the double inequality, for any $\varepsilon$

$$
1- \frac{(\varepsilon+c)^2}{s_n^2} \leq \mathbb{P}(max(|S_1|, |S_2|, ...,|S_n|) \geq \varepsilon) \leq \varepsilon^{-2} s_n^2. \ (KM01)
$$

\bigskip \noindent \textbf{Proof}. We follow the proof in \cite{loeve}. Let $\varepsilon>0$ and put 

$$
A_0=\Omega, \ A_1=(|S_1|<\varepsilon), \ A_k=(|S_1|<\varepsilon, ...,  |S_k|<\varepsilon), \ k\geq 2.
$$

\bigskip \noindent We easily see that the sequence $(A_k)_{1\leq k \leq n}$ is non-increasing and we have

$$
B_2= A_{1} \setminus A_2=(|S_1|<\varepsilon, |S_{2}|\geq \varepsilon),
$$

$$
B_k= A_{k-1} \setminus A_k = (|S_1|<\varepsilon, ...,  |S_{k-1}|<\varepsilon, |S_k|\geq \varepsilon), \ k\geq 3.
$$

\bigskip \noindent We also have

$$
A_n^c= \sum_{1\leq k \leq n} B_k.
$$

\bigskip \noindent To see this quickly, say that $A_n^c=\cup_{1\leq k \leq n} C_k$, where $C_k=(|S_k|\geq \varepsilon)$. We are now accustomed to how rendering a union into a sum of sets since the course of Measure Theory and Integration by taking $D_1=C_1$, $D_2=C_1^c \cap C_2$, $D_k=C_1^c \cap ... \cap C_{k-1}^c C_k$, $k\geq 3$ to have

$$
\bigcup_{1\leq k \leq n} C_k = \sum_{1\leq k \leq n} D_k.
$$ 

\bigskip \noindent We have just to check that the $D_k$'s are exactly the $B_k$'s. In the coming developments, we repeatedly use the fact that an indication function is equal to any of its positive power. Now, for any $1\leq k \leq n$, we may see that $S_k 1_{B_k}$ is independent of $S_n-S_k$ (even when k=n with $S_n-S_k=0$). Reminding that the $S_k$'s are centered, we have

\begin{eqnarray*}
\int_{B_k} S_n^2 \ d\mathbb{P}&=&\mathbb{E}(S_n 1_{B_k})^2\\
&=&\mathbb{E} \biggr( S_k 1_{B_k} + (S_n-S_k) 1_{B_k} \biggr)\\
&=&\mathbb{E} (S_k 1_{B_k})^2 + \mathbb{E}((S_n-S_k) 1_{B_k})^2 + 2 \mathbb{E}((S_k 1_{B_k})(S_n - S_k)) \ \ (L02)\\
&=&\mathbb{E} (S_k 1_{B_k})^2 + \mathbb{E}((S_n-S_k) 1_{B_k})^2 \ \ (L03)\\
&\geq &\mathbb{E} (S_k 1_{B_k})^2 \geq \varepsilon^2 \mathbb{P.
}(B_k).
\end{eqnarray*}

\bigskip \noindent Line (L3) derives from Line (L2) by the fact that $(S_k 1_{B_k})$ and $(S_n - S_k)$ are independent and $S_n-S_k$ is centered. Hence, we get for each $1\leq k \leq n$,

$$
\int_{B_k} S_n^2 \ d\mathbb{P} \geq \varepsilon^2 \mathbb{P}(B_k).
$$

\bigskip \noindent By summing both sides over $k \in  \{1,...,k\}$ and by using the decomposition of $A_n^c$ into the $B_k$'s, we get

$$
\int_{A_n^c} S_n^2 \ d\mathbb{P} \geq \varepsilon^2 \mathbb{P}(A_n).
$$

\bigskip \noindent which, by  the simple remark that

$$
\sum_{i=1}^k \sigma_i^2=s_n^2=\int S_n^2 \ d\mathbb{P}
$$

\bigskip \noindent leads to 

$$
\sum_{i=1}^k \sigma_i^2 \geq \varepsilon^2 \mathbb{P}(A_n),
$$

\bigskip \noindent which is the right-side inequality in Formula (KM01).\\

\noindent To prove the left-side inequality, let us start by remarking that for $2\leq k \leq n$,

$$
S_k 1_{A_{k-1}} = S_{k-1} 1_{A_{k-1}} + X_k 1_{A_{k-1}} = S_k 1_{A_{k}} + S_k 1_{B_{k}}.
$$

\bigskip \noindent Now, on one side, we have

\begin{eqnarray*}
\mathbb{E}(S_{k-1} 1_{A_{k-1}} + X_k 1_{A_{k-1}})^2&=&\mathbb{E}(S_{k-1} 1_{A_{k-1}})^2\\
&+& \mathbb{E}(X_k 1_{A_{k-1}})^2 + 2\mathbb{E}((S_{k-1} 1_{A_{k-1}}) X_k) \ \ (L11) \\
&=&\mathbb{E}(S_{k-1} 1_{A_{k-1}})^2 + \mathbb{E}(X_k 1_{A_{k-1}})^2 \ \ (L12) \\
&=&\mathbb{E}(S_{k-1} 1_{A_{k-1}})^2 + \sigma_k^2 \mathbb{P}(A_{k-1}). \  \ (L13) \\
\end{eqnarray*}

\noindent Line (L12) derives from Line (L11) since $S_{k-1} 1_{A_{k-1}}$ and $X_k$ are independent and $X_k$ is centered. Line (L13) derives from Line (L12) since $X_k^2$ is independent of 
$1_{A_{k-1}}$.\\

\noindent On the other side, we have

\begin{eqnarray*}
\mathbb{E}(S_k 1_{A_{k}} + S_k 1_{B_{k}})^2&=&\mathbb{E}(S_k 1_{A_{k}})^2 + \mathbb{E}(S_k 1_{B_{k}})^2 + 2\mathbb{E}((S_k S_k)(1_{A_{k}} 1_{B_{k}}))\\
&=&\mathbb{E}(S_k 1_{A_{k}})^2 + \mathbb{E}(S_k 1_{B_{k}})^2,
\end{eqnarray*}

\bigskip \noindent since the sets $A_k$ and $B_k$ are disjoint [recall that $B_k=A_{k-1}\setminus A_k=A_{k-1}\cap A_k^c\subset A_k^c$].\\

\noindent We get for $2\leq k \leq n$,

$$
\mathbb{E}(S_{k-1} 1_{A_{k-1}})^2 + \sigma_k^2 \mathbb{P}(A_{k-1})=\mathbb{E}(S_k 1_{A_{k}})^2 + \mathbb{E}(S_k 1_{B_{k}})^2. \ \ (KM02)
$$ 

\bigskip \noindent But the expression $S_k 1_{B_{k}}$, which is used in last term in the right-hand member in Formula (KM02) is bounded as follows

$$
|S_k 1_{B_{k}}| \leq |S_{k-1} 1_{B_{k}}| + |X_k 1_{B_{k}}| \leq (\varepsilon +c) 1_{B_k}.
$$

\bigskip \noindent Hence the last term in the right-hand member in Formula (KM02) itself is bounded as follows

$$
\mathbb{E}(S_k 1_{B_{k}})^2 \leq (\varepsilon +c)^2 \mathbb{P}(B_k).
$$

\bigskip \noindent Further, we may bound below the last term in the left-hand member in Formula (KM02) by $s_k^2 \mathbb{P}(A_{n})$, to get for $2 \leq k \leq n$

$$
\mathbb{E}(S_{k-1} 1_{A_{k-1}})^2 + \sigma_k^2 \mathbb{P}(A_{n}) \leq \mathbb{E}(S_k 1_{A_{k}})^2 + (\varepsilon +c)^2 \mathbb{P}(B_k). \ \ (KM03)
$$ 

\bigskip \noindent Now, we may sum over $k \in \{2,...,n\}$ in both sides to get in the left-hand side

$$
\sum_{k=1}^{n-1} \mathbb{E}(S_{k} 1_{A_{k}})^2 + \sum_{k=2}^{n} \sigma_k^2 \mathbb{P}(A_{n}) \ \ (KM03a)
$$

\bigskip \noindent and in the right-hand side, by rigorously handling the ranges of summation and by using the decomposition of $A_n$'s into the $B_k$'s, we have

\begin{eqnarray*}
&&\sum_{k=1}^{n} \mathbb{E}(S_{k} 1_{A_{k}})^2 - \mathbb{E}(S_{1} 1_{A_{1}})^2 + (\varepsilon +c)^2 \mathbb{P}(A_n^c \setminus B_1)\\
&=&\sum_{k=1}^{n} \mathbb{E}(S_{k} 1_{A_{k}})^2 - \mathbb{E}(S_{1} 1_{A_{1}})^2 + (\varepsilon +c)^2 (\mathbb{P}(A_n^c) - \mathbb{P}(B_1))\\
&\leq&\sum_{k=1}^{n} \mathbb{E}(S_{k} 1_{A_{k}})^2 - \mathbb{E}(S_{1} 1_{A_{1}})^2 + (\varepsilon +c)^2 \mathbb{P}(A_n^c) - (\varepsilon +c)^2 \mathbb{P}(B_1). \  \ (KM03b))\\
\end{eqnarray*}

\noindent By moving the first term in (KM03a) to the right-hand member in (KM03) and by moving the terms in (KM03b) which are preceded by a minus sign to the left-hand member in (KM03) and by reminding that $B_1=A^c$ and $S_1=X_1$, we get

$$
a + \sum_{k=2}^{n} \sigma_k^2 \mathbb{P}(A_{n}) \leq \mathbb{E}(S_{n} 1_{A_{n}})^2 + (\varepsilon +c)^2 \mathbb{P}(A_n^c), \ \ (KM03c)
$$

\bigskip \noindent where

$$
a=\mathbb{E}(X_{1} 1_{A_{1}})^2 + (\varepsilon +c)^2 \mathbb{P}(A_1^c).
$$

\bigskip \noindent But, since $|X_1|\leq c$, we have

\begin{eqnarray*}
\sigma_1^2&=&\mathbb{E}(X_{1}^2)=\mathbb{E}(X_{1} 1_{A_{1}})^2 + \mathbb{E}(X_{1} 1_{A_{1}^c})^2\\
&\leq& \mathbb{E}(X_{1} 1_{A_{1}})^2 + c^2 \mathbb{P}(A_1^c)\\
&\leq& \mathbb{E}(X_{1} 1_{A_{1}})^2 + (\varepsilon+c)^2 \mathbb{P}(A_1^c)=a.\\
\end{eqnarray*}

\bigskip \noindent Since $s_1^2 \mathbb{P}(A_n) \leq s_1^2 \leq a$, we may bound below $a$ by $s_1^2 \mathbb{P}(A_n)$ in Formula (KM03c) to set

\begin{eqnarray*}
\left(\sum_{k=1}^{n} \sigma_k^2\right) \mathbb{P}(A_{n}) &\leq& \mathbb{E}(S_{n} 1_{A_{n}})^2 + (\varepsilon +c)^2 \mathbb{P}(A_n^c)\\
&\leq & \varepsilon^2 \mathbb{P}(A_n) + (\varepsilon +c)^2 \mathbb{P}(A_n^c)\\
&\leq & (\varepsilon+c)^2 \mathbb{P}(A_n) + (\varepsilon +c)^2 \mathbb{P}(A_n^c)=(\varepsilon +c)^2,
\end{eqnarray*}

\bigskip \noindent which implies

$$
\left(\sum_{k=1}^{n} \sigma_k^2\right) (1-\mathbb{P}(A_{n}^c)) \leq (\varepsilon +c)^2, 
$$

\bigskip \noindent and hence

$$
\mathbb{P}(A_{n}^c) \geq 1-\frac{(\varepsilon +c)^2}{\sum_{k=1}^{n} \sigma_k^2},
$$

\bigskip \noindent which is the first inequality in Formula (KM01). The proof is complete now. $\blacksquare$\\

\bigskip \noindent \textbf{(15) Maximal inequality for sub-martingales}.\\

\noindent Let $X_1$, ..., $X_n$ be rel-valued integrable random variables. Let us consider the following sub-$\sigma$-algebras : for $1\leq k \leq n$,

$$
\mathcal{B}_k=\sigma(\{X_j^{-1}(B), \ 1\leq j \leq k, \ B \in \mathcal{B}_{\infty}(\overline{\mathbb{R}})\}).
$$

\bigskip \noindent In clear, each $\mathcal{B}_k$ is the smallest $\sigma$-algebra rendering measurable the mapping $X_j$, $1\leq j\leq k$. It is also clear that $(\mathcal{B}_k)_{1\leq k\leq n}$ is an non-decreasing sequence of sub-$\sigma$-algebras of $\mathcal{A}$.\\

\noindent \textbf{Definition}. The sequence $(X_k)_{1\leq k\leq n}$ is a martingale if and only if

$$
\forall 1\leq k_1 \leq k_2 \leq n, \ \forall A \in \mathcal{B}_{k_1}), \ \int_{A} X_{k_2} \ d\mathbb{P} = \int_{A} X_{k_1} \ d\mathbb{P},
$$

\bigskip \noindent and is a sub-martingale if and only if

$$
\forall 1\leq k_1 \leq k_2 \leq n, \ \forall A \in \mathcal{B}_{k_1}), \ \int_{A} X_{k_2} \ d\mathbb{P} \geq  \int_{A} X_{k_1} \ d\mathbb{P}. \ \Diamond.
$$

\bigskip \noindent Let us adopt the notations given in Inequality (11).\\

\noindent If  $(X_k)_{1\leq k\leq n}$ is a sub-martingale, we have

$$
\mathbb{P}(\max(X_1, X_2, ...,X_n) \leq \varepsilon) \leq \varepsilon^{-1} \mathbb{E}(X_n). (IM01)
$$

\bigskip \noindent \textbf{Proof}. It is clear that

$$
C=(\max(X_1, X_2, ...,X_n) \geq \varepsilon)=\bigcup_{1\leq k \leq n} (X_j \geq \varepsilon)=\sum_{1\leq k \leq n} C_k,
$$

\bigskip \noindent with

$$
C_1=(X_1 \geq \varepsilon), \ C_2=(X_1 < \varepsilon, X_2 \geq \varepsilon), \ C_k=(X_1 < \varepsilon, ..., X_{k-1} < \varepsilon, X_k \geq \varepsilon), \ k\geq 3.
$$

\bigskip \noindent We remark that $C_k \in \mathcal{B}_{\ell}$, for all $1\leq k \leq n$, $k\leq \ell$. We have

\begin{eqnarray*}
\mathbb{E}(X_n)&=& \int X_n \ d\mathbb{P} \geq \int_{X} X_n \ d\mathbb{P} \ \ (L51)\\
&=& \sum_{1\leq k \leq n} \int_{C_k} X_n \ d\mathbb{P} \\
&\geq& \sum_{1\leq k \leq n} \int_{C_k} X_k \ d\mathbb{P} \ \ (L53)\\
&\geq& \sum_{1\leq k \leq n} \varepsilon \mathbb{P}(C_k) \ \ (L54)\\
&=& \varepsilon \sum_{1\leq k \leq n}  \mathbb{P}(C_k)= \varepsilon  \mathbb{P}(C). \ \ (L55)\\
\end{eqnarray*}

\noindent In Line (53), we applied the definition of a sub-martingale. In Line (L54), we applied that $X_k\geq \varepsilon$ on $C_k$. Finally, the combination of Lines (L51) and (L55) gives

$$
\mathbb{P}(C) \leq \varepsilon^{-1} \mathbb{E}(X_n),
$$ 

\bigskip \noindent which is Formula (MT01). $\blacksquare$\\

\bigskip \noindent \textbf{(16) - Kolmogorov's Exponential bounds}.\\

\noindent Let us fix an integer $n$ such that $n\geq 1$. Suppose that we have $n$ independent and centered random variables on the same probability space, as previously, which is \textit{a.s.} bounded. As usual $S_n$ is the partial sum at time $n$ with variance $s_n^2$. We fix $n$ such that $s_n>0$. Define

$$
c=\max_{1\leq k\leq n} \frac{X_k}{s_n} < +\infty.
$$

\bigskip \noindent  The following double inequality which is proved in Point (A1) in Chapter \ref{proba_02_appendix}, Section 
\ref{proba_02_appendixFacts} (page \pageref{proba_02_appendixFactsA1}) will be instrumental in our proofs :

$$
\forall t \in \mathbb{R}_{+}, \ e^{t(1-t)}\leq 1+t \leq e^t. \ (EB1)
$$

\bigskip \noindent Now, let us begin by the following Lemma, which is part, of the body of exponential bounds.

\begin{lemma} \label{expoboundlem01}
\noindent Let $X$ be a centered random variable which is bounded, in absolute value, by $c<+\infty$. Let us denote $\mathbb{E}{X^2}=\sigma^2$. Then for any $t>0$ such that $tc\leq 1$, we have

\begin{eqnarray*}
\mathbb{E}e^{tX} <\exp\biggr( \frac{t^2\sigma^2}{2}\left(1 + \frac{tc}{2} \right)\biggr) \ (EB2)
\end{eqnarray*}

\bigskip \noindent and

\begin{eqnarray*}
\mathbb{E}e^{tX} > \exp\biggr( \frac{t^2\sigma^2}{2}\left(1 - tc)\right)\biggr). \ \ (EB3)
\end{eqnarray*}
\end{lemma}

\bigskip \noindent \textbf{Proof of Lemma \ref{expoboundlem01}}. \noindent We begin to remark that the \textit{mgf} $t \mapsto \mathbb{E}e^{tX}$ admits an infinite expansion on the whole real line of the form 
$$
\mathbb{E}e^{tX} = 1 + \frac{t^2}{2!}\mathbb{E}{X^2} + \frac{t^3}{3!}\mathbb{E}{X^3} + \dots.
$$

 \noindent For $t>0$ and $tc\leq 1$, we have $\mathbb{E}{X^{2+\ell}}\leq \sigma^2 c^{\ell}$ for $\ell>0$. Hence

\begin{eqnarray*}
\mathbb{E}e^{tX} &=& 1 + \frac{t^2\sigma^2}{2}\biggr(1 + 2\biggr(\frac{tc}{3!} +\frac{(tc)^2}{4!}+  \dots\biggr)\biggr)\\
&=&1 + \frac{t^2\sigma^2}{2} \biggr( \sum_{k\geq 3} \biggr(\frac{(tc)^{k-2}}{k!}\biggr). \\
\end{eqnarray*}

\noindent Hence, by using the left inequality in Formula (EB1), we have

\begin{eqnarray*}
&&\mathbb{E}e^{tX} - \biggr(1 + \frac{(t^2\sigma^2}{2}\left(1 + \frac{tc}{2} \right)\\
&\leq&  t^2  \sigma^2 \biggr(2 \sum_{k\geq 3} \biggr(\frac{(tc)^{k-3}}{k!} - \frac{tc}{2}\biggr)\\
&=&  2t^2 (tc) \sigma^2 \biggr(\sum_{k\geq 3} \biggr(\frac{(tc)^{k-3}}{k!} - \frac{1}{4}\biggr)\\
&\leq&  2t^2 (tc) \sigma^2 \biggr(\sum_{k\geq 3} \biggr(\frac{1}{k!} - \frac{1}{4}\biggr) \ (L23)\\
&\leq&  2t^2 (tc) \sigma^2 (e-7/4)\leq 0, 
\end{eqnarray*}

\noindent where we used $tc \leq 1$ in Line (L23). Hence

$$
\mathbb{E}e^{tX} \leq \biggr(1 + \frac{(t^2\sigma^2}{2}\left(1 + \frac{tc}{2} \right) \leq 
\exp\biggr( \frac{t^2\sigma^2}{2}\left(1 + \frac{tc}{2} \right)\biggr),
$$

\bigskip \noindent which proves Formula (EB2).\\

\noindent To prove the left-hand inequality, we remark that $\mathbb{E}{X^{2+\ell}}\geq \sigma^2 (-c)^{\ell}$ for $\ell>0$, we also have

\begin{eqnarray*}
\mathbb{E}e^{tX} &\geq & 1 + \frac{t^2\sigma^2}{2}\biggr(1 + 2\biggr(\frac{-tc}{3!} +\frac{(-tc)^2}{4!}+  \dots\biggr)\biggr)\\
&\geq & 1 + \frac{t^2\sigma^2}{2}\biggr(1 - 2\biggr(\frac{tc}{3!} -\frac{(tc)^2}{4!}+  \dots\biggr)\biggr)\\
\end{eqnarray*}

\noindent The same method, word by word, leads to

\begin{eqnarray*}
\mathbb{E}e^{tX} \geq   -2t^2 (tc) \sigma^2 (e-7/4) \geq 0 
\end{eqnarray*}

\bigskip \noindent and next, by using the right inequality in Formula (EB1), we get

\begin{eqnarray*}
\mathbb{E}e^{tX}  \geq \biggr(1 + \frac{(t^2\sigma^2}{2}\left(1 - \frac{tc}{2} \right) \geq  \exp\biggr(\frac{(t^2\sigma^2}{2}\left(1 - 
tc\right) \biggr),
\end{eqnarray*}

\noindent which establishes Formula (EB3). \ $\square$\\

\bigskip \noindent Here the first result concerning the exponential bounds.\\

\begin{theorem} \label{expoboundth01} Let us use the same notations as in Lemma \ref{expoboundlem01}. Then the assertions below hold true, for any $\varepsilon>0$, for any $n\geq 1$.\\

\noindent (i) for $c\varepsilon\leq 1$,

$$
\mathbb{P}\biggr(S_n > \varepsilon s_n)\biggr) < \exp \biggr(- \frac{\varepsilon^2}{2}\left(1- \frac{\varepsilon c}{2} \right)\biggr).
$$

\noindent (ii) and for $c\varepsilon>1$,

$$
\mathbb{P}\biggr(S_n > \varepsilon s_n)\biggr) < \exp\biggr(-\frac{\varepsilon^2}{4c}\biggr).
$$
\end{theorem}

\bigskip \noindent \textbf{Proof of Theorem \ref{expoboundth01}}. To make the notation shorter, we put $S=S_n$ and $s_n=s$ and some times $S^{\ast}=S/s$. Now let us apply Formulas (EB2) and (EB3) in Lemma \ref{expoboundlem01} in the following way : for $t>0$ and $tc\leq 1$, and since 

$$
\mathbb{E}\exp(t S^{\ast})=\prod_{1\leq k \leq n}\mathbb{E}\exp(t X_k/s),
$$

\noindent we have
\begin{eqnarray*}
\prod_{1\leq k \leq n} \exp\biggr( \frac{(t^2\sigma_k^2}{2s^2}\left(1 -tc \right)\biggr) <\mathbb{E}\exp(t S^{\ast})
\prod_{1\leq k \leq n} \exp\biggr( \frac{(t^2\sigma_k^2}{2s^2}\left(1 + \frac{tc}{2} \right)\biggr) \ \ 
\end{eqnarray*}

\bigskip \noindent This obviously leads to

\begin{eqnarray*}
\exp\biggr( \frac{(t^2}{2}\left(1 -tc \right)\biggr) <\mathbb{E}\exp(t S^{\ast})\\
\exp\biggr( \frac{(t^2}{2}\left(1 + \frac{tc}{2} \right)\biggr). \ \ (DE)\\
\end{eqnarray*}

\noindent From this, we are able to handle both Points (i) or (ii).\\

\noindent For (i), we may apply the the Markov inequality and left-hand inequality in Formula (DE) above to $t>0$, $\varepsilon\>0$ such that 
$c\varepsilon\leq 1$ and  $tc\leq 1$, to get

\begin{eqnarray*}
\mathbb{P}\left(S^{\ast} > \varepsilon\right)&=&\mathbb{P}\left(\exp(t S^{\ast}) > \exp(t\varepsilon)\right)\\
&\leq & \exp(-t\varepsilon) \mathbb{E}\exp(t S^{\ast})\\
&\leq & \exp\biggr(-t\varepsilon \frac{(t^2}{2}\left(1 + \frac{tc}{2} \right)\biggr). \ (L23)\\
\end{eqnarray*}

\noindent We point out that the condition $tc\leq 1$ intervenes only in the conclusion in Line (L23). Taking $t=\varepsilon$ in 
in Line (L23) (which is possible since both conditions $tc\leq 1$ and $\varepsilon c\leq 1$ hold) leads to

\begin{eqnarray*}
\mathbb{P}\left(S^{\ast} > \varepsilon\right) &<& \exp \biggr(-\varepsilon^2 + \frac{\varepsilon^2}{2}\left(1 + \frac{\varepsilon c}{2} \right)\biggr)\\
&=& \exp \biggr(-\frac{\varepsilon^2}{2} \left(1- \frac{\varepsilon c}{2} \right) \biggr),
\end{eqnarray*}

\bigskip \noindent which is the announced result for Point (i).\\ 

\noindent To prove Point (ii), let $c\varepsilon>1$, we use the value $t=1/c$ (here again, the condition $tc\leq 1$ holds) to get

\begin{eqnarray*}
\mathbb{P}\left(S^{\ast} > \varepsilon \right) &<& \exp \biggr(-\frac{\varepsilon}{c} + \frac{1}{2c^2}\left(1 + \frac{1}{2} \right)\biggr)\\
&=& \exp \biggr(-\frac{\varepsilon}{4c}\biggr),
\end{eqnarray*}

\bigskip \noindent which is the announced result for Point (ii). $\square$\\

\noindent Finally, the coming exponential bound is very important when dealing with the \textit{Law of iterated logarithm} (LIL). We have :

\begin{theorem} \label{expoboundth02} Let us use the same notation as in Lemma \ref{expoboundlem01}. 

\bigskip \noindent Let us fix $0<\alpha<1/4$, we set $\beta=2\sqrt{\alpha}$ and 

$$
\gamma=\frac{1+2\alpha+\beta^2/2}{(1-\beta)^2}-1>0.
$$

\noindent  Then there exists $t(\alpha)$ large enough such that for $c(\alpha)$ small enough, that is $c(\alpha)<\alpha/t$ and 
$8c(\alpha)t(\alpha)\leq 1$ such that for $\varepsilon=t(\alpha)(1-2\sqrt(\alpha))$ we have
$$
\mathbb{P}\biggr(S_n > \varepsilon s_n)\biggr) < \exp\biggr(-\frac{\varepsilon^2}{4c}(1+\gamma)\biggr). \Diamond
$$
\end{theorem}

\bigskip \noindent \textbf{Proof of Theorem \ref{expoboundth02}}. The proof is so really technical that some authors like \cite{gutt} omitted and explained :  \textit{this one is no pleasure to prove} it. He referred to \cite{stout}.\\

\noindent Here, we will follow the lines of the proof in \cite{loeve}. However, the presentation and the ordering of the arguments have been significantly improved.\\

\noindent From Formula (DE), we may fix $0<\alpha=t_0 c<1$ so that for all $t \leq t_0$, we have

$$
\mathbb{E}\exp(t S^{\ast}) > \exp \biggr( \frac{t^2}{2}(1-\alpha) \biggr). \ (EB4)
$$

\bigskip \noindent  The principle of all the proof is to fix first $t>0$, as large as necessary, and to choose $c$ so that the desired conclusions hold. Then, let us choose $\alpha$ such that $2 \sqrt{\alpha}<1$. Put

$$
(i) \ \ \beta=2 \sqrt{\alpha},
$$

\bigskip \noindent  We have
$$
(ii) \ \ \gamma=\frac{1+2\alpha+\beta^2/2}{(1-\beta)^2}=\frac{(1+\beta)^2+2\beta + 1}{2(1-\beta)}-1>0.
$$

\bigskip \noindent  The positivity of $\gamma$ is clear since $0<\beta<1$ and $\gamma>0$. We first choose $0<\alpha<1/4$ which guarantees that $1+\beta<2$. Formula (EB4) shows that when $\alpha$ is fixed, the following conditions make sense : For $t$ large enough, we have

$$
(iii) \ \ \biggr(8t^2 \exp \biggr(-\frac{\alpha t^2}{4})\biggr)\biggr) <1/4, \ \ (ii)  \ \ \mathbb{E}e^{tS^{\ast}}>8,
$$

$$
(iv) \ \ \frac{1}{4} \mathbb{E}e^{tS^{\ast}} > 2  \ \ and \ \  (v) \ \ \frac{1}{4t^2} \exp \left(\frac{t^2}{2}\alpha\right)>1.
$$

\bigskip \noindent  We choose a value $t>0$ satisfying points (iii), (iv), (v). Next we suppose that

$$
(vi) \ \  c<\alpha/t \ \, \ (vii) \ \ 8tc\leq 1 \ \ and \ (viii) \ \ c\leq 4t/(1-\beta). 
$$ 

\bigskip \noindent  Once these conditions are set, we may proceed to the proof. First of all, Formula (EB4) is justified by condition (vi). Put $q(x)=\mathbb{P}(S^{\ast} > x)$. By Formula (CF), in Chapter \ref{proba_02_intotp}, page \pageref{formulaCF}, we have, for $t>$ and $Z=\exp(tS^{\ast})$,

$$
\mathbb{E}(Z)=\int_{0}^{+\infty} \mathbb{P}(Z>t) dt, \ (EB5)
$$

\bigskip \noindent  We get
\begin{eqnarray*}
\mathbb{E}(\exp(tS^{\ast})&=&\int_{0}^{+\infty} \mathbb{P}(\exp(tS^{\ast}>y) dy\\
&=&\int_{0}^{+\infty} \mathbb{P}\left(S^{\ast}>\frac{\log y}{t}\right) dy\\
&=&\int_{0}^{+\infty} q\left(\frac{\log y}{t}\right) dy\\
&=&t \int_{-\infty}^{+\infty} e^{tx} q(x) dx\  (L44)\\
\end{eqnarray*}

\noindent Now we split the integral in Line (L44) above by decomposing the integration domain  
$I=]-\infty, +\infty[$ using the intervals $I_1=]-\infty, 0]$, $I_2=]0, (1-\beta)]$, $I_3=]t(1-\beta), t(1+\beta)]$, 
$I_4=]t(1+\beta), 8t]$, and $I_5=]8t, +\infty[$, that is

$$
I=I_1+I_2+I_3+I_4+I_5.
$$

\bigskip \noindent  Let us name the integrals over $I_i$ by $J_i$, $i \in \{1,...,5\}$, respectively\\

\noindent Let us begin by $J_5$. Let $s \in I_5$. We have for $0\leq xc\geq 1$, by Formula (DE) and by Condition (vii)
$$
q(x) < \exp\left( -\frac{x}{4c}\right) < \exp(-2tx),
$$

\bigskip \noindent  where we use $1/(4c)=2t/(8tc)\geq 2t$. If $xc<1$, we apply again Formula (DE) and use $1-xc/2 \geq 1/2$ in the middle member, to get 
$$
q(x)< \exp \biggr(-\frac{x^2}{2} \left(1- \frac{xc}{2} \right) \biggr) < \exp\left( -\frac{x}{4c}\right) <\exp(-2tx)
$$

\bigskip \noindent  We get

$$
J_5=t \int_{8t}^{ta_n} e^{tx} q(x) dx \leq \int_{8t}^{+\infty} e^{tx} \exp(-2tx) dx < \leq \int_{0}^{+\infty} e^{-tx} dx=1.
$$

\bigskip \noindent  Since $q$ is a bounded by one, we have 

$$
J_1=J_5=t \int_{-\infty}^{0} e^{tx} q(x) dx=t \int_{0}^{+\infty} e^{-tx} q(x) dx\leq t \int_{0}^{+\infty} e^{-tx} dx=1. 
$$

\noindent  Now we handle $J_2$ and $J_4$ by using a maximization argument. On $I_4$ and $I_4$, we have $x\leq 0$ and $xc\leq 8tc<1$. From Point (ii) of Theorem \ref{expoboundth01}, and by using again $xc \leq 8tc$ in the second inequality below, we arrive at

$$
e^{tx} q(x)< \exp \biggr(tx-\frac{x^2}{2} \left(1- \frac{xc}{2} \right) \biggr)\leq \exp \biggr(tx-\frac{x^2}{2} \left(1- 4tc\right)\biggr)
\equiv g(x),
$$

\bigskip \noindent  where we remind that $4tc<1/4$. On $\mathbb{R}_{+}$, $g^{\prime}(x)=t-x(1-4tc)$ and thus $g$ attains its maximum at $x_0=t/(1-4tc)$.\\

\noindent Where lies $x_0?$ $x_0 >t(1-\beta)$ is equivalent to $-\beta/(1-\beta) < 4tc$ which is true. As well $x_0 \leq t(1+\beta)$ is equivalent to Condition (viii). Thus $x_0 \in J_3$. Hence on $I_2=]0, (1-\beta)]$, $g$ is non-decreasing and thus, for $x \in I_2$,
\begin{eqnarray*}
g(x)&=&g(t(1-\beta))=t^2 (1-\beta)-\frac{t^2(1-\beta)^2}{2} \left(1- 4tc\right)\\
&=&(1-\beta)\biggr(t^2 - \frac{t^2}{2}(1-\beta)+ \frac{t^2}{2}8tc\frac{1-\beta}{2} \biggr)\\
&=&\frac{t^2}{2}(1-\beta)\biggr(2 - (1-\beta)+ \frac{1-\beta}{2} \biggr) \ \ (we \ used \ Condition \ (iii)\\
&=&\frac{t^2}{2}(1-\beta)\biggr((1+\beta)+ \frac{(1-\beta)}{2} \biggr)\\
&=&\frac{t^2}{2}\biggr(1-\beta^2+ \frac{(1-\beta)^2}{2} \biggr)\\
&=&\frac{t^2}{2}\biggr((1-\beta^2/2)-\frac{1}{2} (1+2\beta)\biggr).\\
\end{eqnarray*}

\noindent  It follows that

\begin{eqnarray*}
J_2&=& t \int_{0}^{t(1-\beta)} e^{g(x)} dx\leq  t \int_{0}^{(1+\beta)} e^{g(t(1-\beta))} dx \\
&\leq& t^2(1-\beta) \exp \biggr((1-\beta^2/2)-\frac{1}{2} (1+2\beta)\biggr).
\end{eqnarray*}

\bigskip \noindent  As well on $I_4=]t(1+\beta), 8t]$ or on $I_4^{\ast}=]t(1+\beta), ta_n]$, $g$ is non-decreasing and we have for 
$x \in I_4^{\ast} \cup I_4$, 

\begin{eqnarray*}
g(x)&=&g(t(1+\beta))=t^2 (1+\beta)-\frac{t^2(1+\beta)^2}{2} \left(1- 4tc\right)\\
&=&(1+\beta)\biggr(t^2 - \frac{t^2}{2}(1+\beta)+ \frac{t^2}{2}8tc\frac{1+\beta}{2} \biggr)\\
&=&\frac{t^2}{2}(1+\beta)\biggr(2 - (1+\beta)+ \frac{1+\beta}{2} \biggr) \ \ (we \ used \ Condition \ (iii)\\
&=&\frac{t^2}{2}(1+\beta)\biggr((1-\beta)+ \frac{(1+\beta)}{2} \biggr)\\
&=&\frac{t^2}{2}\biggr(1-\beta^2+ \frac{(1+\beta)^2}{2} \biggr)\\
&=&\frac{t^2}{2}\biggr((1-\beta^2/2)-\frac{2\beta^2-2\beta+1}{2} \biggr),\\
\end{eqnarray*}

\noindent  since the polynomial $2\beta^2-2\beta+1$ has a negative discriminant and thus, is constantly positive. It follows that

\begin{eqnarray*}
\max(J_4,J_4^{\ast}) &\leq&  t \int_{t(1+\beta)}^{8t} e^{g(x)} dx\\
&\leq&  t \int_{0}^{(1+\beta)} e^{g(t(1+\beta))} dx\\
&\leq 7& t^2 \exp \biggr((1-\beta^2/2)-\frac{1}{2} (1+2\beta)\biggr).
\end{eqnarray*}

\bigskip \noindent  So, we have
\begin{eqnarray*}
\max(J_4,J^{\ast}_4) &<& 7t^2 \exp \biggr(\left(1-\beta^2/2\right)\biggr).
\end{eqnarray*}

\noindent Now we remind that $\alpha$ is fixed and $\alpha=\beta^2/4$ and hence

$$
1-\beta^2/2 = (1-\alpha) +\alpha/2.
$$

\bigskip \noindent  Hence Inequality (EB4) gives

$$
\mathbb{E}e^{tX} > \exp \biggr(\frac{t^2}{2}(1-\alpha)\biggr)= \exp \biggr(\frac{t^2}{2}(1-\beta^2/2))\biggr) 
\exp \biggr(\frac{\alpha t^2}{4})\biggr)
$$

\bigskip \noindent  and hence

$$
\exp \biggr(\frac{\alpha t^2}{4}(1-\beta^2/2))\biggr) < \mathbb{E}e^{tS^{\ast}} \biggr(\frac{\alpha t^2}{4}(1-\beta^2/2))\biggr),
$$

\bigskip \noindent  which, by using Conditions (iiia) and (iiib), leads to

$$
J_2+J_4 <\biggr(8t^2 \exp \biggr(\frac{\alpha t^2}{4})\biggr)\biggr)  \mathbb{E}e^{tS^{\ast}}
$$

\bigskip \noindent  and by Condition (iv), we get 
$$
J_1+J_5) <2 <\frac{1}{4} \mathbb{E}e^{tS^{\ast}} \ and \ J_2+J_4, J_2+J_4^{\ast}) <\frac{1}{4} \mathbb{E}e^{tS^{\ast}}.
$$

\bigskip \noindent  Since $\mathbb{E}e^{tS^{\ast}}=J_1+J_2+J_3+J_4+J_5$, it follows

$$
J_3=t \int_{t(1-\beta)}^{t(1+\beta)} e^{tx} q(x) dx >\frac{1}{2} \mathbb{E}e^{tS^{\ast}}.
$$

\noindent  Now using the bound of $\mathbb{E}e^{tS^{\ast}}$ as in Formula (S) and using the non-increasingness of $q$ and the non-decreasingness of $x \mapsto e^{tx}$ for $t>0$, leads to

$$
\frac{1}{2} \exp \biggr(\frac{t^2}{2}(1-\alpha)\biggr) < t q(t(1-\beta)) \int_{t(1-\beta)}^{t(1+\beta)} e^{t^2(1+\beta)}  dx,
$$ 

\bigskip \noindent  that is

$$
\frac{1}{2} \exp \biggr(\frac{t^2}{2}(1-\alpha)\biggr) < 2t^2 q(t(1-\beta))  e^{t^2(1+\beta)}
$$ 

\bigskip \noindent  and next

$$
q(t(1-\beta)) \geq \biggr[\frac{1}{4t^2} \exp \left(\frac{t^2}{2}\alpha\right)\biggr] \biggr[\exp \biggr(-t^2(1+\beta)+\frac{t^2}{2}(1-\alpha)-\frac{t^2}{2}\alpha\biggr)\biggr], \ (EB6)
$$

\bigskip \noindent with

\begin{eqnarray*}
&&\exp \left(\frac{t^2}{2}\alpha\right) \exp\biggr(-t^2(1+\beta)\frac{t^2}{2}(1-\alpha)-\frac{t^2}{2}\alpha\biggr)\\
&=& \exp \left(\frac{t^2}{2}(1-\alpha-\alpha -2(1+\beta) \right)\\
&=& \exp \left(\frac{t^2}{2}(1-2\alpha - 2(1+\beta) \right)\\
&=& \exp \left(-\frac{t^2}{2}(1+2\alpha +2\beta) \right) \ (EB7)\\
\end{eqnarray*}

\noindent We take $t=\varepsilon/(1-\beta)$. The quantity between the big brackets is bounded below by one in virtue of Condition 
(v).  From this the combination of (EB6) and (EB7) gives

$$
q(\varepsilon) > \exp \left(-\frac{\varepsilon^2}{2}\frac{1+2\alpha+2\beta}{(1-\beta^2)}\right)=
\exp \left(-\frac{\varepsilon^2}{2}\frac{1+2\beta+\beta^2/2}{(1-\beta^2)} \right)
$$

\bigskip \noindent and finally, by Condition (ii), we get

$$
q(\varepsilon) > \exp \left(-\frac{\varepsilon^2}{2} (1+\gamma) \right),
$$

\bigskip \noindent  which was the target. $\square$\\

\bigskip \noindent \textbf{(17) - Billingsley's Inequality (See \cite{billingsley}, page 69)}. \label{ineqBillingsley}\\

\noindent Let $(X_n)_{n\geq 0}$ be a sequence of square integrable and centered real-valued random variables defined on the same probability space $(\Omega, \mathcal{A}, \mathbb{P})$. We have for any $\varepsilon> \sqrt{2}$,

$$
\mathbb{P}\left(\max_{1\leq k \leq n} S_n\geq \varepsilon\right) \leq 2 \mathbb{P}\left(S_n \geq \varepsilon -\sqrt{2\mathbb{V}ar(S_n)}\right).
$$

\bigskip \noindent where, as usual, $S_n$, $n\geq 1$, are the partial sums of the studied sequence.\\

\noindent \textbf{Proof}. Put $s_k^2=\mathbb{V}ar(S_k)$ for $k\geq 1$. As usual,

\begin{eqnarray*}
A&=&\left(\max_{1\leq k \leq n} S_n\geq \varepsilon s_n\right)\\
&=&\sum_{1\leq j \leq n} \left(S_1<\varepsilon s_n, \cdots, S_{j-1}<\varepsilon s_n, S_{j}\geq \varepsilon s_n\right)\\
&\equiv& \sum_{1\leq j \leq n} A_j.\
\end{eqnarray*}

\noindent Now we have

\begin{eqnarray*}
\mathbb{P}(A)&=& \mathbb{P}(A \cap (S_n \geq (\varepsilon-\sqrt{2})s_n)+\mathbb{P}(A \cap (S_n< (\varepsilon-\sqrt{2})s_n)\\
&\leq &  \mathbb{P}(S_n\geq (\varepsilon-\sqrt{2})s_n)+\sum_{1\leq \j \leq n} \mathbb{P}(A_j \cap (S_n< (\varepsilon-\sqrt{2})s_n)\\
&=&\mathbb{P}(S_n\geq (\varepsilon-\sqrt{2})s_n)+\sum_{1\leq \j \leq n-1} \mathbb{P}(A_j \cap (S_n< (\varepsilon-\sqrt{2})s_n),
\end{eqnarray*}

\noindent since $\mathbb{P}(A_n \cap (S_n< (\varepsilon-\sqrt{2})s_n)=\emptyset$. We also have for each $1\leq j <n$, 

$$
(S_j \geq \varepsilon_n) \ and \ (S_n \leq \varepsilon-\sqrt{2})s_n) \Rightarrow (S_n-S_j \geq \sqrt{2} s_n) \Rightarrow (|S_n-S_j| \geq \sqrt{2} s_n.)
$$

\bigskip \noindent Since we still have that $S_n-S_j=X_{j+1}+\cdots+X_{j+1}$, for $1\leq j<n$, is independent of $A_j$, we get

$$
\mathbb{P}(A_j \cap (S_n< (\varepsilon-\sqrt{2})s_n)\leq \mathbb{P}(A_j) \mathbb{P}(|S_n-S_j| \geq s_n\sqrt{2}), \ 1\leq j<n.
$$

\bigskip \noindent Now using the Tchebychev inequality, we get

\begin{eqnarray*}
\mathbb{P}(A)&\leq&\mathbb{P}(S_n\geq (\varepsilon-\sqrt{2})s_n)+\sum_{1\leq \j \leq n-1} \frac{s_n^2-s_j^2}{2s_n^2})\\
&\leq&\mathbb{P}(S_n\geq (\varepsilon-\sqrt{2})s_n)+\sum_{1\leq \j \leq n-1} \frac{1}{2} \mathbb{P}(A_j)\\
&\leq&\mathbb{P}(S_n\geq (\varepsilon-\sqrt{2})s_n)+\frac{1}{2} \sum_{1\leq \j \leq n}  \mathbb{P}(A_j) \\
&=&\mathbb{P}(S_n\geq (\varepsilon-\sqrt{2})s_n)+\frac{1}{2} \mathbb{P}(A)\\
\end{eqnarray*}

\bigskip \noindent which leads to the desired result. $\square$\\

\noindent \textbf{18 - Etemadi's Inequality}. Let  $X_{1}, \cdots, \leq X_{n}$ be $n$ independent real-valued random variables such that the partial sums $S_{k}=X_{1}+...+X_{k}$, $1\leq k\leq n$, are definied. Then for any $\alpha \geq 0$, we have

\begin{equation*}
\mathbb{P}\left(\max_{1\leq k\leq n}\left\vert S_{k}\right\vert \geq 3\alpha \right)\leq 3 \max_{1\leq k\leq n}\mathbb{P}\left(\left\vert S_{k}\right\vert \geq \alpha\right). \ \Diamond
\end{equation*}

\noindent \textbf{proof}. The formula is obvious for $n=1$. Let $n\geq 2$. As usual, denote $B_1=(|X_1|\leq 3\alpha)$, $B_k=(|S_1|<3\alpha, \cdots , |S_{k-1}|<3\alpha, |S_k|\geq \geq 3\alpha)$, $k\geq 2$. By decomposing $(\max_{1\leq j \leq n} |S_j| \geq 3\alpha)$ over the partition 
$$
(|S_n| \geq \alpha)+(|S_n| < \alpha)=\Omega,
$$

\noindent we have

$$
(\max_{1\leq j \leq n} |S_j| \geq 3\alpha) \subset (|S_n| \geq \alpha) \cup (|S_n|<\alpha, \max_{1\leq j \leq n} |S_j| \geq 3\alpha)
$$

\noindent And by the principle of the construction of the $B_j$,

$$
(\max_{1\leq j \leq n} |S_j| \geq 3\alpha)=\sum_{1\leq j \ n} B_j
$$ 

\noindent and hence

$$
(\max_{1\leq j \leq n} |S_j| \geq 3\alpha) \subset (|S_n| \geq \alpha) \cup \sum_{1\leq j \ n-1} (|S_n| <\alpha 3\alpha) \cap B_j
$$

\noindent where the summation is restricted to $j \in \{1,...,n-1\}$ since the event $(|S_n| <\alpha 3\alpha) \cap B_n$ is empty. Further, on $(|S_n| < \alpha) \cup B_j$, we have $(|S_n| < \alpha)$ and $(|S_j| < 3\alpha)$ and the second triangle inequality $|S_n-S_j|\geq |S_j|-|S_n|\geq 3\alpha-\alpha=2\alpha$, that is

$$
(|S_n| < \alpha) \cap B_j \subset B_j \cap (|S_n|<\alpha) \cap (|S_n-S_j|\geq 2\alpha) \subset B_j \cap (|S_n-S_j|\geq 2\alpha).
$$

\noindent Now, we remind that and $B_j$ and $S_n-S_j$ are independent. Translating all this into probabilities gives

\begin{eqnarray*}
\mathbb{P}(\max_{1\leq j \leq n} |S_j| \geq 3\alpha)&\leq&  \mathbb{P}(|S_n| \geq \alpha) + \sum_{1\leq j \ n-1} \mathbb{P}(B_j) \mathbb{P}(|S_n-S_j|\geq 2\alpha)\\
&\leq&  \mathbb{P}(|S_n| \geq \alpha) + \sum_{1\leq j \ n-1} \mathbb{P}(B_j) \biggr(\mathbb{P}(|S_n|\geq \alpha)+\mathbb{P}(|S_j|\geq \alpha)\biggr).\\
\end{eqnarray*} 

\noindent But $(|S_n|\geq 2\alpha)$, $(|S_n| \geq \alpha)$ and $(|S_j|\geq 2\alpha)$ are subsets of 
$$
(\max_{1\leq j \leq n} |S_j| \geq \alpha)
$$ 

\noindent and hence, we may conclude that

\begin{eqnarray*}
\mathbb{P}\left(\max_{1\leq j \leq n} |S_j| \geq 3\alpha\right)&\leq&  \mathbb{P}\left(\max_{1\leq j \leq n} |S_j| \geq \alpha\right) \left(1+ 2 \sum_{1\leq j \leq n}\mathbb{P}(B_j)\right)\\
&\leq&  \mathbb{P}\left(\max_{1\leq j \leq n} |S_j| \geq 3\alpha\right) \left(1+ 2 \mathbb{P}\left(\sum_{1\leq j n}B_j\right)\right)\\
&\leq&  3\mathbb{P}(\max_{1\leq j \leq n} |S_j| \geq 3\alpha). \ \square
\end{eqnarray*}

%% file: proba_02_07_ang.tex
\chapter[Asymptotics Theorems of Independent Random variables]{Introduction to Classical Asymptotic Theorems of Independent Random variables} \label{probab_02_indep}

\section{Easy Introduction} \label{probab_02_indep_sec_01} 

\noindent We are going to quickly discover three classical types of well-known convergences which are related to sequences of independent random variables. In the sequel :\\

\noindent $(X_n)_{n\geq 0}$ is a sequence of centered real-valued random variables defined on the same probability space $(\Omega, \mathcal{A}, \mathbb{P})$. If the expectations $\mu_n=\mathbb{E}X_n$'s exist, we usually center the $X_n$'s at their expectations by taking $X_n-\mu_n$ in order to have centered random variables. If the variances exist, we denote $\sigma_n^2=\mathbb{V}ar(X_n)$ and

$$
s_0^2=0, \ s_1^2=\sigma_1^2, \ s_n^2=\sigma_1^2+...+\sigma_n^2, \ n\geq 2.
$$

\bigskip \noindent The laws we will deal with in this chapter are related to the partial sums

$$
S_0=0, \ S_n=X_1, \ S_n=X_1+...+X_n.
$$

\bigskip \bigskip \noindent \textbf{(a) Discovering the simplest Weak Law of Large Numbers (WLLN)}.\\ 

\noindent Suppose that the random variables $X_n$ are independent and are identically distributed (\textit{iid}) and have the common mathematical expectation $\mu$. We are going to find the limit in probability of the sequence

$$
\overline{X}_n=\frac{S_n}{n}, \ \neq 0.
$$

\bigskip \noindent By Proposition \ref{proba_02_conv_sec_06_prop_01} in Section \ref{proba_02_conv_sec_06} in Chapter \ref{proba_02_conv}, a non-random weak limit is also a limit in probability and vice-versa. So we may directly try to show that $\overline{X}_n$ converges to a non-random limit (which is supposed to be $\mu$). To do this, we have many choices through the Portmanteau Theorem
\ref{proba_02_conv_sec_04_thportmanteau} in Section \ref{proba_02_conv_sec_04} in Chapter \ref{proba_02_conv}. Let us use the characteristic function tool $\Phi_{X_j}=\Phi$ for all $j\geq 1$. Since we have, by Proposition \ref{proba_02_rv_sec_06_propMoments} in Section \ref{proba_02_rv_sec_06} in Chapter \ref{proba_02_rv}, \\

$$
\Phi^{\prime}(0)=i \mu \ \ and \ \ \Phi(0)=1, 
$$

\bigskip \noindent (where $i$ is the normed pure complex number with a positive angle), we may use a one order Taylor expansion of $\Phi$ at zero to have

$$
\Phi(u)=1 + i \mu u + O(u^2), \ as \ u\rightarrow 0. \ (EX)
$$

\bigskip \noindent By the properties of the characteristic function and by taking into account the fact that the variables are \textit{iid}, we have

$$
\Phi_{S_n/n}(u)=\Phi_{X_1+...+X_n}(u/n)=\Phi(u/n)^n, \ u \in \mathbb{R}.
$$

\noindent Now, for $u$ fixed, we have $u/n\rightarrow 0$ as $n\rightarrow \infty$, and we may apply Formula (EX) to have, as $n\rightarrow +\infty$,

$$
\Phi_{S_n/n}(u)=\exp\biggr(n \log(1 + i \mu u/n + O(n^{-2})) \biggr) \rightarrow \exp(i \mu u)=\Phi_{\mu}(u).
$$

\bigskip \noindent Here, we skipped the computations that lead to $n \log(1 + i \mu u/n + O(n^{-2}) \rightarrow i \mu u$. In previous books as\cite{ips-probelem-ang} and \cite{ips-wcrv-ang}, such techniques based on expansions of the logarithm function have been given in details.\\

\noindent We just show that $S_n/n \rightsquigarrow \mu$, hence $S_n/n \rightarrow_{\mathbb{P}} \mu$. This gives us the first law.\\

\begin{theorem} (Kintchine) \label{probab_02_indep_sec_01_theoWLLN} If $(X_n)_{n\geq 0}$ is a sequence of independent and are identically distributed (\textit{iid}) random variables with a finite common mathematical expectation $\mu$, we have the following Weak Law of Large Numbers (WLLN) :

$$
S_n/n \rightarrow_{\mathbb{P}} \mu, \ as \ n \rightarrow +\infty.
$$
\end{theorem}

\bigskip \noindent \textbf{(b) Discovering the Strong Law of Large Numbers (SLLN)}.\\ 

\noindent Before we proceed further, let us state a result of measure theory and integration (See \cite{ips-mestuto-ang}, given in Exercise 3 in Doc 04-05, and its solution in Doc 04-08) in Chapter 5 in \cite{ips-mestuto-ang} in the following famous lemma.

\begin{lemma}  (\textbf{Borel-Cantelli Lemma})  Let $(A_n)_{n\geq 0} \subset \mathcal{A}$.\\

\noindent (i) If the series $\sum_{n\geq 0} \mathbb{P}(A_n)<+\infty$ is convergent, then

$$
\mathbb{P}\left(\limsup_{n \rightarrow +\infty} A_n\right)=0.
$$ 
\label{probab_02_indep_sec_01_lem01}
\bigskip \noindent (ii) If the events $A_n$ are independent and if the series diverges, that is $\sum_{n\geq 0} \mathbb{P}(A_n)=\infty$, then

$$
\mathbb{P}\left(\limsup_{n \rightarrow +\infty} A_n\right)=1.
$$ 
\end{lemma}

\bigskip 

\noindent This lemma is the classical basis of the simple \textit{SLLN}. But before we continue, let us give the following consequence.\\

\begin{corollary} \label{probab_02_indep_sec_01_cor01}  Let $(X_n)_{n\geq 0}$ be a sequence of independent \textit{a.e.} finite real-valued random variables such that $X_n\rightarrow 0$ \textit{a.s} as $n\rightarrow +\infty$. Then for any finite real number $c>0$,

$$
\sum_{n\geq 0} \mathbb{P}(|X_n|\geq c) <+\infty. 
$$ 
\end{corollary}

\bigskip \noindent \textbf{Proof}. Given the assumptions of the corollary, the events $A_n$'s are independent. By the Borel-Cantelli Lemma, $\sum_{n\geq 0} \mathbb{P}(|X_n|\leq c)=+\infty$ would imply $\mathbb{P}(|X_n|>c, i.o)=1$ and hence $(X_n\rightarrow 0) \ \ a.e.$ would be false. The proof is complete with this last remark. $\square$\\

\noindent Let us expose the simple strong law of large number.\\

\begin{theorem} (Simple Strong Law of Large Numbers) \label{probab_02_indep_sec_01_theoSLLN}  Let $(X_n)_{n\geq 0}$ be a sequence of independent centered and square integrable random variables with variance one, that is $\mathbb{E}X_n^2=1$ for all $n\geq 1$. Then
$$
\frac{1}{n} \sum_{1\leq k \leq n} X_k \rightarrow 0 \ a.s. \ as \ n\rightarrow +\infty.
$$
\end{theorem}

\noindent * \textbf{Proof}. Suppose that the assumption of the theorem hold. We are going to use the perfect square method. Put

\begin{equation*}
Y_{n}=S_{n^{2}}/n^{2}, \ n\geq 1,
\end{equation*}

\bigskip \noindent that is,  we only consider the elements of the sequence $(S_k/k)_{k\geq 1}$ corresponding to a square index $k=n^2$. Remark that $\mathbb{V}ar(Y_{n})=n^{-2}$, $n\geq 2$. Fix $0<\beta <1/2$. By Chebychev's inequality, we have
 
\begin{equation*}
\mathbb{P}(\left\vert Y_{n}\right\vert >n^{-\beta })\leq n^{2(1-\beta )}
\end{equation*}

\bigskip \noindent and thus,

\begin{equation*}
\sum_{n}\mathbb{P}(\left\vert Y_{n}\right\vert >n^{-\beta })\leq \sum n^{2(1-\beta )}<\infty.
\end{equation*}

\bigskip \noindent By Borel-Cantelli's Lemma, we conclude that

\begin{equation*}
\mathbb{P}(\lim_{n}\inf (\left\vert Y_{n}\right\vert \leq n^{-\beta })=1.
\end{equation*}%

\bigskip \noindent Let us remind that  

\begin{equation*}
\Omega_0=\lim_{n}\inf (\left\vert Y_{n}\right\vert \leq n^{-\beta
})=\bigcup_{n\geq 0}\bigcap_{r\geq n}(\left\vert Y_{r}\right\vert \leq r^{-\beta }).
\end{equation*}

\noindent Hence, for all $\omega \in \Omega_0$, there exists  $n(\omega)\geq 0$ such that for any $r \geq n$, 

\begin{equation*}
\left\vert Y_{r}\right\vert \leq r^{-\beta}.
\end{equation*}

\noindent * By the sandwich's rule, we conclude that, for any $\omega \in \Omega_0$, we have

\begin{equation*}
Y_{m}(\omega )\rightarrow 0.
\end{equation*}
 
\bigskip \noindent This means that

\begin{equation*}
\Omega_0 \subset (Y_{n}\rightarrow 0).
\end{equation*}

\bigskip \noindent We conclude that $\mathbb{P}(Y_{n}\rightarrow 0)=1$ and hence $Y_n \rightarrow 0$, \textsl{a.s.}.\\

\noindent To extend this result to the whole sequence, we use the decomposition of $\mathbb{N}$ by segments with perfect squares bounds. We have 

\begin{equation*}
\forall (n\geq 0),\exists m\geq 0,k(n)=m^{2}\leq n\leq (\sqrt{k(n)}+1)^{2}.
\end{equation*}

\bigskip \noindent We have  

\begin{equation*}
\mathbb{E}(\frac{1}{n}(S_{n}-S_{k(n)})=0
\end{equation*}

\bigskip \noindent and 

\begin{equation*}
\mathbb{V}ar(\frac{1}{n}(S_{n}-S_{k(n)})=\frac{1}{n^{2}}E\sum_{i=k(n)+1}^{n}X_{i}^{2}%
\leq \frac{1}{n^{2}}(2\sqrt{k(n)}+1)\leq \frac{3\sqrt{n}}{2}=3n^{-3/2}.
\end{equation*}

\bigskip \noindent Hence,
 
\begin{equation*}
\sum_{n}\mathbb{P}(\left\vert \frac{1}{n}(S_{n}-S_{k(n)})\right\vert
>n^{-\beta })\leq 3\sum n^{-(\frac{3}{2}-2\beta )}<\infty
\end{equation*}

\bigskip 
\noindent whenever $\beta <3/4$. We conclude as previously that 

\begin{equation*}
\frac{1}{n}(S_{n}-S_{k(n)})\rightarrow 0,\text{ }a.s.
\end{equation*}

\bigskip \noindent Finally we have 

\begin{equation*}
\frac{S_{n}}{n}=\frac{S_{n}-S_{k(n)}}{n}+\frac{k(n)}{k(n)}\times \frac{S_{k(n)}%
}{n}\rightarrow 0\text{ }a.s.,
\end{equation*}

\bigskip \noindent since
 
\begin{equation*}
1\leq \frac{n}{k(n)}<1+\frac{2}{\sqrt{k(n)}}+\frac{1}{k(n)}
\end{equation*}

\bigskip \noindent and 
\begin{equation*}
\frac{k(n)}{n}\rightarrow 1.
\end{equation*}

\bigskip \noindent We just finished to prove that  
\begin{equation*}
\frac{S_{n}}{n}\rightarrow 0\text{ }a.s.\  \blacksquare
\end{equation*}
 
\bigskip \noindent In a more general case of random variables with common variance, we may center and normalize them to be able to use the result above as in 

\begin{corollary} \label{probab_02_indep_sec_01_corSLLN}  Let $(X_n)_{n\geq 0}$ be a sequence of independent and square integrable random variables with equal variance $\sigma^2>0$,that is 
$\mathbb{V}ar(X_n)=\sigma^2$ for all $n\geq 1$. Then
$$
\frac{1}{n\sigma} \sum_{1\leq k \leq n} (X_k-\mathbb{E}(X_k)) \rightarrow 0 \ a.s. \ as \ n\rightarrow +\infty.
$$
\end{corollary}

\noindent * We may also derive the

\begin{proposition} (Kolmogorov) \label{probab_02_indep_sec_01_prop02} If $(X_n)_{n\geq 0}$ is a sequence of independent random variables with mathematical expectations $\mu_n$ and variances $0<\sigma_n^2<+\infty$, we have

$$
\frac{1}{n} \sum_{1\leq j \leq n} \frac{X_j-\mu_j}{\sigma_j} \rightarrow 0, \ a.s. \ as \ \rightarrow 0.
$$

\bigskip \noindent If the expectations are zero's that is $\mu_n=0$, $n\geq 0$ and if the variances are equal, that $\sigma_n^2=\sigma^2$, $n\geq 0$, we have the simple \textit{SLLN} :

$$
\frac{S_n}{n\sigma} \rightarrow 0, \ a.s. \ as \ \rightarrow 0.
$$
\end{proposition}

\bigskip \noindent \textbf{(c) Discovering the Central limit Theorem}.\\ 

\noindent The Central Limit Theorem in Probability Theory turns around finding conditions under which the sequence of partials sums $S_n$, $n\geq 1$, when appropriately centered and normalized, weakly converges to some random variable. Generally, the probability law of the limiting random variable is Gaussian.\\

\noindent Actually, we already encountered the \textit{CTL} in our series, through Theorem 4 in Chapter 7 in \cite{ips-probelem-ang} in the following way.\\

\noindent If the $X_n$'s are \textit{iid} according to a Bernoulli probability law $\mathcal{B}(p)$, $0<p<1$, $S_n$ follows a Biniomial laws of parameters $p=1-q$ and $n\geq 1$ and we have

$$
Z_n=\frac{S_n-\mathbb{E}S_n}{\mathbb{V}ar(S_n)^{1/2}}=\frac{S_n-npq}{\sqrt{npq}}, \ n\geq 1.
$$

\bigskip \noindent The invoked theorem (in \cite{ips-probelem-ang}) states that, as $n\rightarrow +\infty$,

$$
\forall x \in \mathbb{R}, \ F_{S_n}(x) \rightarrow N(x)= \frac{1}{\sqrt{2\pi}} \int_{-\infty}^{x} \exp(-t^2/2) \ d\lambda(t).
$$

\bigskip \noindent At the light of the Portmanteau Theorem \ref{proba_02_conv_sec_04_thportmanteau} in Section \ref{proba_02_conv_sec_04} in Chapter \ref{proba_02_conv}, we have

\begin{proposition} \label{probab_02_indep_sec_01_prop03} Let $(X_n)_{n\geq 1}$, be a sequence of independent random variables identically distributed as a Bernoulli probability law $\mathcal{B}(p)$, $0<p<1$. Then we have the following Central limit Theorem (CLT)
$$
\frac{S_n-npq}{\sqrt{npq}} \rightsquigarrow \mathcal{N}(0,1), \ as \ n\rightarrow +\infty.
$$
\end{proposition}

\bigskip
\noindent We are going to see that result is a particular case the following one.

\begin{proposition} \label{probab_02_indep_sec_01_prop04} (CLT for an iid sequence with finite variance). Let $(X_n)_{n\geq 1}$, be a sequence of centered and \textit{iid} random variables with common finite variance $\sigma^2>0$. Then, we have the following \textit{CLT}

$$
\frac{S_n}{s_n} \rightsquigarrow \mathcal{N}(0,1), \ as \ n\rightarrow +\infty. \ \ (CLTG)
$$

\bigskip \noindent If the common expectation is $\mu$, we may write

$$
\frac{S_n-n\mu}{\sigma \sqrt{n}} \rightsquigarrow \mathcal{N}(0,1), \ as \ n\rightarrow +\infty.
$$
\end{proposition}

\bigskip \noindent \textbf{Proof}.  The Portmanteau Theorem \ref{proba_02_conv_sec_04_thportmanteau} in Section \ref{proba_02_conv_sec_04} in Chapter \ref{proba_02_conv} offers us a wide set of tools for establishing weak laws. In on dimensional problems, the characteristic method is the favored one. Here, we have $\Phi_{X_j}=\Phi$ for all $i\geq 1$. Let us give the proof for $\sigma=1$. By Proposition \ref{proba_02_rv_sec_06_propMoments} in Section \ref{proba_02_rv_sec_06} in Chapter \ref{proba_02_rv}, we have

$$
\Phi(0)=1, \ \ \Phi^{\prime}(0)=0 \ \ and  \ \ \Phi^{\prime\prime}(0)=-1. 
$$

\bigskip \noindent Let us use two-order Taylor expansion of $\Phi$ in the neighborhood of $0$ to have :

$$
\Phi(u)=1 - u^2/2 + O(u^2), \ as \ u\rightarrow 0. \ (EX2)
$$

\bigskip \noindent By the properties of the characteristic function and by taking into account that the variables are iid, we have

$$
\Phi_{S_n/\sqrt{n}}(u)=\Phi_{X_1+...+X_n}(u/\sqrt{n})=\Phi(u/\sqrt{n})^n, \ u \in \mathbb{R}.
$$

\bigskip \noindent Now for $u$ fixed, we have $u/\sqrt{n}\rightarrow 0$ as $n\rightarrow \infty$, and we may apply Formula (EX2) to have, for $n\rightarrow +\infty$,

$$
\Phi_{S_n/\sqrt{n}}(u)=\exp\biggr(n \log(1 - u^2/(2n) + O(n^{-3/2})) \biggr) \rightarrow \exp(-u^2/2),
$$

\bigskip \noindent where again we skipped details on the expansions of the logarithm function. So we have just proved that

$$
S_n/\sqrt{n} \rightsquigarrow \mathcal{N}(0,1).
$$

\bigskip \noindent If the common expectation is $\mu$, we may transform the sequence to $\left((X_n-\mu)/\sigma \right)_{n\geq 1}$, which is an \textit{iid} sequence of centered random variables with variance one. By applying the result above, we get

$$
\frac{S_n}{\sigma \sqrt{n}} \rightsquigarrow \mathcal{N}(0,1).
$$

\bigskip \noindent We finish the proof by noticing that : $s_n^2=n\sigma^2$, $n\geq 1$. $\square$\\

\bigskip \noindent \textbf{(d) A remark leading the Berry-Essen Bounds}.\\

\noindent Once we have a \textsl{CLT} in the form of Formula (CLTG), the Portmanteau theorem implies that for any fixed $x\in \mathbb{R}$

$$
\biggr| \mathbb{P}\left(\frac{S_n}{s_n} \leq x\right)-N(x)\biggr| \rightarrow 0 \ as \ n \rightarrow +\infty.
$$

\bigskip \noindent Actually, the formula above holds uniformly (See Fact 4 in Chapter 4 in \cite{ips-wcrv-ang}), that is 

$$
B_n=\sup_{x \in \mathbb{R}}\biggr| \mathbb{P}\left(\frac{S_n}{s_n} \leq x\right)-N(x)\biggr| \rightarrow 0 \ as \ n \rightarrow +\infty.
$$

\bigskip \noindent A Berry-Bound is any bound of $B_n$. We will see later in this chapter a Berry-Essen bound for sequence of independent random variables with third finite moments.\\

\bigskip \noindent \textbf{Conclusion}.\\

\noindent Through Theorem \ref{probab_02_indep_sec_01_theoWLLN} and Propositions \ref{probab_02_indep_sec_01_prop02} and \ref{probab_02_indep_sec_01_prop04}, we discovered simple forms of three of the most important asymptotic laws in Probability Theory.\\

\noindent Establishing \textit{WLLN}'s, \textit{SLLN}'s, \textit{CLT}'s, Berrey-Essen bounds, etc. is still a wide and important part in Probability Theory research under a variety of dependence type and in abstract spaces.\\

\noindent For example, the extensions of such results to set-valued random variables constitute an active research field.\\

\noindent The results in this section are meaningful and are indeed applied. But we will give important more general cases in next sections. The coming results represent advanced forms for sequence of independent random variables.\\

%% file: proba_02_07_kol_ang.tex
\section[Kolmogorov's Strong Laws of Large Numbers]{Tail events and Kolmogorov's zero-one law and strong laws of Large Numbers} \label{probab_02_sllnKol}

\noindent This chapter will be an opportunity to revise generated $\sigma$-algebras and to deepen our knowledge on independence.\\

\noindent \textbf{(A) Introduction and statement of the zero-one law}.\\

\noindent At the beginning, let $X=(X_t)_{t\in T}$ be an non-empty of mappings from $(\Omega, \mathcal{A})$ to some measure spaces $(F_t,\mathcal{F}_t)$. The $\sigma$-algebra on $\Omega$ generated by this family is

$$
\mathcal{A}_{X}=\sigma\{X_{t_i}^{-1}(B_{t_i}), \ B_{t_i} \in \mathcal{F}_{s_i}, (t_1, ..., t_p) \in T^p, \ p\geq 1\}.
$$ 

\bigskip \noindent It is left as an exercise to check that $\mathcal{A}_{X}$ is also generated by the class of finite intersections of the form

$$
\mathcal{C}_{X}=\{\bigcap_{1\leq k \leq n} X_{t_i}^{-1}(B_{t_i}), \ B_{t_i} \in \mathcal{F}_{s_i}, (t_1, ..., t_p) \in T^p, \ p\geq 1\}, \ \ (P01)
$$

\bigskip \noindent which is a $\pi$-system.\\

\noindent Coming to our topic on the zero-one law, we already saw from the Borel-Cantelli Lemma \ref{probab_02_indep_sec_01_lem01} that : for a sequence of independent events $(A_n)_{n\geq 0}\subset \mathcal{A}$ on the probability space $(\Omega, \mathcal{A},\mathbb{P})$ such that $A_n \rightarrow A$ as $n\rightarrow +\infty$, then $\mathbb{P}(A) \in \{0,1\}$.\\

\noindent We are going to see that this is a more general law called the Kolmogorov zero-one law. Let $(X_n)_{n\geq 1}$ be a sequence of measurable mappings from $(\Omega, \mathcal{A})$ to some measure spaces $(F_i,\mathcal{F}_i)$, $i\geq 1$. For each $n\geq 0$, the smallest $\sigma$-algebra on $\Omega$ rendering measurable all the mapping $X_k$, $k\geq n$, with respect to $\mathcal{F}$ is

$$
\mathcal{A}_{tail,n}=\sigma\{X_k^{-1}(B_k), \ B_k \in \mathcal{F}_k, \ k\geq n\}.
$$ 

\bigskip \noindent It is usually denoted as $\mathcal{A}_{tail,n}=\sigma(X_k,\ k \geq n)$ and quoted as the $\sigma$-algebra generated by the mappings $X_k$, $k \geq n$.\\

\noindent 

\noindent \textbf{Definition}. The tail $\sigma$-algebra generated by the sequence $(X_n)_{n\geq 0}$, relatively to $\mathcal{F}$, is the intersection 

$$
\mathcal{A}_{tail}= \bigcap_{n\geq 0} \mathcal{A}_{tail,n}.
$$

\noindent The elements of $\mathcal{A}_{tail}$ are, by definition, the \textit{tail events} with respect to the sequence $(X_n)_{n\geq 0}$. $\Diamond$\\

 \noindent Let us give an example. Let $B_k\in \mathcal{F}_k$, for $k\geq 1$. We have that

$$
\liminf_{n\rightarrow +\infty} (X_n \in B_k) \in \mathcal{A}_{tail} \text{ and } \limsup_{n\rightarrow +\infty} (X_n \in B_k) \in \mathcal{A}_{tail}.
$$  

\bigskip \noindent Here is why. Because of the increasingness of the $(\cap_{p\geq n} A_p)_{n\geq 0}$ (in $n$), we have for any fixed $n_0\geq 0$,

$$
\liminf_{n\rightarrow +\infty} (X_n \in B)=\bigcup_{n\geq 1} \bigcap_{p\geq n }(X_p \in B) =\bigcup_{n\geq n_0} \bigcap_{p\geq n }(X_p \in B)
$$

\bigskip \noindent and 

$$
\left\{\bigcup_{k\geq n} (X_n \in B), \ n\geq n_0 \right\} \subset \mathcal{A}_{tail, n_0},
$$

\bigskip \noindent and then $\liminf_{n\rightarrow +\infty} (X_n \in B) \in \mathcal{A}_{tail, n_0}$ for all $n_0\geq 0$ and is in $\mathcal{A}_{tail}$. To get the same conclusion for the superior limit, we applied that conclusion to its complement.\\

\noindent Let us prove a useful result before we proceed further.\\

\noindent \textbf{The zero-one Law}. If the sequence elements of the sequence $(A_n)_{n\geq 0}$ are mutually independent, then any tail $A$ event with respect to that sequence is such that 
$\mathbb{P}(A) \in \{0,1\}$, that is the tail $\sigma$-algebra is behaves as the trivial $\sigma$-algebra. $\Diamond$\\

\bigskip \noindent Before we give the proof, let us get more acquainted with independent $\sigma$-algebras.\\

\noindent \textbf{(B) Independence of $\sigma$-algebras}.\\

\noindent \textbf{Definition}. Two non-empty sub-classes $\mathcal{C}_1$ and $\mathcal{C}_2$ of $\mathcal{A}$ are \textit{mutually independent} if and only if : for any subsets
$\{A_1, ..., A_{\ell_1}\} \subset \mathcal{C}_1$ and $\{B_1, ..., B_{\ell_2}\} \subset \mathcal{C}_2$, for any non-negative, real-valued and measurable functions $h_i$, $i\in \{1,2\}$, defined on a domain containing $\{0,1\}$, $h_1(1_{A_1},...,1_{A_{\ell_1}})$ and $h_2(1_{B_1},...,1_{B_{\ell_2}})$ are independent. $\Diamond$\\

\noindent For easy notation, let us denote by $\mathcal{I}(\mathcal{C}_1)$ the class of all elements of the form $h_1(1_{A_1},...,1_{A_{\ell_1}})$ as described above. $\mathcal{I}(\mathcal{C}_2)$ is
defined similarly.\\

 \noindent \textbf{Example}. Let $\mathcal{F}_1=\{X_t, t \in T\}$ a non-empty family of measurable mappings from $(\Omega, \mathcal{A})$ to some measure space $(F,\mathcal{F})$ and $\mathcal{F}_2=\{Y_s, s \in S\}$ a non-empty family of measurable mappings from $(\Omega, \mathcal{A})$ to some measure space $(G,\mathcal{G})$. Suppose that any finite pairs of sub-families $(X_{t_j})_{1\leq j \leq p}$ ($p\geq 1$)and $(Y_{s_j})_{1\leq j \leq q}$ ($q\geq 1$) of $\mathcal{F}_1$ and $\mathcal{F}_2$ respectively, the random vectors $(X_{t_1}, ..., X_{t_p})^t$ and $(Y_{s_1}, ..., Y_{s_q})^t$ are independent, that is

$$
\mathbb{P}_{(X_{t_1}, ..., X_{t_p},Y_{s_1}, ..., Y_{s_q})}=\mathbb{P}_{(X_{t_1}, ..., X_{t_p})} \otimes \mathbb{P}_{(Y_{s_1}, ..., Y_{s_q})}.  \ \ (DE01)
$$

\bigskip \noindent The classes

$$
\mathcal{C}_{X}=\{X_t^{-1}(B), \ B \in \mathcal{F}, \ t \in T\} \text{ and } \mathcal{C}_{Y}=\{Y_s^{-1}(C), \ C \in \mathcal{G}, \ s \in S\}
$$ 

\bigskip \noindent are independent. To see that, we consider two finite subsets of $\mathcal{C}_X$ and $\mathcal{C}_Y$ of the forms

$$
(X_{t_1}^{-1}(B_1), ..., X_{t_p}^{-1}(B_p)) \text{ and } (Y_{s_1}^{-1}(C_1), ..., S_{t_q}^{-1}(C_q)).
$$

\bigskip \noindent where the $(B_{j})_{1\leq j \leq p} \subset \mathcal{F}$ and $(C_{j})_{1\leq j \leq q} \subset \mathcal{G}$ and, accordingly, two real-valued and measurable functions $h_1$ and $h_1$  of their indicators functions as

$$
H_1=h_1(1_{B_1}(X_{t_1}), ..., 1_{B_p}(X_{t_p}))) \text{ and } H_2=h_1(1_{C_1}(Y_{s_1}), ..., 1_{C_q}(Y_{s_p}))).
$$

\bigskip \noindent So, the functions $H_1$ and $H_2$ are independent because of Formula (DE). $\Diamond$\\

\noindent For now, we need the two results in the next proposition.\\

\begin{proposition} \label{proba_02_slln_01} Let $\mathcal{C}_1$ and $\mathcal{C}_2$  be two mutually independent and non-empty $\pi$-sub-classes of $\mathcal{A}$. Consider the generated $\sigma$-algebra $\mathcal{A}_i=\sigma(\mathcal{C}_i)$, $i\in \{1,2\}$.\\

\noindent (1) Then for any $(A,B) \in \mathcal{A}_1 \times \mathcal{A}_2$, $A$ and $B$ are independent.\\

\noindent (2) $\mathcal{A}_1$ and $\mathcal{A}_2$ are independent.\\

\noindent (3) For any non-negative and real-valued function $Z_i$, $i\in \{1,2\}$, such that each $Z_i$ is $\mathcal{C}_i$-measurable, we have

$$
\mathbb{E}(Z_1 Z_2) =\mathbb{E}(Z_1) \mathbb{E}(Z_2). \ \Diamond
$$
\end{proposition}

\bigskip \noindent \textbf{proof}. We easily see that each $\mathcal{A}_i$, $i\in \{1,2\}$, is also generated by the class of finite intersections of sets which are either elements of $\mathcal{C}_i$ or complements of elements of $\mathcal{A}_i$, denoted

$$
\tilde{\mathcal{C}}_1=\{\bigcap_{1\leq k \leq n} A_i, \ A_i \in \mathcal{C}_{i} \text{ or } A_i^c \in \mathcal{C}_{i}, \ 1 \leq i \leq p, \ p\geq 1\}.
$$

\bigskip \noindent Also, for example, we already learned in Chapter \ref{proba_02_intotp} (Subsection \ref{proba_02_intotp_03} in Section  \ref{proba_02_ssec_indep_ev}) how to choose $h_1$  such that 
$\mathbb{E}h_1(1_{A_1},...,1_{A_{\ell_1}})$ be of the form

$$
\mathbb{P}(A_1^{\prime},...,A_i^{\prime},...,A_{\ell_1}^{\prime}), \ \ (IN02)
$$

\bigskip \noindent where $A_i^{\prime}=A_i$ or $A_i^{\prime}=A_i^{c}$, $1\leq i \leq \ell_1$. In general, for any element $Z_1 \in \mathcal{I}(\mathcal{C}_1)$ of the form 
$h_1(1_{A_1},...,1_{A_{\ell_1}})$ with $\{A_1, ..., A_{\ell_1}\} \subset \mathcal{C}_1$, $\ell_1$, we have

$$
h_1(1_{A_1},...,1_{A_{\ell_1}})=\sum_{\varepsilon \in D_{e\ll_1}} h(\varepsilon) 1_{\prod_{1\leq i \leq \ell_1}} A_i^{(\varepsilon_i)}, \ \ (IN02)
$$

\noindent and $Z_1$ is simply a finite linear combination of elements of $\mathcal{A}_1$. So the independence between $\mathcal{C}_1$ and $\mathcal{C}_2$ is that of
 $\mathcal{A}_1$ and $\mathcal{A}_2$ since the factorization is preserved by finite liner combinations.\\

\bigskip \noindent After these preliminary considerations, we going to prove a first step.\\

\noindent \textbf{Step 1}. We prove that for any $A \in \mathcal{A}_1$, $A$ is independent from $\mathcal{C}_2$.\\

\noindent To see this, define

$$
\mathcal{A}_{0,1}=\{A \in \mathcal{A}_{1}, \forall Z \in \mathcal{I}(\mathcal{C}_2), \ 1_A \text{ independent of } Z \}.
$$

\bigskip \noindent By the assumption we have that $\mathcal{C}_1 \subset \mathcal{A}_{0,1}$. Let us quickly prove that $\mathcal{A}_{0,2}$ is a $\sigma$-algebra. For sure, $\Omega \in \mathcal{A}_{0,2}$. If $A \in \mathcal{A}_{0,1}$, $1_{A^c}$ is a measurable function of $1_A$, and by this is still in $\mathcal{A}_{0,1}$.\\

\noindent Let $(A_1,A_2) \in \mathcal{A}_{0,1}$, $A_2 \subset A_1$. For any non-negative and measurable functions $h(1_{A_2 \setminus A_1})$ and $\ell(Z)$ of $1_{A_2 \setminus A_1}$ and $Z \in \mathcal{I}(\mathcal{C}_2)$, we have

\begin{eqnarray*}
h(1_{A_2 \setminus A_1})&=&h(0) (1_{(A_2 \setminus A_1)^c}) + h(1) (1_{A_2}-1_{A_1})\\
&=&h(0) (1_{A_2^c}+1_{A_1}) + h(1) (1_{A_2}-1_{A_1})
\end{eqnarray*}
 
\bigskip \noindent Hence, by multiplying $h(1_{A_2 \setminus A_1})$ by $\ell(Z)$ and by taking the expectations, we will be able to factorize $\mathbb{E}(\ell(Z) 1_{B})$ for $B \in \{A_1, A_2, A_1^c\}$ in all the terms of the products and, by this, we get   

$$
\mathbb{E} (h(1_{A_2 \setminus A_1})\ell(Z)) = \mathbb{E} (h(1_{A_2 \setminus A_1})) \mathbb{E}(\ell(Z)).
$$

\bigskip \noindent We get that $A_2\setminus A_1 \in \mathcal{A}_{0,1}$.\\

\noindent Finally, let $(A_k)_{leq k \geq 0} \in \mathcal{A}_{0,1}$ be a sequence of pairwise disjoint elements of $\mathcal{A}_{0,1}$. We define 
$$
B_n=\bigcup_{1\leq k \leq n}A_k =\sum_{k \geq 0} A_k \text{ and } B_n=\bigcup_{1\leq k \leq n}A_k =\sum_{1\leq k \leq n} A_k, \ n\geq 0.
$$ 

\bigskip \noindent For any non-negative and measurable functions $h(1_{B_n})$ and $\ell(Z)$ of $1_{B_n}$ and $Z \in \mathcal{I}(\mathcal{C}_2)$, we have

$$
h(1_{B_n}=h(0)\biggr(1-\sum_{1\leq k \leq n} 1_{A_k}\biggr) + h(1) \sum_{1\leq k \leq n} 1_{A_k}.
$$

\bigskip \noindent Here again, by multiplying $h(1_{B_n})$ by $\ell(Z)$ and by taking the expectations, we will be able to factorize any $\mathbb{E}(\ell(Z) 1_{B})$ for $B \in \{B_k, \ 1\leq k \leq n\}$ in all the terms of the product and, by this, we get also that :   

$$
\forall n\geq 0, \ \mathbb{E} (h(1_{B_n})\ell(Z)) = \mathbb{E} (h(1_{B_n})) \mathbb{E}(\ell(Z)).
$$

\bigskip \noindent Next by letting $n\uparrow +\infty$, we get by the Monotone Convergence Theorem that

$$
\mathbb{E} (h(1_{B})\ell(Z)) = \mathbb{E} (h(1_{B})) \mathbb{E}(\ell(Z)),
$$

\bigskip \noindent any non-negative and measurable functions $h(1_{B_n})$ and $\ell(Z)$ of $1_{B_n}$ and of $Z \in \mathcal{I}(\mathcal{C}_2)$. This proves that $B \in \mathcal{A}_{0,1}$. In summary
$\mathcal{A}_{0,1}$ is a Dynkin system containing the $\pi$-system. So by the $\lambda-\pi$-Lemma (See  \cite{ips-mestuto-ang}, Doc 04-02, Chapter 5), it contains  $\mathcal{A}_1$. We conclude that 
$\mathcal{A}_1=\mathcal{A}_{0,1}$ and we get that :\\

\noindent \textit{Any element of $\mathcal{A}_1$ is independent of $\mathcal{I}(\mathcal{C}_2)$}.\\

\bigskip \noindent \textbf{Step 2}. For any $Z_1 \in \mathcal{I}(\mathcal{A}_1)$ of the form $h_1(1_{A_1},...,1_{A_{\ell_1}})$ with $\{A_1, ..., A_{\ell_1}\} \subset \mathcal{A}_1$, $\ell_1$, we have $Z_1$ independent of $\mathcal{C}_2$.\\

\noindent This is an easy consequence of Formula (IN01) and the previous result.\\

\bigskip \noindent \textbf{Final Step 3}. Put

$$
\mathcal{A}_{0,2}=\{B \in \mathcal{A}_{2}, \forall Z \in \mathcal{I}(\mathcal{A}_1), \ 1_B \text{ independent of } Z \}.
$$
 
\bigskip \noindent By the previous steps, $\mathcal{A}_{0,2}$ includes $\mathcal{C}_{2}$. We use the same techniques as in Step 1 to prove that $\mathcal{A}_{0,2}$ is a Dynkin-system and get that 
$\mathcal{A}_{0,2}=\mathcal{A}_{2}$ by the classical methods. Next, we proceed to the same extension as in Step 2 to conclude that any elements of $\mathcal{I}(\mathcal{A}_1)$ is independent of any
other element of $\mathcal{I}(\mathcal{A}_2)$. $\blacksquare$\\

\noindent Now, we may go back the proof of the Kolmogorov law.\\

\bigskip \noindent \textbf{(C) Proof of the zero-one law}. Define the $\sigma$-algebras 

$$
\mathcal{A}_{part, \ n}=\sigma(\{X_k^{-1}(B), \ B \in \mathcal{F}, \ 0\leq k \leq n\}).
$$ 

\noindent If $A$ is a tail event, hence for each $n\geq 1$, $A \in \mathcal{A}_{part, \ n}$. Hence, by the principle underlying Formula (P01) at the beginning at the section and by Proposition above,
we get that $A$ is independent to any $\mathcal{A}_{part, \ n}$, $n\geq 1$. Since these latter sub-classes are $\pi$-system (being $\sigma$-algebras),  $A$ is also independent of 

$$
\sigma\biggr( \mathcal{A}_{part, \ n} \biggr)=\sigma(\{X_k^{-1}(F), \ F \in \mathcal{F}, \ k \geq 0\})=\mathcal{A}_{part, \ 0}.
$$

\bigskip \noindent Since $A \in \mathcal{A}_{part, \ 0}$, we get that $A$ is independent to itself, that is $\mathbb{P}(A)=\mathbb{P}(A \cap A)=\mathbb{P}(A)\mathbb{P}(A)$. The equation 
$\mathbb{P}(A)^2=P(A)$ has only two solutions $0$ or $1$ in $[0,1]$. $\square$\\

\noindent \textbf{(D) Limits Laws for independent random variables}.\\

\noindent We are going to derive series of a three interesting asymptotic laws from the Kolmogorov Inequality (Inequality 14 in Chapter \ref{proba_02_ineg}), the last of them being the celebrated
Three-series law of Kolmogorov.\\

\noindent Let $X_1$, $X_2$, ... be independent centered and square integrable random variables. We denote $\mathbb{V}ar(X_i)=\sigma_i^2$, $1\leq i \leq n$. Define

$$
C_{\infty}=\inf\{C>0, \forall k\geq 0, \ |X_k| \leq C \ \textit{a.s.}\}.
$$

\noindent  Define the partial sums by

$$
S_0=0, \ S_k=\sum_{i=1}^k X_i, \ k\geq 1  \text{ and } s_0=0, \ s_k^2=\sum_{i=1}^k \sigma_i^2.
$$

\noindent We have :

\begin{proposition} \label{proba_02_conv_propd01} The following statements hold.\\

\noindent (1) If $s_n^2$ converges $\sigma^2 \in \mathbb{R}$ as $n\rightarrow +\infty$, then $(S_n)_{n\geq 0}$ converges \textit{a.s} to a \textit{a.s.} finite (possibly constant) random variable.\\

\noindent (2) If the sequence $(X_n)_{n\geq 0}$ is uniformly bounded, that is $C_{\infty}$ is finite, $s_n^2$ converges in $\mathbb{R}$ as $n\rightarrow +\infty$ if and only if $(S_n)_{n\geq 0}$ converges \textit{a.s.}\\

\noindent More precisely, if $s_n^2$ diverges as $n\rightarrow +\infty$ and if $C_{\infty}$ is finite, then $(S_n)_{n\geq 0}$ diverges on any measurable subset of $\Omega$ with a positive probability, that is $(S_n)_{n\geq 0}$ non-where converges.
\end{proposition}

\bigskip \noindent  \textbf{Proof}. Since the random variables are centered and independent, we have for any $0\leq k\leq n$,

$$
\mathbb{V}(S_{n}-S_k)=\sum_{k<j\leq n} \mathbb{V}(X_j)= s_{n}^2-s_{k}^2.
$$

\bigskip \noindent  Let us apply the right-hand Inequality 14 in Chapter \ref{proba_02_ineg}, to get for any $\varepsilon>0$, for any $0\leq k\geq n$,

$$
\mathbb{P}(\max(|S_{k+1}-S_k|,...,|S_{n}-S_k|) \geq \varepsilon ) \leq \frac{|s_{n}^2-s_{k}^2|}{\varepsilon^2},
$$

\bigskip \noindent  in other words, for any $k\geq 0$, $n\geq 0$, for any $\varepsilon>0$

$$
\mathbb{P}\biggr( \bigcup_{1\leq j \leq n} (|S_{k+j}-S_{k}|\geq \varepsilon) \biggr) \leq \frac{s_{n}^2-s_{k}^2}{\varepsilon^2}.
$$

\bigskip \noindent  Let us suppose that $s_n^2$ converges in $\mathbb{R}$ as $n\rightarrow +\infty$, that is $(s_n^2)_{n\geq 0}$ is a Cauchy sequence. By applying the Monotone convergence Theorem, we have for any $k\geq 0$ for any $\varepsilon>0$, 

$$
\mathbb{P}\biggr(\bigcup_{n\geq k} (|S_{n}-S_{k}|\geq \varepsilon)=0\biggr).
$$

\noindent and next any $\varepsilon>0$

$$
\mathbb{P}( \bigcap_{k\geq 0} \biggr( \bigcup_{n\geq k} (|S_{n}-S_{k}|\geq \varepsilon)=0\biggr).
$$

\bigskip \noindent  To conclude, set

$$
\Omega_{p}=\bigcap_{k\geq 0} \biggr(\bigcup_{n\geq k} (|S_{n}-S_{k}|\geq 1/p\biggr), \ p\geq 1 \text{ and } \Omega_{\infty}=\bigcup_{p\geq 1} \Omega_{p}.
$$

\noindent We still have $\mathbb{P}(\Omega_{\infty})=0$ and for any $\omega \in \Omega_{\infty}$, for any $p\geq 1$, $\exists k_0\geq 0$, for all $n\geq k_0$,

$$
|S_{n}-S_{k}|(\omega) < 1/p.
$$

\bigskip \noindent  We conclude that $(S_n)_{n\geq 0}$ is Cauchy on $\Omega_{\infty}$ and then converges on $\Omega_{\infty}$, and simply converges \textbf{a.s.}\\

\noindent It remains to prove that if $C_{\infty}$ is finite and if $s_n^2\rightarrow +\infty$, $(S_n)_{n\geq 0}$ diverges \textit{a.s.}. By Inequality 14 in Chapter \ref{proba_02_ineg}, we also have for any $\varepsilon>0$, for any $0\leq k\leq n$,

$$
\mathbb{P}\biggr( \bigcup_{1\leq j \leq n} \biggr(|S_{k+j}-S_{k}|\geq \varepsilon \biggr) \biggr) \geq 1-\frac{(\varepsilon+C_{\infty})^2}{s_{n}^2-s_{k}^2}.
$$

\bigskip \noindent  For $k$ fixed and $n\rightarrow +\infty$, we get for any $\varepsilon>0$,

$$
\mathbb{P}(\bigcup_{j\geq 0} (|S_{k+j}-S_{k}|\geq \varepsilon)=1.
$$

\bigskip \noindent  and next, for any $\varepsilon>0$,

$$
\mathbb{P}(\bigcap_{k\geq 0} \bigcup_{j\geq 0} (|S_{k+j}-S_{k}|\geq \varepsilon)=1.
$$

\noindent  Denote 

$$
\Omega_0=\biggr(\bigcap_{k\geq 0} \bigcup_{j\geq 0} (|S_{k+j}-S_{k}|\geq 1\biggr).
$$ 

\bigskip \noindent It is clear that $\mathbb{P}(\Omega_0)=1$ and $(S_n)_{n\geq 0}$ is not Cauchy on $\Omega_0$. This proves the two last statements of the proposition. $\square$\\

\begin{proposition} \label{proba_02_conv_propd02} Suppose that the assumptions in Proposition \ref{proba_02_conv_propd02} hold, except we assume that the $X_k$'s are not necessarily centered. Then, if  $(S_n)_{n\geq 0}$ converges \textit{a.s} as $n\rightarrow +\infty$ and the sequence $(X_n)_{n\geq 0}$ is uniformly bounded (that is $C_{\infty}<+\infty$), then the two sequences $(s_n^2)_{n\geq 0}$) and
$\left(\sum_{1\leq k \leq n} \mathbb{E}X_k\right)_{n\geq 0}$ both converge to finite numbers as $n\rightarrow +\infty$.
\end{proposition}   

\bigskip \noindent  \textbf{Proof}. It uses the Kolmogorov construction of probability spaces. At this stage, we know this result only in finite distribution (See Chapter \ref{proba_02_rv}, Section \ref{proba_02_rv_sec_04}, Point (c5)). Here, we anticipate and use Theorem  \ref{KolmConst} (see page  \pageref{KolmConst}) in Chapter \ref{thfondamentalKolm}, and say:\\

\noindent There exists a probability space holding independent random variables $X_k$, $X_k^{\prime}$, $k\geq 0$ such that for each $k$, $X_k=_{d}X_k^{\prime}$.\\

\noindent \label{symmetrizezMzthod} Let us suppose that the assumptions hold and let us define the symmetrized sequence $X^{(s)}_k=X_k-X_k^{\prime}$, $k\geq 0$. Then the sequence $(X^{(s)}_k)_{k\geq 0}$ is centered and uniformly bounded by $2C_{\infty}$. Now, if $\left(\sum_{1\leq k \leq n} X_k\right)_{n\geq 0}$ converges \textit{a.s.}, so does $\left(\sum_{1\leq k \leq n} X_k^{\prime}\right)_{n\geq 0}$ by the equality in law. Hence $\left(\sum_{1\leq k \leq n} X_k^{(s)}\right)_{n\geq 0}$ also converges. Next, by applying Point (2) of Proposition \ref{proba_02_conv_propd02}, the sequence 
$\left(\sum_{1\leq k \leq n} \mathbb{V}(X_k^{(s)})\right)_{n\geq 0}$ converges. Since $\sum_{1\leq k \leq n} \mathbb{V}(X_k^{(s)})=2s_n^2$, we have the first conclusion.\\

\noindent It remains to prove that $\left(\sum_{1\leq k \leq n} \mathbb{E}X_k\right)_{n\geq 0}$ converges. But we have for all $n\geq 0$,

$$
\sum_{1\leq k \leq n} \mathbb{E}(X_k) = \sum_{1\leq k \leq n} X_k - \sum_{1\leq k \leq n} (X_k-\mathbb{E}(X_k))
$$

\bigskip \noindent From this and from that assumption that $(S_n)_{n\geq 0}$ converges \textit{a.s}, we may apply Point (1) of Proposition \ref{proba_02_conv_propd02} to see that the second series in the right-hand of the formula above converges and get our last conclusion. $\square$\\

\noindent \textbf{Remark}. To fully understand this proof, the reader should seriously know the Kolmogorov construction Theorem and its consequences. For example, because of the independence, the vectors
$(X_0,...,X_k)$ and $(X_1^{\prime},...,X_k^{\prime})$ have the same law of $k\geq 0$ and by this, the sequences $(X_k)_{k\geq 0}$ and $(X_k^{\prime})_{k\geq 0}$ have the same law as stochastic processes. So the \textit{a.s.} depends only on the probability law of $(X_k)_{k\geq 0}$, the proved results remain valid on any other probability space for a sequence of the same probability law. We advice the reader to come back to this proof after reading Chapter \ref{thfondamentalKolm}. \\
 
\noindent Before we continue, let us denote for any real-valued random variable $X$ and a real-number $c>0$, the truncation of $X$ at $c$ by 
$$
X^{(c)}=X 1_{(|X|\leq  c)}
$$ 

\bigskip \noindent which is bounded by $c$.\\ 

\begin{proposition} \label{proba_02_conv_propd03} Suppose that the $X_n$'s are are square integrable, centered and independent. If $s_n^2$ converges $\sigma^2 \in \mathbb{R}$ as $n\rightarrow +\infty$, then $(S_n)_{n\geq 0}$ converges \textit{a.s} to a \textit{a.s.} finite . The series $(S_n)_{n\geq 0}$ converges \textit{a.s.} if and only any of the three series below converges :

$$
(i)  \forall c\in \mathbb{R}_+\setminus \{0\}, \ \sum_{k\geq 0}\mathbb{P}(|X_k|\geq c), \ \ \ (ii) \  \sum_{k\geq 0} \mathbb{V}ar\left(X_k^{(c)}\right) \ \ \text{ and }(iii)  \sum_{k\geq 0}\mathbb{E}(X_k^{(c)})
$$
\end{proposition}   

\bigskip \noindent  \textbf{Proof}. \label{proofSeqEquiv} Suppose the three Conditions (i), (ii) and (iii) hold. From (ii), $\left(\sum_{1\leq k \leq n} \biggr( X_k^{(c)}-\mathbb{E}(X_k^{(c)}) \biggr)\right)_{n\geq 0}$ converges. This combined with Condition (iii) implies that $\left(\sum_{1\leq k \leq n} X_k^{(c)}\right)_{n\geq 0}$ converges \textit{a.s.}, based on the remark that $X_k^{(c)}=\biggr( X_k^{(c)}-\mathbb{E}(X_k^{(c)}) \biggr) + 
\biggr(\mathbb{E}(X_k^{(c)})\biggr)$ for all $k\geq 0$.\\

\noindent Next for all $c>0$, the event $(X_k \neq X_k^{(c)})$ occurs only if $(|X_k|\geq c)$ and hence

$$
\mathbb{P}( X_n \neq X_n^{(c)}, \ i.o)=\lim_{n\uparrow +\infty} \mathbb{P}\biggr( \bigcup_{k\geq n} (X_k \neq X_k^{(c)}) \biggr)\leq
\lim_{n\uparrow +\infty} \sum_{k\geq n}\mathbb{P} (|X_k|\geq c).
$$

\noindent Hence, by Condition (i), $\mathbb{P}( X_n \neq X_n^{(c)}, \ i.o)=0$ and next, the series $\sum_{k\geq 0} X_k$ and $\sum_{k\geq 0} X_k^{(c)}$ converge or diverge \textit{a.s.} simultaneously. We get that $\sum_{k\geq 0} X_k$ converges \textit{a.s.}.\\

\noindent Conversely, if $(S_n)_{n\geq 0}$ converges \textit{a.s.}, it follows that $(X_n)_{n\geq 0}$ converges \textit{a.s.} to zero, by Corollary  \ref{probab_02_indep_sec_01_cor01} in Section
 \ref{probab_02_indep_sec_01} below, Condition (i) holds. The latter, by the \textit{a.s.} equivalence between $\sum_{k\geq 0} X_k$ and $\sum_{k\geq 0} X_k^{(c)}$, ensures that 
$\sum_{k\geq 0} X_k^{(c)}$ converges, which by Proposition \ref{proba_02_conv_propd02}, yields Conditions (ii) and (iii).\\

\noindent  It remains to prove that none of the three conditions cannot fail, whenever the series $\sum_{k\geq 0} X_k$ converges \textit{a.s.}. First, by Corollary  \ref{probab_02_indep_sec_01_cor01} in Section  \ref{probab_02_indep_sec_01}, Condition (i) cannot fail and hence the \textit{a.s.} convergence of $\sum_{k\geq 0} X_k^{(c)}$ also cannot fail and this bears Conditions (ii) and (ii). $\square$\\

\noindent Now let us close this introduction to these following important Kolmogorov's Theorems.\\

\noindent \textbf{(E) Strong Law of Large numbers of Limits Laws of Kolmororov}.\\

\noindent Before we state the Kolmogorov laws, we state the following :\\

\noindent \textbf{Kronecker Lemma}. If $(b_n)_{n\geq 0}$ is an increasing sequence of positive numbers and $(x_n)_{n\geq 0}$ is a sequence of finite real numbers such that
 $\left(\sum_{1\leq k \leq n} x_k\right)_{n\geq 0}$ converges to a finite real number $s$, then 

$$
\frac{\sum_{1\leq k \leq n} b_k x_k}{b_n} \rightarrow 0 \ as \ n\rightarrow \infty. \ \Diamond
$$

\bigskip \noindent This Lemma is proved in the Appendix, where it is derived from the Toeplitz Lemma.\\

\noindent Let us begin by \\

\noindent \textbf{(E-a) The Strong Law of Large Numbers for Square integrable and independent random variables}.\\

\begin{theorem} \label{sllnGenKol} Let $(X_n)_{n\geq 0}$ be a sequence of square integrable and independent random variables and let $(b_n)_{n\geq 1}$ be an increasing sequence of finite real numbers. If
$$
\sum_{n\geq 0} \frac{\mathbb{V}(X_n)}{b_n^2} <+\infty, \ \ (CK01)
$$

\bigskip \noindent then we have the following SLLN

$$
\frac{S_n-\mathbb{E}(S_n)}{b_n} \rightarrow 0 \ as \ n\rightarrow +\infty. \ \ (SK01)
$$
\end{theorem}

\bigskip \noindent \textbf{Proof}. The proof comes as the conclusion of the previous developments. Suppose that assumptions of the theorem hold and Condition (CK01) is true. By Proposition \ref{proba_02_conv_propd01}, we have

$$
\sum_{n\geq 0} \frac{X_n-\mathbb{E}(X_n)}{b_n} <+\infty, \ a.s.\ \ (SK02)
$$

\bigskip \noindent  Applying the Kronecker Lemma with $x_k=X_k-\mathbb{E}(X_k)$ and the same sequence $(b_n)_{n\geq 0}$ leads to (SK01). $\square$\\

\noindent \textbf{Example}. If the $X_n$'s have the same variance $\sigma^2$, we may take $b_n=n$, $n\geq 1$ and see that Condition (SK01) is verified since

$$
\sum_{n\geq 1} \frac{\mathbb{V}(X_n)}{b_n^2}=\sigma^2  \sum_{n\geq 1} \frac{1}{n^2}<+\infty,
$$

\bigskip \noindent and next,

$$
\frac{S_n-\mathbb{E}(S_n)}{n} \rightarrow 0 \ as \ n\rightarrow +\infty. \ \ (SK03)
$$

\bigskip \noindent We find again the simple SLLN as in Corollary \ref{probab_02_indep_sec_01_corSLLN} above.\\

\noindent Now, what happens if the first moments of the $X_n$ exist but we do not have information about the second moments of the $X_n$'s. We already saw in Kintchine's Theorem \ref{probab_02_indep_sec_01_theoWLLN} that we have a WLLN if the $X_n$ has the same Law. Here again, the Kolmogorov theory goes far and establishes the SLLN  even if the common second moment is infinite. We are going to see this in the next part.\\

\bigskip \noindent \textbf{(E-b) The Strong Law of Large Numbers for independent and identically random variables with finite mean}.\\

\noindent We will need the following simplified Toeplitz lemma which is proved in the Appendix in its integrability.\\

\noindent \textbf{Simple Toeplitz Lemma}. Suppose that $k(n)=n$ for all $n\geq 1$. Let $(c_k)_{k\geq 0}$ be sequence such that the sequence $(b_n)_{n\geq 0}=(\sum_{1\leq k\leq n}|c_k|_{n\geq 0}$ is non-decreasing and $b_n\rightarrow \infty$. If $x_n\rightarrow x \in \mathbb{R}$ as $n\rightarrow +\infty$, then
$$
\frac{1}{b_n} \sum_{1\leq k \leq n} c_k x_k\rightarrow x\ \ as \ \ n\rightarrow +\infty. \ \ \Diamond\\
$$

\begin{theorem} \label{sllnEPKol} Let $(X_n)_{n\geq 0}$ be a sequence of independent and identically distributed random variables having the same law as $X$. Then\\

\noindent (a) $\mathbb{E}|X|<+\infty$\\

\noindent if and only if \\

\noindent (b) $Sn/n$ converges \textbf{a.s.} to a finite number $c$, which is necessarily $\mathbb{E}(X)$.
\end{theorem}

\bigskip \noindent \textbf{Proof}. Set $A_n=(|X|\geq n)$, $n\geq 0$ with $A_0=\Omega$ clearly.\\

\bigskip \noindent Now suppose that Point (b) holds. We have to prove that of $\mathbb{E}|X|$ is finite. If $X$ is bounded, there is nothing to prove. If not, the upper endpoint of 
$uep(X)$ is infinite and Formula (DF3) in  Proposition \ref{proba_02_sec_01_esp} (See Chapter \ref{proba_02_intotp}, page \pageref{proba_02_sec_01_esp})

$$
-1 + \sum_{n \in [0, \ [uep(X)]^{+}]} \mathbb{P}(|X|\geq n) \leq \mathbb{E}|X|\leq \sum_{n \in [0, \ [uep(X)]^{+}]} \mathbb{P}(|X|\geq n)
$$

\bigskip \noindent  becomes

$$
-1 + \sum_{n \geq 0} \mathbb{P}(|X|\geq n) \leq \mathbb{E}|X|\leq \sum_{n \neq 0} \mathbb{P}(|X|\geq n). \ \ (DF4)
$$

\bigskip \noindent  Then we have 

$$
\frac{X_n}{n}=\frac{S_{n}-S_{n-1}}{n}=\frac{S_{n}}{n}-\frac{n-1}{n}\frac{S_{n-1}}{n-1}\rightarrow c-c=0 \ a.s.
$$

\bigskip \noindent By the Borel-Cantelli Corollary \ref{probab_02_indep_sec_01_cor01}, the serie $\sum_{n\geq 0} \mathbb{P}(|X_n/c|\geq 1)$ is convergent, that $\sum_{n\geq 0} \mathbb{P}(A_n)$ is finite and by Formula (DF4), $\mathbb{E}|X|$ is finite.\\

\noindent Now suppose that $\mathbb{E}|X|$ is finite. If $X$ is bounded,  we are in the case of the last example above with $\mathbb{V}ar(X)=\sigma^2$ is stationary and we have that $S_n/n$ converges to $\mathbb{E}(X)$ \textit{a.s.}. If not, we use the truncated random variables $X_k^{(t)}=X_k 1_{(|X_k|< k)}$, $k\geq 1$ and let 
$S_0^{(t)}=0$, $S_1^{(t)}=X_1^{(t)}$, $S_n^{(t)}=X_1^{(t)}+...+X_n^{(t)}$, $n\geq 2$.\\

\noindent We already explained in page \pageref{proofSeqEquiv} that $S_n^{(t)}/n$ and $S_n/n$ have the same \textit{a.s.} limit or diverge \textit{a.s.} together whenever $\sum_{n\geq 1} \mathbb{P}(|X_n|\geq n)$. But since $\mathbb{E}|X|$ is finite, the series $\sum_{n\geq 1} \mathbb{P}(|X_n|\geq n)$ converges by Formula (DF4). Hence we only have to prove that $S_n^{(t)}/n \rightarrow \mathbb{E}(X)$ \textit{a.s.}. Now, by the Dominated Convergence Theorem, we have

$$
\mathbb{E}(X_n^{(t)}= \int X 1_{(|X|< n)} \ d\mathbb{P} \rightarrow \mathbb{E}{X} \ as \ n\rightarrow +\infty.
$$

\bigskip \noindent By applying the simple Toeplitz Lemma with $c_k=1$ and $x_k=\mathbb{E}(X_n^{(t)}$, we get

$$
\frac{\mathbb{E}(S_n^{(t)}}{n}) \rightarrow \mathbb{E}{X}.
$$

\bigskip \noindent So, our task is to prove that

$$
\frac{S_n^{(t)}-\mathbb{E}(S_n^{(t)}}{n}) \rightarrow 0 \ a.s.
$$

\bigskip \noindent But this derives form Theorem \ref{sllnGenKol} whenever we have

$$
\sum_{n\geq 1} \frac{\mathbb{V}ar(X_n^{(t)})}{n^2} <+\infty.
$$

\bigskip \noindent But we have 

\begin{eqnarray*} 
\sum_{n\geq 1} \frac{\mathbb{V}ar(X_n^{(t)})}{n^2}&=&\sum_{n\geq 1} \frac{\mathbb{E}(X_n^{(t)})^2-(\mathbb{E}(X_n^{(t)}))^2}{n^2}\\
&\leq & \sum_{n\geq 1} \frac{\mathbb{E}(X_n^{(t)})^2}{n^2}\\
&=&\sum_{n\geq 1} \mathbb{E}\biggr(\frac{X^2}{n^2} 1_{(|X|<n)}\biggr).
\end{eqnarray*}

\bigskip \noindent Next, define $B_m=(m-1 \leq |X| <m)$, $m\geq 1$. For $m\geq 1$ fixed, we have

$$
\frac{X^2}{n^2} 1_{(|X|<n)\cap B_m}=\frac{X^2}{n^2} 1_{\emptyset}=0 \ for \ n<m,
$$

\bigskip \noindent and for $n\geq m$
$$
\frac{X^2}{n^2} 1_{(|X|<n)\cap B_m}=\frac{X^2}{n^2} 1_{B_m} \leq \frac{m^2}{n^2} 1_{B_m},
$$

\bigskip \noindent so that

\begin{eqnarray*}
\sum_{n\geq 1}\frac{X^2}{n^2} 1_{(|X|<n)\cap B_m}&=&\sum_{n\geq m}\frac{X^2}{n^2} 1_{(|X|<n)\cap B_m}\\
&\leq & \biggr(m^2 \sum_{n\geq m} \frac{1}{n^2}\biggr) 1_{B_m}.
\end{eqnarray*}

\bigskip \noindent By comparing the series of the form $\sum_{n\geq m} f(n)$ with the integral $\int_{x\geq m} f(x) dx$ for a non-decreasing and continuous function $f(x)=x^{-2}$, we have  

$$
\sum_{n\geq m} \frac{1}{n^2} \leq \int_{m}^{+\infty} x^{-2} \ dx =1/m.
$$

\bigskip \noindent Hence, we have

\begin{eqnarray*}
\sum_{n\geq 1}\frac{X^2}{n^2} 1_{(|X|<n)\cap B_m}\leq m 1_{B_m}=(1+(m-1)) 1_{B_m} \leq (1+|X|) 1_{B_m}.
\end{eqnarray*}

\bigskip \noindent Since, we obviously have $\sum_{m\geq 1} B_m=\Omega$, we may sum over $m$ in the previous formula to have

\begin{eqnarray*}
\sum_{n\geq 1}\frac{X^2}{n^2} 1_{(|X|<n)} \leq (1+|X|).
\end{eqnarray*}

\noindent We arrive at

\begin{eqnarray*}
\sum_{n\geq 1} \frac{\mathbb{V}ar(X_n^{(t)})}{n^2} \leq \mathbb{E}\sum_{n\geq 1}\frac{X^2}{n^2} 1_{(|X|<n)}\leq (1+\mathbb{E}|X|)<+\infty.
\end{eqnarray*}

\bigskip \noindent We reached the desired condition which allows to conclude the proof. $\blacksquare$\\

\noindent This nice theory of Kolmogorov opens the wide field of SLLN's. The first step for the generalization will be the H\'aj\`ek-R\'enyi approach we will see soon in special monograph reserved
to limits laws for sequences of random variables of arbitrary probability laws.

%% file: proba_02_07_Gauss_series.tex
\section{Convergence of Partial sums of independent Gaussian random variables} \label{probab_02_conv_gauss}

\noindent Let us give the following interesting equivalences between different types of convergences for partial sums of independent Gaussian real-valued random variables.\\

\noindent we have

\begin{theorem} Let $(X_n)_{n\geq 1}$ be sequence of independent and centered Gaussian real-valued random variables defined on the same probability space $(\Omega, \mathcal{A}, \mathbb{P})$. Let us define their partial sums $S_n=\sum_{1\leq k \leq n} X_n$, $n\geq 1$ and $s_n^2=\sum_{1\leq k \leq n} \mathcal{V}ar(X_n)$, $n\geq 1$. The the following convergences are equivalent, as $n\rightarrow +\infty$,\\

\noindent (1) $(S_n)_{n\geq 1}$ converges \textit{a.s.} to an \textit{a.s.} finite random variable $Z$.\\

\noindent (2) $(S_n)_{n\geq 1}$ converges in probability to an \textit{a.s.} finite random variable $Z$.\\

\noindent (3) $(S_n)_{n\geq 1}$ weakly converges to an \textit{a.s.} finite random variable $Z$.\\

\noindent (4) $(s_n)_{n\geq 1}$ converges in $\mathbb{R}$.\\

\noindent (5) $(S_n)_{n\geq 1}$ converges in $L^2$.\\

\end{theorem}

\bigskip \noindent \textit{Proof}. The proof is based on the compa rison between type of convergences in Chapter \ref{proba_02_conv}. The implication $(1) \rightarrow (2)$. Next $(2) \rightarrow (3)$ by Point (a) of Proposition  \ref{proba_02_conv_sec_06_prop_01} in Section \ref{proba_02_conv}. Further (3) that for each $t \in \mathbb{R}$,

$$
\mathbb{E} \exp(itS_n) \rightarrow \mathbb{E} \exp(X_{(p, \ +\infty)}).
$$

\bigskip \noindent Since the $X_n$ are independent, we have $S_n\sim \mathcal{N}(0,s_n^2)$. Hence for all $t \in \mathbb{R}$,

$$
 \exp(-t s_n^2/2) \rightarrow \mathbb{E} \exp(X_{(p, \ +\infty)}),
$$

\bigskip \noindent This is possible only if $s_n^2$ converge in $\mathbb{R}$, where we took into account the fact that $Z$ is \textit{a.s.} finite. Now, by Proposition \ref{proba_02_conv_propd01} 
(in Section \ref{probab_02_sllnKol} in Section \ref{proba_02_conv}, page \pageref{probab_02_sllnKol}) (4) implies (1). By this circular argument, the assertions (1) to (4) are equivalent.\\

\noindent Let us handle Assertion (5). Suppose (5) holds with $s_n^2\rightarrow s^2$. Let us denote $S=\sum_{n\geq 1} X_n$. We have for all $n\geq 1$

$$
\mathbb{E}(S-S_n)^=\|S-S_n\|_2^2= (s^2-s_n^2) \rightarrow 0.
$$

\bigskip \noindent So (5) implies (5). Finally, suppose that (2) holds. We have that by Point (3) that $S_n$ converge to $S$ in probability and $s_n^2$ converge to $s^2$, and next 

$$
\|S_n\|_2^2=s_n^2 \rightarrow s^2=\|S\|_2^2.
$$

\bigskip \noindent Thus by Point(c4) of Theorem \ref{proba_02_con_sec_th_01}  (in Section \ref{probab_02_sllnKol} in Section \ref{proba_02_conv}, page \pageref{proba_02_con_sec_th_01}), (2) implies (5). The proof is complete now.\\

%% file: proba_02_07_clt_ang.tex
\section{The Lindenberg-Lyapounov-Levy-Feller Central Limit Theorem} \label{probab_02_indep_sec_03}

\noindent We do not treat the Central limit Theorem on $\mathbb{R}^d$, $d\geq 2$, which is addressed in \cite{ips-wcrv-ang} in its simplest form.\\

\noindent We already described the \textit{CLT} question on $\mathbb{R}$ with \textit{iid} sequences. The current section will give the most finest results for independent random variables. Researchers are trying to export the Linderberg-Levy-Feller Central Limit Theorem to abstract spaces under dependance conditions. In that generalization process, mastering the techniques which are used in the independence case significantly help.\\

\bigskip \noindent Let us begin by the key result of Lyapounov.\newline

\noindent \textbf{(A) Lyapounov Theorem}.\\

\begin{theorem} \label{probab_02_indep_sec_03_th01} Let $X_{1},$ $X_{2},$ ... a sequence of real and independent random variables centered at expectations, with finite $(n+\delta )-$moment, $\delta >0.$ Put for each $n\geq 1$, $S_{n}=X_{1}+...+X_{n}$ and $s_{n}^{2}=$ $\mathbb{E}X_{1}^{2}+\mathbb{E}X_{2}^{2}+...+\mathbb{E}X_{n}^{2}$. We denote $\sigma _{k}^{2}=\mathbb{E}X_{k}^{2},k\geq 1$ and $F_{k}$\ denotes the probability distribution function of $X_{k}$. Suppose that

\begin{equation}
\frac{1}{s_{n}^{2+\delta }}\sum_{k=1}^{n}\mathbb{E}\left\vert
X_{k}\right\vert ^{2+\delta }\rightarrow 0\text{ as }n\rightarrow \infty .
\label{LYAP}
\end{equation}

\bigskip \noindent Then, we have as $n\rightarrow +\infty$, 
\begin{equation*}
S_{n}/s_{n}\rightsquigarrow N(0,1). 
\end{equation*}
\end{theorem}

\bigskip \noindent \textbf{Proof of Theorem \ref{probab_02_indep_sec_03_th01}}. According to Lemma 3 below, if (\ref{LYAP}) holds for $%
\delta >1,$ then it holds for $\delta =1.$ So it is enough to prove the
theorem for $0<\delta \leq 1.$ \ By lemma 4 below, the assumption (\ref{LYAP}%
) implies $s_{n}\rightarrow +\infty $ and%
\begin{equation*}
\max_{1\leq k\leq n}\left( \frac{\sigma _{k}}{s_{n}}\right) ^{2+\delta }\leq
\max_{1\leq k\leq n}\frac{\mathbb{E}\left\vert X_{k}\right\vert ^{2+\delta }%
}{s_{n}^{2+\delta }}\leq \frac{1}{s_{n}^{2+\delta }}\sum_{k=1}^{n}\mathbb{E}%
\left\vert X_{k}\right\vert ^{2+\delta }=:A_{n}(\delta )\rightarrow 0.
\end{equation*}

\bigskip \noindent Let us use the expansion of the characteristic functions 
\begin{equation*}
f_{k}(u)=\int e^{iux}dF_{k}(x)
\end{equation*}

\bigskip \noindent at the order two to get for each $k,1\leq k\leq n$ as
given in Lemma 1 below 
\begin{equation}
f_{k}(u/s_{n})=1-\frac{u^{2}}{2}.\frac{\sigma _{k}^{2}}{s_{n}^{2}}+\theta
_{nk}\frac{\left\vert u\right\vert ^{2+\delta }\mathbb{E}\left\vert
X_{k}\right\vert ^{2+\delta }}{s_{n}^{2+\delta }}.  \label{expan-eachk}
\end{equation}

\bigskip \noindent Now the characteristic function of $S_{n}/s_{n}$ is, for $%
u\in R,$
\begin{equation*}
f_{S_{n}/s_{n}}(u)=\prod\limits_{k=1}^{n}f_{k}(u/s_{n})
\end{equation*}

\noindent S that is
\begin{equation*}
\log f_{S_{n}/s_{n}}(u)=\sum_{k=1}^{n}\log f_{k}(u/s_{n}).
\end{equation*}

\bigskip \noindent Now, we use the uniform expansion of $\log (1+u)$ at the
neighborhood at $1$, that is

\begin{equation}
\sup_{\left\vert u\right\vert \leq z}\left\vert \frac{\log (1+u)}{u}\right\vert =\varepsilon (z)\rightarrow 0.  \label{expans-log}
\end{equation}

\noindent \noindent For each $k$ in (\ref{expan-eachk}), we have

\begin{equation}
f_{k}(u/s_{n})=1-u_{kn}  \label{expans-unif}
\end{equation}

\bigskip \noindent with the uniform bound

\begin{eqnarray*}
\left\vert u_{kn}\right\vert &\leq& \sum_{j=1}^{n}\frac{\left\vert
u\right\vert ^{2}}{2}.\frac{\sigma _{k}^{2}}{s_{n}^{2}}+\frac{\left\vert
u\right\vert ^{2+\delta }\mathbb{E}\left\vert X_{k}\right\vert ^{2+\delta }}{%
s_{n}^{2+\delta }}\\
&=&\frac{\left\vert u\right\vert ^{2}}{2}.\max_{1\leq k\leq
n}\left( \frac{\sigma _{k}^{2}}{s_{n}}\right) +\frac{\left\vert u\right\vert
^{2+\delta }\sum_{j=1}^{n}\mathbb{E}\left\vert X_{k}\right\vert ^{2+\delta}
}{s_{n}^{2+\delta }}=u_{n}.
\end{eqnarray*}

\bigskip \noindent By applying (\ref{expans-log}) to (\ref{expans-unif}), we get
\begin{equation*}
\log f_{k}(u/s_{n})=-u_{kn}+\theta _{n}u_{kn}\varepsilon (u_{n})
\end{equation*}

\bigskip \noindent and next 
\begin{eqnarray*}
\log f_{S_{n}/s_{n}}(u)&=&\sum_{k=1}^{n}\log f_{k}(u/s_{n})\\
&=&-\frac{u^{2}}{2}+\left\vert u\right\vert ^{2+\delta }\theta _{n}A_{n}(\delta )+(\frac{u^{2}}{%
2}+\left\vert u\right\vert ^{2+\delta }\theta _{n}A_{n}(\delta ))\varepsilon
(u_{n})\\
&\rightarrow& -u^{2}/2.
\end{eqnarray*}

\bigskip \noindent We get for u fixed,%
\begin{equation*}
f_{S_{n}/s_{n}}(u)\rightarrow \exp (-u^{2}/2).
\end{equation*}

\bigskip \noindent This completes the proof. $\square$ \newline

\bigskip \noindent \textbf{An expression of Lyapounov Theorem using
triangular arrays}.\newline

\bigskip \noindent Since the proof is based on the distribution of $\{X_{k},1\leq k\leq n\}$ for each $n\geq 1$, it may be extended to triangular array to the following corollary.\newline

\bigskip

\begin{corollary} \label{probab_02_indep_sec_03_cor01} Consider the triangular array $\{X_{nk},1\leq k\leq n,n\geq 1\}$. Put for each $n\geq 1,$ $S_{nn}=X_{n1}+...+X_{nn}$ and $s_{nn}^{2}=$ $\mathbb{E}X_{n1}^{2}+\mathbb{E}X_{n2}^{2}+...+\mathbb{E}X_{nn}^{2}$. Suppose that for each $n\geq 1$, the random variables $X_{nk},1\leq k\leq n,$ \ are centered and independent such
that 

\begin{equation}
\frac{1}{s_{nn}^{2+\delta }}\sum_{k=1}^{n}\mathbb{E}\left\vert
X_{nk}\right\vert ^{2+\delta }\rightarrow 0\text{ as }n\rightarrow \infty .
\end{equation}

\bigskip \noindent Then
\begin{equation*}
S_{nn}/s_{nn}\rightsquigarrow N(0,1).
\end{equation*}
\end{corollary}

\bigskip \bigskip \noindent Now, we are able to prove the Lyapounov-Feller-Levy Theorem (see Lecam for an important historical note with
the contribution of each author in this final result).\newline

\bigskip
\noindent \textbf{(B) The General Central Limit Theorem on $\mathbb{R}$}.\\

\begin{theorem} \label{probab_02_indep_sec_03_th02} Let $X_{1}$, $X_{2},$ ... a sequence of real and independent random variables centered at expectations, with finite $(n+\delta )-$moment, $\delta >0.$ Put for each $n\geq 1,$ $S_{n}=X_{1}+...+X_{n}$ and $s_{n}^{2}=$ $\mathbb{E}X_{1}^{2}+\mathbb{E}X_{2}^{2}+...+\mathbb{E}X_{n}^{2}.$ We denote $\sigma _{k}^{2}=\mathbb{E}X_{k}^{2},k\geq 1$ and  $F_{k}$\ denotes the probability distribution function of $X_{k}$. We have the equivalence between 

\begin{equation}
\max_{1\leq k\leq n}\left( \frac{\sigma _{k}^{2}}{s_{n}}\right) \rightarrow 0%
\text{ \ \ \ and\ \ \ \ }S_{n}/s_{n}\rightsquigarrow N(0,1)  \label{WF1}
\end{equation}

\bigskip \noindent and

\begin{equation}
\forall \varepsilon >0,\text{ \ }\frac{1}{s_{n}^{2}}\sum_{k=1}^{n}\int_{%
\left\vert x\right\vert \geq \varepsilon s_{n}}x^{2}dF_{k}(x)\rightarrow 0%
\text{ as }n\rightarrow \infty .  \label{WF2}
\end{equation}
\end{theorem}

\bigskip  \noindent \textbf{Proof of Theorem \ref{probab_02_indep_sec_03_th02}}.\newline

\bigskip \noindent The proof follows the lines of the proof in \cite{loeve}. But they are extended by more details and adapted and changed in some parts. Much details were omitted. We get them back for making the proof understandable for students who just finished the measure and probability course.\newline

\bigskip \noindent Before we begin, let us establish an important property of

\begin{equation*}
g_{n}(\varepsilon )=\frac{1}{s_{n}^{2}}\sum_{k=1}^{n}\int_{\left\vert x\right\vert \geq \varepsilon s_{n}}x^{2}dF_{k}(x),
\end{equation*}

\bigskip \noindent when (\ref{WF2}) holds. Suppose that this latter holds. We want to show that there exists a sequence $\varepsilon _{n}\rightarrow 0$ such that $\varepsilon _{n}^{-2}g_{n}(\varepsilon _{n})\rightarrow 0$ (this implying also that $\varepsilon _{n}^{-1}g_{n}(\varepsilon
_{n})=o(\varepsilon _{n})\rightarrow 0$ and that $g_{n}(\varepsilon
_{n})=o(\varepsilon _{n}^{2})\rightarrow 0).$ To this end, let $k\geq 1$
fixed$.$ Since $g_{n}(1/k)\rightarrow 0$ as $n\rightarrow \infty ,$ we have $%
0\leq g_{n}(1/k)\leq k^{-3}$ for $n$ large enough.\newline

\bigskip \noindent We will get what we want from an induction on this property. Fix $k=1$ and denote $n_{1}$ an integer such that $0\leq g_{n}(1)\leq 1^{-3}$ for $n\geq n_{1}.$ Now we apply the same property on the sequence $\{g_{n}(\circ ),n_{1}+1\}$ with $k=2$. We find a $n_{2}>n_{1}$ such that $0\leq g_{n}(1/2)\leq 2^{-3}$ for $n\geq n_{2}$. Next we apply the same property on the sequence $\{g_{n}(\circ ),n_{2}+1\}$ with $k=2$. We find a $n_{3}>n_{2}$ such that $0\leq g_{n}(1/3)\leq 3^{-3}$ for $n\geq n_{3}$. Finally, an infinite sequence of integers $n_{1}<n_{2}<..<n_{k}<n_{k+1}<..$ such that for each $k\geq 1,$ one has $0\leq g_{n}(1/k)\leq k^{-3}$ for $n\geq n_{k}$. Put 
\begin{equation*}
\varepsilon _{n}=1/k\text{ on }n_{k}\leq n<n_{k+1}.
\end{equation*}

\bigskip \noindent We surely have $\varepsilon _{n}\rightarrow 0$ and $%
\varepsilon _{n}^{-2}g_{n}(\varepsilon _{n}).$ This is clear from%
\begin{equation*}
\left\{ 
\begin{tabular}{lll}
$\varepsilon _{n}=1/k$ & on & $n_{k}\leq n<n_{k+1}$ \\ 
$\varepsilon _{n}^{-2}g_{n}(\varepsilon _{n})=k^{2}(1/k^{3})\leq (1/k)$ & on
& $n_{k}\leq n<n_{k+1}$.
\end{tabular}
\right. 
\end{equation*}

\bigskip \noindent Now we arge going to use
\begin{equation}
\varepsilon _{n}\rightarrow 0\text{ and }\varepsilon
_{n}^{-2}g_{n}(\varepsilon _{n})\rightarrow 0.  \label{seqTravail}
\end{equation}

\bigskip \noindent \textbf{Proof of (\ref{WF2}) }$\Longrightarrow $\textbf{(\ref{WF1}).} Suppose (\ref{WF2}) holds$.$ Thus there exists a sequence $(\varepsilon _{n})_{n\geq 0}$ of positive numbers such that (\ref{seqTravail}) prevails. First, we see that, for each $j,1\leq j\leq n$,

\begin{eqnarray*}
\frac{\sigma _{j}^{2}}{s_{n}^{2}} &=&\frac{1}{s_{n}^{2}}\int x^{2}dF_{j}=%
\frac{1}{s_{n}^{2}}\left\{ \int_{\left\vert x\right\vert \leq \varepsilon
_{n}s_{n}}x^{2}dF_{j}+\int_{\left\vert x\right\vert \leq \varepsilon
_{n}s_{n}}x^{2}dF_{j}\right\}  \\
&\leq &\frac{1}{s_{n}^{2}}\int_{\left\vert x\right\vert \leq \varepsilon
_{n}s_{n}}x^{2}dF_{k}+\varepsilon _{n}^{2} \\
&\leq &\frac{1}{s_{n}^{2}}\sum_{k=1}^{k}\int_{\left\vert x\right\vert \leq
\varepsilon _{n}s_{n}}x^{2}dF_{k}=g(\varepsilon _{n})+\varepsilon _{n}^{2}.
\end{eqnarray*}

\bigskip \noindent It follows that
\begin{equation*}
\max_{1\leq j\leq n}\frac{\sigma _{j}^{2}}{s_{n}^{2}}\leq g(\varepsilon
_{n})+\varepsilon _{n}^{2}\rightarrow \text{0.}
\end{equation*}

\bigskip \noindent Its remains to prove that $S_{n}/s_{n} \rightarrow \mathbb{N}(0,1).$ To this end we are going to use this array of truncated random variables $\{X_{nk},1\leq k\leq n,n\geq 1\}$ defined as follows. For each fixed $n\geq 1$, we set
\begin{equation*}
X_{nk}=\left\{ 
\begin{tabular}{lll}
$X_{k}$ & if & $\left\vert X_{k}\right\vert \leq \varepsilon _{n}s_{n}$ \\ 
$0$ & if & $\left\vert X_{k}\right\vert >\varepsilon _{n}s_{n}$%
\end{tabular}%
\right. ,1\leq k\leq n.
\end{equation*}

\bigskip \noindent Now, we consider summands $S_{nn}$ as defined in Corollary \ref{probab_02_indep_sec_03_cor01}. Weremark that for any $\eta >0,$

\begin{eqnarray*}
P\left( \left\vert \frac{S_{nn}}{s_{n}}-\frac{S_{n}}{s_{n}}\right\vert >\eta
\right) \leq P\left( \frac{S_{nn}}{s_{n}}\neq \frac{S_{n}}{s_{n}}\right)
=P\left( \frac{S_{nn}}{s_{n}}\neq \frac{S_{n}}{s_{n}}\right)
\end{eqnarray*}

\bigskip \noindent and also, 
\begin{eqnarray*}
\left( \frac{S_{nn}}{s_{n}}\neq \frac{S_{n}}{s_{n}}\right) &=&\left( (\exists
1\leq k\leq n),X_{nk}\neq X_{k}\right)\\
&=&\left( \exists (1\leq k\leq
n),\left\vert X_{k}\right\vert >\varepsilon _{n}s_{n}\right)
=\bigcup\limits_{k=1}(\left\vert X_{k}\right\vert >\varepsilon _{n}s_{n}).
\end{eqnarray*}

\bigskip \noindent We get

\begin{eqnarray*}
P\left( \left\vert \frac{S_{nn}}{s_{n}}-\frac{S_{n}}{s_{n}}\right\vert >\eta
\right)  &\leq &\sum_{k=1}^{n}P(\left\vert X_{k}\right\vert >\varepsilon
_{n}s_{n}) \\
&\leq &\sum_{k=1}^{n}\int_{\left\vert x\right\vert \leq \varepsilon
_{n}s_{n}}dF_{k}=\sum_{k=1}^{n}\int_{\left\vert x\right\vert \leq
\varepsilon _{n}s_{n}}\left\{ \frac{1}{x^{2}}\right\} x^{2}dF_{k} \\
&\leq &\left\{ \frac{1}{(\varepsilon _{n}s_{n})^{2}}\right\}
\sum_{k=1}^{n}\int_{\left\vert x\right\vert \leq \varepsilon
_{n}s_{n}}x^{2}dF_{k} \\
&\leq &\frac{1}{\varepsilon _{n}^{2}}g_{n}(\varepsilon _{n})\rightarrow 0.
\end{eqnarray*}

\noindent Thus $S_{nn}/s_{n}$ and $S_{n}/s_{n}$ are equivalent in probability. This implies that they have the same limit law or do not have a limit law together. So to prove that $S_{n}/s_{n}$ has a limit law, we may prove that $S_{nn}/s_{n}$ has a limit law. Next by Slutsky lemma, it will suffice to establish the limiting law of $S_{nn}/s_{nn}$ whenever $s_{nn}/s_{n}\rightarrow 1$. We focus on this. We begin to remark that, since $\mathbb{E}(X_{k})=0$, we have the decomposition
 
\begin{equation*}
0=\mathbb{E}(X_{k})=\int xdF_{k}=\int_{\left\vert x\right\vert \leq \varepsilon
_{n}s_{n}}xdF_{k}+\int_{\left\vert x\right\vert >\varepsilon
_{n}s_{n}}xdF_{k}
\end{equation*}

\noindent to get that 
\begin{equation*}
\left\vert \int_{\left\vert x\right\vert \leq \varepsilon
_{n}s_{n}}xdF_{k}\right\vert =\left\vert \int_{\left\vert x\right\vert
>\varepsilon _{n}s_{n}}xdF_{k}\right\vert .
\end{equation*}

\noindent We remark also that

\begin{eqnarray}
\mathbb{E}(X_{nk})&=&\int_{\left\vert X_{k}\right\vert \leq \varepsilon
_{n}s_{n}}X_{nk}d\mathbb{P}+\int_{\left\vert X_{k}\right\vert >\varepsilon
_{n}s_{n}}X_{nk}d\mathbb{P} \label{trunc00} \\
&=&\int_{\left\vert X_{k}\right\vert \leq \varepsilon
_{n}s_{n}}X_{k}d\mathbb{P}+\int_{\left\vert X_{k}\right\vert >\varepsilon
_{n}s_{n}}0 \ d\mathbb{P}.  \notag
\end{eqnarray}

\bigskip \noindent Combining all what precedes leads to
\begin{eqnarray*}
\left\vert E(X_{nk})\right\vert &=&\left\vert \int_{\left\vert
X_{k}\right\vert \leq \varepsilon _{n}s_{n}}X_{k}dP\right\vert\\
&=&\left\vert
\int_{\left\vert x\right\vert >\varepsilon _{n}s_{n}}xdF_{k}\right\vert
=\left\vert \int_{\left\vert x\right\vert >\varepsilon _{n}s_{n}}\left\{ 
\frac{1}{x}\right\} x^{2}dF_{k}\right\vert \\
&\leq &\left\vert \int_{\left\vert x\right\vert >\varepsilon _{n}s_{n}}\frac{%
1}{\left\vert x\right\vert }x^{2}dF_{k}\right\vert \leq \frac{1}{\varepsilon
_{n}s_{n}}\left\vert \int_{\left\vert x\right\vert >\varepsilon
_{n}s_{n}}x^{2}dF_{k}\right\vert.
\end{eqnarray*}

\bigskip \noindent Therefore,
\begin{equation}
\frac{1}{s_{n}}\sum_{k=1}^{n}\left\vert E(X_{nk})\right\vert \leq
\varepsilon _{n}^{-1}g(\varepsilon _{n})\rightarrow 0.  \label{expect}
\end{equation}

\bigskip \noindent Based on this,  let us evaluate $s_{nn}/s_{n}$. Notice that for each fixed $n\geq 1$, the $X_{nk}$ are still independent. The technique used in \ref{trunc00} may be summarized as follows : any any measurable function $g(\circ )$ such that $g(0)=0,$

\begin{equation*}
Eg(X_{nk})=\int_{\left\vert X_{k}\right\vert \leq \varepsilon
_{n}s_{n}}g(X_{nk})dP+\int_{\left\vert X_{k}\right\vert >\varepsilon
_{n}s_{n}}g(0)dP=\int_{\left\vert X_{k}\right\vert \leq \varepsilon
_{n}s_{n}}g(X_{nk})dP. \label{esg}
\end{equation*}

\bigskip \noindent By putting these remarks together, we obtain \newpage

\begin{eqnarray*}
1-\frac{s_{nn}^{2}}{s_{n}^{2}} &=&\frac{s_{n}^{2}-s_{nn}^{2}}{s_{n^{2}}} \\
&=&\frac{1}{s_{n^{2}}}\left\{
\sum\limits_{k=1}^{n}EX_{k}^{2}-\sum\limits_{k=1}^{n}E(X_{nk}-E(X_{nk}))^{2}%
\right\}  \\
&=&\frac{1}{s_{n^{2}}}\left\{ \sum\limits_{k=1}^{n}EX_{k}^{2}-\left(
\sum\limits_{k=1}^{n}E(X_{nk}^{2})-E(X_{nk})^{2}\right) \right\}  \\
&=&\frac{1}{s_{n^{2}}}\left\{
\sum\limits_{k=1}^{n}EX_{k}^{2}-\sum\limits_{k=1}^{n}EX_{nk}^{2}+\sum%
\limits_{k=1}^{n}\left( EX_{nk}\right) ^{2}\right\}  \\
\end{eqnarray*}


\begin{eqnarray*}
&=&\frac{1}{s_{n^{2}}}\left\{ \sum\limits_{k=1}^{n}\int
X_{k}^{2}dP-\sum\limits_{k=1}^{n}\int_{\left\vert X_{k}\right\vert \leq
\varepsilon _{n}s_{n}}X_{k}^{2}dP+\sum\limits_{k=1}^{n}\left( EX_{nk}\right)
^{2}\right\}  \\
&=&\frac{1}{s_{n^{2}}}\left\{ \sum\limits_{k=1}^{n}\int_{\left\vert
X_{k}\right\vert >\varepsilon
_{n}s_{n}}X_{k}^{2}dP+\sum\limits_{k=1}^{n}\left( EX_{nk}\right)
^{2}\right\}  \\
&\leq &\frac{1}{s_{n^{2}}}\left\{ \sum\limits_{k=1}^{n}\int_{\left\vert
X_{k}\right\vert >\varepsilon
_{n}s_{n}}X_{k}^{2}dP+\sum\limits_{k=1}^{n}\left( E\left\vert
X_{nk}\right\vert \right) ^{2}\right\} 
\end{eqnarray*}

\bigskip \noindent Finally, we use the simple inequality of real numbers $\left( \sum \left\vert a_{i}\right\vert \right) ^{2}=\sum \left\vert
a_{i}\right\vert ^{2}+\sum_{i\neq j}^{2}\left\vert a_{i}\right\vert \left\vert a_{j}\right\vert \geq \sum \left\vert a_{i}\right\vert ^{2}$ and
conclude from the last inequality that

\begin{eqnarray*}
\left\vert 1-\frac{s_{nn}^{2}}{s_{n}^{2}}\right\vert  &\leq &\frac{1}{%
s_{n^{2}}}\left\{ \sum\limits_{k=1}^{n}\int_{\left\vert X_{k}\right\vert
>\varepsilon _{n}s_{n}}X_{k}^{2}dP+\sum\limits_{k=1}^{n}\left( E\left\vert
X_{nk}\right\vert \right) ^{2}\right\}  \\
&\leq &\frac{1}{s_{n^{2}}}\left\{ \sum\limits_{k=1}^{n}\int_{\left\vert
X_{k}\right\vert >\varepsilon _{n}s_{n}}X_{k}^{2}dP+\left(
\sum\limits_{k=1}^{n}E\left\vert X_{nk}\right\vert \right) ^{2}\right\}\\
&=&g(\varepsilon _{n})+\left( \frac{1}{s_{n}}\sum\limits_{k=1}^{n}E\left\vert
X_{nk}\right\vert \right) .
\end{eqnarray*}

\bigskip \noindent By (\ref{expect}) above, we arrive at

\begin{equation*}
\left\vert 1-\frac{s_{nn}^{2}}{s_{n}^{2}}\right\vert \leq g(\varepsilon
_{n})+\varepsilon _{n}^{-1}g(\varepsilon _{n})\rightarrow 0.
\end{equation*}

\bigskip \noindent It comes that $s_{nn}/s_{n}\rightarrow 1.$ Finally, the proof of this part will be derived  from $S_{nn}/s_{nn}.$ We center the $X_{nk}$
at their expectations. To prove that the sequence of the new summands $T_{nn}/s_{nn}$ converges to $\mathcal{N}(0,1)$, we use Corollary 1 by checking the Lyapounov's condition
\begin{equation*}
\frac{1}{s_{nn}^{3}}\sum_{k=1}^{n}\mathbb{E}\left\vert
X_{nk}-EX_{nk}\right\vert ^{3}\rightarrow 0\text{ as }n\rightarrow \infty .
\end{equation*}

\noindent By

\begin{eqnarray*}
\frac{1}{s_{nn}^{3}}\sum_{k=1}^{n}\mathbb{E}\left\vert
X_{nk}-EX_{nk}\right\vert ^{3} &=&\frac{1}{s_{nn}^{3}}\sum_{k=1}^{n}\mathbb{E%
}\left\vert X_{nk}-EX_{nk}\right\vert \times \mathbb{E}\left\vert
X_{nk}-EX_{nk}\right\vert ^{2} \\
&\leq &\frac{2}{s_{nn}^{3}}\sum_{k=1}^{n}\mathbb{E}\left\vert
X_{nk}\right\vert \times \mathbb{E}\left\vert X_{nk}-EX_{nk}\right\vert ^{2}.
\end{eqnarray*}

\bigskip \noindent Take $g(\cdot )=\left\vert \cdot \right\vert $ in (\ref{esg}) to see again that%
\begin{equation*}
\mathbb{E}\left\vert X_{nk}\right\vert =\int_{\left\vert X_{k}\right\vert
\leq \varepsilon _{n}s_{n}}\left\vert X_{nk}\right\vert dP\leq \varepsilon
_{n}s_{n}.
\end{equation*}

\bigskip \noindent The last two formula yied
\begin{eqnarray*}
\frac{1}{s_{nn}^{3}}\sum_{k=1}^{n}\mathbb{E}\left\vert
X_{nk}-EX_{nk}\right\vert ^{3}&\leq& \frac{2\varepsilon _{n}s_{n}}{s_{nn}^{3}}%
\sum_{k=1}^{n}\mathbb{E}\left\vert X_{nk}-EX_{nk}\right\vert ^{2}\\
&=&\frac{2\varepsilon _{n}s_{n}s_{nn}^{2}}{s_{nn}}=\varepsilon _{n}\frac{2s_{n}}{%
s_{nn}}\rightarrow 0.
\end{eqnarray*}

\noindent It comes by Corrolary \ref{probab_02_indep_sec_03_cor01} that
\begin{equation*}
\frac{T_{nn}}{s_{nn}}=\frac{S_{nn}-\sum_{k=1}^{n}E(X_{nk})}{s_{nn}}%
\rightarrow N(0,1).
\end{equation*}

\bigskip \noindent Since $s_{nn}/s_{n}\rightarrow 1$ and by \ref{expect}%
\begin{equation*}
\left\vert \frac{\sum_{k=1}^{n}E(X_{nk})}{s_{nn}}\right\vert \leq \frac{s_{n}%
}{s_{nn}}\left\{ \frac{1}{s_{n}}\sum_{k=1}^{n}E\left\vert X_{nk}\right\vert
)\right\} \rightarrow 0.
\end{equation*}

\bigskip \noindent We conclude that $S_{nn}/s_{nn}$ converges to $\mathcal{N}(0,1)$.\newline

\noindent \textbf{Proof of : (\ref{WF1})$\Longrightarrow$ (\ref{WF2})}. The convergence to $\mathcal{N}(0,1)$ implies that for any fixed $t\in \mathcal{R}$, we have
\begin{equation}
\prod\limits_{k=1}^{n}f_{k}(u/s_{n})\rightarrow \exp (-u^{2}/2).
\label{limit}
\end{equation}

\bigskip \noindent We are going to use uniform expansions of $\log (1+z)$.
We have 
\begin{equation*}
\lim_{z\rightarrow 0}\left\vert \frac{\log (1+z)-z}{z^{2}}\right\vert =\frac{%
1}{2}
\end{equation*}

\bigskip \noindent this implies
\begin{equation}
\sup_{z\leq u}\left\vert \frac{\log (1+z)-z}{z^{2}}\right\vert =\varepsilon
(u)\rightarrow 1/2\text{ as }u\rightarrow 0.  \label{II_02}
\end{equation}

\bigskip \noindent Now, use the expansion
\begin{equation*}
f_{k}(u/s_{n})=1+\theta _{k}\frac{u^{2}\sigma _{k}^{2}}{2s_{n}^{2}}.
\end{equation*}

\bigskip \noindent This implies that
\begin{equation}
\max_{1\leq k\leq n}\leq \left\vert f_{k}(u/s_{n})-1\right\vert \leq \frac{%
u^{2}}{2}\max_{1\leq k\leq n}\frac{\sigma _{k}^{2}}{s_{n}^{2}}%
=u_{n}\rightarrow 0.  \label{II_01}
\end{equation}

\noindent and next 
\begin{equation*}
\left\vert f_{k}(u/s_{n})-1\right\vert =\theta _{k}^{2}\frac{u^{2}\sigma
_{k}^{2}}{2s_{n}^{2}}\times \frac{u^{2}\sigma _{k}^{2}}{2s_{n}^{2}}\leq %
\left[ \frac{u^{4}}{4}\max_{1\leq k\leq n}\frac{\sigma _{k}^{2}}{s_{n}^{2}}%
\right] \times \frac{\sigma _{k}^{2}}{s_{n}^{2}}.
\end{equation*}

\bigskip \noindent This latter implies

\bigskip 
\begin{equation*}
\sum\limits_{k=1}^{n}\left\vert f_{k}(u/s_{n})-1\right\vert \leq \left[ 
\frac{u^{4}}{4}\max_{1\leq k\leq n}\frac{\sigma _{k}^{2}}{s_{n}^{2}}\right]
=B_{n}(u)\rightarrow 0.
\end{equation*}

\bigskip \noindent By (\ref{II_01}), we see that $\log f_{k}(u/s_{n}$) is uniformly defined in $1\leq k\leq n$ for $n$ large enough and (\ref{limit}) becomes
\begin{equation*}
\sum\limits_{k=1}^{n}\log f_{k}(u/s_{n})\rightarrow -u^{2}/2,
\end{equation*}

\bigskip \noindent that is 
\begin{equation*}
\frac{u^{2}}{2}+\sum\limits_{k=1}^{n}\log f_{k}(u/s_{n})\rightarrow 0.
\end{equation*}

\bigskip \noindent Now using the uniform bound of $\left\vert f_{k}(u/s_{n})-1\right\vert $ by $u_{n}$ to get

\begin{equation*}
\log (f_{k}(u/s_{n}))=f_{k}(u/s_{n})-1+(f_{k}(u/s_{n})-1)^{2}\varepsilon
(u_{n})
\end{equation*}

\bigskip \noindent and then

\begin{eqnarray*}
\frac{u^{2}}{2}+\sum\limits_{k=1}^{n}\log f_{k}(u/s_{n}) &=&\frac{u^{2}}{2}%
+\sum\limits_{k=1}^{n}f_{k}(u/s_{n})-1+(f_{k}(u/s_{n})-1)^{2}\varepsilon
(u_{n}) \\
&=&\left\{ \frac{u^{2}}{2}-\sum\limits_{k=1}^{n}1-f_{k}(u/s_{n})\right\}\\
&+&\left\{ \sum\limits_{k=1}^{n}(f_{k}(u/s_{n})-1)^{2}\right\} \varepsilon (u_{n}),
\end{eqnarray*}

\bigskip \noindent with
\begin{equation*}
\left\vert \left\{ \sum\limits_{k=1}^{n}(f_{k}(u/s_{n})-1)^{2}\right\}
\varepsilon (u_{n})\right\vert \leq B_{n}(u)\left\vert \varepsilon
(u_{n})\right\vert =o(1).
\end{equation*}

\bigskip \noindent We arrive at%
\begin{equation*}
\frac{u^{2}}{2}=\sum\limits_{k=1}^{n}1-f_{k}(u/s_{n})+o(1).
\end{equation*}

\bigskip \noindent If we take the real parts, we have for any fixed $%
\varepsilon >0$,

\begin{eqnarray*}
\frac{u^{2}}{2} &=&\sum\limits_{k=1}^{n}\int (1-\cos \frac{ux}{s_{n}}%
)dF_{k}(x)+o(1) \\
&=&\sum\limits_{k=1}^{n}\int_{\left\vert x\right\vert <\varepsilon
s_{n}}(1-\cos \frac{ux}{s_{n}})dF_{k}(x)\\
&+&\sum\limits_{k=1}^{n}\int_{\left\vert x\right\vert \geq \varepsilon s_{n}}(1-\cos \frac{ux}{s_{n}})dF_{k}(x)+o(1),
\end{eqnarray*}

\bigskip \noindent that is
\begin{equation*}
\frac{u^{2}}{2}-\sum\limits_{k=1}^{n}\int_{\left\vert x\right\vert
<\varepsilon s_{n}}(1-\cos \frac{ux}{s_{n}})dF_{k}(x)=\sum\limits_{k=1}^{n}%
\int_{\left\vert x\right\vert \geq \varepsilon s_{n}}(1-\cos \frac{ux}{s_{n}}%
)dF_{k}(x)+o(1),
\end{equation*}

\bigskip \noindent We have by Fact 2 below that $\sqrt{2(1-\cos a)}\leq
2\left\vert a/2\right\vert ^{\delta }$ for all $\delta ,0<\delta \leq 1.$
Apply this for $\delta =1$ to have%
\begin{eqnarray*}
&&\sum\limits_{k=1}^{n}\int_{\left\vert x\right\vert <\varepsilon
s_{n}}(1-\cos \frac{ux}{s_{n}})dF_{k}(x)\\
 &\leq &\frac{u^{2}}{2s_{n}^{2}}%
\sum\limits_{k=1}^{n}\int_{\left\vert x\right\vert <\varepsilon
s_{n}}x^{2}dF_{k}(x) \\
&=&\frac{u^{2}}{2s_{n}^{2}}(\sum\limits_{k=1}^{n}\int
x^{2}dF_{k}(x)-\sum\limits_{k=1}^{n}\int_{\left\vert x\right\vert \geq
\varepsilon s_{n}}x^{2}dF_{k}(x)) \\
&=&\frac{u^{2}}{2s_{n}^{2}}(s_{n}^{2}-\sum\limits_{k=1}^{n}\int_{\left\vert
x\right\vert \geq \varepsilon s_{n}}x^{2}dF_{k}(x))=\frac{u^{2}}{2}%
(1-g_{n}(\varepsilon )).
\end{eqnarray*}

\bigskip \noindent On the other hand

\begin{eqnarray*}
\sum\limits_{k=1}^{n}\int_{\left\vert x\right\vert \geq \varepsilon
s_{n}}(1-\cos \frac{ux}{s_{n}})dF_{k}(x) &\leq
&2\sum\limits_{k=1}^{n}\int_{\left\vert x\right\vert \geq \varepsilon
s_{n}}dF_{k}(x) \\
&=&2\sum\limits_{k=1}^{n}\int_{\left\vert x\right\vert \geq \varepsilon
s_{n}}\left\{ \frac{1}{x^{2}}\right\} x^{2}dF_{k}(x) \\
&\leq &\frac{2}{\varepsilon _{{}}^{2}s_{n}^{2}}\sum\limits_{k=1}^{n}\int_{\left\vert x\right\vert \geq \varepsilon s_{n}}x^{2}dF_{k}(x)\leq \frac{2}{\varepsilon ^{2}}.
\end{eqnarray*}

\bigskip \noindent By putting all this together, we have

\begin{equation*}
\frac{u^{2}}{2}\leq \frac{u^{2}}{2}(1-g_{n}(\varepsilon ))+\frac{2}{%
\varepsilon ^{2}}+o(1)
\end{equation*}

\bigskip \noindent which leads 
\begin{equation*}
\frac{u^{2}}{2}g_{n}(\varepsilon )\leq \frac{2}{\varepsilon ^{2}}+o(1)
\end{equation*}

\bigskip \noindent which in turns implies

\begin{equation*}
g_{n}(\varepsilon )\leq \frac{2}{u^{2}}(\frac{2}{\varepsilon ^{2}}+o(1)).
\end{equation*}

\bigskip \noindent By letting first $n\rightarrow +\infty $ and secondly $%
u\rightarrow 0,$ we get 
\begin{equation*}
g_{n}(\varepsilon )\rightarrow 0.
\end{equation*}

\bigskip \noindent This concludes the proof. $\square$\\

\newpage
\noindent \textbf{(C) APPENDIX : TOOLS, FACTS AND LEMMAS}

\noindent \textbf{1 - A useful development for the characteristic function}.\\

\bigskip \noindent Consider the characteristic function associated with the real probability distribution function $F$ that is%
\begin{equation*}
\mathbb{R}\ni x\hookrightarrow f(x)=\int e^{itx}dF(x)
\end{equation*}

\bigskip \noindent Suppose that the $n^{th}$ moment exists, that is%
\begin{equation*}
m_{n}=\int x^{n}dF(x).
\end{equation*}

\bigskip \noindent In the following, we also denote%
\begin{equation*}
\mu _{n}=\int \left\vert x\right\vert ^{n}dF(x)
\end{equation*}

\bigskip 
\begin{lemma} \label{probab_02_indep_sec_03_lem01} Let $0<\delta \leq 1.$ If $\mu _{n+2}$ is finite, then we have the following expansion

\begin{equation}
f(u)=1+\sum_{k=1}^{n}\frac{(iu)^{k}mk}{k!}+\theta 2^{1-\delta }\mu
^{n+\delta }\frac{\left\vert u\right\vert ^{n+\delta }}{(1+\delta )(2+\delta
)...(n+\delta )},\left\vert \theta \right\vert \leq 1.  \label{expan-ch}
\end{equation}
\end{lemma}

\bigskip \noindent \textbf{Proof of Lemma \ref{probab_02_indep_sec_03_lem01}}. By using the Lebesgue Dominated
Theorem, we get the $f$ is $n$-times differentiable and the $k$-th
derivative is%
\begin{equation}
f^{(k)}(0)=i^{k}m_{k}=\int x^{k}dF(x),1\leq k\leq n.  \label{moment1}
\end{equation}

\noindent We may use the Taylor-Mac-Laurin formula \ expansion%
\begin{equation*}
f(u)=1+\sum_{k=1}^{n-1}\frac{(iu)^{k}mk}{k!}+\int_{0}^{u}\frac{(u-x)^{n-1}}{%
n!}f^{(n)}(x)dx.
\end{equation*}

\bigskip \noindent We are going to handle $\rho _{n}(u)=\int_{0}^{u}\frac{x^{n-1}}{n!}f^{(n)}(x)dx$. Let us make  the change variable $t=x/u$ and use \ref%
{moment1} to get

\begin{eqnarray*}
\rho _{n}(u) &=&u^{n}\int_{0}^{1}\frac{(1-t)^{n-1}}{n!}%
f^{(n)}(tu)dt\\
&=&(iu)^{n}\int_{0}^{1}\int \frac{(1-t)^{n-1}}{(n-1)!}%
x^{n}e^{itux}dF(x)dt  \label{rhon} \\
&=&(iu)^{n}\int_{0}^{1}\int \frac{(1-t)^{n-1}}{(n-1)!}%
x^{n}(e^{itux}-1+1)dF(x)dt  \notag \\
&=&(iu)^{n}\int_{0}^{1}\int \frac{(1-t)^{n-1}}{(n-1)!}x^{n}dF(x)dt\\
&+&(iu)^{n}\int_{0}^{1}\int \frac{(1-t)^{n-1}}{(n-1)!}x^{n}(e^{itux}-1)dF(x)dt  \notag
\end{eqnarray*}

\bigskip \noindent The first term is
\begin{equation*}
\rho _{n}(1,u)=(iu)^{n}\int_{0}^{1}\frac{(1-t)^{n-1}}{(n-1)!}dt\int
x^{n}dF(x)=(iu)^{n}m_{n}\left[ -\frac{(1-t)^{n}}{n!}\right] _{t=0}^{t=1}=%
\frac{(iu)^{n}m_{n}}{n!}.
\end{equation*}

\bigskip \noindent To handle the second term, we remark that,%
\begin{equation*}
\left\vert e^{ia}-1\right\vert =\sqrt{2-1-\cos (a)}=2\left\vert \sin
(a/2)\right\vert .
\end{equation*}%
Let $0<\delta \leq 1.$ If $\left\vert a/2\right\vert \geq 1,$ we have 
\begin{equation*}
\left\vert e^{ia}-1\right\vert =2\left\vert \sin (a/2)\right\vert \leq
\left\vert a\right\vert \leq 2\left\vert a/2\right\vert ^{\delta }
\end{equation*}

\bigskip \noindent by the decreasingness in $\delta $\ of the function $%
\left\vert a/2\right\vert ^{\delta }$. If $\left\vert a/2\right\vert \leq 1,$
we get by Fact 1 below that also $\left\vert e^{ia}-1\right\vert
=2\left\vert \sin (a/2)\right\vert \leq 2\left\vert a/2\right\vert ^{\delta
} $. We have for all $a\in R,$ for all $0\leq \delta \leq 1,$ 
\begin{equation*}
\left\vert e^{ia}-1\right\vert \leq 2\left\vert a/2\right\vert ^{\delta }
\end{equation*}

\bigskip \noindent Applying this to (\ref{rhon}) yields 
\begin{eqnarray*}
\left\vert \rho _{n}(2,u)\right\vert  &\leq &\left\vert u\right\vert
\int_{0}^{1}\int \frac{(1-t)^{n-1}}{(n-1)!}\left\vert x\right\vert
^{n}\left\vert e^{itux}-1\right\vert dF(x)dt \\
&\leq &2^{1-\delta }\left\vert u\right\vert ^{n+\delta }\int_{0}^{1}\frac{%
(1-t)^{n-1}t^{\delta }}{(n-1)!}dt\int \left\vert x\right\vert ^{n+\delta
}dF(x) \\
&\leq &2^{1-\delta }\left\vert u\right\vert ^{n+\delta }\mu ^{n+\delta
}\int_{0}^{1}\frac{(1-t)^{n-1}t^{\delta }}{(n-1)!}dt.
\end{eqnarray*}

\bigskip \noindent Since by Fact 2 below, 
\begin{equation*}
\int_{0}^{1}\frac{(1-t)^{n-1}t^{\delta }}{(n-1)!}dt=\frac{1}{(1+\delta
)(2+\delta )...(n+\delta )},
\end{equation*}

\bigskip \noindent we get 
\begin{equation*}
\rho _{n}(2,u)=\theta 2^{1-\delta }\mu ^{n+\delta }\frac{\left\vert
u\right\vert ^{n+\delta }}{(1+\delta )(2+\delta )...(n+\delta )},
\end{equation*}

\bigskip \noindent with $\left\vert \theta \right\vert \leq 1.$ By getting
together all these pieces, we get (\ref{expan-ch}). This concludes the proof
of Lemma \ref{probab_02_indep_sec_03_lem01}.\newline

\bigskip

\bigskip \noindent \textbf{FACT 1}. For any $a\in \mathbb{R},$ 
\begin{equation*}
\left\vert e^{ia}-1\right\vert =\sqrt{2(1-\cos a)}\leq 2\left\vert \sin
(a/2)\right\vert \leq 2\left\vert a/2\right\vert ^{\delta }.
\end{equation*}

\bigskip \noindent This is easy for $\left\vert a/2\right\vert >1.$ Indeed
for $\delta >0,\left\vert a/2\right\vert ^{\delta }>0$ and%
\begin{equation*}
2\left\vert \sin (a/2)\right\vert \leq 2\leq 2\left\vert a/2\right\vert
^{\delta }
\end{equation*}

\bigskip \noindent Now for $\left\vert a/2\right\vert >1,$ we have the
expansion%
\begin{eqnarray*}
2(1-\cos a) &=&a^{2}-\sum\limits_{k=2}^{\infty }(-1)^{2}\frac{a^{2k}}{(2k)!}%
=x^{2}-2\sum\limits_{k\geq 2,k\text{ }even}^{\infty }\frac{a^{2k}}{(2k)!}-%
\frac{a^{2(k+1)}}{(2(k+1))!} \\
&=&a^{2}-2x^{2(k+1)}\sum\limits_{k\geq 2,k\text{ }even}^{\infty }\frac{1}{%
(2k)!}\left\{ \frac{1}{a^{2}}-\frac{1}{(2k+1)((2k+2)...(2k+k)}\right\} .
\end{eqnarray*}

\bigskip \noindent For each $k\geq 2,$ for $\left\vert a/2\right\vert <1,$%
\begin{equation*}
\left\{ \frac{1}{a^{2}}-\frac{1}{(2k+1)((2k+2)...(2k+k)}\right\} \geq
\left\{ \frac{1}{4}-\frac{1}{(2k+1)((2k+2)...(2k+k)}\right\} \geq 0.
\end{equation*}

\bigskip \noindent Hence 
\begin{equation*}
2(1-\cos a)\leq a^{2}.
\end{equation*}

\bigskip \noindent But for $\left\vert a/2\right\vert ,$ the function $%
\delta \hookrightarrow \left\vert a/2\right\vert ^{\delta }$ is
non-increasing $\delta ,0\leq \delta \leq 1$. Then%
\begin{equation*}
\sqrt{2(1-\cos a)}\leq \left\vert a\right\vert =2\left\vert a/2\right\vert
^{1}\leq 2\left\vert a/2\right\vert ^{\delta }.
\end{equation*}

\bigskip \noindent \textbf{FACT 2}. For any $1<\delta \leq 1$, for any $%
n\geq 1$ 
\begin{equation*}
\int_{0}^{1}\frac{(1-t)^{n-1}t^{\delta }}{(n-1)!}dt=\frac{1}{(1+\delta
)(2+\delta )...(n+\delta )}.
\end{equation*}

\bigskip \noindent \textbf{Proof}. By integrating by parts, we get%
\begin{equation*}
\int_{0}^{1}\frac{(1-t)^{n-1}t^{\delta }}{(n-1)!}dt=\frac{1}{\delta +1}\left[
\frac{(1-t)^{n-1}t^{\delta }}{(n-1)!}\right] _{t=0}^{t=1}+\frac{1}{\delta +1}%
\int_{0}^{1}\frac{(1-t)^{n-2}t^{\delta +1}}{(n-2)!}dt,
\end{equation*}

\bigskip \noindent that is 
\begin{equation*}
\int_{0}^{1}\frac{(1-t)^{n-1}t^{\delta }}{(n-1)!}dt=\frac{1}{\delta +1}%
\int_{0}^{1}\frac{(1-t)^{n-2}t^{\delta +1}}{(n-2)!}dt.
\end{equation*}

\bigskip \noindent From there, we easily get by induction that, \ for $1\leq \ell \leq n-1$,

\begin{equation*}
\int_{0}^{1}\frac{(1-t)^{n-1}t^{\delta }}{(n-1)!}dt=\frac{1}{\left( \delta
+1\right) \left( \delta +2\right) ...\left( \delta +\ell \right) }%
\int_{0}^{1}\frac{(1-t)^{n-\ell -1}t^{\delta +\ell }}{(n-2)!}dt.
\end{equation*}%
For $\ell =n-1,$ we have%
\begin{eqnarray*}
\int_{0}^{1}\frac{(1-t)^{n-1}t^{\delta }}{(n-1)!}dt &=&\frac{1}{\left(
\delta +1\right) \left( \delta +2\right) ...\left( \delta +n-1\right) }%
\int_{0}^{1}t^{\delta +n-1}dt \\
&=&\frac{1}{\left( \delta +1\right) \left( \delta +2\right) ...\left( \delta
+n\right) }.
\end{eqnarray*}

\bigskip \noindent This finishes the proof. $\square$\\

\bigskip 
\begin{lemma} \label{probab_02_indep_sec_03_lem02} Let $Y$ a random variable with $r_{0}$-th finite moment, $r_{0}>0.$ Then the function $g(x)=\log \mathbb{E}\left\vert
Y\right\vert ^{x}$, $0\leq x\leq r_{0}$, is convex.
\end{lemma}

\bigskip \noindent \textbf{Proof of Lemma \ref{probab_02_indep_sec_03_lem02}}. Let $0\leq r_{1}<r_{2}\leq r_{0}.$ Use
the Cauchy-Scharwz inequality to $\left\vert Y\right\vert ^{(r_{1}+r_{2})/2}$
and $\left\vert Y\right\vert ^{(r_{2}-r_{1})/2}$ to have

\begin{equation*}
\left( E\left\vert Y\right\vert ^{r_{1}}\right) ^{2}\leq E\left\vert
Y\right\vert ^{(r_{1}+r_{2})}\times E\left\vert Y\right\vert ^{(r_{2}-r_{1})}
\end{equation*}

\bigskip \noindent which implies 
\begin{equation*}
2\log E\left\vert Y\right\vert ^{r_{1}}\leq \log \mathbb{E}\left\vert Y\right\vert
^{(r_{1}+r_{2})}+\log \mathbb{E}\left\vert Y\right\vert ^{(r_{2}-r_{1})}
\end{equation*}

\bigskip \noindent that is, since $g$ is continuous,
 
\begin{equation}
g(r_{1})\leq \frac{1}{2}(g(r_{1}+r_{2})+g(r_{2}-r_{1})).  \label{convex1}
\end{equation}

\bigskip \noindent Now, set $x=r_{1}+r_{2}$ and $y=r_{2}-r_{1}$ and (\ref{convex1}) becomes
\begin{equation}
g(\frac{x+y}{2})\leq \frac{1}{2}(g(x)+g(y))  \label{convex2}
\end{equation}

\bigskip \noindent for $0\leq x\leq r_{0}$. Now, the Dominated Convergence
Theorem, the fonction $g(\cdot )$ is continuous. So (\ref{convex2}) implies
the convexity of $g(\cdot)$. $\square$\\

\bigskip 
\begin{lemma} \label{probab_02_indep_sec_03_lem03} Let $X_{1},$ $X_{2},$ ... a
sequence of real and independent random variables centered at expectations,
with finite $(n+\delta )-$moment, $\delta >0.$ Put for each $n\geq 1,$ $%
S_{n}=X_{1}+...+X_{n}$ et $s_{n}^{2}=$ $\mathbb{E}X_{1}^{2}+\mathbb{E}%
X_{2}^{2}+...+\mathbb{E}X_{n}^{2}.$ We denote $\sigma _{k}^{2}=\mathbb{E}%
X_{k}^{2},k\geq 1$ et $F_{k}$\ denotes the probability distribution function
of $X_{k}$. If $\delta >1,$ then any fixed $n\geq 1,$

\begin{equation}
\frac{1}{s_{n}^{2+\delta }}\sum_{k=1}^{n}\mathbb{E}\left\vert
X_{k}\right\vert ^{2+\delta }\leq \left( \frac{1}{s_{n}^{3}}\sum_{k=1}^{n}%
\mathbb{E}\left\vert X_{k}\right\vert ^{3}\right) ^{(\delta -2)/\delta }.
\end{equation}
\end{lemma}

\bigskip \noindent \textbf{Proof of Lemma \ref{probab_02_indep_sec_03_lem03}}. Let $n\geq 1$ be fixed. Let $(\pi
_{1},...,\pi _{n})$ following a multinomial law of $n$ issues having all the
probability $n$ or occuring but only on repeatition. This means that only
one of the $\pi _{k}^{\prime }s$ is one, the remaining being zero. Sey 
\begin{equation*}
Y=\sum_{k=1}^{n}\pi _{k}X_{k}.
\end{equation*}

\bigskip \noindent The meaning of this expression is the following :
\begin{equation*}
Y=X_{k}\text{ on (}\pi _{k}=1).
\end{equation*}

\bigskip \noindent So we have, for $r\geq 0.$%
\begin{equation*}
\left\vert Y\right\vert ^{r}=\sum_{k=1}^{n}\pi _{k}\left\vert
X_{k}\right\vert ^{r}.
\end{equation*}

\bigskip \noindent Hence 
\begin{eqnarray*}
E\left\vert Y\right\vert ^{r} &=&E\sum_{j=1}^{n}\pi _{j}\left\vert
X_{j}\right\vert ^{r}=\sum_{k=1}^{n}P(\pi _{k}=1)E\left( \sum_{j=1}^{n}\pi
_{j}\left\vert X_{j}\right\vert ^{r}\text{ \ \ }\left\vert \text{ }\pi
_{k}=1\right. \right) \\
&=&\frac{1}{n}\sum_{k=1}^{n}E\left\vert X_{k}\right\vert ^{r}.
\end{eqnarray*}

\bigskip \noindent Use now the convexity of $g(r)=\log E\left\vert
Y\right\vert ^{r}$ for $\delta >1$ like that : 
\begin{equation*}
\frac{\delta -1}{\delta }\times 2+\frac{1}{\delta }\times (2+\delta )=3
\end{equation*}

\bigskip \noindent and convexity implies 
\begin{equation*}
g(\frac{\delta -1}{\delta }\times 2+\frac{1}{\delta }\times (2+\delta ))\leq 
\frac{\delta -1}{\delta }g(2)+\frac{1}{\delta }g(2+\delta )).
\end{equation*}%
This implies%
\begin{equation*}
\delta \log E\left\vert Y\right\vert ^{3}\leq (\delta -1)\log E\left\vert
Y\right\vert ^{2}+\log E\left\vert Y\right\vert ^{2+\delta }
\end{equation*}

\bigskip \noindent and by taking exponentials, we get%
\begin{equation*}
\left( E\left\vert Y\right\vert ^{3}\right) ^{\delta }\leq \left(
E\left\vert Y\right\vert ^{2}\right) ^{\delta -1}E\left\vert Y\right\vert
^{2+\delta }\Longrightarrow E\left\vert Y\right\vert ^{3}\leq \left(
E\left\vert Y\right\vert ^{2}\right) ^{\delta -1}\left( E\left\vert
Y\right\vert ^{2+\delta }\right) ^{1/\delta }.
\end{equation*}

\bigskip \noindent Replacing by the values of $E\left\vert Y\right\vert ^{r},
$ we get%
\begin{equation*}
\frac{1}{n}\sum_{k=1}^{n}E\left\vert X_{k}\right\vert ^{3}\leq \frac{1}{%
n^{(\delta -1)/\delta }}s_{n}^{2(\delta -1)/\delta }\left( \frac{1}{n}%
\sum_{k=1}^{n}E\left\vert X_{k}\right\vert ^{2+\delta }\right) ^{1/\delta }.
\end{equation*}

\bigskip \noindent From there, easy computations lead to%
\begin{equation*}
\frac{1}{s_{n}^{3}}\sum_{k=1}^{n}E\left\vert X_{k}\right\vert ^{3}\leq
\left( \frac{1}{s_{n}^{2+\delta }}\sum_{k=1}^{n}E\left\vert X_{k}\right\vert
^{2+\delta }\right) ^{1/\delta }.
\end{equation*}

\bigskip \noindent \textbf{LEMMA 4}. Let $\delta >0$ and let $X$ be a real
random variable such that $\left\vert X\right\vert ^{2+\delta }$ is
integrable. Then
\begin{equation*}
\left( EX^{2}\right) ^{(2+\delta )/2}\leq E\left\vert X\right\vert
^{2+\delta }.
\end{equation*}

\bigskip \noindent \textbf{PROOF}. Use Lemma 2 and the convexity of $%
g(x)=\log E\left\vert X\right\vert ^{x},0<x\leq 2+\delta $ to the convex
combination 
\begin{equation*}
2=\frac{2}{2+\delta }\times \left( 2+\delta \right) +\frac{2}{2+\delta }%
\times 0
\end{equation*}

\bigskip \noindent to get 
\begin{equation*}
g(2)\leq \frac{2}{2+\delta }g(2+\delta )+\frac{2}{2+\delta }g(0).
\end{equation*}%
Since $g(0)=0,$ we have
\begin{equation*}
\log E\left\vert X\right\vert ^{2}\leq \frac{2}{2+\delta }\log E\left\vert
X\right\vert ^{2+\delta },
\end{equation*}

\bigskip \noindent which gives the desired results upon taking the exponentials.
 

%% file: proba_02_07_be_ang.tex
\section{Berry-Essen approximation} \label{probab_02_indep_sec_04}

\noindent Once the central theorem holds, the convergence of the distribution functions of $S_{n}/s_{n}$, denoted $F_{n}$, $n \geq 1$, to that of a standard Gaussian
random variable denoted by $G $ holds uniformly, by a known result of weak convergence (See \cite{ips-wcrv-ang}, chapter 4, Fact 5), that is

\begin{equation*}
\sup_{x \in \mathbb{R}} |F_{n}(x)-G(x)| \rightarrow , \text{ as } n \rightarrow +\infty.
\end{equation*}

\noindent The Berry-Essen inequality is the most important result on the rate of convergence of $F_{n}$ to $G$. Here is a classical form of it.%
\newline

\subsection{Statement of the Berry-Essen Inequality}

\begin{theorem} \label{probab_02_indep_sec_04_th01}
(Berry-Essen) Let $X_{1},X_{2},...$ be independent random variables with zero mean and with partial sums $\{S_{n},n\geq 1\}$. Suppose that $\gamma
_{k}^{3}=\mathbb{E}\left\vert X_{k}\right\vert ^{3}<+\infty $ for all $k\geq 1,$ and set $\sigma _{k}^{2}=Var(X_{k}),$ $s_{n}^{2}=\sum_{1\leq j\leq
k}\sigma _{j}^{2}$ and $\beta _{n}^{3}=\sum_{1\leq j\leq k}\gamma _{j}^{3}$. Then 
\begin{equation*}
\sup_{x\in R}\left\vert \mathbb{P}\left( \frac{S_{n}}{s_{n}}\leq x\right) - \mathbb{P}(N(0,1)\leq x)\right\vert \leq C\frac{\beta _{n}^{3}}{s_{n}^{3}}.
\end{equation*}
\end{theorem}

\bigskip \noindent \textbf{Remarks} This result may be extended to some dependent data. Generally, one seeks to get a Berry-Essen type results each time a Central limit Theorem is obtained.\newline

\noindent The value of $C$ may be of interest and one seeks to have it the lowest possible. In the proof below, $C$ will be equal to $36$. 

\bigskip \noindent \textbf{PROOF} The proof is very technical. But, it is important to do it at least one time, since, it may give ideas when no longer prevails the independence.\newline

\noindent The proof itself depends on two interesting lemmas. We suggest to the reader who wants to develop an expertise in this field, to do the
following.\newline

\noindent 1) The reader who wishes to master this very technical proof is recommended to read the statement and the proof of the Essen Lemma \ref{essen.lem}. This lemma gives the important formula \eqref{essen.lem.form}. It is based on the inversion formula that expresses the density probability function with respect to the characteristic function. It also uses a characterization the supremum of bounded and right-continuous with left-limits (rcll) of real-valued functions vanishing at $\pm \infty$ given in Lemma \ref{essen.lem4}.\newline

\noindent 2) Next, read the statement of Lemma \ref{essen.approx} which gives the approximation of the characteristic function of $S_{n}/s_{n}$ to that of a standard normal random variable which is $exp(-t^{2}/2)$. The proof of this Lemma uses a special expansion of the characteristic function in the neighborhood of zero given in Lemma \ref{essen.dev.charact}.\newline

\noindent From these two points, the proof of the Theorem of Berry-Essen comes out naturally in the following lines by plugging the results of Lemma \ref{essen.approx} in the formula \eqref{essen.lem.form} of Lemma \ref{essen.lem}. And we say :\newline

\noindent By Lemma \ref{essen.lem.form},

\begin{equation}
\sup_{x}\left\vert F_{S_{n}/s_{n}}(x)-G (x)\right\vert \leq \frac{1}{\pi }\int_{-T}^{T}\left\vert \frac{\psi _{S_{n}/s_{n}}(t)-\exp (-t^{2}/2)}{t}\right\vert dt+24A/(\pi T),  \label{essen.berry1}
\end{equation}

\noindent where $A$ is an upper bound of the derivative the standard gaussian distribution function $G $ whose infimum is $1/\sqrt{2\pi }$. Take $A=1/\sqrt{2\pi }$ and $T=T_{n}=s_{n}^{3}/(4\beta _{n}^{3})$. We use Formula (\ref{essen.berry1}) and the following inequality
\begin{equation*}
\left\Vert \psi _{S_{n}/s_{n}}(t)-\exp (-t^{2}/2)\right\Vert \leq 16\exp
(-t^{2}/2)\frac{\beta _{n}^{3}\left\vert t\right\vert ^{3}}{s_{n}^{3}}.
\end{equation*}

\noindent to grap on
\begin{eqnarray*}
\sup_{x}\left\vert F_{S_{n}/s_{n}}(x)-G (x)\right\vert &\leq &\frac{16}{%
\pi }\frac{\beta _{n}^{3}}{s_{n}^{3}}\int_{-s_{n}^{3}/(4\beta
_{n}^{3})}^{s_{n}^{3}/(4\beta _{n}^{3})}t^{2}\exp (-t^{2}/2)dt+\frac{24(1/)}{%
\pi s_{n}^{3}/(4\beta _{n}^{3})} \\
&\leq &\frac{16}{\pi }\frac{\beta _{n}^{3}}{s_{n}^{3}}\int_{s_{n}^{3}/(4%
\beta _{n}^{3})}^{s_{n}^{3}/(4\beta _{n}^{3})}t^{2}\exp (-t^{2}/2)dt+\frac{%
96\beta _{n}^{3}}{\pi \sqrt{2\pi }s_{n}^{3}}.
\end{eqnarray*}

\noindent The integral $\int_{s_{n}^{3}/(4\beta_{n}^{3})}^{s_{n}^{3}/(4\beta _{n}^{3})}t^{2}\exp (-t^{2}/2)dt$ is bounded by
\begin{equation*}
\int_{s_{n}^{3}/(4\beta _{n}^{3})}^{s_{n}^{3}/(4\beta _{n}^{3})}t^{2}\exp
(-t^{2}/2)dt=\frac{3}{2}\sqrt{2\pi }.
\end{equation*}

\noindent We get 
\begin{eqnarray*}
\sup_{x}\left\vert F_{S_{n}/s_{n}}(x)-G (x)\right\vert &\leq &\frac{24\sqrt{3}}{\sqrt{\pi }}\frac{\beta _{n}^{3}}{s_{n}^{3}}+\frac{96\beta _{n}^{3}}{\pi \sqrt{2\pi }s_{n}^{3}} \\
&\leq &\left( \frac{24\sqrt{3}}{\sqrt{\pi }}+\frac{96}{\pi \sqrt{2\pi }}\right) \frac{\beta _{n}^{3}}{s_{n}^{3}} \\
&\leq &36\frac{\beta _{n}^{3}}{s_{n}^{3}}.
\end{eqnarray*}

\noindent This concludes the proof.

\subsection{Tools, Facts and Lemmas}

\begin{lemma}
\label{essen.lem.1} Define the triangle probability density function $pdf$, with parameter $T$ as following%
\begin{equation}
f_{tri}(x)=\frac{1}{T}(1-\frac{\left\vert x\right\vert }{T})1_{(\left\vert
x\right\vert \leq T)}.  \label{tripdf}
\end{equation}

\noindent (i) Then its characteristic function is 
\begin{equation*}
\Phi _{tri(T)}(t)=\sin ^{2}(tT/2)/(tT)^{2}.
\end{equation*}

\noindent (ii) The function 
\begin{equation*}
g(x)=\frac{1-\cos xT}{\pi x^{2}T}, x \in \mathbb{R},
\end{equation*}

\noindent defines a density distribution function and its characteristic function is $1-\left\vert t\right\vert /T$ .
\end{lemma}

\bigskip

\noindent \textbf{Proof}. We have 
\begin{eqnarray*}
\Phi _{tri(T)}(t) &=&\frac{1}{T}\int_{-T}^{T}e^{itx}(1-\frac{\left\vert
x\right\vert }{T})dx \\
&=&\frac{1}{T}\left\{ \int_{-T}^{0}e^{itx}(1+\frac{x}{T})dx+\int_{0}^{T}e^{itx}(1-\frac{x}{T})dx\right\} .
\end{eqnarray*}

\noindent Next, we have  
\begin{equation*}
\int_{-T}^{0}e^{itx}(1+\frac{x}{T})dx=\left[ \frac{e^{itx}}{it}\right]_{-T}^{0}+\frac{1}{T}\int_{-T}^{0}xe^{itx}dx.
\end{equation*}

\noindent By integrating by parts, we get
\begin{eqnarray*}
\int_{-T}^{0}e^{itx}(1+\frac{x}{T})dx &=&\left[ \frac{e^{itx}}{it}\right]_{-T}^{0}+\frac{1}{T}\left[ \frac{xe^{itx}}{it}\right] _{-T}^{0}+\frac{1}{itT}\int_{-T}^{0}e^{itx}dx \\
&=&\left[ \frac{e^{itx}}{it}\right] _{-T}^{0}+\frac{1}{T}\left[ \frac{xe^{itx}}{it}\right] _{-T}^{0}-\frac{1}{itT}\left[ \frac{e^{itx}}{it}\right]_{-T}^{0} \\
&=&\frac{1}{it}(1-e^{-itT})+\frac{1}{it}e^{-itT}+\frac{1}{t^{2}T}(1-e^{-itT}).
\end{eqnarray*}

\noindent Likewise, we get

\begin{eqnarray*}
\int_{0}^{T}e^{itx}(1-\frac{x}{T})dx &=&\left[ \frac{e^{itx}}{it}\right]_{0}^{T}-\frac{1}{T}\int_{0}^{T}xe^{itx}dx. \\
&=&\left[ \frac{e^{itx}}{it}\right] _{0}^{T}-\frac{1}{T}\left[ \frac{xe^{itx}}{it}\right] _{0}^{T}+\frac{1}{itT}\int_{0}^{T}e^{itx}dx \\
&=&\left[ \frac{e^{itx}}{it}\right] _{0}^{T}-\frac{1}{T}\left[ \frac{xe^{itx}}{it}\right] _{0}^{T}+\frac{1}{itT}\left[ \frac{e^{itx}}{it}\right] _{0}^{T}\\
&=&\frac{1}{it}\left(e^{itT}-1)-\frac{1}{it}e^{itT}-\frac{1}{t^{2}T}(e^{itT}-1\right).
\end{eqnarray*}

\noindent By putting all this together, and by adding term by term, we get
\begin{eqnarray*}
\Phi _{tri(T)}(t) &=&\frac{1}{T}\left\{ \frac{2\sin tT}{t}-\frac{2\sin tT}{t}-\frac{2\cos tT-2}{t^{2}T}\right\} \\
&=&\frac{2(1-\cos tT)}{t^{2}T^{2}} \\
&=&\frac{\sin ^{2}tT/2}{t^{2}T^{2}}.
\end{eqnarray*}

\bigskip \noindent Remark that $\Phi _{tri(T)}(t)$ is well defined for $t=0$. From now, we may use the inversion theorem for an absolutely continuous distribution function :%
\begin{equation*}
f_{tri(T)}(t)=\frac{1}{2\pi }\int e^{-itx}\Phi _{tri(T)}(x)dx.
\end{equation*}

\noindent Then for $\left\vert t\right\vert \leq T,$%
\begin{equation*}
\frac{1}{T}\left(1-\frac{\left\vert t\right\vert }{T}\right)=\frac{1}{2\pi }\int
e^{-itx}\frac{\sin ^{2}xT/2}{x^{2}T^{2}}dx,
\end{equation*}

\noindent which gives

\begin{eqnarray*}
1-\frac{\left\vert t\right\vert }{T} &=&\frac{1}{2\pi }\int_{\mathbb{R}} e^{-itx}\frac{\sin ^{2}xT/2}{2\pi x^{2}T^{2}}dx \\
&=&\frac{1}{\pi }\int_{\mathbb{R}} e^{-itx}\frac{1-\cos xT}{x^{2}T^{2}}dx \\
&=&\int_{\mathbb{R}} e^{itx}\frac{1-\cos xT}{\pi x^{2}T^{2}}dx.
\end{eqnarray*}

\bigskip \noindent Taking $t=0$ in that formula proves that 

\begin{equation*}
\frac{1-\cos xT}{\pi x^{2}T^{2}},x\in \mathbb{R}
\end{equation*}

\noindent is a density probability on $\mathbb{R}$, and its characteristic function is $1-\left\vert t\right\vert /T$. This gives

\bigskip

\begin{lemma}
\label{essen.lem.2} The following function 
\begin{equation*}
\frac{1-\cos xT}{\pi x^{2}T^{2}},x\in \mathbb{R}
\end{equation*}

\noindent is a probability density function with characteristic function $1-\left\vert t\right\vert /T.$
\end{lemma}

\bigskip \noindent The following lemma uses the inverse formulas in Proposition \ref{proba02_rv.inv01} (see section Part V, Section \ref{proba_02_rv_sec_06}, Chapter \ref{proba_02_rv}).   

\begin{lemma} \label{essen.lem} Let $U$ and $V$ be two random variables, and suppose that 
\begin{equation}
\sup_{x\in \mathbb{R}}F^{\prime}_{V}(x)\leq A.  \label{essen.lem.cond}
\end{equation}

\noindent Then
\begin{equation}
\sup_{x}\left\vert F_{U}(x)-F_{V}(x)\right\vert \leq \frac{1}{\pi}\int_{-T}^{T}\left\vert \frac{\psi _{U}(t)-\psi _{V}(t)}{t}\right\vert dt+24A/(\pi T).  \label{essen.lem.form}
\end{equation}
\end{lemma}

\bigskip

\noindent{Proof}. Suppose that 
\begin{equation*}
\int_{-T}^{T}\left\vert \frac{\psi _{U}(t)-\psi _{V}(t)}{t}\right\vert dt<+\infty,
\end{equation*}

\noindent for $T>0$, otherwise (\ref{essen.lem.form}) is obvious. We may consider, by using Kolmogorov Theorem, that we are on a probability space holding the ordered pair $(U,V)$ and an absolutely continuous random variable $Z_{T}$ with characteristic function $(1-\left\vert t\right\vert /T)1_{(\left\vert t\right\vert \leq T)}$ as allowed by Lemma \ref{essen.lem.2} such that $Z_T$ is independent from $U$ and $V$. Given the \textit{cdf} $F_{(U,V)}$  of $(U,V)$ and the \textit{cdf} $F_{Z_T}$, the \textit{cdf} of $(U,V,Z_T)$ is given by

$$
F_{(U,V,Z_T)}(u,v,z)=F_{(U,V)}(u,v) \times F_{Z_T}(z), \ \ (u,v,z) \in \mathbb{R}^3.
$$

\noindent The probability space by using the Lebesgue-Stieljes measure of $F_{(U,V,Z_T)}$ following Point (c5), Section \label{proba_02_rv_sec_04} (page \pageref{kolmconst_02}).\\

\noindent Now, we recall the convolution formula on $U$ and $Z_T$ : 
\begin{equation*}
F_{U+Z_T}(x)=\int F_{U}(x-y)f_{Z_{T}}(y) \ dy, \ x \in \mathbb{R}.
\end{equation*}

\noindent Define $F_{U+Z_T}$ likewise. Set 
\begin{equation*}
\Delta (x)=F_{U}(x)-F_{V}(x), \ x \in \mathbb{R}.
\end{equation*}

\noindent and 
\begin{equation}
\Delta _{T}(x)=\int \Delta (x-y)f_{Z_{T}}(y)dy=F_{U+Z_{T}}(x)-F_{V+Z_{T}}(x). \label{essen.deltaT}
\end{equation}

\noindent We remark that for any fixed $t$, $\psi _{Z_{T}}(t)=(1-\left\vert t\right\vert /T)1_{(\left\vert t\right\vert \leq T)}\rightarrow 1$, which is the characteristic function of $0$. Then  $Z_{T}$ weakly converges to $0$, that is equivalent to $Z_{T}\rightarrow _{P}0$. Using results of weak theory (see for example Chapter 5, Subsection 3.2.3, Proposition 21 in \cite{ips-wcrv-ang}) implies that $U+Z_{T}\rightsquigarrow U$ and $V+Z_{T}\rightsquigarrow V$. By returning back to the distribution functions that are continuous, we have from (\ref{essen.deltaT})

\begin{equation*}
\text{For any x, }\Delta _{T}(x)\rightarrow \Delta (x)\text{ as }%
T\rightarrow \infty .
\end{equation*}

\noindent By applying Proposition \ref{proba02_rv.inv01} (see section Part V, Section \ref{proba_02_rv_sec_06}, Chapter \ref{proba_02_rv}), we have for continuity points $x$ and $b$ of both
$F_{U+Z_T}$ and $F_{V+Z_T}$, with $b<x$,

\begin{eqnarray*}
&&(F_{U+Z_T}(x)-F_{U+Z_T}(b)-((F_{U+Z_T}(x)-F_{U+Z_T}(b))\\
&=&\lim_{U \rightarrow +\infty} \frac{1}{2\pi }\int_{-U}^{U}\frac{e^{-ixt}-e^{-ibt}}{it}\left(\psi_U(t)-\psi_V(t)\right)\psi_{Z_T}(t)dt\\
&=&\frac{1}{2\pi }\int_{-T}^{T}\frac{e^{-ixt}-e^{-ibt}}{it}\left(\psi_U(t)-\psi_V(t)\right)\psi_{Z_T}(t)dt,\\
&=&\frac{1}{2\pi }\int_{-T}^{T}\left(\psi_U(t)-\psi_V(t)\right)\psi_{Z_T}(t)\biggr(-\int_{x}^{b} e^{-itv} dv\biggr)dt,
\end{eqnarray*}

\bigskip \noindent since $\psi _{Z_{T}}(t)$ vanishes outside $[-T,T]$. By letting $b\downarrow -\infty$ over the set of continuity points of both $F_{U+Z_T}$ and $F_{V+Z_T}$, and by using the Fatou-Lebesgue convergence theorem at right, we get

\begin{eqnarray*}
&&F_{U+Z_T}(x)-F_{U+Z_T}(x)\\
&=&\frac{1}{2\pi }\int_{-T}^{T}\left(\psi_U(t)-\psi_V(t)\right)\psi_{Z_T}(t)\biggr(-\int_{x}^{-\infty} e^{-itv} dv\biggr)dt\\
&=&\frac{1}{2\pi }\int_{-T}^{T} -e^{-itx} \left(\psi_U(t)-\psi_V(t)\right)\psi_{Z_T}(t)\biggr(-\int_{x}^{-\infty} e^{-itv} dv\biggr)dt\\
\end{eqnarray*}

\noindent which gives, for any continuity point $x$  of both $F_{U+Z_T}$ and $F_{V+Z_T}$. By taking the supremum of those continuity point $x$  of both $F_{U+Z_T}$ and $F_{V+Z_T}$ (which amounts to taking the supremum over $\mathbb{R}$ by right-continuity), we finally get

$$
\|\Delta_T\|_{+\infty} \leq \frac{1}{2\pi}\int_{-T}^{T}\left\vert \frac{\psi _{U}(t)-\psi_{V}(t)}{t}\right\vert dt.
$$

\bigskip \noindent Since we want to prove Formula (\ref{essen.lem.form}), the last formula says it will be enough to prove

\begin{equation}
\left\Vert \Delta \right\Vert _{\infty }\leq 2\left\Vert \Delta
_{T}\right\Vert _{\infty }+24A/(\pi T)  \label{essen.lem.formA}
\end{equation}

\bigskip \noindent We remark that $\Delta$ is bounded and is right-continuous with left-limits at each point of $\mathbb{R}$ and $\Delta (+\infty )=\Delta(-\infty )=0$. Then by Lemma \ref{essen.lem4} below, there exists a $x_{0}\in R$ such that $\left\Vert \Delta \right\Vert _{\infty }=\left\vert\Delta (x_{0})\right\vert $ or $\left\Vert \Delta \right\Vert _{\infty}=\left\vert f(x_{0}-)\right\vert$. We continue with the case where $\left\Vert \Delta \right\Vert _{\infty }=\left\vert \Delta (x_{0})\right\vert =\Delta (x_{0}).$ Handling the other cases is similar. We have for any $s>0$ : 

\begin{equation*}
\Delta (x_{0}+s)-\Delta (x_{0})=\left\{ F_{U}(x_{0}+s)-F_{U}(x_{0})\right\}
-\left\{ F_{V}(x_{0}+s)-F_{V}(x_{0})\right\}
\end{equation*}

\noindent and, by (\ref{essen.lem.cond}), 

\begin{equation*}
F_{V}(x_{0}+s)-F_{V}(x_{0})=\int_{x_{0}}^{x_{0}+s}F_{V}^{\prime }(t)dt\leq
As.
\end{equation*}

\noindent Next 
\begin{equation*}
\Delta (x_{0}+s)-\Delta (x_{0})=\left\{ F_{U}(x_{0}+s)-F_{U}(x_{0})\right\}
-As\geq -As
\end{equation*}

\noindent since $\left\{ F_{U}(x_{0}+s)-F_{U}(x_{0})\right\} \geq 0$ $(F_{U}$
increasing). This gives for any $s\geq 0$%
\begin{equation*}
\Delta (x_{0}+s)\geq \left\Vert \Delta \right\Vert _{\infty }-As.
\end{equation*}

\noindent By applying this to $s=\left\Vert \Delta \right\Vert _{\infty }/(2A))+y$ for $\left\vert y\right\vert \leq \left\Vert \Delta \right\Vert _{\infty}/(2A)$, we get 
\begin{equation}
\Delta \left(x_{0}+\frac{\left\Vert \Delta \right\Vert _{\infty }}{2A}+y\right)\geq \frac{\left\Vert \Delta \right\Vert _{\infty }}{2}-Ay  \label{essen.lem5a}
\end{equation}

\noindent We going to apply this to $\Delta _{T}$, while reminding the definition, to see that 
\begin{eqnarray*}
\Delta _{T}\left(x_{0}+\frac{\left\Vert \Delta \right\Vert _{\infty }}{2A}\right) &=&\int \Delta \left(x_{0}+\frac{\left\Vert \Delta \right\Vert _{\infty}}{2A}-y\right)f_{Z_{T}}(y)dy \\
&=&\int_{\left\{ \left\vert y\right\vert \leq \left\Vert \Delta \right\Vert_{\infty }/(2A)\right\} }\Delta \left(x_{0}+\frac{\left\Vert \Delta \right\Vert_{\infty }}{2A}-y\right)f_{Z_{T}}(y)dy \\
&+&\int_{\left\{ \left\vert y\right\vert >\left\Vert \Delta \right\Vert _{\infty }/(2A)\right\} }\Delta \left(x_{0}+\frac{\left\Vert \Delta \right\Vert _{\infty }}{2A}-y\right)f_{Z_{T}}(y)dy.
\end{eqnarray*}

\noindent On one hand, by (\ref{essen.lem5a}), we have
\begin{eqnarray*}
&&\int_{\left\{ \left\vert y\right\vert \leq \left\Vert \Delta \right\Vert/(2A)\right\} }\Delta (x_{0}+\frac{\left\Vert \Delta \right\Vert _{\infty }}{2A}-y)f_{Z_{T}}(y)dy\\
 &\geq &\int_{\left\{ \left\vert y\right\vert \leq \left\Vert \Delta \right\Vert _{\infty }/(2A)\right\} }\left( \frac{\left\Vert \Delta \right\Vert _{\infty }}{2}-Ay\right) f_{Z_{T}}(y)dy \\
&\geq &\frac{\left\Vert \Delta \right\Vert _{\infty }}{2}\int_{\left\{\left\vert y\right\vert \leq \left\Vert \Delta \right\Vert /(2A)\right\}}f_{Z_{T}}(y)dy \\
&=&\frac{\left\Vert \Delta \right\Vert _{\infty }}{2}P\left( \left\vert Z_{T}\right\vert \leq \frac{\left\Vert \Delta \right\Vert _{\infty }}{2A} \right) \\
&=&\frac{\left\Vert \Delta \right\Vert _{\infty }}{2}\left\{ 1-P\left(\left\vert Z_{T}\right\vert >\frac{\left\Vert \Delta \right\Vert _{\infty }}{2A}\right) \right\}
\end{eqnarray*}

\noindent and for the other term, we use the following trivial inequality 
$$
\Delta(\cdot )\geq -\sup_{x\in R}\left\vert \Delta (x)\right\vert =-\left\Vert\Delta \right\Vert _{\infty}
$$ 

\noindent to have 
\begin{eqnarray*}
&&\int_{\left\{ \left\vert y\right\vert >\left\Vert \Delta \right\Vert
/(2A)\right\} }\Delta \left(x_{0}+\frac{\left\Vert \Delta \right\Vert _{\infty }}{2A}-y\right)f_{Z_{T}}(y)dy\\
 &\geq &-\left\Vert \Delta \right\Vert _{\infty
}\int_{\left\{ \left\vert y\right\vert >\left\Vert \Delta \right\Vert
/(2A)\right\} }f_{Z_{T}}(y)dy \\
&\geq &-\Delta \left\Vert \Delta \right\Vert _{\infty }\int_{\left\{
\left\vert y\right\vert >\left\Vert \Delta \right\Vert /(2A)\right\}
}f_{Z_{T}}(y)dy \\
&=&-\Delta \left\Vert \Delta \right\Vert _{\infty }P\left( \left\vert
Z_{T}\right\vert >\frac{\left\Vert \Delta \right\Vert _{\infty }}{2A}\right).
\end{eqnarray*}

\noindent The two last formulas lead to
\begin{equation*}
\Delta _{T}\left( x_{0}+\frac{\left\Vert \Delta \right\Vert _{\infty }}{2A}%
\right) \geq \frac{\left\Vert \Delta \right\Vert _{\infty }}{2}\left\{
1-3P\left( \left\vert Z_{T}\right\vert >\frac{\left\Vert \Delta \right\Vert
_{\infty }}{2A}\right) \right\}
\end{equation*}

\noindent and next 
\begin{equation}
\left\Vert \Delta _{T}\right\Vert _{\infty }\geq \frac{\left\Vert \Delta
\right\Vert _{\infty }}{2}\left\{ 1-3P\left( \left\vert Z_{T}\right\vert >%
\frac{\left\Vert \Delta \right\Vert _{\infty }}{2A}\right) \right\} .
\label{essen.lem5c}
\end{equation}

\noindent In this last step, we have
\begin{eqnarray*}
P\left( \left\vert Z_{T}\right\vert >\frac{\left\Vert \Delta \right\Vert
_{\infty }}{2A}\right) &=&\int_{\left\{ \left\vert y\right\vert >\left\Vert
\Delta \right\Vert _{\infty }/(2A)\right\} }\frac{1-\cos yT}{\pi Ty^{2}}dy \\
&=&\int_{\Delta \left\Vert \Delta \right\Vert _{\infty }/(2A)}^{+\infty }%
\frac{1-\cos yT}{\pi Ty^{2}}dy \\
&=&\frac{1}{\pi }\int_{\Delta \left\Vert \Delta \right\Vert _{\infty
}T/(4A)}^{+\infty }\frac{1-\cos 2y}{y^{2}}dy \\
&=&\frac{2}{\pi }\int_{\Delta \left\Vert \Delta \right\Vert _{\infty
}T/(4A)}^{+\infty }\frac{\sin ^{2}y}{y^{2}}dy \\
&\leq &\frac{2}{\pi }\int_{\left\Vert \Delta \right\Vert _{\infty
}T/(4A)}^{+\infty }\frac{1}{y^{2}}dy \\
&=&\frac{8A}{\pi T\left\Vert \Delta \right\Vert _{\infty }}.
\end{eqnarray*}

\noindent This and (\ref{essen.lem5c}) yield%
\begin{eqnarray*}
2\left\Vert \Delta _{T}\right\Vert _{\infty } &\geq &\left\{ \left\Vert
\Delta \right\Vert _{\infty }-3\left\Vert \Delta \right\Vert _{\infty
}P\left( \left\vert Z_{T}\right\vert >\frac{\left\Vert \Delta \right\Vert
_{\infty }}{2A}\right) \right\} \\
&\geq &\left\Vert \Delta \right\Vert _{\infty }-\frac{24A}{T\pi },
\end{eqnarray*}

\noindent which implies 
\begin{equation*}
\left\Vert \Delta \right\Vert _{\infty }\leq 2\left\Vert \Delta
_{T}\right\Vert _{\infty }+\frac{24A}{T\pi }.
\end{equation*}

\bigskip
\noindent This was the target, that is Formula (\ref{essen.lem.formA}), which is enough to have the final result (\ref{essen.lem.form}). $\blacksquare$\newline

\noindent \textbf{Technical Lemmas used by the proof}.\\

\begin{lemma}
\label{essen.dev.charact}. Let $X$ be a real random variable with $n+1$
finite moments. Then for any $t\in \mathbb{R},$%
\begin{equation*}
\left\vert \psi _{X}(t)-\sum_{k=0}^{n}\frac{(it)^{k}}{k!}\mathbb{E}%
\left\vert X\right\vert ^{k}\right\vert \leq \min \left( \frac{2\left\vert
t\right\vert ^{n}}{n!}\mathbb{E}\left\vert X\right\vert ^{n},\frac{%
\left\vert t\right\vert ^{n+1}}{(n+1)!}\mathbb{E}\left\vert X\right\vert
^{n+1}\right).
\end{equation*}
\end{lemma}

\bigskip \noindent \textbf{Proof}. We may use the Taylor-Mac-Laurin expansion formula, 
\begin{equation*}
f(y)=\sum_{k=0}^{n}\frac{y^{k}}{k!}f^{(k)}(0)+\int_{0}^{y}\frac{(y-x)^{n}}{n!%
}f^{(n)}(x)dx
\end{equation*}

\bigskip \noindent for $f(y)=e^{iy}.$ We have $f^{(k)}(y)=i^{k}f(y)$ and then 
\begin{equation}
e^{iy}=\sum_{k=0}^{n}\frac{(iy)^{k}}{k!}+i^{n+1}\int_{0}^{y}e^{ix}\frac{%
(y-x)^{n}}{n!}dx.  \label{devM3}
\end{equation}

\bigskip \noindent Then
\begin{equation*}
\left\vert e^{iy}-\sum_{k=0}^{n}\frac{(iy)^{k}}{k!}\right\vert \leq \frac{1}{n!}\int_{0}^{\left\vert y\right\vert }\left\vert y-x\right\vert ^{n}dx\leq 
\frac{\left\vert y\right\vert ^{n+1}}{(n+1)!}.
\end{equation*}

\bigskip  \noindent We apply (\ref{devM3}) for $n-1,$ that is
\begin{equation*}
e^{iy}=\sum_{k=0}^{n-1}\frac{(iy)^{k}}{k!}+i^{n}\int_{0}^{y}e^{ix}\frac{%
(y-x)^{n-1}}{(n-1)!}dx
\end{equation*}

\bigskip \noindent and we use the decomposition $e^{ix}=1+(e^{ix}-1)$ to get
\begin{equation*}
e^{iy}=\sum_{k=0}^{n-1}\frac{(iy)^{k}}{k!}+i^{n}\int_{0}^{y}\frac{(y-x)^{n-1}%
}{(n-1)!}dx+i^{n}\int_{0}^{y}(e^{ix}-1)\frac{(y-x)^{n-1}}{(n-1)!}dx.
\end{equation*}

\bigskip \noindent We have 
\begin{equation*}
i^{n}\int_{0}^{y}\frac{(y-x)^{n-1}}{(n-1)!}dx=(iy)^{n},
\end{equation*}

\bigskip \noindent which leads to
\begin{equation*}
e^{iy}=\sum_{k=0}^{n}\frac{(iy)^{k}}{k!}+i^{n}\int_{0}^{y}(e^{ix}-1)\frac{%
(y-x)^{n-1}}{(n-1)!}dx
\end{equation*}

\bigskip \noindent and next, since $\left\vert (e^{ix}-1)\right\vert \leq 2$,

\begin{equation*}
\left\vert e^{iy}-\sum_{k=0}^{n}\frac{(iy)^{k}}{k!}\right\vert \leq \frac{2}{%
(n-1)!}\int_{0}^{\left\vert y\right\vert }\left\vert y-x\right\vert
^{n}dx\leq \frac{2\left\vert y\right\vert ^{n}}{n!}.
\end{equation*}

\bigskip \noindent We then get
\begin{equation*}
\left\vert e^{iy}-\sum_{k=0}^{n}\frac{(iy)^{k}}{k!}\right\vert \leq \min
\left( \frac{2\left\vert y\right\vert ^{n}}{n!},\frac{\left\vert
y\right\vert ^{n+1}}{(n+1)!}\right) .
\end{equation*}

\bigskip \noindent We apply this to a random real variable $X$ with enough finite moments to get 
\begin{equation*}
\left\vert \mathbb{E}e^{itX}-\mathbb{E}\sum_{k=0}^{n}\frac{(itX)^{k}}{k!}%
\right\vert \leq \mathbb{E}\left\vert e^{itX}-\sum_{k=0}^{n}\frac{(itX)^{k}}{%
k!}\right\vert \leq \mathbb{E}\min \left( \frac{2\left\vert tX\right\vert
^{n}}{n!},\frac{\left\vert tX\right\vert ^{n+1}}{(n+1)!}\right) ,
\end{equation*}

\bigskip \noindent and then
\begin{equation*}
\left\vert \psi _{X}(t)-\sum_{k=0}^{n}\frac{(it)^{k}}{k!}\mathbb{E}%
\left\vert X\right\vert ^{k}\right\vert \leq \min \left( \frac{2\left\vert
t\right\vert ^{n}}{n!}\mathbb{E}\left\vert X\right\vert ^{n},\frac{%
\left\vert t\right\vert ^{n+1}}{(n+1)!}\mathbb{E}\left\vert X\right\vert
^{n+1}\right)
\end{equation*}

\bigskip \begin{lemma} \label{essen.approx} With the notations and assumptions of the Theorem, we
have :\newline

\noindent \textbf{Part 1}. 
\begin{equation*}
\left\vert \exp (itS_{n}/s_{n})-\exp (-t^{2}/2)\right\vert \leq 0.4466464\frac{\beta _{n}^{3}}{s_{n}^{3}}\left\vert t\right\vert ^{3}\exp (-\left\vert t\right\vert ^{2}/2)\text{ for }\left\vert t\right\vert \leq s_{n}/(2\beta _{n}).
\end{equation*}

\bigskip \noindent \textbf{Part 2}. 
\begin{equation*}
\left\vert \exp (itS_{n}/s_{n})-\exp (-t^{2}/2)\right\vert \leq 16\frac{%
\beta _{n}^{3}}{s_{n}^{3}}\left\vert t\right\vert ^{3}\exp (-\left\vert
t\right\vert ^{2}/3)\text{ for }\left\vert t\right\vert \leq
s_{n}^{3}/(4\beta _{n}^{3}).
\end{equation*}
\end{lemma}

\bigskip \noindent \textbf{Proof}.\newline

\noindent \textbf{Proof of Part 1}. Let us prove that%
\begin{equation*}
\left\vert \exp (itS_{n}/s_{n})-\exp (-t^{2}/2)\right\vert \leq 0.5\frac{%
\beta _{n}^{3}}{s_{n}^{3}}\left\vert t\right\vert ^{3}\exp (-\left\vert
t\right\vert ^{3}/3)\text{ \ for } \left\vert t\right\vert \leq
s_{n}/(2\beta _{n})
\end{equation*}

\bigskip \noindent To this end, we use the following expansion 
\begin{equation*}
\exp (itX_{k}/s_{n})=1-\frac{t^{2}\sigma _{k}^{2}}{2s_{n}^{2}}+\theta \frac{%
\left\vert t\right\vert ^{3}\gamma _{k}^{3}}{6s_{n}^{3}},\text{ }\left\vert
\theta \right\vert \leq 1.
\end{equation*}

\noindent For $\left\vert t\right\vert \leq s_{n}/(2\beta _{n}),$

\bigskip 
\begin{equation*}
r_{1,k}=\left\vert \theta \frac{\left\vert t\right\vert ^{3}\gamma _{k}^{3}}{%
6s_{n}^{3}}\right\vert \leq \frac{\left\vert t\right\vert ^{3}\gamma _{k}^{3}%
}{6s_{n}^{3}}.
\end{equation*}

\noindent Next
\begin{eqnarray*}
\log \exp (itX_{k}/s_{n}) &=&\log \left( 1-\frac{t^{2}\sigma _{k}^{2}}{%
2s_{n}^{2}}+r_{1,k}\right) \\
&=&-\frac{t^{2}\sigma _{k}^{2}}{2s_{n}^{2}}+r_{1,k}+r_{2,k},
\end{eqnarray*}

\noindent with, after having used the $c_{r}$-inequality, 
\begin{eqnarray*}
r_{2,k} &\leq &\left\vert -\frac{t^{2}\sigma _{k}^{2}}{2s_{n}^{2}}%
+r_{1}\right\vert ^{2}\leq 2\left\vert -\frac{t^{2}\sigma _{k}^{2}}{%
2s_{n}^{2}}\right\vert ^{2}+2\left\vert r_{1}\right\vert ^{2} \\
&\leq &\frac{1}{2}\left( \frac{\left\vert t\right\vert ^{3}\sigma _{k}^{3}}{%
s_{n}^{3}}\right) \left( \frac{\left\vert t\right\vert \sigma _{k}}{s_{n}}%
\right) +\frac{1}{18}\left( \frac{\left\vert t\right\vert ^{3}\gamma _{k}^{3}%
}{s_{n}^{3}}\right) \left\{ \frac{\left\vert t\right\vert ^{3}\gamma _{k}^{3}%
}{s_{n}^{3}}\right\} .
\end{eqnarray*}

\bigskip \noindent Now $\left\vert t\right\vert \leq s_{n}/(2\beta _{n})$ implies%
\begin{equation*}
\left( \frac{\left\vert t\right\vert \sigma _{k}}{s_{n}}\right) \leq \left( 
\frac{\sigma _{k}}{s_{n}}\times \frac{s_{n}}{2\beta _{n}}\right) =\frac{1}{2}%
\frac{\sigma _{k}}{\beta _{n}}=\frac{1}{2}\left( \frac{\sigma _{k}^{3}}{%
\beta _{n}^{3}}\right) ^{1/3}\leq \frac{1}{2}\left( \frac{\sigma _{k}^{3}}{%
\gamma _{n}^{3}}\right) ^{1/3}=\frac{1}{2}\frac{\sigma _{k}}{\gamma _{k}}%
\leq \frac{1}{2},
\end{equation*}

\noindent by Lyapounov's inequality, that is for $1\leq p \leq q$, 
\begin{equation*}
\left( E\left\vert X_{k}\right\vert ^{p}\right) ^{1/p}\leq \left(
E\left\vert X_{k}\right\vert ^{q}\right) ^{1/q}
\end{equation*}

\noindent and next%
\begin{equation*}
\left\{ \frac{\left\vert t\right\vert ^{3}\gamma _{k}^{3}}{s_{n}^{3}}%
\right\} \leq \left\{ \frac{\gamma _{k}^{3}}{s_{n}^{3}}\times \frac{s_{n}^{3}%
}{8\beta _{n}^{3}}\right\} \leq \frac{1}{8}.
\end{equation*}

\noindent We arrive, after applying again Lyapounov's inequality, at 
\begin{equation*}
r_{2,k}\leq \frac{1}{2}\left( \frac{\left\vert t\right\vert ^{3}\sigma
_{k}^{3}}{s_{n}^{3}}\right) \times \frac{1}{2}+\frac{1}{18}\left( \frac{%
\left\vert t\right\vert ^{3}\gamma _{k}^{3}}{s_{n}^{3}}\right) \frac{1}{8}=%
\frac{37}{144}\frac{\left\vert t\right\vert ^{3}\gamma _{k}^{3}}{s_{n}^{3}}.
\end{equation*}

\noindent Next 
\begin{eqnarray*}
\log \exp (itS_{n}/s_{n}) &=&\sum_{k=1}^{n}\log \exp (itX_{k}/s_{n}) \\
&=&-t^{2}/2+r_{n} \\
&=&-t^{2}/2+\sum_{k=1}^{n}\frac{\left\vert t\right\vert ^{3}\gamma _{k}^{3}}{%
6s_{n}^{3}}+\frac{37}{144}\frac{\left\vert t\right\vert ^{3}\gamma _{k}^{3}}{%
s_{n}^{3}} \\
&=&-t^{2}/2+\frac{61}{144}\frac{\left\vert t\right\vert ^{3}\beta _{n}^{3}}{%
s_{n}^{3}} \\
&\leq &-t^{2}/2+\frac{61}{144}\times \frac{1}{8},
\end{eqnarray*}

\noindent where we used $\left\vert t\right\vert \leq s_{n}/(2\beta _{n})$ at the last step. We already have
\begin{equation*}
\exp (itS_{n}/s_{n})=\exp (-t^{2}/2+r_{n})
\end{equation*}

\noindent so that 
\begin{equation*}
\left\vert \exp (itS_{n}/s_{n})-\exp (-t^{2}/2)\right\vert =\exp (-t^{2}/2)\left\Vert e^{r_{n}}-1\right\Vert .
\end{equation*}

\bigskip \noindent We use the formula $\left\Vert e^{z}-1\right\Vert \leq \left\Vert z\right\Vert e^{\left\Vert z\right\Vert}$ to see that

\begin{eqnarray*}
\left\vert \exp (itS_{n}/s_{n})-\exp (-t^{2}/2)\right\vert &=&\exp
(-t^{2}/2)\left\Vert e^{r_{n}}-1\right\Vert \\
&\leq &\exp (-t^{2}/2)\left\Vert r_{n}\right\Vert e^{\left\Vert
r_{n}\right\Vert } \\
&\leq &\exp (-t^{2}/2)\frac{61}{144}\frac{\left\vert t\right\vert ^{3}\beta
_{n}^{3}}{s_{n}^{3}}e^{61/(8\ast 144)} \\
&\leq &0.4466464\times \exp (-t^{2}/2)\frac{\left\vert t\right\vert
^{3}\beta _{n}^{3}}{s_{n}^{3}} \\
&\leq &16\times \exp (-t^{2}/2)\frac{\left\vert t\right\vert ^{3}\beta
_{n}^{3}}{s_{n}^{3}} \\
&\leq &16\times \exp (-t^{2}/2)\frac{\left\vert t\right\vert ^{3}\beta
_{n}^{3}}{s_{n}^{3}},
\end{eqnarray*}

\noindent since $\exp (-t^{2}/2)\leq \exp (-t^{2}/3)$.\newline

\noindent \textbf{Proof of Part 2}. This is proved as follows. If $s_{n}^{3}/(4\beta _{n}^{3})\leq s_{n}/(2\beta _{n})$, Part 2 is implied by Part 1. Then, we only need to prove Part 2 in the case%
\begin{equation*}
s_{n}/(2\beta _{n})<s_{n}^{3}/(4\beta _{n}^{3}),
\end{equation*}

\noindent and only for $t$ satisfying
\begin{equation*}
s_{n}/(2\beta _{n})<t\leq s_{n}^{3}/(4\beta _{n}^{3}).
\end{equation*}

\noindent Let us proceed by considering the symmetrized form of $X_{k}$, denoted by $X_{k}^{s}$, and defined by \ \ 
\begin{equation*}
X_{k}^{s}=X_{k}-X_{k}^{\prime },
\end{equation*}

\bigskip \noindent where $X_{k}^{\prime }$ is a random variable with the same law than $X_{k}$ and independent of $X_{k}$. Then, obviously, $\mathbb{E}X_{k}^{s}=\mathbb{E}X_{k}-\mathbb{E}X_{k}^{\prime }=\mathbb{E}X_{k}-\mathbb{E}X_{k}=0$ and 
\begin{equation*}
\mathbb{V}ar(X_{k}^{s})=\mathbb{V}ar(X_{k})+\mathbb{V}ar(X_{k}^{\prime })+2\text{ }\mathbb{C}ov(X_{k},X_{k}^{\prime })=\sigma _{k}^{2}+\sigma _{k}^{2}+0=2\sigma _{k}^{2}
\end{equation*}

\bigskip \noindent and finally, by the $C_{r}$-inequality 
\begin{equation*}
\mathbb{E}\left\vert X_{k}^{s}\right\vert ^{r}=\mathbb{E}\left\vert
X_{k}-X_{k}\right\vert ^{r}\leq c_{r}(\mathbb{E}\left\vert X_{k}\right\vert
^{r}+\mathbb{E}\left\vert X_{k}^{\prime }\right\vert ^{r}),
\end{equation*}

\bigskip \noindent with $c_{r}=2^{r-1},$ $r\geq 1.$ Apply it to $r=3$ to get

\begin{equation*}
\mathbb{E}\left\vert X_{k}^{s}\right\vert ^{3}\leq 4(\mathbb{E}\left\vert
X_{k}\right\vert ^{3}+\mathbb{E}\left\vert X_{k}^{\prime }\right\vert
^{3})=8\gamma _{k}^{3}.
\end{equation*}

\bigskip \noindent Now, we remark that we have for any real random variable
\begin{equation*}
\psi _{X}(t)=\int \cos (tx)\text{ }d\mathbb{P}_{X}(x)+i\int \sin (tx)\text{ }d\mathbb{P}_{X}(x).
\end{equation*}

\noindent and 
\begin{equation*}
\psi _{-X}(t)=x\psi _{-X}(t)=\int \cos (tx)\text{ }d\mathbb{P}_{X}(x)-i\int
\sin (tx)\text{ }d\mathbb{P}_{X}(x)=\overline{\psi _{X}(t)},
\end{equation*}

\noindent where $\overline{\psi _{X}(t)}$ is the conjugate of $\psi _{X}(t)$. Next, from this and by independence, we have 
\begin{eqnarray*}
\psi _{X_{k}^{s}}(t) &=&\psi _{X_{k}-X_{k}^{\prime }}(t)=\psi
_{X_{k}}(t)\psi _{-X_{k}^{\prime }}(t)=\psi _{X_{k}}(t)\psi _{-X_{k}}(t) \\
&=&\psi _{X_{k}}(t)\overline{\psi _{X_{k}}(t)}=\left\Vert \psi
_{X_{k}}(t)\right\Vert ^{2},
\end{eqnarray*}

\noindent where, here, $\left\Vert \circ \right\Vert $ denotes the norm in the complex space. Next, we apply Lemma \ref{essen.dev.charact} to $X_{k}^{s} $ at the order $n=2$ to get%
\begin{equation*}
\left\vert \psi _{X_{k}^{s}}(t)-(1-2\sigma _{k}^{2}t^{2})\right\vert \leq 
\frac{8\left\vert t\right\vert ^{3}}{6}\gamma _{k}^{3}=\frac{4\left\vert
t\right\vert ^{3}}{3}\gamma _{k}^{3}.
\end{equation*}

\noindent The triangle inequality leads to
\begin{equation*}
\psi _{X_{k}^{s}}(t)\leq 1-2\sigma _{k}^{2}t^{2}+\frac{4\left\vert
t\right\vert ^{3}}{3}\gamma _{k}^{3},
\end{equation*}

\noindent which gives ($\bowtie$)
\begin{equation*}
\psi _{X_{k}^{s}/s_{n}}(t)\leq \exp \{-2\frac{\sigma _{k}^{2}t^{2}}{s_{n}^{2}%
}+\frac{4\left\vert t\right\vert ^{3}}{3s_{n}^{3}}\gamma _{k}^{3}\}.
\end{equation*}

\bigskip \noindent Denote also $S_{n}^{s}=X_{1}^{s}+...+X_{n}^{s}.$ Then, by reminding that $\psi _{X_{k}^{s}}(t)$ is real and non-negative and that $\left\Vert \psi _{X_{k}^{s}/s_{n}}(t)\right\Vert $ is an absolute value, we have 
\begin{eqnarray*}
\psi _{S_{n}^{s}/s_{n}}(t) &=&\prod\limits_{k=1}^{n}\psi _{X_{k}^{s}}(t)\leq
\prod\limits_{k=1}^{n}\exp \{-\frac{\sigma _{k}^{2}t^{2}}{s_{n}^{2}}+\frac{%
4\left\vert t\right\vert ^{3}}{3s_{n}^{3}}\gamma _{k}^{3}\} \\
&=&\exp \{-\sum\limits_{k=1}^{n}\frac{\sigma _{k}^{2}t^{2}}{s_{n}^{2}}%
+\sum\limits_{k=1}^{n}\frac{4\left\vert t\right\vert ^{3}}{3s_{n}^{3}}\gamma
_{k}^{3}\} \\
&=&\exp \{-t^{2}+\frac{4\left\vert t\right\vert ^{3}\beta _{n}^{3}}{%
3s_{n}^{3}}\}.
\end{eqnarray*}

\noindent Now, for $\left\vert t\right\vert \leq s_{n}^{3}/(4\beta_{n}^{3}), $ it comes that 
\begin{eqnarray*}
\psi _{S_{n}^{s}/s_{n}}(t) &\leq &\exp \{-t^{2}+\frac{4\left\vert
t\right\vert ^{3}\beta _{n}^{3}}{3s_{n}^{3}}\} \\
&\leq &\exp \{-t^{2}+\frac{t^{2}}{3}\}=\exp \{-2t^{2}/3\}.
\end{eqnarray*}

\noindent Since $\psi _{S_{n}^{s}/s_{n}}(t)=\left\Vert \psi
_{S_{n}/s_{n}}(t)\right\Vert ^{2},$ we have%
\begin{equation*}
\left\Vert \psi _{S_{n}/s_{n}}(t)\right\Vert \leq \exp \{-t^{2}/3\}.
\end{equation*}

\noindent Now, since, $s_{n}/(2\beta _{n})<\left\vert t\right\vert$,
\begin{eqnarray*}
1 &\leq &\frac{2\beta _{n}\left\vert t\right\vert }{s_{n}}=\frac{2\beta
_{n}\left\vert t\right\vert ^{3}}{s_{n}}\times \left( \frac{1}{\left\vert
t\right\vert ^{2}}\right) \\
&\leq &\frac{2\beta _{n}\left\vert t\right\vert ^{3}}{s_{n}}\times \left( 
\frac{4\beta _{n}^{2}}{s_{n}^{2}}\right) =\frac{8\beta _{n}^{3}\left\vert
t\right\vert ^{3}}{s_{n}^{3}}.
\end{eqnarray*}

\noindent To conclude, we say that

\begin{eqnarray*}
\left\Vert \psi _{S_{n}/s_{n}}(t)-\exp (-t^{2}/2)\right\Vert &\leq &\exp
(-t^{2}/2)+\left\Vert \psi _{S_{n}/s_{n}}(t)\right\Vert \\
&\leq &\exp (-t^{2}/2)+\exp \{-t^{2}/3\} \\
&\leq &2\exp (-t^{2}/2).
\end{eqnarray*}

\bigskip \noindent We conclude by using the following stuff :
\begin{eqnarray*}
\left\Vert \psi _{S_{n}/s_{n}}(t)-\exp (-t^{2}/2)\right\Vert &\leq &2\exp
(-t^{2}/2)\times (1) \\
&\leq &2\exp (-t^{2}/2)\frac{8\beta _{n}^{3}\left\vert t\right\vert ^{3}}{%
s_{n}^{3}} \\
&=&16\exp (-t^{2}/2)\frac{\beta _{n}^{3}\left\vert t\right\vert ^{3}}{%
s_{n}^{3}}.
\end{eqnarray*}

\bigskip

\noindent The following lemma on elementary real analysis has been used in the proof the Essen Lemma \ref{essen.lem}.

\begin{lemma}
\label{essen.lem4} Let $f$ be a bounded and non-constant right-continuous mapping from $\mathbb{R}$ to $\mathbb{R}$ with left-limits at each point of $\mathbb{R}$, such that 
\begin{equation*}
\lim_{x\rightarrow -\infty }f(x)=0\text{ and }\lim_{x\rightarrow
+\infty}f(x)=0.
\end{equation*}

\noindent Then there exists some $x_0 \in \mathbb{R}$ such that

\begin{equation*}
0<c=\sup_{x\in \mathbb{R}}|f(x)|=|f(x_{0})|\text{ or }c=|f(x_{0}-)|
\end{equation*}

\noindent where $f(x-)$ stands for the left-limit of $f$ at $x$.
\end{lemma}

\bigskip \noindent \textbf{Proof}. Let $c=\sup_{x{\in \mathbb{R}}}|f(x)|$. The number $c$ is strictly positive, otherwise $f$ would be constant and equal to zero, which would be contrary to the assumption. Now since $\lim_{{x\rightarrow -\infty }}f(x)=0$ and $\lim_{{x\rightarrow +\infty }}f(x)=0$, we can find $A>0$ such that

\begin{equation*}
\forall x,(\left\vert x\right\vert >A)\Longrightarrow (\left\vert
f(x)\right\vert <c/2).
\end{equation*}

\noindent So, we get
\begin{equation*}
c=\sup_{x\in {[-A,A]}}|f(x)|.
\end{equation*}

\noindent We remark that $c$ is finite since $f$ is bounded. Now consider a sequence $(x_{n})_{n\geq 0}$ $\subset \lbrack -A,A]$ such that $|f(x_{n})|\rightarrow c.$ Since $(x_{n})_{n\geq 0}$ $\subset \lbrack -A,A]$, by the Bolzano-Weierstrass property, there exists a subsequence $(x_{n(k)})_{k\geq 0}\subset (x_{n})_{n\geq 0}$ converging to some $x_{0}\in \lbrack -A,A]$. Consider

\begin{equation*}
I(\ell )=\{k\geq 1,x_{n(k)}\geq x_{0}\}\text{ \ and }I(r)=\{k\geq
1,x_{n(k)}<x_{0}\}.
\end{equation*}

\bigskip \noindent One of these two set is infinite. If $I(\ell )$ is infinite, we can find a sub-sequence $\left( x_{n(k_{j})}\right) _{j\geq 1}$ such that $x_{n(k_{j})}\geq x_{0}$ for any $j\geq 1$ and $x_{n(k_{j})}\rightarrow x_{0}$ as $j\rightarrow \infty .$ Then by right-continuity, \ $|f(x_{n(k_{j})})|\rightarrow \left\vert f(x_{0})\right\vert $ and as a sub-sequence of $|f(x_{n})|$ which converges to $c,$ we also have $|f(x_{n(k_{j})})|\rightarrow c$ as $j\rightarrow \infty$. Then 

\begin{equation*}
c=\left\vert f(x_{0})\right\vert
\end{equation*}

\bigskip  \noindent If $I(r)$ is infinite, we can find a sub-sequence $\left(x_{n(k_{j})}\right) _{j\geq 1}$ such that $x_{n(k_{j})}<x_{0}$ for any $j\geq 1$ and $x_{n(k_{j})}\rightarrow x_{0}$ as $j\rightarrow \infty$. Then by the existence of the left-limit of $f$ at $x_{0}$, \ $|f(x_{n(k_{j})})|\rightarrow \left\vert f(x_{0}-)\right\vert $ and as a sub-sequence of $|f(x_{n})|$ which converges to $c,$ we also have $|f(x_{n(k_{j})})|\rightarrow c$ as $j\rightarrow \infty$. Then
 
\begin{equation*}
c=\left\vert f(x_{0}-)\right\vert .
\end{equation*}

\bigskip

\begin{lemma}
We have the following inequality, for any complex number $z$
\begin{equation*}
\left\Vert e^{z}-1\right\Vert \leq \left\Vert z\right\Vert e^{\left\Vert
z\right\Vert }
\end{equation*}
\end{lemma}

%% file: proba_02_07_LIL_ang.tex
\section{Law of the Iterated Logarithm} \label{probab_02_lil_sec_06}

\noindent The Law of the Iterated Logarithm, abbreviated \textit{LIL} is one of the classical results in Probability Theory. As usual, it  was discovered for a sequence of \textit{iid} real-valued random variables. From a quick tour of the introduction on the question in \cite{loeve}, in \cite{gutt} and in \cite{feller2}, we may say that the \textit{LIL} goes back to Kintchine, and to Levy in the binary case and finally to Kolmogorov and to Cantelli in the general case for independent random variables. Other important contributors in the stationary case are Hartman \& Wintner, and Strassen. Here, we present the Kolmogorov Theorem as cited by \cite{loeve}.\\

\noindent Throughout this section, the iterated logarithm function $\log(\log(x))$, $x>e$, is denoted by $\log_2(x)$.\\

\noindent Let us ive the statement of \textit{LIL} law, by using the notation introduced above. A sequence of square integrable and centered
real-valued random variables $(X_n)_{n\geq 0}$ defined on the same probability space $(\Omega, \mathcal{A}, \mathbb{P})$ satisfies the \textit{LIL} if we have 

$$
\limsup_{n\rightarrow +\infty} \frac{S_n}{\sqrt{2s_n log_2 s_n^2}}=1, \ a.s.
$$

\bigskip \noindent If the  $(-X_n)_{n\geq 1}$ also satisfies the \textit{LIL}, we also have

$$
\liminf_{n\rightarrow +\infty} \frac{S_n}{\sqrt{2s_n log_2 s_n^2}}=-1, \ a.s.
$$

\bigskip \noindent The two conditions which required in the independent scheme to have the \textit{LIL} are :

$$
s_n\rightarrow +\infty \  as \ n\rightarrow +\infty, \ \ (C1).
$$ 

\bigskip \noindent and

$$
|X_n/s_n|=o((log_2 s_n^2)^{-1}), \ as \ n\rightarrow +\infty. \ \ (C2)
$$ 

\bigskip \noindent The conditions (C2) is used to ensure the following one :\\

\noindent (C3) : For all $\mathbb{R}_{+}\setminus \{0\} \ni c >1$, there exists a sub-sequence $(s_{n_k})_{k\geq 1}$ of $(s_{n})_{n\geq 1}$ such that 

$$
s_{n_{k+1}}/s_{n_k} \sim c \ as \ k\rightarrow +\infty,
$$

\bigskip \noindent which is ensured if $s_{n_k} \sim \beta c^{k}$, where $\beta>0$ is a real constant.\\

\noindent  \textbf{But it is important} that the proof below is based only (C1) and (C3). In the \textit{iid}, we have $s_n=\sigma \sqrt{n}$, $n\geq 1$. For any $c>0$, we may take 
$n_k=\sigma [c^{k}]$, $k\geq 1$ to have (C2).\\
  
\noindent Let us state the Kolmogorov Theorem.\\

\begin{theorem} Let $(X_n)_{n\geq 0}$ be a sequence of square integrable and centered real-valued random variables defined on the same probability space $(\Omega, \mathcal{A}, \mathbb{P})$ such that Condition (C1) and (C3) hold. Then the sequence satisfies the \textit{LIL}, that is 

$$
\limsup_{n\rightarrow +\infty} \frac{S_n}{\sqrt{2s_n log_2 s_n^2}}=1, \ a.s.
$$

\bigskip \noindent and, by replacing $X_n$ by $-X_n$ (which replacement does not change the variances), we have

$$
\liminf_{n\rightarrow +\infty} \frac{S_n}{\sqrt{2s_n log s_n}}=-1, \ a.s.
$$
\end{theorem}

\bigskip \noindent \textbf{Proof}. Let $\delta>0$. By applying (C3), let $(s_{n_k})_{k\geq 1}$ be a sub-sequence of $(s_{n})_{n\geq 1}$ such that $s_{n_k} \sim c^k$, as $k\rightarrow +\infty$, with $1<c <1+\delta$ so that 
$$
2\delta^{\prime}=1 -\frac{1+\delta}{c}>0 \ \ (S1)
$$

\bigskip \noindent and 
$$
(s_{n_{k}} (2 \log_2 s_{n_{k}})^{1/2})/ (s_{n_{k-1}} (2\log_2 s_{n_{k-1}})^{1/2}) \rightarrow c. \ \ (S2)
$$

\bigskip \noindent Now, since the the following class of integers intervals 
$$
\{[1, n_1[, [n_{k-1}, n_k[, k\geq 1\}
$$

\bigskip \noindent is a partition of $\mathbb{N}\setminus \{0\}$, we have for all each $n\geq 1$, there exists a unique $k\geq 1$ such that $n \geq [n_{k-1}, n_k[$ and so 
$$
S_n > (1+\delta) s_n (2\log_2 s_{n})^{1/2},
$$

\bigskip \noindent implies that

$$
S_ {n_{k}}^{\ast}=\sup_{\ell \leq n_k} S_{\ell} \geq S_n > (1+\delta) s_n (2\log_2 s_n^2)^{1/2} \geq (1+\delta) s_{n_{k-1}} 
(2\log_2 s_{n_{k-1}})^{1/2}, 
$$

\bigskip \noindent which by (S1) and (S2), implies, for large values of $k$, that

$$
S_ {n_{k}}^{\ast} > (1+\delta^{\prime}_1)^2 s_{n_{k}} (2\log_2 s_{n_{k}})^{1/2}.
$$

\bigskip \noindent When put together, these formulas above prove that

$$
\biggr(S_n > (1+\delta) s_n (2\log_2 s_n)^{1/2}, \ i.o.\biggr) \subset \biggr(S_ {n_{k}}^{\ast}  > (1+\delta^{\prime}) s_{n_k} (2\log_2 s_{n_{k}}^2)^{1/2}, \ i.o.\biggr)
$$

\bigskip \noindent Hence by Inequality (12) in Chapter \ref{proba_02_ineg} (See page \pageref{ineqBillingsley}), we have

\begin{eqnarray*}
\mathbb{P}\biggr(S_{n_{k}}^{\ast} &>& (1+\delta^{\prime}_2) s_{n_k} (2\log_2 s_{n_{k}})^{1/2}\biggr)\\
&\leq 2& \mathbb{P}\biggr(S_{n_{k}}^{\ast}> \left(1+\delta^{\prime}-\frac{\sqrt{2}}{(2 \log_2 s_{n_{k-1}})^{1/2}}\right) s_{n_k} (2 \log_2 s_{n_{k}}^2)^{1/2}\biggr).
\end{eqnarray*}

\bigskip \noindent So, for any $0<\delta^{\prime \prime}<\delta^{\prime}$, we have for large values of $k$,

$$
\mathbb{P}\biggr(S_{n_{k}}^{\ast} > (1+\delta^{\prime}) (2\log_2 s_{n_{k}}^2)^{1/2}\biggr) \leq 2 
\mathbb{P}\biggr(S_n^{\ast} > (1+\delta^{\prime \prime}) s_{n_k} (2\log_2 s_{n_{k-1}}^2)^{1/2}\biggr).
$$

\noindent At this step, let us apply the exponential inequality, Statement (i) in Theorem \ref{expoboundth02}, to have, with $\varepsilon_{n_k}=(1+\delta^{\prime \prime})(2\log_2 s_{n_{k}})^{1/2}$ and $c_{n_k}$. Since, by assumption,  $c_{n_k}\varepsilon_{n_k} \rightarrow 0$ and for $k$ large to ensure based on $c_{n_k}\varepsilon_{n_k} <1$, the last formula yields

\begin{eqnarray*}
\mathbb{P}\biggr(S_{n_{k}}^{\ast}/s_{n_k} > \varepsilon_{n_k}) &<&\exp(-\frac{\varepsilon_{n_k}^2}{2} (1+\varepsilon_{n_k} c_{n_k}/2) (2\log_2 s_{n_{k-1}})^{1/2}\biggr)\\
&\leq & \exp\biggr(-\frac{\varepsilon_{n_k}^2}{2}\biggr)\\ 
&=& \exp\biggr(-\frac{(1+\delta^{\prime \prime})^2 (2\log_2 s_{n_{k}}}{2})\biggr)\\ 
&\leq & \exp\biggr(- (1+\delta^{\prime \prime}) \log_2 s_{n_{k}}^2\biggr)\\
&=& \frac{1}{(1+\delta)\log s_{n_{k}}^2} \sim \frac{1}{\biggr(2k \log (1+\delta^{\prime \prime})\biggr)^{(1+\delta^{\prime \prime})}}.\\
\end{eqnarray*}

\bigskip \noindent Since the last term in the group of formulas above is the general term of a converging series, we also see that the series of general term

$$
\mathbb{P}\biggr(S_n > (1+\delta) s_n (2\log_2 s_n^2)^{1/2}, \ i.o.\biggr)
$$

\bigskip \noindent also converges. By Point (i) of Borel Cantelli's Lemma \ref{probab_02_indep_sec_01_lem01}, we have

$$
\mathbb{P}\biggr(S_{n_{k}}^{\ast} > (1+\delta^{\prime \prime}) s_{n_k} (2\log_2 s_{n_{k}}^2)^{1/2}, \ i.o \ (in \ k) \ \biggr)=0,
$$

\bigskip \noindent and by the bounds and inclusions that are proved above, we have
 
$$
\mathbb{P}\biggr(S_n > (1+\delta) s_n (2\log_2 s_n^2)^{1/2}, \ i.o.\biggr)=0,
$$

\bigskip \noindent that is, for any arbitrary $\delta>0$, we have
$$
\limsup_{n\rightarrow +\infty} \frac{S_n}{\sqrt{2s_n^2 log_2 s_n^2}}\leq 1+\delta, \ a.s.
$$

\bigskip \noindent which proves that

$$
\limsup_{n\rightarrow +\infty} \frac{S_n}{\sqrt{2s_n log_2 s_n^2}}\leq 1. \ a.s.
$$

\bigskip \noindent To prove that this superior limit is also greater than one, we just remark that the first part of the proof applied to 
the $-X_n$'s with $\delta=1$ leads to 

$$
\mathbb{P}(\{\omega \in \Omega, \ \exists N(\omega), \ \forall n\geq N, \ -S_n \leq 2(2s_n^2 log_2 s_n^2)^{1/2})\})=1
$$

\bigskip \noindent Set

$$
\Omega_{1}^{\ast}=\{\omega \in \Omega, \ \exists N(\omega), \ \forall n\geq N, \ S_n \leq 2(2s_n^2 log_2 s_n^2)^{1/2})\}.
$$

\bigskip \noindent Now, let $\delta>0$ and let $c>1$ such that

$$
\frac{\delta^{\prime}}{2}=1-\biggr((1-\delta) \left(1 -\frac{1}{c^2}\right)^{1/2} - \frac{2}{c}\biggr)>0
$$

\bigskip \noindent so that

$$
\biggr((1-\delta) \left(1 -\frac{1}{c^2}\right)^{1/2} - \frac{2}{c}\biggr)>1-\delta^{\prime}. \ (K1)
$$

\bigskip \noindent Now we select, By Formula (C3), a sub-sequence $(s_{n_k})_{k\geq 1}$  of $(s_{n})_{n\geq 1}$ such that $s_{n_k} \sim c^k$, as $k\rightarrow +\infty$. Next, we wish to apply Theorem \ref{expoboundth02} to the non-overlapping spacings $Y_{n_k}=S_{n_k}-S_{n_{k-1}}$'s of variances $u_{n_k}=s_n^2-s_{n-1}^2$. We immediately check that
$$
u_{n_k}=s_n^2-s_{n-1}^2 \sim \left(1 - \frac{1}{c^2}\right) \ and \ v_{n_k}^2=(2 \log_2 u_{n_k}^2) \sim (2 \log_2 u_k^2). \ (K2)
$$

\bigskip \noindent We have to remark that we still have, as $k+\infty$,

\begin{eqnarray*}
b_{{n_k}}&=&\frac{1}{u_{n_k}}\max_{n_{k-1} <n\leq n_k} |X_{n_k}|\\
&=&\max_{n_{k-1} <n<n_k} |Y_{n_k}|\\
&\leq& \frac{1}{u_{n_k}}\max_{1 \leq n_k} |X_{n_k}|\\
&=&O\biggr((2\log_2 s_{n_{k}}^{2})^{-1/2})\biggr) \rightarrow 0.
\end{eqnarray*}

\bigskip \noindent because of the first part of Formula (K2) above. Now, we are in the position to re-conduct the same method to the $Y_{n_k}$. So for $0<\delta$, $\varepsilon_{n_k}=(1-\delta)v_{n_k} \rightarrow +\infty$. We take $\gamma=(1-\delta)^{-1}$. Applying Theorem \ref{expoboundth02} with that value of $\gamma$ for large values of $k$ so that  we have $b_{{n_k}}$ is small enough and 
$(1-\delta) (2\log_2 u_{n_{k}}^2)^{1/2}$ is large enough, leads to

\begin{eqnarray*}
\mathbb{P}\biggr(Y_{n_{k}}/u_{n_k} > (1+\delta) (2\log_2 u_{n_{k}}^2)^{1/2}\biggr)&>& \exp\left(-\frac{\varepsilon_{n_k}^2}{2}\right)\\
&>&\exp\left((1-\delta) \log_2(u_{n_{k}}^2)\right)\\
&=&\frac{1}{(2k \log u_{n_k})^(1-\delta)}\\
&\sim& \frac{1}{\biggr(2k \log (1+\delta)\biggr)^{(1-\delta)}}.
\end{eqnarray*}

\bigskip \noindent Since the series of general term

$$
\mathbb{P}\biggr(Y_{n_{k}}/u_{n_k} > (1-\delta_1) (2\log_2 u_{n_{k}}^2)^{1/2}\biggr)
$$

\bigskip \noindent diverges and the $Y_{n_{k}}$ are independent, we have by Point (i) of Borel-Cantelli Lemma \ref{probab_02_indep_sec_01_lem01}.

$$
\mathbb{P}\biggr(Y_{n_{k}}/u_{n_k} > (1-\delta) (2\log_2 u_{n_{k}}^2)^{1/2}, \ i.o.\biggr)=1.
$$ 

\bigskip \noindent Set

$$
\Omega_{1}^{\ast}=(Y_{n_{k}}/u_{n_k} > (1-\delta) (2\log_2 u_{n_{k}}^2)^{1/2}, \ i.o.\biggr).
$$

\bigskip \noindent On $\Omega^{\ast}=\Omega_{1}^{\ast} \cap \Omega_{2}^{\ast}$, we have

\begin{eqnarray*}
S_{n_{k}}&=&S_{n_{k}}-S_{n_{k-1}} +S_{n_{k-1}} > (1-\delta) u_{k} (2\log_2 u_{{k}}^2)^{1/2}+S_{n_{k-1}}, \ i.o\\
&\Rightarrow& (S_{n_{k}}  >(1-\delta) s_{n_k} (2\log_2 s_{s_{k}}^2)^{1/2} - 2  s_{{n_{k-1}}} (2\log_2 s_{n_{k-1}}^2)^{1/2}), \ i.o\\
&\Rightarrow& (S_{n_{k}} >(1-\delta)(s_{{n_k}} (2\log_2 s_{n_{k}}^2)^{1/2})\\
&\times& \biggr(\frac{u_{k} (2\log_2 u_{{k}}^2)^{1/2}}{s_{{n_k}} (2\log_2 s_{n_{k}}^2)^{1/2})}-
\frac{s_{{n_{k-1}}} (2\log_2 s_{n_{k-1}}^2)^{1/2})}{s_{{n_k}} (2\log_2 s_{n_{k}}^2)^{1/2}}\biggr), \ i.o. \ \ (L45)
\end{eqnarray*}

\bigskip \noindent But, as $k\rightarrow +\infty$
$$
\biggr(\frac{u_{k} (2\log_2 u_{{k}}^2)^{1/2}}{s_{{n_k}} (2\log_2 s_{n_{k}}^2)^{1/2})}-
\frac{2s_{{n_{k-1}}} (2\log_2 s_{n_{k-1}}^2)^{1/2})}{s_{{n_k}} (2\log_2 s_{n_{k}}^2)^{1/2}}\biggr)
$$

\bigskip \noindent converges to
$$
\biggr((1-\delta) \left(1 -\frac{1}{c^2}\right)^{1/2} - \frac{2}{c}\biggr)>1-\delta^{\prime}.
$$

\bigskip \noindent We conclude that Line (L41) above that

$$
\Omega^{\ast} \subset (S_{n_{k}} >(1-\delta) s_{{n_k}} (2\log_2 s_{n_{k}}^2)^{1/2}, \ i.o).
$$

\bigskip \noindent Since $\mathbb{P}(\Omega^{\ast})=1$, we get that for any $\delta>0$, we have

$$
\mathbb{P}((S_{n_{k}} >(1-\delta) s_{{n_k}} (2\log_2 s_{n_{k}}^2)^{1/2}, \ i.o)=1.
$$

\bigskip \noindent The proof of the theorem is now complete. $\square$\\

%% file: proba_02_08_ang.tex
\chapter{Conditional Expectation} \label{proba_02_ce}

\section{Introduction and definition} \label{proba_02_ce_sec_01}

\noindent We already saw in Chapter \ref{probab_02_indep} the key role played independence in Probability Theory. But a very great part, even the greatest part, among studies in Probability Theory rely 
on some kind on dependence rather that on independence.  However, the notion of independence, in most situations, is used as a theoretical modeling tool or as an approximation method. Actually, many methods which are used to handle dependence are based transformation of independent objects or based on some \textit{nearness} measure from the independence frame. So, the better one masters methods based on independence, the better one understands methods for dependence studies.\\

\noindent However there is a universal tool to directly handle dependence, precisely the \textit{Conditional Mathematical Expectation} tool. This chapter which is devoted to it, is the door  for the study of arbitrary sequences or family of random objects.\\

\noindent The most general way to deal and to introduce to this tool relies on the Radon-Nikodym Theorem as stated in Doc 08-01 in Chapter 9 in \cite{ips-mestuto-ang}. We already spoke a little on it in the first lines in Chapter \ref{proba_02_ineg}.\\

\begin{definition} \label{proba_02_ce_def01} Let $X$ be real-valued random variable $(\Omega, \mathcal{A}, \mathbb{P})$ which is quasi-integrable, that is $\int X^{+}d\mathbb{P}<\infty$ for example. Let $\mathcal{B}$ be a $\sigma$-sub-alg\`{e}bra of $\mathcal{A}$, meaning that $\mathcal{B}$ is a $\sigma$-algebra of subsets of $\Omega$ and $\mathcal{B}\subset \mathcal{A}$. The mapping 
\begin{equation*}
\begin{array}{cccc}
\phi : & \mathcal{B} & \longrightarrow & \overline{\mathbb{R}} \\ 
& A & \hookrightarrow & \phi (B)=\int_{B}Xd\mathbb{P}
\end{array}
\end{equation*}

\bigskip \noindent is well-defined and continuous with respect to $\mathbb{P}$ in the following sense 
\begin{equation*}
(\forall B\in \mathcal{B}), \ (\mathbb{P}(B)=0\Longrightarrow \phi (B)=0).
\end{equation*}

\bigskip \noindent By Radon-Nikodym's Theorem (See Doc 08-01 in Chapter 9 in \cite{ips-mestuto-ang}), there exists a random variable, uniquely defined \textit{a.s.}, 
$$
Z \ : \  (\Omega ,\mathcal{B}) \mapsto \overline{\mathbb{R}}
$$

\bigskip \noindent which is $\mathcal{B}$ -measurable such that 

\begin{equation}
(\forall B\in \mathcal{B}), \left(\int_{B} \ X \ d\mathbb{P=}\int_{B}\ Z \ d\mathbb{P}\right). \label{1.1}
\end{equation}

\bigskip \noindent This random variable $Z$, is defined as the mathematical expectation of $X$ with respect to $\mathcal{B}$ and denoted by by
\begin{equation*}
Z=\mathbb{E}(X/\mathcal{B})=\mathbb{E}^{\mathcal{B}}(X) \ a.s.,
\end{equation*}

\noindent * and the mathematical expectation is uniquely $\mathcal{B}$-\textit{almost surely}.\\

\noindent Let $Y$ be a measurable mapping $Y (\Omega, \mathcal{A})\rightarrow (G,\mathcal{D})$, where $(G,\mathcal{D})$ is an arbitrary measurable space. As previously explained in
the first lines in Section  \ref{probab_02_sllnKol} in Chapter \ref{proba_02_ineg};
$$
\mathcal{B}_{Y} =\{Y^{-1}(H),H\in \mathcal{D}\}
$$ 

\noindent is the $\sigma$-algebra generated by $Y$, the smallest one rendering $Y$ measurable. It is a $\sigma$-sub-algebra of $\mathcal{A}$. The mathematical expectation with respect to $Y$, denoted by $\mathbb{E}(X/Y)$, is the mathematical expectation with respect to $\mathcal{B}_{Y}$ that is

\begin{equation*}
\mathbb{E}(X/Y)=\mathbb{E}(X/\mathcal{B}_{Y}).
\end{equation*}
\end{definition}

\bigskip \noindent \textbf{Extension}. Later we will define the mathematical expectation with respect to a family a measurable mappings similarly to the one with respect to one mapping $Y$ as in the definition.\\

\noindent This definition is one of the most general ones. In the special case where we work with square integrable real-valued random variables, a specific definition based on the orthogonal projection
on the closed linear space $H=L^2(\Omega,\mathcal{B}\mathbb)$ of $\mathcal{B}$-measurable random variables is possible. And the mathematical expectation  of $X$ is its orthogonal projection on $H$. We will see this in Section  \ref{proba_02_ce_sec_01}. But as we will see it, even in the general case, the mathematical expectation is still a linear projection $L^1$ as explained in Remark \ref{rem1} below.

\section{The operator of the mathematical expectation}  \label{proba_02_ce_sec_02}

\noindent We already knew in Doc 08-01 in Chapter 9 in \cite{ips-mestuto-ang}, that $\mathbb{E}(X/\mathcal{B})$ is \textit{a.s.} finite (in the frame of Probability Theory) and is integrable whenever $X$ is. In stating the properties below, we fix $\mathbb{B}$ and we do not need write the mention of \textit{with respect to $\mathcal{B}$}.\\    
 
\bigskip \noindent We have

\begin{proposition} \label{proba_02_ce_prop01}
Considered as an operator from $L^{1}(\Omega ,A,\mathbb{P})$ to $L^{1}(\Omega,\mathcal{B},\mathbb{P})$, the mathematical expectation mapping
\begin{equation*}
\begin{array}{cccc}
L & L^{1}(\Omega ,A,\mathbb{P}) & \longrightarrow & L^{1}(\Omega ,\mathcal{B%
},\mathbb{P}) \\ 
& X & \hookrightarrow & \mathbb{E}(X/\mathcal{B)}%
\end{array}%
\end{equation*}

\bigskip \noindent is linear and satisfies
\begin{equation*}
L^{2}=L\ and \left\Vert L\right\Vert =1.
\end{equation*}

\bigskip \noindent It is non-negative in the following sense
\begin{equation*}
X\geq 0\Longrightarrow L(X)\geq 0.
\end{equation*}

\bigskip \noindent The operator $L$ is non-decreasing in the sense that, for $(X,Y) \in (L^{1})^2$,

$$
X \leq Y \ a.s. \ \Rightarrow  L(X)\leq L(Y), \ a.s.
$$
\end{proposition}

\bigskip \noindent \textbf{Proof}. Let $X$ and $Y$ be two integrable random variables, defined both on $(\Omega,\mathcal{A})$,  
$\alpha$ and $\beta$ two real numbers. Then $\alpha X+\beta Y$ and $\alpha \mathbb{E}(X/\mathcal{B})+\beta \mathbb{E(}Y/\mathcal{B})$
are \textit{a.s.} defined. Then for any $B \in \mathcal{B}$, we have
\begin{equation*}
\int_{B}(\alpha X+\beta Y)d\mathbb{P}=\alpha \int_{B}Xd\mathbb{P}+\beta
\int_{B}Yd\mathbb{P},
\end{equation*}

\bigskip \noindent which, by the the definition of the mathematical expectation,  implies
\begin{equation*}
\alpha \int_{B}\mathbb{E}(X/\mathcal{B})d\mathbb{P}+\beta \int_{B}\mathbb{%
E}(Y/\mathcal{B})d\mathbb{P}=\int_{B}\{\alpha \mathbb{E}(X/\mathcal{B}%
)+\beta \mathbb{E(}Y/\mathcal{B})\}d\mathbb{P}.
\end{equation*}

\bigskip \noindent Since $\alpha \mathbb{E}(X/\mathcal{B})+\beta \mathbb{E(}Y/\mathcal{B})$ is $\mathcal{B}$-measurable, the equality entails that

\begin{equation*}
\mathbb{E}((\alpha X+\beta Y)/\mathcal{B})=\alpha \mathbb{E}(X/\mathcal{B})+\mathbb{E}(Y/\mathcal{B}), \ a.s.
\end{equation*}

\bigskip \noindent Next, as an immediate consequence of the definition, $\mathbb{E}(X/\mathcal{B})=X$ \textit{a.s.} whenever $X$ is $\mathcal{B}$-measurable. 
Since $L(X)=\mathbb{E}(X/\mathcal{B})$ is $\mathcal{B}$-measurable, it comes that 
\begin{equation*}
L^{2}(X)=\mathbb{E}(L(X)/\mathcal{B})=L(X),
\end{equation*}

\bigskip \noindent which implies that $L^{2}=L$. The non-negativity comes from Radon-Nikodym's Theorem which says the if $X$ is non-negative, the mapping
\begin{equation*}
\mathcal{B}\ni B\longmapsto \phi(B)=\int_{B}Xd\mathbb{P}=\int_{B}\ Z \ d\mathbb{P}
\end{equation*}

\bigskip \noindent is non-negative. Hence the its Radon-Nikodym derivative $Z$ is non-negative. Here is an easy proof. Indeed $B_a=(Z<a)\in \mathcal{B}$, $a<0$. Fix $a<0$. If $B_a$ is not a null-set, we would have 
\begin{equation*}
\phi(B)=\int_{B_a}\ X \ dP=\int_{B_a}\ Z\ dP\leq a \mathbb{P}(B_a)<0.
\end{equation*}

\bigskip \noindent This is impossible since we have $\phi (B_a)\geq 0$. So for $k\geq 1$, all the events $B_{-1/k}$ are null-sets. Since 
$$
(Z\leq 0)=\bigcap_{k\geq 1} B_{-1/k},
$$ 

\bigskip \noindent $(Z\leq 0)$ is a null-event and thus, $Z\geq 0$ \textit{a.s.}.\\

\noindent The non-decreasingness is immediate from the combination of the linearity and the non-negativity.\\

\noindent Let us determine the norm of $L$ defined by 
\begin{equation*}
\left\Vert L\right\Vert =\sup \{\left\Vert L(X)\right\Vert_{1}/\left\Vert X\right\Vert _{1}, \ X \in L^{1}, \ \left\Vert X\right\Vert_{1} \neq 0\}.
\end{equation*}

\bigskip \noindent First let us show that $L$ is contracting, that is

\begin{equation*}
|\mathbb{E}(X/\mathcal{B})|\leq \mathbb{E}(|X|/\mathcal{B}).
\end{equation*}

\bigskip \noindent We have $X\leq |X|$ and $-X\leq |X|$. We have $\mathbb{E}(X/\mathcal{B})\leq 
\mathbb{E}(\left\vert X\right\vert /\mathcal{B})$ and $-\mathbb{E}(X/\mathcal{%
B})=\mathbb{E}(-X/\mathcal{B})\leq \mathbb{E}(|X|/\mathcal{B})\geq 0$.\\

\noindent We conclude that  
\begin{equation*}
|\mathbb{E}(X/\mathcal{B})|\leq \mathbb{E}(|X|/\mathcal{B}).
\end{equation*}

\bigskip \noindent By applying the definition of the mathematical expectation, we have

\begin{equation*}
\left\Vert \mathbb{E}(X/\mathcal{B})\right\Vert _{1}=\mathbb{E}|\mathbb{E}(X/%
\mathcal{B})|\leq \mathbb{EE}(|X|/\mathcal{B})=\mathbb{E}(|X|)=\left\Vert
X\right\Vert _{1}.
\end{equation*}

\bigskip \noindent Hence we have 
\begin{equation*}
\left\Vert L(X)\right\Vert _{1}\leq \left\Vert X\right\Vert _{1}.
\end{equation*}

\bigskip \noindent Next we have
\begin{equation*}
\left\Vert L\right\Vert \leq 1.
\end{equation*}

\bigskip \noindent But if $\left\Vert X\right\Vert _{1}\neq 0$ and $X$ is $\mathcal{B}$-measurable, we have $L(X)=X$ and for such random variables, we have$\left\Vert
L(X)\right\Vert _{1}/\left\Vert X\right\Vert _{1}=1$. We conclude that
\begin{equation*}
\left\Vert L\right\Vert =1.
\end{equation*}

\bigskip \noindent The proof of the proposition is complete. $\blacksquare$\\

\begin{remark} \label{rem1} A linear operator $L$ such that $L^2=L$ is called a projection. Hence the mathematical operator is a projection
of $L^1$  to the sub-space of $\mathcal{B}$-measurable functions.
\end{remark}

\section{Other Important Properties}  \label{proba_02_ce_sec_03}

\begin{proposition} \label{probab_02_propAllME}
We have the following facts.\\

\noindent (1) If  $X$ is $\mathcal{B}$-measurable, then  $\mathbb{E}(X/\mathcal{B})=X$ a.s.\\

\noindent (2) The mathematical conditional expectation is anon-negative linear and non-decreasing operator.\\

\noindent (3) The mathematical conditional expectation is a contracting operator, that is, whenever the expressions make sense,

\begin{equation*}
|\mathbb{E}(X/\mathcal{B})|\leq \mathbb{E}(|X|/\mathcal{B}).
\end{equation*}

\bigskip \noindent (4) Let $\mathcal{B}_{1}$ and $\mathcal{B}_{2}$ be two $\sigma$-sub-algebras of $\mathcal{A}$ with $\mathcal{B}_{1}\subset \mathcal{B}_{2}\subset \mathcal{A}$. We have

\begin{equation}
\mathbb{E}(\mathbb{E}(X/\mathcal{B}_{2})/\mathcal{B}_{1})=\mathbb{E}(X/%
\mathcal{B}_{1})  \label{6a}
\end{equation}

\noindent and \\

\begin{equation}
\mathbb{E}(\mathbb{E}(X/\mathcal{B}_{1})/\mathcal{B}_{2})=\mathbb{E}(X/%
\mathcal{B}_{1}).  \label{6b}
\end{equation}

\noindent (5) Let $X$ be a random variable independent of  $\mathcal{B}$ in the following sense : for any mapping $\mathcal{B}$-measurable mapping $Z:\Omega \longmapsto \mathbb{R}$ and for any measurable application $h \ : \ \mathbb{R}\longmapsto \mathbb{R}$, 

\begin{equation*}
\mathbb{E}(Z\times \ h(X))=\mathbb{E}(Z)\times \ \mathbb{E}(h(X)).
\end{equation*}

\bigskip \noindent Then,  if $\mathbb{E}(X)$ exists, we have

\begin{equation*}
\mathbb{E}(X/\mathcal{B})=\mathbb{E}(X).
\end{equation*}

\bigskip \noindent (6) (Monotone Convergence Theorem for Mathematical expectation). Let $(X_{n})_{n\geq 0}$ be a non-decreasing sequence of integrable random variables which are all non-negative or all integrable. Then we have

\begin{equation*}
\mathbb{E}(\lim_{n\rightarrow \infty }X_{n}/\mathcal{B})=\lim_{n\rightarrow
\infty }\mathbb{E}(X_{n}/\mathcal{B})
\end{equation*}

\bigskip \noindent (7) (Fatou-Lebesgue Theorems). Let $(X_{n})_{n\geq 0}$ be quasi-integrable real-valued random variables which is \textit{a.s.} bounded below by an integrable random variable, then

$$
\mathbb{E}(\liminf_{n\rightarrow \infty }X_{n}/\mathcal{B})\leq \liminf_{n\rightarrow \infty }\mathbb{E}(X_{n}/\mathcal{B}).
$$

\bigskip \noindent If the sequence is \textit{a.s.} bounded above by an integrable random variable, then

$$
\mathbb{E}(\limsup_{n\rightarrow \infty }X_{n}/\mathcal{B})\geq \limsup_{n\rightarrow \infty }\mathbb{E}(X_{n}/\mathcal{B}).
$$

\bigskip \noindent If the sequence is uniformly \textit{a.s.} bounded by an integrable random variable $Z$ and converges \textit{a.s.} to $X$, then

$$
\lim_{n\rightarrow \infty }\mathbb{E}(X_{n}/\mathcal{B})= \mathbb{E}(X /\mathcal{B})
$$

\bigskip \noindent and

$$
\mathbb{E}(|X|/\mathcal{B}) \leq \mathbb{E}(|Z|/\mathcal{B}).
$$

\bigskip \noindent (8) Let $X$  be a quasi-integrable random variable. Let $Z$ be $\mathcal{B}$-measurable and non-negative or integrable. Then we have 
\begin{equation*}
\mathbb{E}(ZX/\mathcal{B})=Z\times \mathbb{E}(X/\mathcal{B}).
\end{equation*}

\end{proposition}

\bigskip \noindent \textbf{Proof}.\\

\noindent Points from (1) to (3) are already proved in the first proposition.\\

\noindent \textbf{Proof of Point (4)}. First we know that $\mathbb{E}(X/\mathcal{B}_{1})$ is $\mathcal{B}_{1}$-measurable  and thus $\mathcal{B}_{2}$-measurable. By Point (1), we have

\begin{equation*}
\mathbb{E}(\mathbb{E}(X/\mathcal{B}_{1})/\mathcal{B}_{2})=\mathbb{E}(X/\mathcal{B}_{1}).
\end{equation*}

\bigskip \noindent Formula (\ref{6a}) is proved. Next, for any $B\in \mathcal{B}_{1}\subset 
\mathcal{B}_{2}$, we have

\begin{equation*}
\int_{B}Xd\mathbb{P}=\int_{B}\mathbb{E}(X/\mathcal{B}_{2})d\mathbb{P},
\end{equation*}

\bigskip \noindent since $B$ is also in $\mathcal{B}_{2}$. Now we apply the definition of the mathematical expectation with respect to 
$\mathcal{B}_{1}$ in the right-hand member to have

\begin{equation*}
\int_{B}Xd\mathbb{P}=\int_{B}\mathbb{E}(X/\mathcal{B}_{2})d\mathbb{P}=\int_{B}\mathbb{E}(\mathbb{E}(X/\mathcal{B}_{2})/\mathcal{B}_{1})d\mathbb{P}.
\end{equation*}

\bigskip \noindent Since $\mathbb{E}(\mathbb{E}(X/\mathcal{B}_{2})/\mathcal{B}_{1})$ is $\mathcal{B}_{1}$-measurable, we conclude that 
$\mathbb{E}(X/\mathcal{B}_{1})=\mathbb{E}(\mathbb{E}(X/\mathcal{B}_{2})/\mathcal{B}_{1})$ \textit{a.s.} Thus, we reach Formula (\ref{6b}).\\

\noindent \textbf{Proof of (5)}. It is clear that the constant mapping $\omega \longmapsto \mathbb{E}(X)$ is $\mathcal{B}$-measurable. Hence for any $B \in \mathcal{B}$,
\begin{equation*}
\int_{B} \ X \ d\mathbb{P}=\mathbb{E}(1_{B}\times \ X)=\mathbb{E}(1_{B})\times \mathbb{E}(X)=\mathbb{E}(X)\int_{B}d\mathbb{P}=\int_{B}\mathbb{E}(X)d\mathbb{P},
\end{equation*}

\bigskip \noindent which proves that
  
\begin{equation*}
\mathbb{E}(X/\mathcal{B})=\mathbb{E}(X), \ a.s.
\end{equation*}

\noindent \textbf{Proof of Point (6)}. Let $(X_{n})_{n\geq 0}$ be a non-decreasing sequences of random variables such that $\mathbb{E}(X_{n}^{+})<\infty$ pour tout $n\geq 0$. For any $B\in \mathcal{B}$, we have
\begin{equation*}
\int_{B}X_{n}d\mathbb{P}=\int_{B}\mathbb{E}(X_{n}/\mathcal{B})d\mathbb{P}.
\end{equation*}

\bigskip \noindent Since the sequences $(X_{n})_{n\geq 0}$ and $(\mathbb{E}(X_{n}/\mathcal{B}))_{n\geq 0}$ are non-decreasing of integrable random variable, we may apply the Monotone Convergence Theorem
to get 

\begin{equation*}
\int_{B}\lim_{n\rightarrow \infty }X_{n}d\mathbb{P}=\lim_{n\rightarrow
\infty }\int_{B}X_{n}d\mathbb{P}=\lim_{n\rightarrow \infty }\int_{B}\mathbb{E%
}(X_{n}/\mathcal{B})d\mathbb{P=}\int_{B}\lim_{n\rightarrow \infty }\mathbb{E}%
(X_{n}/\mathcal{B})d\mathbb{P}.
\end{equation*}

\bigskip \noindent Since $\lim_{n\rightarrow \infty }\mathbb{E}(X_{n}/\mathcal{B})d\mathbb{P}$ is $\mathcal{B}$-measurable, we have

\begin{equation*}
\mathbb{E}(\lim_{n\rightarrow \infty }X_{n}/\mathcal{B})=\lim_{n\rightarrow
\infty }\mathbb{E}(X_{n}/\mathcal{B}).
\end{equation*}

\bigskip \noindent \textbf{Proof of Point (7)}. Based on the Monotone convergence Theorem for conditional expectation in the previous Point (6), the Fatou-Lebesgue Theorem and the Lebesgue Dominated Theorem are proved as in the unconditional case, as done in Chapter 6 in \cite{ips-mestuto-ang}. $\blacksquare$\\

\bigskip \noindent \textbf{Proof of Point (8)}. Let $Z$ be a $\mathcal{B}$-measurable random variable non-negative or integrable. Thus $Z\mathbb{E}(X/\mathcal{B})$ is  $\mathcal{B}$-measurable. Suppose that  $Z=1_{C}$, $C\in \mathcal{B}$, that is $Z$ is $\mathcal{B}$-measurable indicator function. We have for any $B\in \mathcal{B}$, 

\begin{eqnarray*}
\int_{B} \ Z\ X \ d\mathbb{P}&=&\int_{B}1_{C}\ X \ d\mathbb{P}=\int_{BC}\ X\ d\mathbb{P}\\
&=&\int_{BC}\mathbb{E}(X/\mathcal{B})d\mathbb{P}=\int_{B}1_{C}\mathbb{E}(X/\mathcal{B})d\mathbb{P}.
\end{eqnarray*}

\bigskip \noindent Thus, we have
\begin{equation*}
\int_{B}\ Z \ X\ d\mathbb{P=}\int_{B} \ Z \ \mathbb{E}(X/\mathcal{B})d\mathbb{P}.
\end{equation*}

\bigskip \noindent Since $Z\times \mathbb{E}(X/\mathcal{B})$ is $\mathcal{B}-$measurable, we get  
\begin{equation*}
\mathbb{E}(ZX/\mathcal{B})=Z\mathbb{E}(X/\mathcal{B}).
\end{equation*}

\noindent To finish the proof, we follow the famous three steps method by extending the last formula to elementary functions based on $\mathcal{B}$-measurable sets, next to non-negative random variables using Point (6) above and finally to an arbitrary random variable $Z$ using the additivity of both the expectation and the conditional expectation.

\bigskip * We may and do have the same theory by using non-negative $\mathcal{B}$-measurable random variables in place of the elements of $\mathcal{B}$ in the definition of the mathematical expectation.

\section{Generalization of the definition}  \label{proba_02_ce_sec_04}

Let us define by $\mathbb{E}(X/ \mathcal{B})$ a $\mathcal{B}$-measurable random variable such that for all $\mathcal{B}$-measurable mapping $h : (\Omega,\mathcal{B})\longmapsto (\mathbb{R},\mathcal{B})$, we have  
\begin{equation}
\int h\ X \ d\mathbb{P=}\int h \mathbb{E}(X/ \mathcal{B}) \ d\mathbb{P}  \label{4.2}
\end{equation}

\noindent We are going to quickly show that the definitions based on Formulas \eqref{4.2} and \eqref{1.1} respectively are the same.\\

\noindent Before we do it, let us just say that Formula \eqref{4.2} usually offers a more comfortable handling of the mathematical expectation.\\

\noindent \textbf{Proof of the equivalence between Formulas \eqref{4.2} and \eqref{1.1}}. The implication 
$(\ref{4.2})\Longrightarrow (\ref{1.1})$ by taking $h=1_{B}$ for $B \in \mathcal{B}$. To prove the converse implication, we use the classical three steps methods. Suppose that Formula \eqref{1.1} holds.\\

\noindent \textit{Step 1}. If $h=1_{B}$ for $B\in \mathcal{B}$, Formula \eqref{4.2} is obvious.\\

\noindent \textit{Step 2}. $h$ is an elementary function of the form
\begin{equation*}
h=\sum_{i=1}^{p}\alpha _{i}1_{B_{i}},
\end{equation*}

\bigskip \noindent where $B_{i}\in \mathcal{B}$ et $\alpha _{i}$ $\in \mathbb{R}$. By using the linearity, we have

\begin{eqnarray*}
\int h \ X \ d\mathbb{P}&=&\int \left( \sum_{i=1}^{p}\alpha _{i}1_{B_{i}}\right) X \ d\mathbb{P}\\
&=&\sum_{i=1}^{p}\alpha _{i}\left( \int 1_{B_{i}}X \ d\mathbb{P}\right)\\
&=&\sum_{i=1}^{p}\alpha _{i}\left( \int 1_{B_{i}}\mathbb{E}(X/\mathcal{B}) \ d\mathbb{P}\right)\\
&=&\int \left( \sum_{i=1}^{p}\alpha _{i}1_{B_{i}}\right) \mathbb{E}(X/\mathcal{B})d\mathbb{P=}\int h\mathbb{E}(X/\mathcal{B})\ d\mathbb{P}.
\end{eqnarray*}

\bigskip \noindent \textit{Step 3}. $h$ is $\mathcal{B}$-measurable and non-negative There exists a sequence of elementary  $(h_{n})_{n\geq 0}$ based on elements of $\mathcal{B}$ such that
 $h_{n}\nearrow h$ and thus, 
\begin{equation*}
h_{n}X^{+}\nearrow hX^{+}\text{ and }h_{n}\mathbb{E}(X^{+}/\mathcal{B}%
)\nearrow h\mathbb{E}(X^{+}/\mathcal{B}).
\end{equation*}

\bigskip \noindent By the monotone convergence Theorem, we have

\begin{equation*}
\int h_{n}X^{+}d\mathbb{P}=\int h_{n}\mathbb{E}(X^{+}/\mathcal{B})d\mathbb{P} \nearrow \int hX^{+}d\mathbb{P}=\int h\mathbb{E}(X^{+}/\mathcal{B})d\mathbb{P}.
\end{equation*}

\bigskip \noindent We similarly get 
\begin{equation*}
\int h \ X^{-} \ d\mathbb{P}=\int h\mathbb{E}(X^{-}/\mathcal{B})d\mathbb{P}.
\end{equation*}

\noindent Thus by quasi-integrability, we have 
\begin{eqnarray*}
\int h \ X \ d\mathbb{P}&=&\int h\mathbb{E}(X^{+}/\mathcal{B})d\mathbb{P}-\int h \ X^{-}\ d\mathbb{P}\\
&=&\int h\mathbb{E}(X^{+}/\mathcal{B})d\mathbb{P-}\int h \ \mathbb{E}(X^{-}/\mathcal{B})d\mathbb{P}=\int h\mathbb{E}(X/\mathcal{B}) \ d\mathbb{P}.
\end{eqnarray*}

\noindent The proof is over. $\square$\\

\noindent With the second definition, some properties are easier to prove as the following one.

\begin{proposition}
if $Z$ is $\mathcal{B}$-measurable either non-negative or integrable, we have for any quasi-integrable random variable, 
\begin{equation*}
\mathbb{E}(Z\times \ X/\mathcal{B})=Z\times \ \mathbb{E}(X/\mathcal{B}).
\end{equation*}
\end{proposition}

\bigskip \noindent \textbf{Proof}. Let $h$ be any non-negative and real valued $\mathcal{B})$-measurable function. We have

\begin{equation*}
\int h\{ZX\}d\mathbb{P}=\int \{hZ\}Xd\mathbb{P}=\int \{hZ\}\mathbb{E}(X/%
\mathcal{B})d\mathbb{P}=\int h\left\{ Z\mathbb{E}(X/\mathcal{B})\right\} d%
\mathbb{P}.
\end{equation*}

\noindent Since $Z\mathbb{E}(X/\mathcal{B})$ is $\mathcal{B}$-measurable, we get that $\mathbb{E}(ZX/\mathcal{B})=Z\mathbb{E}(X/\mathcal{B})$ \textit{a.s.} $\square$\\

\section{Mathematical expectation with respect to a random variable}  \label{proba_02_ce_sec_05}

\noindent Let us consider that $\mathcal{B}$ is generated by a measurable mapping $Y \ : \ (\Omega,,\mathcal{A})\rightarrow (G,\mathcal{D})$, where $(G,\mathcal{D})$ is an arbitrary leasure space, that is

\begin{equation*}
\mathcal{B}=\mathcal{B}_{Y}=\{Y^{-1}(H),H\in \mathcal{D}\}
\end{equation*}

\bigskip \noindent Par definition, we denote

$$
\mathbb{E}(X/\mathcal{B}_{Y})=\mathbb{E}(X/Y).
$$

\bigskip \noindent Let us prove that any real valued and $\mathcal{B}_{Y}$-measurable mapping $h$ is of the form $g(Y)$,  where $g$ is a measurable mapping defined on $(G,\mathcal{D})$ and takes its values in 
$\overline{\mathbb{R}}$.\\

\noindent To see this, let us use again the four steps method. In the first step, let us suppose is an indicator function of an element of $\mathcal{B}_Y$. So there exists
$C\in \mathcal{D}$, such that

\begin{equation*}
h=1_{Y^{-1}(C)}=1_{C}(Y)
\end{equation*}

\bigskip \noindent Clearly $g=1_C$ is a real-valued measurable mapping defined on $G$ such that $h=g(Y)$.\\

\noindent In a second step, let  $h$ be a of the form
\begin{equation*}
h=\sum_{i=1}^{p}\alpha _{i}1_{Y^{-}1}(B_{i})=\sum_{i=1}^{p}\alpha
_{i}1_{B_{i}}(Y)=\left( \sum_{i=1}^{p}\alpha _{i}1_{B_{i}}\right) (Y)=g(Y),
\end{equation*}

\noindent where $g$ is clearly $\mathcal{D}$-measurable. We easily move to non-negative $\mathcal{B}_{Y}$ -measurable function
by Point (6) of Proposition ... and the classical fact that any measurable and non-negative function is a non-decreasing limit
of a sequence of non-negative elementary function. The extension to an arbitrary quasi-integrable $\mathcal{B}$-measurable function is done by using the positive and negative parts.\\

\noindent In summary, whenever it exists, $\mathbb{E}(X/Y)$ is has the form 
\begin{equation*}
\mathbb{E}(X/Y)=g(Y),
\end{equation*}

\noindent $g$ is a measurable mapping defined on $(G,\mathcal{D})$ and takes its values in $\overline{\mathbb{R}}$.\\

\bigskip \noindent This function $g$ is also called the \textbf{regression function} of $X$ in $Y$ denoted as

\begin{equation*}
E(X/Y=y)=g(y).
\end{equation*}

\noindent It is very interesting to see a discrete version of that formula, which is very commonly used. Suppose that $Y$ takes a countable number of values denoted by $(y_j)_{j\in J}$,
$j \subset \mathbb{N}$. We recall that we have

\begin{equation*}
\mathbb{E}(X)=\mathbb{E}(\mathbb{E}(X/Y)),
\end{equation*}

\bigskip \noindent and next by using the regression function, we have

\begin{eqnarray}
\mathbb{E}(X)=\mathbb{E}(\mathbb{E}(X/Y))&=&\mathbb{E}g(Y) \notag \\
&=&\sum_{j\in J}g(y_{j})\mathbb{P}(Y=y_{j})=\sum_{j\in J}\mathbb{E}(X/Y=y_{j})\mathbb{P}(Y=y_{j}).  \label{ecd1}
\end{eqnarray}

\bigskip \noindent This gives 

\begin{equation*}
\mathbb{E}(X)=\sum_{j\in J}\mathbb{E}(X/Y=y_{j})\mathbb{P}(Y=y_{j}).
\end{equation*}

\bigskip \noindent If $X$ itself is  discrete and takes the values $(x_{i})_{i\in I}$, we have 

\begin{equation}
\mathbb{E}(X)=\sum_{j\in J}x_{i}\mathbb{P}(X=x_{i}/Y=y_{j})\mathbb{P}(Y=y_{j}).  \label{ecd2}
\end{equation}

\noindent Let us study the Jensen's inequality for the mathematical expectation.  We keep the same notations. For some details on convex function in our series, on may consult Exercise 6 and its solution in Doc 03.09 in Chapter 4 in \cite{ips-mestuto-ang}.\\

\section{Jensen's Inequality for Mathematical Expectation}  \label{proba_02_ce_sec_06}

\begin{theorem} \label{proba_02_jensenME} Let $X$ be random variable supported by an interval $I$ on which is defined a real-valued convex function $\phi$. Suppose that $X$ and $\phi(X)$ are integrable. Then for any $\sigma$-sub-algebra $\mathcal{B}$ of $\mathcal{A}$, we have

$$
\phi(\mathbb{E}(X/\mathcal{B}) \leq \mathbb{E}(\phi(X)/\mathcal{B}).
$$  
\end{theorem}

\noindent \textbf{Proof}. Let us follows \cite{chung} in the first proof therein. Let us proceed by step.\\

\noindent \textit{Step 1}. Let us suppose that $X$ takes a finite number of distinct values $(x_j), \ j\in J$, $J$ finite. Let us denote
$B_j=(X=x_j)$, $j \in J$ so that

$$
X=\sum_{j\in J} x_j 1_{B_j}, \ and \  \sum_{j\in J} 1_{B_j} =1_{\Omega}=1. \ \ (F1)
$$

\bigskip \noindent and hence

$$
\phi(X)=\sum_{j\in J} \phi(x_j) 1_{B_j}.
$$

\bigskip \noindent By the linearity of the mathematical expectation, we have

$$
\mathbb{E}(\phi(X)/\mathcal{B})= \sum_{j\in J} \phi(x_j) \mathbb{E}(1_{B_j}/\mathcal{B}). \ \ (F2)
$$ 

\bigskip \noindent But the real numbers $\mathbb{E}(1_{B_j}/\mathcal{B})$ add up to one since, because of Formula (F1), we get

$$
\sum_{j\in J} \mathbb{E}(1_{B_j}/\mathcal{B})=  \mathbb{E}\left(\sum_{j\in J} 1_{B_j}/\mathcal{B}\right)=\mathbb{E}( 1_{\Omega}/ \mathcal{B})=1.
$$

\bigskip \noindent Hence by the convexity of $\phi$, the right-hand member of Formula (F2) satisfies 

$$
\sum_{j\in J} \phi(x_j) \mathbb{E}(1_{B_j}/\mathcal{B}) \geq \phi\biggr( \biggr(\sum_{j\in J} x_j \mathbb{E}(1_{B_j}\biggr) /\mathcal{B}\biggr), 
$$

\bigskip \noindent and, surely, the right-hand member is  

$$
\phi\bigg( \mathbb{E}\biggr(\sum_{j\in J} x_j 1_{B_j}\biggr)/\mathcal{B}\biggr)=\phi(\mathbb{E}(X/\mathbb{B}). \ \ (F3)
$$

\bigskip \noindent By comparing the left-hand term of Formula (F2) and the right-hand term of Formula (F3), we get the desired result.\\

\noindent \textit{Step 2}. For a \textit{a.s.} finite general random variable, we already know from Measure Theory and Integration that $X$ is limit of a sequence elementary functions $(X_p)_{p\geq 1}$
with $|X_p|\leq |X|$ for all $p\geq 1$.\\

\noindent If $X$ is bounded \textit{a.s}, say $|X|\leq A<+\infty$ \textit{.a.s}, then by the continuity $\phi$, we have

$$
\max\biggr(|\phi(X)|, \sup_{p\geq 1} |\phi(X_p)|\biggr) \leq \|\phi\|_{[-A,A]} <+\infty.
$$ 

\bigskip \noindent By applying the result of Step 1, we have for all $p\geq 1$

$$
\phi(\mathbb{E}(X_p/\mathcal{B}) \leq \mathbb{E}(\phi(X_p)/\mathcal{B}).
$$  

\bigskip \noindent By applying the Dominated Convergence Theorem in both sides, we get the desired result.\\

\noindent \textit{Step 3}. Now suppose that $X$ is not bounded above. By Proposition 17.6 in \cite{choquet}, each point of $(a, \phi(a))$ of the Graph $\Gamma$ of the convex function $\phi$ has a supporting line, that is a straight line which passes through $(a, \phi(a))$ and is below $\Gamma$. A quick drawing may help to catch the meaning of this. For each $n\geq 1$, consider a supporting line at the point $(n, \phi(n))$ with equation $f_n(x)=A_nx+B_n$.\\

\noindent If $X$ is not bounded below, we consider, for each $n\geq 1$, a supporting line at the point $(-n, \phi(-n))$ with equation $g_(x)=C_nx+D_n$.\\

\noindent We may have $X$ bounded below and not bounded above, $X$ not bounded below and bounded above or $X$ neither bounded below and nor bounded above. In all these situations, we will have similar  way to handle the situation. Let us take the last case. We define

$$
\phi_n=g_n 1_{]-\infty,-n[} + \phi 1_{[-n,n]} + f_n 1_{]n,+\infty[}
$$

\bigskip \noindent We may check quickly that each $\phi_n$ is convex, $\phi_n \leq \phi$ and $\phi_n \uparrow \phi$ as $n\uparrow +\infty$. By denoting, for each $n\geq 1$,

$$
E_n=\|\phi\|_{[-n,n]}
$$

\bigskip \noindent and $a_n=|A_n|+|C_n|+|E_n|$ and $b_n=|B_n|+|D_n|$, we have for all $n\geq 1$, for all $x\in \mathbb{R}$

$$
|\phi_n(x)|\leq a_n|x|+b_n,
$$

\bigskip \noindent and next for all $n\geq 1$, for all $p\geq 1$

$$
|\phi_n(X_p)|\leq a_n|X_p|+b_n \leq a_n|X|+b_n, \ (F4)
$$

\bigskip \noindent Since for each $n\geq 1$, $\phi_n$ is convex, the result of Step 1 gives for all $p\geq 1$,

$$
\phi_n(\mathbb{E}(X_p/\mathcal{B}) \leq \mathbb{E}(\phi_n(X_p)/\mathcal{B}).
$$  

\bigskip \noindent By fixing $n\geq 1$, by letting $p\rightarrow +\infty$ and by applying the Dominated Convergence Theorem in both sides on the account of Formula (F4), we get 

$$
\phi_n(\mathbb{E}(X/\mathcal{B}) \leq \mathbb{E}(\phi_n(X)/\mathcal{B}).
$$  

\bigskip \noindent By letting $n\uparrow +\infty$, and by applying the Monotone convergence Theorem of the integrable functions in the right-hand member, we get the general conclusion. $\blacksquare$

\section[Orthogonal Projection in $L^2$]{The Mathematical Expectation as an Orthogonal Projection in $L^2$}  \label{proba_02_ce_sec_07}

\noindent Let us suppose that $X \in E=L^2(\Omega, \mathcal{A},\mathbb{P})$. For any $\sigma$-sub-algebra $\mathcal{B}$ of $\mathcal{A}$, let us consider $H=L^2(\Omega, \mathcal{B},\mathbb{P})$ the square integrable and real-valued $\mathcal{B}$-measurable functions. At least $1=1_{\Omega}$ and $0=1_{\emptyset}$ are elements of $H$.\\

\noindent We already know that $L^2(\Omega, \mathcal{A},\mathbb{P})$ is a Hilbert space endowed with inner product

$$
L^2(\Omega, \mathcal{A},\mathbb{P})^2 \ni (X,Y) \mapsto \langle X,Y \rangle=\mathbb{E}(XY).
$$

\bigskip \noindent  We have the following projection theorem in Hilbert spaces (See for example Theorem 6.26 in \cite{chidumefa}, page 109).\\

\begin{proposition} \label{probab_02_project} Suppose that $E$ is a Hilbert space and $H$ a closed sub-linear space. Fix $x \in E$. We have the following facts.\\

\noindent (1) There exists a unique element $p_H(x) \in H$ such that

$$
d(x,H)=\inf \{\|x-h\|, \ h \in H\}=\|x=p_H(x)\|.
$$

\bigskip \noindent (2) $p_H(x)$ is also the unique element of $H$ such that $x-p_H(x)$ is orthogonal all elements of $H$. 
\end{proposition}

\noindent We are going to apply it in order to characterize $\mathbb{E}(X/\mathcal{B})$. We have

\begin{theorem} \label{probab_02_MEL2} For any $\sigma$-sub-algebra $\mathcal{B}$ of $\mathcal{A}$, $H(\mathcal{B})=L^2(\Omega, \mathcal{B},\mathbb{P})$ is a closed linear space, and for any
$X \in E=L^2(\Omega, \mathcal{A},\mathbb{P})$,

$$
\mathbb{E}(X/\mathcal{B})=p_{H(\mathcal{B})}(X). \ a.s.
$$
\end{theorem}

\bigskip \noindent  \textbf{Proof}. Let us begin to show that $H(\mathcal{B})$ is closed. Let $Z$ a limit of a sequence of elements of $H(\mathcal{B})$. Since $E$ is a Hilbert space, we know that
$Z$ is still in $E$, thus is square integrable. By Theorem \ref{proba_02_con_sec_th_01}  in Chapter \ref{proba_02_conv}, page \pageref{proba_02_con_sec_th_01}, the concerned sequence converges to $Z$ in probability and next, by the relation between weak and strong limits seen in the same chapter, a sub-sequence of the sequence converges \textit{a.s.} to $Z$. Finally $Z$ being an \textit{a.s.} limit of
a sequence $\mathcal{B}$-measurable functions is $\mathcal{B}$-measurable. In total $Z \in H(\mathcal{B})$. thus $H(\mathcal{B})$ is closed in $E$.\\

\noindent Now for any $X \in E$, Point (2) of Proposition \ref{probab_02_project} characterizes $Z=p_{H(\mathcal{B})}(X)$. Thus for any $B \in \mathcal{B}$, $h=1_B \in H(\mathcal{B})$, $\langle X-Z, h \rangle=0$, that is $\langle X, Z \rangle=\langle X, h \rangle$

$$
\int_B X \ d\mathbb{P}= \int_B Z \ d\mathbb{P}.
$$

\bigskip \noindent  We conclude that $Z=\mathbb{E}(X/\mathcal{B})$. $\square$

\section{Useful Techniques}  \label{proba_02_ce_sec_08}

\noindent In a great number of situations, we need to compute the mathematical expectation of a real-valued function $h(X)$ of $X$ and we have to use a conditioning based on another random variables $Y$
where of course $X$ and $Y$ are defined on the same probability space $(\Omega, \mathcal{A}, \mathbb{P})$, even if they may have their values in different measure spaces  $(E_1, \mathcal{G}_{E_1})$ and
$(E_2, \mathcal{G}_{E_2})$. Suppose $h : (E_1, \mathcal{G}_{E_1}) \rightarrow \mathbb{R}$ is measurable and that $h$ is non-negative or $h(X)$ is integrable. We already know that there exists a measurable
function $g : (E_2, \mathcal{G}_{E_2}) \rightarrow \mathbb{R}$ such that

$$
\mathbb{E}(h(X)/Y)=g(Y) \ \ (CG01)
$$

\bigskip \noindent and by this, we have

$$
\mathbb{E}(h(X))=\mathbb{E}(g(Y)). \ \ (CG02)
$$ 

\noindent Now suppose that $Y$ is continuous with respect to a measure $\nu$ that is given on the measure space $(E_2, \mathcal{G}_{E_2})$. By reminding that
$g(y)=\mathbb{E}(h(X)/(Y=y))$ for $y \in E_2$, we have from formula (CG02)

\begin{eqnarray*}
\mathbb{E}(h(X))&=&\mathbb{E}(g(Y))\\
&=&\int_{E_2} g(y) \ d\mu(y)\\
&=&\int_{E_2} \mathbb{E}(h(X)/(Y=y)) \ d\nu(y)\\
\end{eqnarray*}

\bigskip \noindent which is, in short, \label{cg03}

\begin{eqnarray*}
\mathbb{E}(h(X))=\int_{E_2} g(y) \ d\nu(y). \ \ (CG03)
\end{eqnarray*}

\bigskip \noindent This result takes the following particular forms.\\

\noindent \textbf{(I) Conditioning by a discrete random variable}.\\

\noindent Suppose that $Y$ is discrete, that is the values set of $Y$ is countable and is written as $D=\{y_j, \in J\}$, $J \subset \mathbb{N}$. Hence the probability law of $Y$ is continuous with respect to the counting measure $\nu$ with support  $D$. Thus by the Discrete Integration Formula (DIF1) (see page \pageref{dif1}) applied to Formula (CG03), we get \label{cd}

$$
\mathbb{E}(h(X))=\sum_{j\in J} \mathbb{E}(h(X)/(Y=y_j)) \mathbb{P}(X=j_j). \ \ (CD)
$$ 

\bigskip \noindent \textbf{(II) Conditioning by an absolutely continuous real random vector}.\\

\noindent Suppose that $E_2=\mathbb{R}^r$, $r\geq 1$. If $Y$ is continuous with respect to the Lebesgue measure $\lambda_r$, we get the formula \label{cc01}

$$
\mathbb{E}(h(X))=\int_{\mathbb{R}^r} \mathbb{E}(h(X)/(Y=y)) \ f_Y(y) \ d\lambda_r(y). \ \ (CC01)
$$ 

\bigskip \noindent \textbf{(III) Conditional probability density function}.\\

\noindent On top of the assumptions in Part  (II) above, let us suppose also that $E_2=\mathbb{R}^s$, $s\geq 1$, with $d=r+s$. Let us suppose that $Z=(X^t, Y^t)^t$ has a \textit{pdf} 
$f_Z \equiv f_{(X,Y)}$ with respect to the Lebesgue measure on $\mathbb{R}^d$. Thus the marginal \textsl{pdf}'s of $X$ and $Y$ are defined by

$$
f_X(x)=\int_{\mathbb{R}^r} f_{(X,Y)}(x,y) \ d\lambda_r(y), \ \ x \in \mathbb{R}^s
$$

\bigskip \noindent and 

$$
f_Y(y)=\int_{\mathbb{R}^s} f_{(X,Y)}(x,y) \ d\lambda_s(x), \ \ y \in \mathbb{R}^r.
$$

\bigskip \noindent Let us define

$$
f_{X|Y=y}(x)=\frac{f_{(X,Y)}(x,y)}{f_Y(y)}, \ \ f_Y(y)>0, 
$$

\bigskip \noindent as the conditional \textit{pdf} of $X$ given $Y=y$. The justification of such a definition relies in the fact that replacing $\mathbb{E}(h(X)/(Y=y))$ 

$$
\int_{\mathbb{R}^s} h(x) f_{X|Y=y}(x) \ d\lambda_s(x)
$$

\bigskip \noindent and using Tonelli's Theorem (when $h$ in non-negative) or Fubini's Theorem (when $h(X)$ is integrable) leads to

\begin{eqnarray*}
&&\int_{\mathbb{R}^r} \biggr(\int_{\mathbb{R}^s} h(x) f_{X|Y=y}(x) \ d\lambda_s(x) \biggr) \ f_Y(y) \ d\lambda_t(y)\\
&=&\int_{\mathbb{R}^r} \biggr(\int_{\mathbb{R}^s} h(x) \frac{f_{(X,Y)}(x,y)}{f_Y(y)} \ d\lambda_s(x) \biggr) \ f_Y(y) \ d\lambda_t(y)\\
&=&\int_{\mathbb{R}^s} h(x) \biggr(\int_{\mathbb{R}^s}  \frac{f_{(X,Y)}(x,y)}{f_Y(y)} \ d\lambda_t(y) \biggr)  \ d\lambda_s(x)\\
&=&\int_{\mathbb{R}^s} h(x) f_X(x) d\lambda_s(x)\\
&=&\mathbb{E}(h(X)).
\end{eqnarray*}

\bigskip \noindent This leads to the frequent use the following formula \label{cc02}

$$
\mathbb{E}(h(X))=\int_{\mathbb{R}^r} \biggr(\int_{\mathbb{R}^s} h(x) f_{X|Y=y}(x) \ d\lambda_s(x) \biggr) \ f_Y(y) \ d\lambda_r(y). \ \ (CC02)
$$

%% file: proba_02_09_ang.tex
\chapter{Probability Laws of family of Random Variables} \label{thfondamentalKolm} 
\Large
\section{Introduction}

\noindent We have already studied the finite product measure in Chapter 8 in \cite{ips-mestuto-ang} for $\sigma$-finite measures defined on arbitrary measure spaces. In this chapter we give the Theorem of Kolmogorov which establishes the arbitrary product of probability measures, but in special measure spaces.\\

\noindent This theorem of Kolmogorov is the foundation of the modern theory of probability. There is nothing above it, in term of probability laws. On this basis, the modern theory of random analysis, which extends Real Analysis (paths smoothness, differentiability, integration according to different types, etc) is built on.\\

\noindent We recommend the reader, especially the beginner, to read it as many as possible and to often and repeatedly come back to it in order to see its deepness and to understand its consequences.\\

\noindent Among special spaces on which the construction is made, we count Polish spaces. A Polish space is a complete and separable metric space $(E,d)$ like $(\mathbb{R}^s, \rho)$, $s\geq 1$, where $\rho$ is one of its three classical metrics. An interesting remark is that the finite product of Polish spaces is a Polish space. The finite Borel product $\sigma$ of Polish spaces is generated by the product of open balls.\\

\noindent Here, the level of abstraction is moderately high. Once again, we recommend the beginner to go slow and to give himself the needed time to understand the definitions and the notation. This chapter may be considered as a continuation of Chapter 8 in \cite{ips-mestuto-ang}.\\

\noindent We already encountered this Kolmogorov construction in finite dimensions in Chapter \ref{proba_02_rv} in pages \pageref{kolmconst_01} and \pageref{kolmconst_02}. The results in this chapter will be the most general extension of this kind of result.\\

\noindent In the first section, we state and prove the existence of the product probability measure. Next, we will see how to  state a number of particular forms involving Lebesgue-Stieljes measures.\\

\bigskip \section{Arbitrary Product Probability Measure}

\noindent Let $(E_{t},\mathcal{B}_{t},\mathbb{P}_{t})$, $t\in T\neq \emptyset$, be a family of probability spaces. We define the product space by
\begin{equation*}
E=\prod\limits_{t\in T}E_{t}.
\end{equation*}

\bigskip \noindent If $T$ is finite, even countable, we may use the classical notation : $T=\{t_j, j\geq 0\}$. It make senses to speak about the first factor $E_{t_1}$, the second $E_{t_2}$, etc. The elements of

\begin{equation*}
E=\prod_{j \geq 0} E_{t_j},
\end{equation*}

\bigskip \noindent may be denoted by $x=(x_{t_1}, x_{t_2}, ...)$  as an ordered set.\\

\noindent But, the index set $T$ may arbitrary and uncountable. For example $T$ be may a set of functions. If the functions are real-valued, $T$ is uncountable and has a partial order. Sometimes we may not have an order at all. So the general appropriate way to study $E$ seems to consider $E$ as a set of functions. Thus, an element $x$ of $E$, written as 

\begin{equation*}
x=(x_{t})_{t\in T}=(x(t),t\in T),
\end{equation*}

\bigskip \noindent is perceived as a function $x$ which corresponds to each $t \in T$ a value $x(t)=x_t \in E_t$.\\

\noindent Let us begin to introduce the projections and give relevant notation. We denote by $\mathcal{P}_{f}$ the class of finite and non-empty subsets of $T$. Given an element  $S=\{s_{1},...,s_{k}\}$ of $\mathcal{P}_{f}$, we may write in any order of the subscripts. This leads to the class of ordered and non-empty subsets denoted by $\mathcal{P}_{of}$. Elements of $\mathcal{P}_{of}$ are written as $k$-tuples $S=(s_{1},...,s_{k})$, $k\geq 1$. For any $S=(s_{1},...,s_{k}) \in \mathcal{P}_{of}$, we have the finite product space

$$
E_{S}=\prod\limits_{j=1}^{k}E_{s_{j}}
$$

\bigskip \noindent which is endowed with the finite product $\sigma$-algebra 
\begin{equation*}
\mathcal{B}_{S}=\bigotimes\limits_{j=1}^{k}\mathcal{B}_{s_{j}}.
\end{equation*}

\bigskip \noindent The projection of this space $E_{S}$ is defined by

\begin{equation}
\begin{array}{cccc}
\Pi _{S}: & (E,\mathcal{B}) & \longrightarrow & (E_{S},\mathcal{B}_{S}) \\ 
& x=(x_{t})_{t\in T} & \longmapsto & \Pi _{S}(x)=(x_{s_{1}},...,x_{s_{k}})
\end{array}
.  \label{intro1}
\end{equation}

\bigskip \noindent We name $S$ as the index support of the projection.\\

\noindent Our first objective is to define a $\sigma $-algebra on $E$, which renders measurable all the projections on finite sub-products spaces.

\subsection{The Product $\sigma$-algebra on the product space}

Let us begin by the definition 

\begin{definition} \label{probab_02_sigmaAlgebraEngProj} The product $\sigma$-algebra on $E$, denoted by $\mathcal{B}$, is the smallest $\sigma$-algebra rendering measurable all the projections of finite index support.
\end{definition}

\bigskip \noindent We already know that such a $\sigma$-algebra exists. Compared to the finite product $\sigma$-algebra, there is nothing new yet.\\

\noindent In the sequel, we have to change the order of elements of $V \in \mathcal{P}_{of} $. So the following recall may be useful. Indeed, by permuting the elements of $V=(v_{1},...,v_{k}) \in \mathcal{P}_{of}$, $k\geq 1$, by means of a permutation $s$ of $\{1,2,...,n\}$, the correspondence

\begin{equation}
\begin{array}{cccc}
& (E,\mathcal{B}) & \longrightarrow & (E_{s(S)},\mathcal{B}_{s(S)}) \\ 
& (x_{v_{1}},...,x_{v_{k}}) & \longmapsto & (x_{s(v_{1})},...,x_{s(v_{k})})%
\end{array}%
.  \label{intro2}
\end{equation}

\bigskip \noindent is a one-to-one mapping. Also, in parallel of the notation of $E_V$, we may and do adopt the following notation

$$
x_V=(x_{s(v_{1})},...,x_{s(v_{k})}).
$$

\bigskip  \noindent As well, the space $\mathcal{S}_V$ denote the class of measurable rectangles in $E_V$.\\

\bigskip \noindent Now, as in the finite product case, we have to see how to generate $\mathcal{B}$ by what should correspond to the class of measurable rectangles. Here, we use the phrasing in \cite{loeve} of measurable cylinders. Let $S=(s_{1},...,s_{k}) \in \mathcal{P}_{of}$, $k\geq 1$. A finite measurable rectangle in $E_S$ is generally denoted by
 
\begin{equation}
A_{S}=\prod\limits_{j=1}^{k}A_{s_{j}},\text{ }(A_{s_{j}}\in \mathcal{B}%
_{s_{j}},1\leq j\leq k)  \label{intro3}
\end{equation}

\bigskip \noindent It is clear that $\Pi _{S}^{-1}(A_{S})$ is the set of all  $x=(x_{t})_{t\in T}$ such that  
\begin{equation*}
x_{s_{j}}\in A_{s_{j}},\text{ }1\leq j\leq k.
\end{equation*}

\bigskip \noindent We write the above fact as

$$
\Pi _{S}^{-1}(A_{S})=A_S \times \prod_{t\notin S} E_t. \ \ (FP)
$$

\bigskip

\begin{definition} The class $\mathcal{S}$ of measurable cylinders of $E$ is the class of subsets of $E$ which are of the form $\Pi _{S}^{-1}(A_{S})$, $S\in \mathcal{P}_{of}$.\\

\noindent In other words, a measurable cylinder of $E$ is a product of measurable subsets $A_t \in \mathcal{B}_t$ of the form

$$
\prod_{t \in T} A_t, \ \ (SP01)
$$

\bigskip \noindent such that at most a finite number of the $A_t$, $t\in T$, are non-empty. If $V=\{v_1, ..., v_k\} \subset T$ is such that $A_t=\emptyset$ for $t\notin V$, then the product

$$
A_V=\prod_{v \in V} A_v \ \ (SP02)
$$

\bigskip \noindent is called a finite support of the cylinder and the cylinder is written as

$$
c(A_V)= A_V \times \prod_{t\notin V} E_t. \ \ (SP03)
$$
\end{definition}

\noindent \textbf{Remarks}. The following remarks are important.\\

\noindent (1) In the definition of the support of the cylinder in Formula (SP02), the order of $V$ is not relevant in the writing of the cylinder $c(A_V)$, but it really counts in the writing of the support $A_V$.\\

\noindent (2) A support is not unique. For example if one of $E_v$, $v \in V$, is equal to $E_v$, may may drop it from the support. As well, we may add to $V$ any other $w \notin V$ such that $A_w=E_w$ : we may drop full spaces from the support and add full spaces to it.\\

\noindent (3) Formula (SP03) means that $x=(x_t)_{t \in T}$ is in the cylinder only if $x_v \in E_v$, $v\in V$, and we do not care about where are the $x_t$, $t\notin V$. The only knowledge about them is that they remain in their full space $E_t$, $t\notin V$.\\

\noindent (4) The notation $c(A_V)$ introduced in Formula (SP03) stands for cylinder of support $A_V$ where $V$ is non-empty set of $T$.\\

\noindent (5) For the sake of shorter notation, we may write the formula in (SP03) in the form

$$
c(A_V)= A_V \times E^{\prime}_S \ \ where \ \ E^{\prime}_S=\prod_{t\in V} E_t
$$ 

\bigskip  \noindent The coming concept of coherence, which is so important to the theory of Kolmogorov, depends on the understanding of the remarks above and the next remark.\\

\noindent \textbf{(6) Common index support of two cylinders}. Consider two cylinders

\begin{equation}
c(A_V)=A_{V}\times \prod\limits_{t\notin V}E_{t}\text{ \ and \ \ } c(B_{W})=B_{W} \times \prod\limits_{t\notin W}E_{t}.  \label{intro6a}
\end{equation}

\bigskip \noindent of respective supports $V=(v_{1},...,v_{q}) \in \mathcal{P}_{of}$ and $W=(w_{1},...,w_{p}) \in \in \mathcal{P}_{of}$, $p\geq 1$, $q \geq 1$.\\

\noindent We want to find a common support for both $c(A_V)$ and $c(A_W)$. We proceed as follows. Let us form an ordered set $U$ by selecting first all the elements of $V$ in the ascendent order of the subscripts. Next we complete by adding the elements of $W$ which were not already in $U$, still in the ascendent order of the subscripts. At the arrival, the elements of $W$, corresponding to the common elements of $V$ and $W$ if they exists, may not be present if $U$ in the original order of the their subscripts in $W$. Rather, they are present in $U$ in the subscripts order of some permutation $s(W)$ of $W$.\\

\noindent \textbf{Example}. Suppose that $V=(v_1, v_2, v_3, v_4, v_5)(1,4,7,2, 5)$ and $W=(w_1,w_2,w_3,w_4)=(5, 2, 10, 8)$. We have
$$
U=(1, 4, 7, 2, 5, 10, 8)
$$ 

\noindent * So the elements of $W$ are given in order in $(w_2,w_1,w_3,w_4)$ which is $s(W)$ with $s(1)=2$, $s(2)=2$, $s(3)=3$ and $s(4)=4$.\\

\noindent But we already saw that the order of the subscripts of $V$, or $W$ or $U$ does not alter the cylinders $c(A_V)$, $c(B_W)$, $c(A_U)$ or $c(B_U)$. We have

\begin{equation}
c(A_V)=A_{U}\times E^{\prime}_U \text{ \ and \ \ } c(B_{W})=B_{U}\times E^{\prime}_U.  \label{commonB}
\end{equation}

\bigskip  \noindent Actually, we formed $A_{U}$ (resp. $A_{U}$) by adding full spaces $E_t$ to the support $A_V$  (resp. $B_W$) for $t \in U \setminus V$) (resp. for $t \in U \setminus W$) .\\

\noindent We say that we have written $c(A_V)$ and $c(B_W)$ with a common index support $U$. This consideration will be back soon.\\

\noindent (5) In the definition of a cylinder, the finite support, say $A_V$, is a measurable rectangle. But in general, $A_V$ may be a measurable subset of $E_V$, that is $A_V \in \mathcal{B}_V$ and we still have

$$
\Pi^{-1}(A_V)=A_V \times E^{\prime}_V,
$$

\bigskip \noindent which is to be interpreted as

$$
x \in A_V \times E^{\prime}_V \Leftrightarrow x_V \in A_V.
$$

\bigskip \noindent Now, we are ready to go further and to give important properties of $\mathcal{S}$.\\

\begin{proposition} \label{probab_02_semiAlgebra} 
$\mathcal{S}$ is a semi-algebra.

\end{proposition}

\noindent \textbf{Proof}. (i) Let us see that $E \in \mathcal{S}$. If we need to prove it, we consider a point $t_{0}\in T$, put $A_{t_{0}}=E_{t_{0}}$ and get that
\begin{equation*}
E=A_{t_{0}}\times \prod\limits_{t\neq t_{0}}E_{t}=c\left(A_{\{t_0\}}\right).
\end{equation*}

\bigskip \noindent (ii) Next, by the definition of a cylinder of support $S \in \mathcal{P}_{of}$, checking that $x$ belongs to $c(A_S)$ or not depends only of $x_S \in A_S$ or not. So we have

\begin{equation}
c(A_{S})^{c}=A_{S}^{c}\times \prod\limits_{t\notin S}E_{t},\label{intro5}
\end{equation}

\noindent Next, let us check that the complement of any element of $\mathcal{S}$ is a finite sum of elements of $S$. We already knew that the class of measurable rectangles in $E_V$ is a semi-algebra, so that $A_{S}^{c}$ is a finite sum of elements of measurable rectangles $A^{(j)}_S$, $1\leq j \leq r$, $r\geq 1$, of $E_S$. And it becomes obvious that

\begin{eqnarray*}
c(A_{S})^{c}&=&A_{S}^{c}\times \prod\limits_{t\notin S}E_{t}\\
&=&\biggr( \sum_{1 \leq j \leq r } A^{(j)}_S \biggr) \times \prod\limits_{t\notin S}E_{t}\\
&=& \sum_{1 \leq j \leq r } \biggr( A^{(j)}_S \times \prod\limits_{t\notin S}E_{t}\biggr) \\
&=& \sum_{1 \leq j \leq r } c(A^{(j)}_S).
\end{eqnarray*}

\bigskip \noindent Our checking is successful.\\

\noindent (iii) Finally, let us check that $\mathcal{S}$ is stable under finite intersection. To do so, let us consider two cylinders

\begin{equation}
c(A_V)=A_{V}\times \prod\limits_{t\notin V}E_{t}\text{ \ and \ \ } c(B_{W})=B_{W} \times \prod\limits_{t\notin W}E_{t},
\end{equation}

\noindent and next their expressions using a common index support as explained earlier, we have

\begin{equation}
c(A_V)=A_{U}\times E^{\prime}_U \text{ \ and \ \ } c(B_{W})=B_{U}\times E^{\prime}_U.
\end{equation}

\bigskip \noindent It becomes clear that we have

\begin{equation}
c(A_V) \cap c(B_{W})=c(A_V)=\biggr(A_{U}\cap B_{U}\biggr) \times E^{\prime}_U,
\end{equation}

\bigskip \noindent which is element of $\mathcal{S}$ since

$$
A_{U}\cap B_{U} =\prod_{t \in U} A_t \cap B_t.
$$

\noindent * We also have

\begin{theorem} \label{probab_02_sigmaAlgebraEngSA} The $\sigma$-algebra on $E$ generated by the projections, denoted $\mathcal{B}$ in Definition \ref{probab_02_sigmaAlgebraEngProj}, is also generated by the class of cylinders of finite support $\mathcal{S}$, called the product $\sigma$-algebra and denoted as

\begin{equation*}
\mathcal{B}=\bigotimes\limits_{t\in T}\mathcal{B}_{t}.
\end{equation*}

\end{theorem}

\noindent \textbf{Proof}. Let us denote by $\mathcal{B}$ the $\sigma$-algebra on $E$ generated by the projections with finite support and by $\mathcal{B}_0$ the one generated by $\mathcal{S}$.\\

\noindent (1) Let us prove that $\mathcal{B}\subset \mathcal{B}_0$. Let us fix $V \in \mathcal{P}_{of}$. For any measurable rectangle $A_V$ which in $E_V$, that is 
$A_V\in \sigma(\mathcal{S}_V)$, we already now, since Formula (FP), that

$$
\pi^{-1}(A_V)=A_V \times E^{\prime}_V 
$$ 

\bigskip \noindent and next

$$
\pi^{-1}_V(A_V)=A_V \times E^{\prime}_V  \in \mathcal{S} \subset \mathcal{B}_0.
$$ 

\bigskip \noindent So each projection $\pi^{-1}_V$ of finite support is $\mathcal{B}_0$-measurable. We conclude $\mathcal{B}\subset \mathcal{B}_0$ by the definition of $\mathcal{B}$.\\

\noindent 2) Let us prove that $\mathcal{B}_0\subset \mathcal{B}$. It is enough to prove that  $\mathcal{S}\subset \mathcal{B}$. But any element $A$ of $\mathcal(S)$ can be written as

$$
A=\prod_{t \in V} A_t \times \prod_{t \notin V} E_t =:A_V \times E^{\prime}_V, \ A_t \in \mathcal{B}_t
$$

\bigskip  \noindent which is

$$
A=\Pi_{V}^{-1}(A_{V}),
$$

\bigskip \noindent and then to $\mathcal{B}$, since $A_{V} \in \mathcal{B}_V$ and $\Pi_{V}$ is $\mathcal{B}$-measurable. $\square$\\

\section[Measurability for a family of Random Variables]{Stochastic Process, Measurability for a family of Random Variables} \label{probab_02_sec_procSto}

\noindent \textbf{I - General case : family of random variable}.\\

\noindent Now we have the product space 
$$
E=\prod_{t \in T} E_t,
$$

\noindent endowed with the product $\sigma$-algebra

\begin{equation*}
\mathcal{B}=\bigotimes\limits_{t\in T}\mathcal{B}_{t}.
\end{equation*}

\noindent We may study the measurability of mappings $X \ :  \  (\Omega, \mathcal{A}) \rightarrow (E, \mathcal{B})$. According to the notation above, we denote 

$$
\forall \omega \in \Omega, \ X(\omega)=(X_{t}(\omega ))_{t\in T}.
$$

\bigskip \noindent For all $t \in T$, the mapping $\omega \longmapsto X_{t}(\omega )$ taking its values in $E_{t}$ is called the $t$-th component or margin.We immediately have that for each $t \in T$,
\begin{equation*}
X_{t}=\Pi _{t}\circ X.
\end{equation*}

\bigskip \noindent It become clear that if $X$ is measurable, thus each margin $X_{t}$, $t\in T$, is also  measurable. Actually, this is a characterization of the measurability of such mappings.

\begin{proposition} \label{probab_02_margins_01} A mapping  $ X \: \ (\Omega ,\mathcal{A}) \longrightarrow (E,\mathcal{E})$ is measurable if and only if each margin $X_{t}$, $t \in T$, is measurable. Indeed, we have
\end{proposition}

\bigskip \noindent \textbf{Proof}. We only need to prove the implication that if all the margins $X_{t}$, $t \in T$, are measurable, then $X$ is. Suppose that all the margins $X_{t}$, $t \in T$, are measurable.
It will be enough to show that $X^{-1}(cl(A_V))$ is measurable whenever $cl(A_V)\in \mathcal{S}$. By using the notation above, we have

\begin{eqnarray*}
\omega \in X^{-1}(c(A_V)) \ \ &\leftrightarrow& \ \  X(\omega) \in c(A_V)\\
 &\leftrightarrow& \ \  (\forall v \in V, X_v(\omega ) \in A_{v})\\
&\leftrightarrow& \ \  (\forall v \in V, \omega \in X_v^{-1}(A_{v}))\\
&\leftrightarrow& \ \  \omega \in \bigcap_{v \in A} X_v^{-1}(A_{v}),\\
\end{eqnarray*}

\noindent which, by the measurability of the $X_v$'s, gives

$$
X^{-1}(c(A_V))=\bigcap_{v \in A} X_v^{-1}(A_{v}) \in \mathcal{A}. \ \square
$$

\bigskip \noindent \noindent \textbf{II - Stochastic Processes.}\\

\noindent Let consider the special case where  all the $E_t$ are equal to one space $E_0$ on which is defined a $\sigma$-algebra $\mathcal{B}_0$. The product space is denoted by

$$
E=E_0^{T}
$$

\bigskip \noindent and is the class of all mappings defined from $T$ to $E_0$. As in the general context,  elements of $E$ are denoted $x=(x_t)_{t \in T}$, where for all $t\in T$, $x_t \in E_0$. The product
$\sigma$-algebra is denoted by

$$
\mathcal{B}=\mathcal{B}_0^{\otimes T}.
$$

\bigskip \noindent We have the general terminology :\\

\noindent (1) A measurable application $X \ : \ (\Omega, \mathcal{A}) \rightarrow \left(E_0^{T},\mathcal{B}_0^{\otimes T}\right)$ is called a \textbf{stochastic process}.\\

\noindent (2) $E_0$ is called the states space of the stochastic process.\\

\noindent (3) $T$ is called the time space in a broad sense.\\

\noindent (4) If $T=\{1\}$ is a singleton, the stochastic process is called a simple random variable.\\

\noindent (5) If $T=\{1, ..., k\}$ is finite with $2\leq k \in \mathbb{N}$, the stochastic process is  called a random vector.\\

\noindent (6) If $T=\mathbb{N}$, the stochastic process is  a sequence of random variables.\\

\noindent (7) If $T=\mathbb{Z}$, the stochastic process is called a time series.\\

\noindent (8) If $T=\mathbb{R}_+$, the terminology of time space is meant in the real-life case.\\

\noindent (9) If $T$ in endowed with a partial order, we generally speak of a random field or random net.\\

\noindent (10) For any $\omega \in \Omega$, the mapping

$$
T \ni t \mapsto X_t(\omega).
$$ 

\bigskip \noindent is called a path of the stochastic process on $E_0$.\\

\section{Probability Laws of Families of random Variables}

\noindent \textbf{(I) The concept of Coherence}.\\

\noindent Consider a probability measure $\mathbb{P}$ on the product measure space $(E,\mathcal{B})$. The image measure on a sub-product $(E_{S},\mathcal{B}_{S})$ by the projection $\Pi _{S}$, where $S=(s_{1},...,s_{k}) \in \mathcal{P}_{of}$ is
\begin{equation*}
\mathbb{P}_{S}=\mathbb{P}\Pi _{S}^{-1},
\end{equation*}

\bigskip \noindent We recall that for any $B_{S} \in \mathcal{B}_S$, we have 

\begin{equation*}
\mathbb{P}_{S}(B_{S})=\mathbb{P}(\Pi _{S}^{-1}(B_{S})).
\end{equation*}

\bigskip \noindent We get the family of probability measures 
$$
\{\mathbb{P}_{S},S\in \mathcal{P}_{of}(T)\},
$$ 

\bigskip \noindent which we called \textit{the family of marginal probability measures with finite index support}. By a language abuse, we also use the phrase of  family of \textit{finite-dimensional} marginal probability measures of $\mathcal{P}$.\\

\noindent We are going to discover some important relations between the \textit{finite-dimensional} marginal probability measures. But we should also keep in mind that, for any $S \in \mathcal{P}_{of}(T)$, a probability measure $\mathbb{P}_{S}$ on $(E_S, \mathcal{B}_S)$ is characterized by its values on $\mathcal{S}_S$, which the class of measurable rectangles on $E_S$.\\

\noindent First, let us consider $(S_1,S_2) \in \mathcal{P}_{of}(T)^2$, such that one of them is a permutation of the other, that is $S_{1}=s(S_{2})$, where $S_{1}=(s_1,...,s_k)$, and $s$ is a permutation $\{1,2,...,k\}$. Consider any $A_{S_{1}} \in \mathcal{S}_{S_1}$. We have

\begin{equation*}
s(A_{S_{1}})=s(A_{s_{1}}\times ...\times A_{s_{k}})=A_{s(s_{1})}\times
...\times A_{s(s_{k})}=A_{S_{2}}.
\end{equation*}

\bigskip \noindent Furthermore, the projection on $E_{S_{2}}$ is the composition of the projection on $E_{S_{1}}$ and the permutation of that projection by $s$, which gives

\begin{equation*}
\Pi _{S_{2}}=s\circ \Pi _{S_{1}},
\end{equation*}

\noindent from which, by the characterization of a finite product probability by its values on the measurable rectangles, we have
\begin{equation*}
\Pi _{S_{2}}^{-1}=\Pi _{S_{1}}^{-1}\circ s^{-1},
\end{equation*}

\bigskip \noindent and similarly,

\begin{equation*}
\Pi _{S_{1}}^{-1}=\Pi _{S_{2}}^{-1}\circ s.
\end{equation*}

\noindent Hence, for any $B_{S_{1}}\in \mathcal{B}_{S_{1}}$, we have
\begin{equation*}
\mathbb{P}\Pi _{S_{1}}^{-1}(B_{S_{1}})=\mathbb{P}\Pi _{S_{2}}^{-1}\circ
s(B_{S_{1}}),
\end{equation*}

\noindent which leads to

\begin{equation*}
\mathbb{P}_{S_{1}}(B_{S_{1}})=\mathbb{P}_{S_{2}}(s(B_{S_{1}})),
\end{equation*}%

\noindent and 
\begin{equation*}
\mathbb{P}_{S_{1}}(\cdot )=\mathbb{P}_{s(S_{1})}(s(\cdot )).
\end{equation*}

\noindent * We already reached a first coherence (or consistency) condition. Let us discover a second one. Let $U=(u_{1},..,u_{r})\subset S=(s_{1},...,s_{k})$, where the inclusion holds with the preservation of the ascendent order of the subscripts. Then the projection on $E_{U}$ is obtained by the projection on $E$ on $E_{S}$ first, and next by the projection of $E_{S}$ on $E_{U}$ denoted $\Pi _{S,U}$. Accordingly to the notation above, we have for any $B_U \subset E_U$,

$$
\Pi _{S,U}^{-1}(B_U)=B_U \times E^{\prime}_{S\setminus U},
$$

\noindent which is interpreted as

$$
x_S=(x_U, \ x_{{S\setminus U}}) \in \Pi _{S,U}^{-1}(B_U) \Leftrightarrow x_U \in B_U.
$$

\noindent Going back to the considerations which were made above about the projection on $S$, we have 
 
\begin{equation*}
\Pi _{U}=\Pi _{S,U}\circ \Pi _{S},
\end{equation*}

\noindent and next, 
\begin{equation*}
\mathbb{P}\Pi _{U}^{-1}=\mathbb{P}\Pi _{S}^{-1}\Pi _{S,U}^{-1},
\end{equation*}

\bigskip \noindent and finally,

\begin{equation*}
\mathbb{P}_{U}=\mathbb{P}_{S}\Pi _{S,U}^{-1}.
\end{equation*}

\bigskip \noindent We get a second relation between the marginal probability measures. Based on the previous developments, we may define

\begin{definition} A family of finite-dimensional probability measures $\{\mathbb{P}_{S},S\in \mathcal{P}_{of}(T)\}$ is said to be coherent if and only if we have the following two conditions, called coherence coherent or consistency conditions :\\

\noindent (CH1a) \label{CH1a} For any ordered and finite subsets $U$ and  $S$ of $T$ such that $U$ is subset of $S$ with the preservation of the ascendent ordering of the subscripts of $U$ in $S$, we have

\begin{equation*}
\mathbb{P}_{U}=\mathbb{P}_{S}\Pi _{S,U}^{-1}.
\end{equation*}

\noindent (CH2) \label{CH2}  For any ordered and finite subset $S$ of $T$ and for any permutation $s$ of $E_S$, for any $B_S \in \mathcal{B}_S$,

\begin{equation*}
\mathbb{P}_{S}(B_{S})=\mathbb{P}_{s(S)}(s(B_{S})).
\end{equation*}
\end{definition}

\bigskip \noindent \textbf{Important Remarks}.\\

\noindent (a) The condition (CH2) is useless when $T$ is endowed with a total ordering. In that case, we may and do write the finite subsets of $T$ always in that total order.\\

\noindent (b) The main coherence condition (CH1a) may have different equivalent forms.\\

\noindent (b1) First we may write (CH1a) when $S$ has only one point more than $U$. From the new condition, we have the general one by simple induction.\\

\noindent (b2) We may also consider $V$ and $W$ two finite ordered subsets of $T$ such that $U=V\cap W$ is not empty, and as usual, we suppose that $V\cap W$ is in $V$ and in $W$ with the same ascendent order of the subscripts. Condition (CH1a) gives

$$
\mathbb{P}_{S}\Pi _{V,U}^{-1} = \mathbb{P}_{S}\Pi _{W,U}^{-1}
$$

\bigskip \noindent and this, in turn, implies (CH1a) for $U=V \subset W=S$. So we have the following new coherence condition :\\

\noindent (CH1b) \label{CH1b} For any $U=(u_{1},..,u_{r})\in \mathcal{P}_{of}(T)$ and $S=(u_{1},...,u_{r},u_{r+1})$ with  $u_{r+1}\notin U$, 
\begin{equation*}
\mathbb{P}_{U}=\mathbb{P}_{S}\Pi _{S,U}^{-1}
\end{equation*}

\bigskip \noindent which is equivalent to saying that for $B \in \mathcal{B}_{(u_{1},...,u_{r})}$,

$$
\mathbb{P}_{(u_{1},...,u_{r})}(B)=\mathbb{P}_{(u_{1},...,u_{r},u_{r+1})}(B \times E_{u_{r+1}}).
$$

\bigskip \noindent (CH1c) \label{CH1c} For any two finite and ordered subsets of $V$ and $W$ of $T$ such that $U=V\cap W$ is not empty and is  in $V$ and in $W$ with the same ascendent order of the subscripts, we have

$$
\mathbb{P}_{S}\Pi _{V,U}^{-1} = \mathbb{P}_{S}\Pi _{W,U}^{-1}.
$$

\bigskip \noindent \textbf{(II) Towards the construction of a probability law of a coherent family of marginal probability}.\\

\noindent In this part, we try to solve the following problems.\\

\noindent (i) Given a coherent (or consistence) family of real-valued, non-negative, normed and additive mappings $\mathbb{L}_{V}$, $V\in \mathcal{P}_{of}(T)$ defined on $\mathcal{B}_V$, and denoted

$$
\mathcal{F}=\{\mathbb{L}_{V},\ V\in \mathcal{P}_{of}(T)\},
$$

\bigskip \noindent \textbf{does it exists a real-valued, normed and non-negative and additive mapping}  $\mathbb{L}$ on $\mathcal{B}$ such that the elements of $\mathcal{F}$ are the finite-dimensional margins of $\mathbb{L}$, that is for any $V \in \mathcal{P}_{of}(T)$,

$$
\mathbb{L}_{V}=\mathbb{L} \Pi_{V}^{-1}?
$$

\bigskip \noindent (ii) Given a coherent family of finite-dimensional probability measures $\mathbb{P}_{V}$, $S \in \mathcal{P}_{of}(T)$, defined on $\mathcal{B}_V$ and denoted

$$
\mathcal{F}=\{\mathbb{P}_{V},\ V\in \mathcal{P}_{of}(T)\},
$$

\bigskip \noindent \textbf{does it exists a probability measure} $\mathbb{P}$ on $\mathcal{B}$ such that the elements of $\mathcal{F}$ are the finite-dimensional marginal probability measures of $\mathbb{P}$, that is for any ,

$$
\mathbb{P}_{V}=\mathbb{P} \Pi _{V}^{-1}?
$$

\bigskip \noindent Of course, if Problem (ii) is solved, Problem (i) is also solved, by taking $\mathbb{L}=\mathbb{P}$. On the other side, the solution of Problem (i) is the first step to the solution of Problem (ii).\\

\noindent We are going to see that Problem (i) has a solution with no supplementary conditions. We have

\begin{theorem} \label{preKolm} Given a coherent family $\mathcal{F}=\{\mathbb{L}_{V}, V\in \mathcal{V}_{of}(T)\}$ of non-negative, normed and additive applications, as described above, it exists \textbf{a normed and non-negative and additive application}  $\mathbb{L}$ on $\mathcal{B}$ such that the elements of $\mathcal{F}$ are the finite-dimensional margins of $\mathbb{L}$, that is for any $V \in \mathcal{P}_{of}(T)$,

$$
\mathbb{L}_{V}=\mathbb{L} \Pi_{V}^{-1}.
$$
\end{theorem}

\bigskip \noindent \textbf{Proof}. We adopt the notation introduced before to go faster. Let us suppose we are given a coherent family of 
$\{\mathbb{L}_{V}$, $S\in \mathcal{P}_{of}(T)\}$. Let us define $\mathcal{S} \subset \mathcal{P}(E)$ the class of cylinders of finite support, the following mapping

\begin{equation}
A_{V}\times E^{\prime}_V \mapsto \mathbb{L}(A_{V}\times E^{\prime}_V)=\mathbb{L}_{V}(A_{V})
\label{kolgA}
\end{equation}

\bigskip \noindent for all $V \in \mathcal{P}_{of}(T)$ and $A_V \in \mathcal{S}_V$, or in an other notation

\begin{equation}
\Pi _{V}^{-1}(A_{V}) \mapsto \mathbb{L}(\Pi _{V}^{-1}(A_{V}))=\mathbb{L}_{V}(A_{V}) \label{kolgB}
\end{equation}

\bigskip \noindent The first thing to do is to show that $\mathbb{L}$ is well-defined. Indeed, the support of $A=A_{V}\times E^{\prime}_V \in \mathcal{S}$ (with $A_V \in \mathcal{S}_V$) is not unique. Let us consider an other represent of $A$ : $A=A_{W}\times E^{\prime}_W$,  
$A_W \in \mathcal{S}_W$). If $U = V\cap W$ is empty, it means that all the factors of $A$ are full spaces and and so $A=E$ and, and since all the $\mathbb{L}_t$'s are normed, we have

$$
\mathbb{L}_{V}(A_{V})=\mathbb{L}_{W}(A_{W})=1=\mathbb{L}(A).
$$  

\bigskip \noindent If $U$ is not empty, what ever how it is ordered, it is present in $V$ according to a certain order corresponding to a permutation of $r$ of it. Also, there exists a permutation $s$ of $E_W$ such that $r(U)$ is in $W$ with the preservation of the ascendent order of the script. So we may denote $r(U)=(u_1,...,u_k)$, $V=(v_1,...,v_p)$, $s(W)=(w_1,...,w_q)$, $p\geq k$, $q\geq k$. We have

$$
A_V= A_{r(U)} \times E_{V\setminus r(U)} \ \ and \ \ A_{s(W)}=A_{r(U)} \times E_{s(W)\setminus r(U)}
$$

\noindent and it is clear that

\begin{eqnarray*}
\mathbb{L}_V(A_V)&=&\mathbb{L}_V(A_{r(U)} \times E_{V\setminus r(U)}) \ \ (L11)\\
&=&\mathbb{L}_V\left(\Pi^{-1}_{V,r(U)}(A_{r(U)})\right) \ \ (L12)\\
&=&\mathbb{L}_{r(U)}(A_{r(U)}) \ \ (L13)\\
&=&\mathbb{L}_{r(U)}(r(A_{U})) \ \ (L14)\\
&=&\mathbb{L}_{U}(A_{U}) \ \ (L15)\\
\end{eqnarray*}

\noindent In Lines (L11)-(L13), we used the coherence condition (CH1a) while (CH2) was used in Lines (L14) and (L15).\\

\noindent At the arrival, using any writing of $A \in \mathcal{S}$ leads to the same value. Then the mapping $\mathbb{L}$ is well-defined and normed.\\

\noindent In the next step, we have to show that $\mathbb{L}$ is additive of $\mathcal{S}$. For this, let us consider an element of $\mathcal{S}$ that is split into two disjoint elements of $\mathcal{S}$. Suppose

\begin{equation*}
A=B+C
\end{equation*}

\noindent with 

\begin{equation*}
A=A_{U} \times E^{\prime}_U, \ \ B=A_{V} \times E^{\prime}_V \ \ and \ \
C=C_{W} \times E^{\prime}_W.
\end{equation*}

\bigskip \noindent Let consider $Z=U \cup V \cup W$ given in some order of the subscripts. There exist permutations $s$, $r$ and $p$ of $E_U$, $E_V$ and $E_W$ respectively such that  $r(U)$, $s(V)$ and $p(W)$ are given in $Z$ with the preservation of the ascendent order of the subscripts and we have :

\begin{equation*}
A=\biggr( A_{r(U)} \times E_{Z\setminus r(U)}\biggr) \times E^{\prime}_Z \equiv A^{\ast}_{Z} \times E^{\prime}_Z,
\end{equation*}

\begin{equation*}
B=\biggr( B_{s(V)} \times E_{Z\setminus s(V)}\biggr) \times E^{\prime}_Z \equiv B^{\ast}_{Z} \times E^{\prime}_Z,
\end{equation*}

\bigskip \noindent and

\begin{equation*}
C=\biggr( C_{p(Z)} \times E_{Z\setminus p(W)}\biggr) \times E^{\prime}_Z \equiv C^{\ast}_{Z} \times E^{\prime}_Z,
\end{equation*}

\bigskip \noindent with

$$
A^{\ast}_{Z} \times E^{\prime}_Z = \biggr(B^{\ast}_{Z} \times E^{\prime}_Z \biggr)+ \biggr(C^{\ast}_{Z} \times E^{\prime}_Z\biggr)
$$

\bigskip \noindent This is possible only if we have

$$
A^{\ast}_{Z} = B^{\ast}_{Z} +C^{\ast}_{Z},
$$

\bigskip \noindent with

$$
\mathbb{L}_Z(A^{\ast}_{Z})=\mathbb{L}_U(A_U), \ \ \mathbb{L}_Z(B^{\ast}_{Z})=\mathbb{L}_V(B_V), \ \ and \ \
\mathbb{L}_Z(C^{\ast}_{Z})=\mathbb{L}_W(C_Z)
$$

\bigskip \noindent Using the coherence conditions, we have

\begin{eqnarray*}
\mathbb{L}(B)&=&\mathbb{L}_{Z}(B^{\ast}_Z)) \ \ \ (L31)\\
&=&\mathbb{L}_{Z}\left(\Pi_{r(V),Z}^{-1}\left(B_{r(V)}\right)\right) \ \  (L32)\\
&=&\mathbb{L}_{r(V)}\left(B_{r(V)}\right) \  (L33)\\
&=&\mathbb{L}_{V}(B_V).
\end{eqnarray*}

 \bigskip \noindent By doing the same for $A$ and $C$, we have

\begin{eqnarray*}
&&\mathbb{L}(A)=\mathbb{L}_{U}(A_U)=\mathbb{L}_{Z}(A^{\ast}_Z)),\\
&&\mathbb{L}(B)=\mathbb{L}_{V}(B_V)=\mathbb{L}_{Z}(B^{\ast}_Z)),\\
&&\mathbb{L}(C)=\mathbb{L}_{W}(C_W)=\mathbb{L}_{Z}(C^{\ast}_Z))\\
\end{eqnarray*}

\noindent Bu using the additivity of $\mathbb{L}_{Z}$, we conclude that

$$
\mathbb{L}(A)=\mathbb{L}(B)+\mathbb{L}(C).
$$

\bigskip \noindent The mapping $\mathbb{L}$ is normed and additive on the semi-algebra. From Measure Theory and Integration (See Doc 04-02, Exercise 15, in \cite{ips-mestuto-ang}),  $\mathbb{L}$ is automatically extended to a normed and additive mapping on the algebra $\mathcal{C}$ generated by $\mathcal{S}$.\\

\bigskip \noindent Now, we face Problem (ii). Surely, the assumptions and the solution Problem (i) ensure that there exists a normed and additive mapping $\mathbb{L}$ whose margins are the $\mathbb{P}_V$, $V \in \mathcal{P}_{of}(T)$. Let us call it $\mathbb{P}$. All we have to do is to get an extension of $\mathbb{P}$ to $\sigma(\mathcal{C})=\mathcal{B}$.\\

\noindent A way to do it is to use Caratheodory's Theorem (See \textbf{Doc 04-03 in \cite{ips-mestuto-ang}} for a general revision). But, unfortunately, we need special spaces. Suppose that each $(E_t,d_t)$, $t \in T$, is Polish space, that is a metric separable and complete separable space. The following facts are known in Topology. For $V \in \mathcal{P}_{of}(T)$, the space $E_V$ is also a Polish space. In such spaces, the extension of $\mathbb{P}$ to probability measure is possible. The proof heavily depends on topological notions, among them a characterization of compact sets.\\

\noindent We give the proof in the last section as an Appendix. In the body of the text, we focus on probability theory notions. However, we strongly recommend the learners to read the proof in small groups. We have the following Theorem.

\begin{theorem} \label{KolmConst} (Fundamental Theorem of Kolmogorov) Let us suppose that each $(E_t,d_t)$, $t \in T$ is Polish Space. For any $V \in \mathcal{P}_{of}(T)$, $E_V$ is endowed with the Borel $\sigma$-algebra associated with the product metric of the metrics of its factors.\\

\noindent For $T \neq \emptyset$, given a coherent family $\mathcal{F}=\{\mathbb{P}_{V}, \ V\in \mathcal{P}_{of}(T)\}$ of finite-dimensional probability measures, there exists a unique probability measure on $\mathcal{B}$ such that the elements of $\mathcal{F}$ are the finite-dimensional margins of $\mathbb{L}$, that is for any $V \in \mathcal{P}_{of}(T)$,

$$
\mathbb{P}_{V}=\mathbb{P} \Pi_{V}^{-1}.
$$
\end{theorem}

\bigskip \noindent Now, we are going to derive different versions of that important basis of Probability Theory and provide applications and examples.\\

\noindent To begin, let us see how to get the most general forms the Kolmogorov construction in finite dimensions (See Chapter \ref{proba_02_rv} in pages \pageref{kolmconst_01} and \pageref{kolmconst_02}). Let us repeat a terminology we already encountered. For any mapping
$$
X \ : \ (\Omega, \mathcal{A}, \mathbb{P}) \rightarrow (E, \mathcal{B})
$$

\noindent we have $X(\omega)=(X_t(\omega))_{t \in T}$. For any $V=(v_1, ...,v_k) \in \mathcal{P}_{of}(T)$, $k\geq 1$,

$$
X_V\equiv (X_{v_1}, ..., X_{v_k})=\Pi_V(X),
$$

\noindent is called a finite-dimensional (and ordered) margin of $X$. We have

\begin{theorem} \label{KolmConstTraj} Let us suppose that each $(E_t,d_t)$, $t \in T$ is a Polish Space. For any $V \in \mathcal{P}_{of}(T)$, $E_S$ is endowed with the Borel $\sigma$-algebra associated with the product metric of the metrics of  its factors.\\

\noindent For $T \neq \emptyset$, given a coherent family $\mathcal{F}=\{\mathbb{P}_{V}, \ V\in \mathcal{P}_{of}(T)\}$ of probability measures, there exists a probability space $(\Omega, \mathcal{A}, \mathbb{P})$ and a measurable mapping $X$ defined on $\Omega$ with values in
$(E, \mathcal{B})$ such that for any $V \in \mathcal{P}_{of}(T)$, $k\geq 1$,

$$
\mathbb{P}_{V}=\mathbb{P}_{X_V}=\mathbb{P} X_{V}^{-1}.
$$

\bigskip \noindent In other words, there exists a probability space $(\Omega, \mathcal{A}, \mathbb{P})$ holding a measurable mapping $X$ with values in
$(E, \mathcal{B})$ such that the finite-dimensional marginal probability measures $\mathbb{P}_{V}$, $V\in \mathcal{P}_{of}(T)$ are the probability laws of the finite-dimensional (and ordered) margins $X_V$ of $X$.\\

\bigskip \noindent Furthermore, the probability laws of the finite-dimensional (and ordered) margins $X_V$ of $X$ determine the probability law of $X$.
\end{theorem}

\noindent \textbf{Proof}.  We apply Theorem \ref{KolmConst} above to get the probability measure on $\mathbb{P}$ on $(E, \mathcal{B})$ whose finite-dimensional marginal probabilities are the $\mathbb{P}_{V}$, $V\in \mathcal{P}_{of}(T)$. Now we take

$$
(\Omega, \mathcal{A}, \mathbb{P})=(E, \mathcal{B}, \mathbb{P})
$$

\noindent * and set $X$ as the identity mapping. We have

\begin{eqnarray*}
\mathbb{P}_{X_V}&=&\mathbb{P} X_{V}^{-1}= \mathbb{P} \Pi_V(X)^{-1}\\
&=&\mathbb{P} \Pi_V^{-1} X^{-1}\\
&=&\mathbb{P}_V X^{-1}.
\end{eqnarray*}

\noindent Now for any $B_V \in \mathcal{B}_V$,

$$
\mathbb{P}_{X_V}(B_V)=\mathbb{P}_V(\{\omega \in \Omega, X(\omega)=\omega \in B_V\})=\mathbb{P}_V(B_V).
$$

\bigskip \noindent We get the desired result : $\mathbb{P}_{X_V}=\mathbb{P}_{V}$, $V\in \mathcal{P}_{of}(T)$. To finish, the probability law of $X$ is given by

$$
\mathbb{P}_X(B)=\mathbb{P}X^{-1}(B), \ B\in \mathcal{B}.
$$ 

\bigskip \noindent By the uniqueness of the Caratheodory's extension for a $\sigma$-additive and proper mapping of from an algebra to the $\sigma$-algebra generated, the probability measure $\mathbb{P}X^{-1}$ on $\mathcal{B}$ is characterized by its values on $\mathcal{S}$. But an element of $\mathcal{S}$ is of the form 
$$ 
B=B_V \times E^{\prime}_V, B_V \in \mathcal{B}, \ V\in \mathcal{P}_{of}(T).
$$

\bigskip \noindent We have have 

$$
X \in B=B_V \times E^{\prime}_V \Leftrightarrow X_V \in B_V
$$

\bigskip \noindent that is

$$
X^{-1}(B)=X_V^{-1}(B_V).
$$

\bigskip \noindent Hence, by applying $\mathbb{P}$ at both sides, we get

$$
\mathbb{P}_X(B_V \times E^{\prime}_V)=\mathbb{P}_{X_V}(B_V).
$$

\noindent * Since the values of $\mathbb{P}_X$ are functions only of the values of the probability laws of the finite-dimensional (and ordered) margins $X_V$ of $X$, these latter finally determine $\mathbb{P}_X$. $\square$\\

\noindent We are continuing to see developments of the Kolmogorov Theorem in special sections.

\section[Skorohod's Construction]{Skorohod's Construction of real vector-valued stochastic processes} \label{proba_02_skoro}

\noindent \textbf{(I) - The General Theorem}.\\

\noindent Let $T$ be an non-empty index set. For each $t\in T$, let be given $E_t=\mathbb{R}^{d(t)}$, where $d(t)$ is positive integer number. Let us consider a family of probability distribution functions described as follows : for $V=(v_1,..., v_k)\in \mathcal{P}_{of}(T)$, $k\geq 1$, we set $d(V)=d(v_1)+...+d(v_k)$. The probability distribution function associated to $V$ is defined for $x_{v_j} \in \mathbb{R}^{d(v_j)}$, $1\leq j \leq k$, by

$$
\mathbb{R}^{d(V)} \ni (x_{v_1}, ..., x_{v_k}) \mapsto F_V(x_{v_1}, ..., x_{v_k}).
$$

\bigskip \noindent The family of $\{F_V, \ V=(v_1,..., v_k)\in \mathcal{P}_{of}(T)\}$ is coherent if and only if :\\

\noindent \label{CHS1} (CHS1), for $V \in \mathcal{P}_{of}(T)$, for any $(x_{v_1}, ..., x_{v_k}) \in \mathbb{R}^{d(V)}$, for any permutation of 
$E_V$,

$$
F_V(x_{v_1}, ..., x_{v_k})=F_{s(V)}(s(x_{v_1}, ..., x_{v_k}))
$$

\bigskip \noindent and \\

\noindent \label{CHS2} (CHS2) for $V \in \mathcal{P}_{of}(T)$, for any $(x_{v_1}, ..., x_{v_k}) \in \mathbb{R}^{d(V)}$, for any $u \in T\setminus V$,

$$
F_V(x_{v_1}, ..., x_{v_k})=\lim_{u \uparrow \{+\infty\}^{d(u)}}F_{V}(x_{v_1}, ..., x_{v_k},u).
$$

\noindent In all this part, by writing $(x_{v_1}, ..., x_{v_k}) \in \mathbb{R}^{d(V)}$, we also mean that  $x_{v_j} \in \mathbb{R}^{d(v_j)}$, for all $j \in \{1,...,k\}$.\\

\noindent Here is the Skorohod Theorem as follows.\\

\begin{theorem} \label{proba_02_skorohod01} Given a coherent family of probability distribution functions $\{F_V, \ V=(v_1,..., v_k)\in \mathcal{P}_{of}(T)\}$, there exists a probability space $(\Omega, \mathcal{A}, \mathbb{P})$ holding a measurable mapping $X$ with values in
$(E, \mathcal{B})$ such that each finite-dimensional marginal probability distribution function $F_{X_V}$, $V\in \mathcal{P}_{of}(T)$ is $F_V$, that is for  of $V=(v_1,..., v_k)$, $k\geq 1$, for $(x_{v_1}, ..., x_{v_k}) \in \mathbb{R}^{d(V)}$,

$$
F_V(x_{v_1}, ..., x_{v_k})=\mathbb{P}(X_{v_1}\leq x_{v_1}, ...,X_{v_k}\leq x_{v_k}).
$$
\end{theorem}

\noindent \textbf{Proof}. The proof results from the application of The Kolmogorov Theorem and a smart use of the Lebesgue-Stieljes measures. We remind first that for any $\ell \geq 1$, a finite measure on $\mathbb{R}^{\ell}$ is characterized by its values of the elements of the form

$$
]-\infty,a]= \prod_{1\leq j \leq \ell} ]-\infty, a_j], \ a=(a_1,...,a_{\ell})
$$   

\bigskip \noindent which form a $\pi$-system denoted $\mathcal{D}_{\ell}$, which in turn, generates $\mathcal{B}(\mathbb{R}^{\ell})$. Consider the unique Lebesgue-Stieljes probability measure $\mathbb{P}_{V}$ on $\mathbb{R}^{d(V)}$ associated with $F_V$, $V=(v_1,..., v_k)\in \mathcal{P}_{of}(T)\}$, $k\geq 1$. By keeping the previous notation, we have that for any $(x_{v_1}, ..., x_{v_k}) \in \mathbb{R}^{d(V)}$ and for any permutation of $E_V$

\begin{eqnarray*}
\mathbb{P}_{V}(v)&=&F_V(x_{v_1}, ..., x_{v_k})\\
&=&F_{s(V)}(s(x_{v_1}, ..., x_{v_k}))\\
&=&\mathbb{P}_{s(V)}\left(s \left(\prod_{1\leq j \leq k} ]-\infty, x_{v_j}]\right)\right)\\
&=&\mathbb{P}_{s(V)} s^{-1} \left(\prod_{1\leq j \leq k} ]-\infty, x_{v_j}]\right).\\
\end{eqnarray*}

\noindent Since the probability measures $\mathbb{P}_{V}$ and $\mathbb{P}_{s(V)} s^{-1}$ coincide on $\mathcal{D}_{d(V)}$, they are equal and the first coherence condition is proved. To prove the second, we have

\begin{eqnarray*}
\mathbb{P}_{V}\left(\prod_{1\leq j \leq k} ]-\infty, x_{v_j}]\right)&=&F_V(x_{v_1}, ..., x_{v_k})\\
&=& \lim_{\mathbb{R}^{d(u)} \uparrow \{+\infty\}^{d(u)}}F_{V \cup \{u\}}(s(x_{v_1}, ..., x_{v_k}), u)\\
&=& \lim_{\mathbb{R}^{d(u)} \uparrow \{+\infty\}^{d(u)}} \mathbb{P}_{V \cup \{u\}}\left( \prod_{1\leq j \leq k} ]-\infty, x_{v_j}]
\times ]-\infty, u]\right)\\
&=&\mathbb{P}_{V \cup \{u\}}\left( \prod_{1\leq j \leq k} ]-\infty, x_{v_j}]
\times E_u \right).
\end{eqnarray*}

\bigskip \noindent But the two probability measures on $\mathbb{R}^{d(V)}$ : $\mathbb{P}_{V}(B)$ and $\mathbb{P}_{V \cup \{u\}}\left(B \times E_u\right)$, $B \in \mathcal{B}_V$, coincide on $\mathcal{D}_{d(V)}$. Hence for any $B \in \mathcal{B}_V$, we have

$$
\mathbb{P}_{V}(B)=\mathbb{P}_{V \cup \{u\}}\left(B \times E_u\right).
$$

\noindent Thus, the coherence condition (CH1b) holds. Finally there exists a probability measure $\mathbb{P}$ on $(E, \mathcal{B})$ whose finite-dimensional marginal probability measures are the elements of $\{\mathbb{P}_V, \ V\in \mathcal{P}_{of}\}$. Let us take

$$
(\Omega, \mathcal{A}, \mathbb{P})=(E, \mathcal{B}, \mathbb{P})
$$

\noindent and set $X$ as the identity mapping. We have, for any $(x_{v_1}, ..., x_{v_k}) \in \mathbb{R}^{d(V)}$,

\begin{eqnarray*}
\mathbb{P}(X_{v_1}\leq x_{v_1}, ...,X_{v_k}\leq x_{v_k})&=&\mathbb{P}(\{\omega \in E, \omega_{v_1}\leq x_{v_1}, ..., \omega_{v_k}\leq x_{v_k}\})\\
&=&\mathbb{P}(\prod_{1\leq j \leq k} ]-\infty, x_{v_j}] \times E^{\prime}_V)\\
&=&\mathbb{P}_V\left(\prod_{1\leq j \leq k} ]-\infty, x_{v_j}]\right)\\
&=&F_V(x_{v_1}, ...,X_{v_k}).
\end{eqnarray*}

\bigskip \noindent The proof is finished. $\square$

\bigskip \noindent \textbf{Other forms of the Skorohod Theorem using densities of probability}.\\

\noindent Suppose we have the similar following situation as earlier. For each $t\in T$, let be given $E_t=E_0^{d(t)}$, where $d(t)$ is positive integer number. Let $\nu$ be a $\sigma$-finite measure on $E_0$. On each $E_V=E_0^V$, $V \in \mathcal{P}_{of}(T)$, we have a the finite product probability :

$$
\nu_V= \nu^{\otimes d(v_1)}  \nu\otimes^{\otimes d(v_2)} \otimes \cdots \otimes \nu^{\otimes d(v_k)}=\nu^{\otimes d(V)}.  
$$

\noindent * Now a family of marginal probability density functions (\textit{pdf}) $\{f_V, \ V=(v_1,..., v_k)\in \mathcal{P}_{of}(T)\}$, each $f_V$ is \textit{pdf} with respect to $\nu^{}\otimes d(V)$ on $E_0^{d(V)}$, is said to be coherent if the two conditions hold :\\

\noindent \label{CHSD1} (CHSD1) For $V \in \mathcal{P}_{of}(T)$, for any $(x_{v_1}, ..., x_{v_k}) \in E_0^{d(V)}$, for any permutation of 
$E_V$, we have

$$
f_V(x_{v_1}, ..., x_{v_k})=f_{s(V)}(s(x_{v_1}, ..., x_{v_k}))
$$

\bigskip \noindent and, \\

\noindent (CHSD2) \label{CHSD2} for $V \in \mathcal{P}_{of}(T)$, for any $(x_{v_1}, ..., x_{v_k}) \in \mathbb{R}^{d(V)}$, for any $u \in T\setminus V$,

$$
f_V(x_{v_1}, ..., x_{v_k})=\int_{E_{v_u}} f_{V}(x_{v_1}, ..., x_{v_k},u) d\nu^{\otimes d(v_u)}(u),
$$

\bigskip \noindent meaning that $f_V(x_{v_1}, ..., x_{v_k})$ is a marginal \textit{pdf} of $f_{V}(x_{v_1}, ..., x_{v_k},u)$.\\

\bigskip \noindent Let us consider the finite distribution probability measure on $E_V$ defined by

$$
E_V (=E_0^{d(V)}) \ni B \mapsto \mathbb{P}_V(B)= \int_{E_{V}} f_{V}(x_{v_1}, ..., x_{v_k}) d\nu^{\otimes d(V)}(x_{v_1}, ..., x_{v_k}),
$$

\bigskip \noindent It is easy to see that (CHSD1) and (CHSD2) both ensure that the $\mathbb{P}_V$ form a coherent family of finite dimensional probability measures. We apply theorem to conclude that :\\

\noindent For any coherent family of marginal probability density functions (\textit{pdf}) $\{f_V, \ V=(v_1,..., v_k)\in \mathcal{P}_{of}(T)\}$, each $f_V$ is \textit{pdf} with respect $\nu^{\otimes d(V)}$ on $E_0^{d(V)}$, there exists a probability space $(\Omega, \mathcal{A}, \mathbb{P})$ holding a measurable mapping $X$ with values in $(E, \mathcal{B})$ such that each finite-dimensional \textit{pdf} $f_{X_V}$, $V\in \mathcal{P}_{of}(T)$ is $F_V$, that is for  of $V=(v_1,..., v_k)$, $k\geq 1$, for $(x_{v_1}, ..., x_{v_k}) \in \mathbb{R}^{d(V)}$,

$$
d\mathbb{P}(X_{v_1}, \cdots , X_{v_k})=f_V d\nu^{\otimes d(V)}.
$$

\bigskip \noindent In general, this is used in the context of $\mathbb{R}^d$, $d\geq 1$, with $\nu$ being the Lebesgue measure or a counting measure on $\mathbb{R}$. But it goes far beyond as a general law.

\section{Examples} \label{thfondamentalKolm_examples}

\noindent To make it simple, let $T=\mathbb{R}_+$ or $T=\mathbb{N}$. So we do not need to care about the first coherence condition since we have a natural order. Let be given $E_t=E_0^{d(t)}$, where $E_0$ is a polish space and $d(t)$ is positive integer number.\\

\bigskip \noindent \textbf{Problem 1}. Given a family of Probability measures $\mathbb{P}_{t}$ on each $E_t$, of dimension $d(t)$. Does-it exist a probability space $(\Omega, \mathcal{A}, \mathbb{P})$ holding a stochastic process $(X_t)_{t\in T}$ with independent margins such that each margin $X_t$ follows the probability law $\mathbb{P}_{t}$.\\

\noindent \textbf{Solution}. We can easily that the family of finite dimensional probability measure, 

$$
\mathbb{P}_{(t_1,t_2,...,t_k}=\bigotimes_{j=1}^{k} \mathbb{P}_{t_j};
$$

\bigskip  \noindent for $t_1<...<t_k$, defined by, for any $B_j \in \mathcal{B}(E_t)$, $j \in \{1,...,k\}$,

$$
\mathbb{P}_{(t_1,t_2,...,t_k}\left(\prod_{1\leq j \leq k} B_j\right)=\prod_{j=1}^{k} \mathbb{P}_{t_j}(B_j).
$$

\bigskip  \noindent is coherent since for $t_{j+1}>t_j$,

\begin{eqnarray*}
\mathbb{P}_{(t_1<t_2<...<t_k}\left(\prod_{1\leq j \leq k} B_j \times E_{j+1}\right)&=&\prod_{j=1}^{k} \mathbb{P}_{t_j}(B_j) \times \mathbb{P}_{t_{j+1}}(E_{j+1})\\
&=& \prod_{j=1}^{k} \mathbb{P}_{t_j}(B_j).
\end{eqnarray*}

\noindent Thus, the answer is positive.\\

\bigskip \noindent \textbf{Problem 2}. Many techniques are based on the symmetrization method as in the proof of Proposition \ref{proba_02_conv_propd02} (See page \pageref{symmetrizezMzthod}). We need to have two sequences $(X_n)_{n\geq 1}$ and $(Y_n)_{n\geq 1}$ on the same probability space and having their values on $\mathbb{R}^{d}$ such that $X_n=_d Y_n$ for each $n\geq 1$. Is it possible?\\

\noindent Here is the statement of the problem for independent margins.\\

\noindent Given a family of Probability measures $\mathbb{P}_{t}$ on each $E_t$, of dimension $d(t)$. Does-it exist a
 probability space $(\Omega, \mathcal{A}, \mathbb{P})$ holding a stochastic process $(X_t)_{t\in T}$ with independent margins such that :\\

\noindent (a) $X_t \in E_t^2$, that is  $X_t=(X_t^{(1)},X_t^{(2)})^{t}$\\

\noindent (b) For each $t\in T$, for each $i\in \{1,2\}$, $\mathbb{P}_{X_t^{(i)}}=\mathbb{P}_{t}$.\\

\noindent If this problem is solved and if $E_0$ is a linear space, we may form the symmetrized form $X^{(s)}=X_t^{(1)}-X_t^{(2)}$ with
$X_t^{(1)}=_{d}X_t^{(2)}$.\\

\noindent \textbf{Solution}. Let us apply the solution of Problem 1 for the case where $E_0=\mathbb{R}$ in the context of independent margins. We notice that nothing is said about the dependence between $X_t^{(1)}$ and $X_t^{(2)}$. So we may take, for any $t \in T$, an arbitrary probability distribution function $F_t$ on $\mathbb{R}^{2d(t)}$ such that the margins 

$$
F_{t}(\underset{d(t)\text{ times}}{\underbrace{x_{1},\cdots, x_{d(t)}}}, \underset{d(t)\text{ times}}{\underbrace{+\infty,\cdots,+\infty}})
$$   

\bigskip \noindent and

$$
F_{t}(\underset{d(t)\text{ times}}{\underbrace{+\infty,\cdots,+\infty}}, \underset{d(t)\text{ times}}{\underbrace{x_{d(t)+1},\cdots, x_{2d(t)}}} 
),
$$   

\noindent are equal both the probability distribution function of $\mathbb{P}_t$. This is possible by the use of copulas. With such a frame, we apply again the Skorohod Theorem to get our solution.\\

\noindent \textbf{Problem 3}. Existence of the Poisson Process. Given $\theta>0$, by the solution of Problem 1, there exists a probability space $(\Omega, \mathcal{A}, \mathbb{P})$ holding a sequence independent random variables identically distributed as the standard exponential law $\mathcal{E}$ denoted $X_1$, $X_2$, etc.\\

\noindent Let us call them the independent and exponential inter-arrival times.\\

\noindent Let us define the arrival times $Z_0$, $Z_j=X_1+\cdots+X_j$, $j\geq 1$, so that we have
$$
Z_0 < Z_1 < \cdots Z_j \cdots
$$  

\bigskip \noindent If we suppose that the $Z_j$ are the arrival times of clients at a desk (say a bank desk) and $Z_0=0$ is the opening time of the desk, we may wish to know the probability law of the number at arrived clients at a time $t>0$,

$$
N(]0,t])=N_t=\sum_{j\geq 1} 1_{(Z_j\leq t)}.
$$ 

\bigskip \noindent Here, we say that we have a standard Poisson Process (\textit{SPP}) of intensity $\theta$. Sometimes, authors mean $(N_t)_{t\geq 0}$ which is the counting function of the \textit{SPP}, others mean $(X_n)_{n\geq 1}$ which is the sequence of arrival times or $(Z_n)_{n\geq 1}$ which is the sequence of inter-arrival times.\\

\bigskip \noindent \textbf{Problem 4}. Existence of Brownian Movement by the exercise. Let $0 = t_{0} <t_{1} <... <t_{n}$ be $n$ real numbers and consider $Y_{1}, Y_{2}, ..., Y_{n}$,  $n$ non-centered Gaussian with respective variances $t_{1}, t_{2} -t_{1}, ..., t_{n}-t_{n-1}$. Set

\begin{equation*}
X=(X_{1}, X_{2}, ..., X_{n}) = (Y_{1}, Y_{1} + Y_{2}, ..., Y_{1} + {Y_2} + \cdots + Y_{n}).
\end{equation*}

\bigskip  \noindent a) Find the density of $X$.\\

\noindent b) Give the distribution function of $X$.\\

\noindent c) Now, consider that family of distribution functions indexed by the ordered and finite subsets of $\mathbb{R}_{+}$ : for any
$(x_{1}, x_{2}, ..., x_{k}) \in \mathbb{R}^k$, $k\geq 1$,
\begin{eqnarray*}
&&F_{(t_{1}, t_{2}, ... <t_{k})} (x_{1}, x_{2}, ..., x_{k}) \ \ \ \ \ \ \ \ \ (BR01)\\
&=&\int _{-\infty}^{x_{1}} dy_{1} \int_{-\infty}^{x_{2}} dy_{2} \cdots \int_{-\infty}^{x_{k}} \prod_{i=1}^{k} 
\frac{1}{\sqrt {2 \pi(t_{i} -t_{i-1})}} \exp \biggr(-\frac{1}{2} \frac{(y_{i} -y_{i-1)}^{2}}{t {i} -t {{-1}}}\biggr) \ dy_{k},
\end{eqnarray*}

\bigskip  \noindent where $(t_{1} <t_{2}, ... <t_{k}) $ is an ordered and finite subset of $\mathbb{R}_{+}$, with $y_{0}=t_{0}=0$.\\

\noindent c1) Say on the basis of questions (a) and (b), why do we have, for all $(t_{1} <t_{2}, ... <t_{k})$, $k \geq 2$; for all 
$(x_{1}, x_{2}, ..., x_{k-1}) \in \mathbb{R}^{k-1}$

\begin{equation*}
\lim_{x_{k} \uparrow \infty} {F (t_{1}, t_{2}, ..., t_{k})} (x_{1}, x_{2}, ..., x_{k}) = F_{(t_{1}, t_{2}, ..., t_{k-1})}(x_{1}, x_{2}, ..., x_{k-1})
\end{equation*}

\bigskip \noindent c2) Show this property  directly from the definition (BR01).\\

\noindent c3) Conclude by the Kolmogorov-Skorohod, that there is a stochastic process 
$(\Omega, \mathcal{A}, \mathbb{P}, (B_{t})_{t \in \mathbb{R} _{+}})$ for which we have

\begin{equation*}
F _{(B (t_{1}), B (t_{2}), ..., B (t_{k}))} (x_{1}, x_{2}, ..., x_{k }) = F _{(t_{1}, t_{2}, ... <t_{k})} (x_{1}, x_{2}, ..., x_{k})
\end{equation*}

\bigskip \noindent for any  finite and ordered subset $(t_{1}, t_{2}, ... <t_{k})$ of $\mathbb{R}_{+}$.\\

\noindent Alternatively, use the system of \textit{pdf}'s indexed by the ordered and finite subsets of $\mathbb{R}_{+}$ : for any
$(x_{1}, x_{2}, ..., x_{k}) \in \mathbb{R}^k$, $k\geq 1$,

\begin{eqnarray*}
&&f_{(t_{1}, t_{2}, ... <t_{k})} (x_{1}, x_{2}, ..., x_{k}) \ (BR02)\\
&=&\prod_{i=1}^{k} \frac{1}{\sqrt {2 \pi(t_{i} -t_{i-1})}} \exp \biggr(-\frac{1}{2} \frac{(y_{i} -y_{i-1)}^{2}}{t {i} -t {{-1}}}\biggr)
\end{eqnarray*}

\bigskip  \noindent where $(t_{1} <t_{2}, ... <t_{k}) $ is an ordered and finite subset of $\mathbb{R}_{+}$, with $y_{0}=t_{0}=0$, and
the coherence condition (CHSD2) (page \pageref{CHSD2}) to justify the existence of such a process.\\

\noindent d) Such a stochastic process $(B_{t})_{t \in \mathbb {R}_{+}}$ is called \textit{Brownian motion} in Probability Theory and 
\textit{Wiener Process} in Statistics.\\

\noindent Show or state the following facts.\\

\noindent d1) Its finite distributions are non-centered Gaussian vectors.\\

\noindent d2) For all $ 0\leq s <t$, $B_{t} -B_{s}$ and $B_{s}$ are independent.\\

\noindent d3) For all $t\geq 0$, 

$$
B_{t} -B_{s} \sim B (t-s) \sim \mathcal{N}(0, t-s).
$$

\bigskip  \noindent d4) $\Gamma(s, t)=\mathbb{C}ov(B_{t}, B_{s}) = \min(s, t)$, \ $(t,s) \in \mathbb{R}_{+}^2$.\\

\noindent (e) A stochastic process $(X_t)_{t\geq 0}$ is said to be Gaussian if and only if its finite margins are Gaussian vectors.\\

\noindent Show the following points :\\

\bigskip \noindent (e1) Show that the probability law of a Gaussian Process is entirely determined by its mean function

$$
m(t)=\mathbb{E}(X_t), \ t\in \mathbb{R}_{+}.
$$

\bigskip \noindent and by its variance-covariance function

$$
\Gamma(s, t)=\mathbb{C}ov(X_{t}, X_{s}), \ (t,s) \in \mathbb{R}_{+}^2.
$$

\bigskip \noindent (e2) Deduce from this that the probability law of the Brownian Process is entirely the variance-covariance function

$$
\Gamma(s, t)=\mathbb{C}ov(X_{t}, X_{s})=\min(s,t), \ (t,s) \in \mathbb{R}_{+}^2.
$$


%% file: proba_02_09_proof_ang.tex
\section[Proof of the Fundamental Theorem of Kolmogorov]{Caratheodory's Extension and Proof of the Fundamental Theorem of Kolmogorov} \label{proba_02_caraKolm}

\noindent Our departure point is the end of the Proof of Theorem \ref{preKolm}. The construction mapping $\mathbb{L}$; we denote now as 
$\mathbb{P}$ is additive and normed on the algebra $\mathcal{C}=a(\mathcal{S})$ generated by $\mathcal{S}$. By Carath\'{e}odory Theorem (Doc 04-03 in \cite{ips-mestuto-ang}), $\mathbb{P}$ is uniquely extensible to a probability measure whenever it is continuous at  $\emptyset$, that is, as $n\rightarrow \infty $,

\begin{equation*}
\biggr(\mathcal{C\ni } \ A_{n} \downarrow \emptyset\biggr) \implies \biggr(\mathbb{P}(A_{n})\downarrow 0\biggr).
\end{equation*}

\bigskip \noindent Actually, we are going to use an \textit{ab contrario} reason. Suppose that there exists a non-increasing sequence
$(A_n)_{n\geq 0} \subset \mathcal{C}$ and $\mathbb{P}(A_{n})$ does not converges to zero. Since the sequence $(\mathbb{P}(A_{n}))_{n\geq 1}$ is non-increasing, its non-convergence to zero is equivalent to 
\begin{equation*}
\exists \varepsilon >0), \ (\forall n\geq 1,\mathbb{P}(A_{n})>\varepsilon).
\end{equation*}

\bigskip \noindent At the beginning let us remark that $\mathcal{C}=a(\mathcal{S})$ is formed by the finite sum of elements of $\mathcal{S}$, we rely on the above considerations on finite sums of elements of $\mathcal{S}$, and easily get that any element of $\mathcal{C}$, and then any $A_n$ is of the form
$$
A_n=B_{V_{n}} \times E^{\prime}_{V_{n}}, \ n\geq 1,
$$

\bigskip \noindent where $B_{V_{n}} \in \mathcal{B}_{V_{n}}$, and ${V_{n}} \in \mathcal{P}_{of}(T)$. Hence the whole sequence does involve only a countable spaces $E_t$, $t \in T_0$, where

\begin{equation*}
T_{0}=\bigcup\limits_{j=1}^{\infty }V_{n}
\end{equation*}

\bigskip \noindent and denote, accordingly,  
\begin{equation*}
E_{T_{0}}=\prod\limits_{t\in T_{0}}E_{t}.
\end{equation*}

\bigskip \noindent So we may ignore all the other factors in $E$.

\begin{equation*}
A_{n}=B_{V_{n}}\times \prod\limits_{t\notin T_{0}\backslash V_{n}}E_{t}.
\end{equation*}

\bigskip \noindent As well, we may use the natural order of integers and only consider supports index sets of the form $V_n=(1,...,m(n))$, $n\geq 1$. Now we are going to use the following key topological property : in a Polish space, for any Borel set $B$, for any finite measure $\mu$, for any $\varepsilon >0$, there exists a compact set $K(\varepsilon)$ such
$$
\mu(B\setminus K(\varepsilon)) < \varepsilon. 
$$

\bigskip \noindent Then for any $\varepsilon >0$, for any $n\geq 1$, there exists a compact set $B^{\prime}_{V_{n}} \subset B_{V_{n}}$ such that

\begin{equation*}
P_{V_{n}}(B_{V_{n}}-B_{V_{n}}^{\prime })\leq \varepsilon 2^{-(n+1)}.
\end{equation*}

\bigskip \noindent Let us denote 
\begin{equation*}
A_{n}^{\prime}=B_{V_{n}}^{\prime }\times E^{\prime}_{T_0 \setminus V_{n}}
\end{equation*}

\bigskip \noindent Hence for each $n\geq 1$,
\begin{equation*}
P(A_{n}-A_{n}^{\prime })=P_{V_{n}}(B_{V_{n}}-B_{V_{n}}^{\prime})<\varepsilon 2^{-(n+1)}.
\end{equation*}

\bigskip \noindent Let us set

\begin{equation*}
C_{n}=A_{1}^{\prime }\cap ...\cap A_{n}^{\prime }.
\end{equation*}

\bigskip \noindent We have for each $n\geq 1$,

\begin{eqnarray*}
P(A_{n}-C_{n})&=&P(A_{n}\cap (\bigcup\limits_{j=1}^{n}(A_{j}^{\prime})^{c}))\leq \sum\limits_{j=1}^{n}P(A_{n}\cap (A_{j}^{\prime })^{c}))\\
&\leq& \sum\limits_{j=1}^{n}P(A_{n}-A_{j}^{\prime })\\
&\leq&\sum\limits_{j=1}^{n}P(A_{j}-A_{j}^{\prime })\\
&\leq&\sum\limits_{j=1}^{n}\varepsilon 2^{-(j+1)}<\varepsilon /2.
\end{eqnarray*}

\noindent But $C_{n}\subset A_{n}^{\prime}\subset A_{n}$, we have
\begin{equation*}
P(C_{n})=P(A_{n})-P(A_{n}-C_{n})>P(A_{n})-\varepsilon /2>\varepsilon /2.
\end{equation*}

\bigskip \noindent We conclude that for all $n\geq 1$, $C_{n}$ is non-empty. So, by the axiom of Choice, we may choose, $n\geq 1$,
$x^{(n)}=(x_1^{(n)}, \ x_2^{(n)}, ....) \in E_{T_0}$ such that $x^{(n)}$. By the non-decreasingness of the sequence $(C_n)_{n\geq 1}$, the sequence $(x^{(n)})_{n\geq 1}$ is in $C_1$, which we recall is such that

$$
C_{1}\subset A_{1}^{\prime }=B_{V_{1}}^{\prime }\times E^{\prime}_{T_0 \setminus V_{1}}
$$

\begin{equation*}
\forall n\geq 1, \ (x_{1}^{(n)},x_{2}^{(n)},...,x_{m_{1}}^{(n)})) \in E_{V_{1}}^{\prime }.
\end{equation*}

\bigskip \noindent Since $B^{\prime}_{V_{1}}$ est compact, there exists a sub-sequence $(x^{(n_{1,k})})$ of $(x^{(n)})$ such that
$$
(x_{1}^{(n_{1,k})},x_{2}^{(n_{1,k})},...,x_{m_{1}}^{(n_{1,k})})) \rightarrow (x^{\ast}_{1}, ..., x^{\ast}_{m(1)}) \in B^{\prime}_{V_{1}}\subset B_{V_{1}}.
$$

\bigskip \noindent But the sub-sequence Now by the nature $(x^{(n_{1,k})}) \in C_2$ whenever $n_{1,k}\geq 2$ (which happens for from some value 
$k_2>0$ since the sequence $(n_{1,k})_{k\geq 1}$ is an increasing sequence of non-negative integers). We thus have
\begin{equation*}
\forall k>k_2\geq 1, \ (x_{1}^{(n_{1,k})},x_{2}^{(n_{1,k})},...,x_{m_{1}}^{(n_{1,k})})) \in E_{V_{2}}^{\prime }.
\end{equation*}

\bigskip \noindent We conclude similarly that there exists a sub-sequence $(x^{(n_{2,k})})_{k\geq 1}$ of $(x^{(n_{1,k})})_{k\geq 1}$ such that

$$
(x_{1}^{(n_{2,k})},x_{2}^{(n_{2,k})},...,x_{m_{1}}^{(n_{2,k})})) \rightarrow (x^{\ast}_{1}, ..., x^{\ast}_{m(2)}) \in  B_{V_{2}}.
$$

\bigskip \noindent It is important the for a common factor $j$  between $A_{V_1}$ and $A_{V_2}$, the limit $x^{\ast}_{j}$ remains unchanged as the limit of a sub-sequence of converging sequence. We may go so-forth and consider the diagonal sub-sequence

$$
(x^{(n_{k,k})})_{k\geq 1}.
$$

\bigskip \noindent We have that for each $n\geq 1$, there exists $K(n)>0$ such that $(x^{(n_{k,k})})_{k\geq K(n)} \subset C_n$. Hence, for each 
$n\geq 1$,

$$
(x^{\ast}_{1}, ..., x^{\ast}_{m(n)}) \in  B_{V_{n}}.
$$

\bigskip \noindent So by denoting

$$
x^{\ast}=(x^{\ast}_{1}, x^{\ast}_{2}, ... ),
$$

\bigskip \noindent we get that $x^{\ast}$ belongs to each $A_n$, $n\geq 0$. Hence $A$ is not empty.

%% file: proba_02_appendix.tex
\chapter{Appendix} \label{proba_02_appendix}

\section{Some Elements of Topology}. \label{proba_02_appendix_topology}

\bigskip
\noindent \textbf{I - Stone-Weierstrass Theorem}.\\

\noindent Here are two forms of Stone-Weierstrass Theorem. The second is more general and is the one we use in this text.

\begin{proposition} \label{proba02_rv_append_sw_prop1}
Let $(S, d)$ be a compact metric space and $H$ a non-void subclass of the class $\mathcal{C}(S, \mathbb{R})$ of all real-valued continuous functions defined on $S$. Suppose that $H$ satisfies the following conditions.\\

\bigskip \noindent (i) $H$ is \textit{lattice}, that is, for any couple $(f,g)$ of elements of $H$,  $f\wedge g$\ et $f\vee g$ are in $H$
\newline

\bigskip \noindent (ii) For any couple $(x,y)$ of elements of $S$ and for any couple $(a,b)$ of real numbers such that $a=b$ if $x=y$, there exists a couple $(h,k)$ of elements of $H$ such that
\begin{equation*}
h(x)=a\text{ and }k(y)=b.
\end{equation*}

\bigskip \noindent Then $H$ is dense in $\mathcal{C}(S, \mathbb{R})$ endowed with the uniform topology, that is each continuous
function  from $S$ to $\mathbb{R}$ is the uniform limit of a sequence of elements in $H$.
\end{proposition}

\begin{theorem} \label{proba02_rv_append_sw_prop2}
Let $(S, d)$ be a compact metric space and $H$ a non-void subclass of the class $\mathcal{C}(S, \mathbb{C})$ of all real-valued continuous functions defined on $S$. Suppose that $H$ satisfies the following conditions.\\

\noindent (i) $H$ contains all the constant functions.\\

\noindent (ii) For all $(h,k)\in H^{2}$, $h+k\in H,h\times k\in H,\overline{u}\in H$.\\

\noindent (iii) $H$ separates the points of $S$, i.e., for two distinct elements of $S$, $x$ and $y$, that is $x\neq y$, there exists $h\in H$ such that 
\begin{equation*}
h(x)\neq h(y).
\end{equation*}

\noindent Then $H$ is dense in $\mathcal{C}(S, \mathbb{C})$ endowed with the uniform topology, that is each continuous
function  from $S$ to $\mathbb{C}$ is the uniform limit of a sequence of elements in $H$.
\end{theorem}

\bigskip \noindent \textbf{Remark}.\newline

\noindent If we work in $\mathbb{R}$, the condition on the conjugates - $\overline{u}\in H$\ - becomes needless.\\

\noindent But here, these two classical versions do not apply. We use the following extension.

\begin{corollary} \label{sec_EF_cor_05} Let $K$ be a non-singleton compact space and $\mathcal{A}$ be a non-empty sub-algebra of $C(K,\mathbb{C})$. Let $f \in C(K,\mathbb{C})$. Suppose that there exists $K_0 \subset K$ such that $k\setminus K_0$ has at least two elements and $f$ is constant on $K_0$. Suppose that the following assumption hold.

\noindent (1) $\mathcal{A}$ separates the points of $K\setminus K_0$ and separates any point of $K_0$ from any point of $K\setminus K_0$.\\

\noindent (2) $A$ contains all the constant functions.\\

\noindent (3) For all $f \in \mathcal{A}$, its conjugate function $\bar{f}=\mathcal{R}(f) - i \mathcal{Im}(f) \in \mathcal{A}$,\\

\noindent  Then 
$$
f \in \overline{\mathcal{A}}.
$$
\end{corollary}

\noindent A proof if it available in \cite{loSW2018}.\\

\noindent \textbf{II- Approximations of indicator functions of open sets by Lipschitz function}.\\

\noindent We have the

\begin{lemma} \label{proba02_rv_charac_opensets} Let $(S,d)$ be an arbitrary metric space and $G$ be an open set of $S$. Then there exists a non-decreasing sequence $(f_{k})_{k\geq 1}$ of non-negative real-valued and Lipschitz functions defined on $S$ converging to $1_G$.\\
\end{lemma}

\bigskip \noindent \textbf{Proof}. Let $G$ be an open set of $S$. For any integer number $k\geq 1$, set the function $f_{k}(x)=\min
(kd(x,G^{c}),1)$, $x\in S$. We may see that for any $k\geq 1$, $f_{m}$ has values in $[0,1]$, and is bounded. Since $G^{c}$ is closed, we have
\begin{equation*}
d(x,G^{c})=\left\{ 
\begin{array}{c}
>0\text{ if }x\in G \\ 
0\text{ if }x\in G^{c}\text{ }%
\end{array}%
\right. .
\end{equation*}

\bigskip \noindent Let us show that $f_{k}$ is a Lipschitz function. Let us handle $\left\vert f_{k}(x)-f_{k}(y)\right\vert $ through
three cases.\\

\noindent Case 1. $(x,y)\in (G^{c})^{2}$. Then 
\begin{equation*}
\left\vert f_{k}(x)-f_{k}(y)\right\vert =0\leq k\text{ }d(x,y).
\end{equation*}

\bigskip \noindent Case 2. $x\in G$\ and  $y\in G^{c}$\ (including also the case where the roles of $x$ and $y$ are switched). We have 
\begin{equation*}
\left\vert f_{k}(x)-f_{k}(y)\right\vert =\left\vert \min
(kd(x,G^{c}),1)\right\vert \leq k\text{ }d(x,G^{c})\leq k\text{ }d(x,y),
\end{equation*}

\bigskip \noindent by the very definition of $d(x,G^{c})=\inf \{d(x,z),$\ z$\in G^{c}\}.$\newline

\noindent Case 3. $(x,y) \in G^{2}$. We use Lemma \ref{proba02_rv_append_annexe2} in this section, to get 
\begin{equation*}
\left\vert f_{k}(x)-f_{k}(y)\right\vert =\left\vert \min
(kd(x,G^{c}),1)-\min (kd(y,G^{c}),1)\right\vert \leq \left\vert
kd(x,G^{c})-kd(y,G^{c})\right\vert ,
\end{equation*}%
\begin{equation*}
\leq k\text{ }d(x,y)
\end{equation*}

\bigskip \noindent by the second triangle inequality. Then  $f_{k}$ is a Lipschitz function with coefficient $k$. Now, let us show that
\begin{equation*}
f_{k}\uparrow 1_{G}\text{ as k}\uparrow \infty .
\end{equation*}

\bigskip \noindent Indeed, if $x\in G^{c}$, we obviously have $f_{k}(x)=0\uparrow
0=1_{G}(x)$. If $x\in G$, that $d(x,G^{c})>0$ and $kd(x,G^{c})\uparrow \infty$ as $k\uparrow \infty$. Then for $k$ large enough, 
\begin{equation}
f_{k}(x)=1\uparrow 1_{G}(x)=1 \text{ as} k\uparrow \infty.  \label{limfm}
\end{equation}

\bigskip \noindent $\blacksquare$.\\

\noindent \textbf{III -Lipschitz property of finite maximum or minimum}.\\

\noindent We have the

\begin{lemma} \label{proba02_rv_append_annexe2} For any real numbers $x$, $y$, $X$, and $Y$, 
\begin{equation}
\left\vert \min (x,y)-\min (X,Y)\right\vert \leq \left\vert x-X\right\vert
+\left\vert y-Y\right\vert.  \label{proba02_rv_append_annexe2_for}
\end{equation}
\end{lemma}

\bigskip \noindent \textbf{Proof}. Let us have a look at the four possibles case.\\

\noindent Case 1 : $min(x,y)=x$ and $min(X,Y)=X$. We have
\begin{equation*}
\left\vert \min (x,y)-\min (X,Y)\right\vert \leq \left\vert x-X\right\vert
\end{equation*}

\noindent Case 2 : $min(x,y)=x$ and $min(X,Y)=Y$. If $x\leq Y$, we have $Y\geq X$, we have
\begin{equation*}
0 \leq \min (X,Y) - \min (x,y)=Y-x \leq X-x
\end{equation*}

\noindent If $x>Y$, we have $X\geq Y$, we have
\begin{equation*}
0 \leq \min (x,y)-\min (X,Y)=x-Y\leq y-Y
\end{equation*}

\noindent Case 3 : $min(x,y)=y$ and $min(X,Y)=Y$. We have
\begin{equation*}
\left\vert \min (x,y)-\min (X,Y)\right\vert \leq \left\vert y-Y\right\vert
\end{equation*}

\noindent Case 4 : $min(x,y)=y$ and $min(X,Y)=X$. This case id handled as for Case 2 by permuting the roles of $(x,y)$ and $(X,Y)$.\\

\newpage
\section[Orthogonal matrix]{Orthogonal Matrices, Diagonalization of Real Symmetrical Matrices and Quadratic forms}\label{proba_02_appendix_ortho}

\bigskip
\noindent \textbf{I - Orthogonal matrices}.\\

\noindent We begin by this result.

\begin{proposition} \label{proba_02_appen_prop_01} For any square $d$-matrix $T$, we have the equivalence between the following assertions.\\

\noindent (1) $T$ is invertible and the inverse matrix $T^{-1}$ of $T$ is its transpose matrix, that is

$$
TT^t=T^t T=I_{d},
$$

\bigskip \noindent where $I_{d}$ is the identical matrix of dimension $d$.\\

\noindent (2) $T$ is an isometry, that is $T$ preserves the norm : For any $x\in \mathbb{R}^d$

$$
\left\Vert Tx\right\Vert =\left\Vert x\right\Vert.
$$

\bigskip \noindent (3) The columns $\biggr( T^{(1) },T^{(2)},..,T^{(d)} \biggr)$ form an orthonormal basis of $\mathbb{R}^{d}$.\\

\bigskip \noindent (4) The transposes of the lines $\biggr(T_{1}^t,...,T_{d}^t \biggr)$ form an orthonormal basis of $\mathbb{R}^{d}$.\\

\noindent Besides, if $T$ is orthogonal, its transpose is also orthogonal and satisfies
$$
det(T)=\pm 1.
$$
\end{proposition}

\bigskip \noindent Before we give the proof, we provide the definition of an orthogonal matrix.\\

\begin{definition} $ $\\
\noindent A square $d$-matrix is orthogonal if and only if one of the equivalent assertions of Proposition \ref{proba_02_appen_prop_01} holds.
\end{definition}

\bigskip \noindent Now we may concentrate of the \\

\noindent \textbf{Proof of Proposition \ref{proba_02_appen_prop_01}}.\\

\noindent Recall that, in finite dimension linear theory, the $d$-matrix $B$ is the inverse of the $d$-matrix $A$ if and only if $AB=I_d$ if and only if $BA=I_d$. (See the reminder at the end of the proof).\\

\noindent Let us show the following implications or equivalences.\\

\noindent (i)  $(1) \Leftrightarrow (3)$. By definition, for any  $(i,j) \in \{1,...,d\}^d$,

$$
\left(T^t T\right) _{ij}=\left( T^t\right) _{i}T^{\left( j\right)}=\left(T^{\left( j\right)}\right)^t T^{\left( j\right)}=\langle T^{\left( i\right)},T^{\left( j\right)}\rangle. \ (I01)
$$

\noindent and
$$
\left(TT^t\right)_{ij}=T_i (T^t)_j=T_i (T_j)^t=\biggr(T_i^t\biggr)^t (T_j)^t=\langle T_i^t, T_j^t\rangle. \ (I02)
$$

\bigskip \noindent By Formula (I01), we have the equivalence between $T^t T=I_d$ and (3), and thus, (1) and (3) are equivalent.\\

\noindent (i) : $(1) \Leftrightarrow (4)$. The same conclusion is immediate by using Formuka (I02) instead of Formula (I01).\\

\noindent (ii) $(3) \Longleftrightarrow (2)$. We have for all $x\in \mathbb{R}^{d}$,
 
\begin{eqnarray*}
\left\Vert Tx\right\Vert ^{2} &=&\langle Tx,Tx\rangle =^{t}\left( Tx\right)
\left( Tx\right) =^{t}x^{t}TTx \\
&=&\sum_{i=1}^{n}\sum_{j=1}^{n}\left( ^{t}TT\right)
_{ij}x_{i}x_{j}=\sum_{i=1}^{n}\sum_{j=1}^{n}\langle T^{\left( i\right)
},T^{\left( j\right) }\rangle x_{i}x_{j}. \ (IS01)
\end{eqnarray*}

\bigskip \noindent Hence, (3) implies that for all $x\in \mathbb{R}^{d}$,

\begin{equation*}
\ \left\Vert Tx\right\Vert ^{2}=\sum_{i=1}^{n}x_{i}^{2}=\ \left\Vert
x\right\Vert ^{2},
\end{equation*}

\bigskip \noindent which is the definition of an isometry.\\

\noindent (iii) $(2) \Longleftrightarrow (3)$. Let us suppose that (2) holds.\\

\noindent To show that each $T^{(i_0)}$, for a fixed  $i_0\in \{1,...,d\}$ is normed, we apply (ISO) to the  vector $x$ whose coordinates are zero except $x_{i_0}=1$. We surely have $\|x\|=1$ and all the terms of

$$
\sum_{i=1}^{n}\sum_{j=1}^{n}\langle T^{\left( i\right)},T^{\left( j\right) }\rangle x_{i}x_{j}
$$

\bigskip \noindent are zero except for $i=j=i_0$, and the summation reduces to  $\langle T^{\left( i_0\right)},T^{\left( i_0\right) }\rangle x_{i_0}^2=\langle T^{\left( i_0\right)},T^{\left( i_0\right) }\rangle$. Equating the summation with  $\|x\|^2$ gives that 
$$
\langle T^{\left( i_0\right)},T^{\left( i_0\right) }\rangle=1.
$$

\bigskip \noindent So,  $T^{(i_0)}$ is normed.\\

\noindent To show that two different columns $T^{(i_0)}$ and $T^{(j_0)}$, for a fixed ordered pair $(i_0,j_0) \in \{1,...,d\}^2$, are orthogonal, we apply (ISO) the vector $x$ whose coordinates are zero except $x_{i_0}=1$ and $x_{j_0}=1$. We surely have $\|x\|^2=2$ and all the terms of

$$
\sum_{i=1}^{n}\sum_{j=1}^{n}\langle T^{\left( i\right)},T^{\left( j\right) }\rangle x_{i}x_{j}
$$

\bigskip \noindent are zero except for $i=j=i_0$, $i=j=j_0$ and $(i,j)=(i_0,j_0)$, and the summation reduces to  

$$
\langle T^{\left( i_0\right)},T^{\left( i_0\right) }\rangle x_{i_0}^2 + \langle T^{\left( j_0\right)},T^{\left( j_0\right) }\rangle x_{i_0}^2 + 2 \langle T^{\left( i_0\right)},T^{\left( j_0\right) }\rangle x_{i_0} x_{j_0}^2
$$

\bigskip 
\noindent By equating with the summation with $\|x\|^2$, we get

$$
2 = 2 + 2 \langle T^{\left( i_0\right)},T^{\left( j_0\right)}\rangle.
$$

\bigskip 
\noindent This implies that $\langle T^{\left( i_0\right)},T^{\left( j_0\right)}\rangle=0$.\\

\noindent We conclude that (2) holds whenever (3) does.\\

\noindent We obtained the following equivalences

\begin{table}[htbp]
	\centering
		\begin{tabular}{c c c}
		(1) & $\Leftrightarrow$  & (3) \\
		  $\Updownarrow$  &   &  $\Updownarrow$\\
		(4) &   & (2), \\
			
		\end{tabular}
\end{table}

\noindent from which we deive the equivalence between the four assertions. \\

\noindent It remains the two last points. That the transpose of $T$ is orthogonal with $T$, is a direct consequence of the equivalence between assertions (3) and (4). Since a square matrix and its transpose
have the same determinant and since $TT^t=I_d$, we get that $1=det(TT^t)=det(T)det(T^t)=det(T)^2$. $\blacksquare$.\\

\bigskip 
\noindent \textbf{A useful reminder}. \\

\noindent In finite dimension linear theory, the $d$-matrix $B$ is the inverse of the $d$-matrix $A$ if and only if $AB=I_d$ if and only if $BA=I_d$. But in an arbitrary algebraic structure $(E,\star)$ endowed with an internal operation $\star$ having a unit element $e$, that is an element of $e$ satisfying $x\star e=e \star x=x$ for all $x\in E$, an inverse $y$ of $x$ should should fulfills : 
$x\star y=y \star x=e$. The definition may be restricted to $x\star y=e$ or to $y \star x=e$ if the operation  $e$ is commutative.\\

\noindent In the case of $d$-matrices, the operation is not commutative. So using only one of the two conditions $AB=I_d$ and $BA=I_d$ to define the inverse of a matrix $A$ is an important result of linear algebra in finite dimensions.\\

\bigskip \noindent \textbf{II - Diagonalization of symmetrical matrices}.\\

\noindent \textbf{Statement and proof}.\\

\noindent We have the important of theorem.\\

\begin{theorem} For any real and symmetrical $d$-matrix $A$, there exists an orthogonal $d$-matrix $T$ such that $TAT^t$ is a diagonal matrix $diag(\delta_1, ..., \delta_d)$, that is

$$
TAT^t =diag(\delta_1, ..., \delta_d),
$$ 

\bigskip \noindent where $\delta_i$, $1\leq i \leq d$, are finite real numbers. 
\end{theorem}

\bigskip \noindent \textbf{Remark}. In other words, any real and symmetrical $d$-matrix $A$ admits $d$ real eigen-values (not necessarily distinct) $\delta_i$, $1\leq i \leq d$ and the passage matrix may be chosen to be an orthogonal matrix.\\

\noindent \textbf{Proof}. \noindent Let us suppose that $A$ is symmetrical, which means that for any $u \in \mathbb{R}^d$, we have

$$
\langle Au,v \rangle=\langle u, Av \rangle. \ (S)
$$.

\noindent In a first step, let us borrow tools from Analysis. The linear application
$\mathbb{R}^d \ni u \rightarrow Au$ is continuous so that

\begin{equation*}
\left\Vert A\right\Vert =\sup_{u \in \mathbb{R}^d, \ \left\Vert u\right\Vert \leq 1} \left\Vert Au\right\Vert
=\sup_{u \in \mathbb{R}^d, \ \left\Vert X\right\Vert u 1}\left\Vert
Au\right\Vert <+\infty.
\end{equation*}

\bigskip \noindent Since the closed ball is closed a compact set in $\mathbb{R}^d$, there exists, at least, $u_{0}$ such that
$\left\Vert u_{0}\right\Vert =1$ and
 
\begin{equation*}
\sup_{\left\Vert u\right\Vert \leq 1}\left\Vert Au\right\Vert =\left\Vert
Au_{0}\right\Vert.
\end{equation*}

\bigskip \noindent In a second step, let us assume that $A$ has two eigen-vectors $u$ and $v$ associated to two distinct real eigen-valeues $\mu$ and $\lambda$. By Formula (S) above, we have

\begin{eqnarray*}
&&\mu \langle u,v \rangle=\lambda \langle u, v \rangle\\
&\Rightarrow& (\mu-\lambda) \langle u, v \rangle=0.
\end{eqnarray*}

\bigskip \noindent We get that $u$ and $v$ are orthogonal. We get the rule : two eigen-vectors of a symmetric square matrix which are associated to two distinct real eigen-values are orthogonal.\\

\bigskip \noindent In a third step, let us show that if a linear sub-space $F$ in invariant by $A$, that is for all $u\in F$, $Au \in F$ 
(denoted $AF\subset F$), then the orthogonal $F^{\perp}$ of $F$ is also invariant by $A$. Indeed, if $F$ is $A$-invariant and $v \in F^{\perp}$,
we have

$$
£\forall u \in F, \ \langle Av, u \rangle=\langle v, Au \rangle=0;
$$

\bigskip \noindent since $Au \in F$.\\

\bigskip \noindent Finally, in the last and fourth step, we have for $u \in \mathbb{R}^d$, such that $\|Au\|=1$, by applying the Cauchy-Schwartz Inequality

$$
\|Au\|=^2 =\langle Au, Au \rangle=\langle u, A^{2}u \rangle\leq \|A^{2}u\|. \ (S1)
$$

\noindent The equality is reached for some $u_1$ (with $\|u_1\|=1$) only if  $u_1$ and $A^{2}u_1$ are linearly dependent, that is exists $\lambda$ such that

$$
A^{2}u_1=\lambda u_1,
$$

\bigskip \noindent meaning that $u_1$ is an eigen-vector of $A^2$ where, by taking the norms, we have 

$$
\lambda=\|A^{2}u_1\|.
$$

\noindent Now, we have all the tools to solve the problem by induction. By definition, we have

$$
\|A\|^2=\|Au_0\|^2 \leq \|A^{2}u_0\| \leq \|A\| \|Au_0\|\leq \|A\| \|A\| \|u_0\|,
$$

\bigskip \noindent and hence Formula (S1) becomes an equality for $u_0$. The conclusion of the fourth step says that $u_0$ is an eigen-vector of $A^2$ associated to $\lambda=\|A\|^2$. For $\lambda=\mu^2$, this leads to $A^2u_0=u^2 u_0$, that is

$$
(A-\mu I_d)(A+\mu I_d)u_0=0.
$$

\noindent Now, either $(A+\mu I_d)u_0=0$ and $u_0$ is an eigen-vector of $A$ associated to $-\mu$, or $v_0=(A+\mu I_d)u_0\neq 0$ and $v_0$ is an eigen-vector of $A$ associated to $\mu$. In both case, the eigen-value is $\pm \|A\|$.\\

\noindent We proved that $A$ has at least on real eigen-vector we denote by $e_1$ associated to $\lambda_1=\pm \|A\|$. In a next step, let us denote $F_{1}=Lin(\{e_{1}\})$ and $G_{2}=F_{1}^{\perp}$. It is clear that $F_1$ is invariant by $A$, so is $G_{2}$. We consider the restriction of $A$ on $G_{2}$. We also have that $A_{2}$ symmetrical and clearly  $\left\Vert A_{2}\right\Vert
\leq \left\Vert A\right\Vert $. We  find an eigen-vector $e_2$ of $A_2$, thus of $A$, associated with $\lambda_2=\pm \|A_2\|$ and 
$|\lambda_1|\geq |\lambda_2|$ and $e_1$ and $e_2$ are orthogonal. We do the same for $F_2=Lin({e_1,e_2})$, $G_3=F_2^{\perp}$ and $A_3$ the restriction of $A_2$ (and hence of $A$) on $G_3$. We will find an eigen-vector $e_3$ of $A_2$, thus of $A$, associated with $\lambda_3=\pm \|A_3\|$ and $|\lambda_1|\geq |\lambda_2|\geq |\lambda_3|$ with  $\{e_1, e_2,e_3\}$ orthonormal. We proceed similarly to get exactly $d$ normed eigen-vectors orthogonal associated to a decreasing sequence of eigen-values in absolute values. $\blacksquare$\\

\noindent \textbf{(b) Some consequences}.\\

\noindent \textbf{(b1) Determinant}.\\

\noindent If $TAT^t =diag(\delta_1, ..., \delta_d)$, where $T$ is orthogonal, we have

$$
det(TAT^t)=det(T) det(T^t) det(A)=det(T)^2 det(A)=det(A)
$$ 

\bigskip \noindent and next
 
$$
det(A)=det\bigg(diag(\delta_1, ..., \delta_d)\biggr)
$$

\bigskip \noindent which leads to

$$
det(A)=\prod_{1}^{d} \delta_j.
$$

\bigskip \noindent \textbf{(b2) A useful identity}.\\

\noindent If $TAT^t =diag(\delta_1, ..., \delta_d)$, where $T$ is orthogonal, we have, for any $(i,j) \in \{1,...,d\}^2$, 

$$
\sum_{j=1}^{d} \delta_h \biggr(T^{(h)}(T^{(h)})^t\biggr)_{ij}=a_{ij}. \ (UID)
$$

\bigskip \noindent \textbf{Proof}. Let us denote $D=diag(\delta_1, ..., \delta_d)$ and suppose that $TAT^t=D$, $T$ being orthogonal. We get A=$T^tDT$. Hence, for any $(i,j) \in \{1,...,d\}^2$, we have

$$
(A)_{ij}=\biggr(T^tD\biggr)_{i}T^{(j)}.
$$

\bigskip \noindent But the $h$ elements, $1\leq h \leq d$, of the line $\biggr(T^tA\biggr)_{i}$ are $\biggr(T^t\biggr)_{i}D^{(h))}$, which are 

$$
\sum_{1\leq r \leq d} t_{ti} \delta_h \delta_{rh}=\delta_h t_{hi}.
$$
 
\bigskip \noindent Thus, we have

$$
(A)_{ij}=\biggr(T^tD\biggr)_{i}T^{(j)}=\sum_{1\leq h \leq d} \delta_h t_{hi} t_{hj}=\sum_{1\leq h \leq d} \delta_h  \biggr(T^{(h)} (T^{(h)})T\biggr)_{ij}.
$$

\bigskip \noindent \textbf{III - Elements from Bi-linear Forms and Quadratic Forms Theory}.\\

\noindent Before we begin, let us remind Formula (ACBT), seen in the proof of (P5) in Points (b)-(b2) in Section \ref{proba_02_rv_sec_04} in Chapter \ref{proba_02_rv} : for any $(p\times d)$-matrix $A$, any $(d \times s)$-matrix $C$ and any $(q \times s)$-matrix $B$, the $ij$-element of $ACB^t$ is given by

$$
\sum_{1\leq k\ \leq s} \sum_{1\leq p\ \leq p} a_{ih} c_{hk} b_{jk}. \ \textbf{}
$$
 
\bigskip \noindent Let us apply this to vectors $u=A^t \in \mathbb{R}^d$, $v=B^t \in \mathbb{R}^k$ and to a matrix $(d \times k)$-matrix $C$. The unique element of the $(1 \times 1)$-matrix  $u^t C v$ is

$$
u^t C v = \sum_{1\leq i\ \leq d} \sum_{1\leq j \leq k} c_{ij} u_{i} v_{j}.
$$

\bigskip \noindent This formula plays a  key role in bi-linear forms studies in finite dimensions.\\

\noindent If $d$-matrix $C$ is diagonal, that is $C_{ij}=0$ for $i\neq j$, we have

$$
u^t C v = \sum_{1\leq i\ \leq d} \delta_{j}  v_{j}^2. \ \ (C0)
$$

\bigskip \noindent where $\delta_j=c_{jj}$, $j \in \{1,...,d\}$. We may use the Kronecker's symbol defined by

\begin{equation*}
\delta_{ij}=\left\{ 
\begin{tabular}{lll}
$1$ & if & $i=j$ \\ 
$0$ & if  & $i\neq j$
\end{tabular}
\right. .
\end{equation*}

\bigskip \noindent to get the following notation of a diagonal matrix. A $d$-diagonal matrix $D$ whose diagonal elements are denoted by $\delta_j$, $j \in \{1,...,d\}$, respectively, may be written as follows :

$$
D=diag(\delta_1,...,\delta_2)=\biggr(\delta_i \delta_{ij}\biggr)_{1\leq i,j \leq d}.
$$

\bigskip \noindent \textbf{(a) - Bi-linear Forms}.\\

\noindent By definition, a function 

\begin{equation*}
f \ : \ \mathbb{R}^{d} \times \mathbb{R}^{k}\rightarrow \mathbb{R}
\end{equation*}

\bigskip \noindent is bi-linear if and only if :\\

\noindent (i) for any fixed $u \in \mathbb{R}^{d}$, the partial application $v \mapsto f\left( u,v\right) $ is linear\\

\noindent and \\

\noindent (ii) for any fixed $v \in \mathbb{R}^{k}$, the partial application $u \mapsto f\left( u,v\right) $ is linear.\\

\noindent The link with matrices theory is the following. Let $\left( e_{1},e_{2},...e_{n}\right)$ be an orthonormal basis of $\mathbb{R}^{d}$ and
$\left( \varepsilon _{1},\varepsilon _{2},...,\varepsilon _{k}\right)$ an orthonormal basis of $\mathbb{R}^{k}$. Let us define the $(d \times k)$-matrix $A$ by

$$
a_{ij}=f(e_{i},\varepsilon _{j})
$$

\bigskip \noindent and denote the coordinates of $u\in \mathbb{R}^d$ and $v\in \mathbb{R}^k$ in those bases by 
\begin{equation*}
u=\sum_{i=1}^{n}u_{i}e_{i} \text{ and } v=\sum_{j=1}^{m}v_{j}\varepsilon _{j}.
\end{equation*}

\noindent We have the following expression of the bi-linear form

$$
f(u,v)= u^t A v.
$$

\bigskip \noindent The proof is the following :

\begin{eqnarray*}
f(u,v) &=&f\left( \sum_{i=1}^{n}u_{i}e_{i},\sum_{j=1}^{m}v_{j}\varepsilon
_{j}\right) \\
&=&\sum_{i=1}^{n}u_{i}f\left( e_{i},\sum_{j=1}^{m}v_{j}\varepsilon
_{j}\right) \\
&=&\sum_{i=1}^{n}u_{i}\sum_{j=1}^{m}v_{j}f\left( e_{i},\varepsilon
_{j}\right) \\
&=&\sum_{i=1}^{n}\sum_{j=1}^{m}u_{i}a_{ij}v_{j}
\end{eqnarray*}

\bigskip \noindent Thus, we may conclude with the help of Formula  (uTCv) above.\\

\bigskip \noindent \textbf{(b) - Quadratic forms}.\\

\noindent For any bi-linear form $f \ : \ (\mathbb{R}^{d})^2 \rightarrow \mathbb{R}$, the mapping

$$
x \ni \mathbb{R}^{d} \mapsto Q_f(x)=f(x,x) 
$$

\bigskip \noindent is called the quadratic form associated with $f$.\\

\noindent The quadratic form is said to be \textbf{semi-positive} if and only if $Q_f(x)\geq 0$, for all $x \in \mathbb{R}^{d}$.\\

\noindent It is said to be \textbf{positive} if and only if  $Q_f(x)>0 x$, for all $0 \neq x \in \mathbb{R}^{d}$.\\

\noindent We already know that $f$ may be represented by a $d$-matrix $A$ and thus, $Q_f$ may be represented as

$$
Q_f(u)=u^t A u, \ u \in \mathbb{R}^{d}.
$$

\bigskip \noindent But, since $Q_f(u)^t=u^t A^t u=Q_f(u)$, we also have for all $u \in \mathbb{R}^{d}$ that

$$
Q_f(u)=u^t \frac{A+A^t}{2} u, \ u \in \mathbb{R}^{d}.
$$

\bigskip \noindent The matrix $B=(A+A^t)/1$ is symmetrical and we have

$$
Q_f(u)=u^t B u, \ u \in \mathbb{R}^{d},
$$

\bigskip \noindent which leads to the : \\

\begin{proposition} Any quadratic form $Q$ on $\mathbb{R}^{d}$ is of the form.

$$
Q(u)=u^t B u, \ u \in \mathbb{R}^{d},
$$

\bigskip \noindent where $B$ is a $d$-symmetrical form.
\end{proposition}

\bigskip \noindent \textbf{(c) - Canonical reduction of a Quadratic forms}.\\

\noindent Reducing a quadratic form $Q$ on $\mathbb{R}^{d}$ to a canonical form consists in finding an invertible linear change of variable $v=Tu=(v_1,...,v_d)^d$ such that $Q_0(v)=Q(Tv)$ is of the form

$$
Q_0(v)=\sum_{1\leq j \leq d} \delta_j v_j^2.
$$

\bigskip \noindent This may be achieved in finite dimension in the following ways. Let $B$ be a symmetrical matrix associated to the quadratic form $B$. According to Part II of this section, we can find an orthogonal matrix $T$ such that $T B T^t$ is a diagonal matrix $D=diag(\delta_1, ..., \delta_d)$. For $v=Tu=(v_1,...,v_d)^d$, we have

$$
Q(vT)=v^t (T^t AT) v=v^t D v=\sum_{1\leq i \leq d, \ 1\leq j \leq d} d_{ij} v_i v_j=\sum_{1\leq j \leq d} \delta_{j} v_j^2.
$$

\bigskip \noindent This leads to the

\begin{proposition} Any quadratic form $Q$ on $\mathbb{R}^{d}$ the form.

$$
Q(u)=u^t B u, \ u \in \mathbb{R}^{d},
$$

\bigskip \noindent where $B$ is a $d$-symmetrical matrix, may be reduced to the canonical for

$$
Q(u)=\sum_{1\leq j \leq d} \delta_j v_j^2,\ \ (CF)
$$

\bigskip \noindent where $v=Tu$ and the columns of $T^t$ form an orthonormal basis of $\mathbb{R}^d$ and are eigen-vector of $B$ respectively associated to the eigen-values $\delta_j$, $1\leq j\leq d$.\\
\end{proposition}

\bigskip \noindent \textbf{Consequences}. From the canonical form (CF), we may draw the straightforward following facts based on the facts that $T$ is invertible and its determinant is $\pm 1$. Hence each element $u \in \mathbb{R}^{d}$ is of the from $u=Tv$. Hence, Formula (FC) holds for all $u \in \mathbb{R}^{d}$ with $u=Tv$.\\ 

\noindent (1) If all the eigen-values are non-negative, the quadratic form $Q$ is semi-positive.\\

\noindent (2) If all the eigen-values are positive, the quadratic form $Q$ is semi-positive.\\

\noindent (3) If the quadratic form $Q$ is semi-positive and $B$ is invertible or has a non-zero determinant, that it is positive.\\

\bigskip \noindent Before we close the current section, let us remind that a canonical form  as in (CF) is not unique. But the three numbers of positive terms $\eta_{(+)}$, of negative terms
($\eta_{-}$) and zero terms ($\eta_{0}$) are unique. The triplet $(\eta_{(-)}, \eta_{(0)},\eta_{(+)})$ is called the signature of the quadratic form.\\

====

\section[Limits on $\mathbb{R}$]{What should not be ignored on limits in $\overline{\mathbb{R}}$ - Exercises with Solutions} \label{proba_02_appendix_limits}

\bigskip

\noindent \textbf{Definition} $\ell \in \overline{\mathbb{R}}$ is an accumulation point of a sequence 
 $(x_{n})_{n\geq 0}$ of real numbers finite or infinite, in $\overline{\mathbb{R}}$, if and only if there exists a sub-sequence $(x_{n(k)})_{k\geq 0}$ of
 $(x_{n})_{n\geq 0}$ such that $%
x_{n(k)}$ converges to $\ell $, as $k\rightarrow +\infty $.\newline

\noindent \textbf{Exercise 1}.\\

\noindent  Set $y_{n}=\inf_{p\geq n}x_{p}$ and $z_{n}=\sup_{p\geq n}x_{p} $ for all $n\geq 0$. Show that :\newline

\noindent \textbf{(1)} $\forall n\geq 0,y_{n}\leq x_{n}\leq z_{n}$.\newline

\noindent \textbf{(2)} Justify the existence of the limit of $y_{n}$ called limit inferior of the sequence $(x_{n})_{n\geq 0}$, denoted by $%
\liminf x_{n}$ or $\underline{\lim }$ $x_{n},$ and that it is equal to the following%
\begin{equation*}
\underline{\lim }\text{ }x_{n}=\lim \inf x_{n}=\sup_{n\geq 0}\inf_{p\geq
n}x_{p}.
\end{equation*}

\noindent \textbf{(3)} Justify the existence of the limit of $z_{n}$ called limit superior of the sequence $(x_{n})_{n\geq 0}$ denoted by $%
\lim \sup x_{n}$ or $\overline{\lim }$ $x_{n},$ and that it is equal%
\begin{equation*}
\overline{\lim }\text{ }x_{n}=\lim \sup x_{n}=\inf_{n\geq 0}\sup_{p\geq
n}x_{p}x_{p}.
\end{equation*}

\bigskip

\noindent \textbf{(4)} Establish that 
\begin{equation*}
-\liminf x_{n}=\limsup (-x_{n})\noindent \text{ \ \ and \ }-\limsup
x_{n}=\liminf (-x_{n}).
\end{equation*}

\newpage \noindent \textbf{(5)} Show that the limit superior is sub-additive and the limit inferior is super-additive, i.e. :  for two sequences
$(s_{n})_{n\geq 0}$ and $(t_{n})_{n\geq 0}$ 
\begin{equation*}
\limsup (s_{n}+t_{n})\leq \limsup s_{n}+\limsup t_{n}
\end{equation*}

\noindent and
\begin{equation*}
\lim \inf (s_{n}+t_{n})\geq \lim \inf s_{n}+\lim \inf t_{n}.
\end{equation*}

\noindent \textbf{(6)} Deduce from (1) that if%
\begin{equation*}
\lim \inf x_{n}=\lim \sup x_{n},
\end{equation*}%
then $(x_{n})_{n\geq 0}$ has a limit and 
\begin{equation*}
\lim x_{n}=\lim \inf x_{n}=\lim \sup x_{n}
\end{equation*}

\bigskip

\noindent \textbf{Exercise 2.} Accumulation points of $(x_{n})_{n\geq 0}$.\newline

\noindent \textbf{(a)} Show that if $\ell _{1}$=$\lim \inf x_{n}$ and $\ell
_{2}=\lim \sup x_{n}$ are accumulation points of $(x_{n})_{n\geq 0}.
$ Show one case and deduce the second one and by using Point (3) of Exercise 1.\newline

\noindent \textbf{(b)} Show that $\ell _{1}$ is the smallest accumulation point of $(x_{n})_{n\geq 0}$ and $\ell _{2}$ is the biggest.
(Similarly, show one case and deduce the second one and by using Point (3) of Exercise 1).\newline

\noindent \textbf{(c)} Deduce from (a) that if $(x_{n})_{n\geq 0}$ has a limit $\ell$,  then it is equal to the unique accumulation point and so,%
\begin{equation*}
\ell =\overline{\lim }\text{ }x_{n}=\lim \sup x_{n}=\inf_{n\geq
0}\sup_{p\geq n}x_{p}.
\end{equation*}

\noindent \textbf{(d)} Combine this result with Point \textbf{(6)} of Exercise 1 to show that a sequence $(x_{n})_{n\geq 0}$ of $\overline{\mathbb{R}}
$ has a limit $\ell $ in $\overline{\mathbb{R}}$ if and only if\ $\lim \inf
x_{n}=\lim \sup x_{n}$ and then%
\begin{equation*}
\ell =\lim x_{n}=\lim \inf x_{n}=\lim \sup x_{n}.
\end{equation*}

\newpage

\noindent \textbf{Exercise 3. } Let $(x_{n})_{n\geq 0}$ be a non-decreasing sequence
of $\overline{\mathbb{R}}$. Study its limit superior and its limit inferior and deduce that%
\begin{equation*}
\lim x_{n}=\sup_{n\geq 0}x_{n}.
\end{equation*}

\noindent Deduce that for a non-increasing sequence $(x_{n})_{n\geq 0}$
of $\overline{\mathbb{R}},$%
\begin{equation*}
\lim x_{n}=\inf_{n\geq 0}x_{n}.
\end{equation*}

\bigskip

\noindent \textbf{Exercise 4.} (Convergence criteria)\newline

\noindent \textbf{Prohorov Criterion} Let $(x_{n})_{n\geq 0}$ be a sequence of $\overline{%
\mathbb{R}}$ and a real number $\ell \in \overline{\mathbb{R}}$ such that: Every subsequence of $(x_{n})_{n\geq 0}$ 
also has a subsequence ( that is a subssubsequence of $(x_{n})_{n\geq 0}$ ) that converges to $\ell .$
Then, the limit of $(x_{n})_{n\geq 0}$ exists and is equal $\ell .$\newline

\noindent \textbf{Upcrossing or Downcrossing Criterion}. \newline

\noindent Let $(x_{n})_{n\geq 0}$ be a sequence in $\overline{\mathbb{R}}$ and two real numbers $a$ and $b$ such that $a<b.$
We define%
\begin{equation*}
\nu _{1}=\left\{ 
\begin{array}{cc}
\inf  & \{n\geq 0,x_{n}<a\} \\ 
+\infty  & \text{if (}\forall n\geq 0,x_{n}\geq a\text{)}%
\end{array}%
\right. .
\end{equation*}%
If $\nu _{1}$ is finite, let%
\begin{equation*}
\nu _{2}=\left\{ 
\begin{array}{cc}
\inf  & \{n>\nu _{1},x_{n}>b\} \\ 
+\infty  & \text{if (}n>\nu _{1},x_{n}\leq b\text{)}%
\end{array}%
\right. .
\end{equation*}%
.

\noindent As long as the $\nu _{j}'s$ are finite, we can define for $\nu
_{2k-2}(k\geq 2)$

\begin{equation*}
\nu _{2k-1}=\left\{ 
\begin{array}{cc}
\inf  & \{n>\nu _{2k-2},x_{n}<a\} \\ 
+\infty  & \text{if (}\forall n>\nu _{2k-2},x_{n}\geq a\text{)}%
\end{array}%
\right. .
\end{equation*}%
and for $\nu _{2k-1}$ finite, 
\begin{equation*}
\nu _{2k}=\left\{ 
\begin{array}{cc}
\inf  & \{n>\nu _{2k-1},x_{n}>b\} \\ 
+\infty  & \text{if (}n>\nu _{2k-1},x_{n}\leq b\text{)}%
\end{array}%
\right. .
\end{equation*}

\noindent We stop once one $\nu _{j}$ is $+\infty$. If $\nu
_{2j}$ is finite, then 
\begin{equation*}
x_{\nu _{2j}}-x_{\nu _{2j-1}}>b-a. 
\end{equation*}

\noindent We then say : by that moving from $x_{\nu _{2j-1}}$ to $x_{\nu
_{2j}},$ we have accomplished a crossing (toward the up) of the segment $[a,b]$
called \textit{up-crossings}. Similarly, if one $\nu _{2j+1}$
is finite, then the segment $[x_{\nu _{2j}},x_{\nu _{2j+1}}]$ is a crossing downward (down-crossing) of the segment $[a,b].$ Let%
\begin{equation*}
D(a,b)=\text{ number of up-crossings of the sequence of the segment }[a,b]\text{.}
\end{equation*}

\bigskip

\noindent \textbf{(a)} What is the value of $D(a,b)$ if \ $\nu _{2k}$ is finite and $\nu
_{2k+1}$ infinite.\newline

\noindent \textbf{(b)} What is the value of $D(a,b)$ if \ $\nu _{2k+1}$ is finite and $\nu
_{2k+2}$ infinite.\newline

\noindent \textbf{(c)} What is the value of $D(a,b)$ if \ all the $\nu _{j}'s$ are finite.%
\newline

\noindent \textbf{(d)} Show that $(x_{n})_{n\geq 0}$ has a limit iff
for all $a<b,$ $D(a,b)<\infty.$\newline

\noindent \textbf{(e)} Show that $(x_{n})_{n\geq 0}$ has a limit iff
for all $a<b,$ $(a,b)\in \mathbb{Q}^{2},D(a,b)<\infty .$\newline

\bigskip

\noindent \textbf{Exercise 5. } (Cauchy Criterion). Let $%
(x_{n})_{n\geq 0}$ $\mathbb{R}$ be a sequence of (\textbf{real numbers}).\newline

\noindent \textbf{(a)} Show that if $(x_{n})_{n\geq 0}$ is Cauchy,
then it has a unique accumulation point $\ell \in 
\mathbb{R}$ which is its limit.\newline

\noindent \textbf{(b)} Show that if a sequence $(x_{n})_{n\geq 0}\subset 
\mathbb{R}$ \ converges to $\ell \in \mathbb{R},$ then, it is Cauchy.%
\newline

\noindent \textbf{(c)} Deduce the Cauchy criterion for sequences of real numbers.

\newpage

\begin{center}
\textbf{SOLUTIONS}
\end{center}

\noindent \textbf{Exercise 1}.\newline

\noindent \textbf{Question (1)}. It is obvious that :%
\begin{equation*}
\underset{p\geq n}{\inf }x_{p}\leq x_{n}\leq \underset{p\geq n}{\sup }x_{p},
\end{equation*}

\noindent since $x_{n}$ is an element of $\left\{
x_{n},x_{n+1},...\right\} $ on which we take the supremum or the infimum.%
\newline

\noindent \textbf{Question (2)}. Let $y_{n}=\underset{p\geq 0}{\inf }%
x_{p}=\underset{p\geq n}{\inf }A_{n},$ where $A_{n}=\left\{
x_{n},x_{n+1},...\right\} $ is a non-increasing sequence of sets : $\forall n\geq 0$,
\begin{equation*}
A_{n+1}\subset A_{n}.
\end{equation*}

\noindent So the infimum on $A_{n}$ increases. If $y_{n}$ increases in $%
\overline{\mathbb{R}},$ its limit is its upper bound, finite or infinite. So%
\begin{equation*}
y_{n}\nearrow \underline{\lim }\text{ }x_{n},
\end{equation*}

\noindent is a finite or infinite number.\newline

\noindent \textbf{Question (3)}. We also show that $z_{n}=\sup A_{n}$ decreases and $z_{n}\downarrow \overline{\lim }$ $x_{n}$.\newline

\noindent \textbf{Question (4) \label{qst4}}. We recall that 
\begin{equation*}
-\sup \left\{ x,x\in A\right\} =\inf \left\{ -x,x\in A\right\}, 
\end{equation*}

\noindent which we write 
\begin{equation*}
-\sup A=\inf (-A).
\end{equation*}

\noindent Thus,

\begin{equation*}
-z_{n}=-\sup A_{n}=\inf (-A_{n}) = \inf \left\{-x_{p},p\geq n\right\}.
\end{equation*}

\noindent The right hand term tends to $-\overline{\lim}\ x_{n}$ and the left hand to $\underline{\lim} (-x_{n})$ and so 

\begin{equation*}
-\overline{\lim}\ x_{n}=\underline{\lim }\ (-x_{n}).
\end{equation*}

\bigskip \noindent Similarly, we show:
\begin{equation*}
-\underline{\lim } \ (x_{n})=\overline{\lim} \ (-x_{n}).
\end{equation*}

\noindent 

\noindent \textbf{Question (5)}. These properties come from the formulas, where $A\subseteq \mathbb{R},B\subseteq \mathbb{R}$ :%
\begin{equation*}
\sup \left\{ x+y,A\subseteq \mathbb{R},B\subseteq \mathbb{R}\right\} \leq
\sup A+\sup B.
\end{equation*}

\noindent In fact : 
\begin{equation*}
\forall x\in \mathbb{R},x\leq \sup A
\end{equation*}

\noindent and
\begin{equation*}
\forall y\in \mathbb{R},y\leq \sup B.
\end{equation*}

\noindent Thus 
\begin{equation*}
x+y\leq \sup A+\sup B,
\end{equation*}

\noindent where 
\begin{equation*}
\underset{x\in A,y\in B}{\sup }x+y\leq \sup A+\sup B.
\end{equation*}%
Similarly,%
\begin{equation*}
\inf (A+B\geq \inf A+\inf B.
\end{equation*}

\noindent In fact :

\begin{equation*}
\forall (x,y)\in A\times B,x\geq \inf A\text{ and }y\geq \inf B.
\end{equation*}

\noindent Thus 
\begin{equation*}
x+y\geq \inf A+\inf B,
\end{equation*}

\noindent and so
\begin{equation*}
\underset{x\in A,y\in B}{\inf }(x+y)\geq \inf A+\inf B
\end{equation*}

\noindent \textbf{Application}.\newline

\begin{equation*}
\underset{p\geq n}{\sup } \ (x_{p}+y_{p})\leq \underset{p\geq n}{\sup } \ x_{p}+\underset{p\geq n}{\sup } \ y_{p}.
\end{equation*}

\noindent All these sequences are non-increasing. By taking the infimum, we obtain the limits superior :

\begin{equation*}
\overline{\lim }\text{ }(x_{n}+y_{n})\leq \overline{\lim }\text{ }x_{n}+%
\overline{\lim }\text{ }x_{n}.
\end{equation*}

\bigskip

\noindent \textbf{Question (6)}. Set

\begin{equation*}
\underline{\lim } \ x_{n}=\overline{\lim } \ x_{n}.
\end{equation*}

\noindent Since : 
\begin{equation*}
\forall x\geq 1,\text{ }y_{n}\leq x_{n}\leq z_{n},
\end{equation*}%

\begin{equation*}
y_{n}\rightarrow \underline{\lim} \ x_{n}
\end{equation*}%

\noindent and 

\begin{equation*}
z_{n}\rightarrow \overline{\lim } \ x_{n},
\end{equation*}

\noindent we apply the Sandwich Theorem to conclude that the limit of $x_{n}$ exists and :

\begin{equation*}
\lim \text{ }x_{n}=\underline{\lim }\text{ }x_{n}=\overline{\lim }\text{ }%
x_{n}.
\end{equation*}

\bigskip 
\noindent \textbf{Exercice 2}.\newline

\noindent \textbf{Question (a).}\\

\noindent Thanks to Question (4) of Exercise 1, it suffices to show this property for one of the limits. Consider the limit superior and the three cases:\\

\noindent \textbf{The case of a finite limit superior} :

\begin{equation*}
\underline{\lim} x_{n}=\ell \text{ finite.}
\end{equation*}

\noindent By definition, 
\begin{equation*}
z_{n}=\underset{p\geq n}{\sup }x_{p}\downarrow \ell .
\end{equation*}

\noindent So: 
\begin{equation*}
\forall \varepsilon >0,\exists (N(\varepsilon )\geq 1),\forall p\geq
N(\varepsilon ),\ell -\varepsilon <x_{p}\leq \ell +\varepsilon .
\end{equation*}

\noindent Take less than that:

\begin{equation*}
\forall \varepsilon >0,\exists n_{\varepsilon }\geq 1:\ell -\varepsilon
<x_{n_{\varepsilon }}\leq \ell +\varepsilon.
\end{equation*}

\noindent We shall construct a sub-sequence converging to $\ell$.\\

\noindent Let $\varepsilon =1:$%
\begin{equation*}
\exists N_{1}:\ell -1<x_{N_{1}}=\underset{p\geq n}{\sup }x_{p}\leq \ell +1.
\end{equation*}

\noindent But if 
\begin{equation}
z_{N_{1}}=\underset{p\geq n}{\sup }x_{p}>\ell -1, \label{cc}
\end{equation}

\noindent there surely exists an $n_{1}\geq N_{1}$ such that%
\begin{equation*}
x_{n_{1}}>\ell -1.
\end{equation*}

\noindent If not, we would have 
\begin{equation*}
( \forall p\geq N_{1},x_{p}\leq \ell -1\ ) \Longrightarrow \sup \left\{
x_{p},p\geq N_{1}\right\} =z_{N_{1}}\geq \ell -1,
\end{equation*}

\noindent which is contradictory with (\ref{cc}). So, there exists $n_{1}\geq N_{1}$ such that
\begin{equation*}
\ell -1<x_{n_{1}}\leq \underset{p\geq N_{1}}{\sup }x_{p}\leq \ell -1.
\end{equation*}

\noindent i.e.

\begin{equation*}
\ell -1<x_{n_{1}}\leq \ell +1.
\end{equation*}

\noindent We move to step $\varepsilon =\frac{1}{2}$ and we consider the sequence%
 $(z_{n})_{n\geq n_{1}}$ whose limit remains $\ell$. So, there exists $N_{2}>n_{1}:$%
\begin{equation*}
\ell -\frac{1}{2}<z_{N_{2}}\leq \ell -\frac{1}{2}.
\end{equation*}

\noindent We deduce like previously that $n_{2}\geq N_{2}$ such that%
\begin{equation*}
\ell -\frac{1}{2}<x_{n_{2}}\leq \ell +\frac{1}{2}
\end{equation*}

\noindent with $n_{2}\geq N_{1}>n_{1}$.\\

\noindent Next, we set $\varepsilon =1/3,$ there will exist $N_{3}>n_{2}$ such that%
\begin{equation*}
\ell -\frac{1}{3}<z_{N_{3}}\leq \ell -\frac{1}{3}
\end{equation*}

\noindent and we could find an $n_{3}\geq N_{3}$ such that%

\begin{equation*}
\ell -\frac{1}{3}<x_{n_{3}}\leq \ell -\frac{1}{3}.
\end{equation*}

\noindent Step by step, we deduce the existence of $%
x_{n_{1}},x_{n_{2}},x_{n_{3}},...,x_{n_{k}},...$ with $n_{1}<n_{2}<n_{3}%
\,<...<n_{k}<n_{k+1}<...$ such that

$$
\forall k\geq 1, \ell -\frac{1}{k}<x_{n_{k}}\leq \ell -\frac{1}{k},
$$

\bigskip \noindent i.e.

\begin{equation*}
\left\vert \ell -x_{n_{k}}\right\vert \leq \frac{1}{k},
\end{equation*}

\bigskip \noindent which will imply: 
\begin{equation*}
x_{n_{k}}\rightarrow \ell 
\end{equation*}

\bigskip \noindent Conclusion : $(x_{n_{k}})_{k\geq 1}$ is very well a subsequence since $n_{k}<n_{k+1}$ for all $k \geq 1$ 
and it converges to $\ell$, which is then an accumulation point.\\

\noindent \textbf{Case of the limit superior equal $+\infty$} : 
$$
\overline{\lim} \text{ } x_{n}=+\infty.
$$

\bigskip  \noindent Since $z_{n}\uparrow +\infty ,$ we have : $\forall k\geq 1,\exists
N_{k}\geq 1,$ 
\begin{equation*}
z_{N_{k}}\geq k+1.
\end{equation*}

\noindent For $k=1$, let $z_{N_{1}}=\underset{p\geq N_{1}}{\inf }%
x_{p}\geq 1+1=2.$ So there exists 
\begin{equation*}
n_{1}\geq N_{1}
\end{equation*}

\noindent such that :
\begin{equation*}
x_{n_{1}}\geq 1.
\end{equation*}

\noindent For $k=2$, consider the sequence $(z_{n})_{n\geq n_{1}+1}.$
We find in the same manner 
\begin{equation*}
n_2 \geq n_{1}+1
\end{equation*}%
\noindent and 
\begin{equation*}
x_{n_{2}}\geq 2.
\end{equation*}

\noindent Step by step, we find for all $k\geq 3$, an $n_{k}\geq n_{k-1}+1$ such that
\begin{equation*}
x_{n_{k}}\geq k,
\end{equation*}

\noindent which leads to $x_{n_{k}}\rightarrow +\infty $ as $k\rightarrow +\infty $.\\

\noindent \textbf{Case of the limit superior equal $-\infty$} : 

$$
\overline{\lim }x_{n}=-\infty.
$$

\bigskip  \noindent This implies : $\forall k\geq 1,\exists N_{k}\geq 1,$ such that%
\begin{equation*}
z_{n_{k}}\leq -k.
\end{equation*}

\noindent For $k=1$, there exists $n_{1}$ such that%
\begin{equation*}
z_{n_{1}}\leq -1.
\end{equation*}
But 
\begin{equation*}
x_{n_{1}}\leq z_{n_{1}}\leq -1.
\end{equation*}

\noindent Let $k=2$. Consider $\left( z_{n}\right) _{n\geq
n_{1}+1}\downarrow -\infty .$ There will exist $n_{2}\geq n_{1}+1:$%
\begin{equation*}
x_{n_{2}}\leq z_{n_{2}}\leq -2
\end{equation*}

\noindent Step by step, we find $n_{k1}<n_{k+1}$ in such a way that $x_{n_{k}}<-k$ for all $k$ bigger than $1$. So
\begin{equation*}
x_{n_{k}}\rightarrow +\infty 
\end{equation*}

\bigskip

\noindent \textbf{Question (b).}\\

\noindent Let $\ell$ be an accumulation point of $(x_n)_{n \geq 1}$, the limit of one of its sub-sequences $(x_{n_{k}})_{k \geq 1}$. We have

$$
y_{n_{k}}=\inf_{p\geq n_k} \ x_p \leq x_{n_{k}} \leq  \sup_{p\geq n_k} \ x_p=z_{n_{k}}.
$$

\bigskip \noindent The left hand side term is a sub-sequence of $(y_n)$ tending to the limit inferior and the right hand side is a 
sub-sequence of $(z_n)$ tending to the limit superior. So we will have:

$$
\underline{\lim} \ x_{n} \leq \ell \leq \overline{\lim } \ x_{n},
$$

\bigskip
\noindent which shows that $\underline{\lim} \ x_{n}$ is the smallest accumulation point and $\overline{\lim } \ x_{n}$ is the largest.\\

\noindent \textbf{Question (c).} If the sequence $(x_n)_{n \geq 1}$ has a limit $\ell$, it is the limit of all its sub-sequences,
so subsequences tending to the limits superior and inferior. Which answers question (b).\\

\noindent \textbf{Question (d).} We answer this question by combining point (d) of this exercise and Point 6) of the Exercise 1.\\

\noindent \textbf{Exercise 3}. Let $(x_{n})_{n\geq 0}$ be a non-decreasing sequence, we have:%
\begin{equation*}
z_{n}=\underset{p\geq n}{\sup} \ x_{p}=\underset{p\geq 0}{\sup} \ x_{p},\forall
n\geq 0.
\end{equation*}

\noindent Why? Because by increasingness,%
\begin{equation*}
\left\{ x_{p},p\geq 0\right\} =\left\{ x_{p},0\leq p\leq n-1\right\} \cup
\left\{ x_{p},p\geq n\right\}.
\end{equation*}

\bigskip

\noindent Since all the elements of $\left\{ x_{p},0\leq p\leq
n-1\right\} $ are smaller than than those of $\left\{ x_{p},p\geq n\right\} ,$
the supremum is achieved on $\left\{ x_{p},p\geq n\right\} $ and so 
\begin{equation*}
\ell =\underset{p\geq 0}{\sup } \ x_{p}=\underset{p\geq n}{\sup }x_{p}=z_{n}.
\end{equation*}

\noindent Thus
\begin{equation*}
z_{n}=\ell \rightarrow \ell .
\end{equation*}

\noindent We also have $y_n=\inf \left\{ x_{p},0\leq p\leq n\right\}=x_n$, which is a non-decreasing sequence and so converges to
$\ell =\underset{p\geq 0}{\sup } \ x_{p}$. \\

\bigskip

\noindent \textbf{Exercise 4}.\\

\noindent Let $\ell \in \overline{\mathbb{R}}$ having the indicated property. Let $\ell ^{\prime }$ be a given accumulation point.%
\begin{equation*}
 \left( x_{n_{k}}\right)_{k\geq 1} \subseteq \left( x_{n}\right) _{n\geq 0}%
\text{ such that }x_{n_{K}}\rightarrow \ell ^{\prime}.
\end{equation*}

\noindent By hypothesis this sub-sequence $\left( x_{n_{K}}\right) $
has in turn a sub-sub-sequence $\left( x_{n_{\left( k(p)\right) }}\right)_{p\geq 1} $ such that $x_{n_{\left( k(p)\right) }}\rightarrow
\ell $ as $p\rightarrow +\infty $.\newline

\noindent But as a sub-sequence of $\left( x_{n_{\left( k\right)
}}\right) ,$ 
\begin{equation*}
x_{n_{\left( k(\ell )\right) }}\rightarrow \ell ^{\prime }.
\end{equation*}%
Thus
\begin{equation*}
\ell =\ell ^{\prime}.
\end{equation*}

\noindent Applying that to the limit superior and limit inferior, we have:%
\begin{equation*}
\overline{\lim} \ x_{n}=\underline{\lim}\ x_{n}=\ell.
\end{equation*}

\noindent And so $\lim x_{n}$ exists and equals $\ell$.\\

\noindent \textbf{Exercise 5}.\\

\noindent \textbf{Question (a)}. If $\nu _{2k}$ is finite and if $\nu _{2k+1}$ is infinite, then there ate exactly $k$ up-crossings : 
$[x_{\nu_{2j-1}},x_{\nu _{2j}}]$, $j=1,...,k$, that is, we have $D(a,b)=k$.\\

\noindent \textbf{Question (b)}. If $\nu _{2k+1}$ is finite and $\nu _{2k+2}$ is infinite, then there are exactly $k$ up-crossings:
$[x_{\nu_{2j-1}},x_{\nu_{2j}}]$, $j=1,...,k$, that is we have $D(a,b)=k$.\\

\noindent \textbf{Question (c)}. If all the $\nu_{j}'s$ are finite, then there are an infinite number of up-crossings : 
$[x_{\nu_{2j-1}},x_{\nu_{2j}}]$, $j\geq 1k$ : $D(a,b)=+\infty$.\\

\noindent \textbf{Question (d)}. Suppose that there exist $a < b$ rationals such that $D(a,b)=+\infty$. 
Then all the $\nu _{j}'s$ are finite. The subsequence $x_{\nu_{2j-1}}$ is strictly below $a$. 
So its limit inferior is below $a$. This limit inferior is an accumulation point of the sequence $(x_n)_{n\geq 1}$, 
so is more than $\underline{\lim}\ x_{n}$, which is below $a$.\\

\noindent Similarly, the subsequence $x_{\nu_{2j}}$ is strictly below $b$. So the limit superior is above $a$. 
This limit superior is an accumulation point of the sequence $(x_n)_{n\geq 1}$, so it is below $\overline{\lim}\ x_{n}$, 
which is directly above $b$. This leads to :

$$
\underline{\lim}\ x_{n} \leq a < b \leq \overline{\lim}\ x_{n}. 
$$

\bigskip
\noindent That implies that the limit of $(x_n)$ does not exist. In contrary, we just proved that the limit of $(x_n)$ exists, 
meanwhile for all the real numbers $a$ and $b$ such that $a<b$, $D(a,b)$ is finite.\\

\noindent Now, suppose that the limit of $(x_n)$ does not exist. Then,

$$
\underline{\lim}\ x_{n} < \overline{\lim}\ x_{n}. 
$$

\bigskip
\noindent We can then find two rationals $a$ and $b$ such that $a<b$ and a number $\epsilon$ such that $0<\epsilon$, such that 

$$
\underline{\lim}\ x_{n} < a-\epsilon < a < b < b+\epsilon <  \overline{\lim}\ x_{n}. 
$$

\bigskip
\noindent If $\underline{\lim}\ x_{n} < a-\epsilon$, we can return to Question \textbf{(a)} of Exercise 2 and construct a sub-sequence of $(x_n)$
which tends to $\underline{\lim}\ x_{n}$ while remaining below $a-\epsilon$. Similarly, if $b+\epsilon < \overline{\lim}\ x_{n}$, 
we can create a sub-sequence of $(x_n)$ which tends to $\overline{\lim}\ x_{n}$ while staying above $b+\epsilon$. 
It is evident with these two sequences that we could define with these two sequences all $\nu_{j}$ finite and so $D(a,b)=+\infty$.\\

\noindent We have just shown by contradiction that if all the $D(a,b)$ are finite for all rationals $a$ and $b$ such that $a<b$, 
then, the limit of $(x_n)_{n\geq 0}$ exists.\\

\noindent \textbf{Exercise 5}. Cauchy criterion in $\mathbb{R}$.\\

\noindent Suppose that the sequence is Cauchy, $i.e.$,
$$
\lim_{(p,q)\rightarrow (+\infty,+\infty)} \ (x_p-x_q)=0.
$$

\bigskip
\noindent Then let $x_{n_{k,1}}$ and $x_{n_{k,2}}$ be two sub-sequences converging respectively to $\ell_1=\underline{\lim}\ x_{n}$ and $\ell_2=\overline{\lim}\ x_{n}$. So

$$
\lim_{(p,q)\rightarrow (+\infty,+\infty)} \ (x_{n_{p,1}}-x_{n_{q,2}})=0.
$$

\noindent, By first letting $p\rightarrow +\infty$, we have

$$
\lim_{q\rightarrow +\infty} \ \ell_1-x_{n_{q,2}}=0,
$$

\bigskip
\noindent which shows that $\ell_1$ is finite, else $\ell_1-x_{n_{q,2}}$ would remain infinite and would not tend to $0$. 
By interchanging the roles of $p$ and $q$, we also have that $\ell_2$ is finite.\\

\noindent Finally, by letting $q\rightarrow +\infty$, in the last equation, we obtain
$$
\ell_1=\underline{\lim}\ x_{n}=\overline{\lim}\ x_{n}=\ell_2.
$$

\bigskip
\noindent which proves the existence of the finite limit of the sequence $(x_n)$.\\

\noindent Now suppose that the finite limit $\ell$ of $(x_n)$ exists. Then

$$
\lim_{(p,q)\rightarrow (+\infty,+\infty)} \ (x_p-x_q)=\ell-\ell=0,
$$0

\bigskip
\noindent which shows that the sequence is Cauchy.\\

\bigskip
\noindent \textbf{Improper Riemann integral of an odd function on $\mathbb{R}$}. \label{proba_02_oddfunction} Consider

$$
f(x)=\frac{1}{\pi(1+x^2)}, \ x \in \mathbb{R}.
$$

\bigskip \noindent We have

$$
\int_{-\infty}^{+\infty} \frac{1}{\pi(1+x^2)} \ dx = \int_{-\infty}{+\infty} d(\tan x)= \biggr[\tan x\biggr]_{-\infty}^{+\infty}=\pi.
$$

\bigskip \noindent Hence $f$ is a \textit{pdf}. Let $X$ be a random variable associated to the \textit{pdf} $f$. Set $g(x)=xf(x)$, $x\in \mathbb{R}$. Since $g^+$ and $g^-$ are non-negative and locally bounded and Riemann integrable, we have

$$
\int_{\mathbb{R}} g^+(x) d\lambda(x)=\int_{\mathbb{R}} g^-(x) d=+\infty \text{ and } \int_{\mathbb{R}} g^+(x) d\lambda(x)=\int_{\mathbb{R}} g^+(x) d=+\infty,
$$

\bigskip \noindent by using for example the D'Alembert criterion. Hence $\mathbb{E}(X)$ does not exist.\\

\section{Important Lemmas when dealing with limits on limits in $\overline{\mathbb{R}}$} \label{proba_02_appendix_02}

\noindent \textbf{(1) - Cesaro generalized Limit}.\\

\noindent The following result is often quoted as the Cesaro lemma.\\

\begin{lemma} \label{fourrier_sec_01_lem02} Let $(x_n)_{n\geq 1} \subset \mathbb{R}$ be a sequence of finite real numbers converging to $x \in \mathbb{R}$, then sequence of arithmetic means

$$
y_n=\frac{x_1+...+x_n}{n}, \ n\geq 1
$$

\noindent S also converges to $x$.
\end{lemma}

\bigskip \noindent  \textbf{Proof}. Suppose that $(x_n)_{n\geq 1} \subset \mathbb{R}$ converge to $x \in \mathbb{R}$ as $n\rightarrow +\infty$. Fix $\varepsilon>0$. Thus, there exists $N\geq 1$ such that
$|x_n-x|<\varepsilon$ for all $n\geq N$. Now, for any $n\geq N$, we have

\begin{eqnarray*}
|y_n-x|&=&\biggr|\frac{x_1+...+x_n}{n} - \frac{y+...+y}{n}\bigg|\\
&=&\bigg|\frac{(x_1-x)+...+(x_n-x)}{n}\biggr|\\
&=&\bigg|\frac{(x_1-x)+...+(x_N-x)}{n}\bigg|+\bigg|\frac{(x_{N+1}-x)+...+(x_n-x)}{n}\biggr|\\
&=& A_N +\bigg|\frac{(x_{N+1}-x)+...+(x_n-x)}{n}\biggr|.
\end{eqnarray*}

\bigskip \noindent  with $A_N=|(x_1-x)+...+(x_N-x)|$, which is constant with $N$. Hence for for any $n\geq N$, we have

\begin{eqnarray*}
|y_n-x|<&\leq& A_N +\frac{1}{n } (|x_{N+1}-x|+...+|x_n-x|)\\
&\leq& \frac{A_N}{n} +\frac{(n-N)}{n } \varepsilon\\
\end{eqnarray*}

\bigskip \noindent  we conclude that, for all $\varepsilon>0$,

$$
\limsup_{n\rightarrow +\infty} |y_n-x|\leq \varepsilon,
$$

\bigskip \noindent  and this achieves the proof. $\square$\\

\noindent \textbf{Remark}. The limit of sequence of arithmetic means $(x_1+...+x_n)/n$, $n\geq 1$, may exists and that of $(x_n)_{n\geq 1}$ does not. In that sense the limit of the arithmetic means, whenever it exists, is called the \textit{Cesaro generalized limit} of the sequence of $(x_n)_{n\geq 1}$.\\

\noindent \textbf{(2) - Toeplitz Lemma}. Let $(a_{n,k})_{(n\geq 1, \ 1\leq k\leq k(n)}$ be an array of real numbers such that \\

\noindent (i) For any fixed $k\geq 1$, $a_{n,k}\rightarrow 0$ as $n\rightarrow +\infty$,\\

\noindent (ii) there exists a finite real number $c$ such that $\sup_{n\geq 1} \sum_{1\leq h \leq k(n)} |a_{n,k}| \leq c$.\\

\noindent Let $(x_n)_{n\geq 1}$ be a sequence of real number and define $y_n=\sum_{1\leq h \leq k(n)} x_k a_{n,k}$ and $b_n=\sup_{1\leq h \leq k(n)} a_{n,k}$, $n\geq 1$. We have the following facts.\\

\noindent (1) If $x_n\rightarrow 0$ as $n\rightarrow +\infty$, then $y_n\rightarrow 0$ as $n\rightarrow +\infty$.

\noindent (2) If $b_n\rightarrow 1$ and $x_n\rightarrow x \in \mathbb{R}$ as $n\rightarrow +\infty$, then  $y_n\rightarrow x$ as $n\rightarrow +\infty$.\\

\noindent (3) Suppose that $k(n)=n$ for all $n\geq 1$. Let $(c_k)_{k\geq 0}$ be sequence such that the sequence $(b_n)_{n\geq 0}=(\sum_{1\leq k\leq n}|c_k|_{n\geq 0}$ is non-decreasing and $b_n\rightarrow \infty$. If $x_n\rightarrow x \in \mathbb{R}$ as $n\rightarrow +\infty$, then
$$
\frac{1}{b_n} \sum_{1\leq k \leq n} c_k x_k\rightarrow x\ \ as \ \ n\rightarrow +\infty. \ \ \Diamond\\
$$

\bigskip \noindent \textbf{Proof}. All the convergence below are meant as $n\rightarrow +\infty$.\\

\noindent \textbf{Proof of (1)}. Since $x_n$ converges to $0$, we can find for any fixed $\eta>0$ a number $k_0=k_0(\eta)>0$ such that for any $k\geq k_0$, $|x_k|\leq \eta/c$ and (by this), we have for any $n\geq 0$

\begin{eqnarray*}
|y_n|&\leq &\max(\sum_{1\leq h \leq k_0} |x_k| |a_{n,k}|, \sum_{1\leq h \leq k_0} |x_k| |a_{n,k}| + \sum_{k_0 \leq h \leq k(n)} |x_k| |a_{n,k}|) \ (L1)\\
&\leq &\sum_{1\leq h \leq k_0} |x_k| |a_{n,k}| + (\eta/c) \sum_{1 \leq h \leq k(n)} |a_{n,k}|\\
&\leq&\sum_{1\leq h \leq k_0} |x_k| |a_{n,k}| + \eta,
\end{eqnarray*}

\bigskip \noindent  in short

\begin{eqnarray*}
|y_n| \leq \sum_{1\leq h \leq k_0} |x_k| |a_{n,k}| + \eta.
\end{eqnarray*}

\bigskip \noindent  In the Line (L1) above, it is not sure that $k(n)$ might exceed $k_0$, so we bound by the first argument of the \textit{max} if $k(n)\leq k_0$. The the last equation, we let $n$ go to infinity to have, for all $\eta>0$,

\begin{eqnarray*}
\limsup_{n\rightarrow +\infty} |y_n| \leq  \eta,
\end{eqnarray*}

\bigskip \noindent  since that finite number of $k_0$ sequences $|x_k| |a_{n,k}|$ (in $n$) converge to zero, which implies that $y_n$ converges to zero.

\noindent \textbf{Proof of 2}. We have
$$
y_n= x \sum_{1\leq h \leq k_0} a_{n,k} + \sum_{1\leq h \leq k_0} a_{n,k}(x_k-x)\\
$$

\noindent which implies

$$
|y_n-x|\leq |x||\sum_{1\leq h \leq k_0} a_{n,k}-1| + c|x_k-x|,
$$

\noindent * which by the assumptions lead to $y_n\rightarrow 0$.\\

\noindent \textbf{Proof of (3)}. By setting $a_{n,k}=c_k/b_{n}$, $1\leq k n$, we inherit the assumption is the former points with $c=1$ and we may conclude by applying Point (2). $\square$\\

\noindent \textbf{(3) - Kronecker Lemma}. If $(b_n)_{n\geq 0}$ is an increasing sequence of positive numbers and $(x_n)_{n\geq 0}$ is a sequence of finite real numbers such that 
$\left( \sum_{1\leq k \leq n} x_k\right)_{n\geq 0}$ converges to a finite real number $s$, then 

$$
\frac{\sum_{1\leq k \leq n} b_k x_k}{b_n} \rightarrow 0 \ as \ n\rightarrow \infty. \ \Diamond
$$

\noindent \textbf{Proof}. Set $b_0=0$, $a_k=b_{k+1}-b_k$, $k\geq 0$, $s_1=0$, $s_{n+1}=x_1+..+s_n$, $n\geq 2$. We have

\begin{eqnarray*}
\frac{\sum_{1\leq k \leq n} b_k x_k}{b_n}&=&\frac{\sum_{1\leq k \leq n} b_k (s_{k+1}-s_k)}{b_n}\\
&=& s_{n+1} - \frac{1}{b_n}\sum_{1\leq k \leq n} b_k s_{k}. \ \ (L2)\\
\end{eqnarray*}

\bigskip \noindent  To see how to get Line, we just have to develop the summation and to make the needed factorizations as in

\begin{eqnarray*}
&&b_n(s_{n+1}-s_{n})+b_{n-1}(s_{n}-s_{n-1})+b_{n-2}(s_{n-1}-s_{n-2})+ \cdots b_{3}(s_{4}-s_{3})+b_{2}(s_{3}-s_{2})+b_{1}(s_{2}-s_{1})\\
&&b_n s_{n+1} - s_{n})(b_n-b_{n-1}) - s_{n-1} (s_{n-1}-b_{n-2})+ \cdots s_{2}(b_{2}-b_{1})+ s_{1} b_1.
\end{eqnarray*}

\noindent From Line (L2), we may apply Point (3) of the Toeplitz's Lemma above, since $b_n=a_1+...+a_n$, to conclude that the expression in Line (L2) converges to zero as $n\rightarrow +\infty$. $\square$

\section{Miscellaneous Results and facts} \label{proba_02_appendixFacts}

\noindent \textbf{A - Technical formulas}.\\ 

\noindent \textbf{A1}. We have for all $t>0$, we have \label{proba_02_appendixFactsA1}

$$
\forall t \in \mathbb{R}_{+}, \ e^{t(1-t)}\leq 1+t \leq e^t. \  \Diamond
$$

\noindent \noindent \textbf{Proof}. Put $g(t)=(1+t)-e^{t}$ $t\geq 0$. It is clear that $g^{\prime}(t)=1-e^{t}$ is non-positive and hence $g$
is non-decreasing on $\mathbb{R}_{+}$, and thus : for any $t \in \mathbb{R}_{+}$, $g(t)\leq g(0)$, which leads to the right-hand. To deal with the left-hand one, we put $g(t)=e^{t(t-1)}-t(t-1))$, $t \in \mathbb{R}$. The first two derivatives of $g$ are

$$
g^{\prime}(t)=(-2t+1)e^{t(t-1)}-1 \ \ and \ \ g^{\prime\prime}(t)=(4t^2-4t-1)e^{t(t-1)}, \  \ t\in \mathbb{R}.
$$  

\bigskip \noindent The zeros of $4t^2-4t-1)e^{t(t-1)}$ are $t_1=(1-\sqrt{2})/2$ and $t_2=(1+\sqrt{2})/2$. Since $t_1\leq 0$, $g^{\prime\prime}$ is negative on $[0,t_2]$ vanishes on $t_2$ and positive on $]t_2, \ +\infty[$. This means that $g^{\prime}$ which vanishes at $0$ and tends to $-1$ at $+\infty$, decreases on $[0,t_2]$, reaches its minimum at value at $t_2$ and increases to $-1$ on $]t_2, \ +\infty[$. So we have proved that $g^{\prime}$ is non-positive on $t \in \mathbb{R}_{+1}$. We conclude that for any $t \in \mathbb{R}_{+}$, $g(t)\leq g(0)$, which is the right-hand member of the inequality.\\

\noindent \textbf{Acknowledgement}. This proof is due to Cherif Mamadou Moctar Traoré, University of Bamako, Mali.\\

 \section{Quick and powerfull algorithms for Gaussian probabilities}

\noindent \textbf{Visual Basic$^{\texttrademark}$ codes to compute $F(z)=\mathbb{P}(\mathcal{N}(0,1)\leq z)$}.\\

\begin{lstlisting}
Function ProbaNormale(z As Double) As Double
Dim a1 As Double, a2 As Double, a3 As Double, A4 As Double
Dim A5 As Double, w As Double, W1 As Double, P0 As Double

a1 = 0.31938153
a2 = -0.356563782
a3 = 1.781477937
A4 = -1.821255978
A5 = 1.330274429

W1 = Abs(z)
w = 1 / (1 + 0.2316419 * W1)
W1 = 0.39894228 * Exp(-0.5 * W1 * W1)
P0 = (a3 + w * (A4 + A5 * w))
P0 = w * (a1 + w * (a2 + w * P0))

P0 = W1 * P0

If z <= 0 Then
        P0 = 1 - P0
End If

ProbaNormale = 1 - P0
End Function
\end{lstlisting}

\bigskip \noindent \textbf{Quantile Function}.\\

\noindent The quantile function or inverse function of $F(z)=\mathbb{P}(\mathcal{N}(0,1)\leq z)$ is computed by :\\

\begin{lstlisting}
Public Function inverseLoiNormal(z As Double) As Double
Dim a1 As Double, a2 As Double, a3 As Double, A4 As Double, A5 As Double
Dim A6 As Double
Dim W1 As Double, w As Double, W2 As Double, Q As Double
 
a1 = 2.515517: a2 = 0.802853: a3 = 0.010328
A4 = 1.432788: A5 = 0.189269: A6 = 0.001308

	If z <= 0 Then
		inverseLoiNormal = -4
		Exit Function
	ElseIf z >= 1 Then
    inverseLoiNormal = 4
    Exit Function
	End If

Q = 0.5 - Abs(z - 0.5)
w = Sqr(-2 * Log(Q))
W1 = a1 + w * (a2 + a3 * w): W2 = 1 + w * (A4 + w * (A5 + A6 * w))
inverseLoiNormal = (w - W1 / W2) * Sgn(z - 0.5)
End Function
\end{lstlisting}